\documentclass[10pt]{amsart}
\usepackage{graphicx}
\usepackage{caption}
\usepackage{subcaption}
\usepackage{amssymb}
\usepackage{epstopdf}
\usepackage{nicefrac}
\usepackage{enumerate}
\usepackage{float}
\usepackage{color}
\usepackage{pst-all}
\usepackage{tabularx}
\usepackage{hyperref}
\newtheorem{thm}{THEOREM}[section]

\newtheorem{cor}[thm]{COROLLARY}
\newtheorem{defn}[thm]{DEFINITION}

\newtheorem{case}[thm]{CASES}
\newtheorem{hyp}[thm]{HYPOTHESIS}
\newtheorem{lemma}[thm]{LEMMA}

\newtheorem{prop}[thm]{PROPOSITION}
\newtheorem{quest}[thm]{QUESTION}
\newtheorem{remark}[thm]{REMARK}

\newcommand{\ds}{\displaystyle}

\newcommand{\F}{{\mathcal F}}
\newcommand{\G}{\Gamma}
\newcommand{\e}{{\epsilon}}
\newcommand{\ve}{{\varepsilon}}
\newcommand{\g}{{\gamma}}
\newcommand{\wmW}{\widehat{\mW}}
\newcommand{\bRt}{{\bf R}_{0}}
\newcommand{\bRa}{{\bf R}_{a}}

\newcommand{\TP}{\bT_{\Phi}}
\newcommand{\bR}{{\bf R}}

\newcommand{\bT}{{\bf T}}

\newcommand{\mK}{{\mathbb K}}
\newcommand{\mN}{{\mathbb N}}
\newcommand{\mR}{{\mathbb R}}
\newcommand{\mS}{{\mathbb S}}
\newcommand{\mT}{{\mathbb T}}
\newcommand{\mZ}{{\mathbb Z}}
\newcommand{\mW}{{\mathbb W}}
\newcommand{\cA}{{\mathcal A}}
\newcommand{\cB}{{\mathcal B}}
\newcommand{\cC}{{\mathcal C}}
\newcommand{\cD}{{\mathcal D}}
\newcommand{\cE}{{\mathcal E}}
\newcommand{\cG}{{\mathcal G}}
\newcommand{\cI}{{\mathcal I}}
\newcommand{\cJ}{{\mathcal J}}
\newcommand{\cK}{{\mathcal K}}
\newcommand{\cL}{{\mathcal L}}
\newcommand{\cM}{{\mathcal M}}

\newcommand{\cO}{{\mathcal O}}
\newcommand{\cP}{{\mathcal P}}
\newcommand{\cR}{{\mathcal R}}
\newcommand{\cS}{{\mathcal S}}
\newcommand{\cT}{{\mathcal T}}
\newcommand{\cU}{{\mathcal U}}
\newcommand{\cV}{{\mathcal V}}
\newcommand{\cW}{{\mathcal W}}
\newcommand{\cX}{{\mathcal X}}
\newcommand{\cZ}{{\mathcal Z}}

\newcommand{\cGK}{\cG_{K}}
\newcommand{\cGM}{\cG_{\fM}}
\newcommand{\wcGK}{\widehat{\cG}_{K}}

\newcommand{\fMR}{\fM_{\bRt}}
\newcommand{\fMC}{\fM_{\fC}}

\newcommand{\whPsi}{{\widehat \Psi}}
\newcommand{\whPhi}{{\widehat \Phi}}
\newcommand{\whE}{{\widehat E}}
\newcommand{\whF}{{\widehat F}}

\newcommand{\fC}{{\mathfrak{C}}}
\newcommand{\fD}{{\mathfrak{D}}}
\newcommand{\fG}{{\mathfrak{G}}}
\newcommand{\fL}{{\mathfrak{L}}}
\newcommand{\fM}{{\mathfrak{M}}}
\newcommand{\fN}{{\mathfrak{N}}}

\newcommand{\fT}{{\mathfrak{T}}}
\newcommand{\fU}{{\mathfrak{U}}}
\newcommand{\fV}{{\mathfrak{V}}}
\newcommand{\fW}{{\mathfrak{W}}}
\newcommand{\fX}{{\mathfrak{X}}}
\newcommand{\fY}{{\mathfrak{Y}}}
\newcommand{\fZ}{{\mathfrak{Z}}}
\newcommand{\oU}{\overline{U}}

\newcommand{\ox}{\overline{x}}
\newcommand{\oy}{\overline{y}}

\newcommand{\ow}{\overline{w}}

\newcommand{\opsi}{\overline{\psi}}
\newcommand{\ophi}{\overline{\phi}}
\newcommand{\ovg}{\overline{\gamma}}
\newcommand{\ovl}{\overline{\lambda}}
\newcommand{\ovp}{\overline{\varphi}}

 \newcommand{\vp}{{\varphi}}

\newcommand{\A}{{\rm Area}}
\newcommand{\interior}{{\rm int}}

\newcommand{\whzeta}{\widehat{\zeta}}
\newcommand{\psg}{{\rm pseudo}{\star}{\rm group}}
\newcommand{\wtPhi}{\widetilde{\Phi}}

\begin{document}

\begin{abstract}
  In this work, we study  the dynamical properties of   Krystyna Kuperberg's aperiodic flows on $3$-manifolds. We introduce the notion of a ``zippered lamination'', and with suitable generic hypotheses, show that the unique minimal set for such a flow is  an invariant zippered lamination. We obtain    a  precise description of the topology and dynamical properties of the minimal set, including   the presence of non-zero entropy-type invariants and chaotic behavior. Moreover, we show that   the minimal set does not have stable shape, yet  satisfies the Mittag-Leffler condition for homology groups. 
 \end{abstract}

\title{The dynamics of generic  Kuperberg flows}


\author{Steven Hurder}	
\address{Steven Hurder, Department of Mathematics, University of Illinois at Chicago, 322 SEO (m/c 249), 851 S. Morgan Street, Chicago, IL 60607-7045}
\email{hurder@uic.edu}
 \thanks{Preprint date: June 19, 2013; Revised October 14, 2015}

\author{Ana Rechtman}
\thanks{AR supported in part by   CONACyT  Postdoctoral Research Fellowship}\address{Ana Rechtman, Institut de Recherche Math\'ematique Avanc\'ee,
Universit\'e de Strasbourg,
7 rue Ren\'e Descartes,
67084 Strasbourg, France}
\email{rechtman@math.unistra.fr}

\date{}

\maketitle

 \vfill
 \eject

    \tableofcontents

  \vfill
 \eject

 \listoffigures

\vfill
\eject

\section{Introduction} \label{sec-intro}

The ``Seifert Conjecture'', as originally formulated in 1950 by Seifert in \cite{Seifert1950},  asked: ``Does   every non-singular vector field on the $3$-sphere $\mS^3$ have a periodic orbit?''  
The partial answers to this question have a long history.
 F.W.~Wilson    constructed in the 1966  work \cite{Wilson1966}   a smooth flow on a  \emph{plug}  with exactly two periodic orbits, which was used to modify a given flow on a $3$-manifold to obtain one with only isolated periodic orbits. 
Paul Schweitzer   showed in the  1974 work \cite{Schweitzer1974}, that for any closed   $3$-manifold $M$, there exists  a non-singular $C^1$-vector field on $M$ without periodic orbits. Schweitzer's result      suggested a modified version of   Seifert's question:  ``Does every  non-vanishing  $C^\infty$  vector field on a closed $3$-manifold   have a periodic orbit?''  
Krystyna Kuperberg  showed in her celebrated 1994 work \cite{Kuperberg1994}, that the smooth Seifert Conjecture is  also false by inventing a construction  of aperiodic plugs which is renowned for its simplicity,  beauty and subtlety.  

\begin{thm}[K. Kuperberg] \label{thm-mainK}
On every closed oriented 3-manifold $M$, there exists  a $C^\infty$  non-vanishing vector field without periodic orbits.
\end{thm}
The goal of this work is to understand the dynamical properties of   such ``Kuperberg flows'', especially the structure of their minimal sets. 
In the exploration of  these properties, we  reveal their beauty and   discover the hidden complexity of the Kuperberg dynamical  systems.

 Let us recall the strategy of the proofs of the results cited above. 
A \emph{plug} is a compact 3-manifold with boundary in $\mR^3$,
equipped with a flow   satisfying additional conditions. The flow in a
plug  is assumed to be parallel to the ``vertical'' part of the boundary, so that it may be inserted in any coordinate chart of a   $3$-manifold $M$ to modify the given flow on $M$, and    changes only those orbits entering and leaving the ``horizontal'' faces of the plug.  
 Another  assumption  on the flow in a plug is that there are orbits, which are said to be \emph{trapped},  which  enter the plug and never exit. 
The closure of such an orbit limits to a compact invariant set
contained entirely within the interior of the plug, thus the 
plug must contain at least one minimal set. In the case of the   Wilson Plug, the two periodic orbits   are the minimal sets.

A plug is said to be \emph{aperiodic} if it contains no closed
orbits. Schweitzer   observed in  \cite{Schweitzer1974}  
that the role of the periodic orbits in a Wilson Plug could be
replaced by   Denjoy minimal sets,  resulting in an aperiodic  plug,
which could then be used to ``open up'' the isolated closed orbits
provided by Wilson's result. The   flow in the Schweitzer Plug is only $C^1$, due to the topology of the minimal set contained in the plug, around which all trapped orbits for the flow must     accumulate.   Harrison  constructed in  \cite{Harrison1988} a modified  ``non-flat''   embedding of the Denjoy continuum into a $3$-ball, which she used to construct  an aperiodic  plug with   a $C^2$-flow.
In contrast,    Handel  showed in \cite{Handel1980}   that if the
trapped orbits of a plug accumulate on a     minimal set whose topological
dimension is one and is the only \emph{invariant} set for the flow in the plug, then the
minimal set is \emph{surface-like}: the flow restricted to the minimal
set is topologically conjugated to the minimal set of a flow on a surface.

Kuperberg's construction in \cite{Kuperberg1994} of   aperiodic  smooth flows on plugs   introduced a fundamental new idea,   that of  ``geometric surgery'' on a modified version of the   Wilson Plug $\mW$, to obtain  the \emph{Kuperberg Plug} $\mK$ as a quotient space,  $\tau \colon \mW \to \mK$. The Wilson vector field $\cW$ on $\mW$  is modified to provide    a smooth   vector field $\cK$ on the quotient. The   flow of $\cK$ is denoted by $\Phi_t$. This is said to be a \emph{Kuperberg  flow} on $\mK$.

The periodic orbits for the Wilson flow on $\mW$ get ``cut-open'' when they are  mapped to $\mK$, and there they become trapped orbits for $\Phi_t$.     The essence  of the novel strategy behind the aperiodic property of $\Phi_t$  is perhaps best described by a quote from the paper by Matsumoto \cite{Matsumoto1995}:
\begin{quotation}
  We therefore must demolish the two closed orbits in the Wilson Plug beforehand. But producing a new plug will take us back to the starting line. The idea of Kuperberg is to \emph{let closed orbits demolish themselves}. We set up a trap within enemy lines and watch them settle their dispute while we take no active part. 
\end{quotation}
The images in $\mK$ of the cut-open periodic orbits from the Wilson flow $\Psi_t$ on $\mW$,  generate two orbits for the Kuperberg flow $\Phi_t$ on $\mK$, which are called  the \emph{special     orbits} for $\Phi_t$. These two special orbits play an absolutely central role in the study of the dynamics of a Kuperberg flow. 

  There followed after Kuperberg's seminal work,  a collection of  three works explaining in further detail the proof of the aperiodicity for the   Kuperberg flow, and investigating its dynamical properties: 
  \begin{itemize}
\item the S{\'e}minaire Bourbaki lecture  \cite{Ghys1995} by \'{E}tienne Ghys; 
\item  the notes  by  Shigenori Matsumoto \cite{Matsumoto1995} in Japanese,     later translated into English; 
\item  the joint paper     \cite{Kuperbergs1996} by Greg Kuperberg and Krystyna Kuperberg. 
\end{itemize} 
It was observed in these works that the special orbits in $\mK$ each limit to the other, 
and that a Kuperberg flow has   a unique minimal  set, which we denote by $\Sigma$.
  The topological and dynamical   properties of the minimal sets for the  Wilson, Schweitzer
and Harrison Plugs are   fundamental aspects  of the constructions of the flows in these plugs. For the Kuperberg flow, the minimal set $\Sigma$ is not specified by the construction, but rather its topological properties   are a consequence of the strong interaction of the special orbits.  We will see that with the proper geometric assumptions  in the construction of $\mK$,   we are able to  obtain a detailed understanding of the  topological and dynamical   properties of the minimal set  $\Sigma$ as a result.

 The   \emph{Radius Inequality},  stated as hypothesis (K8) in Section~\ref{sec-kuperberg},  is a topological property of the insertion maps  used to construct the quotient space $\mK$ from the Wilson Plug $\mW$. It is an absolutely remarkable aspect of Kuperberg's construction, that the Radius Inequality  is   essentially all  that is required to show that the quotient flow $\Phi_t$ is aperiodic. 
Moreover,     the smooth  insertion maps which satisfy hypothesis (K8) admit    many variations in their local behavior near the special orbits for the Wilson flow, with each choice  yielding  an aperiodic  ``Kuperberg flow''.

In order to describe the \emph{generic   hypotheses} that we introduce, we require some notions which  are described more precisely  in Sections~\ref{sec-wilson}, \ref{sec-kuperberg} and \ref{sec-radius}.
  The \emph{modified Wilson Plug} $\mW$, as defined in Section~\ref{sec-wilson}, contains a cylinder set $\cR \subset \mW$ which is invariant under the Wilson flow $\Psi_t$ on $\mW$, and the boundary of $\cR$ consists of the periodic orbits for the flow.  There is a ``notched'' subset $\cR' \subset \cR$, illustrated in Figure~\ref{fig:notches}, which maps to a closed subset $\tau(\cR') \subset \mK$ by the quotient map  $\tau \colon \mW \to \mK$, as illustrated in Figure~\ref{fig:notched8}.  The $\Phi_t$-flow of  $\tau(\cR')$  is   a non-compact, embedded surface,   $\fM_0 \subset \mK$,   with boundary  consisting of the special orbits   in $\mK$.   
  Thus, the closure  $\fM = \overline{\fM_0}$ is a flow invariant, compact connected subset of $\mK$, which contains the closure of the special orbits, hence as observed in   \cite{Ghys1995,Kuperbergs1996,Matsumoto1995},   the minimal set $\Sigma \subset \fM$.  The existence of this compact  subset $\fM$ which is invariant for the Kuperberg flow $\Phi_t$ is a remarkable consequence  of the construction, and is   the key to a deeper understanding of the dynamical properties of the flow $\Phi_t$.

 This work introduces several new concepts and   techniques which are used in  the study of the space $\fM_0$ and its closure $\fM$. The
first is the notion of \emph{propellers}, which are surfaces with boundary,
possibly minus a point at infinity, embedded in $\mW$ so that they wrap around the core cylinder $\cR$ in the Wilson Plug $\mW$, as illustrated in Figure~\ref{fig:propeller}. 
The projections of these surfaces to $\mK$ are assembled   according to the dynamics of the
Kuperberg flow, to yield the embedded surface $\fM_0$ as partially illustrated in Figure~\ref{fig:choufleur}.

The   study of   the topological structure of $\fM_0$  reveals the fundamental role played by the local dynamics of the flow $\Phi_t$ in a small open neighborhood of the core cylinder $\tau(\cR') \subset \mK$, and especially in sufficiently small  open $\e$-neighborhoods of  the  special orbits along the boundary of $\tau(\cR')$. In fact, we show in Proposition~\ref{prop-syndeticgk*} that the return times of a Kuperberg flow to these $\e$-neighborhoods form a syndetic set, so that the global dynamics of $\Phi_t$ is essentially determined by its local dynamics   in  $\e$-neighborhoods of  $\tau(\cR')$. 
The local dynamics of $\Phi_t$  depend on the choices made   constructing the Wilson vector field $\cW$ on $\mW$ and the Kuperberg vector field $\cK$ on $\mK$. We formulate  generic conditions    on the choice  of $\cW$ and the construction of $\cK$,  in order to eliminate pathologies in the dynamics of the flows  for $\cW$ and $\cK$  that may otherwise be possible.

The first type of generic condition is   given by Hypothesis~\ref{hyp-genericW}. In brief, this states      that  the vertical component of the Wilson vector field $\cW$  has  quadratic vanishing along the periodic  orbits.   The anti-symmetry imposed on the modified Wilson vector field $\cW$ forces   the generic case to be  a quadratic vanishing condition.

The second type of generic condition  is  imposed on the insertion maps $\sigma_i$ for $i=1,2$ which are introduced in Section~\ref{sec-kuperberg} as part of  the construction of the space $\mK$. These maps are required to satisfy the Radius Inequality, in order to obtain an aperiodic flow on $\mK$, but the smooth properties of these maps also control many other aspects of the global dynamics of $\Phi_t$. For example, 
the propellers associated to  the flow  $\Phi_t$ which form $\fM_0$   are generated by the images of curves under the   embedding maps $\sigma_i$ and their inverses,   and without assumptions on the geometries of these curves, it seems impossible to control the geometry of the space $\fM_0$. 
The precise statements  of the    ``generic hypotheses'' we impose   requires various preliminary notations, and are formally   given in Definition~\ref{def-generic}. In brief, they   assert that   the insertion maps $\sigma_i$  are  ``uniformly quadratic'' in an open neighborhood of each periodic orbit.

 We say   that a Kuperberg  flow $\Phi_t$   is \emph{generic} if it satisfies the  conditions of Definition~\ref{def-generic},  which includes the quadratic hypotheses on the Wilson flow in Hypothesis~\ref{hyp-genericW}. The goal of this work is then to make a complete investigation of the dynamical properties of   \emph{generic Kuperberg flows}.  
 
 We next discuss the results of this work. First, two natural problems  concerning the topological dynamics of the flow $\Phi_t$ are to identify its wandering set $\fW$ and its non-wandering set $\Omega$, which are  defined in Section~\ref{sec-minimalset}. The minimal set $\Sigma$ is always contained in the non-wandering set, and we show in   Theorem~\ref{thm-wandering} that: 
\begin{thm}\label{thm-wander2}
Let $\Phi_t$ be a Kuperberg flow of $\mK$ which satisfies Hypothesis~\ref{hyp-genericW}. Then the non-wandering set $\Omega$ is a subset of the closure  $\fM$ of the embedded manifold $\fM_0$, and thus the complement of $\fM$ consists of wandering points for the flow.
\end{thm}

With the additional hypothesis that the flow is generic, then Theorem~\ref{thm-density} shows that the minimal set $\Sigma$ equals the non-wandering set $\Omega$, and we have:
 \begin{thm} \label{thm1}
Let $\Phi_t$ be a generic Kuperberg flow of $\mK$,  then   
$\Sigma =   \Omega = \fM$.
\end{thm}

  The papers \cite{Ghys1995} and \cite{Kuperbergs1996} gave examples where the   minimal set $\Sigma$  equals the space $\fM$, and suggested that the inclusion   $\Sigma \subset \fM$ may be an equality   in more generality than the examples they gave. The proof of Theorem~\ref{thm1} is inspired by the examples and related remarks in these   works.

It was   remarked  above that the topological properties  of the minimal set $\Sigma$ are a consequence of the ``strong interaction'' of the special orbits. The identification  $\Sigma = \fM$ for a generic flow $\Phi_t$ allows to make this remark precise, using the structure theory for the space $\fM$ that we develop in this work. In fact,  a notable aspect of this work is the precise    description of $\fM$ that is developed, which is possibly the first detailed description of an exceptional minimal set of topological dimension 2 for a flow.

Consider the submanifold $\fM_0 \subset \mK$ as   a stratified space,
with  the interior being a stratum of topological dimension $2$, and
the boundary curves  (the \emph{special orbits}) being a stratum of
topological dimension 1.  Thus, the closure $\fM$  of $\fM_0$ inherits   a type of stratified structure, where the $2$-dimensional stratum of $\fM$ are the leaves of a   laminated structure on a subset of $\fM$, obtained from the closure in $\mK$ of the interior of $\fM_0$. The $1$-dimensional stratum of $\fM$ can be described as the  boundaries of the ``leaves'' of the ``lamination''  $\fM$.  Difficulties arise from    this description though, as the special orbits for $\Phi_t$ are the two boundary components of the manifold $\fM_0$, which are dense in $\fM$  for a generic flow.  That is, the $1$-stratum of $\fM$ is dense in the $2$-stratum, which is not a normal property for a lamination.   

In fact, we show in Theorem~\ref{thm-zippered} that for a generic Kuperberg flow, the space $\fM$ has a stronger property than these informal observations, in that it satisfies   the conditions in Definition~\ref{def-zl} which gives it a local product structure, which is  analogous to that  of a lamination.
We call this type of structure a    \emph{zippered lamination}, as it is ``sewn together'' along the special boundary orbits.
 \begin{thm} \label{thm-ZL}
For a generic Kuperberg flow,  the  space $\fM$ has the structure of a zippered  lamination with $2$-dimensional  leaves.  
\end{thm}

The   basis of many of our results concerning $\fM$, including the proof of Theorem~\ref{thm-ZL} above, is a detailed structure theory for $\fM_0$ based on the properties of \emph{finite} or \emph{infinite}, \emph{simple} or \emph{double propellers}, arising from  either \emph{boundary} or \emph{interior  notches} in $\fM_0$. The amazing complexity of these aspects of the embedding of $\fM_0$ in $\mK$ are organized using the     level function on $\fM_0$, which is introduced in this work.

The notion of the  \emph{relative level} of two points along the same  orbit of the Kuperberg flow was introduced by Kuperberg in \cite{Kuperberg1994}, and is a key technique for showing that the flow obtained is aperiodic.  It is defined by equations \eqref{def-level+} and \eqref{def-level-} in Section~\ref{sec-radius} below.  Further properties of the concept of relative level  were  developed by Ghys and Matsumoto in \cite{Ghys1995, Matsumoto1995}.    In  Section~\ref{sec-proplevels}, we show that    the relative level function
 along orbits   induces a well-defined level function on $\fM_0$, for which the Reeb cylinder $\tau(\cR')$ is at level $0$. 
The   level function on $\fM_0$ is used to decompose this space into
an infinite    union of propellers $\fM_{\ell}$ for $\ell
\geq 1$, where $\fM_{\ell}$ consists of sets at level $\ell$: for
$\ell=1$ this set is composed by two non-compact propellers attached to $\tau(\cR')$, while
$\fM_\ell$ for $\ell\geq 2$ is formed by  $2^{\ell}$  
families of compact  propellers. 
The level decomposition of $\fM_0$ is used to analyze its geometry and dynamics, and consequently also that of    $\fM$.

One of the original motivations for this work was a question posed by  Krystyna Kuperberg concerning the \emph{topological shape} of the minimal set $\Sigma$ for the flow $\Phi_t$.   
Section~\ref{sec-shape} gives a very brief introduction to shape theory, including the definitions of ``stable shape'' in Definition~\ref{def-stableshape} and for a continua to be ``movable''    in Definition~\ref{def-movable}. It is a simple observation that the subspaces $\fM$ and its dense subset $\fM_0$ have the same topological shape, and thus for a generic flow, we can use the level decomposition of $\fM_0$ to study the shape properties of the minimal set $\Sigma$. 
As a cumulation of   the   results in this work, we show in Section~\ref{sec-shape} the following results:

\begin{thm} \label{thm-stable}
For a generic Kuperberg flow,  the   minimal set $\Sigma$ does not have stable shape.
\end{thm}

 The Mittag-Leffler condition for homology groups, as introduced in Proposition~\ref{prop-MLmove}, is a homology version of the movable condition.  Proposition~\ref{thm-MLhomology} yields the following consequence:
\begin{thm}\label{thm-ML}
For a generic Kuperberg flow,  the   minimal set $\Sigma$ satisfies the Mittag-Leffler condition for homology groups.
\end{thm}

 Theorems~\ref{thm-stable} and \ref{thm-ML}   follow from three fundamental  properties of the minimal set $\Sigma$. First, that for the generic Kuperberg flow, the minimal set $\Sigma$ equals $\fM$, which  is the closure of the non-compact, embedded surface $\fM_0$ obtained from the flow of the Reeb cylinder $\tau(\cR') \subset \mK$. Second, that the topology of a sufficiently fine open neighborhood approximation $\fM \subset U_k$ -- as occurs in a shape approximation of $\fM$ -- has increasing topological complexity as the index $k$ tends to infinity. Third,  for any finite time, the trace of the flow of the   cylinder $\tau(\cR')$ retracts to the cylinder $\tau(\cR')$, which implies that for $k$ sufficiently large, the topological complexity of the approximating spaces $U_{k}$ are contained in an arbitrarily small open neighborhood of a space $\fM_0^1$ introduced here,   from which the Mittag-Leffler condition is shown to follow. 
 
 Thus, our study of the shape properties of $\Sigma$ uses and combines almost all of the results of this paper.
The details of these arguments are often quite subtle and tedious, and begin with the construction in Section~\ref{sec-shape} of   a   decreasing nested sequence of  compact   domains  $\fN_{\ell} \subset \mK$ for $\ell \geq 0$, such that each $\fN_{\ell}$ satisfies $\fM \subset \fN_{\ell}$ and has the homotopy type of a finite $CW$-complex. The spaces $\fN_{\ell}$ are constructed using  the level function on  $\fM_0$ and the observation in Section~\ref{sec-doublepropellers} that   the double propellers introduced there are nested. That is, the double   propellers at level $\ell +1$ are contained in the closures of the interiors of   double propellers at level $\ell$. We show that   the first homology groups of the spaces $\fN_{\ell}$ have ranks which grow without bound  with $\ell$, and yet for an appropriate choice of a subsequence of these spaces, the ranks of the maps induced on homology from the bonding maps have constant rank $3$. We use this to conclude   that the shape of $\Sigma$ is not stable, but does satisfy the Mittag-Leffler condition.

Our final collection of results concern the entropy invariants   and dynamical complexity  of a   Kuperberg flow $\Phi_t$.
The topological entropy of a Kuperberg flow is zero, for as noted by Ghys   in \cite{Ghys1995}, this follows as a consequence of a   result of Katok on $C^2$-flows on $3$-manifolds in \cite{Katok1980}, which  implies that the entropy of an aperiodic smooth flow on a $3$-manifold must be zero.  Katok's proof in \cite{Katok1980}   uses the  \emph{Pesin Theory} for smooth flows, for example as given in \cite{BP2002}, to conclude from  $h_{top}(\Phi_t) > 0$ that $\Phi_t$ must have a periodic orbit.

Note that  Theorem~\ref{thm-wander2} above shows that the non-wandering set of $\Phi_t$ is contained in $\fM$, and hence the flow entropy $h_{top}(\Phi_t) = h_{top}(\Phi_t | \fM)$, where the latter is the entropy for the flow restricted to $\fM$.  This naturally suggests the question, whether some geometric property of   $\fM$ may  directly imply that $h_{top}(\Phi_t | \fM) =0$? 

  It is a standard technique for the study of the dynamics of flows, to choose a section to the flow and study the dynamics of the induced return maps.  There is a natural choice of ``section'' for the Kuperberg flow $\Phi_t$, given by the rectangle  $\bRt \subset \mK$  defined by \eqref{eq-goodsection}  in Section~\ref{sec-pseudogroup}. 
  There is a well-defined ``return map'' $\whPhi$  to $\bRt$  for the flow $\Phi_t$.    However, the vector field $\cK$   is tangent to $\bRt$ along the center horizontal line, which complicates the study of the   dynamics for $\whPhi$. For this reason, we introduce in Section~\ref{sec-pseudogroup} the pseudogroup $\cGK$ generated by the return map $\whPhi$ restricted to   open subsets of $\bRt$ on which the map is continuous.  

  We show in Section~\ref{sec-pseudogroup} that   the global dynamics of the Kuperberg flow is determined by  the actions of the generators for $\cGK$. 
 In particular, the compact set $\fMR = \fM \cap \bRt$ is locally invariant under the return map $\whPhi$, and so is invariant under the local actions of elements of $\cGK$.  The space $\fMR$ is the closure of the set $\fM_0 \cap \bRt$, and it turns out that for the   collection  of special generating maps   $\ds \cGK^{(1)} \subset \cGK$  as defined in \eqref{eq-symmetricgenset}, the action of these generators on $\fM_0 \cap \bRt$ has a systematic description in terms of the level decomposition  for the propellers in $\fM_0$. This is discussed  in Section~\ref{sec-doublepropellers}. 
 
 Even more is true. The action of the set $\ds \cGK^{(1)}$ on the space $\fM_0 \cap \bRt$  has a description in terms of a Cayley graph for the pseudogroup, as described in Section~\ref{sec-normal}.
 Choosing a
 Riemannian metric on
 the Kuperberg plug, induces a  metric on $\fM_0$ for which this space
 is quasi-isometric to a tree $\TP$ with valence at most 4. The
 notation system for the level decomposition of $\fM_0$  also 
 labels the vertices of the embedded tree $\TP \subset \fM_0$. 
The properties of the tree model for the dynamics of $\cGK$ are discussed further in Sections~\ref{sec-normal} and \ref{sec-growth}.

 Section~\ref{sec-entropyflow}   introduces the entropy   associated to  a  \emph{ finite symmetric set of generators} for  a pseudogroup, following the ideas introduced by Ghys, Langevin and Walczak in \cite{GLW1988}.   To obtain entropy invariants of $\Phi_t$, it is necessary to choose a  {finite symmetric generating set} in $\cGK$.  We work with two such choices:
 \begin{itemize}
\item The   collection $\ds \cG_{\whPhi}^{(1)}$   defined in \eqref{eq-gensetforPwhPsi}, with  associated    entropy $\ds h_{GLW}(\cG_{\whPhi}^* | \fMR)$.
\item The   collection $\ds \cGK^{(1)}$  defined in \eqref{eq-symmetricgenset}, with associated   entropy $\ds h_{GLW}(\cGK^* | \fMR)$. 
\end{itemize}
The results in Section~\ref{sec-entropyflow}   yield the implications:
\begin{equation}
h_{top}(\Phi_t | \fM)  > 0 ~ \Longrightarrow ~ h_{GLW}(\cG_{\whPhi}^* | \fMR) > 0 ~ \Longrightarrow ~ h_{GLW}(\cGK^* | \fMR) > 0 ~ .
\end{equation}
The structure theory of $\fM_0$ and   the subexponential growth estimate Corollary~\ref{cor-subexponentialgrowth},  are   used in the proof of 
  Theorem~\ref{thm-entropyvanish} to show that $\ds h_{GLW}(\cGK^* | \fMR) = 0$. 
We thus conclude:
 \begin{thm}\label{thm-entropyvanishintro}
Let $\Phi_t$ be a generic Kuperberg flow, then $h_{top}(\Phi_t | \fM) =0$.
\end{thm}

One of the intriguing aspects of the dynamics of the $\cGK$-action associated to a Kuperberg flow, is the presence of    families of ``horseshoe-like'' structures, as illustrated in Figure~\ref{fig:transform2},  which show that the $\Phi_t$-orbits have a  form of chaotic behavior. However,
the  rate  of contraction for the return maps of the   flow $\Phi_t$ near the special orbits is ``too slow'' for this chaotic behavior to yield positive    entropy, as is seen in the calculations   we make in Section~\ref{sec-entropyflow} for the proof of Theorem~\ref{thm-entropyvanishintro}.  We say that the flow $\Phi_t$ has  ``slow chaos'' near $\fM$. 

We show   in the work    \cite{HR2014a}, that by varying the embeddings $\sigma_i$ for $i=1,2$ so they no longer satisfy the \emph{Radius Inequality}, then the ``slow chaos'' for the Kuperberg flow, becomes   ``rapid chaos'' associated with a hyperbolic attractor   for the perturbed flow.  These observations imply that  the construction of the Kuperberg  flow $\Phi_t$ places it at ``the boundary of hyperbolicity'', in the manner of partially hyperbolic systems \cite{BDV2005}.

Another pseudogroup model for the dynamics of $\Phi_t$ is introduced in Section~\ref{sec-entropylamination}, which is the pseudogroup $\cGM$ acting on the \emph{simple curves} in $\fMR$. This action on curves induces an action on a Cantor set $\fC \subset \fMR$ which is  transverse to the leaves of the lamination   $\fM$. The action of $\cGM$ on $\fC$ can be thought of as the ``essential model'' for the chaotic behavior on the flow $\Phi_t$.

The ``slow chaos'' property of $\Phi_t$ is quantified  by introducing    the \emph{slow lamination entropy} for $\cGM$, denoted by  $h_{GLW}^{\alpha}(\cGM)$  for $0 < \alpha < 1$.   This   invariant for a pseudogroup action,  is the analog of the \emph{slow flow entropy}  introduced in the works of Katok and Thouvenot \cite{KatokThouvenot1997} and Cheng and Li \cite{ChengLi2010}. Our main result in this section is Theorem~\ref{thm-slowentropy}, which relates the growth of orbits for the pseudogroup with the dynamics of the insertion maps used in the construction of $\mK$.
  
\begin{thm} \label{thm-chaos}
Let $\Phi_t$ be a generic Kuperberg flow. If the insertion maps $\sigma_j$  have ``slow growth''   in the sense of  Definition~\ref{def-slowgrowth}, then $h_{GLW}^{1/2}(\cGM) > 0$. \end{thm}

\begin{remark}\label{rmk-slowentropy}
{\rm 
Let $\Phi_t$ be a generic Kuperberg flow  satisfying Definition~\ref{def-slowgrowth}. The works  \cite{deCarvalho1997,DHP2011}    suggests that $h_{GLW}^{1/2}(\cGM) > 0$ implies  the Cantor set $\fC$ has Hausdorff dimension at least $1/2$.
}
\end{remark}

One of the conclusions of the calculations in Sections~\ref{sec-entropyflow} and \ref{sec-entropylamination} of the entropy for the pseudogroups associated to $\Phi_t$,  is that  the quadratic vanishing of the vertical component of the Wilson vector field, as specified  in Hypothesis~\ref{hyp-genericW}, is a key to showing    that $h_{top}(\Phi_t | \fM) = 0$.    Remark~\ref{rmk-slowentropy} is a speculation that it may be possible to show a more direct relationship between the rate of vanishing for the Wilson field $\cW$ along its periodic orbits and the Hausdorff dimension of the closure $\fM$ of $\fM_0$.
 
If  the Kuperberg vector field $\cK$ is constructed using   a \emph{piecewise smooth} Wilson vector field,  with hyperbolic contracting singularities along its periodic orbits, then   Theorem~\ref{thm-hyperbolicentropy} states that $h_{GLW}(\cGM^*) > 0$, as a consequence of the preceding calculations, suitably adapted.  Question~\ref{quest-entropyflow1} then poses the problem of showing that    $h_{top}(\Phi_t | \fM) > 0$ for such flows.

 It was remarked in  \cite{Ghys1995,Kuperberg1994} that   Kuperberg Plugs can also be constructed for which the manifold $\mK$ and its flow $\cK$ are real analytic. Details of this construction are given in the  Ph.D. Thesis \cite{Rechtman2009}  of the second author.  
The authors expect that for real analytic flows, many of  the results of this work remain valid   without the   generic hypotheses in Definition~\ref{def-generic}. 
 
 In general, there are many further questions about the dynamics of flows formed by a ``Kuperberg surgery'', which is the colloquial name for the construction described in Section~\ref{sec-kuperberg}.

 To conclude this introduction to our work, we explain how the paper
 is organized,  and at the same time we
summarize the   properties   of the dynamics of a Kuperberg Plug,
making the distinction between the results which were known for a general Kuperberg flow, and those results obtained in this work.  
 
 Of course,  the fundamental  result  is  Kuperberg's  theorem, which is Theorem~\ref{thm-mainK} above:
 For any choice of modified Wilson flow on $\mW$ as constructed in Section~\ref{sec-wilson}, and any pair of insertion maps $\sigma_i$ 
 which satisfy the Radius Inequality from Section~\ref{sec-kuperberg}, the Kuperberg flow $\Phi_t$ constructed on $\mK$ is aperiodic. Sections~\ref{sec-wilson} to \ref{sec-minimalset} give a self-contained and very detailed proof of this result, which is based on a synthesis of the results of the papers \cite{Kuperberg1994},
\cite{Ghys1995} and \cite{Matsumoto1995}.
 
The   papers by Ghys \cite{Ghys1995} and Matsumoto \cite{Matsumoto1995} include further results on the dynamics of the Kuperberg flow $\Phi_t$.  In particular, they show  that it has a unique minimal set $\Sigma$
  contained in the interior of $\mK$. Also, Matsumoto showed that the
  Kuperberg Plug traps a set with non-empty interior. These results are also discussed and proved in   Sections~\ref{sec-wilson} to \ref{sec-minimalset} of this work, which also establish    notations  and  fundamental techniques required for  the remaining parts of the work.   

The pseudogroup $\cGK$ acting on a rectangle $\bRt\subset \mK$ is introduced in
Section~\ref{sec-pseudogroup}. We choose five maps among the
generators of this pseudogroup that reflect the flow dynamics near the
minimal set $\Sigma$: the pseudogroup orbits are syndetic in the flow orbits of
points in $\Sigma$. The choice of the rectangle is not arbitrary,  as
it takes advantages of the symmetries in the construction of the plugs
$\mW$ and $\mK$. One unavoidable consequence of the anti-symmetry property  of the Kuperberg flow, is that   there must be discontinuities  for the return map of the flow. The nature of these  discontinuities is described in   Section~\ref{sec-pseudogroup},  and
again later in   Section~\ref{sec-entropyflow} where they enter into  the calculation of the entropy for various pseudogroups.

The set $\fM_0$ that is dense in $\fM$ is described in
Sections~\ref{sec-intropropellers} and \ref{sec-proplevels}, giving
the definition of propeller and the decomposition of $\fM_0$ in level
sets. The relation of this decomposition with the pseudogroup $\cGK$
is explained in
Section~\ref{sec-doublepropellers}, with the introduction of the
notion of double propellers. Double propellers play a fundamental role
in the proofs of Theorem~\ref{thm-stable}. 

To complete the description of $\fM_0$ we introduce in Section~\ref{sec-normal}  a graph $\TP\subset
\fM_0$, which is a tree with an additional loop added at the root point $\omega_0 \in \TP$.   The  vertices of $\TP$ are in fact defined by the action of the five special   generators of
$\cGK$ acting on the root point $\omega_0$. We call $\cGK^*$ the set of words obtained by composition of
the special generators. The geometric interpretation suggests an algebraic
decomposition of words in $\cGK^*$ as {\it normal} words: a
composition of two level
monotone words, that is words along which the level is only increasing
or only decreasing. This simplification is described  in
Section~\ref{sec-normal} and is a key tool for the entropy calculations in Sections~\ref{sec-entropyflow} and \ref{sec-entropylamination}. It  is also used to calculate in Section~\ref{sec-growth} the area growth rate of the special leaf $\fM_0$ in $\fM$.

In the approach to
the set $\fM_0$ described in Section~\ref{sec-proplevels} we introduce and discuss  
two types of irregular behaviors  that can arise in the description of the structure of $\fM_0$, and that depend on the choices made in the construction of the Kuperberg
Plug (Remark~\ref{rmk-bubbles}). Each type of  irregularity results  in the existence of ``double propellers''
attached to $\fM_0$, that change the embedding of $\fM_0$ in a
uniformly bounded way, without affecting the dynamical
invariants of $\fM$. One type of irregularity is studied in
Section~\ref{sec-bubbles}. The second type is more subtle and the
generic hypothesis is needed to control it, hence the description and
boundedness of this phenomena is discussed later in the text, in   
Section~\ref{sec-geometry}. Neither of these types of irregularities have been considered previously in the study of the Kuperberg flow. 

Theorem~\ref{thm1} is proved in Section~\ref{sec-generic}.   The   key to  the proof are the  estimates established in Section~\ref{sec-generic} for the orbit behavior of the return map $\whPhi$  near the critical line $r=2$ in $\bRt$.
The proofs of these inequalities  are based on the generic assumptions  on the insertion maps.   It is to be expected that some form of strong regularity hypotheses is required to prove Theorem~\ref{thm1}, as its conclusion is a type of ``Denjoy Theorem'' for a smooth flow on a surface lamination.
 In fact,  \cite[Theorem~3]{Kuperbergs1996} gives an example of a PL flow for which the minimal set is $1$-dimensional, hence not all of $\fM$. It  seems reasonable to conjecture that a smooth example, not satisfying the generic hypotheses, can also be constructed for which the inclusion $\Sigma \subset \fM$ is proper.

Section~\ref{sec-wandering}     uses   the notion of double
propellers to prove that the set $\fM$ is the non-wandering set of the
Kuperberg Plug. This result allows us to restrict the entropy
calculation for the flow 
to $\fM$, as discussed   in Section~\ref{sec-entropyflow}.  The calculations in this section use almost all the results of the study of the geometry and dynamical properties of the flow established in previous sections.  In particular, the calculations show the usefulness of    the $\psg$ $\cGK^*$ defined by the 
  nice set of generators for $\cGK$ mentioned above. The computation of
the entropy that we present uses the results and techniques from
Sections~\ref{sec-pseudogroup}, \ref{sec-normal} and
\ref{sec-wandering}, but avoids the use of Pesin Theory for flows.

The   definition of a zippered lamination is given in     Definition~\ref{def-zl}, and  Theorem~\ref{thm-zippered} shows 
 that for a generic flow, the space $\fM$ is   a zippered lamination. The holonomy pseudogroup $\cGM$ of the zippered lamination   $\fM$ is   introduced in Section~\ref{sec-entropylamination}, and its generators are related to the five special generators of the pseudogroup $\cGK$ induced by the return map of the flow $\Phi_t$ to the section $\bRt$. The notable properties of $\cGM$ is that it collapses the symmetry built into the Kuperberg flow, and also that the holonomy induced by a leafwise path can give a more efficient representation of the action of words in $\cGK$. 
 These two properties  suggest the possibility of the entropy of the lamination
  being positive, even if the flow has zero topological
entropy. In Section~\ref{sec-entropylamination} we show that the entropy of the action of $\cGM$ vanishes as well, but under extra
hypothesis on the insertion maps,  we prove Theorem~\ref{thm-chaos},
implying that some chaotic behavior exist in its orbit structure anyway. In
Section~\ref{sec-growth},  we use normal forms of words in $\cGM$ to
establish the growth type of the leaves of $\fM$.

Finally, Theorems~\ref{thm-stable} and \ref{thm-ML} are proved in Section~\ref{sec-shape},
using the double propellers of Section~\ref{sec-doublepropellers} to
build a suitable sequence of neighborhoods of $\fM$ which are used to study its shape properties.

The reader will quickly observe   one of the significant contributions of this work, which is an extensive collection of illustrations which accompany the text.  Many of the dynamical properties that we discuss here are difficult to grasp without these illustrations, which we hope will aid the reader to a full understanding of the beauty and complexity of the dynamics of Kuperberg flows. The approach in this paper invokes four perspectives on the dynamics of the Kuperberg flow, and the  illustrations help to understand the relationship between the dynamics as analyzed using each of these viewpoints.  The models of the dynamical behavior illustrated in each case  are related by non-linear transformations, so that it often requires some effort to visualize the correspondence between each approach.

This paper owes a profound debt to the authors of the works \cite{Ghys1995,Kuperberg1994,Kuperbergs1996,Matsumoto1995,Rechtman2009} whose text and pictures provided many insights to the Kuperberg dynamics during the development of this work, and inspired many of the illustrations in this paper.

 This work was made possible by the help of many colleagues and
 institutions. We  thank our colleagues Alex Clark,  \'{E}tienne Ghys
 and Krystyna Kuperberg for their insights and suggestions during the
 development of this work. We also thank the funding agency CONACyT of
 Mexico for its postdoctoral Research Fellowship support,     the
 Mathematics Department of UIC, and the Mathematics Department of
 University of Chicago for its welcoming support of the second
 author, and the University of Strasbourg for its welcoming support of the
 first author. We gave numerous talks on this work during the
 preparation of this manuscript, we thank the colleagues that
 contributed with comments and questions. 

 Finally, the authors owe many thanks to the referee, who provided many and repeated helpful comments to improve the presentation, and to improve and correct the proofs of various results in this text through the multiple revisions of this manuscript.

  \vfill
 \eject
 \bigskip
 
\section{The modified Wilson Plug} \label{sec-wilson}

A $3$-dimensional plug is a manifold $P$ endowed with a vector field $\cX$ satisfying the following characteristics: The manifold $P$ is of the form $D \times [-2,2]$, where $D$ is a compact 2-manifold with boundary $\partial D$. Set 
$$\partial_v P = \partial D \times [-2,2] \quad , \quad \partial_h^- P = D \times \{-2\} \quad , \quad \partial_h^+ P = D \times \{2\}$$
Then   the boundary  (with corners) of $P$ has a decomposition
$$\partial P ~ = ~  \partial_v P \cup \partial_h P ~ = ~  \partial_v P \cup \partial_h^-P \cup \partial_h^+ P ~ .$$

Let $\frac{\partial}{\partial z}$ be the \emph{vertical} vector field on $P$, where $z$ is the coordinate on  $[-2,2]$.

The vector field $\cX$ on $P$ must satisfy the conditions:
\begin{itemize}
\item[(P1)] \emph{vertical at the boundary}: $\cX=\frac{\partial}{\partial z}$ in a neighborhood of $\partial P$; thus, $\partial_h^- P$ and $\partial_h^+ P$  are the entry and exit regions of $P$ for the flow of $\cX$, respectively; 
\item[(P2)]   \emph{entry-exit condition}: if a point $(x,-2)$ is in the same trajectory as $(y,2)$, then $x=y$. That is,  an orbit that traverses $P$, exits just in front of its entry point;
\item[(P3)] \emph{trapped orbit}: there is at least one entry point whose entire \emph{forward} orbit is contained in $P$; we will say that its orbit is \emph{trapped} by $P$;
\item[(P4)] \emph{tameness}: there is an embedding $i \colon P\to \mR^3$ that preserves the vertical direction on the boundary $\partial P$.
\end{itemize}

   A plug is \emph{aperiodic} if there is no  closed orbit for $\cX$. 
   
   Note that conditions (P2) and (P3) imply that if the forward orbit of a point $(x,-2)$ is trapped, then the backward orbit of $(x,2)$ is also trapped. 

A {\it semi-plug} is  a manifold $P$ endowed with a vector field $\cX$ as above, satisfying conditions (P1), (P3) and (P4), but not necessarily (P2).  The concatenation of a semi-plug with an inverted copy of itself, that is a copy where the direction of the flow is inverted, is then a plug. Note that we can generalize the above definition to higher dimensions: just take the manifold $D$ to have dimension $n-1$, where $n$ is the dimension of the ambient manifold of the flow.

Condition   (P4)  implies that given any open ball  $B({x},\e) \subset
\mR^3$ with $\e > 0$ and $x$ a point, there exists a modified embedding $i' \colon P \to B({x},\e)$ which preserves the vertical direction again.   Thus, a plug can be used to change a vector field $\cZ$ on any $3$-manifold $M$ inside a flowbox, as follows. Let $\vp \colon U_x \to (-1,1)^3$ be a coordinate chart which maps the vector field $\cZ$ on $M$ to the vertical vector field $\frac{\partial}{\partial z}$. Choose a modified embedding $i' \colon P \to B({x},\e) \subset (-1,1)^3$, and then  replace the flow $\frac{\partial}{\partial z}$ in the interior of $i'(P)$ with the image of $\cX$. This results in a   flow $\cZ'$ on $M$.

 The   entry-exit condition implies that    a  periodic orbit of $\cZ$ which meets $\partial_h P$ in a non-trapped point,   will remain periodic after this modification. An orbit of $\cZ$ which meets $\partial_h P$ in a  trapped point  never exits the plug $P$, hence  after modification,   limits to a closed invariant   set contained in $P$.  A closed invariant set contains a minimal set for the flow, and thus, a plug  serves as a device    to insert a minimal set  into a flow.

We   next introduce the    ``modified Wilson Plug'', which is the first step in the construction   of the Kuperberg Plug.    
Consider the rectangle 
$${\bR} = [1,3]\times[-2,2] = \{(r,z) \mid 1 \leq r \leq 3 ~ \& -2 \leq z \leq 2\} .$$ 
Choose a $C^\infty$-function $g \colon \bR \to [0,g_0]$ for $g_0 > 0$,  which satisfies the ``vertical'' symmetry condition $g(r,z) = g(r,-z)$. 
The value of $g_0 > 0$ is arbitrary; to be definite, we fix $g_0 = 1$ throughout this work. 
Also, require that   $g(2,-1) = g(2,1) = 0$,   that $g(r,z) = 1$ for $(r,z)$ near the boundary of $\bR$, and that $g(r,z) > 0$ otherwise.  
Later, in \eqref{eq-generic1} we will specify that $g(r,z) = 1$ for all points outside of an $\e_0$-neighborhood of the vanishing points $(2,-1)$ and $(2,1)$.

\begin{figure}[!htbp]
\centering
 \includegraphics[width=40mm]{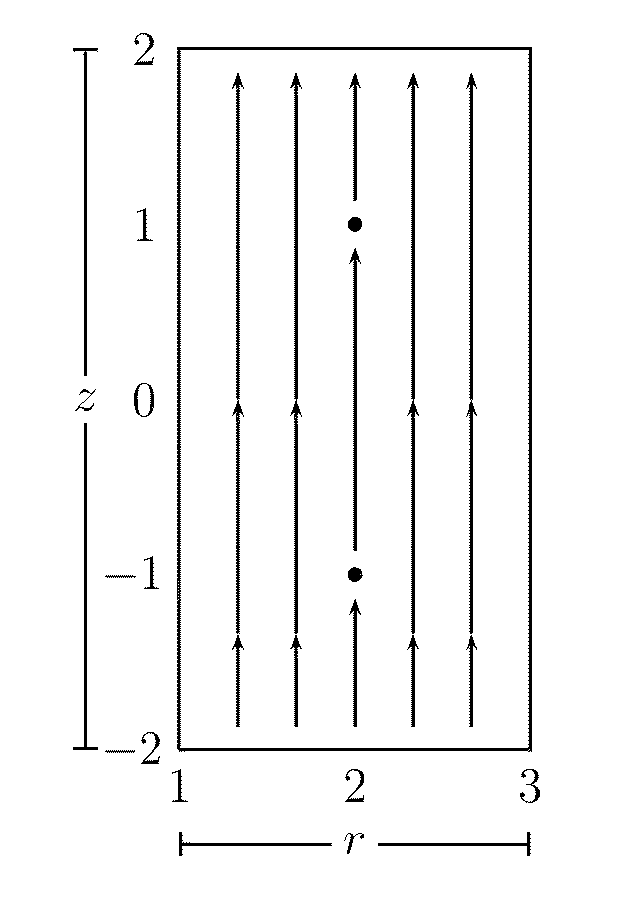}
\caption{Vector field $\cW_v$ \label{fig:wilson1}}
\vspace{-10pt}
\end{figure}

Define the vector field $\cW_v = g \cdot \frac{\partial}{\partial  z}$
which has two singularities, $(2,\pm 1)$ and   is otherwise
everywhere vertical, as  illustrated in Figure~\ref{fig:wilson1}. 

Next, choose a $C^\infty$-function $f \colon \bR \to [-1,1]$ which satisfies the following conditions:
\begin{enumerate}
\item[(W1)] $f(r,-z) = -f(r, z)$ ~ [\emph{anti-symmetry in z}].
\item[(W2)] $f(\xi) = 0$ for $\xi$ near the boundary of $\bR$.
\item[(W3)] $f(r,z) \geq 0$ for  $-2 \leq z \leq 0$, and  $f(r,z) > 0$ for $5/4 \leq r \leq 11/4$ ~ and ~  $-7/4 \leq z < 0$. 
\item[(W4)] $f(r,z) \leq 0$ for  $0  \leq z \leq 2$, and $f(r,z) < 0$ for $5/4 \leq r \leq 11/4$ ~ and ~  $0  < z \leq 7/4$.
\item[(W5)] $f(r,z) =1$ for $5/4 \leq r \leq 11/4$  ~ \text{and} ~  $-7/4 \leq z \leq -1/4$. 
\item[(W6)] $f(r,z) = -1$ for $5/4 \leq r \leq 11/4$ ~ \text{and}~  $1/4 \leq z \leq 7/4$. 
\end{enumerate}
Condition (W1) implies that $f(r,0) =0$ for all $1 \leq r \leq 3$,  and that  
  Conditions (W5) and (W6) are equivalent. Note that Conditions (W5) and (W6)  are stated more precisely than in the works \cite{Kuperberg1994,Ghys1995,Matsumoto1995}, as it is convenient  to specify the values of $f$ on the specified domain in later considerations.

  Next, define the manifold with boundary
\begin{equation}\label{eq-wilsoncylinder}
\mW=[1,3] \times \mS^1\times[-2,2] \cong {\mathbf R} \times \mS^1
\end{equation}
with cylindrical coordinates  $x = (r, \theta,z)$.  That is, 
 $\mW$ is a      solid cylinder with an open core removed, obtained by rotating  the rectangle $\bR$,  considered as embedded in $\mR^3$, around the $z$-axis. 
 
\begin{figure}[!htbp]
\centering
\begin{subfigure}[c]{0.3\textwidth}{\includegraphics[width=35mm]{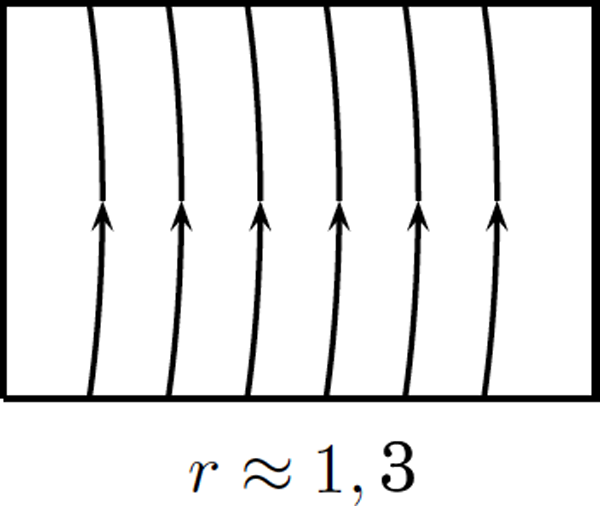}}\end{subfigure}
\begin{subfigure}[c]{0.3\textwidth}{\includegraphics[width=35mm]{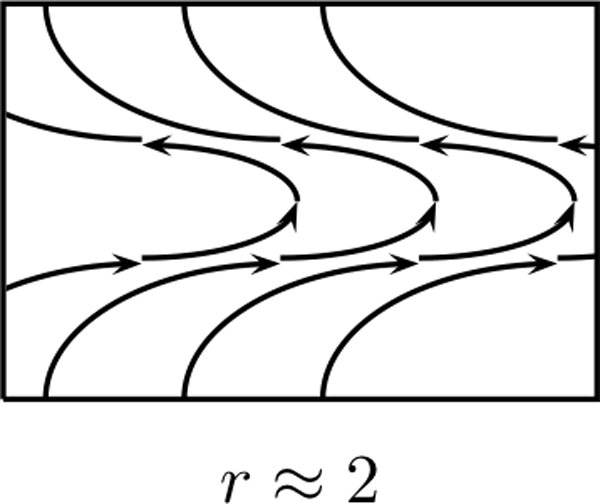}}\end{subfigure}
\begin{subfigure}[c]{0.3\textwidth}{\includegraphics[width=35mm]{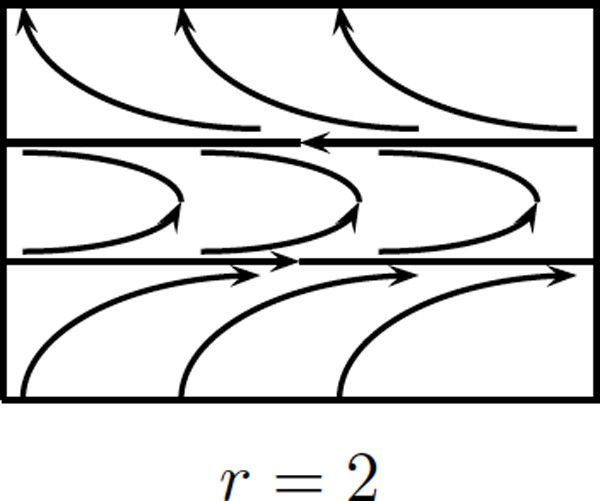}}\end{subfigure}
\caption{$\cW$-orbits on the cylinders $\{r=const.\}$ \label{fig:flujocilin}}
\vspace{-6pt}
\end{figure}

\begin{figure}[!htbp]
\centering
\begin{subfigure}[c]{0.4\textwidth}{\includegraphics[height=60mm]{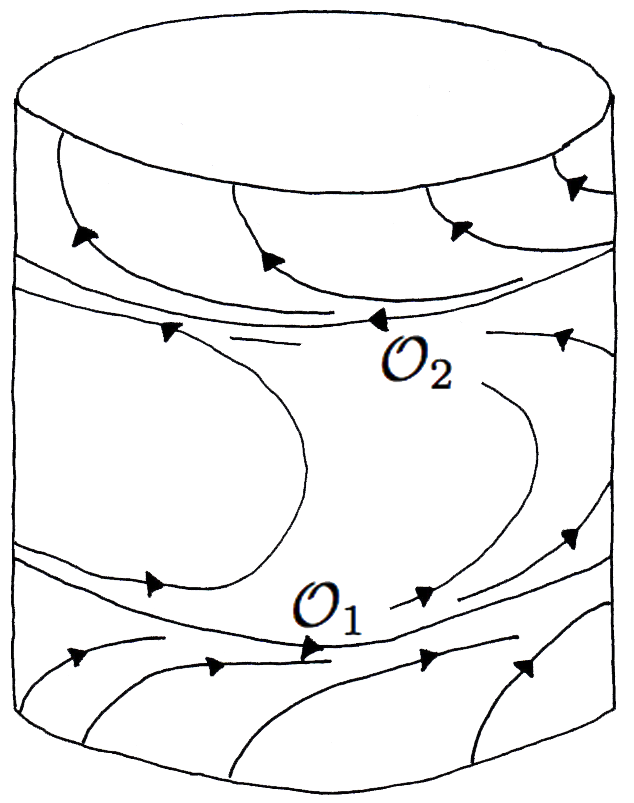}}\end{subfigure}
\begin{subfigure}[c]{0.4\textwidth}{\includegraphics[height=60mm]{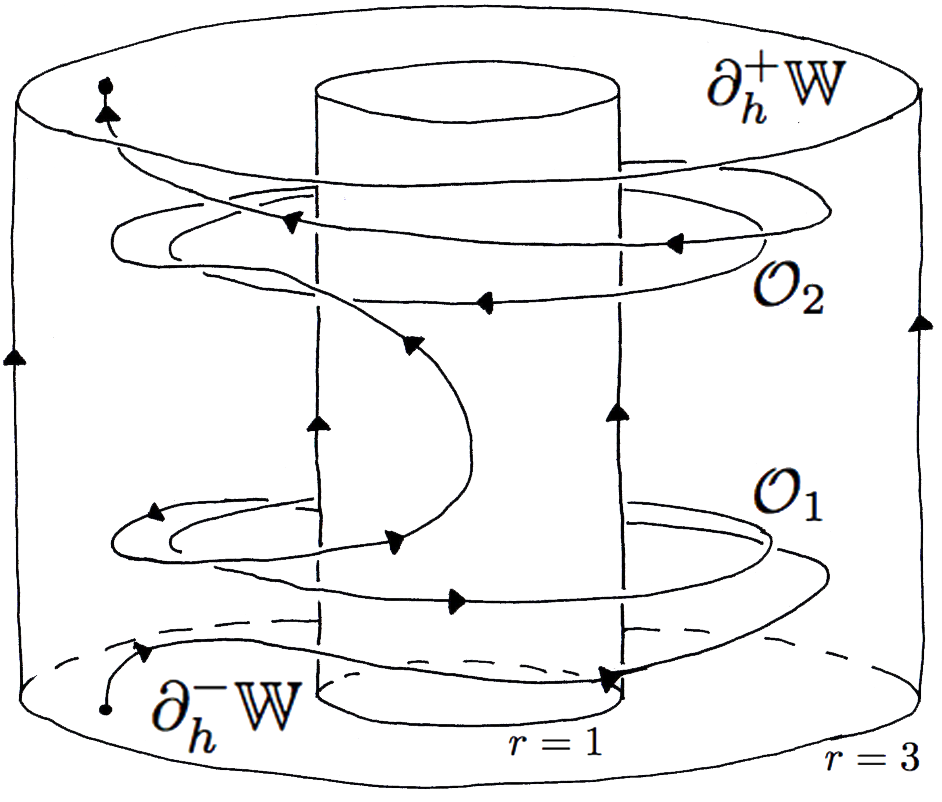}}\end{subfigure}
 \caption[{$\cW$-orbits in the cylinder  $\cC=\{r=2\}$}]{$\cW$-orbits in the cylinder  $\cC=\{r=2\}$ and in $\mW$ \label{fig:Reebcyl}}  
\vspace{-6pt}
\end{figure}

Extend the functions $f$ and $g$ above to $\mW$ by setting $f(r, \theta, z) = f(r, z)$ and $g(r, \theta, z) = g(r, z)$, so that they are   invariant under rotations around the $z$-axis. 
The \emph{modified Wilson vector field} $\cW$ on $\mW$  is defined by 
\begin{equation}\label{eq-wilsonvector}
\cW =g(r, \theta, z)  \frac{\partial}{\partial  z} + f(r, \theta, z)  \frac{\partial}{\partial  \theta} ~ .
\end{equation}
Let $\Psi_t$ denote the flow of $\cW$ on $\mW$.  
Observe that the vector field $\cW$ is vertical near the boundary of $\mW$,  and is horizontal for the points $(r, \theta, z) = (2,\theta, \pm 1)$. 
Also, $\cW$ is tangent to the cylinders $\{r=const.\}$. The flow $\Psi_t$ on the cylinders $\{r=const.\}$ is illustrated (in cylindrical coordinate slices) by   Figures~\ref{fig:flujocilin} and \ref{fig:Reebcyl}.

We give some of the basic properties of the Wilson flow.  
Let $R_{\varphi} \colon \mW \to \mW$ be rotation by the angle $\varphi$. That is, $R_{\varphi}(r,\theta,z) = (r, \theta + \varphi, z)$. 
Define the closed subsets:
\begin{eqnarray*}
\cC ~ & \equiv & ~  \{r=2\}   \quad \text{[\emph{The ~ Full ~ Cylinder}]}\\
\cR ~ & \equiv & ~  \{(2,\theta,z) \mid  -1 \leq z \leq 1\} \quad \text{[\emph{The ~ Reeb ~ Cylinder}]}\\
\cA ~ & \equiv & ~  \{z=0\} \quad \text{[\emph{The ~ Center ~ Annulus}]}\\
\cO_i ~ & \equiv & ~  \{(2,\theta,(-1)^i)  \} \quad \text{[\emph{Periodic Orbits, i=1,2}]}
\end{eqnarray*}
Then $\cO_1$ is the lower   boundary circle  of the Reeb cylinder $\cR$, and $\cO_2$ is the upper boundary circle. 

\begin{prop}\label{prop-wilsonproperties}
Let $\Psi_t$ be the flow on $\mW$ defined above, then:
\begin{enumerate}
\item $R_{\varphi} \circ \Psi_t = \Psi_t \circ R_{\varphi}$ for all $\varphi$ and $t$.
\item The flow $\Psi_t$ preserves the cylinders $\{r=const.\}$   and  in particular 
    preserves the cylinders $\cR$ and $\cC$. 
\item $\cO_i$ for $i=1,2$ are  the periodic orbits for $\Psi_t$.
\item For $x = (2,\theta,-2)$, the forward orbit   $\Psi_t(x)$ for $t > 0$ is trapped.
\item For $x = (2,\theta,2)$, the backward orbit  $\Psi_t(x)$ for $t < 0$ is trapped.
\item For $x = (r,\theta,z)$ with $r \ne 2$, the orbit $\Psi_t(x)$ terminates in the top face $\partial_h^+ \mW$ for some $t \geq 0$, and terminates in $\partial_h^- \mW$ for some $t \leq 0$.
\item The flow $\Psi_t$  satisfies   the entry-exit condition (P2) for plugs.
\end{enumerate}
\end{prop}
\proof The only assertion that needs a comment is the last, which
follows by (W1) and the symmetry condition imposed on the functions
$g$ and $f$.
\endproof
 
Observe that for the choice of $g_0 =1$ for the maximum value of the
function $g$, the typical orbit of $\cW$ rises from $z=-2$ to $z=2$ in
less than one revolution around the circle parameter $\theta$. A
smaller choice of $g_0$ much closer to $0$ will result in a flow which
climbs much slower, and so will possibly make many more revolutions
before transiting the cylinder. This observation will be used in the discussion of the pseudogroup dynamics in Section~\ref{sec-pseudogroup}, and also in the
discussion of   propellers in Sections~\ref{sec-intropropellers}, \ref{sec-doublepropellers}.

 \vfill
 \eject
 
\section{The Kuperberg plug}\label{sec-kuperberg}

\begin{figure}[!htbp]
\centering
{\includegraphics[width=80mm]{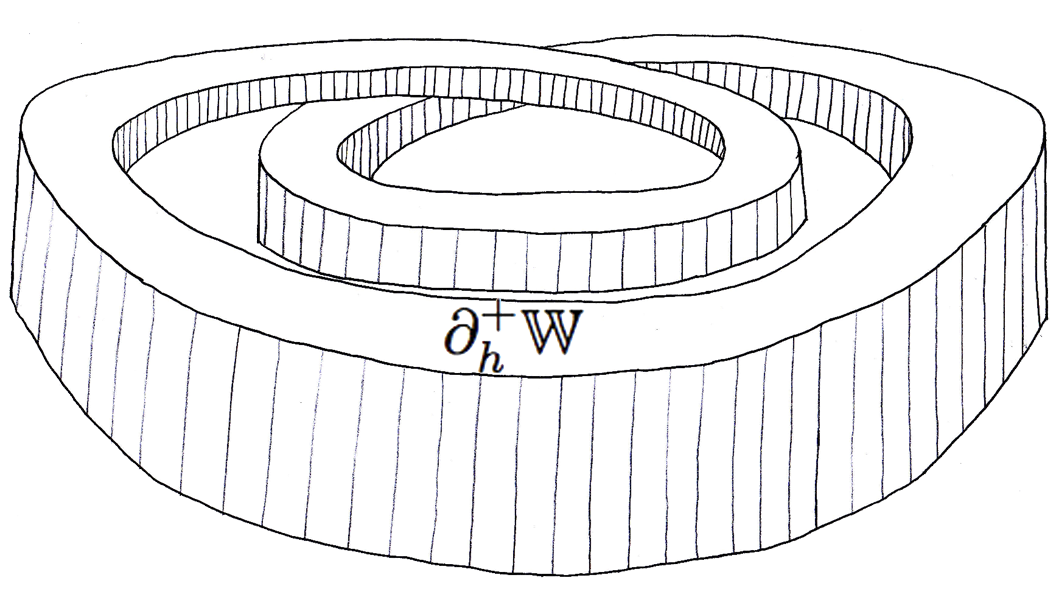}}
\caption{Embedding of Wilson Plug $\mW$ as a {\it folded figure-eight} \label{fig:8doblado} }
\vspace{-6pt}
\end{figure}

The construction of the Kuperberg Plug begins with    the modified Wilson Plug   $\mW$ with vector field $\cW$. 
The first step is to re-embed the manifold $\mW$   in $\mR^3$ as a {\it folded
 figure-eight}, as shown in Figure~\ref{fig:8doblado}, preserving the
vertical direction. A simple but basic point is that the embeddings of the faces of the plug are not   ``planar''.

The fundamental   idea of the Kuperberg Plug is to   construct  two  insertions  of $\mW$ into itself,  in
such a way that the two periodic orbits will be trapped by these
self-insertions. Moreover, the insertions are made so that the resulting space $\mK$ is again embedded in $\mR^3$. A key subtlety of the construction arises in the precise requirements on these self-insertions. As with the construction of the modified Wilson Plug, the description of this construction in   the works \cite{Kuperberg1994,Ghys1995,Matsumoto1995} is qualitative, as this suffices to prove the aperiodicity of the resulting flow. As will be seen in later sections of this work, 
other  properties of the dynamics of the flow $\Phi_t$ in the resulting plug $\mK$ are strongly influenced by the precise nature of these maps, so   we specify the definitions and some properties of the insertion maps more carefully in this section. Later in the manuscript, we formulate       additional ``generic''  requirements on the insertion maps,      in order to obtain  further properties about the dynamics of the flow $\Phi_t$.

\begin{figure}[!htbp]
\centering
{\includegraphics[width=60mm]{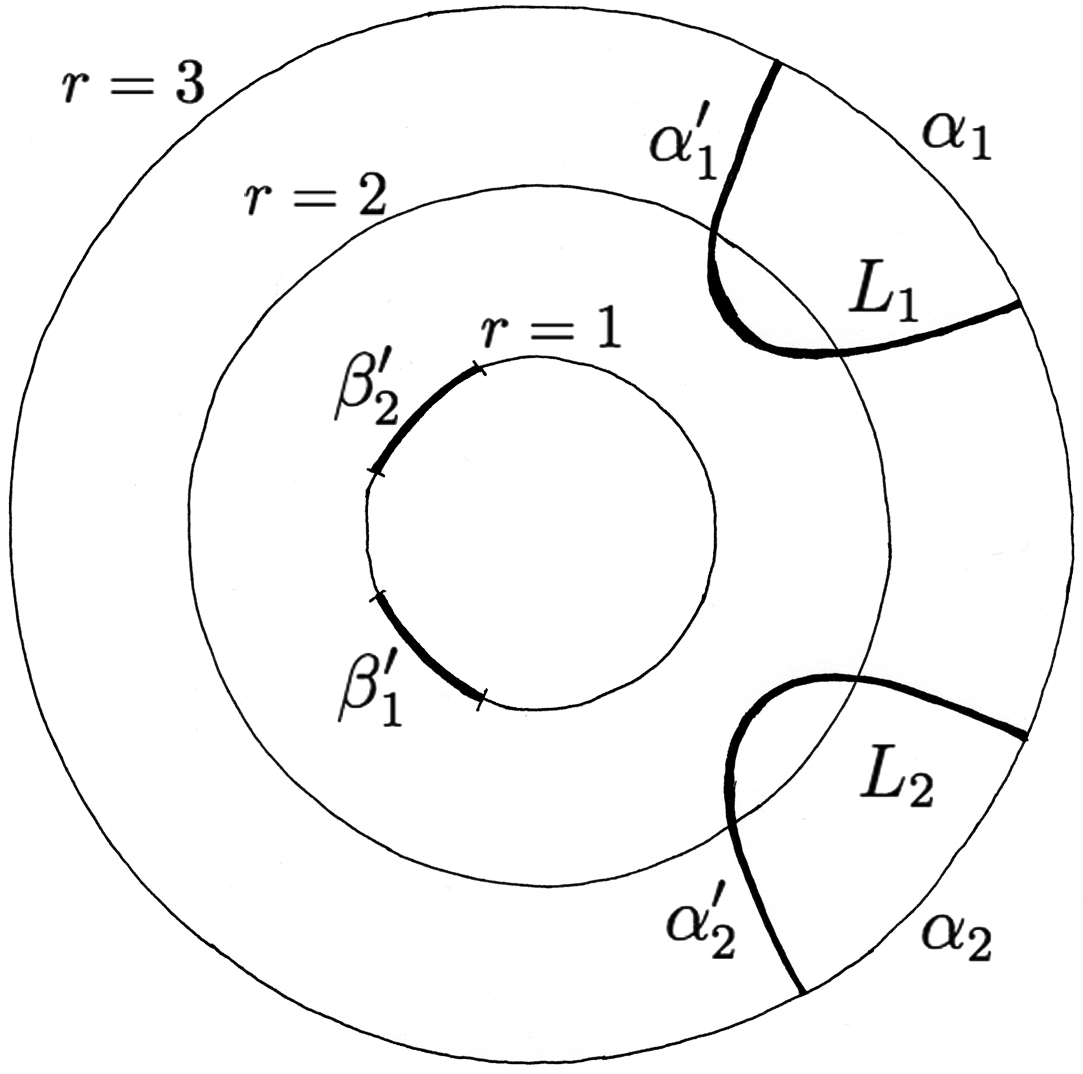}}
\caption{ The disks $L_1$ and $L_2$ \label{fig:insertiondisks}}
\vspace{-6pt}
\end{figure}

The construction begins with the choice    in the annulus $[1,3] \times \mS^1$  of  two   closed regions $L_i$, for $i=1,2$, which are topological disks. Each region has boundary defined by   two arcs: for $i=1,2$, $\alpha^\prime_i$ is the boundary contained in  the
interior of $[1,3] \times \mS^1$  and $\alpha_i$ in the outer boundary contained in the circle
$\{r=3\}$, as depicted in Figure~\ref{fig:insertiondisks}. 

We fix these curves precisely. Let $\zeta_1 = \pi/4$ and $\zeta_2 = - \pi/4$, then  define the arcs
$$
\alpha_1   ~ = ~ \{(3, \theta) \mid   ~  |\theta - \zeta_1| \leq 1/10\} \quad , \quad 
\alpha_2 ~ = ~ \{(3, \theta) \mid   ~ |\theta - \zeta_2| \leq 1/10\} ~ . 
$$
Let $\alpha_i'$ be the curves in the interior of $L_i$ which in polar coordinates $(r,\theta)$ are parabolas with minimum values $r = 3/2$ and base the line segment $\alpha_i$, as depicted in Figure~\ref{fig:insertiondisks}. We choose an explicit form for the embedded curves, for example,     given by
$\ds \alpha_i' \equiv \{ r =  3/2 + 300/2 \cdot (\theta - \zeta_i)^2  \}$.

Consider the closed sets $D_i \equiv L_i \times[-2,2] \subset \mW$, for $i = 1,2$.  
Note that  each $D_i$ is homeomorphic to a closed $3$-ball, that $D_1 \cap D_2 = \emptyset$  and each $D_i$ intersects the  cylinder $\{r=2\}$ in   a rectangle. 
Label   the top and bottom faces of these regions by
\begin{equation}\label{eq-regions}
L_1^{\pm}=  L_1\times \{\pm 2\}  ~ , ~  
L_2^{\pm} =    L_2\times \{\pm 2\} ~.
\end{equation}

The next step is to define insertion maps   $\sigma_i \colon D_i \to \mW $, for $i=1,2$, in such
a way that the periodic orbits $\cO_{1}$ and $\cO_{2}$ for the $\Psi_t$-flow  intersect $\sigma_i(L_i^-)$ in   points
corresponding to   $\cW$-trapped entry points for the modified Wilson Plug. 
Consider the two disjoint arcs $\beta_i'$ in the inner boundary circle $\{r=1\}$,
\begin{eqnarray*}
\beta_1' & = & \{(1, \theta) \mid   ~  |\theta - (\zeta_1 + \pi) | \leq 1/10\}\\
\beta_2' & = & \{(1, \theta) \mid   ~ |\theta - (\zeta_2 + \pi) | \leq 1/10\} ~ .
\end{eqnarray*}

\begin{figure}[!htbp]
\centering
{\includegraphics[width=72mm]{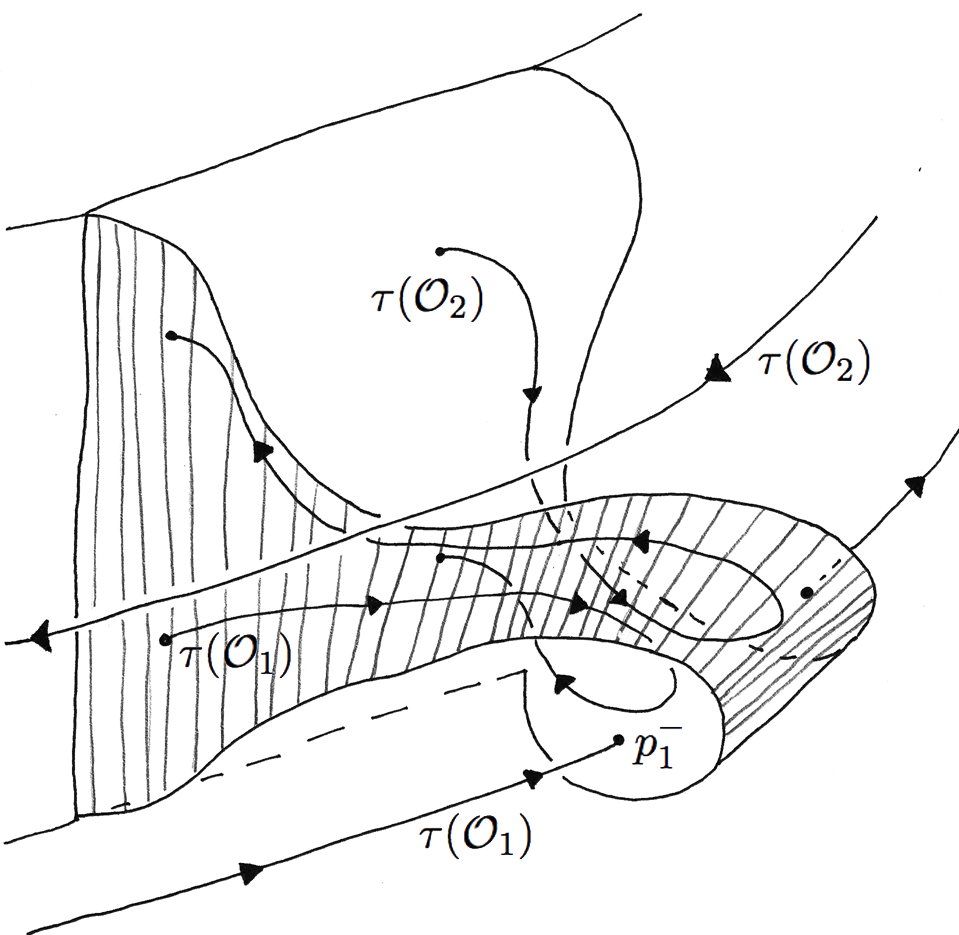}}
\caption{The image of $L_1\times [-2,2]$ under $\sigma_1$ \label{fig:twisted} }
\vspace{-6pt}
\end{figure}

 For $i=1,2$, choose orientation preserving  diffeomorphisms $\sigma_i \colon \alpha_i' \to \beta_i'$, and extend these maps to   smooth embeddings   $\sigma_i \colon D_i \to  \mW $, as illustrated in  Figure~\ref{fig:twisted}, which   satisfy the conditions:   
\begin{itemize}
\item[(K1)]   $\sigma_i(\alpha_i'\times z)=\beta_i'\times z$ for   $z\in [-2,2]$, and  the interior arc $\alpha_i^\prime$ is mapped to a boundary arc $\beta_i'$;  
\item[(K2)]     $\cD_i = \sigma_i(D_i) \subset  \{(r, \theta, z) \mid 1 \leq r \leq 5/2 ~ , ~ |\theta - (\zeta_i + \pi) | \leq 1/10\}$,    thus $\cD_1 \cap \cD_2 =\emptyset$;
\item[(K3)]  $\sigma_1(L_1^-) \cap \{r \geq 2\} \subset \{z < 0\}$  and $\sigma_2(L_2^-) \cap \{r \geq 2\} \subset \{z > 0\}$; 
\item[(K4)]    For every $x \in L_i$, the image   $\cI_{i,x} \equiv \sigma_i(x \times
  [-2,2])$ is an arc contained in a trajectory of $\cW$;
\item[(K5)]  Each slice $\sigma_i(L_i\times\{z\})$ is transverse to the vector field $\cW$, for all $-2\leq z \leq 2$; 
\item[(K6)]   $\cD_i$ intersects the periodic orbit $\cO_i$ and not $\cO_j$, for   $i \ne j$.
\end{itemize}

The ``horizontal faces'' of the embedded regions   $\cD_i = \sigma_i(D_i) \subset \mW$ are labeled by
\begin{equation}\label{eq-tongues}
\cL_1^{\pm}= \sigma_1(L_1^\pm)  ~ , ~   
\cL_2^{\pm} =   \sigma_2(L_2^\pm) ~.
\end{equation}

Note that the arcs $\cI_{i,x}$ in condition (K4) are  line segments from $\sigma_i(x \times \{-2\})$   to $\sigma_i(x \times \{2\})$ which follow the $\cW$-trajectory and traverse the insertion from the bottom face to the top face. 
Since $\cW$ is vertical near the boundary of $\mW$, and horizontal in the two periodic orbits,  we have that the arcs $\cI_{i,x}$    are vertical near the inserted curve $\sigma_i(\alpha_i')$ and horizontal at the intersection of the insertion with the periodic orbit $\cO_i$. 
Thus,  the embeddings of the surfaces $\sigma_i(L_i\times\{z\})$ make a {\it half   turn} upon insertion, for each   $-2 \leq z \leq 2$,  
  as depicted in Figure~\ref{fig:twisted}.  The turning is clockwise for the bottom insertion $i=1$  as illustrated  in Figure~\ref{fig:twisted} and counter-clockwise for the upper insertion $i=2$, which is not illustrated.

The image of the first insertion $\sigma_1(D_1)$ in Figure~\ref{fig:twisted}   
intersects the first periodic orbit of $\cW$ and is disjoint of the
second periodic orbit.  
  The image of the second insertion
$\sigma_2(D_2)$ is disjoint from the first
insertion and the first periodic orbit, and  intersects
the second periodic orbit.

The embeddings $\sigma_i$ are also required to satisfy two further conditions, which are the key to showing that the resulting Kuperberg flow   $\Phi_t$ is \emph{aperiodic}:   
\begin{itemize}
\item[(K7)] For $i=1,2$, the disk $L_i$ contains a point $(2,\theta_i)$ such that
  the image under $\sigma_i$ of the vertical segment
  $(2,\theta_i)\times[-2,2] \subset D_i \subset \mW$ is an arc
  $\{r=2\} \cap \{\theta_i^- \leq \theta \leq \theta_i^+\} \cap
  \{z=(-1)^i\}$  of the periodic orbit $\cO_i$ of $\cW$.   \\
\item[(K8)] {\it Radius Inequality}: For all $x = (r, \theta,  z) \in  L_i \times [-2,2]$, let $x' = (r', \theta',z') = \sigma_i(r, \theta,  z) \in \cD_i$,  then $r' < r$ unless  $x =  (2,\theta_i, z)$.
\end{itemize}

The Radius Inequality (K8) is one of the most fundamental  concepts of Kuperberg's construction. The condition (K4) and the fact that  the flow of the vector field $\cW$ on $\mW$ preserves the radius coordinate on $\mW$, allows restating (K8) in the more concise form for points in the faces $\cL_i^-$ of the insertion regions $\cD_i$. For $x' = (r', \theta',z') = \sigma_i(r, \theta,  z) \in \cD_i$ we have
 \begin{equation}\label{eq-radius}
r(\sigma_i^{-1}(x')) \geq r' ~ {\rm for} ~ x' \in \cL_i^- ~, ~ {\rm with }~ r(\sigma_i^{-1}(x')) = r' ~{\rm if ~ and ~ only ~ if} ~ x' = \sigma_i(2, \theta_i,  -2) ~. 
\end{equation}
The Radius Inequality (K8) is illustrated in Figure~\ref{fig:radius} below. This is an ``idealized'' case, as it implicitly assumes that the relation between the values of $r$ and $r'$ is ``quadratic'' in a neighborhood of the special points $(2,\theta_i)$, which is not required in order that  (K8) be satisfied. Later in this work, this ``quadratic condition'' will be added as part of the generic hypotheses on the construction.
\begin{figure}[!htbp]
\centering
{\includegraphics[height=42mm]{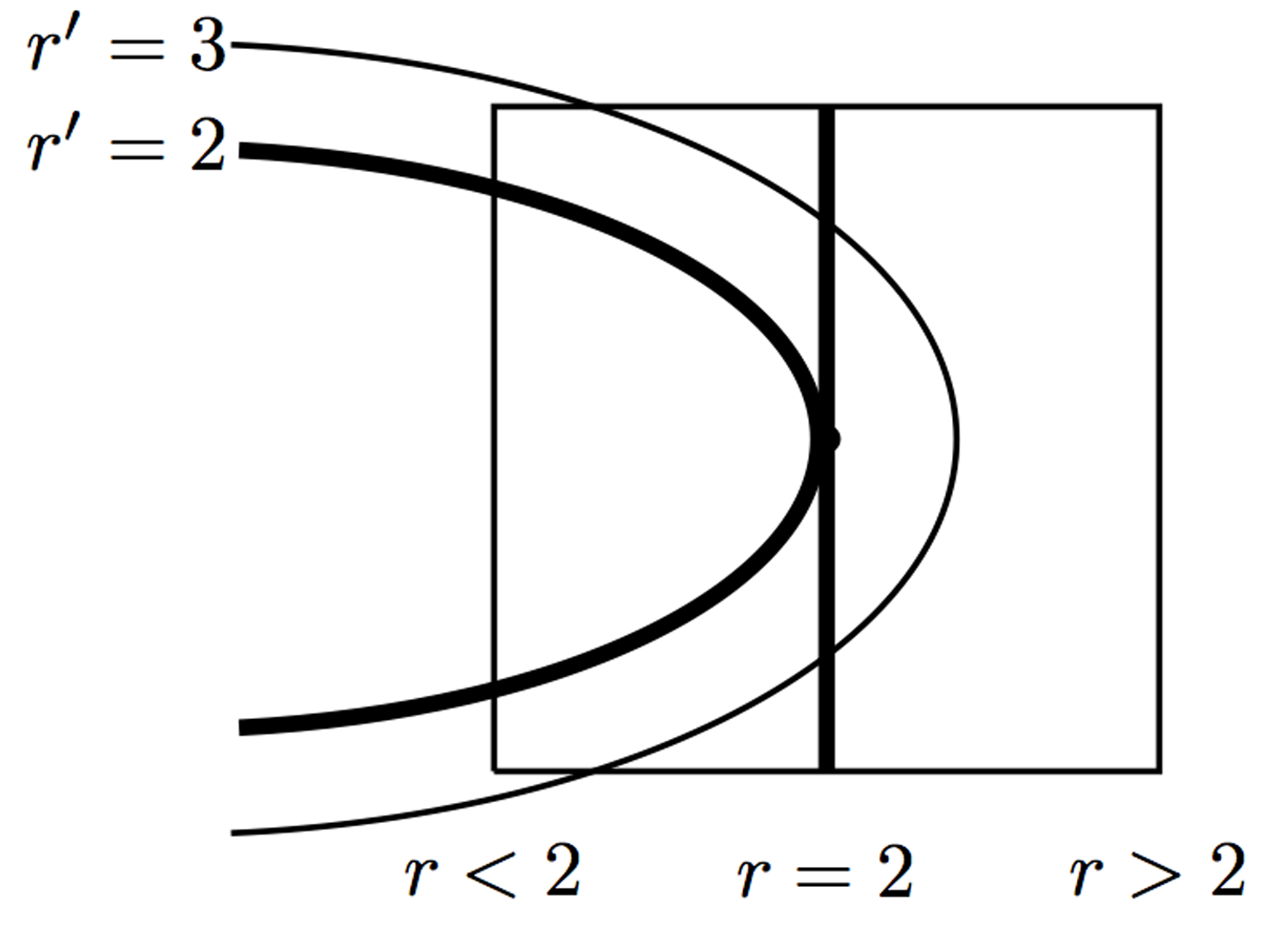}}
\caption{The radius inequality illustrated \label{fig:radius}}
\vspace{-6pt}
\end{figure}

 The  embeddings $\sigma_i \colon L_i \times [-2,2] \to \mW$, for $i=1,2$,  can be  constructed by first choosing smooth  embeddings of the faces  $\ds  \sigma_i \colon L_i^- \times \{-2\} \to \mW$ so that the image surfaces    are transverse to the vector field  $\cW$ on $\mW$, and satisfy the conditions (K1), (K3), (K7) and (K8). Then we extend the   embeddings of the faces $L_i^- \times \{-2\}$ to the cylinder sets $L_i \times [-2,2]$ by flowing the images using a reparametrization of the flow of $\cW$, so that we obtain embeddings of $L_i \times [-2,2]$ satisfying conditions (K1) to (K8), and as pictured in Figure~\ref{fig:twisted}.

\begin{figure}[!htbp]
\centering
{\includegraphics[width=120mm]{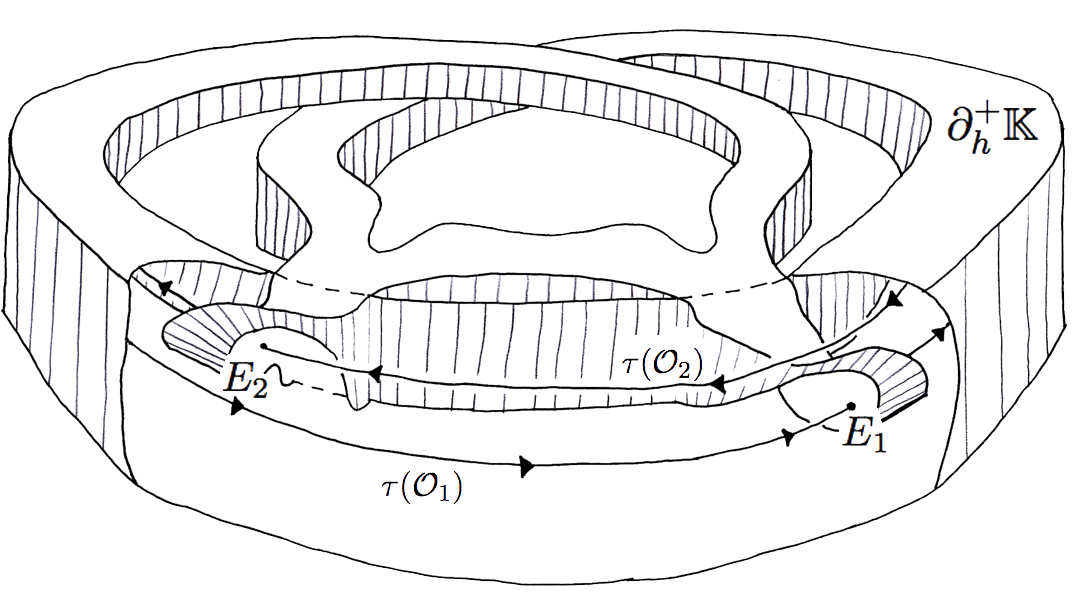}}
 \caption{The Kuperberg Plug $\mK$ \label{fig:K}}
 \vspace{-6pt}
\end{figure}

Finally,  define $\mK$ to be the quotient manifold obtained from $\mW$ by
identifying the sets $D_i$ with $\cD_i$. That is, for each point $x
\in D_i$ identify $x$ with $\sigma_i(x) \in \mW$, for $i = 1,2$. 
 This is    illustrated  in Figure~\ref{fig:K}.

 The restricted $\Psi_t$-flow on the inserted disk $\cD_i = \sigma_i(D_i)$ is not compatible with the image of the restricted $\Psi_t$-flow on $D_i$.  Thus, to obtain a smooth vector field $\cX$ from this construction,  it is necessary to modify $\cW$ on each insertion $\cD_i$. 
 The idea is to  replace  the vector field $\cW$  on the interior of each  region $\cD_i$ with the image vector field, so that the dynamics of $\Phi_t$ in the interior of each insertion region $\cD_i$ reverts back to the Wilson dynamics on $D_i$. This requires a minor technical step first. 

Smoothly  reparametrize   the image  of $\cW|_{D_i}$ under $\sigma_i$  on an open neighborhood of the boundary of $\cD_i$ so that it agrees with the restriction of $\cW$ to the same neighborhood. This is possible since the vector field $\cW$ is vertical on a sufficiently small open  neighborhood of the boundary of $D_i$ so is mapped by $\sigma_i$ to an orbit segment of $\cW$ by (K4). We obtain a vector field $\cW_i'$ on $\cD_i$ with the same orbits as the image of $\cW|_{D_i}$. 
The case of $\cD_1$ is illustrated in Figure~\ref{fig:twisted}.

 Then modify $\cW$ on each insertion $\cD_i$, replacing it with the modified image   $\cW_i'$.   Let $\cW'$ denote the vector field on  $\mW$ after these modifications and note that $\cW'$ is smooth.
By the modifications made above, the 
vector field $\cW^\prime$ descends to a smooth vector field on $\mK$ denoted by $\cK$. 
Let $\Phi_t$ denote the flow of the vector field $\cK$ on $\mK$.   
The  \emph{Kuperberg Plug} is the resulting space, $\mK \subset \mR^3$.
 
 \bigskip
 
\section{Transition points and the radius function}\label{sec-radius}
  
 In this section, we introduce notations that will be used throughout this work and also some basic concepts which are fundamental for relating the dynamics of the two vector fields $\cW$ and $\cK$. These results are contained in the literature \cite{Ghys1995, Kuperbergs1996, Kuperberg1994, Matsumoto1995}, though in a variety of differing notations and presentations.

Recall that   $\cD_i = \sigma_i(D_i)$ for $i =1,2$   are   solid $3$-disks embedded in $\mW$.
 Introduce  the sets:
\begin{equation}\label{eq-notchedW}
\mW' ~   \equiv   ~ \mW - \left \{ \cD_1  \cup \cD_2 \right\} \quad , \quad 
\wmW~   \equiv   ~ \overline{\mW - \left \{ \cD_1  \cup \cD_2 \right\}} ~ .
\end{equation}
The  closure $\wmW$ of $\mW'$  is the  \emph{pi\`{e}ge de Wilson  creus\'{e}}  as defined in \cite[page 292]{Ghys1995}.
The compact space $\wmW \subset \mW$ is the result of  ``drilling out''  the interiors of $\cD_1$ and $ \cD_2$, as the   terminology \emph{creus\'{e}} suggests.

For $x, y \in \mK$, we say that $x \prec_{\cK} y$ if   there exists $t \geq 0$ such that $\Phi_t(x) = y$.
Likewise, for $x',y' \in \mW$, we say that  $x' \prec_{\cW} y'$  if   there exists $t \geq 0$ such that $\Psi_t(x') = y'$.

  Let   $\tau \colon \mW  \to \mK$ denote the quotient map,
   which for $i = 1,2$,    identifies   a point $x \in D_i$ with its image $\sigma_i(x)   \in \cD_i$.
Then   the  restriction $\tau' \colon \mW' \to \mK$ is injective and onto. Let   $(\tau')^{-1}  \colon \mK \to \mW'$ denote the inverse map, which followed  by the inclusion $\mW' \subset \mW$, yields  the (discontinuous) map  $\tau^{-1} \colon \mK \to \mW$, where  $i=1,2$, we have:
\begin{equation}\label{eq-radiusdef}
\tau^{-1}(\tau(x)) =x ~ {\rm for} ~ x \in D_i  ~ ,~ {\rm and} ~   \sigma_i(\tau^{-1}(\tau(x))) = x ~ {\rm for} ~ x \in \cD_i~.
\end{equation}
For $x \in \mK$,  let $x=(r,\theta,z)$  be  defined as  the $\mW$-coordinates   of  $\tau^{-1}(x) \subset \mW'$. 
In this way, we obtain   (discontinuous)     coordinates $(r,\theta,z)$ on $\mK$. 
In particular, let $r \colon \mW' \to [1,3]$ be 
 the restriction of the   radius coordinate on $\mW$, then the function is extended to   the \emph{radius function}  of $\mK$, again denoted by $r$, where for $x \in \mK$ set $r(x) = r(\tau^{-1}(x))$.  
 
     The flow of the vector field $\cW$ on $\mW$ preserves the radius function on $\mW$, so   $x' \prec_{\cW} y'$ implies that      $r(x') = r(y')$.  
However, $x \prec_{\cK} y$ need not imply that $r(x) = r(y)$. 
The    points of discontinuity for the function $t \mapsto r(\Phi_t(x))$ play a fundamental role in the study of the dynamics of   Kuperberg flows.

Let $\partial_h^-\mK \equiv  \tau(\partial_h^- \mW\setminus (L_1^-\cup L_2^-))$ and $\partial_h^+\mK \equiv \tau(\partial_h^+ \mW\setminus (L_1^+\cup L_2^+))$  denote  the   bottom and top horizontal faces of $\mK$, respectively.   Note that each of the surfaces $\partial_h^-\mK$ and $\partial_h^+\mK$ are closed homeomorphic to  a twice-punctured torus, as can be seen in Figure~\ref{fig:K} and have boundary  which is the union of two circles. By the choice of $\cK$, the vertical boundary component  $\partial_v\mK \equiv \tau(\wmW \cap \partial_v \mW)$   is  tangent to $\cK$.

Points $x' \in \partial_h^- \mW$ and $y' \in \partial_h^+ \mW$ are said to be \emph{facing}, and we  write $x' \equiv y'$,  
if      $x' = (r, \theta, -2)$ and $y' = (r, \theta, 2)$ for some $r$ and $\theta$.  The entry/exit property of the Wilson flow is then equivalent to the property that  $x' \equiv y'$  if $[x',y']_{\cW}$ is an orbit from $\partial_h^- \mW$ to $\partial_h^+ \mW$ whenever $r(x') \ne 2$. 
There is also a notion of facing points for $x, y \in \mK$,   if either of   two cases are satisfied: 
\begin{itemize}
\item For  $x = \tau(x') \in \partial_h^-\mK$ and $y = \tau(y') \in \partial_h^+\mK$, if  $x' \equiv y'$ then $x \equiv y$.
\item For    $i=1,2$, with $x = \sigma_i(x')$ and $y = \sigma_i(y')$,     if  $x' \equiv y'$ then $x \equiv y$.
\end{itemize} 
 The context in which the notation $x \equiv y$ is used dictates which usage applies.

Consider the embedded disks  $\cL_i^{\pm} \subset \mW$ defined by \eqref{eq-tongues}, which appear  as the   faces of the
insertions  in $\mW$.  Their images in the quotient manifold $\mK$ are denoted by:
\begin{equation}\label{eq-sections}
E_1= \tau(\cL_1^-) ~ , ~ S_1= \tau(\cL_1^+) ~ , ~   
E_2= \tau(\cL_2^-) ~ , ~ S_2= \tau(\cL_2^+)   ~ .
\end{equation}
Note that $\tau^{-1}(E_i) = L_i^-$, while $\tau^{-1}(S_i) = L_i^+$.  
Also, introduce the notation:
\begin{equation}\label{eq-transversal}
\cT_{\cK}  ~ = ~  E_1 \cup E_2 \cup S_1 \cup S_2 ~. 
\end{equation}
Then by the formulation \eqref{eq-radius} of the Radius Inequality,
the points of discontinuity along $\cK$-orbits for the radius function $r \colon \mK \to [1,3]$ are contained in the set $\cT_{\cK}$.

The set   $\cT_{\cK} \subset \mK$ is transverse to the flow $\cK$. The
{\it transition   points}  of an orbit of $\cK$ are those points where the orbit intersects
the transversal $\cT_{\cK}$, or terminate in a boundary component  $\partial_h^-\mK$ or $\partial_h^+\mK$.
 They    are then either \emph{primary} or \emph{secondary} transition
 points, where $x\in \mK$ is: 
\begin{itemize}
\item a \emph{primary entry  point} if $x\in \partial_h^-\mK$;
\item  a \emph{primary exit  point} if $x \in \partial_h^+\mK$;
\item a \emph{secondary entry  point} if $x \in E_1 \cup E_2$;
\item a \emph{secondary exit  point}  $x \in S_1 \cup S_2$.
\end{itemize}
If a $\cK$-orbit contains no transition points, then it lifts to a $\cW$-orbit in $\mW$   flowing   from    $\partial_h^- \mW$ to     $\partial_h^+ \mW$.

The  \emph{special points} for the flow $\Phi_t$ are the images, for $i=1,2$, 
\begin{equation}\label{eq-special}
p_i^{-} = \tau(\cO_i \cap \cL_i^{-}) \in E_i ~ , ~ p_i^{+} = \tau(\cO_i \cap \cL_i^{+}) \in S_i ~ .
\end{equation}
Note that by definitions and the Radius Inequality, we have $r(p_i^{\pm}) = 2$ for $i=1,2$. 
  
\begin{figure}[!htbp]
\centering
{\includegraphics[width=80mm]{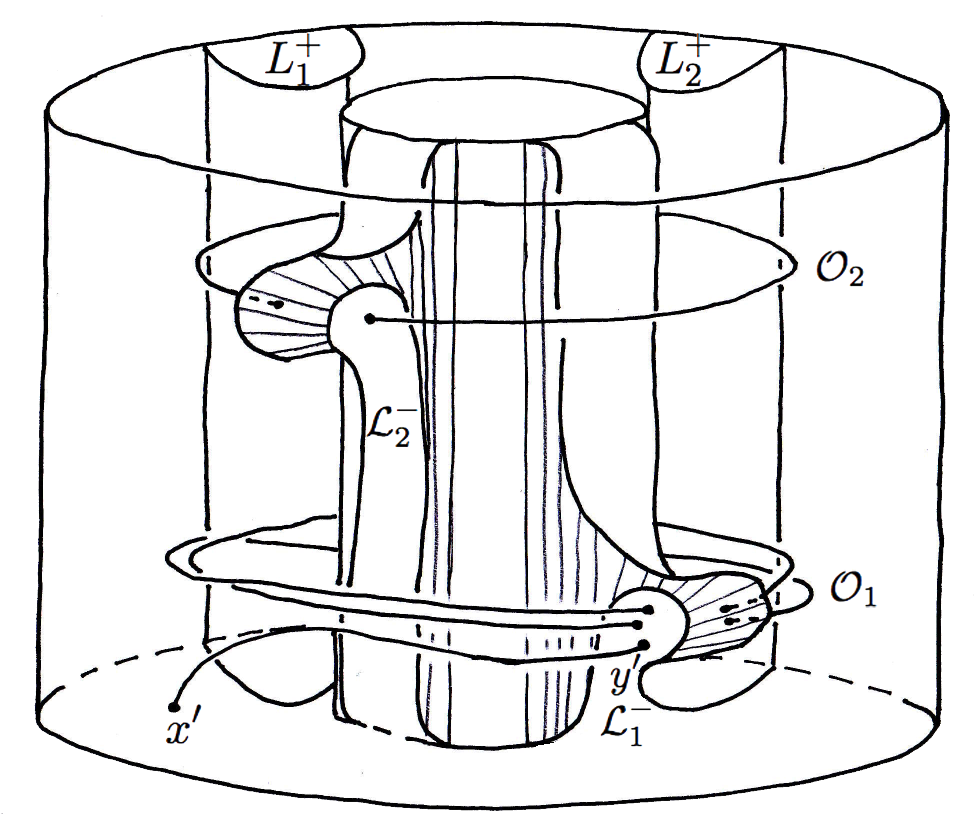}}
\caption{$\cW$-arcs lifted to $\wmW$ \label{fig:cWarcs}}
\vspace{-6pt}
\end{figure}

 A \emph{$\cW$-arc} is a closed segment $[x,y]_{\cK} \subset \mK$ of the flow of $\cK$ whose endpoints $\{x,y\}$ lie in $\cT_{\cK}$,  while the   interior $(x,y)_{\cK}$ of the arc is disjoint from $\cT_{\cK}$. The open interval $(x,y)_{\cK}$  is then the image under $\tau$ of a unique $\cW$-orbit segment in $\mW'$, denoted by $(x',y')_{\cW}$  where $\tau(x') = x$ and $\tau(y') = y$ (see Figure~\ref{fig:cWarcs}.) 
 Let $[x', y']_{\cW}$ denote the closure of $(x',y')_{\cW}$   in $\wmW$, then  we say that $[x', y']_{\cW}$ is the \emph{lift} of $[x,y]_{\cK}$.  Note that   the radius function $r$ is constant along   $[x', y']_{\cW}$.

 The properties of the Wilson flow $\cW$ on $\wmW$ determine   the endpoints of lifts   $[x',y']_{\cW}$.  We state the six cases which arise explicitly, as they will be cited in  later arguments.    Figure~\ref{fig:cWarcs} helps in visualizing these cases.  
 \begin{lemma}\label{lem-cases1}
Let  $[x,y]_{\cK} \subset \mK$ be a $\cW$-arc, and let $[x',y']_{\cW} \subset \wmW$ denote its lift.
 \begin{enumerate}
 
\item {\bf (p-entry/entry)}   If $x$ is a primary entry point, 
then $\ds x' \in    \partial_h^-  \mW \setminus (L_1^-   \cup L_2^- )$ 
and  if $y$ is   a secondary   entry point, we have   $y' \in \cL_i^-$ for  $i =1$ or $2$.  
  
\item   {\bf (p-entry/exit)} If $x$ is a primary entry point, then $x' \in   \partial_h^-
  \mW \setminus ( L_1^-   \cup L_2^- )$ and  if  $y$ is an exit point, 
  then   $\ds y' \in \partial_h^+\mW$ is a primary exit point, and by  the entry/exit condition on $\mW$ we have   $x\equiv y$.     

\item  {\bf   (s-entry/entry)} If $x$ is a secondary entry point, then $x' \in    L_i^-$  for  $i =1$ or $2$, and  if $y$ is   an entry point,  then we have   $y'  \in    \cL_j^- $ where $j =1,2$   is not necessarily    equal to $i$.

\item   {\bf  (s-entry/exit)} If $x$ is a secondary entry point, then $x' \in    L_i^-$  for  $i =1$ or $2$, and  if  $y$ is an exit point, then $y'
  \in    L_i^+$ and $x \equiv y$ by the entry/exit condition  of $\mW$.  

\item  {\bf  (s-exit/entry)} If $x$ is a secondary exit point, then $x' \in   \cL_i^+$  for  $i =1$ or $2$, and  if $y$ is   an entry point, so that $y' \in   \cL_j^-$ then $j=2$ if  $i=2$, and $j =1,2$  if $i=1$.

\item {\bf  (s-exit/exit)} If $x$ is a secondary exit point, then $x' \in   \cL_i^+$  for  $i =1$ or $2$, and  if   $y$ is a  primary
  exit point, $y' \in  \{ \partial_h^+ \mW \setminus ( L_1^+ \cup L_2^+) \}$. If $y$ is a secondary exit point, then $y'  \in    L_j^+$, where  $j =1$ or $2$ is not necessarily    equal to $i$. 
\end{enumerate}
\end{lemma}

 Figure~\ref{fig:cWarcs} illustrates some of the notions discussed in this section. 
 The disks $L_1^-$ and $L_2^-$ contained in  $\partial_h^-\mW$ are
 drawn in the bottom face, though they are partially obscured by the
 cylinder $\{r=1\}$.
The  image of $L_1^-$ under $\sigma_1$ is  the entry face of the insertion region in the lower half of the cylinder, 
while the  image of $L_2^-$ under $\sigma_2$ is  the entry face of the insertion region in the upper half of the cylinder. 
Analogously, the   disks $L_1^+$ and $L_2^+$ in $\partial_h^+\mW$  are mapped to the exit region of the insertion regions. 
The intersection of the $\cW$-periodic orbits $\cO_1$ and $\cO_2$ with $\wmW$ is illustrated,
as well as two  $\cW$-arcs in $\mW'$ that belong to the same orbit. One $\cW$-arc  goes from
$\partial_h^-\mW$ to $\cL_1^-$, hence from a principal entry point to
a secondary entry point (as in Lemma~\ref{lem-cases1}.1).  The second $\cW$-arc   goes from $\cL_1^+$ to
$\cL_1^-$, thus from a secondary exit point to a secondary entry point (as in Lemma~\ref{lem-cases1}.5).

    Introduce the radius coordinate function along $\cK$-orbits. For $x \in \mK$  set     $\rho_x(t) \equiv r(\Phi_t(x))$.  Note that if $\Phi_{t}(x) \not\in \cT_{\cK}$ then the    function $\rho_x(t)$  is locally constant at $t$, and thus we have: 
\begin{lemma} \label{lem-constant}
If the $\cK$-arc $\{\Phi_t(x) \mid t_0 \leq t \leq t_{1}\}$  contains no transition point, then $\rho_x(t)=\rho_x(t_0)$ for all $t_0 \leq  t \leq  t_1$. \hfill $\Box$
\end{lemma}

The  \emph{level function}  along an orbit  indexes the discontinuities of the  radius function. 
Given $x \in \mK$, set      $n_x(0) = 0$   and for $t > 0$ define 
\begin{equation}\label{def-level+}
n_x(t) = \#   \left\{ (E_1 \cup E_2) \cap \Phi_s(x) \mid  0 < s \leq t \right\} - \#   \left\{ (S_1 \cup S_2) \cap \Phi_s(x) \mid  0 < s \leq t  \right\} .
\end{equation}
 That is,  $n_x(t)$ is the total number of secondary entry points, minus
 the total number of secondary exit points,  traversed by the flow of $x$ over the interval 
 $0 < s \leq  t$.  For example, suppose that   $x_1 = \Phi_{t_1}(x)$ is the first transition point for $t > 0$.  If $x_1$ is a secondary entry point, then $n_x(t_1) = 1$, while $n_x(t_1) = -1$ if $x_1$ is a secondary exit point. 
Thereafter,     $n_x(t)$ changes value by $\pm 1$  at  each $t > t_1$ such that $\Phi_t(x) \in \cT_{\cK}$ and whether the value increases or decreases, indicates whether the transition point is an entry or exit point.

The function can be extended to negative time by setting, for $t < 0$,
 \begin{equation}\label{def-level-}
n_x(t) =  \#   \left\{ (S_1 \cup S_2) \cap \Phi_s(x) \mid  t < s \leq 0  \right\}    -   \#   \left\{ (E_1 \cup E_2) \cap \Phi_s(x) \mid  t < s \leq 0 \right\}.
\end{equation}

 Endow the manifold  $\mK$ with a Riemannian metric, induced from the natural coordinates on $\wmW$. Along the boundaries $\partial_h^\pm \wmW$ it is necessary to make a small adjustment to the Riemannian metric so that it becomes smooth on the quotient space $\mK$. We can assume that the vector field $\cK$ on $\mK$ has unit length, and so the adjustments are made to the metric transverse to $\cK$ in a small open neighborhood of these boundaries.
 Let    $d_{\mK}$ denote the resulting path length metric on $\mK$.

Let   $d_{\mW}$ be the path length metric in $\mW$. By the assumption that $\cK$ is unit length, the metric $d_{\mW}$ is just the same as that derived from the time coordinate along the flow $\Phi_t$. For example,  by this convention,  the length of the circles $\{(r,\theta,0) \mid 0 \leq \theta \leq  2\pi\} \subset \mK$ is $2\pi  r$.

 For $x' \prec_{\cW} y'$ in $\mW$,  let   $d_{\cW}(x',y')$ denote the
 path length of the $\cW$-orbit segment $[x' ,y']_{\cW}$ between
 them. Similarly, for $x \prec_{\cK} y$ in $\mK$,  let
 $d_{\cK}(x,y)$ denote the path length of the $\cK$-orbit segment $[x,y]_{\cK}$.  Note that if $[x,y]_{\cK}$ is a $\cW$-arc with lift $[x',y']_{\cW}$, then   $d_{\cK}(x,y) = d_{\cW}(x',y')$.

We establish some basic length estimates which are used in later sections.

\begin{lemma}\label{lem-wilsonlengths}
Let $0 < \e < 1$. There exists $L(\e) > 0$ such that for any $\xi \in \mW$ with $|r(\xi) -2| \geq \e$,    the total $\cW$-orbit segment  $[x',y']_{\cW}$ through $\xi$ has   length bounded above by $L(\e)$.
\end{lemma}
\proof
Since $r(\xi) \ne 2$, the orbit of $\cW$ containing $\xi$ is finite, hence there exists $x_{\xi}' \in   \partial_h^- \mW$ such that 
$x_{\xi}' = \Phi_{t}(\xi)$ for some $t \leq 0$. Likewise, there exists    $s \geq  0$ such that $y_{\xi}' = \Phi_{s}(\xi) \in \partial_h^+ \mW$. Then $[x_{\xi}',y_{\xi}']_{\cW}$ is the complete $\cW$-orbit containing $\xi$. In particular, $d_{\cW}(x_{\xi}',y_{\xi}') < \infty$.

 The function $\xi \mapsto d_{\cW}(x_{\xi}',y_{\xi}')$ is continuous
 on the open domain $r(\xi) \ne 2$. Since the flow in $\mW$ is rotationally
 invariant, this function depends only on $r(\xi)$. Let $L(\e)$ denote the maximum of this length function on the compact domain   $|r(\xi) -2| \geq \e$. 
\endproof

\begin{lemma}\label{lem-segmentlengths} 
There exists $0 < d_{min} <
  d_{max}$ such that  if  $[x',y']_{\cW} \subset \wmW$ is the lift of
  a $\cW$-arc $[x,y]_\cK$, then we have the uniform estimate
\begin{equation}\label{eq-segmentlengths}
d_{min} \leq  d_{\cW}(x',y') \leq d_{max} ~ .
\end{equation}
\end{lemma}
\proof
  
 First, suppose that $x' \in \partial_h^- \mW$ and so either $y' \in \cL_i^-$ for $i = 1,2$  or $y' \in \partial_h^+ \mW$. The set of points  $x' \in \partial_h^- \mW$ whose forward orbit   has first transition point in $\cL_i^-$ is a compact set containing the circle $\{(2, \theta, -2) \mid 0 \leq \theta \leq 2\pi \}   \subset \mW$ in its interior. Thus, there is   a lower and upper bound for  the length  $d_{\cW}(x',y')$. 
 
 Away from the core circle in $\partial_h^- \mW$,
the $\cW$-orbit for a point   $x' \in \partial_h^- \mW$ with $r(x') \ne 2$ is not trapped, 
so terminates in $y' \in \partial_h^+ \mW$. Moreover, the set of $x'$ whose   $\cW$-orbit    does not intersect the compact set $\cL_i^-$ has $r(x')$ value bounded away from $2$ by the previous case, and thus there exists $\e > 0$ so  $|r(x') - 2| \geq \e$. Then by Lemma~\ref{lem-wilsonlengths} the length $d_{\cW}(x',y')$  is then bounded above by $L(\e)$. 
The lower bound on $d_{\cW}(x',y')$ follows from the compactness of the set 
$\ds \{x' \in \partial_h^- \mW \mid ~ |r(x') - 2| \geq \e\}$.

The analysis of the    cases where $x' \in \cL_i^+$, for $i=1,2$, proceeds similarly.   
 \endproof

\begin{cor}\label{cor-segmentlengths} 
Let $[x,y]_{\cK} \subset \mK$ be a   \emph{$\cW$-arc}. Then there is a uniform length estimate
 \begin{equation}\label{eq-segmentlengths2}
d_{min} \leq  d_{\cK}(x,y) \leq d_{max} ~ . 
\end{equation}
\end{cor}
These results then combine to give the observation:
\begin{cor}\label{cor-orbitnumber} 
The Wilson orbit through $x' \in \mW$ with $r(x') \ne 2$ contains at
most a finite number of lifts of distinct  $\cW$-arcs, partially-ordered by the relation $\prec_{\cW}$.
\end{cor}

 We summarize the dynamical properties of orbits for  $\Phi_t$ which follow from the  above results. 
 
 For each $x \in \mK$ the $\cK$-orbit $\{\Phi_t(x) \mid a \leq t \leq b\}$, 
 where $-\infty \leq a < b \leq \infty$,  can be decomposed in a
 finite or infinite
 family of $\cW$-arcs $\ds \left\{[x_i', y_{i+1}']_{\cW}
 \right\}\subset \mW$.
 These $\cW$-arcs are indexed relative to the initial point $x$ by the
 level function $n_x(t)$ and the radius function $\rho_x(t)$, both
 functions are constant along $(x_i, y_{i+1})_{\cK}$. 
  Thus, the $\cW$-arcs can be grouped according to the values $\rho_x(t)$ and so grouped, lie in orbits of the  Wilson flow in   disjoint cylinders $\{r=const.\}$.

 The $\cK$-orbit of $x$ is then determined by understanding when two
 segments  $[x_i', y_{i+1}']_{\cW}$ and $[x_j', y_{j+1}']_{\cW}$  have
 the same $r$-value and so lie in a  common cylinder. Then whether they are contained in a  common $\cW$-orbit   within that cylinder, so that $x_i' \prec_{\cW} x_j'$ or  $x_j' \prec_{\cW} x_i'$ . 
 By Corollary~\ref{cor-orbitnumber}  there can be only a finite number of   segments in each finite $\cW$-orbit. 
 
 Once the arcs have been grouped according to their ``$\cW$-orbit type'', it then remains to understand how these $\cW$-arcs are assembled to give the full $\cK$-orbit, while satisfying the constraints imposed by the Radius Inequality (K8).  
 
 \bigskip
 
\section{Semi-local dynamics}\label{sec-semilocal}

In this section, we establish  a variety of properties of the level function and how it relates the   dynamical  behavior of the Wilson flow $\Psi_t$ with the much more complicated dynamical behavior of the Kuperberg flow $\Phi_t$. This relationship  is one of the main themes in the understanding of the dynamics of the flow $\Phi_t$ as developed in this work. 

We first consider     the properties of finite length segments of $\cK$-orbits.
The results presented below  are formulations, in our notation, of results which are contained in the works 
   \cite{Ghys1995, Kuperbergs1996, Kuperberg1994, Matsumoto1995} and are the most basic techniques used in  the analysis of the dynamics of the Kuperberg flow.

We say that the $\cK$-orbit of  $x \in \mK$ is \emph{trapped in forward time} if  the forward orbit is defined for all $t \geq 0$,  so that  the segment $\{\Phi_t(x) \mid t \geq 0\} \subset \mK$. In particular, this forward orbit never intersects $\partial_h^+ \mK$.  Likewise,  the $\cK$-orbit of  $x \in \mK$ is \emph{trapped in backward time} if the segment $\{\Phi_t(x) \mid t \leq 0\} \subset \mK$. When it is clear whether forward or backward time is meant, such as for points on the faces   $\partial_h^{\pm} \mK$, we simply refer to a \emph{trapped orbit}. 

The $\cK$-orbit of $x$ is \emph{infinite} if it is trapped in both forward and backward time.

 For a $\cK$-orbit segment  $[x,y]_{\cK}$, we define its ``lift'' to $\wmW$, which is a union of lifts of $\cW$-orbit segments. This is a fundamental technical construction and introduces another useful notational convention.
 Let   $0 \leq t_0 < t_1 < \cdots < t_n$ be such that $x_{\ell} = \Phi_{t_{\ell}}(x)$ are  the successive transition points in  $[x,y]_{\cK}$ so that   
 \begin{equation}\label{eq-decomposition}
[x, y]_{\cK} = [x,x_0]_{\cK} \cup [x_0,x_1]_{\cK} \cup \cdots \cup [x_{\ell},x_{\ell + 1}]_{\cK} \cup \cdots \cup [x_{n-1},x_n]_{\cK} \cup [x_n,y]_{\cK} \subset \mK ~.
\end{equation}
By convention,  $[x,x_0]_{\cK}$ is defined to be empty if $x = x_0$, and otherwise has a well-defined lift $[x',x_0']_{\cW} \subset \wmW$.  The case for $y$ is analogous. 
 Each $\cW$-arc $[x_{\ell},x_{\ell + 1}]_{\cK}$ for $0 \leq \ell < n$
 lifts to a $\cW$-orbit segment $[x_{\ell}' ,y_{\ell + 1}']_{\cW} \subset \wmW$, where  $\tau(x_{\ell +1}') =  \tau(y_{\ell +1}') = x_{\ell +1}$ though $y_{\ell +1}' \ne x_{\ell +1}'$.  The  lift of $[x,y]_{\cK}$ to $\mW$ is defined  to be the collection of orbit segments in $\mW$, 
 \begin{equation}\label{eq-decompositionW}
 \left\{ [x',y'_0]_{\cW} ,  [x'_0,y'_1]_{\cW} , \cdots , [x'_{\ell},y'_{\ell + 1}]_{\cW} , \cdots , [x'_{n-1},y'_n]_{\cW} , [x'_n,y']_{\cW} \right\} ~.
\end{equation}
The lift of a $\cK$-orbit segment  is closely related to the behavior of the  level function $n_x(t)$ along an orbit. 
Given  a $\cK$-orbit segment $[x_0 , x_n]_{\cK}$, where $x_0 = \Phi_{0}(x)$,   decompose it into $\cW$-arcs as in   \eqref{eq-decomposition}.
The  orbit $\Phi_t(x)$   for $t > 0$ flows from $x$  until it hits a
transition point $x_1$. This initial $\cW$-arc
$[x_0,x_1]_{\cK}$ 
is   the image under $\tau$ of  a $\cW$-orbit segment   $[x_0',y_1']_{\cW} \subset \wmW$ as in   \eqref{eq-decompositionW}. 
After this first    transition point, it then follows the image of the
vector field $\cW$ under an insertion   map $\sigma_i$, which is again
the image of a $\cW$-orbit segment  $[x_1',y_2']_{\cW} \subset \wmW$. At the
next transition point, the flow $\Phi_t(x)$ either exits the insertion
and reverts back to the same orbit of $\cW$ in $\wmW$; or, it flows
into another insertion, and thus the orbit  is the image of a
$\cW$-arc into an insertion within  the first insertion. This process continues along the entire $\cK$-orbit segment. The level function $n_x(t)$ counts the number of such insertions within insertions along the $\cK$-segment $[x_0,x_n]_{\cK}$, so  measures  the ``depth of penetration'' of the orbit into the self-insertion process used to construct the plug $\mK$.

\begin{figure}[!htbp]
\centering
{\includegraphics[width=80mm]{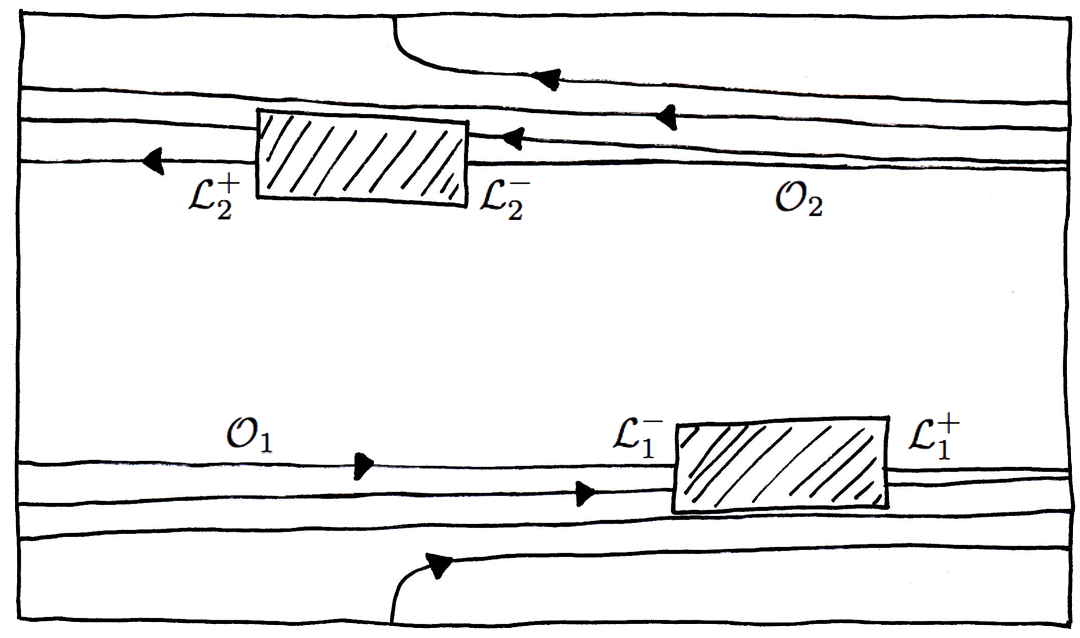}}
\caption[{Decomposition into $\cW$-arcs}]{Decomposition into $\cW$-arcs in the cylinders $\{r=2\}\subset \wmW$ \label{fig:cW-decomposition}}
\vspace{-6pt}
\end{figure}

We illustrate this orbit decomposition in
Figure~\ref{fig:cW-decomposition}, that represents the cylinder
$\{r=2\}$ in $\wmW$. Consider the first periodic orbit and its
intersection with $\wmW$, which is pictured as the lower
horizontal line in the drawing.
Flowing in forward time, the lift of the $\cK$-orbit first intersects   the
  entry face  $\cL_1^-$ of the lower insertion in a ``special point''. The flow then continues   to
 a point in the lower boundary of the same cylinder
$\{r=2\}$ where it enters. Thus the next $\cW$-arc in this orbit is in the same
cylinder. It starts climbing and turning, until hits the
first insertion and  exits the cylinder $\{r=2\}$ through  $\cL_1^-$. 
We prove, in Proposition~\ref{prop-propexist}, that after a certain
time this orbit comes back to the cylinder $\{r=2\}$ in the facing
point. Thus the $\cK$-orbit will eventually pass the insertion and continue in the
same $\cW$-orbit for another turn before hitting the insertion
again. This pattern then repeats infinitely many times. Flowing in
backward time, the lift of the $\cK$-orbit containing the arc
$\cO_1\cap \wmW$, first intersects the exit face $\cL_1^+$ in a
special point. The orbit then continues to a point in the upper bound
of the cylinder. It thus descends turning until it hits the exit face
$\cL_2^+$. Proposition~\ref{prop-propexist} applies, implying that
after a certain time this backward orbit comes back to the cylinder
and continues descending and turning. This process repeats infinitely
many times. Remark that the $\cK$-orbit containing the arc
$\cO_1\cap\wmW$ accumulates in forward time on $\cO_1\cap\wmW$ and in
backward time on $\cO_2\cap\wmW$. Applying the same analysis to the
$\cK$-orbit containing $\cO_2\cap\wmW$, we obtain the same
conclusion. This argument is explained again in the proof of Proposition~\ref{prop-r=2any}.

 We next give a sequence of technical results concerning   the properties of the lift \eqref{eq-decompositionW} of a path as in \eqref{eq-decomposition}, beginning with the very simple     concept  of a \emph{short-cut}, which   was introduced in \cite{Kuperberg1994}, and developed in more detail in \cite{Ghys1995, Kuperbergs1996, Matsumoto1995}.

 \begin{lemma}[Short-cut]\label{lem-shortcut}
 Suppose that    $x \in E_i$ and $y \in S_i$ are facing. Then there exists $x', y' \in \mW$ such that $\tau(x') = x$,  $\tau(y') = y$ and $x' \prec_{\cW} y'$. 
That is, there is a  $\cW$-orbit segment $[x', y']_{\cW} \subset \mW$ (the short-cut) between  $x'$ and $y'$.
\end{lemma}
 \proof Let $x' \in \cL_i^-$ with $\tau(x') = x$ and $y' \in \cL_i^+$
 with $\tau(y') = y$. Consider $x'' \in L_i^-$ with $\sigma_i(x'') = x'$ and $y'' \in L_i^+$ with $\sigma_i(y'') = y'$. Then $x \equiv y$ implies $x'' \equiv y''$ by definition. By (K4) the image $\sigma_i(\cI_{i, x''})$ is a $\cW$-orbit segment in $\cD_i$ with endpoints $x', y'$. 
  Thus, $x' \prec_{\cW} y'$.
 \endproof
 
Note that the short-cut path $[x',y']_{\cW} \subset \mW$ has endpoints
in $\cL_i^{\pm}$ and ``bridges the gap'' between these two faces in $\mW$. 
In   Figure~\ref{fig:cW-decomposition}, this may be viewed as filling
in with a flow line  across any of the insertion boxes, from an entry point to an exit point.

Note that for a short-cut, the length of the path segment  $[x',y']_{\cW}$ is
bounded above by Lemma~\ref{lem-segmentlengths}, independent of the
choice of the points $x,y$. On the other hand if we are only given that $x\prec_\cK y$, then no such a priori bound exists. Thus, the $\cW$-path  $[x',y']_{\cW}$ is truly a ``short-cut'' between $x'$ and $y'$ when compared to the length of the $\cK$-path between $x$ and $y$.

  The next result shows that given a pair of $\cW$-arcs in $\mK$ whose \emph{inner endpoints} are facing,   there is a $\cW$-orbit segment in $\mW$ containing the lifts of both $\cW$-arcs.

\begin{lemma} \label{lem-shortcut1}
Let $x, y, z, u \in \mK$   be successive transition points on a $\cK$-orbit, such that $x'
\prec_{\cW} y'$, $z' \prec_{\cW} u'$,  $y$ is an   entry point, $z$ is an   exit point and $y \equiv z$. Then $x' \prec_{\cW} u'$ and hence    $r(x') = r(u')$.
\end{lemma}
\proof
Observe that $x' \prec_{\cW} y'$ implies $y$ must be a secondary entry
point, thus by Lemma~\ref{lem-cases1}, the endpoint $y'$ of the
$\cW$-orbit segment  $[x',y']_{\cW}$ must lie in $\cL_i^-$ for some $i
=1,2$. Similarly, $z$ must be a secondary exit point, thus the
endpoint  $z'$ of the $\cW$-orbit segment  $[z',u']_{\cW}$ must lie in $\cL_j^+$ for some $j$. The assumption $y \equiv z$ implies that $i=j$.   Then  by Lemma~\ref{lem-shortcut},    there is a $\cW$-arc  in $\mW$  between $y'$ and $z'$, a short-cut. Thus,  $y' \prec_{\cW} z'$, and so 
  $x' \prec_{\cW} y' \prec_{\cW} z' \prec_{\cW} u'$ which implies $x' \prec_{\cW} u'$ and so  $r(x') = r(u')$. 
\endproof

 The following result   gives a criteria for when  a pair of $\cW$-arcs, whose inner endpoints lift to   points which can be joined by a $\cW$-orbit segment in $\mW$,  are themselves contained in a $\cW$-orbit segment and have facing endpoints.

\begin{lemma} \label{lem-shortcut2}   Suppose that $x  \in \mK$ is an entry point,  $u \in \mK$ is an exit point,
$[x,y]_{\cK}$ and $[z,u]_{\cK}$ are  $\cW$-arcs with lifts
  $[x',y']_{\cW}$ and $[z',u']_{\cW}$. If  $y' \prec_{\cW} z'$,  then $x' \prec_{\cW} u'$ and    $x \equiv u$. 
\end{lemma}

\proof
Let $[x',y']_{\cW} \subset \wmW$ be the lift of $[x,y]_{\cK}$. The fact that   
$x$ is an entry point implies that $x' \in \partial_h^- \mW$ and so either $y' \in \cL_i^-$ for $i = 1,2$, or $y' \in  \partial_h^+ \mW$. 
The assumption that $y' \prec_{\cW} z'$  implies there is a
$\cW$-orbit segment from $y'$ to $z'$ so the case $y' \in  \partial_h^+ \mW$ is impossible.
Similarly, let $[z',u']_{\cW} \subset \mW$ be the lift of
$[z,u]_{\cK}$, then $u$ is an exit point implies that $u'
\in \partial_h^+ \mW$ and  so $z' \in \cL_j^+$ for $j =1,2$ as $y'
\prec_{\cW} z'$ is given.  Thus, concatenating the $\cW$-orbits segments from $x'$ to $y'$ to $z'$ to $u'$ yields   $x' \prec_{\cW} u'$. The entry/exit property  for the Wilson flow on $\mW$ implies that  $x'$ and $u'$     are facing, hence  $x \equiv u$.
\endproof

The next result   gives a criteria using   the level function, 
for when  a sequence of three $\cK$-arcs admit a lift to a segment of the $\Psi_t$-flow. 

 \begin{lemma}\cite[Corollary 4.2]{Matsumoto1995} \label{lem-shortcut3}
Let $x   \in \mK$. Given    successive transition points   $x_{\ell} = \Phi_{t_{\ell}}(x)$ where $0 = t_0 < t_1 < t_2 < t_3$, suppose that $n_{x_0}(t) \geq 0$ for all $0 \leq  t < t_3$, and that $n_{x_0}(t_{2}) = 0$. 
Then    $x_1 \equiv x_{2}$, $x_0'  \prec_{\cW} y_3'$ and hence $r(x_0') = r(y_3')$.
\end{lemma}
\proof
For $0\leq \ell <  3$,  let  $[x_\ell', y_{\ell+1}']_{\cW} \subset \wmW$ be the     lift of the $\cW$-arc  $[x_\ell,x_{\ell+1}]_{\cK}$.
The fact that $n_{x_0}(t_0) = 0$ and the assumption    $n_{x_0}(t_1) \geq 0$ implies that $n_{x_0}(t_1) =1$.
 Then $x_1$ is a secondary entry point,  with $y_1' \in \cL_i^-$ for $i=1$ or $i=2$. 
Similarly, the assumption     $n_{x_0}(t_{2}) = 0$  implies   that $x_{2}$ is a secondary exit point,  with $x_{2}' \in \cL_j^+$ for $j=1$ or $j=2$. 
 
For the middle  $\cW$-arc  $[x_1', y_2']_{\cW}  \subset \wmW$,  we must then have  $x_1' \in L_i^-$ and $y_2' \in L_j^+$.  
By the entry/exit assumption on $\mW$,  we have  $x_1' \equiv y_2'$ so
that $i=j$ and $x_1\equiv x_2$. It follows  that  there exist a short-cut $[y_1', x_2']_{\cW} \subset \cD_i$ between $y_1'$ and $x_2'$. 
Thus we have  $x_0' \prec_{\cW} y_3'$.
\endproof

\begin{figure}[!htbp]
\centering
{\includegraphics[width=70mm]{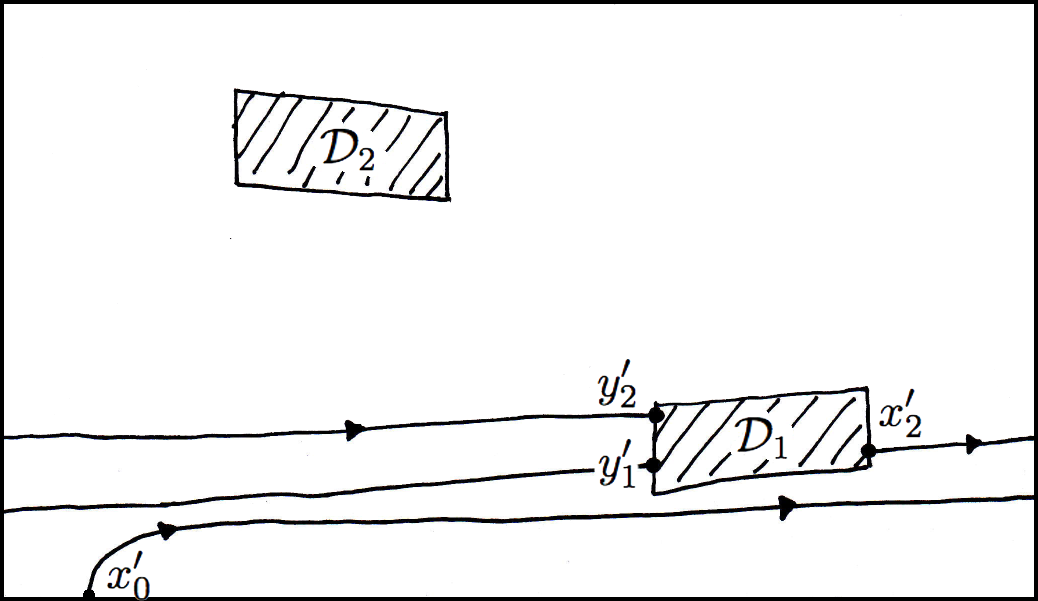}}
\caption{Three step curves \label{fig:3stepcurves}}
\vspace{-6pt}
\end{figure}

In geometric terms, the lifts of the $\cW$-arcs in Lemma~\ref{lem-shortcut3} to $\wmW$ do not form a continuous $\cW$-orbit segment, due to the discontinuity of the map $\tau$, but by replacing the middle segment 
$[x_1', y_2']_{\cW}$ with a short-cut $[y_1', x_2']_{\cW}$ through $\cD_i$ we obtain a $\cW$-orbit segment $[x_0' ,  y_3']_{\cW}$.
This is illustrated in Figure~\ref{fig:3stepcurves}.

\medskip

 Next, we use Lemmas~\ref{lem-shortcut} to \ref{lem-shortcut3} to obtain a general form of Lemma~\ref{lem-shortcut3}, which is one of the fundamental results concerning the finite orbits of the Kuperberg flow. 
 
\begin{prop}\cite[Lemme p.297]{Ghys1995} \label{prop-shortcut4}
Let $x   \in \mK$.  For $n \geq 3$, assume there are given    successive transition points   $x_{\ell} = \Phi_{t_{\ell}}(x)$ for $0 = t_0 < t_1 < \cdots < t_{n-1} < t_n$.
  Suppose that $n_{x_0}(t) \geq 0$ for all $0 \leq  t < t_n$ and that $n_{x_0}(t_{n-1}) = 0$. 
Then    $x_0'  \prec_{\cW} y_n'$, and hence $r(x_0') = r(y_n')$.
Moreover, if $x_0$ is an    entry point and $x_n$ is an exit point, then 
  $x_0\equiv x_n$.
\end{prop}
\proof
 For $0\leq \ell <  n$, let $[x_\ell', y_{\ell+1}']_{\cW} \subset \wmW$ be the     lift of
 the $\cW$-arc  $[x_\ell,x_{\ell+1}]_{\cK}$.
 The fact that $n_{x_0}(t_0) = 0$ and the assumption    $n_{x_0}(t_1) \geq 0$ implies that $n_{x_0}(t_1) =1$.
Then $x_1$ is a secondary entry point,  with $y_1' \in \cL_i^-$ for $i=1$ or $i=2$. 
Similarly, the assumption     $n_{x_0}(t_{n-1}) = 0$ and  $n_{x_0}(t_{n-2}) \geq 0$  implies   that $x_{n-1}$ is an exit point. Then $x_{n-1}$ must be a secondary exit point   and thus  $x_{n-1}' \in \cL_j^+$ for $j=1$ or $j=2$. 
The  case $n=3$ then follows from  Lemmas~\ref{lem-shortcut2}  and \ref{lem-shortcut3}.

For the case $n > 3$, we  proceed by induction.  Assume that  the
result holds for all $\cK$-segments containing at most  $n$ transition
points.   We will now prove the result for a segment with $n+1$
  transition points $x_\ell$ for $0\leq \ell \leq n$. Observe that by
  hypothesis, $n_{x_0}(t_1)=1=n_{x_0}(t_{n-2})$.   If there exists $3 \leq   \ell \leq n-2$ such that $n_{x_0}(t_{\ell -1}) = 0$, then consider the least such $\ell$. By the inductive hypothesis,  we have that 
$x_0' \prec_{\cW} y_{\ell}'$, and as $n_{x_0}(t_{\ell -1}) = 0$, we also have that $x_{\ell-1}' \prec_{\cW} y_n'$. 
Both  $\cW$-segments contain the arc $[x_{\ell-1}', y_{\ell}']_{\cW}
  \subset \wmW$.  Thus,  the last arc of the $\cW$-orbit segment
  $[x_0', y_{\ell}']_{\cW}$ and the first arc of the $\cW$-orbit
  segment $[x_{\ell-1}', y_n']_{\cW}$  must agree, hence   $x_0'
  \prec_{\cW} y_n'$ as claimed. 

  If in addition,   $x_0$ is an entry point and $x_n$ is an exit
  point, then   $x_0'  \prec_{\cW} y_n'$ implies   $x_0\equiv x_n$ by Lemma~\ref{lem-shortcut2}.

If  $n_{x_0}(t) \geq 1$ for all $t_1 \leq t < t_{n-1}$, then apply the
inductive hypothesis to the segment $[x_1, x_{n-1}]_{\cK}$ to obtain
$x_1' \prec_{\cW}  y_{n-1}'$. 
Using the induction, we have that since $x_1$ is an entry point
  and $x_{n-1}$ is an exit point, hence $x_1\equiv x_{n-1}$. Then,
  Lemma~\ref{lem-shortcut1} implies that $x_0'  \prec_{\cW} y_{n}'$ and $r(x_0') =
r(y_{n}')$. If in addition, $x_0$ is an entry point and $x_n$ is an
exit point,  then again, $x_0'  \prec_{\cW} y_n'$ implies that $x_0\equiv x_n$.
\endproof

Proposition~\ref{prop-shortcut4}  implies   that the flow $\Phi_t$ on $\mK$ satisfies the entry/exit condition.  

\begin{prop}\cite[Proposition 5.2]{Matsumoto1995} \label{prop-ee} 
Let   $x \in \partial_h^-\mK$   be a primary entry   point and suppose  $x_{\ell} = \Phi_{t_{\ell}}(x)$ for $0 = t_0
< t_1 < \cdots < t_{n-1} < t_n$   are    successive transition points. 
   If $y = x_n$ is a primary exit point,  then $x \prec_{\cW} y$ and hence $x
\equiv y$.  Moreover, $n_x(t) \geq 0$ for $0 \leq t < t_{n}$.
\end{prop}
\proof
  If $n=1$, then   $[x_0,x_1]_{\cK}$ is a $\cW$-arc, and
  the conclusion follows by the entry/exit condition for $\mW$. The
  case $n=2$ is impossible by Lemma~\ref{lem-cases1}. 
  
If    $n = 3$,  note that $x_1$ must be a secondary entry point,  hence $n_x(t_1) =1$.
The case  $n_x(t_2) = 2$ is impossible, as $x_3$ is a primary exit
point, hence $x_2$ must be a secondary exit point, contrary to
assumption. Thus $n_x(t_2) = 0$. Then by
Lemma~\ref{lem-shortcut3} we have $x_1' \equiv y_2'$ and $x_0'
\prec_{\cW} y_3'$. As $x_0'$ is  a primary entry point, and $y_3'$ is
a primary exit point of $\mW$, they must be facing, or $x \equiv y$.

Now assume that    $n > 3$. As before,  $x_1$ must be a secondary entry point
with $n_x(t_1) = 1$, so it suffices  to prove that
$n_x(t_{n-1})=0$ and $n_x(t)\geq 0$ for all $0 \leq t \leq t_{n-1}$,
and then apply Proposition~\ref{prop-shortcut4}. This will 
imply that $x_0' \prec_{\cW}y_n'$ and as $x_0=x$ is a primary entry
point and $x_n=y$ is a primary exit point, then $x\equiv y$.

We claim that $n_x(t) \geq  0$ for all $0 \leq t < t_n$.
 Suppose not, then there exists a least $2 < \ell < n$ such that
$n_x(t_{\ell}) = -1$. Then $n_x(t) \geq 0$ for all $0 \leq t < t_{\ell}$ and   $n_x(t_{\ell -1})=0$.
By Proposition~\ref{prop-shortcut4} we have 
 $x_0' \prec_{\cW} y_{\ell}'$,  and since $x_{\ell}$ is an
exit point, we have $x_0\equiv x_\ell$. This implies that $x_\ell$ is a
primary exit point, which is a contradiction as $\ell<n$.
    We can thus assume that $n_x(t) \geq  0$ for all $t_0 \leq t \leq t_n$ and that  $n_x(t_1) = 1$.

If $n_x(t_{n-1}) = k >0$, consider the greatest integer $0<\ell<n$ such that $n_x(t_\ell)=k$ and $n_x(t_{\ell-1})=k-1$.  
Then the $\cK$-orbit segment $[x_{\ell}, x_n]_{\cK}$ satisfies the
conditions of Proposition~\ref{prop-shortcut4} with $x_\ell$ an entry point, so that $x_{\ell}'
\prec_{\cW} y_n'$   and $x_\ell\equiv x_n$. As $x_n$ is a primary exit point and $x_{\ell} \equiv x_n$, the entry/exit property of the Wilson flow  implies that $x_{\ell}$ is a primary entry point. This is a contradiction as we have that   $\ell > 0$.
 Then   $n_x(t_{n-1}) = 0$ and the proof is finished.  
\endproof

 In summary,    the proofs of   Propositions~\ref{prop-shortcut4} and \ref{prop-ee}  show that, given a $\cK$-orbit segment   $[x, y]_{\cK}$   
 in $\mK$ with $x$ an entry point, $y$ a facing exit point, the condition  $n_x(t) \geq 0$ for the $\cK$-orbit between the two points,  then there exists  a $\cW$-orbit segment $[x_0',y_n']_{\cW}$ in $\mW$ between lifts $x_0'$ of $x$ and $y_n'$ of $y$. 
 The $\cW$-orbit segment  $[x_0', y_n']_{\cW}$ is obtained by an inductive ``short-cut''  procedure,  which   replaces the ``dynamics'' of the Kuperberg flow with that of the Wilson flow.  
 
  \bigskip
 
 \section{Dynamics and level}\label{sec-global}

In this section,   we consider   when it is possible to perform an inverse of the short-cut reduction used in Section~\ref{sec-semilocal}, which replaces a $\cW$-arc with a suitable $\cK$-orbit segment that spans the gap between the endpoints of a $\cW$-orbit segment.
The solution to this problem  depends on the relations between   the functions $\rho_x(t)$, $n_x(t)$, and  the lengths of $\cW$-orbit segments in $\mW$ as   will be shown  below.  An important application of this analysis   gives criteria    for when the $\cK$-orbit of a point $x \in \mK$ must necessarily escape through a face of $\mK$, either in forward or backward time.

 Recall that $\cD_i = \sigma_i(D_i) \subset \mW$ is the inserted compact   region. Define
 \begin{equation}\label{eq-rmax}
R_* = \max \{R_1 , R_2\} ~ , ~ R_i = \max \{r(x) \mid x \in \cD_i\} > 2
\end{equation}
Note that  for any $x \in \mW$ with $r(x) > R_*$, the $\cW$-orbit of
$x$ in $\mW$  does not intercept the inserted regions $\cD_i$, but it
might intersect the regions $D_i$. Hence
the $\cK$-orbit of   $\tau(x) \in \mK$   escapes through  the top
face  $\partial_h^+ \mK$ of $\mK$ or intersects the secondary exit
regions $S_i$ for $i=1$ or 2 in forward time, and escapes through  the bottom face $\partial_h^- \mK$ or intersects the secondary entry
regions $E_i$ for $i=1$ or 2 in backward time. 
Also, for the case where  $r(x) = R_*$, the orbit of $\tau(x)$ in
$\mK$ behaves in the same way as the flow on the boundary of $\cD_i$ agrees with the Wilson flow on $\mW$.

\begin{figure}[!htbp]
\centering
{\includegraphics[width=50mm]{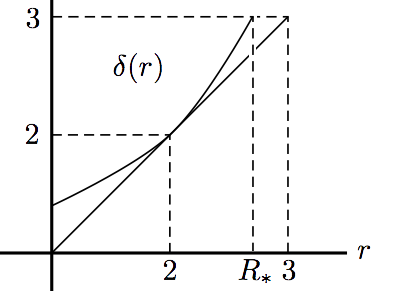}}
 \caption{\label{fig:r(delta)}  The function $\delta(r)$}
 \vspace{-6pt}
\end{figure}

 Let $\cC(r) = \{x \in \mW \mid r(x) = r \}$ denote the cylinder in
 $\mW$ of radius $r$, with $\cC(2)=\cC$ as defined in
 Section~\ref{sec-wilson}. Observe that the entry regions $E_i$, for
 $i=1,2$ intersect the cylinders $\cC(r)$ for $1\leq r\leq R_*$ in 
   lines.
 Then for $1 < r \leq R_*$, the radius inequality and the compactness of $\cC(r)$  imply there exists a  lower bound   
\begin{equation}\label{eq-minradineq}
\delta(r) = \min \{  r(\sigma_i^{-1}(x)) \mid x \in \cC(r) \cap E_i ~, ~ i = 1,2 \} ~ \geq ~ r
\end{equation}
with equality only for $r=2$. Note that  $\delta(2)=2$ and $\delta(R_*)=3$, and 
   that $\delta(r)$ is an increasing function of $r$. The graph of $r \mapsto \delta(r)$ is illustrated in Figure~\ref{fig:r(delta)}.

Fix $r_0 > 2$,  set $r_1 = \delta(r_0)$ and   define  $r_{k} = \delta(r_{k -1})$ recursively, if  $r_{k -1} \leq R_*$. By the compactness of the regions $\cD_i$ there exists $N(r_0) \geq 0$ such that  $r_{k} \leq R_*$  for $k < N(r_0)$, but $r_{k} > R_*$   for $k = N(r_0)$. Note that  $r_{k}$ is not defined for $k > N(r_0)$.

We make an observation concerning  the hypothesis ``$x$ is not a secondary exit point'' which is often imposed in the statements of the following results. 
Consider the level function  $n_{x}(t)$   for the case where $x \in \mK$ is a secondary exit point.
For $\e > 0$   sufficiently small so that the orbit segment $\{\Phi_t(x) \mid 0 < t \leq \e\}$ contains no transition point, then we have that $n_{x}(t) = 0$ for $ 0 \leq t < \e$ by the definition   \eqref{def-level+}.  
Recall that  the radius function $r(x) = r(\tau^{-1}(x))$ on $\mK$ has a discontinuity at a secondary exit point, where  $r(x) = r(x')$ for $x' \in   D_i$ with $\tau(x') = x$.
 As $x$ is a secondary exit point, we have that $\rho_{x}(t) = r(\Phi_t(x)) = r(x'')$ for $0 < t < \e$ where $x'' = \sigma_i(x')$ as in \eqref{eq-radiusdef}. 
By the Radius Inequality (K8), we have $r(x'')  < r(x')$  unless $r(x') = r(x'') = 2$ where $x$ is a special point. 
  Thus,  if $x$ is a secondary exit point with $r(x)>2$, we have that  $\rho_x(t) < r(x)$ for $0 < t < \e$. 
  The assumption $x$ is not a secondary point eliminates this possibility. 
 
 We now give a sequence of results relating the level function with the properties of the orbits of $\Phi_t$.
First, we show that $N(r_0)$ gives a uniform bound on the level
function for  forward $\cK$-orbits. In what follows, when we say ``for
all $t$'' we mean ``for all $t$ such that $\Phi_t(x)$ is defined''.

\begin{lemma}\label{lem-maxstep}
Let $r_0 > 2$. Suppose that  $x \in \mK$ with $r(x) \geq r_0$ and    that $n_x(t) \geq 0$ for all $t \geq 0$.  
If $x$ is not a secondary exit point, then $n_x(t) \leq N(r_0)$ for all $t \geq 0$.
\end{lemma}
\proof
  Suppose there exists  $t_* > 0$ such that $n_x(t_*) > N(r_0)$. We show this leads to a contradiction. 
   Let $x_{\ell} = \Phi_{t_{\ell}}(x)$ with $0  \leq  t_0 < t_1 < \cdots < t_m < t_*$      be the transition points for  $\{\Phi_t(x) \mid 0 \leq t \leq t_*\}$.     
    
    Suppose there exists $k  > 1$  such that $n_x(t_{k-1} ) = n_x(t_0) = 0$.  
        Then the segment 
  $[x_0 , x_{k}]_{\cK}$ satisfies the hypotheses of
  Proposition~\ref{prop-shortcut4}, so we have   $x_0' \prec_{\cW} y_{k}'$  and thus $r(x_{k - 1}') = r(x) \geq r_0$. 
  The length of the $\cW$-orbit containing the segment $[x_0'  ,y_{k}']_{\cW}$ is bounded above by $L(\e)$ for $\e = r_0 -2$,  from 
  Lemma~\ref{lem-wilsonlengths}. Hence there can exist at most a finite number of values $k > 1$ for which $n_x(t_{k-1} ) = n_x(t_0) = 0$.
   
   Let $k  \geq 0$ be the largest index such that $n_x(t_{k-1} ) = n_x(t_0) = 0$.  Set $\ell_1 = k$. 
        Then the segment 
  $[x_0 , x_{\ell_1}]_{\cK}$ satisfies the hypotheses of
  Proposition~\ref{prop-shortcut4}, so we have   $x_0' \prec_{\cW} y_{\ell_1}'$  and thus $r(x_{\ell_1 - 1}') = r(x) \geq r_0$. 
 As $n_x(t) \geq 0$ for all $t \geq 0$, we must have  $1=n_x(t_{\ell_1}) > n_x(t_{\ell_1-1} )$  and thus  $x_{\ell_1}$ is a secondary entry point,  and  so $r(x_{\ell_1}') \geq r_1$. 
  
The   assumption  that  $\ell_1$   is maximal implies that $n_{x}(t) >  0$ for $t_{\ell_1} \leq t \leq t_*$.  
If  $n_x(t_{\ell_1}) = 1 < N(r_0)$, then     repeat the above process
to choose $\ell_2$   to be the largest index such that $n_x(t_{\ell_2-1})=1$. 

  Continuing in this way for $j < N(r_0)$,   choose the sequence $0
  = \ell_0 < \ell_1 < \cdots < \ell_n$ where $n \leq N(r_0)$ as long
  as possible, such that $\ell_j$ is the largest index satisfying 
  $n_x(t_{\ell_j -1}) = j-1$ and $n_x(t_{\ell_j}) = j$, for $0 < j \leq n$. 
  Then $r(x_{\ell_j}') \geq \delta(r(x_{\ell_j -1}')) \geq \delta(r_{j -1}) = r_{j}$ since
  $x_{\ell_j}$ is a secondary entry point  and $r(x_{\ell_j -1}') = r(x_{\ell_{j -1}}') \geq r_{j -1}$.

The assumption  $n_x(t_*) > N(r_0)$ implies  it is possible to choose an  index $\ell_n$ with $n_x(t_{\ell_n}) = N(r_0)$. 
By the choice of $N(r_0)$ we have $r_n > R_*$   and
$n_x(t)\geq N(r_0)$ for all $t\geq t_{\ell_n}$, and hence  the forward $\cK$-orbit  $\{\Phi_t(x) \mid t \geq t_{\ell_n}\}$ does not contain a transition point.   Thus, $\rho_x(t)$ is constant for $t \geq t_{\ell_n}$ and hence the maximum value of $n_x(t)$ is achieved at $t = t_{\ell_n}$, where the value is $N(r_0)$. This  contradicts the assumption that  $n_x(t_*) > N(r_0)$.
\endproof

 \begin{cor}\label{cor-maxstep}
 Suppose that $x \in \mK$ satisfies $\rho_x(t) \geq r_0>2$ for all $t \geq 0$. Then there exists $n(r_0)$ such that 
  $n(r_0)\leq n_x(t) \leq N(r_0)$ for all $t \geq 0$.  
 \end{cor}
 \proof
Suppose no such lower bound on $n_x(t)$ exists. Let  $x_{\ell} = \Phi_{t_{\ell}}(x)$  with $0 \leq t_0 < t_1 < \cdots < t_n < \cdots$ be    the transition points in $\cK$-orbit     $\{\Phi_t(x) \mid  t \geq 0\}$ and $\{\ell_k \mid k =1,2, \ldots\}$ be an increasing subsequence such that $n_x(t_{\ell_k}) = - k$.

Consider the   flow $\Xi^k_t(x) = \Phi_{t_{\ell_k} -t}(x)$ on $\mK$, then $\Xi^k_0(x)   = x_{\ell_k}$ and $\Xi^k_{t_{\ell_k}}(x) = x$. Moreover, $r(\Xi^k_t(x)) \geq r_0$ for all $0 \leq t \leq t_{\ell_k}$. Then  the level  with respect to the flow $\Xi^k_t$ starting at $x_{\ell_k}$    increases to $k$ for $0 \leq t \leq t_{\ell_k}$. We can then apply the method of proof of Lemma~\ref{lem-maxstep} to conclude that the level is bounded above by $N(r_0)$, independent of the choice of $k$, which contradicts the above.  
 \endproof
 
These results combine to yield   an   ``escape'' result for $\cK$-orbits. 
   
  \begin{prop} \label{prop-r>2} 
Let $x \in \mK$ satisfy $r_0 = r(x) > 2$, assume that $x$ is not a secondary exit point and suppose that $n_x(t) \geq 0$ for all $t \geq 0$. Then the forward $\cK$-orbit of $x$ is not trapped.  
  \end{prop}
  \proof 
  Suppose that the forward orbit of $x$ is trapped, then there exists      $0 \leq t_0 < t_1 < \cdots < t_n < \cdots$ so that the points $x_{\ell} = \Phi_{t_{\ell}}(x)$ are the forward transition points for the $\cK$-orbit 
     $\{\Phi_t(x) \mid  t \geq 0\}$.    By  Lemma~\ref{lem-maxstep}  and our assumption, we have   that $0 \leq n_x(t_{\ell}) \leq N(r_0)$ for all $\ell \geq 0$.

     Let $0 \leq n_0 < N(r_0)$ be the least integer such that there exists    
     $0 \leq \ell_0 < \ell_1 < \cdots < \ell_k < \cdots$ such that  $n_x(t_{\ell_j}) = n_0$. That is, 
     $\ds n_0 = \liminf_{\ell > 0} ~ n_x(t_{\ell})$. 
          Since $n_0$ is the least such integer, there exists $k \geq 0$ such that  $n_x(t) \geq n_0$ for all $t \geq t_{\ell_k}$. 
      Then for each integer $\alpha \geq 1$, the segment $[x_{\ell_k}, x_{\ell_{k+\alpha} + 1}]_{\cK}$ 
     satisfies the hypotheses of Proposition~\ref{prop-shortcut4}, so
     we have   $x_{\ell_k}' \prec_{\cW} y_{\ell_{k+\alpha} + 1}'$.  
Then by Lemma~\ref{lem-wilsonlengths} and  Corollary~\ref{cor-segmentlengths}, considering the lengths of the lifted  $\cW$-orbit segments $[x_{\ell_k}', y_{\ell_{k+\alpha} + 1}']_{\cW}$  yields the   estimate  
$\ds  d_{min}\cdot (\ell_{k+\alpha}-\ell_k)\leq L(r_0-2)$ for all $\alpha$. 
   However,   we can choose $\ell_{k+\alpha}$ arbitrarily large and hence also 
      $(\ell_{k+\alpha} - \ell_k)$, which yields a contradiction. 
         \endproof

  \begin{cor} \label{cor-r>2} 
Let $x \in \mK$ satisfy $\rho_x(t) > 2$ for all $t \geq 0$, and suppose that  $\ds \liminf_{t > 0} ~ n_x(t) > -\infty$. 
Then the forward $\cK$-orbit of $x$ is not trapped.
  \end{cor}
\proof
Let $t_0 \geq 0$ be such that $x_0 = x_{t_0} = \Phi_{t_0}(x)$ is a transition point with 
$\ds n_x(t_0) = \liminf_{t > 0} ~ n_x(t)$. Then $n_{x_0}(t) \geq 0$
for all $t \geq t_0$ and we have  $r(x_0) >   2$. Then apply Proposition~\ref{prop-r>2} for $x = x_0$.
\endproof

\begin{figure}[!htbp]
\centering
{\includegraphics[width=80mm]{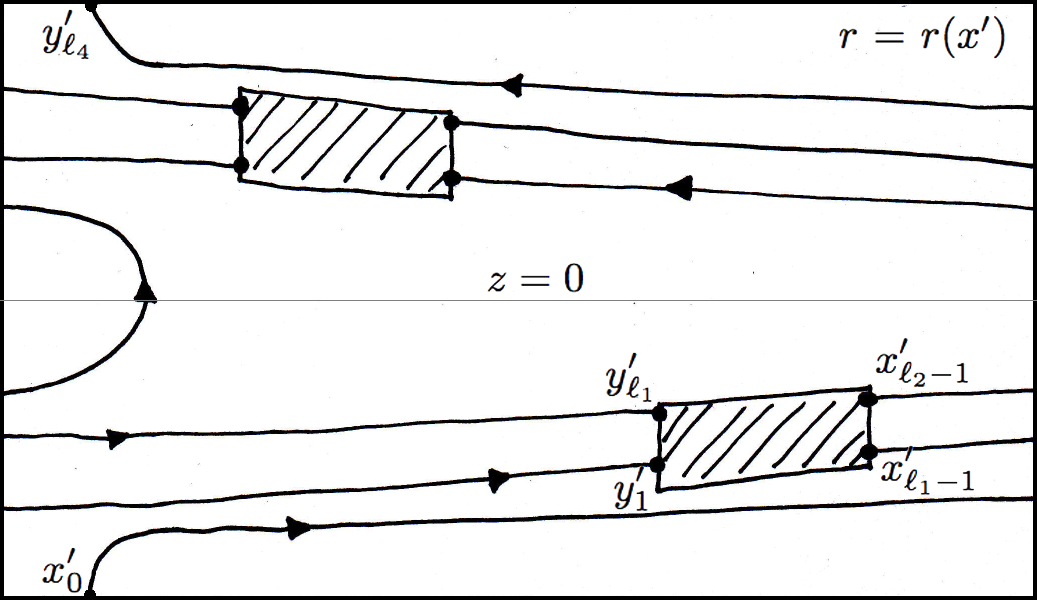}}
 \caption[{The arcs of the lifted $\cK$-segment $[x,y]$} ]{\label{fig:Wsegments} The arcs of the lifted $\cK$-segment $[x,y]$ on the cylinder $r=r(x')$ with $x=x_0$ and $y=x_{\ell_4}$   }
 \vspace{-6pt}
\end{figure}

 \medskip
 Ghys remarks in \cite[Lemme, page 300]{Ghys1995},   that for an entry point $x \in \mK$ which escapes to its facing exit point  $y$, so that  $x \equiv y$, ``the $\cK$-orbit of $x$ contains the image under $\tau$ of all the $\cW$-arcs that lie in $\mW^\prime$ between $x'$ and $y'$, where $y'$ is the   exit point such that $x^\prime\equiv y^\prime$.''  
  Propositions~\ref{prop-min}  and \ref{prop-propexist} below give a precise formulation of  this assertion.   Figure~\ref{fig:Wsegments} illustrates these lifts.

   \begin{prop}  \label{prop-min} 
  Let $x'\in \mW$ such that $x =
\tau(x')$ is a primary entry point in $\mK$ with $r(x) > 2$. Then the
$\cK$-orbit of $x$ escapes from $\mK$ at some primary exit point $w$, we have that $\rho_x(t)\geq r(x)$ for all $t$
and the collection of lifts of the $\cW$-arcs in $[x,w]_{\cK}$
contains all the $\cW$-arcs of the $\cW$-orbit of $x'$ that are in
$\wmW$.
   \end{prop}
   \proof

Suppose that the forward orbit of $x$ is trapped, then there exists      $0 = t_0 < t_1 < \cdots < t_n < \cdots$ with $x_{\ell} = \Phi_{t_{\ell}}(x)$ the transition points for the $\cK$-orbit     $\{\Phi_t(x_0) \mid  t \geq 0\}$.   
          We claim that $n_{x_0}(t) \geq 0$ for all $t \geq 0$, then the result follows from   Proposition~\ref{prop-r>2} and Proposition~\ref{prop-ee}.

  Suppose $n_{x_0}(t) < 0$    for some $t > 0$, then there is a least $k > 0$ such that $n_{x_0}(t_{k}) = -1$, and thus 
  $n_{x_0}(t_{k-1}) = 0$ and $x_{k}$ is an exit point. 
  Proposition~\ref{prop-shortcut4} then implies that $x_0' \prec_{\cW} x_k'$ and  $x_0 \equiv x_k$. 
  Since $x_0$ is a primary entry point, we must have that $x_k$ is a
  primary exit point, which is a contradiction. Thus, $n_{x_0}(t) \geq 0$ for all $t \geq 0$. 

Next, note that as ${x_0}$ is a primary entry point,    $[x_0',y_1']_{\cK} \subset \wmW$ is the initial $\cW$-arc in the
intersection of the $\cW$-orbit of  $x'$ with $\wmW$. Then 
let  $\ell_1, \ell_2, \ldots, \ell_m$ be the collection of  all the indices such that
$n_{x_0}(t_{\ell_k-1})=0$ for all $k=1,2,\ldots,
m$. Proposition~\ref{prop-shortcut4} implies that $x_0' \prec_\cW y_{\ell_1}'$ and $x_1\equiv x_{\ell_1-1}$, hence  the subsequent    lift 
$[x_{\ell_1-1}' ,y_{\ell_1}']_{\cK} \subset \wmW$ is the second arc in the
intersection of $\cW$-orbit of
  $x_0'$ with $\wmW$. Continuing by induction, we obtain that all the $\cW$-arcs for $x_0'$   lying in $\wmW$  are contained in the $\cK$-orbit of $x_0$.
\endproof

There is  an important consequence of this     result.
  \begin{cor} \label{cor-r>2exit} 
Let $z \in \mK$ satisfy $r(z) > 2$ and $n_z(t) \geq 0$
for all $t \geq 0$. Then the forward $\cK$-orbit of $z$ exits $\mK$ in
a point $y = \Phi_{T}(z)$ at which  $n_z(T) =  0$. Furthermore, $z$ is contained in a bounded orbit of $\Phi_t$, starting at a primary entry point $x$ of $\mK$, where $x \equiv y$. 
  \end{cor}
 
  \proof
 Proposition~\ref{prop-r>2}  implies that the flow $\Phi_t$ has a primary exit point $y$ at time $T$.
Note that $n_z(T) \geq 0$ implies $r(y) \geq r(z) > 2$ by Proposition~\ref{prop-shortcut4}. 

 Then apply Proposition~\ref{prop-min} to the reverse flow $\Xi^T_t(\xi) \equiv \Phi_{T -t}(\xi)$ 
  to deduce that there
exists $t_* \geq T$ such that $x =\Xi^T_{t_*}(y)$ is a primary \emph{exit} point for  $\Xi^T_t$ and hence a primary \emph{entry} point for $\Phi_t$.   
Moreover, the proof of Proposition~\ref{prop-min} shows that 
$n_{x}(t)\geq 0$ for all $t \geq 0$. Note that $y = \Phi_{T}(z) = \Phi_{t_*}(x)$ so we also have    $n_{x}(t_*) = 0$.

 The fact $n_z(T)=0$ follows form the equality
  $n_z(t)=n_x(t_*-T+t)-n_x(t_*-T)$ for all $t\geq 0$. Thus
  $n_z(T)=n_x(t_*)-n_x(t_*-T)=-n_x(t_*-T)\leq 0$. Since $T\geq 0$ and
  $n_z(T)\geq 0$ by hypothesis, $n_z(T)=0$ as claimed. 
  \endproof

  The conclusion of Corollary~\ref{cor-r>2exit}  has a simple
  geometric interpretation. Given a \emph{primary} entry point  $x \in
  \mK$ with $r(x) > 2$, the points  $z = \Phi_t(x)$ along the
  $\cK$-orbit of $x$ for which $n_{x}(t) = 0$  are \emph{exactly} those
  points on the Wilson flow of $x'$   which lie in
  $\wmW$. Proposition~\ref{prop-min}  implies that each such $\cW$-arc
  of the $\cW$-orbit of   $x'$ is contained in the $\cK$-orbit of
  $x$. Observe that by taking $z\in \mK$ with $r(z)>2$ and $n_z(t) \geq 0$
for all $t \geq 0$, we assume that $z'\in\mW$ such that $\tau(z')=z$
belongs to a $\cW$-orbit whose entry point is in
$\partial_h^-\mW-(L_1^1\cup L_2^-)$.

 Note that   Proposition~\ref{prop-min}  implies that the endpoints of
 the $\cK$-orbit through $z$ are facing, so that  $x\equiv y$.   
 The following result is valid for secondary entry points, and  is used repeatedly in subsequent arguments.

\begin{prop} \label{prop-propexist}
Let $x \in \mK$ be a secondary entry point with $r(x) > 2$, let $y$ be a secondary exit point  and suppose that $x \equiv y$. 
Then $x \prec_{\cK} y$ and the collection of lifts of the $\cW$-arcs in $[x,y]_{\cK}$
contains all the $\cW$-arcs of the $\cW$-orbit of $x'$ that are in
$\wmW$, where $\tau(x')=x$. Hence, if $[\xi',x']_{\cW}$ and $[y',z']_{\cW}$ are $\cW$-arcs, then their images under $\tau$  are both contained in the $\cK$-orbit through $x$.
\end{prop}
\proof
Corollary~\ref{cor-r>2exit}  implies there exists a maximal $t_* > 0$ such that  $n_x(t) \geq 0$ for $0 \leq t < t_*$ as $x$ is assumed to be a secondary entry point. 
Then $x_* = \Phi_{t_*}(x)$ must be an exit point. 

Let    $0 = t_0 < t_1 < \cdots < t_n = t_*$ with $x_{\ell} = \Phi_{t_{\ell}}(x)$ be the transition points for the $\cK$-orbit     $\{\Phi_t(x) \mid  0 \leq t \leq t_n\}$.  
Then $n_x(t) \geq 0$ for all $0 \leq t < t_n$, and $n_x(t_{n-1}) = 0$.
Hence $x_{n-1}$ and $x_n$ are exit points. By
Proposition~\ref{prop-shortcut4}, $x_0'\prec_\cW y_n'$ and $x_0\equiv
x_n$. The point facing $x$ is unique, so we must have
$y=x_n=\tau(y_n')$ and $x\prec_\cK y$. 

The proof that the segment $[x,y]_\cK$ contains all the $\cW$-arcs on the $\cW$-orbit of $x'$ that are in $\wmW$ follows as in the proof of Proposition~\ref{prop-min}.
The first transition point $x_1$ must be a secondary entry as $n_x(t_1) \geq 0$, thus $n_x(t_1) =1$. 
If $n_x(t_2) = 0$ then $x_2$ is a secondary exit, so by Proposition~\ref{prop-shortcut4} we have that $x' \prec_{\cW} x_3'$ and so $x_1 \equiv x_2$. 
Thus the lift of the $\cW$-arc $[x_2,x_3]_\cK$ is contained in the
$\cW$-orbit of $x'$.
Otherwise, $n_x(t_2) = 2$ and we continue the analysis of cases, until we obtain the first instance where $n_x(t_{\ell -1}) = 0$. 
We can then apply Proposition~\ref{prop-shortcut4} to conclude that 
 $x_0'  \prec_{\cW} y_{\ell}'$. Thus the lift of the $\cW$-arc
 $[x_{\ell-1},x_\ell]_\cK$ is contained in the $\cW$-orbit of $x'$. The rest follows as in the
 proof of Proposition~\ref{prop-min}.
\endproof

  \medskip
  
Let $x' \in \partial_h^- \mW$ with $r(x') > 2$ and let  $y' \in \partial_h^+ \mW$ be the unique facing point.  Let $x = \tau(x')$ and $y = \tau(y')$. 
If  $x$ is a primary entry point, 
then  $x \prec_{\cK} y$ by     Proposition~\ref{prop-min}.    
If $x$ is a secondary entry point, then $x \prec_{\cK} y$ by    Proposition~\ref{prop-propexist}. 
Thus, in both cases there exists  $T_x > 0$ such that $\Phi_{T_x}(x) = y$.
 
 \begin{cor}\label{cor-lengths}
 The function $x' \mapsto T_{x'}$ is well-defined   and continuous  for $x' \in \partial_h^- \mW$ with $r(x') > 2$.
 Moreover, $T_{x'}$ tends to infinity as $r(x')$ tends to $2$. \hfill $\Box$
 \end{cor}

  \bigskip
 
\section{Trapped and infinite orbits}\label{sec-trapped} 

 In this section, we begin the analysis of the trapped orbits of the flow $\Phi_t$.
We first  consider the   orbits of points $x = \tau(x')$ with   $r(x') =2$.  The following result  is a general formulation of 
an  observation of Kuperberg, which provided the original insight
leading to  her construction \cite{Kuperberg2011} and was sketched in
Figure~\ref{fig:cW-decomposition}.

  \begin{prop} \label{prop-r=2any}
Let  $x'  \in \mW'$ with $r(x') = 2$ and set $x = \tau(x')$. Then the $\cK$-orbit of $x$ is trapped in either forward or backward time, or both. Moreover,  the collection of $\cW$-arcs for the $\cK$-orbit  of $x$
contains every $\cW$-arc lying in $\wmW$ of the $\cW$-orbit of
$x'$ and  there is  a subsequence 
$\{ x_{\ell_j}' \mid j=1,2, \ldots \}$ of transition points  with 
$r(x_{\ell_j}') = 2$ such that $\{ x_{\ell_j} = \tau(x_{\ell_j}') \mid j = 1,2, \ldots \}$  converges to a special point $p_i^-$ for $i=1,2$. 
\end{prop}
\proof
Assume without loss of generality that $x$ is not   a transition point and 
let   $x_{\ell} = \Phi_{t_{\ell}}(x)$, for  $0 < t_0 < t_1 < \cdots < t_n < \cdots$,    be  the transition points for the $\cK$-orbit    $\{\Phi_t(x) \mid  t \geq 0\}$.
There are  five cases to consider: 
\begin{enumerate}
\item $\ds z(x') = - 1 $ ~; 
\item $\ds  z(x') =   1 $ ~ ;
\item $\ds -2 \leq z(x')<-1 $ ~ ;
\item  $\ds -1< z(x')<1 $ ~;
\item  $\ds 1< z(x')\leq 2 $ ~ ;
\end{enumerate}
Consider first the case of the forward $\cK$-orbit of a point  $x = \tau(x')$ with $x' \in \mW'$,  $r(x') =2$ and $z(x')=-1$. 
Then $x_0 \in E_1$ is the special point $p_1^-$ of \eqref{eq-special} and   for the $\cW$-arc $[x_0', y_1']_{\cW} \subset \mW'$ we have  $x_0' \in L_1^-  \subset \partial_h^- \mW$.
Then $r(x_0') = 2$ by Condition (K7) so that $r(y_1') = 2$ also. The
$\cW$-arc $[x_0', y_1']_{\cW}$ flows upward from   $z(x_0') = -2$
until it intersects at $y_1' \in \cL_1^-$ with  $z(y_1') < -1$, as in Figure~\ref{fig:cW-decomposition}. 
Let $x_1' \in L_1^-$ satisfy $\tau(x_1') = \tau(y_1')$, then   the Radius Inequality implies  that $r(x_1') > 2$. 
Let $\ox_1 \in S_1$ be the facing point to $x_1$. Then by  Proposition~\ref{prop-propexist},   we have $x_1 \prec_{\cK}  \ox_1$ and  so $\ox_1 = x_{\ell_1-1}$ for some $\ell_1 >2$. Then $r(x_{\ell_1-1}') =2$,
$z(y_1') < z(x_{\ell_1-1}') < -1$ and $[x_{\ell_1-1}', y_{\ell_1}']_{\cW}$ is the subsequent $\cW$-arc. Thus $x_{\ell_1 } =\tau(y_{\ell_1 }')$ must be a secondary entry point again and we can repeat this argument inductively to obtain 
  a subsequence $\{ y_{\ell_i}' \mid i=1,2, \ldots \}$ in $\cL_1^-$ with 
$x_{\ell_i} =\tau(y_{\ell_i}')$ in $E_1$ converging to $p_1^-$, so that the forward orbit of $x$ is trapped. 
The assertion about the collection of lifts of $\cW$-arcs in $\wmW$ follows  from Proposition~\ref{prop-propexist}.
A similar
analysis for the backward orbit of $x$  yields a subsequence of transition
points converging to $p_2^-$.

  For the  case of $x = \tau(x')$ with $x' \in \wmW$,  $r(x') =2$ and
  $z(x')=1$,    note that the first forward transition point $x_0 \in
  E_2$ is the special point $p_2^-$ of \eqref{eq-special}  and so  for
  the $\cW$-arc $[x_0', y_1']_{\cW} \subset \mW'$ we have  $x_0' \in
  L_2^- \subset \partial_h^- \mW'$. Then again $r(x_0') = 2$ by
  Condition (K7) and the rest of the analysis proceeds similarly to the previous case. Note that the forward orbits of  a point  $x = \tau(x')$ in the two cases where $r(x') =2$ and $z(x')= \pm 1$ limit to  $p_1^-$, while their backward orbits tend to $p_2^-$.

There are three remaining  cases: either $-2 \leq z(x') < -1$, $-1 < z(x') < 1$, or $1 < z(x') \leq 2$. 
All three cases proceed in a manner analogous  to  the two cases above. 

If  $-2 < z(x') < -1$, then  the forward $\cW$-orbit of $x'$ is asymptotic to the periodic orbit
$\cO_1$, so in particular is trapped for the Wilson Plug and each
time the $\cW$-orbit enters an insertion region $\cD_i \subset \mW$,  it
subsequently exits through the facing point. 
Thus,   the  forward $\cK$-orbit of the point  $x = \tau(x')$
has first transition point $x_0 \in E_1$. As $x'$ is not in a
periodic orbit, we have $x_0' = \tau^{-1}(x_0) \in L_1^-$ with $r(x_0')
> 2$.  Let $x_{\ell_1-1}\in S_1$ be the facing point, by
Proposition~\ref{prop-propexist} $x_0\prec_\cK x_{\ell_1-1}$. We then proceed as above, to obtain 
  a subsequence $\{ x_{\ell_i} \mid i=0,1,2, \ldots \}$ in $E_1$ with
  $x_{\ell_i}$ converging to $p_1^-$.  
  
  Note that the backward
  $\cW$-orbit of $x'$ is  not   trapped and so exits $\mW$ at a point
  $y'$. The backward $\cK$-orbit of $x$ contains the point
  $y=\tau(y')$. If   $y \not\in L_i^-$ for $i=1,2$, then   the backward $\cK$-orbit of
  $x$ is not trapped.  If $y \in L_i^-$ for $i=1,2$, then the process continues in the backward direction.

For the case $-1 < z(x') < 1$, the forward $\cW$-orbit of $x'$ in $\mW$
is asymptotic to the periodic orbit $\cO_2$,   and each time the
$\cW$-orbit enters an insertion region $\cD_i$ it subsequently exits
the same region. The   rest of
the analysis proceeds similarly to the   case when $-2 \leq z(x') <
-1$, yielding a subsequence of secondary entry points for the forward $\cK$-orbit which converge to $p_2^-$. 
The backward  $\cW$-orbit of $x'$ in $\mW$ is asymptotic to the periodic orbit $\cO_1$, yielding in the same manner  a subsequence of secondary entry points for the reverse  $\cK$-orbit converging to $p_1^-$.

  The last case, for $1 < z(x') \leq 2$, reduces to the   case for $-2 \leq z(x') < -1$ by reversing the flow $\Phi_t$. Thus, the forward $\cK$-orbit of $x$ may escape through $\partial_h^+ \mK$ or may be trapped for all forward time. The backward $\cW$-orbit of $x'$ always converges to $\cO_2$ and so yields a subsequence of secondary entry points for the reverse flow converging to $p_2^-$.   
 \endproof

We observe two additional     consequences of   Proposition~\ref{prop-r=2any} and its proof.
\begin{cor} \label{cor-r=2special}
 Let $x \in  \mK$ and suppose that  either $x$ 
  is a primary entry point with
  $r(x) =2$, or $x = \tau(x')$  where $x' \in \cO_i \cap \mW'$ for $i=1,2$. Then the level function based at $x$ satisfies
  $n_{x}(t) \geq 0$ for all $t \geq 0$. Moreover, if $\Phi_t(x)$ is
  not a transition point for  $t\geq 0$, then $n_{x}(t) = 0$ if and only if $\rho_x(t)= 2$. \hfill $\Box$
\end{cor}

   \medskip

\begin{cor}\label{cor-A_0infinite}
Let $x\in \mK$ and suppose   there exists $t_*\in \mR$ such that
$y=\Phi_{t_*}(x)$ is a secondary entry point with $\rho_x(t_*)=2$. Then  for
all $t>t_*$ we have that $\rho_y(t)\geq 2$. Moreover, the forward orbit of $x$ 
contains an infinite sequence of secondary entry points that limit to
$p_i^-\in E_i^-$ for $i=1$ or $i=2$.  \hfill $\Box$
\end{cor}

\medskip

We next consider the   $\cK$-orbits of points $x \in \mK$ for which $r(x) < 2$. 
We recall results by Ghys \cite{Ghys1995} and Matsumoto  \cite{Matsumoto1995}, which give 
  conditions such that the orbit of $x$ is trapped in forward time and
  that   the closure of the $\cK$-orbit of $x$ contains  a special
  point. Hence the   closure of the $\cK$-orbit of $x$ contains the closures of the orbits of both special  points $p_1^-$ and $p_2^-$.

\begin{figure}[!htbp]
\centering
{\includegraphics[width=60mm]{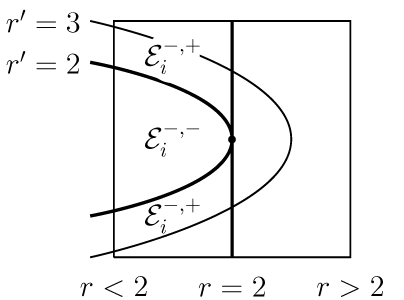}}
 \caption{\label{fig:Upsilon} The regions $\cE_i^{-,-}$ and $\cE_i^{-,+}$ of $\cL_i^-$}
 \vspace{-6pt}
\end{figure}

  Fix $i = 1,2$ and consider the restriction of the insertion map
  $\sigma_i  \colon L_i^-  \to \cD_i \subset \mW$. Express this map in
  polar coordinates $(r',\theta')$ on the domain $L_i^-$ and
  $(r,\theta,z)$ on the image. The image under $\sigma_i$ of the curve $\{r' =2\} \cap L_i^-$   is a ``parabolic curve'' $\Upsilon$ which is tangent to the vertical line $\{r = 2\}$.
For each $i =1,2$, we define two regions  (see Figure~\ref{fig:Upsilon})  contained in the image of the region  $\cL_i^- \cap \{r < 2\}$:
\begin{itemize}
\item  $\cE_i^{-,-} = \sigma_i(\{r' < 2\} \cap L_i^-) \cap \{r < 2\}$    with outer boundary $\Upsilon$
\item $\cE_i^{-,+} = \sigma_i(\{r' > 2\} \cap L_i^-) \cup \{r < 2\}$
  with inner boundary  $\Upsilon$. 
\end{itemize}

  We make two basic observations, as used in   the proof of \cite[Proposition~7.1]{Matsumoto1995}. 
  First,  the tangency of $\Upsilon$ with   the vertical line $\{r = 2\}$ and Taylor's Theorem implies that 
for $\delta > 0$ sufficiently small, the   vertical line segment 
$\cE_i^{-,-} \cap \{r = 2-\delta\}$ has length at least $C' \cdot \sqrt{\delta}$ for some fixed $C' > 0$.

The second observation     is that the Wilson vector field $\cW$ is
horizontal along the periodic orbits $\cO_i$ as the vertical component
$g(r, \theta, z)  \frac{\partial}{\partial  z}$ in
\eqref{eq-wilsonvector} vanishes at the points $(2, \theta, \pm 1)$. Thus, there exists    $C'' > 0$ so that for $\delta > 0$ sufficiently small, the flow $\Psi_t(x')$ of a point $x' \in \mW$ with $r(x') = 2 - \delta$   intersects the face $\cL_i^-$ in a sequence of points whose vertical spacing near the planes $\{z = \pm 1\}$ is bounded above by    $C'' \cdot \delta$.  A more precise estimate can be obtained using the method of proof  for Lemma~\ref{lem-density2}, but these approximate results   suffice for the proof of Proposition~\ref{prop-wander} below.

\begin{defn}\label{def-matsumoto}
We say that  $\delta_M > 0$ is a Matsumoto constant if for all $ 0 < \delta \leq \delta_M$ and $C', C''$   as above,   then  
\begin{equation}\label{eq-matsumoto}
C'' \cdot \delta < C' \cdot \sqrt{\delta}
\end{equation}
The  \emph{Matsumoto region} for $\delta_M$  is the set $\cU(\delta_M) = \tau(\{2-\delta_M < r < 2\}) \subset    \mK$. 
\end{defn}

 \begin{prop}\label{prop-wander} \cite[Proposition~7.1]{Matsumoto1995}
 Let $\delta_M > 0$  be a Matsumoto constant, and       $x = \tau(x') \in \cU(\delta_M)$.  
\begin{enumerate}
\item If $z(x') = -2$, then the \emph{forward orbit} of $x$ is trapped. Moreover,    there exists a subsequence of   transition times, 
$0 < t_0 < t_{\ell_1} < \cdots < t_{\ell_i} < \cdots $, such that
$r(x_{\ell_i}')  \leq 2$     for $i \geq 0$  with $\ds \lim_{i \to  \infty} r(x_{\ell_i}')  =  2$.  
\item If $z(x') = 2$, then    the backward orbit of $x$ is trapped. Moreover, there exists a subsequence of   transition times, 
$0 > t_0 >  t_{\ell_1} >  \cdots > t_{\ell_i} > \cdots $, such that  $r(x_{\ell_i}')  \leq 2$     for $i \geq 0$ with $\ds \lim_{i \to  \infty} r(x_{\ell_i}')  =  2$.
\item If $z(x') = 0$, then the orbit of $x$ is   infinite and there exists a bi-infinite subsequence of   transition times 
  such that  $r(x_{\ell_i}')  \leq 2$     for all $i$, with $\ds \lim_{i \to  \pm \infty} r(x_{\ell_i}')  =  2$.
\end{enumerate}
\end{prop}

\begin{figure}[!htbp]
\centering
{\includegraphics[width=70mm]{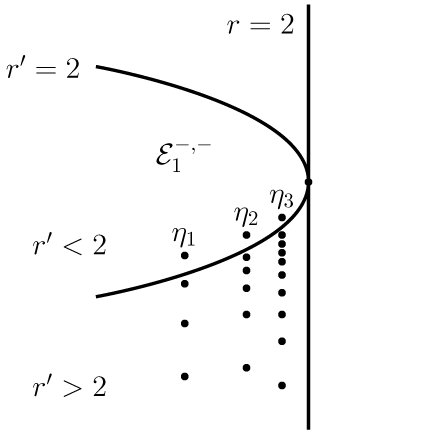}}
 \caption{\label{fig:trapping} Trapped orbits intersecting $\cE_1^{-,-}$ infinitely often}
\vspace{-6pt}
\end{figure}

\proof
We consider first the case $z(x') = -2$. The case $z(x') =2$ follows analogously, by reversing the flow $\Phi_t$. Figure~\ref{fig:trapping} illustrates the following construction.

Given  $x_0' = (r_0,\theta, -2)$ which satisfies $2-\delta_M < r_0 < 2$, the $\cW$-orbit of $x_0'$
has increasing $z$-coordinate in $\mW$ and intersects the face
$\cL_1^-$   in the set $\cE_1^{-,+}$, possibly repeatedly,  until it
first intersects either the curve $\Upsilon$ or the region $\cE_1^{-,-}$. 
If $\Psi_{t_*}(x_0') \in \Upsilon$ for some $t_* > 0$, then $r(\tau(\Psi_{t_*}(x_0'))) = 2$ and the analysis of the forward orbit reduces to the case of Corollary~\ref{cor-A_0infinite}.

Otherwise, the assumption that $2-\delta_M < r_0 < 2$ implies the initial segment of the
$\cK$-orbit of $x_0 = \tau(x_0')$     contains a sequence of  $\cW$-arcs which
coincide with   the initial $\cW$-arcs of the $\cW$-orbit of $x_0'$, up until the initial point of the $\cW$-arc that is contained in
$\cE_1^{-,-}$. Then the $\cK$-orbit jumps to a radius $r_1$
with $\delta_M < r_0 < r_1 < 2$ and this pattern repeats itself.
We   make this   sketch of   proof   precise.

  Let    $x_{\ell} = \Phi_{t_{\ell}}(x_0)$   with   $0 = t_0 < t_1 < \cdots < t_n < \cdots$     be the transition points for the
 $\cK$-orbit    $\{\Phi_t(x_0) \mid  t \geq 0\}$. Denote the
 corresponding lifts to $\mW$ of the $\cW$-arcs  in this forward orbit  by $[x_{\ell}', y_{\ell +1}']_{\cW}$.  As remarked above,  we can assume that   $r(x_{\ell}') \ne 2$ for all $\ell$.
The orbit of $x_0$ is then described by an iterative process, as follows.

  Note that  $x_0' \in \partial_h^- \mW$ implies that  $y_1' \in \cL_1^-$,  so $n_{x_0}(t_1) = 1$.  
 If   $x_1 \in    \cE_1^{-,-}$,  set $\lambda_1 = 1$,  $\eta_{1} =
 x_1$ and $r_1 = r(\eta_1)$. Then  note that $r(y_1') = r(x_0') = r_0$
 while $r_0 < r(\eta_1') < 2$ by the definition of the region  $\cE_1^{-,-}$.

Otherwise, we have $x_1 \in    \cE_1^{-,+}$ and so $r(y_1') = r_0$
while $r(x_1') > 2$ by definition. Proposition~\ref{prop-propexist}
implies that there is a least $\ell_1>1$ such
that $n_{x_0}(t_{\ell_1-1})=0$, so that $x_0$ and $x_{\ell_1}$ satisfy
$x_0'\prec_\cW x_{\ell_1-1}'$ by Proposition~\ref{prop-shortcut4} and $r(x_{\ell_1-1}')=r_0$. Then $x_{\ell_1}$ is a secondary entry point.
 
If $x_{\ell_1} \in \cE_i^{-,-}$ then set $\lambda_1 = \ell_1$,  $\eta_{1} = x_{\ell_1}$ and $r_1 = r(\eta_1)$.
Otherwise, $x_{\ell_1} \in \cE_i^{-,+}$ and so  $r(x_{\ell_1}') > 2$, and we repeat this process.
By Corollary~\ref{cor-orbitnumber} and the choice of $\delta$, this inductive process can be repeated at most a finite number of times, 
until we obtain  $\{0 = \ell_0 < \ell_1 < \cdots < \ell_{k_0}\}$ 
with $r(x_{\ell_i -1}')  =r_0$, $r(x_{\ell_i}')  > 2$ for $1 \leq
i < k_0$ and $r_0 < r(x_{\ell_{k_0}}')  < 2$. Thus
$x_{\ell_{k_0}}\in \cE_i^{-,-}$. Then set $\lambda_1 = \ell_{k_0}$, $\eta_1 = x_{\lambda_1}$ and $r_1 = r(\eta_1)$. 

We thus obtain the  secondary entry points $\{y_{\ell_i}' \mid 0 \leq i < k_0\} \subset \{ r = r_0\} \cap \cL_1^-$ all
contained in a vertical line $r = r_0$, plus a new point 
$\eta_1 \in E_1$ with lifts $y_{\lambda_1}'\in \{r=r_0\}\cap \cL_1^-$ and
$x_{\lambda_1}'=\eta_1'$ and $r_0 < r_1 = r(\eta_1') < 2$.  We can then repeat the above process by considering the $\cK$-orbit of $\eta_1$, 
which yields an ascending finite sequence of secondary entry points $\{x_{\ell_i}\}$ with $\{y_{\ell_i}'\}$ 
 contained in the vertical line segment $\cE_1^{-,+} \cap \{r = r_1\}$, until they reach the point $\eta_2 = x_{\lambda_2} \in \cE_1^{-,-}$ with radius 
$r_0 < r_1 < r_2 = r(\eta_2) < 2$.

 This   process then continues recursively,  so that the sequence of 
 transition points $\{x_i \mid i \geq 0\}$ contains  an infinite
 collection of finite subsequences lying on the lines $\{r=r_j \} \cap
 \cL_1^-$ with $r_0 < r_1 < r_2 < \cdots<2$, where the initial points
 $\eta_j$ for each such     ``finite stack of points'' is defined as the
 secondary entry points where the sequence transitions from the region
 $\cE_1^{-,+}$ to the region $\cE_1^{-,-}$.  Moreover, the sequence of
 radii $r_j \to 2$ by  the Radius Inequality.
 Finally, note that the forward $\cK$-orbit of $x_0$ satisfies $r(\eta_{\ell}) < r(\eta_{\ell +1})$ so the orbit cannot be recurrent.

 As noted above, if $x_0' = (r_0,\theta, 2)$ with $2-\delta_M < r_0 < 2$,  then we obtain a similar conclusion using the reverse time flow and the $z$-symmetry of the flow $\cW$ on $\mW$.
 
 If $x_0' = (r_0,\theta, 0)$ with $2-\delta_M < r_0 < 2$, then the forward $\cW$-orbit of $x_0'$  
has increasing $z$-coordinate in $\mW$, and intersects the face
$\cL_2^-$ in the set $\cE_2^{-,+}$,   possibly repeatedly, until it
first intersects either the curve $\Upsilon$, or the region $\cE_2^{-,-}$.  The analysis then proceeds as above. The backward $\cW$-orbit of $x_0'$  
has decreasing $z$-coordinate in $\mW$, and  we obtain a similar
conclusion using the $z$-symmetry of the Wilson flow. 
 \endproof

Suppose that  $r(x) < 2$ and assume there exists $t_0 > 0$ such that $\rho_x(t_0) = 2$.   Then by Corollary~\ref{cor-A_0infinite}, the forward orbit $\{\Phi_t(x) \mid t \geq t_0\}$  is trapped in the region $\{r \geq 2\}$.

For   $r(x) < 2$, it remains to consider the case where      $\rho_x(t) \ne 2$ for all $t \geq 0$.    The proof  of the next   result  is  adapted from  the proof of the Th\'eor\`eme  in  \cite[page 301]{Ghys1995}. 
\begin{prop}\label{prop-cluster-}
Let $x\in \mK$ satisfy $r(x) < 2$ and assume that   its forward orbit is trapped and satisfies  $\rho_x(t) \ne 2$ for all $t \geq 0$.
Then  there exists an infinite subsequence of transition points 
$\{ x_{\ell_i} \mid i=1,2, \ldots \}$  with  $r(x_{\ell_i}') <  2$ such that $\ds \lim_{i\to \infty} ~ r(x_{\ell_i}') = 2$.
\end{prop}
\proof
Let   $x_{\ell} =
\Phi_{t_{\ell}}(x)$ with  $0 < t_0 < t_1 < \cdots < t_n < \cdots$     be  the transition points for   $\{\Phi_t(x) \mid
t \geq 0\}$, with  associated  lifts of $\cW$-arcs $[x_{\ell}', y_{\ell +1}']_{\cW}$. 

We first show that $r(x_{\ell}')   < 2$ for an infinite number of indices $\ell$. 
If not, then    there exists  $\ell_0 \geq 0$ which is  the least index such that $r(x_{\ell}')   \geq  2$ for all $\ell \geq \ell_0$.
  As $r(x_{\ell}') \ne 2$ for all $\ell \geq 0$,   thus  $r(x_{\ell}') > 2$ for all $\ell \geq \ell_0$.
 In particular, the choice of $\ell_0$ implies that $r(x_{\ell_0 -1}')< 2$. 
Then  $x_{\ell_0} \in E_i$  for $i=1$ or $i=2$,  and let   $\ox_{\ell_0}' \in L_i^+$ be the facing point to $x_{\ell_0}' \in L_i^-$. As $r(x_{\ell_0}')>2$,  Proposition~\ref{prop-propexist} implies that 
$x_{\ell_0}\prec_{\cK} \ox_{\ell_0}$ where $\ox_{\ell_0} =\tau(\ox_{\ell_0}')$. Thus there exists $\ell_1 > \ell_0$ such that $x_{\ell_1} = \ox_{\ell_0}$. Then 
$x_{\ell_0-1}' \prec_{\cW} x_{\ell_1}'$  by Lemma~\ref{lem-shortcut1} and so $r(x_{\ell_1}') = r(y_{\ell_0}')  = r(x_{\ell_0-1}' ) < 2$. This contradicts the choice of $\ell_0$, and thus  the set of indices $\ell$ for which $r(x_{\ell}') < 2$ must be infinite.

Next, we show that there exists a subsequence $\{ x_{\ell_i}' \mid i=1,2, \ldots \}$   with $r(x_{\ell_i}') < 2$ such that $\ds \lim_{i\to \infty} ~ r(x_{\ell_i}') = 2$.
Set  $\ds r_* = \limsup \left\{ r(x_{\ell}' ) \mid  \ell \geq 0 ~{\rm such ~that} ~ r(x_{\ell}') < 2\right\}$ and let 
 $\e_*  = 2-r_* \geq 0$ and  $r_*' = 2 - \e_*/2  \leq 2$. 
Assume  that   $\e_* > 0$,  then we show this yields a contradiction.

Note that  for all $0 < \e  < \e_*$ the collection $\{x_{\ell}'   \mid  2 - \e  \leq r(x_{\ell_i}') < 2  \}$ is finite by  the definition of $r_*$. 
 We can thus choose    $\ell_* \geq 0$ so  that $r(x_{\ell_*}') < 2$ and satisfying the condition that for all $\ell \geq \ell_*$,  if $r(x_{\ell}') < 2$    then  $r(x_{\ell}' ) < r_*'$. 
 Set $t_* = t_{\ell_*}$.

Observe that  Radius Inequality (K8) implies that there exists $\delta_*^+ > 0$ so that if   $r(x_{\ell}') \leq r_*'$ and  $x_{\ell +1}$ is a secondary entry point, then $r(x_{\ell +1}') > r(x_{\ell}') + \delta_*^+$. Conversely, the Radius Inequality also implies there exists $\delta_*^- > 0$ such that   if  $r(x_{\ell}') \leq r_*'$ and  $x_{\ell +1}$ is a secondary exit point, then $r(x_{\ell +1}') < r(x_{\ell}') - \delta_*^-$.

\begin{lemma} \label{lem-negativebound}
 There exists a greatest lower bound $N_*$ such that  $n_x(t) \geq N_*$   for all $t \geq 0$.
\end{lemma}
\proof
The set of values $\{ n_x(t) \mid 0 \leq t \leq t_* \}$ is finite, so it suffices to show   there is a lower bound on $n_x(t)$ when     $t \geq t_*$.

If no such bound exists, then  there exists an increasing subsequence of indices $\{ \ell_i > \ell_* \mid i \geq 1\}$ with $n_x(t_{\ell_1}) = n_x(t_{*}) - 1$ and
$n_x(t_{\ell_i}) = n_x(t_{*}) - i$ for all $i> 1$. Moreover, we can assume that $\ell_{i}$ is the least index $\ell > \ell_{i}$ such that $n_x(t_{\ell}) = n_x(t_{\ell_i}) - 1$. Consequently,  each point $x_{\ell_i}$ must be a secondary exit point.

If $\ell_{i+1} > \ell_i + 1$, then $n_x(t_k) \geq n_x(t_{\ell_i})$ for $\ell_i \leq k < \ell_{i+1}$, so by Proposition~\ref{prop-shortcut4} we have $x_{\ell_i}' \prec_{\cW} y_{\ell_{i+1}}'$ and $r(x_{\ell_i}') = r(y_{\ell_{i+1}}')$. 
For $i=1$,  this yields 
$$r(x_{\ell_{1}}') < r(x_{\ell_{1}-1}')   - \delta_*^- = r(x_{\ell_{*}}')   - \delta_*^- < r_*'$$ 
so that $r(x_{\ell_{1}}') < 2$. 
Then by induction, we have  that $r(x_{\ell_i}') < r(x_{\ell_*}') - i \, \delta_*^- < r_*'- i \, \delta_*^-$ and in particular  $r(x_{\ell_i}') < 2$. 

As $r(x_{\ell_i}') \geq 1$ and $r_*' < 2$, the value of the index $i$ is bounded by $i \leq  1/\delta_*^-$, contradicting our assumption that the subsequence is infinite. 
\endproof

Let $N_* = \min \{ n_x(t)  \mid  t \geq 0\}$,  and let $0 < t_{\ell_1} < t_{\ell_2} < \cdots < t_{\ell_k}$ be the sequence of times such that 
$n_x(t_{\ell_i}) = N_*$. Note that  the sequence $t_{\ell_i} < t_{\ell_i +1} < \cdots < t_{\ell_{i+1} +1}$ satisfies the hypotheses of Proposition~\ref{prop-shortcut4}  for each $i \geq 1$, and thus $x_{\ell_i}' \prec_{\cW} x_{\ell_{i+1}}'$ so $r(x_{\ell_i}') =  r(x_{\ell_{i+1}}')$. Then as in the proof of Proposition~\ref{prop-r>2}, this implies that $x_{\ell_1}' \prec_{\cW} x_{\ell_i}'$ for all $i > 1$. However, $r(x_{\ell_1}') < r_*' < 2$ so by Lemma~\ref{lem-wilsonlengths} the length of the $\cW$-orbit through $x_{\ell_1}'$ has an upper bound. This implies there is an upper bound on the number of $\cW$-segments $[x_{\ell_i}', y_{\ell_{i+1}+1}']_{\cW}$ which implies that $k$ is bounded above by a constant $k_1$ depending only on $r(x_{\ell_1}')$.

Then note that $n_x(t) \geq N_* +1$ for all $t \geq t_{\ell_{k_1} +1}$. Repeat the above  arguments  inductively,  to conclude that for each $i \geq 0$, there exists only a finite number of transition points $\{ t_{\ell} \mid \ell > 0\}$ with $r(x_{\ell}') < 2$ and $n_x(t_{\ell}) = N_*+i$.  It follows that for some index $\ell > \ell_*$, we have that $r(x_{\ell}') \geq r_*'$. However, $r_*' \leq r(x_{\ell}') \leq 2$  is forbidden, so we obtain that $r(x_{\ell}') > 2$ for all $\ell \geq \ell_k$ for $k$ sufficiently large. 
By Corollary~\ref{cor-r>2exit}  it follows that the orbit of $x_{\ell_k}$ escapes from $\mK$, contrary to assumption.

Therefore,    $r_* = 2$ and so there exists an infinite subsequence of transition points 
$\{ x_{\ell_i} \mid i=1,2, \ldots \}$  with  $r(x_{\ell_i}') <  2$ such that $\ds \lim_{i\to \infty} ~ r(x_{\ell_i}') = 2$.

Finally, note that the above proof shows that 
$\ds \liminf_{k \to \infty}  ~ \{ r(x_{\ell}') \mid \ell \geq k  \} = 2$.  
\endproof

 \medskip
We point out a corollary of the proof of Proposition~\ref{prop-cluster-}.

\begin{cor}\label{cor-always<2}
Suppose that  $x \in \mK$ satisfies $\rho_x(t) < 2$ for all $t$, then $x$ is not trapped. \hfill $\Box$
\end{cor}

 Finally, consider the   case  where
 $x \in \mK$ with  $\rho_x(t) > 2$ for all $t\geq 0$. Sections~\ref{sec-semilocal} and \ref{sec-global} developed various criteria for when the orbit of $x$ is not trapped. For example, if the level function $n_x(t) \geq 0$ for all $t \geq 0$, then by Proposition~\ref{prop-r>2} the $\cK$-orbit escapes. The techniques used above yield the following result for such   trapped orbits.  Note that later in this work, in  Proposition~\ref{prop-r>2R0} with the assumption of   Hypotheses~\ref{hyp-SRI} and \ref{hyp-genericW}, we obtain a stronger conclusion about such orbits.

 \begin{prop}\label{prop-cluster+}
Let $x \in \mK$ be trapped in forward time. 
If  $\rho_x(t) > 2$ for all $t \geq 0$, then there exists a special point $p_j^-$  as in \eqref{eq-special} and a subsequence 
$\{ x_{\ell_i} \mid i=1,2, \ldots \}$  of forward transition points   such that $\ds \lim_{i\to \infty} ~ x_{\ell_i}' = p_j^-$. The analogous result holds if $x$ is trapped in backward time. 
\end{prop}
\proof
 By assumption,  $r(x_\ell')>2$ for all transition points $x_{\ell}=\Phi_{t_{\ell}}(x)$ with $t_\ell\geq 0$.
Suppose that  $n_x(t)$ admits a minimum value, then consider  the
least $\ell$ such that $n_x(t_\ell)$ is equal to this minimum. Then,
$r(x_\ell')>2$ and $n_{x_\ell}(t)\geq 0$ for all $t\geq 0$ and   by
Proposition~\ref{prop-r>2} the $\cK$-orbit of $x$ escapes.

Hence,  there is an infinite sequence
$t_{\ell_0}<t_{\ell_1}<\cdots$ where $\ell_i$ is the least index
with $n_x(t_{\ell_i})=-i$, and all   $x_{\ell_i}$ must be secondary
exit points. We claim that $r(x_{\ell_i}')>r(x_{\ell_{i+1}}')>2$ for
all $i$. Observe that $r(x_\ell')\geq r(x_{\ell_0}')$ for all $\ell_0\leq
\ell<\ell_1$ and $r(x_{\ell_1-1}')=r(x_{\ell_0}')$ by
Proposition~\ref{prop-shortcut4}, so that  $2<r(x_{\ell_1}')<r(x_{\ell_1-1}')=r(x_{\ell_0})$
by the Radius Inequality.

 Continue in this way, to obtain a sequence of points $x_{\ell_i}$ for which 
 $$ 2 < \cdots< r(x_{\ell_{i+1}}') < r(x_{\ell_i}') < \cdots < r(x_{\ell_1}') < r(x_{\ell_0}') .$$
 It follows that $\ds \lim_{i \to \infty} ~   \{ r(x_{\ell_{i+1}}') -
 r(x_{\ell_i}')\} \to 0$ and $\ds \lim_{i \to \infty} ~  r(x_{\ell_i}')  = r_*
 \geq 2$. The Radius Inequality implies that for each $j = 1,2$, there
 is a    unique fixed point for the radius coordinate on the range and
 domain of the inverse $\sigma_j^{-1}$ of the insertion map, which is
 the point $r=2$ and $z=(-1)^j$. Thus, $\ds \lim_{i\to \infty} ~
 r(x_{\ell_i}) = 2$ and there is a subsequence of this subsequence
 converging to one of the special points. 
\endproof

  \begin{cor}\label{cor-anyx}
Suppose that the $\cK$-orbit of $x \in \mK$ is forward trapped,  then  there is a sequence of times $0 < t_1 < t_2 < \cdots$ such that
$\rho_x(t_{\ell}) \to 2$ and $z(\Phi_{t_{\ell}}(x)) \to \pm 1$. The analogous conclusion holds when the   $\cK$-orbit of $x \in \mK$ is backward trapped.
\end{cor}
\proof
Assume that the $\cK$-orbit of $x$ is forward trapped. 

If $\rho_x(t) > 2$ for $t \geq 0$, then the claims  follow from Proposition~\ref{prop-cluster+}. 

If $\rho_x(t) =2$ for some $t \geq 0$, then the claims  follow from   Proposition~\ref{prop-r=2any}.

If $r(x) < 2$ and $\rho_x(t)\neq 2$ for $t \geq 0$ then the claims  follow from  Proposition~\ref{prop-cluster-} 

If   the   $\cK$-orbit of $x \in \mK$ is backward trapped,   consider the reverse $\cK$-orbit to obtain the claims.
\endproof

  \bigskip
 
\section{The Kuperberg minimal set}\label{sec-minimalset}

In this section, we give the proof that the   flow $\Phi_t$ has no
periodic orbits in the plug $\mK$, and then discuss some of the
properties of the minimal set $\Sigma$ for the flow. We also state
some results on the non-wandering set of the flow.

\begin{thm}\label{thm-aperiodic}
The flow $\Phi_t$ in the Kuperberg Plug  is aperiodic.
\end{thm}
\proof
 
 Suppose there exist a periodic orbit   $x \in \mK$ for the flow $\Phi_t$. We show that this leads to a contradiction. 
 
Let   $x_{\ell} = \Phi_{t_{\ell}}(x)$ with  $0 \leq t_0 < t_1 < \cdots < t_n < \cdots $    be  the forward transition points in   $\{\Phi_t(x) \mid t \geq 0\}$ and suppose that  $x_n$  is the first   transition point   with $x_n = \Phi_{t_n}(x_0) = x_0$ for $t_n > 0$. Note that we then have  $\Phi_{t + t_n}(x) = \Phi_{t}(x)$ for all $t$.  

Let $[x_{\ell}', y_{\ell +1}']_{\cW} \subset \wmW$ be the lift of the $\cW$-arc $[x_{\ell}, x_{\ell + 1}]_{\cK} \subset \mK$ for $0 \leq \ell \leq n$. Then $[x_{0}', y_{1}']_{\cW} = [x_{n}', y_{n+1}']_{\cW}$. 
Recall that  the   $r$-coordinate is constant on each $\cW$-arc 
$[x_{\ell}', y_{\ell +1}']_{\cW}$ and define $r_0 = \min \{ r([x_{\ell}', y_{\ell +1}']_{\cW}) \mid 0 \leq \ell < n\}$.  We may   assume   this minimum occurs for  $[x_{0}', y_{1}']_{\cW}$.

If $r([x_{\ell}', y_{\ell +1}']_{\cW})  = 2$ for some $0 \leq \ell < n$, 
then the $\cK$-orbit of $x_{\ell}$ is either trapped or infinite   by
Proposition~\ref{prop-r=2any}. The same is true for $x$, which
contradicts the assumption that $x$ is a periodic point. Thus $r(x_\ell')\neq
2$ for all $\ell$. 
If $r([x_{\ell}', y_{\ell +1}']_{\cW})  < 2$ for some $0 \leq \ell < n$, then by Proposition~\ref{prop-cluster-} the orbit of $x$ contains a special orbit in its closure, which contradicts the assumption that the orbit is periodic. Thus $r_0 > 2$. 

As $r([x_{1}', y_{2}']_{\cW}) \geq r_0 > 2$, we have that $r(x_1') >
r(y_1')$ and then Condition  (K8) implies that $x_1$ must be a secondary entry point. Thus, $n_{x_0}(t_1) = 1$.

Next, observe that if $n_{x_0}(t) < 0$ for some $t > 0$, then there is a least $\ell > 0$  such that $n_{x_0}(t_{\ell}) < 0$. Hence   $n_{x_0}(t_{\ell -1}) =0$, and  both $x_{\ell -1}$ and $x_{\ell}$ are    secondary exit points by the minimality of $\ell$. 
Then by Proposition~\ref{prop-shortcut4}, we have $x_0' \prec_{\cW} y_{\ell}'$. Thus   $r(x_{\ell -1}') = r(x_0')=r_0$. 
As  $r(x_{\ell}') \ne 2$ and $x_\ell$ is a secondary exit point, it follows that $r(x_{\ell}') < r(x_{\ell -1}') = r_0$, which is a contradiction. Thus, we have $n_{x_0}(t) \geq 0$ for all $t\geq 0$.

  Next, suppose that $n_{x_0}(t_n) = 0$, we can then apply
Proposition~\ref{prop-shortcut4} to the orbit segment $[x_0 ,
x_{n+1}]_{\cK}$ to conclude that $x_0' \prec_{\cW} y_{n+1}'$. Hence,
$x_0'\prec_{\cW} y_1'\prec_{\cW} x_n'\prec_{\cW}y_{n+1}'$, but since
$x_0=x_n$ we conclude $x_0'\prec_{\cW}y_1'\prec_{\cW}x_0'$. Hence, the
$\cW$-arc $[x_0',y_1']$ is contained in one of the $\cW$-periodic
orbits, in particular $r(x_0)=2$ which is a contradiction.

The only possibility left is that $n_x(t_n) > 0$, but this is impossible by considering the reversed flow and then it becomes the forbidden case $n_x(t_n) < 0$ considered above. 
\endproof

The above proof is essentially that given by Kuperberg in \cite{Kuperberg1994}.   Proposition~\ref{prop-shortcut4} is not formally stated in her paper, though  it is stated explicitly   in the subsequent treatments   \cite{Ghys1995,Kuperbergs1996,Matsumoto1995}.

We   now give a summary of  general properties of the minimal set for the Kuperberg flow, as observed in \cite{Ghys1995,Kuperbergs1996,Matsumoto1995}. These are based on Theorem~\ref{thm-aperiodic} and   results of previous sections.
   Recall that  $p_i^{-} = \tau(\cL_i^- \cap \cO_i)$ for $i=1,2$ are the special entry points.  Define the orbit closures in $\mK$, 
\begin{equation}  \label{eq-minset1}
\Sigma_1 ~   \equiv   ~ \overline{ \{\Phi_t(p_1^-) \mid -\infty < t < \infty \}}     \quad , \quad 
\Sigma_2 ~   \equiv   ~ \overline{ \{\Phi_t(p_2^-) \mid -\infty < t < \infty \}}     ~ .
\end{equation}

\begin{thm} \label{thm-minimal} 
For the closed   sets   $\Sigma_i$ for $i=1,2$ we have:
\begin{enumerate}
\item $\Sigma_i$ is $\Phi_t$-invariant;
\item $r(x) \geq 2$ for all $x \in \Sigma_i$;
\item $\Sigma_1 = \Sigma_2$ and we set  $\Sigma = \Sigma_1 = \Sigma_2$;
\item $\Sigma \subset \cZ$, where $\cZ \subset \mK$ is any closed invariant set for $\Phi_t$.
\item $\Sigma$ is the unique minimal set for $\Phi_t$; 
\end{enumerate}
  \end{thm}
\proof
For 1)   note that the closure of any $\Phi_t$-orbit is a $\Phi_t$-invariant set.

For 2)  note that for any point $x \in \tau(\cO_1 \cap \wmW)$,
Corollary~\ref{cor-r=2special} shows that $n_x(t) \geq 0$ for all $t
\geq 0$, hence $\rho_x(t) \geq r(x) =2$. Thus all points $y \in \Sigma_i$ in the closure of the orbit also satisfy $r(y) \geq 2$.

For 3) first observe that the proof of Proposition~\ref{prop-r=2any}
for $p_1^-$ implies that $\Sigma_1$ contains $p_2^-$, and thus $\Sigma_2 \subset \Sigma_1$. It likewise implies that $\Sigma_1 \subset \Sigma_2$, and thus  $\Sigma_1 = \Sigma_2$. 

For 4) let $x \in \cZ$ then $\cZ$ is invariant implies that $x$ has an
infinite orbit. By Theorem~\ref{thm-aperiodic} the orbit of $x$ is not periodic.

If $\rho_x(t) =2$ for some $t\in \mR$, then  Proposition~\ref{prop-r=2any} shows that the $\Phi_t$-orbit of   $x$   contains a special point in its closure, hence $\Sigma \subset \cZ$. 
If $r(x) < 2$ and $\rho_x(t)\neq 2$ for all $t\in \mR$, then  Proposition~\ref{prop-cluster-} shows that the $\Phi_t$-orbit of   $x$   contains a special point in its closure, hence $\Sigma \subset \cZ$. If $r(x) > 2$ and $\rho_x(t) \leq 2$  for some $t$, then we are reduced to the above cases. Otherwise, if $\rho_x(t) > 2$ for all $t$, then  Proposition~\ref{prop-cluster+} shows that the $\Phi_t$-orbit of  $x$  again contains a special point in its closure. Thus, in all cases, the closure of an infinite orbit must contain $\Sigma$ in its closure. 

For  5) note that the closure of the $\Phi_t$-orbit of any $x \in \Sigma$ contains $\Sigma$, so the set is minimal. Suppose that $\cZ$ is a minimal set for $\Phi_t$ then given any $x \in \cZ$, by 4) we have $\Sigma \subset \cZ$, so they must be equal.  
 \endproof

We conclude this discussion of the Kuperberg minimal set, with some observations concerning other aspects of the topological dynamics of a Kuperberg flow, which follow from the results  of the previous sections. Recall that the orbits of the Kuperberg flow are divided into those which are finite, forward or backward trapped, or trapped in both directions and so infinite. Correspondingly, the asymptotic properties of the orbits must be considered within this restraint, as the asymptotic behavior outside of the plug $\mK$ of an orbit which escapes is not known.

First, we recall some standard definitions from topological dynamics. The \emph{forward limit set} of a forward trapped point $x \in \mK$ is the $\Phi_t$-invariant set
\begin{equation}
\alpha(x)  ~   \equiv   ~ \bigcap_{T \to \infty} ~ \overline{ \{\Phi_t(x) \mid t \geq T \}} ,
\end{equation}
and $x$ is \emph{forward recurrent} if $x \in \alpha(x)$. 
The \emph{backward limit set} of a backward trapped point $x \in \mK$ is the $\Phi_t$-invariant set
\begin{equation}
\omega(x)  ~   \equiv   ~ \bigcap_{T \to -\infty} ~ \overline{ \{\Phi_t(x) \mid t \leq T \}} 
\end{equation}
and $x$ is \emph{backward recurrent} if $x \in \omega(x)$. 

\begin{cor}
Let $x$ be forward trapped, then $\Sigma \subset \alpha(x)$. Likewise, if   $x$ is backward trapped, then   $\Sigma \subset \omega(x)$. 
\end{cor}
 \proof
Theorem~\ref{thm-minimal}.4 yields  the inclusions  $\Sigma \subset \alpha(x)$ and  $\Sigma \subset \omega(x)$. 
 \endproof

Next, consider the opposite extreme from recurrent points. A point $x \in \mK$  is \emph{forward wandering} if there exists an open set $x \in U \subset \mK$ and $T_U > 0$ so that for all $t \geq T_U$ we have $\Phi_t(U) \cap U = \emptyset$. Similarly, $x$ is \emph{backward wandering} if there exists an open set $x \in U \subset \mK$ and $T_U < 0$ so that for all $t \leq T_U$ we have $\Phi_t(U) \cap U = \emptyset$. A point $x$ with infinite orbit   is \emph{wandering} if it is forward and backward wandering. 

Define the following subsets of $\mK$:
\begin{eqnarray*}
\fW^{0} & \equiv & \{ x \in \mK \mid x ~ {\rm orbit ~ is ~ finite}\} \label{eq-wanderingfinite}\\
\fW^+ & \equiv & \{ x \in \mK \mid x ~ {\rm orbit ~ is ~ forward ~   wandering}\} \label{eq-wandering+}\\
\fW^- & \equiv & \{ x \in \mK \mid x ~ {\rm orbit ~ is ~ backward ~   wandering}\} \label{eq-wandering-}\\
\fW^{\infty} & \equiv & \{ x \in \mK \mid x ~ {\rm is ~ wandering}\} \label{eq-wanderinginfinite}
\end{eqnarray*}
Note that $x \in \fW^{0}$ if and only if the orbit of $x$ escapes through $\partial_h^+ \mK$ in forward time, and escapes though $\partial_h^- \mK$ in backward time. 
Define 
\begin{equation}\label{eq-dynamicdecomp}
\fW    ~ = ~  \fW^{0}    \cup \fW^+ \cup \fW^- \cup \fW^{\infty} \quad ; \quad  \Omega = \mK -  \fW .
\end{equation}

 The set $\Omega$ is called the \emph{non-wandering} set for $\Phi_t$. A point $x$ with forward trapped orbit    is characterized by the property: 
  $x \in \Omega$  if   for all $\e > 0$ and $T > 0$, there exists $y$
  and $t > T$ such that $d_{\mK}(x,y) < \e$ and $d_{\mK}(x, \Phi_t(y))
  < \e$, where $d_{\mK}$ is a distance function defined in Section~\ref{sec-radius}.
There are obvious corresponding statements for points which are   backward trapped or infinite. 

We now give some of the properties of the wandering and non-wandering sets.
\begin{lemma} \label{lem-nonwandering}
$\Omega$  is a closed, $\Phi_t$-invariant subset, with $\Sigma \subset \Omega$.
\end{lemma}
\proof
The fact that each of the sets $\fW^{0}$, $\fW^{+}$, $\fW^{-}$ and  $\fW^{\infty}$ is open follows directly from the definition of wandering, as does their invariance under the flow $\Phi_t$. For $x \in \Sigma$, its  orbit is  recurrent so is not wandering, hence $x \in \Omega$.
\endproof

\begin{lemma} \label{lem-wanderingboundary}
 If  $x \in \mK$   is   a primary entry or exit point, then $x \in \fW$.  
\end{lemma}
\proof 
Let $x \in \partial_h^- \mK$ be a primary entry point. Then for some $\e > 0$, the choice of the vector field $\cW$ on $\mW$ implies that the coordinate function $z(\Phi_t(z))$ is strictly increasing for $0 \leq t \leq \e$. For $0 < \delta < \e$, let $B(x,\delta)$ be the closed ball of radius $\delta$ centered at $x$. Then for $\delta$ sufficiently small, the image $\Phi_t(B(x,\delta))$ is disjoint from $B(x,\delta)$ for all $t > \e$. Thus, $x \in \fW^{0}$ if its orbit is not trapped and $x \in \fW^{+}$ if its orbit is  trapped. Similar considerations apply for a primary exit point $x \in \partial_h^+ \mK$, to show that either $x \in \fW^{0}$ or $x \in \fW^{-}$.
\endproof

\begin{cor} \label{cor-wanderinginfinite}
For each $x \in \Omega$, the $\Phi_t$-orbit of $x$  is infinite.
\end{cor}
 
 We can also restrict the radius coordinates of the non-wandering orbits.
 \begin{prop} \label{prop-wanderinginterior}
$\Omega \subset \{x \in \mK \mid r(x) \geq 2\}$. 
\end{prop}
\proof
Let $x \in \Omega$ and suppose that $r(x) = r_0 < 2$.    Note that $x$ is an infinite orbit implies that $r(x) > 1$ and $x$ cannot be contained in the orbit of a special point as $r(x) < 2$. 
Moreover, it follows as in the proof of Lemma~\ref{lem-negativebound},  that the level function $n_x(t)$ has a greatest lower bound $N_* \leq 0$ for all $t \in \mR$. 
Let $x_* = \Phi_{t_*}(x)$ be a transition point with $n_x(t_*) = N_*$. 
As $\Omega$ is $\Phi_t$-invariant, it will suffice to show that $x_*$ is wandering, which is a contradiction.  Thus, we may assume that $x = x_*$ and   the level function satisfies $n_x(t) \geq 0$ for $t \geq 0$.

Let   $x_{\ell} = \Phi_{t_{\ell}}(x)$ with  $0 = t_0 < t_1 < \cdots < t_{\ell} < \cdots $    be  the forward transition points in   $\{\Phi_t(x) \mid t \geq 0\}$. 
Let $[x_{\ell}', y_{\ell +1}']_{\cW} \subset \wmW$ be the lift of the $\cW$-arc $[x_{\ell}, x_{\ell + 1}]_{\cK}$, for   $\ell \geq 0$.

The proof of Proposition~\ref{prop-cluster-} shows that  there are
only a finite number of   indices $\ell \geq 0$  such that $x_{\ell}$
is a transition point with   level $n_x(t_{\ell}) = 0$. Let $\ell_0
\geq 1$ be the greatest index for which $n_x(t_{\ell_0 -1}) = 0$. Then
$n_x(t_{\ell}) \geq  1$ for all $\ell \geq \ell_0$, which implies that both $x_{\ell_0}$ and $x_{\ell_0 +1}$  are  secondary entry points and thus 
   $x_{\ell_0}' \in  L_i^-$  for $i = 1$ or $i = 2$ and $y_{\ell_0 +1}' \in  \cL_j^-$  for $j = 1$ or $j = 2$.

 We have $r_0 = r(x_{0}') = r(x_{\ell_0 -1}')$ and set $r_1 = r(x_{\ell_0}')$. As $x$ is not on a special orbit,  the Radius Inequality (K8) implies that $r_1 > r_0$.  By the geometry of the insertion maps $\sigma_i$ it follows that $x_{\ell_0}' \in  L_i^-$ must be an interior point for the compact region $L_i^-$.  The flow $\Phi_t$   is transverse to the section $\cT_{\cK}$ defined in \eqref{eq-transversal},  so it follows that    there exists $\delta_0 > 0$ such that, for each   $0 \leq \ell \leq \ell_0$, the ball $B(x_{\ell}, \delta_0)$ intersects just one component of $\cT_{\cK}$ and is disjoint from its boundary.

 By the continuity of the flow, we can choose $0 < \e \leq \delta_0$ sufficiently small so that 
  $\Phi_{t_{\ell}}(B(x, \e)) \subset  B(x_{\ell}, \delta_0)$ for each $0 \leq \ell \leq \ell_0$.
  Moreover, for $r_0' = r_0 + (r_1 - r_0)/3$ and $r_1' = r_1 - (r_1 - r_0)/3$ we require that $r(y) \leq r_0'$ for each $y \in B(x, \e) \cap \cT_{\cK}$ and 
  $r(z) \geq r_1'$ for each $z \in B(x_{\ell_0}, \delta_0) \cap \cT_{\cK}$.

 For  $y \in B(x, \e)$, it then   follows that for each $0 \leq \ell \leq \ell_0$, $\Phi_{t_{\ell}'}(y)$ is an interior point of $\cT_{\cK}$ for some $t_{\ell}'$   close to $t_{\ell}$. In particular,   $\Phi_{t_{\ell_0}'}(y)$ is a secondary entry point with   $n_y(t_{\ell_0'}') = n_x(t_{\ell_0}) = 1$. 
As in the proof of Proposition~\ref{prop-cluster-},  it follows   that  $n_y(t) \geq  1$ for all $t \geq t_{\ell_0}$ and $r(\Phi_t(y)) \geq r(\Phi_{t_{\ell_0}'}(y)) \geq r_1'$ for all $t \geq t_{\ell_0}'$. Thus, $\Phi_t(B(x, \e)) \cap B(x, \e) = \emptyset$ for all $t \geq t_{\ell_0}$, so that $x$ is forward wandering, as was to be shown.
\endproof

 \bigskip 

\section{The Kuperberg pseudogroup}\label{sec-pseudogroup}

In this   section, we define  a pseudogroup $\cGK$ acting on a rectangle   $\bRt \subset \mK$ which captures the dynamics of the flow $\Phi_t$.
The study of the action of  $\cGK$ leads to 
a deeper understanding of the geometry and topology of the minimal set $\Sigma$, and gives a framework for the rest of the paper.    
The analysis of the dynamical properties of the action
of $\cGK$ on $\bRt$ uses many of the results and techniques developed
in the previous sections and provides an interpretation of  these   results in a pseudogroup setting.

\begin{figure}[!htbp]
\centering
{\includegraphics[width=120mm]{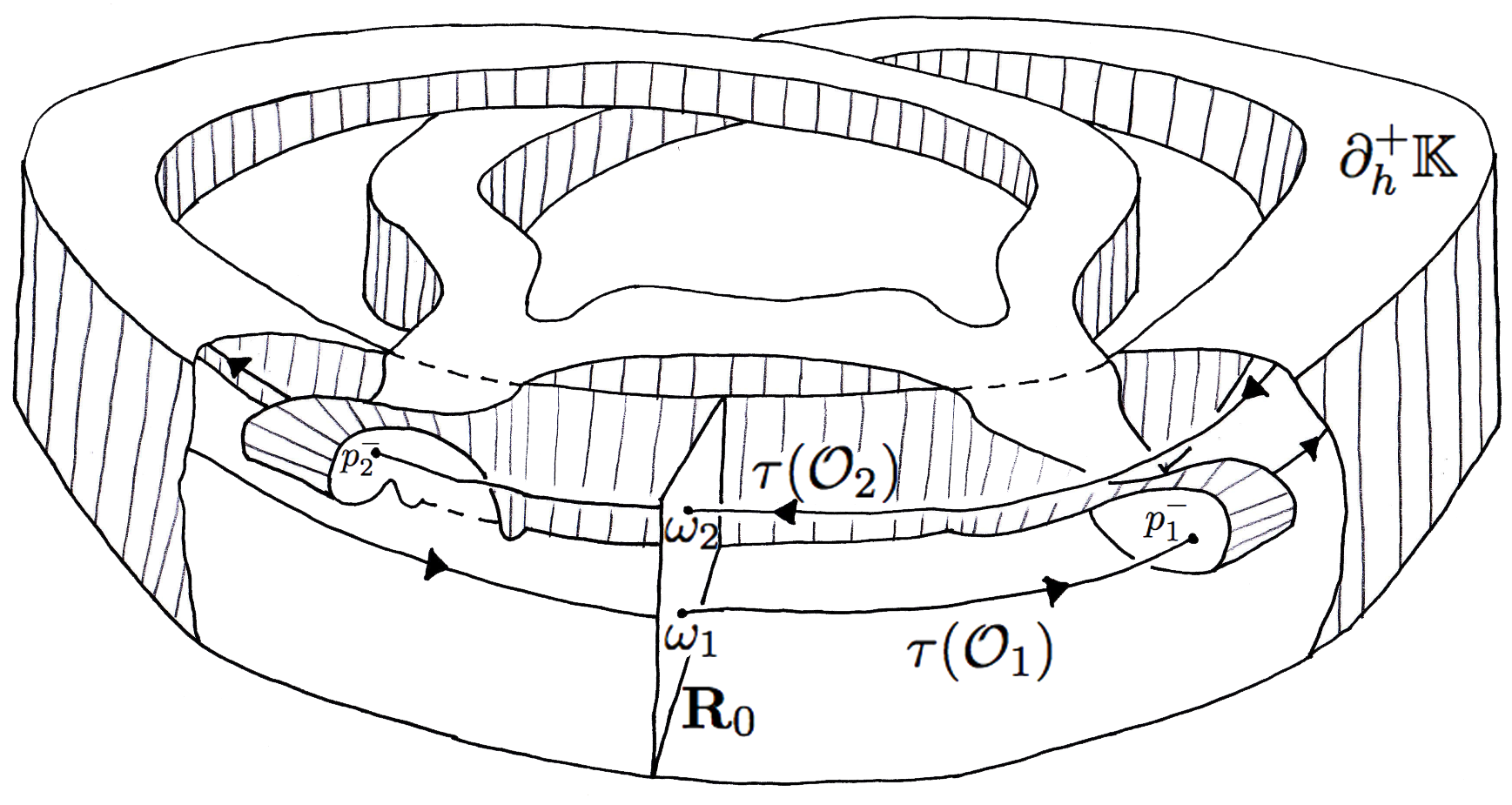}}
 \caption{\label{fig:KR} The rectangle $\bRt$ in the Kuperberg plug $\mK$}
 \vspace{-6pt}
\end{figure}

Choose a value of $\theta_0$ such that the rectangle $\bRt$ as defined in cylindrical coordinates, 
\begin{equation}\label{eq-goodsection}
\bRt \equiv \{\xi = (r, \theta_0, z) \mid ~ 1 \leq r \leq 3 ~,   ~ -2
\leq z \leq 2\}  \,\subset \mW' ~ , 
\end{equation}
is disjoint from both the regions $D_i$ and their insertions $\cD_i$
for $i=1,2$, as defined in Section~\ref{sec-kuperberg}. 
 For example, for the curves  $\alpha_i$ and $\beta_i'$  defined in
Section~\ref{sec-kuperberg}, we can take $\theta_0 = \pi$, so that
$\bRt$ is between the embedded regions  $\cD_i$ for $i=1,2$ as illustrated in Figure~\ref{fig:KR}.

As $\bRt \subset \mW'$,  the quotient map $\tau \colon \mW \to \mK$ is injective on $\bRt$  and   we denote  its image  in $\mK$  also   by $\bRt$ with coordinates   $r = r(\xi)$ and $z = z(\xi)$  for  $\xi \in \bRt$. 
The    metric on   $\bRt$ is given by    
\begin{equation}\label{eq-metricR}
 d_{\bRt}(\xi, \xi') = \sqrt{(r' - r)^2 + (z'-z)^2} \quad  {\rm  for} ~  \xi = (r, \pi, z) ~ , ~ \xi' = (r', \pi, z') ~ .
\end{equation}
 Introduce the special points $\omega_i \in \bRt$, given by 
  \begin{equation}\label{eq-omegas}
 \omega_1  = \cO_1 \cap \bRt = (2,\pi, -1)  \quad , \quad \omega_2 = \cO_2 \cap \bRt = (2,\pi, 1)   ~ .
\end{equation}
The  first transition point of the forward orbit of    $\omega_i$ is the special entry point 
 $p_i^{-} = \tau(\cL_i^- \cap \cO_i) \in E_i$. 
The first transition point of the backward orbit of    $\omega_i$  is the special exit point 
 $p_i^{+} = \tau(\cL_i^+ \cap \cO_i) \in S_i$.

We next give  a basic result concerning the behavior of $\cK$-orbits with respect to $\bRt$.

 \begin{prop}\label{prop-escapeS}
 Let $x \in \mK$ be such that its forward (or backward) $\cK$-orbit is trapped, then the orbit   intersects $\bRt$.    Hence, if the $\cK$-orbit of $x$   does not intersect $\bRt$,   it escapes from $\mK$ in both forward and backward time.  
 \end{prop}
 \proof   Both the  forward and the backward  $\cK$-orbit of a special point $p_i^-$  limits to each special point $p_j^-$  by Proposition~\ref{prop-r=2any}, and thus intersects $\bRt$ repeatedly. 
 It follows that there is  an open neighborhood of each of the special points $\{\omega_1,\omega_2\}$ consisting of  points whose   forward and   backward  $\cK$-orbit intersects $\bRt$.  
 Suppose that $r(x) < 2$, then by Proposition~\ref{prop-cluster-} the assumption that the forward orbit of $x$ is trapped   implies that it contains a special point in its closure, hence must intersect $\bRt$  infinitely often.  
The case where $r(x) \geq 2$   follows similarly by Propositions~\ref{prop-r=2any} and  \ref{prop-cluster+}. The case when  the backward orbit is trapped follows by reversing the flow and applying the above case.
\endproof
 
 \begin{cor}\label{cor-escapeS}
 If $x \in \mK$ belongs to an infinite $\cK$-orbit, then it     intersects $\bRt$ in an infinite sequence of points in both forward and backward time. \hfill $\Box$
  \end{cor}

The  first return map $\whPhi$ on $\bRt$ for the Kuperberg flow $\Phi_t$ is defined at $\xi \in \bRt$ if there is a $\cK$-orbit segment $[\xi, \eta]_{\cK}$ with $\eta \in \bRt$ and its interior $(\xi, \eta)_{\cK}$ is disjoint from $\bRt$. We then set $\whPhi(\xi) = \eta$. 
The domain of $\whPhi$ is the set:
\begin{equation}\label{eq-domainwhPhi}
Dom(\whPhi) \equiv \left\{ \xi \in \bRt \mid \exists ~ t > 0 ~ \text{ such ~ that} ~ \Phi_t(\xi) \in \bRt ~ \text{and} ~  \Phi_s(\xi)\notin \bRt 
  ~ \text{for} ~   0<s<t   \right\} .
\end{equation}
 Corollary~\ref{cor-escapeS} implies that every $x \in \bRt$ with infinite orbit is in the domain of $\whPhi$, and thus the dynamical behavior of the map  $\whPhi$ reflects the dynamical behavior of the infinite orbits for $\Phi_t$.    We also note that the map $\whPhi \colon Dom(\whPhi) \to \bRt$ has many points of discontinuity, as will be shown.

Recall the formal definition of a pseudogroup modeled on a space $X$: 
\begin{defn}\label{def-pseudogroup}
A pseudogroup $\cG$ modeled on a topological space $X$ is a collection of  homeomorphisms between open subsets of $X$ satisfying the following properties:
\begin{enumerate}
\item For every open set $U \subset X$, the identity $Id_U \colon U \to U$ is in $\cG$.
\item For every $\varphi \in \cG$  with  $\varphi \colon U_{\varphi} \to V_{\varphi}$ where $U_{\varphi}, V_{\varphi} \subset X$ are open subsets of $X$, then    also  $\varphi^{-1} \colon V_{\varphi} \to U_{\varphi}$ is in $\cG$.
\item  For every $\varphi \in \cG$  with  $\varphi \colon U_{\varphi} \to V_{\varphi}$  and each open subset $U' \subset U_{\varphi}$, then the restriction $\varphi \mid U'$ is in $\cG$.
\item  For every $\varphi \in \cG$  with  $\varphi \colon U_{\varphi} \to V_{\varphi}$ and   every $\varphi' \in \cG$  with  $\varphi' \colon U_{\varphi'} \to V_{\varphi'}$, 
if $V_{\varphi} \subset U_{\varphi'}$ then the composition $\varphi' \circ \varphi$ is in $\cG$.
\item If $U \subset X$ is an open set,   $\{U_{\alpha} \subset X \mid \alpha \in \cA\}$ are open sets whose union is $U$, $\varphi \colon U \to V$ is a homeomorphism to an open set $V \subset X$ and for each $\alpha \in \cA$ we have $\varphi_{\alpha} = \varphi \mid U_{\alpha} \colon U_{\alpha} \to V_{\alpha}$ is in $\cG$, then $\varphi$ is in $\cG$.
\end{enumerate}
\end{defn}

 We apply this definition to the first return map $\whPhi$ on $\bRt$ to obtain:
   \begin{defn}\label{def-pseudogroupK}
  Let  $\cGK$   be the pseudogroup  generated by the  map    $\whPhi$ acting  on    $X = \bRt$. 
That is, if $U \subset Dom(\whPhi)$ is an open set with $V = \whPhi(U)$ and the restriction $\whPhi
 \mid U$ is a homeomorphism, then both $\whPhi | U$ and $\whPhi^{-1} |
 V$ are in $\cGK$. The remaining elements of $\cGK$ are given by
 adding maps as required so that  the conditions of   Definition~\ref{def-pseudogroup} are satisfied.   
\end{defn}

In the rest of this section, we consider the   basic dynamical properties of  the action of   $\cGK$ on $\bRt$ and begin developing  the relationship between the dynamics of the flow $\Phi_t$ on $\mK$ and that of the induced action of $\whPhi$. We first consider   the relation between the actions of the maps $\whPhi$ and $\whPsi$ on $\bRt$, where $\whPsi$ denotes the return map to $\bRt$ for   the Wilson flow $\Psi_t$. The dynamical properties   of $\Psi_t$ on $\mW$  are  described in Proposition~\ref{prop-wilsonproperties}, and illustrated in  Figures~\ref{fig:flujocilin} and \ref{fig:Reebcyl}.

The  first return map $\whPsi$ on $\bRt$ for the Wilson flow $\Psi_t$ is defined at $\xi \in \bRt$ if there is a $\cW$-orbit segment $[\xi, \eta]_{\cW}$ with $\eta \in \bRt$ and its interior $(\xi, \eta)_{\cW}$ is disjoint from $\bRt$. We then set $\whPsi(\xi) = \eta$. 
The domain of $\whPsi$ is the set:
\begin{equation}\label{eq-domainwhPsi}
 Dom(\whPsi) \equiv \left\{ \xi \in \bRt \mid \exists ~ t > 0 ~
    \text{ such ~ that} ~ \Psi_t(\xi) \in \bRt ~ \text{and}  ~
  \Psi_s(\xi)\notin \bRt ~ \text{for}~   0<s<t  \right\}.
\end{equation}
The radius function is constant  along the orbits of the Wilson flow,  so that $r(\whPsi(\xi)) = r(\xi)$ for all $\xi \in Dom(\whPsi)$. 
Also, note that the points $\omega_i$ for $i=1,2$ defined in \eqref{eq-omegas} are   fixed-points for $\whPsi$.
 For all  other points $\xi \in \bRt$ with $\xi \ne \omega_i$,  it was assumed in Section~\ref{sec-wilson} that the function $g(r,\theta, z) > 0$, so the   $\cW$-orbit of $\xi$   has a ``vertical drift'' arising from the term 
$g(r, \theta, z)  \frac{\partial}{\partial  z}$ in the formula \eqref{eq-wilsonvector} for $\cW$.

For    $\xi \in  \bRt$  sufficiently close to the vertical boundary   $\partial_v \mW$, function $f(r,\theta,z) =0$ and hence the vector field $\cW$ is tangent to $\bRt$. Then    the $\Psi_t$-flow of $\xi$ is contained  in $\bRt$ and hence $\xi \not\in Dom(\whPsi)$.

For $\xi \in  \bRt$  such that $\cW$ is not vertical at $\xi$, then the $\Psi_t$-flow of $\xi$ exits $\bRt$ and has increasing $z$-coordinate, and so 
 flows upward until it either returns to $\bRt$, or exits through $\partial_h^+ \mW$ without
intersecting $\bRt$, in which case it fails to make a complete
revolution around the cylinder. In the latter case,  $\xi \not\in Dom(\whPsi)$. If $z(\xi) \geq 0$, then whether its forward $\Psi_t$-flow returns to $\bRt$ or not, depends strongly on the choices of the functions $g(r,\theta,z)$ and $f(r,\theta,z)$ and the minimum distance from the forward orbit of $\xi$ to the fixed-point   $\omega_2$.

On the other hand, let $\xi = (r,\pi, z) \in \bRt$ with $z< 0$ and suppose that $\cW$ is not vertical at $\xi$. Then as before,    the forward $\Psi_t$-orbit of $\xi$ exits $\bRt$ and has increasing $z$-coordinate, so must intersect the annular region $\cA = \{z=0\}$. If the forward orbit does not intersect $\bRt$ before crossing $\cA$, 
then by the anti-symmetry hypothesis on the Wilson flow, the $\Psi_t$-orbit of $\xi$ returns to $\bRt$ at the point $\eta =  (r,\pi, -z)$, and  so  $\whPsi(\xi) = \eta$. 
 
Note that these remarks imply that for all $\xi \in \bRt$ with $z(\xi) < 0$ and such that $\cW$ is not vertical at $\xi$,  then $\xi \in Dom(\whPsi)$. However, which of the two cases above occurs again depends strongly on the choices of the functions $g(r,\theta,z)$ and $f(r,\theta,z)$ and the minimum distance from the forward orbit of $\xi$ to the fixed-points $\omega_i$.

 Finally, we consider the domain of $\whPsi$ in an open neighborhood of the vertical segment $\{r=2\}$. 
For $\xi \in \bRt$ sufficiently close to one of  the fixed points   $\omega_i$  the function  $g(r,\theta,z)$ is arbitrarily small, so the forward $\cW$-orbit  of $\xi$ must intersect    $\bRt$, and thus  $\xi \in Dom(\whPsi)$. In particular, for each $\omega_i$ the set $ Dom(\whPsi)$   contains a  open neighborhood of  $\omega_i$. The   return map for $\whPsi$ at other points close to $\{r=2\}$ is  considered in three cases, which are illustrated in Figures~\ref{fig:flujocilin} and \ref{fig:Reebcyl}.

First, the $\Psi_t$-flow of the points in   the open   segment   
$$ \cI_{0}  = \left\{(2,\pi ,z) \mid  -1 < z < 1 \right\} \subset \{r=2\} \cap \bRt$$
always return to the same segment, so we have $\cI_{0} \subset  Dom(\whPsi)$ with $\whPsi \colon \cI_0 \to \cI_0$. Moreover,  $\whPsi$  is bijective  when restricted to $\cI_0$. However, $\whPsi$ is not continuous on $\cI_0$ as we will discuss further below.

The conditions (W5) and (W6) in Section~\ref{sec-wilson} imply there exists a least $0 < \e_f < 1/4$ such that the function $f(2, \theta, z) > 0$ for $-2 + \e_f < z < 0$. Thus,   for the open segment 
$$ \cJ_{0}  = \left\{(2,\pi,z) \mid  -2 + \e_f \leq z < -1 \right\}   \subset \{r=2\} \cap \bRt$$
we have $\cJ_0 \subset Dom(\whPsi)$ and $\whPsi(\cJ_0) \subset \cJ_0$. By the anti-symmetry condition (W1) for the function $f$,  we also have that $f(2, \theta, z) < 0$ for $0 < z < 2 - \e_f$. Thus,  for the open segment 
$$  \cK_{0}  = \left\{(2,\pi,z) \mid  1 < z \leq 2 - \e_f \right\} \subset \{r=2\} \cap \bRt$$
we have $\cK_0 \subset Dom(\whPsi^{-1})$ and $\whPsi^{-1}(\cK_0) \subset \cK_0$.  
Thus,   $\cI_0$, $\cJ_0$ and  $\whPsi^{-1}(\cK_0)$ are   in the interior of  $Dom(\whPsi)$. 

\begin{figure}[!htbp]
\centering
\begin{subfigure}[c]{0.3\textwidth}{\includegraphics[height=85mm]{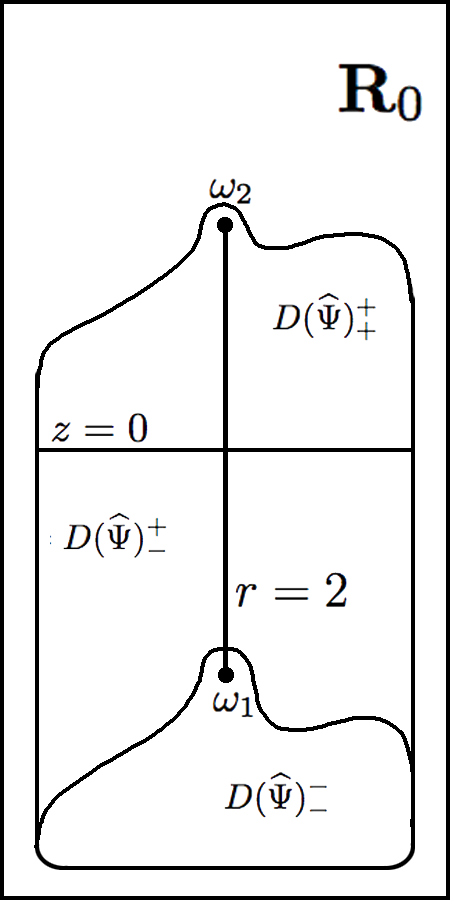}}\caption{Domains of $\whPsi$}\end{subfigure}
\begin{subfigure}[c]{0.3\textwidth}{\includegraphics[height=85mm]{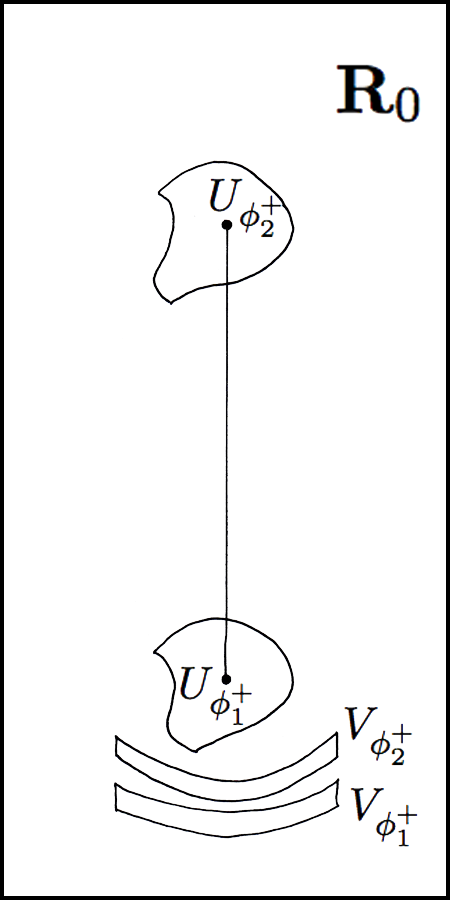}}\caption{Domains of $\phi_1^+, \phi_2^+$}\end{subfigure}
\begin{subfigure}[c]{0.3\textwidth}{\includegraphics[height=85mm]{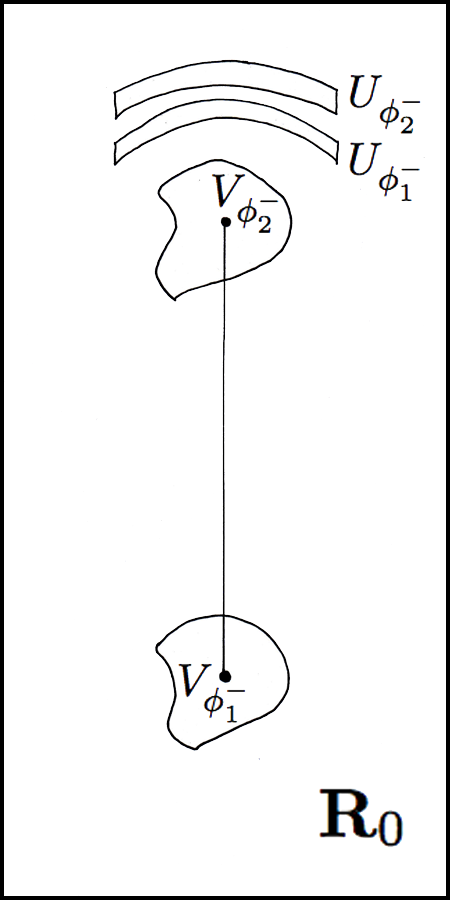}}\caption{Domains of $\phi_1^-, \phi_2^-$}\end{subfigure}
\caption{\label{fig:GKdomains}  Domains and ranges for the maps  $\{\whPsi, \phi_1^+, \phi_2^+, \phi_1^-, \phi_2^-\}$}
\vspace{-6pt}
\end{figure}

In order to illustrate the regions in  the domain of $\whPsi$ as described above, assume that  $0 \leq g(r,\theta,z) \leq 1/10$ and that $g(r,\theta,z) = 1/10$  when allowed. Thus,  the flow of $\Psi_t$ rises at an approximate  rate of $1/10$ in $\mW$, and $z(\whPsi(\xi)) \approx z(\xi) + r(\xi)/10$.
   Figure~\ref{fig:GKdomains}.(A)  pictures the three regions of the domain in this case.

 We next consider the continuity properties of the return map $\whPsi$.      
  Recall that   the Wilson flow  reverses direction at the annulus $\cA = \{z=0\} \subset \mW$, and is anti-symmetric  with respect to the annulus $\cA$ by Condition (W1)  in Section~\ref{sec-wilson}, so that $\cW$   is tangent to $\bRt$ along the line $\cT = \cA \cap \bRt$.

Let $\xi = (r, \pi, z) \in  Dom(\whPsi)$  and   $\eta=\whPsi(\xi) = (r,\pi, 0)$. Then there exists a $\cW$-orbit segment $[\xi , \eta]_{\cW}$   which intersects $\bRt$ only in its endpoints, and the $\Psi_t$ flow of $\xi$ is tangent to $\bRt$ at  $\eta \in \cT$. 
Let $\xi' = (r,\pi, z')$ with $z - \e' < z' < z$, for $\e' > 0$ sufficiently small, then $\eta' = \whPsi(\xi')$ is defined and satisfies $z(\whPsi(\xi')) < 0$, with value depending continuously on $z'$. 
On the other hand, for $\xi'' = (r, \pi, z'')$ with $z < z'' < z+\e''$, with $\e'' > 0$ sufficiently small, then the point $\eta'' = \whPsi(\xi'')$ is again well-defined. 
But the $\cW$-orbit segment $[\xi'' , \eta'']_{\cW}$  is no longer tangent to $\bRt$ near $\eta$, as it traverses $\mW$ in a counter-clockwise direction in the region $\{z< 0\}$, until   it crosses the plane $\{z=0\}$ before reaching the surface $\bRt$, and then afterwards reverses direction and subsequently intersects $\bRt$ from the opposite direction in the region $\{z > 0\}$, at a point close to  $\whPhi^2(\xi)$.   Consequently, the map $\whPsi$ has a discontinuity at $\xi$, and so  $\cL_1 = \whPsi^{-1}(\cT)$ is a curve of discontinuities for $\whPsi$. 

There is another type of discontinuity for $\whPsi$ that arises for $\xi = (r,\pi,0) \in \cT$. For $\e' > 0$ sufficiently small,  then for $\xi' = (r,\pi, -z')$ with $0 < z' < \e'$, we have   $\whPsi(\xi') = (r,\pi, z')$ due to the anti-symmetry of  $\Psi_t$. 
On the other hand,   for $\xi'' = (r, \pi, z'')$ with $0 < z'' < \e'$, suppose that $\eta'' = \whPsi(\xi'')$ is   well-defined then  the $\cW$-orbit segment $[\xi'' , \eta'']_{\cW}$  traverses $\mW$ in a  clockwise direction in the region $\{z> 0\}$, and so the value $z(\whPsi(\xi''))$ is   much larger than $0$. 
Thus, $\cT \cap Dom(\whPsi)$ is   a set of discontinuities for $\whPsi$.

Define the three domains of continuity for the   induced return map $\whPsi$ on $\bRt$, as illustrated in    Figure~\ref{fig:GKdomains}.(A).
\begin{eqnarray*}
D(\whPsi)_-^- ~ & = & ~ \{ \xi \in Dom(\whPsi) \mid z(\xi) < 0 ~ {\rm and}~ z(\whPsi(\xi)) \leq 0 \}\\
D(\whPsi)_-^+ ~ & = & ~ \{ \xi \in Dom(\whPsi) \mid z(\xi) < 0 ~ {\rm and}~ z(\whPsi(\xi)) > 0 \}\\
D(\whPsi)_+^+ ~ & = & ~ \{ \xi \in Dom(\whPsi) \mid z(\xi) \geq 0 ~ {\rm and}~ z(\whPsi(\xi)) > 0 \} ~ ,
\end{eqnarray*}
Restriction of the   map  $\whPsi$  thus yields three   continuous maps with disjoint domains,    denoted by
\begin{equation}\label{defPsi}
\whPsi_- = \whPsi | D(\whPsi)_-^- \quad , \quad \whPsi_0 = \whPsi | D(\whPsi)_-^+ \quad , \quad \whPsi_+ = \whPsi | D(\whPsi)_+^+ ~ .
\end{equation}
 
 A comparable analysis of the domains of continuity for the return map $\whPhi$ for the flow $\Phi_t$  is far more complicated, and will be postponed until Section~\ref{sec-entropyflow}, when techniques have been introduced which are sufficient for describing these domains.  
Our strategy here is to define a   collection  of special  maps in $\cGK$, obtained by the restriction of $\whPhi$ to particular domains of continuity, and  which correspond to   key aspects of the dynamics of the flow $\Phi_t$. We show that these special elements and all their compositions in $\cGK$   capture all of the essential dynamical properties of the flow $\Phi_t$. 
 We first consider    elements of $\cGK$ corresponding to the entry and exit dynamics of the two insertions.

For $i = 1,2$,  let $U_{\phi_i^+} \subset Dom(\whPhi)$ be the
  subset of $\bRt$ consisting of points 
 $\xi \in Dom(\whPhi)$ with $\eta = \whPhi(\xi)$,  such that the $\cK$-arc $[\xi, \eta]_{\cK}$
  contains a single transition point $x$, with  $x \in E_i$. 
Note that for such $\xi$, we see from Figures~\ref{fig:cWarcs} and ~\ref{fig:KR}, that its $\cK$-orbit exits the surface $E_i$ as the $\cW$-orbit of a point  $x' \in L_i^-$ with $\tau(x') = x$, flowing upwards from $\partial_h^- \mW$ until it   intersects $\bRt$ again. If the $\cK$-orbit of $\xi$ enters   $E_i$  but   exits   through $S_i$ before crossing $\bRt$, then it is not considered to be in the domain $U_{\phi_i^+}$ as it contains more than one transition point.

Let  $\phi_i^+ \colon U_{\phi_i^+} \to V_{\phi_i^+}$ denote the element of $\cGK$   defined by the restriction of $\whPhi$. 
As the $\cK$-arcs $[\xi, \eta]_{\cK}$ defining $\phi_i^+$ do not intersect $\cA$, the restricted map $\phi_i^+$ is continuous. 
The inverse map $(\phi_i^+)^{-1} \colon V_{\phi_i^+} \to U_{\phi_i^+}$
is also in $\cGK$ for $i =1,2$. The sets $U_{\phi_i^+}$ and
$V_{\phi_i^+}$ are sketched in the center illustration   in Figure~\ref{fig:GKdomains}.

  For $i = 1,2$,  let $U_{\phi_i^-} \subset Dom(\whPhi)$ be the  
  subset of $\bRt$ consisting of points 
 $\xi \in Dom(\whPhi)$ with $\eta = \whPhi(\xi)$,  such that the $\cK$-arc $[\xi, \eta]_{\cK}$
  contains a single transition point $x$, with  $x \in S_i$.  
Then let  $\phi_i^- \colon U_{\phi_i^-} \to V_{\phi_i^-}$ denote the element of $\cGK$   defined by the restriction of $\whPhi$.
Again, as the $\cK$-arcs $[\xi, \eta]_{\cK}$ defining the maps $\phi_i^-$ do not intersect $\cA$,   the restricted map $\phi_i^-$ is continuous.
The   inverse map $(\phi_i^-)^{-1} \colon V_{\phi_i^-} \to U_{\phi_i^-}$
 is also in $\cGK$ for $i =1,2$. The sets $U_{\phi_i^-}$ and
$V_{\phi_i^-}$ are sketched in the right hand side illustration  in Figure~\ref{fig:GKdomains}.

We comment on some details of the regions in Figures~\ref{fig:GKdomains}.(B) and (C).  For the map 
$\phi_i^+$, $i=1,2$,  the domain contains a neighborhood of the point
$\omega_i$. Flowing the domain $U_{\phi_i^+}$ forward to $E_i$ and
then applying the map $\sigma_i^{-1}$ we obtain a set
$\widetilde{U_{\phi_i^+}}\subset L_i^-$ containing points with
$r$-coordinate equal to 2. Observe that the Radius Inequality implies
that the maximum radius of points in $\widetilde{U_{\phi_i^+}}$ is
bigger than the maximum radius of points in $U_{\phi_i^+}$. The first
intersection of the $\cW$-orbits of points in
$\widetilde{U_{\phi_i^+}}$ with $\bRt$ is thus a region containing
points with $r$-coordinate equal to 2 and since these points climb
slower than other points, the region folds at $\{r=2\}$. By the choice
of $\bRt$, the points coming from $L_2^-$ are farther than the ones
coming from $L_1^-$ and thus $V_{\phi_2^+}$ is above $V_{\phi_1^+}$.
Similar considerations apply to the maps $\phi_i^-$ for $i=1,2$.
 
Observe that each generator $\phi_i^+$ for $i=1,2$ corresponds to a flow through a transition point which increases the level function $n_x(t)$ by $+1$, 
so the inverse of $\phi_i^+$  decreases the level function by $1$.  The map $\phi_i^-$ decreases the level by $-1$ and its inverse increases the level by $+1$. 
  
 Next, we develop a relation between the action of the return map $\whPsi$ for the Wilson flow and the action of the pseudogroup $\cGK$.
Define    subsets of $\bRt$  contained in the domain of  the Wilson map $\whPsi$:
\begin{eqnarray*}
U_{\psi}^*  ~ &   = & ~ \{ (r, \pi, z) \in Dom(\whPsi)  \mid r > 2\}  \subset \bRt \cap \{r > 2\} \\
\bRt^*  ~ &   = & ~ \{ (r, \pi, z) \in Dom(\whPsi)  \mid r \geq 2\} \setminus \{\omega_1 , \omega_2\}  
\end{eqnarray*}
  
 Consider  $\xi, \eta \in \bRt^*$ such that   $\eta$ is
 contained in the forward $\cW$-orbit of $\xi$, then 
 if $r(\xi)>  2$, Propositions~\ref{prop-min} and \ref{prop-propexist} imply that
 $\eta$ is in the forward $\cK$-orbit of $\xi$. Thus,  there exists some $k > 0$ such that $\eta = \whPhi^{k}(\xi)$. The following result is a pseudogroup version of this 
 observation. 
    
 \begin{lemma} \label{lem-psi}
The continuous maps defined by   restriction,
\begin{equation}\label{eq-defpsidomain}
\psi_- = \whPsi_- | \{ U_{\psi}^* \cap  D(\whPsi)_-^- \}  \quad , \quad \psi_0 = \whPsi_0   | \{ U_{\psi}^* \cap  D(\whPsi)_-^+ \} \quad , \quad \psi_+ = \whPsi_+  | \{ U_{\psi}^* \cap D(\whPsi)_+^+ \} ~ , 
\end{equation}
each belong to  $\cGK$ after restricting to the interiors of their domains.
 \end{lemma}
 \proof
Consider first the case of $\psi_-$.  
  We show that for each  $\xi \in U_{\psi}^* \cap \interior(D(\whPsi)_-^-)$ there is an open neighborhood $\xi \in U_{\xi} \subset U_{\psi}^* \cap \interior(D(\whPsi)_-^-)$ such that the restriction $\whPsi \mid U_{\xi}$ is in $\cGK$, then the claim that $\psi_- | \interior(D(\whPsi)_-^-)  \in \cGK$ follows by the gluing condition Condition~\ref{def-pseudogroup}.5.
   Let $[\xi, \eta]_{\cW}$ be the $\cW$-orbit segment from $\xi$ to $\eta = \whPsi(\xi)$. 

  If  $[\xi, \eta]_{\cW}$ contains no transition point,  then it is also a $\cK$-orbit segment, so  $\whPsi(\xi) = \whPhi(\xi)$. 
  As the insertion regions $\cD_i$ are closed,   the continuity of the flows implies    there is an open neighborhood $U_{\xi}\subset U_{\psi}^* \cap  D(\whPsi)_-^-$ of $\xi$ such that  $\whPsi | U_{\xi} = \whPhi | U_{\xi}$. By Condition~\ref{def-pseudogroup}.3,   $\whPsi | U_{\xi}$ is in $\cGK$.

If   $[\xi, \eta]_{\cW} \cap \cD_i$ is non empty, which for this case implies $i=1$, then the first intersection with $\cD_1$  is a point $x' \in  \cL_1^-$. This is  followed   along the  $\cW$-orbit segment $[\xi, \eta]_{\cW}$  by a     point $y'  \in \cL_1^+$ before returning to $\eta \in \bRt$.  Thus, the 
  first transition point in the $\cK$-orbit of $\xi$ must be the secondary entry point $x = \tau(x') \in  E_1$.

As $\xi \prec_{\cW} \eta$ with $x \equiv y$, and $r(\xi) > 2$ is given,
Proposition~\ref{prop-propexist} implies that  $\xi \prec_{\cK} \eta$. Moreover,  from its proof  we   have that  $n_{\xi}(t) \geq 0$ for $0 \leq t \leq t_{\eta}$ where $\eta = \Phi_{t_{\eta}}(\xi)$, and so also  $\rho_{\xi}(t) > 2$ for $0 \leq t \leq t_{\eta}$.
Let $\xi = \xi_0 \prec_{\cK} \xi_1 \prec_{\cK} \cdots \prec_{\cK} \xi_k = \eta$ be the points of intersection of  $\cK$-arc $[\xi,\eta]_{\cK}$   with $\bRt$, then  $\eta = \whPhi^k(\xi)$ and $r(\xi_{\ell}) > 2$ for each $0 \leq \ell \leq k$. 
By the transversality of the flow $\cK$ with the faces  $E_1$ and $S_1$,  it  follows    that 
 there is an open neighborhood $\xi \in U_{\xi} \subset U_{\psi}^*$ such that  $\whPhi^k | U_{\xi}$ is defined, and  thus is the composition of elements in $\cGK$. 
 Thus  $\whPsi | U_{\xi} = \whPhi^k | U_{\xi} \in \cGK$.
 
 The   claim for the    cases $\psi_0$ and $\psi_+$   is shown in the same way.    \endproof

Next we show that the maps in $\cGK$ defined in Lemma~\ref{lem-psi}   admit  continuous extensions to   open neighborhoods of the space $\bRt^*$. This technical result is fundamental for many subsequent applications of the pseudogroup approach to the study of dynamics for the flow $\Phi_t$.

 \begin{lemma} \label{lem-psi-ext}
There exists an open set $U_{\psi} \subset \bRt$  containing $\bRt^*$ such that the restrictions of $\whPsi$ to the domains of continuity in $U_{\psi}$ define    elements of    $\cGK$.
 \end{lemma}
 \proof
 We first consider the case of the restriction $ \whPsi_- | \{U_{\psi}^* \cap  D(\whPsi)_-^-\} $. 
 
For $\xi \in \bRt^* \cap   \interior(D(\whPsi)_-^-)$ with $r(\xi) > 2$,   Lemma~\ref{lem-psi} shows there is an open   neighborhood   $\xi \in U_{\xi} \subset U_{\psi}^*$ such that the restriction   $\whPsi | U_{\xi}$ defines an  element   $\psi \in \cGK$.

For $\xi \in \bRt^* \cap  \interior(D(\whPsi)_-^-)$  with $r(\xi) = 2$, the return map $\whPsi$ is  defined on a sufficiently small  open neighborhood of $\xi$. 
Set $\eta =   \whPsi(\xi)$ and note that $r(\eta) = 2$.  
Let $[\xi, \eta]_{\cW}$ be the $\cW$-orbit segment from $\xi$ to $\eta$. 
  If  $[\xi, \eta]_{\cW}$ contains no transition point,  then this is also a $\cK$-orbit segment, so 
 $\whPsi(\xi) = \whPhi(\xi)$ and there is an open neighborhood $U_{\xi} \subset \bRt$ such that  $\whPsi | U_{\xi} = \whPhi | U_{\xi}$.      Thus,  $\whPsi | U_{\xi} \in \cGK$   defines an extension of the map $\psi$ defined by Lemma~\ref{lem-psi}  to the open neighborhood $U_{\xi}$ of $\xi$.

  Now assume that   $[\xi, \eta]_{\cW} \cap \cD_i$ is non empty, and so $i=1$ in this case. Then the first intersection with $\cD_1$  is a point $x' \in  \cL_1^-$. This is  followed   along the  $\cW$-orbit segment $[\xi, \eta]_{\cW}$  by a     point $y'  \in \cL_1^+$ before returning to $\eta \in \bRt$.  Thus, the 
  first transition point in the $\cK$-orbit of $\xi$ must be the secondary entry point $x = \tau(x') \in  E_1$.
As $\xi \prec_{\cW} \eta$ with $x \equiv y$, and $r(\xi) = 2$ is given, and $\xi \ne \omega_i$, the Radius Inequality implies that $r(x') > 2$.
 Then as above, Proposition~\ref{prop-propexist} implies that  $\xi \prec_{\cK} \eta$.

  Let $\xi = \xi_0 \prec_{\cK} \xi_1 \prec_{\cK} \cdots \prec_{\cK} \xi_k = \eta$ be the points of intersection of  $\cK$-arc $[\xi,\eta]_{\cK}$   with $\bRt$, then $\eta = \whPhi^k(\xi)$. As before, we   have that  $n_{\xi}(t) \geq 0$ for $0 \leq t \leq t_{\eta}$ where $\eta = \Phi_{t_{\eta}}(\xi)$ and so also  $\rho_{\xi}(t) \geq 2$ for $0 \leq t \leq t_{\eta}$.

 Then there is an open neighborhood $\xi \in U_{\xi} \subset \bRt$ such that  $\whPhi^k | U_{\xi}$ is defined and is the composition of elements in $\cGK$. It follows that $\whPsi | U_{\xi} = \whPhi^k | U_{\xi}$, and so  $\whPsi | U_{\xi} \in \cGK$.
The vertical line segments  $\{r=2 ; z \ne \pm 1\} \cap  \interior(D(\whPsi)_-^-)$ thus admit  coverings by open sets  $U_{\xi} \subset \bRt$ such that $\whPhi | U_{\xi}$ is the restriction of an element of $\cGK$. 
Thus by the gluing condition Condition~\ref{def-pseudogroup}.5, there exists an open set   $U_{\psi_-} \subset \bRt  \cap   \interior(D(\whPsi)_-^-)$    containing $\bRt^*  \cap   \interior(D(\whPsi)_-^-)$ for which $\whPsi_- | U_{\psi_-}$ is continuous, and   the restriction of $\psi_- = \whPsi_- | U_{\psi_-}$ defines  a map in  $\cGK$.

The proofs for the other two cases in \eqref{eq-defpsidomain} follow similarly, yielding domains denoted by $U_{\psi_0}$ and $U_{\psi_+}$ such that 
$\whPsi_0 | U_{\psi_0} \in \cGK$ and $\whPsi_+ | U_{\psi_+} \in \cGK$.
  \endproof
 
 We thus obtain the following maps in $\cGK$
 \begin{equation}\label{eq-defpsiextended}
\psi_- = \whPsi_- | U_{\psi_-}  \quad , \quad \psi_0 = \whPsi_0   | U_{\psi_0}  \quad , \quad \psi_+ = \whPsi_+  | U_{\psi_+} ~.
\end{equation}
 
  \begin{remark}\label{rmk-defpsi}
For convenience of notation in subsequent discussions,   we will use the notation $\psi$ for the union of the three maps $\psi_-$, $\psi_0$ and $\psi_+$ in \eqref{eq-defpsiextended}. The domains of these three maps are disjoint, so $\psi \in \cGK$ by   the gluing condition Condition~\ref{def-pseudogroup}.5. Moreover, for  $\xi$ in the domain of $\psi$ the meaning of $\psi(\xi)$ is then clear, as it is defined by exactly one of these maps.
\end{remark}

 We have established the existence of five special elements, $\{\phi_1^+, \phi_1^-, \phi_2^+, \phi_2^-,  \psi\} \subset \cGK$, each of which reflects aspects of the dynamics of the flow $\Phi_t$ which played key roles in    analyzing the dynamics of the flow $\Phi_t$  in the previous sections of this work. Also, let $Id$ denote the identity map on $\bRt$.  We introduce the subset of $\cGK$ generated by the compositions of these words:

  \begin{defn}\label{def-pseudogroup+}
  Let   $\cGK^* \subset \cGK$ denote the  collection of all maps  formed by compositions of the maps $\{Id, \phi_1^+, \phi_1^-, \phi_2^+, \phi_2^-,  \psi\}$ and their restrictions to open subsets in their domains.
\end{defn}
Note that $\cGK^*$ need not satisfy the condition Definition~\ref{def-pseudogroup}.5 on unions of maps, so that $\cGK^*$ is not a sub-pseudogroup of $\cGK$.
The reason for not imposing this condition, is that many dynamical properties for flows admit corresponding versions for local maps defined by compositions of maps in a generating set, but not if we allow for arbitrary unions as in condition (5). This issue is discussed further in 
 Hurder \cite{Hurder2014} and Matsumoto \cite{Matsumoto2010}. Following the convention of \cite{Matsumoto2010}, we say that $\cGK^*$ is a   \emph{$\psg$}, where the ``$\star$'' refers to the definition as the composition of maps.

A set $\cS \subset \mR$ is \emph{syndetic}   if    there exists    $\nu_{\cS} >0$ such that for all $a \in \mR$ the interval $[a, a + \nu_{\cS}]$ satisfies $\cS \cap [a, a + \nu_{\cS}] \ne \emptyset$. The set $\cS$ has \emph{bounded gaps} if there exists   $\mu_{\cS} >0$ so that    $| \cJ | \leq \mu_{\cS}$ for all intervals $\cJ \subset \mR$ with $\cJ \cap \cS = \emptyset$. These two notions are clearly the same.

 Let $x \in \mK$ have  forward trapped orbit and let  $x_{\ell} = \Phi_{t_{\ell}}(x)  \in \cT_{\cK}$, where $0 \leq t_0 < t_1 < \cdots$, be the sequence of transition points. 
    Corollary~\ref{cor-segmentlengths} implies that there is an upper bound on the lengths $|t_{\ell +1} - t_{\ell}|$, so the set of transition times has bounded gaps and  is a syndetic subset of $\{t \geq 0\}$.   The study of the dynamical properties of the $\cK$-orbit in previous sections used this property repeatedly.    
 
 The dynamics of the pseudogroup $\cGK$ acting on $\bRt$ provides a second discrete model for the dynamics of the flow $\Phi_t$. The following result shows that the $\cGK$-pseudogroup dynamics is an accurate model of the $\cK$-orbit dynamics along the infinite orbits in the region $\{ r\geq 2\}$. 
    
     \begin{prop}\label{prop-sigmadyn} 
 Let $\xi \in \bRt$ have  infinite orbit for the flow $\Phi_t$ with $\rho_{\xi}(t) \geq 2$ for all $t$. Then the set  $\ds \cS_{\xi} = \{ s \mid \Phi_s(\xi) \in \bRt \}$ is syndetic, for a constant $\nu_{\cK}$ which is independent of $\xi$. 
     \end{prop}
\proof 

We first make some observations   about transition points and the rectangle $\bRt$.  
   Recall that Lemma~\ref{lem-cases1} lists the   possible cases for a
   $\cW$-arc  $[x,y]_{\cK} \subset \mK$ with lift $[x',y']_{\cW}$ in
   $\wmW$ and transition points $\{x , y\}$.     Figures~\ref{fig:cWarcs} and \ref{fig:KR}    help to visualize   these possible cases. 
 Note that there does not exists a $\cW$-arc  $[x,y]_{\cK}$ for the   case where $x \in S_2$ and $y  \in E_1$. Also,  for $x \in E_1$ and $y \in  S_2$ or   for   $x \in E_2$ and  $y \in  S_1$, again no such $\cW$-arc exists by the facing property of the $\Psi_t$-flow.  
 
 The following   results describe the  cases where a  $\cW$-arc
 $[x,y]_{\cK} \subset \mK$ intersects   the rectangle $\bRt$.

\begin{case}\label{cases2}
 There are   7   cases where a $\cW$-arc  $[x, y]_{\cK}$ always intersects   $\bRt$:
\begin{enumerate}
\item   $x \in E_i$   for $i=1,2$ and $y \in E_1$  
\item $x \in S_1$    and $y  \in E_j$  for $j=1,2$  
\item $x \in S_1$ and $y \in  S_j$   for $j=1,2$  
\item $x \in S_2$ and $y  \in E_2$.  
\end{enumerate}
\end{case}

 \begin{case}\label{cases3}
 There are  6   cases where a $\cW$-arc  $[x, y]_{\cK}$ may not intersect   $\bRt$:
\begin{enumerate}
\item $x \in E_i$ for $i=1,2$ and $y \in E_2$
\item   $x \in E_i$   for $i=1,2$ and $y \in S_i$  (with $x \equiv y$)
\item $x \in S_2$ and $y \in S_i$ for $i=1,2$.
\end{enumerate}
\end{case}

The   cases where $x \in E_i$ for $i=1,2$ and  $y \in E_2$  are notable, for they represent  ``entry/entry'' transitions. 
Similarly, the cases where 
$x \in S_2$ and  $y \in S_i$ for $i=1,2$, represent ``exit/exit'' transitions. 
These $\cW$-arcs   will intersect $\bRt$  if $r(x)$ is sufficiently close to $2$, but may not intersect $\bRt$  for $r(x)$ near the upper limit $R_*$ defined in \eqref{eq-rmax}.

To prove Proposition~\ref{prop-sigmadyn}, let $\xi_0 \in \bRt$ such  that $\xi_1 = \whPhi(\xi_0) = \Phi_{s_1}(\xi_0) \in \bRt$ and   $[\xi_0, \xi_1]_{\cK}$ is a $\cK$-orbit segment containing no interior intersections with $\bRt$. We show that  there is a uniform upper bound, independent of the point $\xi_0$, on the length of  $[\xi_0, \xi_1]_{\cK}$.

   First, suppose that    $[\xi_0, \xi_1]_{\cK}$ contains no
   transition point, then it is a  $\cW$-orbit segment and  so  $\xi_1 = \psi(\xi_0)$. 
    If  $[\xi_0, \xi_1]_{\cK}$ does not contain an interior point with   $z=0$, then the segment $[\xi_0, \xi_1]_{\cK}$ makes one complete revolution in $\mW$ from its start to finish in $\bRt$ and so   admits a uniform upper and lower  bound  on its length.   Otherwise, if  $[\xi_0, \xi_1]_{\cK}$ intersects the annulus $\cA = \{z=0\}$,  the point $\xi_0$ must lie below $\cA$ and flow counter-clockwise until it meets $\cA$. This half of the flow   traverses less than one revolution,  and then flows clockwise to the point $\xi_1$ which must be symmetric with $\xi_0$ with respect to $\cA$. Thus,  there is again a uniform upper bound on the length of the orbit segment between $\xi_0$ and $\xi_1$.  
  Note that, as $\xi_0$ can be arbitrarily close to $\cA$, there is no lower bound on this length.

If $[\xi_0, \xi_1]_{\cK}$ contains transition points, label them  $\{x_1, \ldots , x_k\}$   where  $x_i = \Phi_{t_i}(\xi_0)$ for  $0 < t_1 < \cdots < t_k < s_1$ with $k \geq 1$.   Note that by assumption,       $r(x_i) \geq 2$ for all $1 \leq i \leq k$.   The lengths of the segments  $[\xi_0 , x_1]_{\cK}$ and $[x_k , \xi_1]_{\cK}$ in $\mK$ admit   a uniform upper (and lower) bound as the subsets $\bRt$,  $E_i$ and $S_i$ for $i=1,2$, are compact.

If $[\xi_0, \xi_1]_{\cK}$ contains exactly one transition point $x_1$,
the above discussion shows that its length is bounded, as it is a
union of two $\cK$-orbit segments, both with an endpoint in $\bRt$. Moreover, if $x_1$ is a secondary entry point in $E_i$ then $\xi_1 = \phi_i^+(\xi_0)$ and  if $x_1$ is a secondary exit point in $S_i$ then    $\xi_1 = \phi_i^-(\xi_0)$.

For the remaining cases where $k \geq 2$, note that by Corollary~\ref{cor-segmentlengths}, the lengths of the arcs $[x_i, x_{i+1}]_{\cK}$ for $1 \leq i < k$ admit   uniform upper (and lower) bounds. Thus, it   suffices to show there is an upper bound on the index $k$, independent of the initial point $\xi_0$. 
Also, note that the assumption   $[x_i, x_{i+1}]_{\cK}$ does not intersect $\bRt$   limits the possibilities   for these $\cK$-arcs to the cases listed in Cases~\ref{cases3}. We obtain an upper bound estimate on the number of such $\cK$-arcs which can occur.

\begin{lemma}\label{lem-boundedentries} There exists $N_*$ such that if $[x_1, x_k]_{\cK}$ is a segment of the $\cK$-orbit in $\mK$ with transition points $\{x_1, \ldots , x_k\}$ for $k \geq 2$, each $x_i$ a \emph{secondary entry point}  with $r(x_i) \geq 2$ and $[x_1, x_k]_{\cK} \cap \bRt = \emptyset$, then $k \leq  N_*$. 
\end{lemma}
\proof
By assumption,  $x_i$ is a secondary entry point for each  $1 \leq i < k$, so there exists  $x_i' \in L_i^-$ and $y_i'\neq x_i'$  satisfying $\tau(x_i') =\tau(y_i')= x_i$. 
Moreover,  $y_{i+1}' \in \cL_2^-$ for $1 \leq i < k$,  for     $y_{i+1}' \in \cL_1^-$ implies    $\cW$-arc $[x_i', y_{i+1}']_{\cW} \cap \bRt \ne \emptyset$, contrary to assumption.
Thus, each    $\cW$-arc $[x_i', y_{i+1}']_{\cW}$ flows from the point $x_i' \in \partial_h^- \mW$ to the point $y_{i+1}' \in \cL_2^-$ as so misses the surface $\cL_1^-$.

In fact,  the $\cW$-orbit of each $x_i'$ must rise   sufficiently fast vertically so that it crosses the annulus $\cA$ before intersecting $\bRt$, where it reverses direction and then flows until it terminates at $y_{i+1}' \in \cL_2^-$. We note that this is an \emph{exceptional} condition, in that with some choices of embeddings $\sigma_1$ and $\sigma_2$ it may be impossible to satisfy. 

Note that there exists $\e'' > 0$ such that for all $x'
\in \partial_h^- \mW$ with $2 \leq r(x') \leq 2 + \e''$, the
$\cW$-orbit of $x'$  intersects $\bRt$ before intersecting $\cL_2^-$. Thus,  $r(x_i') \geq 2 + \e''$ for $1 \leq i < k$. In particular, $r(x_1') \geq 2 + \e''$. 

 As $n_{x_1}(t_{i+1}) = n_{x_1}(t_i)+1$ for $1 \leq i < k$, we apply the Radius Inequality recursively to obtain  that $r(x_{i+1}') > r(x_i') > r(x_1') \geq 2 + \e''$. 
 By   Lemma~\ref{lem-maxstep},  the number $k$ of successive secondary
 entry transitions in the $\cK$-orbit segment $[x_1, x_k]_{\cK}$  is bounded above by the constant $N(2+\e'')$ introduced in its proof. Set   $N_* \equiv N(2+\e'')$ and the result follows.
\endproof

\begin{lemma}\label{lem-boundedexits} 
There exists  $N^* > 0$ so that if $[x_1, x_k]_{\cK}$ is a $\cK$-orbit
segment with transition points $\{x_1, \ldots , x_k\}$ for $k \geq 2$,     each $x_i$ is a \emph{secondary exit point}   with $r(x_i) \geq 2$ and $[x_1, x_k]_{\cK} \cap \bRt = \emptyset$ , then $k \leq  N^*$. 
\end{lemma}
\proof
For $1 \leq i < k$, $x_i$ is a secondary exit point so there exists  $x_i' \in \cL_i^+$ satisfying $\tau(x_i') = x_i$.
A $\cW$-arc from $\cL_1^+$ to   $L_j^+$ for $j=1,2$ must intersect $\bRt$ as noted previously, so we must have $x_{i}' \in \cL_2^+$ for $1 \leq i < k$ and thus $y_{i+1}' \in L_2^+$ for $1 \leq i < k-1$ as well. Note that $y_k' \in L_j^+$ for either $j=1$ or $j=2$ is allowed. 

There exists $\e''' > 0$ such that for all $x' \in \cL_2^+$
with $2<r(x')\leq 2+\e'''$, the $\cW$-orbit of $x'$ must intersect
$\bRt$ before it exits $\mW$ through $\partial_h^-\mW$. Thus, we must have    $r(x_i') > 2+\e'''$ for $1 \leq i < k$. 

 Let $0 < t_1 < \cdots < t_k$ be such that $x_i = \Phi_{t_i}(x)$ for $1 \leq i \leq k$.
Then the level function satisfies $n_{x_1}(t_i) = 1-i$ for $0 \leq i \leq k$ and thus $2 < r(x_{i+1}') < r(x_i')$ for all $0\leq i < k$. 
By  the Radius Inequality, using the same argument as in the proof of Lemma~\ref{lem-maxstep}, there exists $N^*$ such that $i \geq N^*$ implies that 
$r(x_i')\leq 2+\e'''$. By the choice of $\e'''$ this implies that $[x_i', y_{i+1}']_{\cK}$ intersects $\bRt$. Thus, we must have $k \leq N^*$ and the result follows. 
\endproof

We now conclude  the proof of Proposition~\ref{prop-sigmadyn},  for  the    cases   $k \geq 2$, 
which involves an analysis of the possible cases which can arise for the  $\cK$-arcs $[x_i, x_{i+1}]_{\cK}$ for $1 \leq i < k$.

  $\bullet$~        If $x_i \in E_1$ then $x_{i+1} \in E_2$ is possible though ``exceptional'', while $x_{i+1} \in S_j$ implies $j=1$ and $x_i \equiv x_{i+1}$.
          
    $\bullet$~    If $x_i \in E_2$ then  $x_{i+1} \in S_1$ is not possible, while  if   $x_{i+1} \in S_2$   then  $x_i \equiv x_{i+1}$.     
 The case   $x_{i+1} \in E_2$ is possible, although   exceptional and then   Lemma~\ref{lem-boundedentries} implies that we can repeat this case  at most $N_*$ times to yield consecutive secondary entry points in $E_2$,  before  there is  a transition point in $S_2$.
 
        Combining these two cases, it follows that if $x_i \in E_j$  then after at most $N_* +1$ transition points, there follows a secondary exit point.

   $\bullet$~     If $x_i \in S_1$ then the $\cK$-arc $[x_i , x_{i+1}]_{\cK}$ must intersect $\bRt$ for all     $x_{i+1}$.

     $\bullet$~       If $x_i \in S_2$ and $x_{i+1} \in S_1$,  then   the $\cK$-arc $[x_i , x_{i+1}]_{\cK}$ must be followed by a $\cK$-arc which intersects $\bRt$.          
  If $x_i \in S_2$ and  $x_{i+1} \in S_2$, then by Lemma~\ref{lem-boundedexits}, this case can      be repeated successively at most $N^*$ times before the path must intersect $\bRt$. 

Combining  all possible cases above, starting with $x_1$, the number of possible cases which can occur is bounded by  $k \leq 3 + N^* + N_*$. 
Thus, there is a uniformly bounded number of transition points
$\{x_1, \ldots , x_k\}$ which can arise in a $\cK$-orbit segment $[x_1, x_k]_{\cK}$ which does not intersect $\bRt$.  
\endproof

The following consequence of Proposition~\ref{prop-sigmadyn} shows that, in essence,  the dynamical properties of $\Phi_t$ restricted to the non-wandering set $\Omega$ are  determined by the action of the pseudogroup $\cGK$ restricted to the compact invariant subset $\ds \Omega \cap  \bRt$.

\begin{cor}\label{cor-min=min}
Let   $x \in \Omega$, then  there exists $-\nu_{\cK} \leq t_x \leq \nu_{\cK}$ such that $\Phi_{t_x}(x) \in \bRt$. 
\end{cor}

 Recall that    the $\psg$ $\cGK^*$  is the subset of $\cGK$  generated by   compositions of  maps in   the collection $\ds \{\phi_1^+, \phi_1^-, \phi_2^+, \phi_2^-, \psi\}$. It is natural to ask if the actions of $\cGK$ and $\cGK^*$ on $\bRt$ are equivalent,  which leads to the consideration of induced maps in $\cGK$ which are not products of these generators. These maps are related to the maps appearing  in the list   Cases~\ref{cases3}.
In  the next result, we show that the dynamics of  $\cGK$ and $\cGK^*$ restricted to the non-wandering set $\Omega$ agree, at least for points $\xi \in \bRt$ with $2 \leq r(\xi) \leq r_e$  where $r_e$ is the ``exceptional radius'' as defined   in \eqref{eq-eradius} below.

Consider the $\cK$-arcs  as in  Cases~\ref{cases3}.1. Given    $x' \in
L_i^-$  with $2<r(x') \leq 2+\e''$ the $\cW$-orbit of $x'$   intersects the
rectangle $\bRt$ before intersecting the surface   $\cL_2^-$, as in
the proof of Lemma~\ref{lem-boundedentries}.  Then this case does not arise    if $2 \leq r(x) < 2+\e''$ for $x = \tau(x')$.

 Consider next the $\cK$-arcs  as in  Cases~\ref{cases3}.3. Recall
 from the proof of Lemma~\ref{lem-boundedexits} that for
 $x'\in \cL_2^+$ with $2<r(x')\leq 2+\e'''$, the $\cW$-orbit of
 $x'$ intersects $\bRt$ before intersecting the exit region
 $\partial_h^-\mW$. 
Then set 
\begin{equation}\label{eq-eradius}
r_e = \min\{2+\e'', 2+\e''' \} .
\end{equation}

\begin{prop}\label{prop-generators}
Let $\xi\in \bRt$   and suppose $\eta=\whPhi(\xi)$, with   $\eta=\Phi_{t_{\eta}}(\xi)$. 
Assume that $2 \leq \rho_{\xi}(t) \leq r_e$ for $0 \leq t \leq t_{\eta}$, 
then  there exists $\phi \in \cG_K^*$ such that $\phi(\xi)=\eta$.
\end{prop}

\proof 
If   the $\cK$-segment $[\xi, \eta]_{\cK}$ contains no transition points, then  $\whPhi(\xi) = \whPsi(\xi) = \psi(\xi)$ and the claim follows. 
So consider the case when the $\cK$-segment $[\xi, \eta]_{\cK}$ contains at least one transition point and let 
      $x_1 = \Phi_{t_1}(\xi)$ with  $0 < t_1  < t_{\eta}$ be the first transition point in the $\cK$-segment $(\xi, \eta)_{\cK}$. Then $x_1$ is a secondary entry or exit point and we consider the possible cases.

Suppose that $x_1$ is a secondary entry point, so the $\cW$-orbit of
$\xi$   intersects one of the surfaces $\cL_i^-$   in a point $y_1'$
with $\tau(y_1') = x_1$. Note that $r(y_1') = r(\xi) \leq r_e $.  Let $x_1' \in L_i^-$ be such that $\tau(x_1') =   x_1$, then
$r(x_1') \leq 2+\e''$. Let $x_2$ be
the next transition point in the $\cK$-orbit of $\xi$, we have the
three following cases:
\begin{itemize}
\item If $x_2\in E_1$, then the $\cW$-arc $[x_1,x_2]_\cK$ intersects
  $\bRt$ at $\eta$ and thus $\eta=\phi_i^+(\xi)$.
\item If $x_2\in E_2$, then $\eta$ belongs also to the $\cW$-arc
  $[x_1,x_2]_\cK$ since $r(x_1')\leq 2+\e''$. Then $\eta=\phi_i^+(\xi)$.
\item If $x_2\in S_i$ for $i=1,2$, then $x_1\equiv x_2$ and  $x_2'$ is
  in the $\cW$-orbit of $\tau^{-1}(\xi)$. In this case it might happen
  that $[x_1,x_2]_\cK$ does not intersects the rectangle $\bRt$. If
  this is the case, we have to
  consider more transitions points. Let  $x_i = \Phi_{t_i}(\xi)$ with
  $0 < t_1 < t_2 <t_3 <
  t_{\eta}$ and $1\leq i\leq 3$, be the first three transition points. 

According to Cases~\ref{cases2} if $x_3$ is a secondary entry point,
the $\cW$-arc $[x_2,x_3]_\cK$ must intersect $\bRt$. Thus $\eta\in
[x_2,x_3]_\cK$ and since $x_2$ is in the Wilson orbit of $\xi$,
Proposition~\ref{prop-propexist} implies that $\eta=\psi(\xi)$. If
$x_3$ is a secondary exit point, the choice of $r_e$ implies that
$[x_2,x_3]_\cK$ must intersect $\bRt$ at $\eta$. Again we conclude
that $\eta=\psi(\xi)$.
\end{itemize}

Next, suppose that $x_1$ is a secondary   exit point,  then the
$\cW$-orbit of $\xi$  exits  through  $\partial_h^+ \mW$ in a point
$y_1' \in L_i^+$   with $\tau(y_1') = x_1$, for $i=1,2$.  Let $x_1'\in
\cL_i^+$ be such that $\tau(x_1')=x_1$.

Assume first that $x_1\in S_1$. Let $x_2$ be the following transition point in the $\cK$-orbit of
$\xi$. Cases~\ref{cases2} imply that the $\cW$-arc $[x_1,x_2]_\cK$
must intersect $\bRt$, and thus $\eta=\phi_1^-(\xi)$.

We are left with the case $x_1\in S_2$. If $x_2$ is a secondary entry
point, then $x_2\in E_2$ and the $\cW$-arc $[x_1,x_2]_\cK$ must
intersect $\bRt$. We conclude that $\eta=\phi_2^-(\xi)$. If $x_2$ is a
secondary exit point, then by the choice or $r_e$ we have that
$r(x_1')<2+e'''$ and thus the $\cW$-arc $[x_1,x_2]_\cK$ must
intersect $\bRt$, implying that $\eta=\phi_2^-(\xi)$.

These   cases exhaust the possibilities for the $\cK$-orbit segment $[\xi, \eta]_{\cK}$ so we have  $\phi(\xi)=\eta$ where $\phi$ is one of the generators of $\cGK^*$. 
 \endproof

\begin{prop}\label{prop-syndeticgk*}
Let $\xi \in \bRt$ have   infinite orbit for the flow $\Phi_t$ with
$r(\xi)<r_e$ and  $\rho_{\xi}(t)\geq 2$ for all $t$. Then the set $\ds
\cS_{\xi}^* = \{ s \mid \Phi_s(\xi) \in \cGK^*(\xi)\subset \bRt ~, ~ r(\Phi_s(\xi)) < r_e\}$ is syndetic, for a constant $\nu_{\cK}^*$ which is independent of $\xi$.
\end{prop}
\proof 
Let   $s_0  \in \cS_\xi^*$ and let $s_1 \in \cS_\xi^*$ satisfy  $s_0 < s_1$, such that $s_1$ is the least such value. We need to show that there exists a value $\nu_{\cK}^*$ independent of $s_0$ such that $s_1 - s_0 \leq \nu_{\cK}^*$.

It is given that  $\xi \in \bRt$ has  infinite orbit with $\rho_{\xi}(t) \geq 2$ for all $t$, so   the set $\ds \cS_{\xi} = \{ s \mid \Phi_s(\xi) \in \bRt \}$  is syndetic in $\mR$ with constant $\nu_{\cK}$ by  Proposition~\ref{prop-sigmadyn}. Also,  by definition  there is an  inclusion $\cS_\xi^*\subset \cS_\xi$.
Consider  successive points $s_0 , s_1 \in \cS_\xi^*$ with $s_0 < s_1$. It suffices to show that there exists $C(r_e) > 0$, depending on $r_e$ but independent of $\xi$, so that $[s_0, s_1] \cap \cS_{\xi}$ contains  at most $C(r_e)$ points. We show this first for the case $s_0 \geq 0$.

Suppose that  for all $s_0 \leq t \leq s_1$ we have $2\leq \rho_{\xi}(t)\leq r_e$. Then by Proposition~\ref{prop-generators}, the intersection $\cS_\xi \cap (s_0,s_1)$ is empty.

If $\rho_{\xi}(t)\leq r_e$ for all $ t\geq 0$, then the above argument shows that $[0, \infty) \cap \cS_\xi = [0, \infty) \cap \cS_\xi^*$.
Otherwise, 
let $t_1>0$ be the least time for which $\rho_{\xi}(t_1)\geq r_e$. Then $x_1=\Phi_{t_1}(\xi)$ is a secondary entry point. 
 By Proposition~\ref{prop-propexist}   there exists $T_{x_1}>0$ such
that $\ox_1 = \Phi_{t_1+T_{x_1}}(\xi)$ is the secondary exit point facing
$x_1$. Thus for $s=T_{x_1}+\e$ and $\e>0$ small,
$\rho_{\xi}(t_1+s)<r_e$. 

Consider the backwards flow of the   secondary entry point $x_1$  to the first  intersection of the orbit with $\bRt$ to obtain a point  $\xi_1 = \Phi_{u_1}(\xi) \in \bRt$ where $0 \leq u_1 < t_1$. By the choice of $t_1$ and the remarks above, we have that   $\xi_1\in \cGK^*(\xi)$ and so   $u_1 \in \cS_\xi^*$. Also, note that $t_1 - u_1 < \nu_{\cK}$ as $\cS_\xi$  is syndetic for the constant $\nu_{\cK}$.

 Next, let $\xi_2 = \Phi_{u_2}(\xi)$  be the first intersection with $\bRt$ of the forward orbit of   $\ox_1 = \Phi_{t_1+T_{x_1}}(\xi)$, so that   $t_1+T_{x_1} < u_2$ and   $u_2 - (t_1+T_{x_1}) < \nu_{\cK}$. 
 Since $\rho_{\xi}(t_1+T_{x_1}+\e )<r_e$ for $\e > 0$ sufficiently small, by the choice of $r_e$ in \eqref{eq-eradius},  there are no transition points in the $\cK$-orbit segment between $\Phi_{t_1+T_{x_1}}(\xi)$ and $\xi_2$. 
Then by    Lemma~\ref{lem-shortcut},  we have that $\xi_1 \prec_{\cW} \xi_2$, which implies that   $\xi_2=\psi(\xi_1)$ for the generator $\psi \in \cGK^*$. Thus,  $\xi_2\in  \cGK^*(\xi)$ with $r(\xi_2) < r_e$. 
Observe that $u_2 - u_1 <    2\nu_{\cK} + T_{x_1}$.

Lemma~\ref{lem-maxstep} and
Corollary~\ref{cor-lengths} imply that there exists $T_{r_e} > 0$ such that    $T_{x_1} \leq T_{r_e}$   for all $x_1$ with $r(x_1) \geq r_e$. 
  
  Define $\nu_{\cK}^* =  2\nu_{\cK} + T_{r_e}$. 
We can then 
apply the above process  recursively along the forward $\cK$-orbit of $\xi$ to obtain that $[0,\infty) \cap \cS_{\xi}^* $ is syndetic in $[0,\infty)$ for the constant $\nu_{\cK}^*$. 
 The conclusion for the backward flow follows by reversing the time parameter as usual. 
\endproof
 
  \bigskip
  
\section{The level decomposition} \label{sec-fM}

We now begin the  study of the dynamics of the Kuperberg flow $\Phi_t$ from a more topological point of view. Many of the   results of the previous sections are   given topological interpretations in this approach, which culminates in Section~\ref{sec-zippered} with a description of    the dynamics of the flow in terms of the structure of the ``zippered lamination'' $\fM$     containing  the minimal set $\Sigma$.

\begin{figure}[!htbp]
\centering
{\includegraphics[width=80mm]{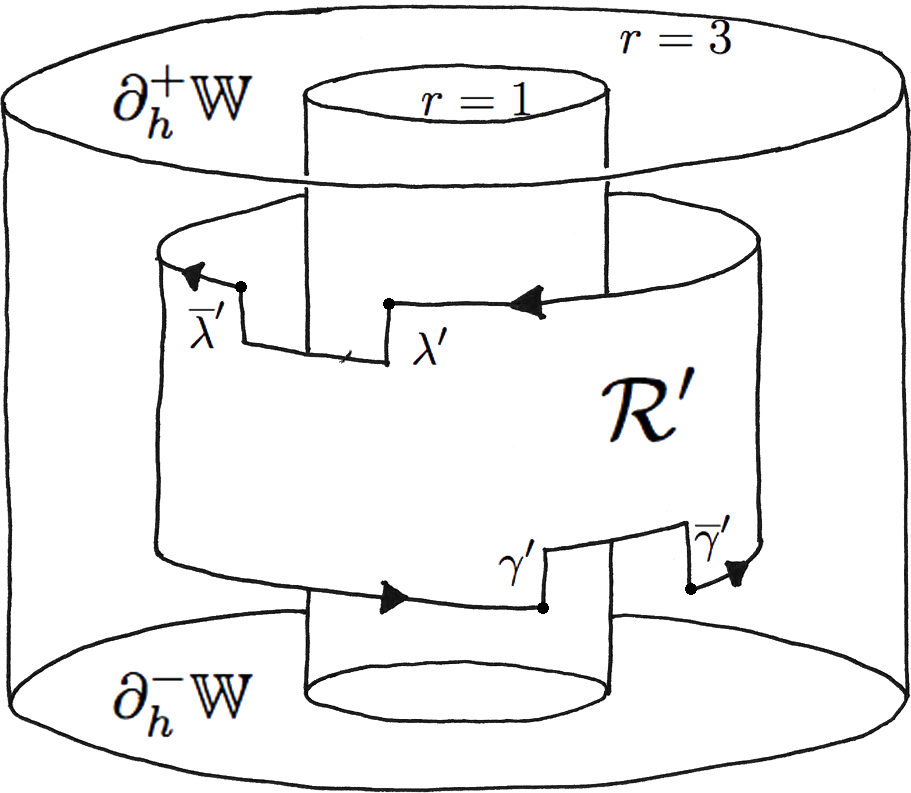}}
\caption{\label{fig:notches} The notched cylinder $\cR'$ embedded in $\mW$}
\vspace{-6pt}
\end{figure}

Recall that the periodic orbits $\cO_i$ of the Wilson flow are the
boundary circles for the Reeb cylinder   $\cR  \subset \mW$.  We introduce the 
 \emph{notched Reeb cylinder},  $\cR' = \cR \cap \mW'$, which has two
closed ``notches''   removed from   $\cR$ where it intersects the
closed insertions $\cD_i \subset \mW$ for $i=1,2$.
Figure~\ref{fig:notches}  illustrates the cylinder $\cR'$ inside
$\mW$.

Consider  the $\cK$-orbit of the image   $\tau(\cR') \subset \mK$ and its closure
 \begin{equation}\label{eq-lamination}
\fM_0 ~ \equiv ~  \{\Phi_t(\tau(\cR')) \mid -\infty < t < \infty \} \quad  , \quad  \fM ~ \equiv ~ \overline{\fM_0}     \subset \mK ~.
\end{equation}
Note that since the special points $p_i^{-} = \tau(\cL_i^- \cap \cO_i)$ for $i=1,2$ are contained in      $\fM_0$,  we have $\Sigma \subset \fM$. Moreover, the flow $\Phi_t$ restricts to a flow on $\fM$, so the dynamics of $\Phi_t$ on $\Sigma$ is defined by the restricted dynamics of the flow on $\fM$.

  The horizontal segments of the boundary of $\cR'$ which are not in
  the notches, are the orbit
  segments $\cO_i\cap \mW$ of $\cW$, while the vertical segments labeled $\g'$, $\ovg'$, $\lambda'$ and $\ovl'$ 
  in Figure~\ref{fig:notches}   are transverse to the flow of $\cW$.  
  If $x \in \tau(\gamma')$ or $x \in \tau(\lambda')$ is on a vertical
  segment,  the $\cK$-orbit of $x$ points into the notch. If in
  addition $z(x) \ne \pm 1$, 
  then by Proposition~\ref{prop-propexist},  the forward $\cK$-orbit of $x$ returns in finite time to the opposite vertical boundary $\tau(\ovg')$,  or $\tau(\ovl')$ respectively, of the notch. 
  For each $0 < \e < 1$, introduce the compact set $\fM_0^\e$ as follows.
Observe that for $x\in \g'$  with $z(x)>-1$ Proposition~\ref{prop-propexist} implies that  there is a finite
time $T_x$ such that $\Phi_{T_x}(\tau(x))\in \tau(\overline{\g}')$. Analogously, for $x\in \lambda'$
with $z(x)<1$ there is a finite
time $T_x$ such that $\Phi_{T_x}(\tau(x))\in \tau(\overline{\lambda}')$. Let
$$ \fM_0^\e      =   \tau(\cR')   \cup    \{\Phi_t(x)\,|\, x\in \tau(\g') ~ , ~  z(x)\geq -1+\e ~ ,~ 0\leq t\leq T_x\} 
     \cup     \{\Phi_t(x)\,|\, x\in \tau(\lambda') ~ , ~ z(x)\leq 1-\e ~ , ~ 0\leq t\leq T_x\}. $$

The regions added to $\tau(\cR')$ 
``fill the gap'' made by the notches in $\tau(\cR')$, and yield the compact
surface with boundary $\fM_0^\e$ embedded in $\mK$ (the embedding $\tau$ of
$\cR'$ in $\mK$ is illustrated in Figure~\ref{fig:notched8}). The time $T_x\to
\infty$ as   $x$ tends to the endpoint of $\g'$ or $\lambda'$,   and the boundary curves of the compact surface $\fM_0^\e$ become increasing complicated, though the surface remains embedded at all times.
 The limit of these compact surfaces gives approximations to the set $\fM_0$.

The properties of the level function $n_x(t)$ for the $\cK$-orbit of   $x \in \mK$ were  fundamental for the analysis   of the orbital dynamics of the Kuperberg flow in previous sections. We next show that the level function along orbits yields a well-defined function on $\fM_0$ and use it to introduce the \emph{level  decomposition} of $\fM_0$. The geometry of $\fM_0$ is quite complicated to visualize. Beginning in   Section~\ref{sec-intropropellers},  we  formulate the structure of each of the components in the level decomposition in terms of   \emph{propellers} and develop a precise description of these sets.

For each $x \in \fM_0$ there exists some   $y' \in  \cR'$ for which $y \in \tau(\cR')$ such that $x = \Phi_{t}(y)$ for some $t \in \mR$. The point $y$ is not unique, but as any such choice satisfies $r(y) = 2$,  the proof of Proposition~\ref{prop-r=2any} implies that  $r(x) \geq 2$. The following result is a consequence of this observation and   previous results.
\begin{prop}\label{prop-levelsM}\label{prop-levels}
There is a well-defined level function 
\begin{equation}
n_0 \colon \fM_0 \to \mN = \{0,1,2,\ldots\} .
\end{equation}
\end{prop}
\proof
First, for $x \in \tau(\cR')$ set $n_0(x) = 0$.

Let $x \in \fM_0$, then there exists $y \in \tau(\cR')$ such that $x =
\Phi_{s_0}(y)$. Define $n_0(x) = n_{y}(s_0)$.
Note that we allow both positive and negative values for $s_0$, using either formulas \eqref{def-level+} or \eqref{def-level-}. 

  Let $w \in \tau(\cR')$ be another point such that $x = \Phi_{s_1}(w)$. Assume without loss of generality that $y\prec_{\cK} w$ so that there exists $T > 0$ with $w = \Phi_T(y)$. Note that $n_{y}(s + T) = n_{w}(s) + n_{y}(T)$, so it suffices to show that $n_y(T) = 0$.
  
  Suppose that  $y \in \tau(\cR')$ where $y = \tau(y')$ with   $z(y') = \pm 1$, so that $y$ is on a special orbit. Then for $w \in \tau(\cR')$ with  $y\prec_{\cK} w$ we must have $w=\tau(w')$ where $z(w') = z(y')$. This implies    there are no transition points in the $\cW$-orbit  between $y'$ and $w'$, thus $n_y(T) = 0$. 
  
Now suppose that $z(y') \ne \pm 1$, and let    $0 \leq  t_0 < t_1 < \cdots < t_k \leq T$ with $x_{\ell} = \Phi_{t_{\ell}}(y)$ be  the transition points for the $\cK$-orbit  between $y$ and $w$, for $k \geq 0$.  The case $k=0$ is again immediate, so assume that $k \geq 1$. We show that $n_{y}(T) = 0$, which will follow using an inductive argument  on the number of transition points $k$ between $y$ and $w$. 

Note that $y' \in \cR'$ and $z(y') \ne \pm 1$ implies that the first transition point $x_0 \in E_i$ for $i=1,2$, so $n_y(t_0) = 1$ and satisfies $r(x_0') > 2$. 
Let $\ell > 0$ be the least integer such that $n_y(t_{\ell}) = 0$, which exists by the proof of Proposition~\ref{prop-propexist}. Then $x_{\ell}$ is an exit point with $x_0 \equiv x_{\ell}$ and   so $x_{\ell} \in \tau(\cR')$. Choose $u = \Phi_{t'}(y)$ on the $\cK$-orbit between $x_{\ell}$ and $x_{\ell +1}$, so that   $n_y(t') = n_y(0) = 0$ and we have reduced the problem to showing that   $n_{u}(T-t') = 0$, where the $\cK$-orbit segment $[u, w]_{\cK}$ now has $k-\ell -1$ transition points. The claim now follows by induction. 
\endproof

Define the  \emph{level sets} of  $\fM_0$ as follows: 
\begin{equation}\label{eq-levelsets}
\fM_0^n ~ = ~ \{ x \in \fM_0 \mid n_0(x) \leq n\} ~, ~ n = 0,1,2, \ldots
\end{equation}
The set $\fM_0^0$  contains the   notched Reeb cylinder $\tau(\cR')$ and   the level $0$ points in the boundary of the notches in $\tau(\cR')$. 
 The descriptions of the sets  $\fM_0^n$ for $n > 0$ is more subtle and will follow from a detailed study of the $\cK$-orbit   of the vertical segments $\gamma'$ and $\lambda'$  in     Section~\ref{sec-proplevels}.

  \bigskip
  
\section{Embedded surfaces and propellers}\label{sec-intropropellers}

The invariant set $\fM_0$  is the union of its level sets,  as defined in \eqref{eq-levelsets}. While the set 
  $\fM_0^0$   as described in the last section is rather simple, the description of the level sets $\fM_0^n$ for $n > 0$ leads to the introduction of one of the main ideas of this work, the notion of \emph{finite} and \emph{infinite}   propellers,  which   are obtained  from the   flow $\Psi_t$  of  selected arcs in $\mW$. These are  defined and studied in this section. For example, $\fM_0^1$ is obtained from $\fM_0^0$ by attaching 
   two non-compact surfaces which  are  infinite propellers,  while the    level sets $\fM_0^n$ for $n > 1$ are obtained
  by attaching  finite propellers to $\fM_0^{n-1}$, where the   complexity of these added finite propellers at level $n$ increases  as $n \to \infty$.

Consider a   curve in the entry region, $\g \subset \partial_h^- \mW$, with a parametrization  $w(s) = (r(s), \theta(s), -2)$ for $0 \leq s \leq 1$, 
 and assume that the map $w \colon [0,1] \to \partial_h^- \mW$ is a homeomorphism onto its image. We use the notation $w_s = w(s)$ when convenient, so that 
  $w_0 = w(0)$ denotes the initial point and $w_1 = w(1)$    the terminal point  of $\g$. 
For    $\e > 0$,  assume    that $r(w_0) = 3$,   $r(w_1)  = 2+\e$ and $2 +  \e < r(w_s) < 3$ for $0 < s < 1$.

The $\cW$-orbits of the
points in $\g$ traverse $\mW$ from $ \partial_h^- \mW$ to $ \partial_h^+ \mW$ and hence the flow of $\g$ generates a compact invariant surface $P_{\g} \subset \mW$.
The surface $P_{\g}$ is  parametrized by $(s,t) \mapsto \Psi_t(w(s))$ 
for $0 \leq s \leq 1$ and $0 \leq t \leq T_s$, where $T_s$ is the exit time for the $\cW$-orbit of $w(s)$. Observe that as  $s \to 1$ and $\e \to 0$,  Corollary~\ref{cor-lengths} implies that the exit time $T_s    \to \infty$.

\begin{figure}[!htbp]
\centering
\begin{subfigure}[c]{0.4\textwidth}{\includegraphics[height=54mm]{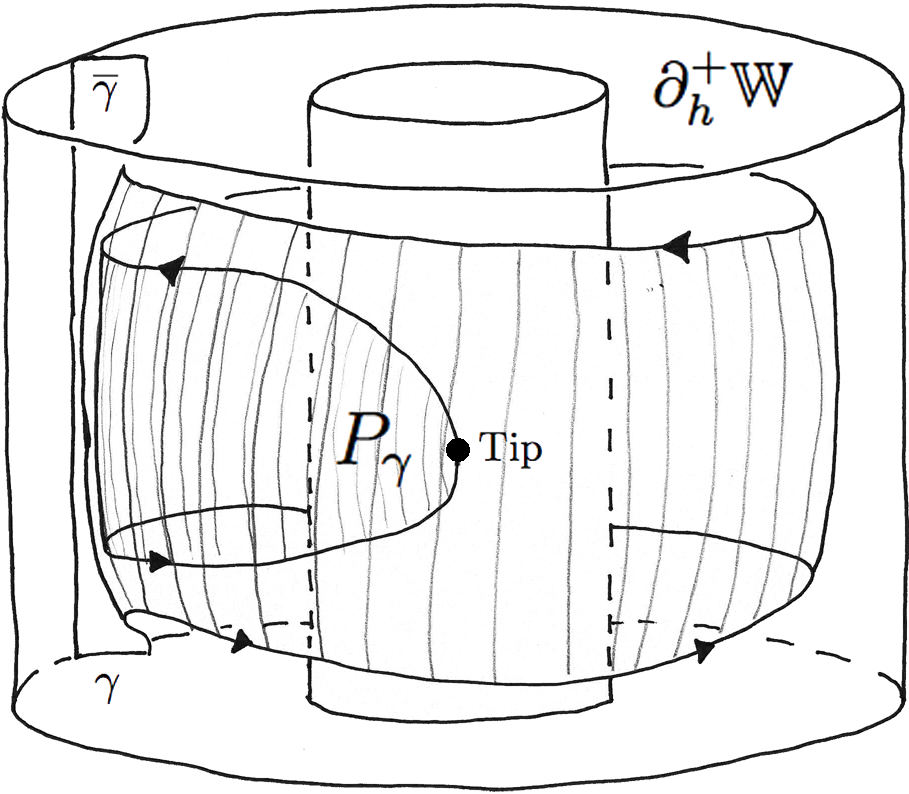}}\end{subfigure}
\begin{subfigure}[c]{0.5\textwidth}{\includegraphics[height=48mm]{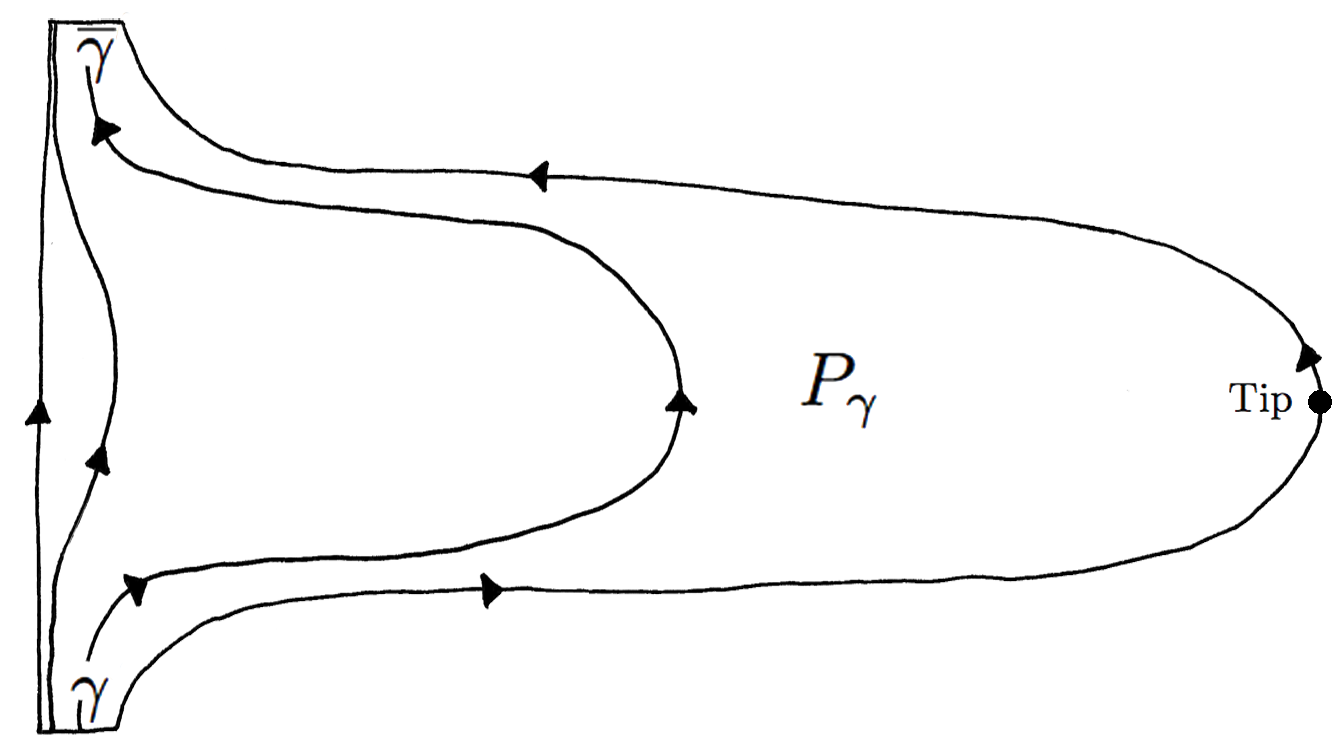}}\end{subfigure}
\caption{\label{fig:propeller}   Embedded and flattened finite propeller}
\vspace{-6pt}
\end{figure}

The surface  $P_{\g}$ is called a {\it propeller}, due to the nature of its shape in $\mR^3$. It takes the form of a ``tongue'' wrapping around the core cylinder $\cC(2 + \e)$ which contains the orbit of $w_1$. To visualize the shape of this surface, consider the   case where $\g$
 is topologically transverse to the cylinders $\cC(r_0) = \{r=r_0\}$ for $2+ \e \leq r_0 \leq 3$. 
The transversality assumption implies that the radius  $r(w_s)$ is \emph{monotone decreasing} as $s$ increases.
   Figure~\ref{fig:propeller}  illustrates the surface $P_{\g}$ as  a ``flattened'' propeller on the right and   its embedding in $\mW$ on the left. As $\e \to 0$ the surface  approaches  the   cylinder $\cC= \{r=2\}$ in an infinite spiraling manner.

We comment on the details in Figure~\ref{fig:propeller}.
The horizontal boundary   $\partial_h  P_{\g}$ is composed of the
initial curve   $\g \subset \partial_h^- \mW$ and its      mirror
image   $\ovg \subset \partial_h^+ \mW$ via the entry/exit condition
on the Wilson Plug. The vertical boundary   $\partial_v P_{\g}$ is
composed of the vertical segment $w_0 \times [-2,2]$ in $\partial_v \mW$ and the
orbit  $\{\Psi_t(w_1) \mid 0 \leq t \leq T_1\}$ which is the inner (or long) edge in the interior of $\mW$.
One way to visualize the surface,   is to consider the product surface $\g \times [-2,2]$ and  then start deforming it by an isotopy which follows the flow lines of $\cW$, as illustrated in Figure~\ref{fig:propeller}.
In the right hand side of the figure,  some of the orbits in the propeller
are presented, while in the left hand side,  just the boundary orbit is
presented.

 Consider the orbit $\{\Psi_t(w_1) \mid 0 \leq t \leq T_1\}$  of the endpoint $w_1$ with $r(w_1) = 2+\e$.
  The path $t \mapsto \Psi_t(w_1)$    makes a certain
number of turns in the positive $\mS^1$-direction before reaching
the core annulus $\cA$ at $z = 0$. The Wilson vector field $\cW$ is vertical on the plane $\cA$,     so the flow of $w_1$ then crosses   $\cA$, after which the orbit $\Psi_t(w_1)$ starts turning in the negative direction and ascending until it reaches $\partial_h^+ \mW$. 
The point where the flow $\Psi_t(w_1)$ intersects $\cA$ is called the \emph{tip} of the propeller $P_{\g}$.

The anti-symmetry of the vector field $\cW$ implies that the number of turns in one direction (considered as a
real number)   equals   the number of turns in the opposite direction. To be precise, for $w_1 = (r_1, \theta_1, -2)$,  let $\Psi_t(w_1) = (r_1(t), \theta_1(t), z_1(t))$ in coordinates. The function $z_1(t)$ is monotone increasing and by the symmetry, we have $z_1(T_1 /2) = 0$. Thus, the tip is the point 
$\Psi_{T_1 /2}(w_1)$.

 Set  $\Theta(2+\e) = \overline{\theta_1}(T_1/2) - \overline{\theta_1}(0)$, where   $\overline{\theta_1}(t)$ is a \emph{continuous} function with $\theta_1(t) = \overline{\theta_1}(t) \mod (2\pi)$. 
 The function $\Theta(2+\e)$ measures the total angle advancement of
 the curve $t \mapsto \Psi_t(w_1)$ between the initial point $w_1$ and
 the tip $\Psi_{T_1 /2}(w_1)$. 
 The number $\Theta(2+\e)$ depends only on the radius $2+\e$ of the   point $w_1$, as the flow $\Psi_t$ is rotationally symmetric by Proposition~\ref{prop-wilsonproperties}.
Note that $\Theta(2+\e) \to \infty$ as $\e \to 0$.
Also, introduce the notation 
\begin{equation}\label{eq-Delta}
\Delta(2+\e) = \lfloor \Theta(2+\e)/2\pi \rfloor ~ .
\end{equation}
 for the   integer part of $\Theta(2+\e)/2\pi$, which is the number of times that the curve  $\Psi_{t}(w_1)$for $0 \leq t \leq T_1/2$ makes a complete  circuit around the cylinder $\cC(2+\e)$.

 Next,  for fixed $0 \leq a < 2\pi$, consider the intersection of $P_{\g}$ with a slice 
 \begin{equation}\label{eq-slices}
\bRa \equiv \{\xi = (r, a, z) \mid ~ 1 \leq r \leq 3 ~,   ~ -2 \leq z \leq 2\} ~  .
\end{equation}
The rectangle $\bRt$ defined by \eqref{eq-goodsection} in Section~\ref{sec-pseudogroup} corresponds to the value   $a=\pi$.
Each rectangle $\bRa$ is tangent to the Wilson flow along the annulus $\cA$ and also near the boundary $\mW$, but is transverse to the flow at all other points. The case when  $a = \Theta(2 + \e) \mod (2\pi)$ is special,   as  the tip of the propeller is tangent to $\bRa$.

\begin{figure}[!htbp]
\centering
\begin{subfigure}[c]{0.4\textwidth}{\includegraphics[width=50mm]{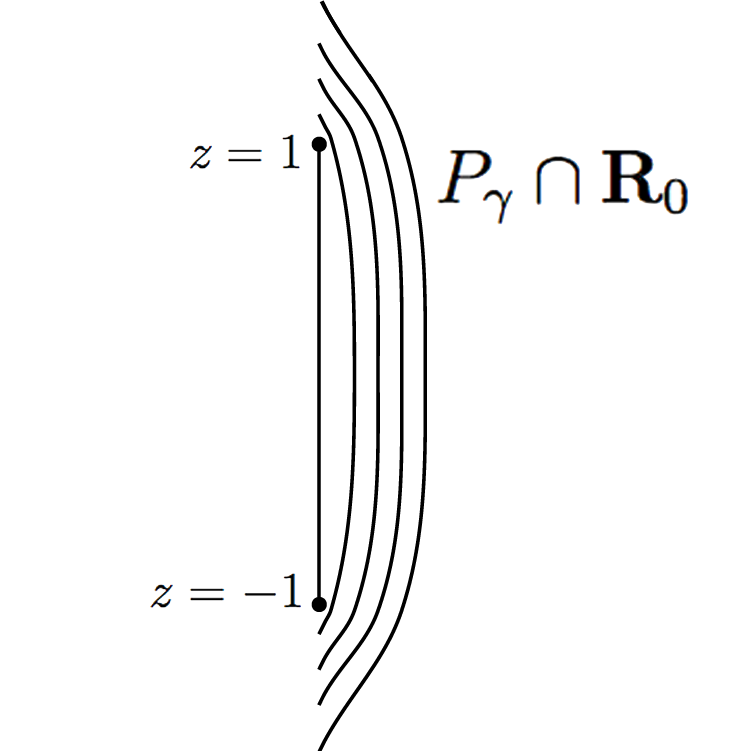}}\caption{Infinite propeller $P_{\gamma}$}\end{subfigure}
\begin{subfigure}[c]{0.4\textwidth}{\includegraphics[width=50mm]{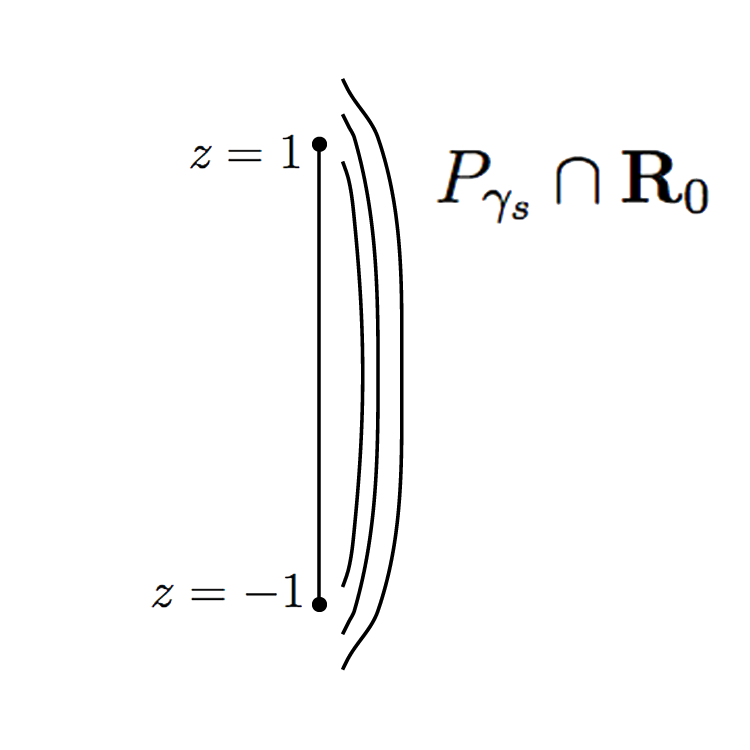}}\caption{Finite propeller $P_{\gamma_s}$}\end{subfigure}
\caption{\label{fig:arcspropeller} Trace of   propellers   in $\bRt$}
\vspace{-6pt}
 \end{figure}

Assume that $a \ne  \Theta(2 + \e) \mod (2\pi)$, then the flow
$\Psi_t(w_1)$ intersects $\bRa$ in a series of points
on the line $\cC(2+\e) \cap \bRa$ that are paired, as illustrated in Figure~\ref{fig:arcspropeller}. 
Moreover, the intersection $P_{\g} \cap \bRa$ consists of a finite sequence of arcs between the symmetrically paired points of $\Psi_t(w_1) \cap \bRa$.
 The number of such arcs is equal to   $\Delta(2+\e) \pm 1$.

 We comment on the details of Figure~\ref{fig:arcspropeller}.
 The vertical line between the points $(2,a,-1)$ and $(2,a,1)$ (marked
 in the figure simply by $z=-1$ and $z=1$, respectively) is the trace
of the Reeb cylinder in $\bRa$. The trace of a propeller in $\bRa$ is a collection of arcs that
have their endpoints in the vertical line $\{r=2 + \e\}$. In the left
hand figure, $r(w_1)=2$ and the propeller in consideration is
infinite. The curves form an infinite family, here just four arcs are shown,
accumulating on the vertical line. The right hand figure illustrates the case  $r(w_1)>2$  and the propeller is
finite.

\begin{figure}[!htbp]
\centering
{\includegraphics[height=72mm]{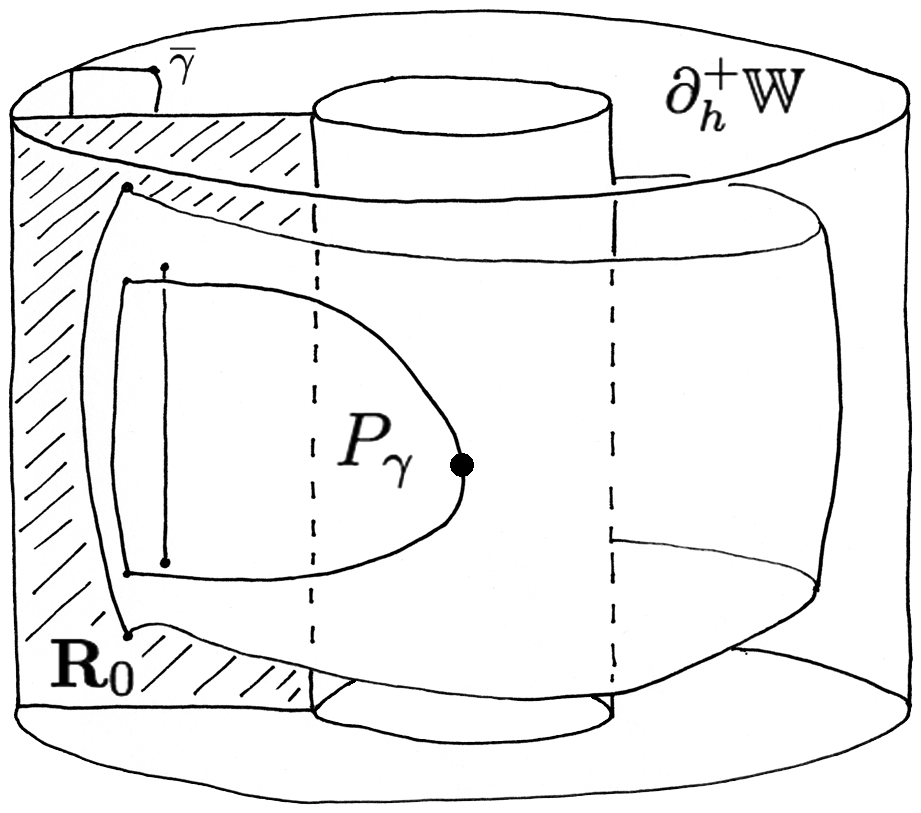}}
\caption{\label{fig:propellerslices} Intersection of finite propeller $P_{\gamma}$ with  $\bRt$}
\vspace{-6pt}
 \end{figure}

\begin{figure}[!htbp]
\centering
{\includegraphics[width=40mm]{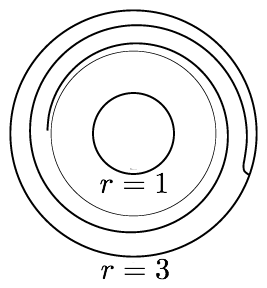}}
\caption{\label{fig:spiralpropeller} Intersection of a finite propeller $P_\g$ with $\cA$}
 \vspace{-6pt}
 \end{figure}

Let us now describe the intersection of $P_{\g}$ with the   annulus $\cA = \{ z=0\}$, which traces out the midpoints of the curves $t \mapsto \Psi_t(w(s))$ for $0 \leq s \leq 1$. That is, the intersection  is the curve $s \mapsto  \Psi_{T_s /2}(w(s)) \in \cA$, which  is a spiral, turning in the positive $\mS^1$ direction $\Delta(2+\e)$ times around the core circle, as in Figure~\ref{fig:spiralpropeller}. 
The point of the curve on the boundary of $\cA$ in Figure~\ref{fig:spiralpropeller} corresponds to the orbit  $\Psi_{t}(w_0) $  and
the point at the end of the spiral    closest to the circle $\{r=2\}$ corresponds to the tip of
$P_{\g}$. If we change  $\cA$ for another annulus 
$\cA_b=\{z=b\}$ where $b \ne 0$, the intersection   $P_{\g}  \cap \cA_b$ will be a spiral
turning in the positive $\mS^1$-direction that is strictly shorter than
the one in $\cA$. By symmetry, the traces of $P_{\g}$ on $\cA_{-b}$ and $\cA_b$ are the same.

Finally,  consider the case where $\e = 0$, so that    the endpoint $w_1$ of   $\g$
lies in  the cylinder $\cC$.  Then for $0 \leq s < 1$,  we have $r(w(s)) > 2$, so the $\cW$-orbit of 
$w_s \in    \partial_h^- \mW$   escapes from $\mW$. Define the curve $\ovg$ in $\partial_h^+ \mW$ to be the trace of these facing endpoints in $\partial_h^+ \mW$, parametrized by $\ow(s)$ for $0 \leq s < 1$, where $w(s) \equiv \ow(s)$. 
Define $\ds \ow_1 = \lim_{s\to 1} ~ \ow(s)$ so that $\ow_1 \equiv w_1$ also. 

Note that  the forward $\Psi_t$-orbit of $w_1$ is   asymptotic to the periodic orbit
$\cO_1$,  while the backward   $\Psi_t$-orbit of $\ow_1$ is   asymptotic to the periodic orbit
$\cO_2$. Introduce their ``pseudo-orbit'', 
\begin{equation}\label{eq-Zset}
\cZ_{\g} = \cZ_{\g}^- ~ \cup ~ \cZ_{\g}^+ ~ , ~ \cZ_{\g}^- =
\{\Psi_t(w_1) \mid t \geq 0 \}  ~ \text{and} ~ \cZ_{\g}^+ =  \{\Psi_t(\ow_1) \mid t \leq 0 \}
\end{equation}
Each curve $\cZ_{\g}^{\pm}$ traces out a semi-infinite ray in $\cC$ which spirals from the bottom or top face to a periodic orbit, and thus 
$\cZ_{\g}$ traces out two semi-infinite curves in $\cC$ spiraling to the periodic orbits $\cO_1 \cup \cO_2$. 
 
For $0 < \e \leq 1$,  denote by $\g^{\e}$ the curve with image $w([0, \e])$, parametrized by 
\begin{equation}\label{eq-sparem}
w^{\e}(s) = w(\e \cdot s) . 
\end{equation}

\begin{defn}\label{def-infpropeller}
Let $\g$ be a curve parametrized by $w \colon [0,1]\to \mW$   as above, with $r(w_0) = 3$ and $r(w_1)=2$.   Introduce the \emph{infinite propeller} and its closure in $\mW$:
\begin{equation}\label{eq-infpropeller}
P_{\g} ~ \equiv ~  \cZ_{\g} ~ \cup ~ \bigcup_{\e > 0} ~ P_{\g^\e}   \quad , \quad \overline{P}_{\g} ~ \equiv ~ \overline{\bigcup_{\e > 0} ~ P_{\g^\e}}
\end{equation}
\end{defn}

\begin{prop} \label{prop-propellerclosure}
The closure $\overline{P}_{\g}$ of an infinite propeller contains the Reeb cylinder $\cR$,    with
\begin{equation}
\overline{P}_{\g} ~ = ~ P_{\g} ~ \cup ~ \cR   ~ .
\end{equation}
In particular, the   periodic orbits $\cO_{1}$ and $\cO_2$ are contained in $\overline{P}_{\g}$. 
\end{prop}
\proof 
The endpoints of the arcs in $P_{\g^\e} \cap \bRa$ tend to points in $\cC \cap \bRa$ as $\e \to 0$, as illustrated in  Figure~\ref{fig:arcspropeller}, so the closure $\overline{P}_{\g}$ contains the set $\cZ_{\g}$. Thus, 
the  intersection  $\{\Psi_t(w_s) \mid 0 \leq t \leq T_s\} \cap \bRa$ for $s \to 1$ contains pairs of points arbitrarily close to the intersections $\cO_i \cap \bRa$ for $i=1,2$. The family of arcs connecting these points  is nested with the vertical line $\cR \cap \bRa$ as the inner boundary. Thus,  $P_{\g}  \cap \bRa$ contains arcs joining these points which are arbitrarily close to  $\cR \cap \bRa$ and so  $\cR \cap \bRa \subset \overline{P}_{\g} \cap \bRa$.
\endproof

Note that the first-half $\cW$-orbit segments $\{\Psi_t(w_s) \mid 0 \leq t \leq T_s/2\}$ for $s \to 1$ have $\cO_1$ in their closure and    the second-half $\cW$-orbit segments $\{\Psi_t(w_s) \mid T_2/2 \leq t \leq T_s\}$ for $s \to 1$ have $\cO_2$ in their closure.

  \bigskip
 
\section{Propellers and the level decomposition}\label{sec-proplevels}

In this section, we analyze the structure of the level decomposition of $\fM_0 \subset \mK$ introduced in Section~\ref{sec-fM}, using the concepts of propellers   introduced in the last section.  A key point is the introduction of an ``intuitive'' labeling system for the collection of propellers generated by the flow of the notched Reeb cylinder. This labeling gives order to the propellers at each level, which grow in number at an exponential rate.  A second key point is the beginning of the study of the dynamics of the $\cK$-orbit in terms of the action of the pseudogroup $\cGK$ on $\bRt$.   In order to eliminate various pathologies that can arise in the study of this action, we first impose ``generic'' properties of the choices made in the construction of   the Wilson flow $\cW$  and on the insertions $\sigma_i$ for $i=1,2$. 
For this purpose, we   formulate    additional   assumptions on the insertions $\sigma_i$ for $i=1,2$, which yield a stronger form of the Radius Inequality (K8) introduced by Kuperberg in \cite{Kuperberg1994}.

 Let $(r', \theta',z')  = \sigma_i(x) \in \cD_i$ for    $i=1,2$, where     $x = (r, \theta,  z) \in  D_i$ is a point in the domain of $\sigma_i$.
Let $\pi_z(r', \theta', z') = (r', \theta',-2)$ denote the projection of $\mW$ along the $z'$-coordinate. 
We assume that $\sigma_i$ restricted to the bottom face, $\sigma_i \colon L_i^- \to \mW$,  has image  transverse to the vertical fibers of $\pi_z$. 
 Then  $\pi_z \circ \sigma_i \colon L_i^- \to \mW$ is a diffeomorphism into the face $\partial_h^- \mW$, with   image  denoted by $\fD_i \subset \partial_h^- \mW$. 
 Given this assumption, let $\ds \vartheta_i = (\pi_z \circ \sigma_i)^{-1} \colon \fD_i \to L_i^-$ denote the inverse map, so we have: 
\begin{equation}\label{eq-coordinatesVT}
\vartheta_i(r', \theta',-2) = (r(\vartheta_i(r',\theta')), \theta(\vartheta_i (r',\theta')), -2)  = (R_{i,r'}(\theta', -2)  , \Theta_{i,r'}(\theta' , -2), -2)  ~.
\end{equation}
 We formalize the   assumptions on   the insertion maps $\sigma_i$ that are intuitively 
implicit in  Figure~\ref{fig:twisted}. Set

 \begin{hyp} [{\it Strong Radius Inequality}]  \label{hyp-SRI}
  For   $i=1,2$, assume that: 
  \begin{enumerate}
\item \label{item-SRI-2} $\sigma_i \colon L_i^- \to \mW$ is transverse to the fibers of $\pi_z$; 
\item \label{item-SRI-1} $r' = r(\sigma_i(r,\theta,z))< r$, except for  $(2,\theta_i, z)$ and  then  $z(\sigma_i(2,\theta_i, z)) = (-1)^i$; 
  \item \label{item-SRI-3} $\Theta_{i,r'}(\theta') = \theta(\vartheta_i(r', \theta',-2))$ is an increasing function of $\theta'$ for each fixed $r'$; 
\item  \label{item-SRI-4} $R_{i,r'}(\theta') = r(\vartheta_i(r',  \theta',-2))$ has non-vanishing derivative for $r'=2$, except for the case of  $\theta_i'$ 
  defined by  $\vartheta_i(2,\theta_i',-2)= (2,\theta_i,-2)$; 
  \item For $r'$ sufficiently close to $2$, we require that the $\theta'$ derivative of $R_{i,r'}(\theta')$ vanish at a unique point denoted by $\theta'(i,r')$.
\end{enumerate}
Consequently, each  surface $\cL_i^-$   is transverse to the coordinate vector fields $\partial/\partial \theta$ and $\partial/\partial z$ on $\mW$.
  \end{hyp}

Recall from \eqref{eq-wilsonvector} that we have
$\cW =g(r, z)  \frac{\partial}{\partial  z} + f(r, z)  \frac{\partial}{\partial  \theta}$,
where    $g \colon \bR \to [0,1]$,  which  satisfies the ``vertical''
symmetry condition $g(r,z) = g(r,-z)$, $g(2,-1) = g(2,1) = 0$. Also, we have    that $g(r,z) = 1$ for $(r,z)$ near the boundary of $\bR$, and that $g(r,z) > 0$ otherwise.  
These assumptions are made more precise by specifying    that  
\begin{equation}\label{eq-generic1}
g(r,z) = 1 \quad \text{for} \quad  (r-2)^2 + (|z| -1)^2 \geq \e_0^2 
\end{equation}
where $0 < \e_0  \leq 1/4$ is sufficiently small so that the closed $\e_0$-neighborhood of each special point $p_i^{\pm}$  intersects the insertion regions $\cL_i^{\pm}$ in the interior of the face. We also require that $2 + \e_0 \leq r_e$ where $r_e$ is the exceptional radius defined in \eqref{eq-eradius}.
Finally, we require that $g(r,z)$ is monotone increasing as a function of the distance $\ds \sqrt{(r-2)^2 + (|z| -1)^2}$ from the special points $(2, -1)$ and $(2,1)$.

As $g(r,z) \geq 0$,  the first derivatives of $g$ must vanish at the
points $(2,\pm 1)$ and the $2 \times 2$ Hessian matrix of second
derivatives at these points must be positive \emph{indefinite}. 
The function $g$ is said to be   \emph{non-degenerate} if its Hessian matrix is   \emph{positive definite} at the points $(2,\pm 1)$, and the value of  $g(r,z)$ is a non-decreasing function of the distance from the points in the $\e_0$-ball  around each.

 \begin{hyp}\label{hyp-genericW} 
 Assume that $f(r, z)$ satisfies the conditions (W1) to (W6) in Section~\ref{sec-wilson}, that the condition \eqref{eq-generic1} holds and that   $g$  is  non-degenerate.   
  \end{hyp}

   Hypotheses~\ref{hyp-SRI} and \ref{hyp-genericW} are not required for the results in   previous sections, though some version of these hypotheses appear to be implicitly   assumed by   certain conclusions stated   in  \cite{Ghys1995,Matsumoto1995}.

\medskip

We next apply these assumptions to the study of the $\cK$-orbit of  the
notched Reeb cylinder $\cR'$.
First, note that   for each $1 \leq r < 3$,  if the intersection    $\cC(r) \cap \cL_i^-$ is non-empty, 
then it is a curve. Hypothesis~\ref{hyp-SRI}.\ref{item-SRI-3} implies that each such arc of intersection in 
$\cC(r)$ has the property that,  as the $z$-coordinate increases
along the curve, the $\theta$-coordinate   decreases for the curve in $\cL_1^-$, or increases for the curve in $\cL_2^-$.

Consider  the intersections of the faces $\cL_i^{\pm}$ with the Reeb cylinder $\cR$, 
\begin{equation}\label{eq-curvesprime}
\g' = \cR \cap \cL_1^- ~ , ~ \ovg' = \cR \cap \cL_1^+ \quad ; \quad \lambda' = \cR \cap \cL_2^- ~ , ~ \ovl' = \cR \cap \cL_2^+
\end{equation}
which  are arcs    transverse to the $\Psi_t$-flow on $\cR$, as illustrated in Figure~\ref{fig:notches}. 
  Label their preimages in $\partial_h^{\pm}\mW$ under the insertion maps $\sigma_i$ by 
\begin{equation}\label{eq-curves}
\g = \sigma_1^{-1}(\g') ~ , ~ \ovg = \sigma_1^{-1}(\ovg' ) \quad ; \quad
 \lambda = \sigma_2^{-1}(\lambda') ~ , ~ \ovl = \sigma_2^{-1}(\ovl') .
 \end{equation}
 
One endpoint  of the curve $\g$ is contained in the boundary $L_1^- \cap \{r=3\}$ and the other endpoint is the special point $\sigma_1^{-1}(p_1^-) \in  L_1^- \cap \{r=2\}$.  Similarly, one endpoint  of the curve $\lambda$ is contained in the boundary $L_2^- \cap \{r=3\}$ and the other endpoint is the special point $\sigma_2^{-1}(p_2^-) \in  L_2^- \cap \{r=2\}$.

\begin{figure}[!htbp]
\centering
\begin{subfigure}[c]{0.4\textwidth}{\includegraphics[width=50mm]{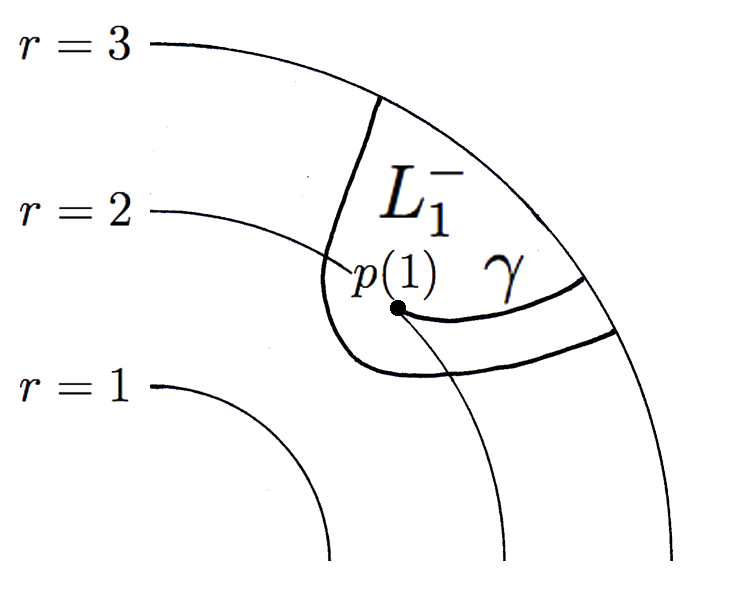}}\end{subfigure}
\begin{subfigure}[c]{0.4\textwidth}{\includegraphics[width=50mm]{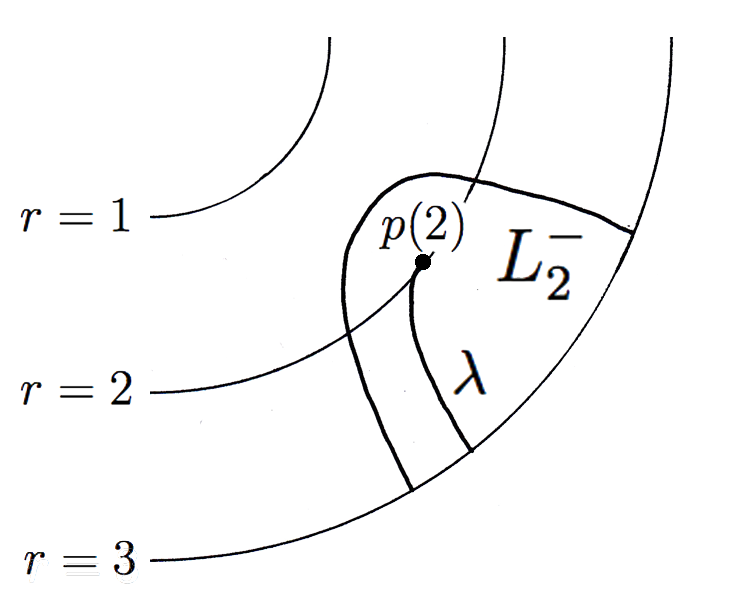}}\end{subfigure}
\caption{\label{fig:gamma'}  The curves $\g$ and $\lambda$ in $\partial_h^-\mW$}
\vspace{-6pt}
 \end{figure}

Hypothesis~\ref{hyp-SRI}   implies  that  both     $\g$ and $\lambda$ are  transverse to the cylinders $\{r=r_0\}$ for $2 < r_0 \leq 3$. As the $z$-coordinate of $x' \in \g'$   decreases towards $z =-1$, the radius coordinate $r$ of the corresponding point in $\gamma$ monotonically decreases to $r=2$,  while the angle coordinate $\theta$ monotonically increases towards $\theta_1$ as defined in (K7). 
  For $\lambda'$, as the $z$-coordinate of $x' \in \lambda'$ increases towards $z =1$, 
the radius of the corresponding point in   $\lambda$ monotonically decreases,   
while the angle coordinate $\theta$ monotonically increases towards $\theta_2$ as defined in (K7). 
Thus,   the graphs of the curves $\g$ and $\lambda$ in  $\partial_h^-\mW$ appear  as in Figure~\ref{fig:gamma'}.

\begin{lemma}\label{lem-facingcurves}
The curve $\ovg \subset L_1^+$ is the facing curve to $\g \subset L_1^-$ and $\ovl \subset L_2^+$ faces $\lambda \subset L_2^-$.
Moreover, for each point $x \in \g$ with $r(x) > 2$, the facing point $\ox$ satisfies $x \prec_{\cK} \ox$. 
Similarly, for each point $y \in \lambda$ with $r(y) > 2$, the facing point $\oy$ satisfies $y \prec_{\cK} \oy$. 
\end{lemma}
\proof
The fact that $\ovg$ is facing to $\g$ and that $\ovl$ is facing to $\lambda$, follows from the construction of the insertions $\sigma_i$. The assertions
$x \prec_{\cK} \ox$ for $r(x) > 2$ and $y \prec_{\cK} \oy$ for $r(y) > 2$, then follow from  Proposition~\ref{prop-propexist}. 
\endproof

Let $\g$ be   parametrized by $w \colon [0,1]\to \mW$, with $r(w(0)) = 3$ and $r(w(1))=2$. 
For $0 \leq \e \leq 1$, define $\g^{\e}$ by $w(\e \cdot s)$ for $0
\leq s \leq 1$, as in \eqref{eq-sparem}.
Then Lemma~\ref{lem-facingcurves} implies that for $0 < \e < 1$ and $0 \leq s \leq 1$, we have  $\g^{\e}(s) \prec_{\cK} \ovg^{\e}(s)$. 
However, the $\cW$-orbits of the points  $w_1 = w(1)$ and $\ow_1 \equiv
w_1$ are both trapped, so cannot satisfy $w_1 \prec_{\cW} \ow_1$.
Thus, the curves $\g$ and $\ovg$ are the horizontal boundary curves in $\partial_h^{\pm} \mW$ of
the infinite propeller  $P_{\g}$  defined by  the $\Phi_t$-flow of $\g$. Similar conclusions hold for  $\lambda$, $\ovl$ and the propeller $P_{\lambda}$.

Note that in  the Definition~\ref{def-infpropeller} of the infinite propeller $P_{\g}$,  we formed   the union  of the propellers $P_{\g^\e}$ defined by the curves $\g^{\e}$ for $0 < \e < 1$,  with  the   two trapped orbits $\cZ_{\g}$.  The $\cW$-orbit given by the ``long edge'' of each $P_{\g^\e}$ is a finite $\cW$-arc, but as $\e \to 0$  their limit converges to the union $\cZ_{\g}$ of infinite orbits.  
This  behavior   is reminiscent of the ``Moving Leaf Lemma'' in \cite{EpsteinVogt1978,Sullivan1976}, which is the key to understanding the orbit behavior  in the counter-examples to the \emph{Periodic Orbit Conjecture} \cite{Epstein1972}.
 
The $\cK$-orbits of the curves  $\g$ and $\lambda$ projected to $\mK$ have a  complicated, hierarchical   structure, which we next describe  using their $\Phi_t$-flows as the starting model. 
Define the \emph{notched  propellers}  of $\g$ and $\lambda$   by $P_{\g}'=P_{\g}\cap \mW'$ and    $P_{\lambda}'=P_{\lambda}\cap \mW'$, respectively. Note that the vertical ``transverse'' boundary curves for the notches are not included in $P_{\g}'$ and $P_{\lambda}'$.

\begin{figure}[!htbp]
\centering
{\includegraphics[width=120mm]{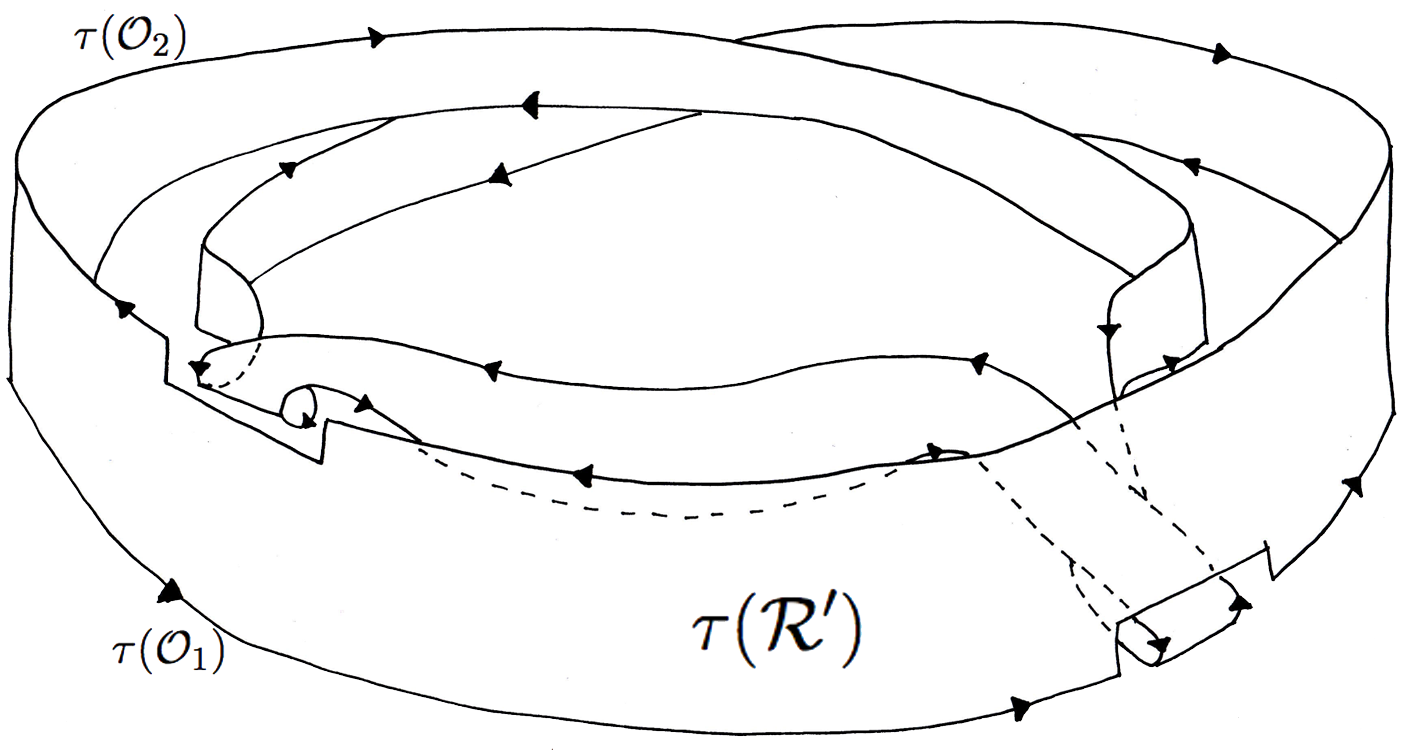}}
\caption{\label{fig:notched8}  Embedding of $\cR'$ in $\mK$}
\vspace{-6pt}
\end{figure}

For $x \in \tau(\cR')$, the level function $n_x(t)$ increases from $0$ to $1$ when the orbit of $x$ intersects either curve $\tau(\g)$ or $\tau(\lambda)$, then drops back to $0$ when it exits through the curves $\tau(\ovg)$ or $\tau(\ovl)$. 
 Thus,  we have
 \begin{equation}\label{eq-propeller0}
\fM_0^0 ~ = ~ \tau(\cR' ~ \cup ~   \ovg  ~ \cup ~ \ovl)
\end{equation}
which is a cylinder, minus two rectangles, embedded in $\mR^3$ as a
folded eight having two parts that are tangent to the boundary of the notches, as in    Figure~\ref{fig:notched8}.

 We comment on the details of Figure~\ref{fig:notched8}.
Note that the embedding $\tau$ of $\cR'$ in $\mK$ differs from its inclusion in $\mW$ as illustrated in Figure~\ref{fig:notches}, as the Reeb cylinder in embedded as a  folded figure-eight.   The insertions $\sigma_i$ fold  parts of the inner loop in the folded eight and make them tangent to the erased rectangles. 
 This process brings two vertical lines in the cylinder to  arcs contained in the erased boundary of the rectangles,
  covering the middle half of it.   
  The remaining portions on each interval correspond to the  cylindrical complement  $\cC - \cR$, 
  and will be considered later in Section~\ref{sec-doublepropellers}.

The set $\fM_0^1$ again has a simple intuitive description. It is obtained by attaching the points in $\fM_0$ of level $1$   to $\fM_0^0$. The curves 
$\tau(\g)$ or $\tau(\lambda)$ were observed to have level $1$, and thus 
the points in their $\cW$-orbits which lie in $\mW'$ are also at level $1$.  
 In particular,  the   images  $\tau(P_{\g}')$ and $\tau(P_{\lambda}')$ are at level $1$, so are contained in $\fM_0^1$.

Each of the infinite propellers $\tau(P_{\g}')$ and
$\tau(P_{\lambda}')$ is obtained by removing  an infinite sequence of notches from the propellers $P_{\g}$ and $P_{\lambda}$ and then applying $\tau$. 
There are two types of notches: 
\begin{enumerate}
\item an infinite collection of notches, each   having a corner containing   a point of level $2$ in the intersection of the $\cK$-orbit of a special point $p_i^{\pm}$ with the surfaces $E_j$ for $j=1,2$;
\item a finite number, possibly zero,  of 
interior  notches   corresponding to the intersections of $P_{\g}$ and $P_{\lambda}$ with $\cD_j$ for $j=1,2$, which do not contain points in the $\Phi_t$-orbit of the special points. 
\end{enumerate}

\begin{remark}\label{rmk-bubbles}
The intersections of the second type give rise to what we call
``bubbles'' in the propeller and are discussed in detail later in this
section and     in Sections~\ref{sec-bubbles} and \ref{sec-geometry}.

The precise shape of propellers depends on the choices made in the
construction of the Kuperberg plug. In this section and
Section~\ref{sec-bubbles} we consider a first type of internal notches
and generated bubbles that can arise.

In Section~\ref{sec-generic} we introduce new hypothesis on the
construction that allow to give a more accurate description of the
shape of propellers. This analysis is thus postponed until
Section~\ref{sec-geometry}. One consequence of it, will be another
possible situation where internal notches and bubbles arise, as
explained at the end of that section.
\end{remark}

  The construction of the set $\fM_0^{2}$ uses ``gluing'' of propellers at higher level to the boundaries of these deletions  from the infinite propellers at level $1$. In order to give a more precise description of $\fM_0^1$, in preparation for describing the nuances of the construction of $\fM_0^n$ for $n \geq 2$, 
we develop an indexing system to label these intersections.

  Recall from \eqref{eq-special} that  for $i=1,2$,  $p_i^{\pm} = \tau(\cO_i \cap \cL_i^{\pm})$ is called a special point.  The corresponding points    $\omega_i \in   \bR_0$ for $i=1,2$, as in \eqref{eq-omegas}, are used to define the level function along the special orbits. That is, 
the level function on $\fM_0$ along the $\cK$-orbit of $\omega_i$  is defined by $n_0(x) = n_{\omega_i}(t_x)$ for  $x = \Phi_{t_x}(\omega_i)$. 
 Note that  to obtain the full $\cK$-orbits of $\omega_i$,    it is     necessary to consider both their forward and backward orbits.  
  
 The first transition point  for the forward $\cK$-orbit of $\omega_1$ is the  entry point  $p_1^- \in E_1$ for which  $n_0(p_1^-) = 1$.   
   Set $p(1) = \tau^{-1}(p_1^-)  \in L_1^-$  and note that    $r(p(1)) = 2$. 
The  forward $\cW$-orbit of  $p(1)$  is trapped in the region $\cC \cap \{z < -1\}$ and thus  intercepts  $\cL_1^- \cap \cC$ in an infinite  sequence of points with   increasing $z$-coordinates between $-2$ and $-1$. 
Label these points   $p'(1; 1 , \ell)$ for   $\ell \geq 0$. 

Note that   $r(p'(1; 1 , \ell)) =2$ for all $\ell \geq 0$ and we have 
\begin{equation}\label{eq-interlaced1}
-2 < z(p'(1; 1 , 0)) < \cdots < z(p'(1; 1 , \ell)) <  z(p'(1; 1 , \ell +1)) < \cdots < -1 
\end{equation}
where $z(p'(1; 1 , \ell)) \to -1$ as $\ell \to \infty$.

   Set  $p(1; 1 , \ell)  = \sigma_1^{-1}(p'(1; 1 , \ell))   \in L_1^-$ for $\ell \geq 0$.
Then  $r(p(1; 1 , \ell))  > 2$  by the Radius Inequality (K8) and the
sequence $p(1; 1 , \ell)$ accumulates  in $L_1^-$  on  $p(1)$ as $\ell
\to \infty$. Hence  $r(p(1; 1 , \ell)) \to 2$ as $\ell \to \infty$. 
Note that  $n_0(\tau(p(1; 1 , \ell))) =2$ for all $\ell \geq 0$.

 Similarly,   the first transition point  for the forward $\cK$-orbit of $\omega_2$ is the  entry point  $p_2^- \in E_2$ with $n_0(p_2^-) = 1$.  
 Set $p(2) = \tau^{-1}(p_2^-) \in   L_2^-$ and note that    $r(p(2)) = 2$.

  The forward $\cW$-orbit of $p(2)$  is also trapped in the region $\cC \cap \{z < -1\}$ and thus  intercepts    $\cL_1^-\cap \cC$ in an infinite  sequence of points with increasing  $z$-coordinates also between $-2$ and $-1$. Label these points  $p'(2; 1 , \ell)$ for  $\ell \geq 0$.  
  
 Note that   $r(p'(2; 1 , \ell)) =2$ for all $\ell \geq 0$, and we have 
 \begin{equation}\label{eq-interlaced2}
-2 < z(p'(2; 1 , 0)) < \cdots < z(p'(2; 1 , \ell)) <  z(p'(2; 1 , \ell +1)) < \cdots < -1 
\end{equation}
where $z(p'(2; 1 , \ell)) \to -1$ as $\ell \to \infty$.  

      Set  $p(2; 1 , \ell)    = \sigma_1^{-1}(p'(2; 1 , \ell)) \in L_1^-$ for $\ell \geq 0$.  
Again,   $r(p(2; 1 , \ell))  > 2$  by the Radius Inequality (K8) 
and the sequence $p(2; 1 , \ell)$ also accumulates  in $L_1^-$  on  $p(1)$ as $\ell \to \infty$. 
Thus, $r(p(2; 1 , \ell)) \to 2$ as $\ell \to \infty$ and   $n_0(\tau(p(2; 1 , \ell))) =2$ for all $\ell \geq 0$.

\begin{lemma}\label{lem-interlaced12}
The sequences $\{p'(1; 1 , \ell) \mid \ell \geq 0\}$ and $\{p'(2; 1 ,
\ell') \mid \ell' \geq 0\}$ are \emph{interlaced} on the line segment
$\cC \cap \cL_1^- \cap \{z < -1\} \subset \mW$. 
If  $z(p'(1;1,0))<z(p'(2;1,0))$ then we have
\begin{equation}\label{eq-interlaced12}
 -2 < \cdots < z(p'(1; 1 , \ell))  < z(p'(2; 1 , \ell)) <  z(p'(1; 1 , \ell +1))  < \cdots < -1 .
\end{equation}
The analogous conclusion holds when $z(p'(2;1,0))<z(p'(1;1,0))$.
\end{lemma}
\proof
Observe that $L_1^-$ follows $L_2^-$ in the direction of the
$\theta$-coordinate in $\mW$. Since  $p'(1; 1 , 0)$ and $p'(2; 1 ,
0)$ are the first intersections of the $\cW$-orbit of $p(1)$ and
$p(2)$ with $\cL_1^-$, we cannot predict which is below. Assuming
that $z(p'(1;1,0))<z(p'(2;1,0))$,  the inequality 
$z(p'(2; 1 , \ell)) <  z(p'(1; 1 , \ell +1))$, for $\ell \geq 1$, follows from the fact that the $\Psi_t$-flow preserves the cylinder $\cC$ and so preserves the height relationship. The case when $z(p'(2;1,0))<z(p'(1;1,0))$ follows in the same way. 
\endproof

  Next, consider the backward orbits of the points  $\omega_i$, for $i=1,2$, which intersect $S_i$ in the points   $p_i^+$ with     $\tau^{-1}(p_i^+) \in L_i^+$ and $r(\tau^{-1}(p_i^+)) = 2$.   
Thus,   the backward Wilson orbits of $\tau^{-1}(p_i^+)$ for $i=1,2$ are trapped in the region $\cC \cap \{z > 1\}$ and so   intercept  $\cL_2^-$ in infinite  sequences of points with $r =2$ and $z$-coordinates between $1$ and $2$.   

For the backward $\cW$-orbit of  $\tau^{-1}(p_1^+)$, label these points   $p'(1; 2 , \ell) \in \cL_2^-$,  for $\ell \geq 0$, with       $r(p'(1; 2 , \ell)) =2$ and $1 < z(p'(1; 2 , \ell)) < 2$, where $z(p'(1; 2 , \ell)) \to 1$ as $\ell \to \infty$.
  
       Set  $p(1; 2 , \ell)   = \sigma_2^{-1}(p'(1; 2 , \ell)) \in L_2^-$.    
We then  have $n_0(\tau(p(1; 2 , \ell)) ) = 2$   for $\ell \geq 0$ by formula \eqref{def-level-}.     
The   Radius Inequality (K8) implies that   $r(p(1; 2 , \ell))   > 2$  and note that the sequence $p(1; 2 , \ell)$   accumulates   on  $p(2)$ as $\ell \to \infty$. Thus, $r(p(1; 2 , \ell)) \to 2$ as $\ell \to \infty$.

Similarly, the   backward  $\cW$-orbit of $\tau^{-1}(p_2^+)$ intercepts  $\cL_2^- \subset \mW$ in a sequence of points   with $r=2$ and  $z$-coordinate between $1$ and $2$. 
Label these points   $p'(2; 2 , \ell) \in \cL_2^-$,   for $\ell \geq 0$.
    Then $r(p'(2; 2 , \ell)) =2$ and $1 < z(p'(2; 2 , \ell)) < 2$, where $z(p'(2; 2 , \ell)) \to 1$ as $\ell \to \infty$.

Set  $p(2; 2 , \ell)   = \sigma_2^{-1}(p'(2; 2 , \ell)) \in L_2^-$.    
     Again, we have $n_0(\tau(p(2; 2 , \ell))) = 2$ for $\ell \geq 0$ by formula \eqref{def-level-}.
  Then  $r(p(2; 2 , \ell))  > 2$  by the Radius Inequality (K8) and the sequence $p(2; 2 , \ell)$   accumulates   on  $p(2)$ as $\ell \to \infty$. Thus, $r(p(2; 2 , \ell)) \to 2$ as $\ell \to \infty$.

The analog of Lemma~\ref{lem-interlaced12} follows by the same arguments.
\begin{lemma}\label{lem-interlaced22}
The sequences $\{p'(1; 2 , \ell) \mid \ell \geq 0\}$ and $\{p'(2; 2 , \ell') \mid \ell' \geq 0\}$ are \emph{interlaced} on the line segment $\cC \cap \cL_2^- \cap \{z > 1\} \subset \mW$. 
If  $z(p'(1;2,0))<z(p'(2;2,0))$ then
\begin{equation}\label{eq-interlaced22}
 1 < \cdots <     z(p'(2; 2 , \ell +1))  < z(p'(1; 2 , \ell))  < z(p'(2; 2 , \ell))   <  \cdots  < 2 .
\end{equation}
The analogous conclusion holds when $z(p'(2;2,0))<z(p'(1;2,0))$.
\end{lemma}

\begin{remark}\label{rmk-labels}
The points $\tau(p(1; \cdot, \cdot))$ belong to the $\cK$-orbit of $\omega_1$,
     the points $\tau(p(2; \cdot, \cdot))$ belong to  the $\cK$-orbit of $\omega_2$, 
 while the points $\tau(p(\cdot; 1, \cdot)) \in E_1$ and  the points $\tau(p(\cdot; 2, \cdot)) \in E_2$.
\end{remark}

Next, we consider in detail the intersections of the propellers
$P_{\g}$ and $P_{\lambda}$ with the surfaces $\cL_i^-$ for
$i=1,2$. Recall that $\cL_i^-$ was chosen  in
Section~\ref{sec-kuperberg} to be transverse to the vector field $\cW$ and the insertion $\cD_i$ is obtained from the $\Psi_t$-flow of its points. It follows
that the surfaces $\cL_i^-$ intersect the propellers  $P_{\g}$ and
$P_{\lambda}$ transversally, thus each such intersection must be a
union of closed line segments whose boundaries are contained in the
boundaries of either $\cL_i^-$ or in the boundary of the propellers. The intersections of $\cD_i$ with a propeller are then obtained as the $\Psi_t$-flow of each such line segment in $\cL_i^-$, until it reaches the facing line segment in $\cL_i^+$.  
We consider the possible cases for the line segments   $\cL_i^- \cap P_{\g}$. The analysis for the propeller $P_{\lambda}$ will be analogous.

Recall that the intersection of $P_\g$ with $\bRt$ is an infinite
collection of arcs whose endpoints are on the vertical line $\{r=2\}$,
the lower ones having $z$-coordinate less than $-1$ and the upper ones
having $z$-coordinate bigger than 1,
as illustrated in Figure~\ref{fig:arcspropeller}. The forward $\Psi_t$-flow of
these arcs gives the intersection $\cL_i^-\cap P_\g$. For $i=1$, an
arc in $\cL_1^-\cap P_\g$ either has both endpoints on the boundary
$\partial \cL_1^-$, either has its upper endpoint on $\partial \cL_1^-$
and its lower endpoint is a point $p'(1;1,\ell)$ for $\ell\geq 0$. For $i=2$, an
arc in $\cL_2^-\cap P_\g$ either has both endpoints on the boundary
$\partial \cL_2^-$, either has its lower endpoint on $\partial \cL_2^-$
and its upper endpoint is a point $p'(1;2,\ell)$ for $\ell\geq 0$. 

\begin{defn}\label{def-notches}
If a segment of $\cL_i^- \cap P_{\g}$ has both endpoints on $\partial
\cL_i^-$, the corresponding intersection of $P_\g$ with $\cD_i$ defines a ``rectangular notch'' in the interior of the propeller,  called an \emph{interior notch} for $P_{\g}$.
If  a segment of $\cL_i^- \cap P_{\g}$ has an endpoint in the interior
of $\cL_i^-$, the corresponding intersection of $P_\g$ with $\cD_i$ defines a \emph{boundary notch} for $P_{\g}$, or sometimes simply as a ``notch''.  
\end{defn}

\begin{figure}[!htbp]
\centering
{\includegraphics[width=120mm]{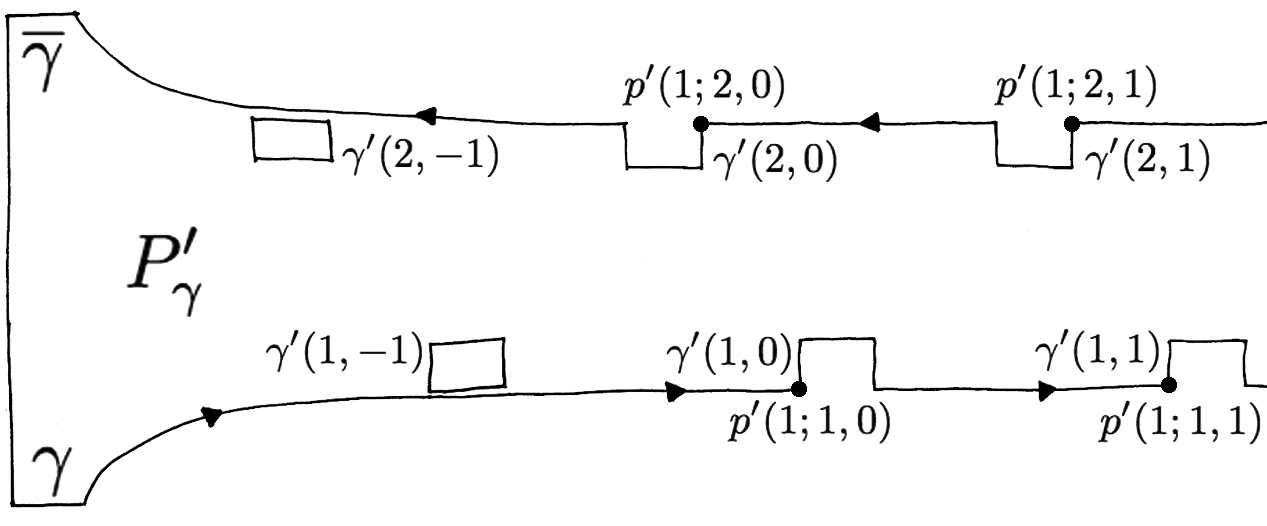}}
\caption{\label{fig:intnotches}  Internal and (boundary) notches of $P'_\g$}
\vspace{-6pt}
\end{figure}

The surface $P_{\g}'$ is   obtained from $P_{\g}$ by deleting the interior and boundary notches for $i=1$ and $i=2$. This is illustrated in Figure~\ref{fig:intnotches}.

For  the propeller $P_{\lambda}$ formed by the $\Psi_t$-flow of the curve $\lambda$, analogous considerations and notations apply to the connected components of $P_{\lambda} \cap \cL_i^-$ for $i=1, 2$.

  Observe that $P'_\g$
has a boundary notch for each
intersection of the $\cW$-orbits  of $p(1)=\tau^{-1}(p_1^-)$ and
$\tau^{-1}(p_1^+)$ with $\cL_1^-\cup \cL_2^-$. Similarly,  $P'_\lambda$
has a boundary notch for each
intersection of the $\cW$-orbits  of $p(2)=\tau^{-1}(p_2^-)$ and
$\tau^{-1}(p_2^+)$ with $\cL_1^-\cup \cL_2^-$. There are an infinite number of such  boundary notches. 

In contrast, there are at most finitely many internal notches, since these
correspond to the intersections of $P'_\g$ with $\cL_1^-\cup \cL_2^-$ for
which the $\cW$-orbits of $p(1)=\tau^{-1}(p_1^-)$ and
$\tau^{-1}(p_1^+)$ are away from the insertion region, or to the intersections of $P'_\lambda$ with $\cL_1^-\cup \cL_2^-$ for
which the $\cW$-orbits of $p(2)=\tau^{-1}(p_2^-)$ and
$\tau^{-1}(p_2^+)$ are away from the insertion region.  Let $-b \geq 0$ denote the number of internal notches, with $b=0$ if there are no internal notches, so that the $\cW$-orbit of $p(1)$  makes at least $|b|$ revolutions in the $\theta$-coordinate, before it intersects the surface $\cL_1^-$.

We now describe the set  $\fM_0^1$. 
The images  $\tau(P_{\g}')$ and $\tau(P_{\lambda}')$ consist of points of level $1$,  and the exit  boundaries of the notches in these propellers   are also at level $1$, so all belong to $\fM_0^1$. 

We index the notches at level $2$ as follows.
The segment  $\g$ is said to be \emph{based} at its inner endpoint  $p(1) = \tau^{-1}(p_1^-) \in L_1^-$ and 
  $\lambda$     to be \emph{based} at its inner endpoint  $p(2) = \tau^{-1}(p_2^-) \in   L_2^-$.

  Next, introduce labels for the four families   of boundary notches which arise. It may help to consider the illustration Figure~\ref{fig:intnotches}  to keep track of the following definitions.

 Let  $\g'(1 , \ell) \subset \cL_1^-$ for $b   \leq \ell < \infty$ denote the curves    corresponding to  the intersection of $P_{\g}$ with $\cL_1^-$, where    $\gamma'(1,\ell)$ for $b \leq \ell <0$ denotes an interior curve of the intersection with   $\cL_1^-$, assuming that such exists.   
Let $\g'(1,\ell)$ for $\ell\geq 0$ denote the boundary curve  with lower  endpoint $p'(1; 1 , \ell)$.    
Set  $\g(1 , \ell) = \sigma_1^{-1}(\g'(1 , \ell)) \subset L_1^-$.  We
say that the curve  $\g(1 , \ell)$ for $\ell\geq 0$ is  based at the point  $p(1; 1 , \ell)$. For all $x \in \tau(\g(1 , \ell))$, we have $n_0(x) = 2$ and $r(x) > 2$.

 Let  $\lambda'(1 , \ell) \subset \cL_1^-$ for $b \leq \ell < \infty$  denote the curves    corresponding to the intersection of $P_{\lambda}$ with  $\cL_1^-$, where   
    $\lambda'(1,\ell)$ for $b\leq \ell <0$ denotes an interior curve of the intersection with
$\cL_1^-$, assuming that such exists.   Let $\lambda'(1,\ell)$ for $\ell\geq 0$ denote the boundary curve  with lower endpoint $p'(2; 1 , \ell)$.  
 Set  $\lambda(1 , \ell) = \sigma_1^{-1}(\lambda'(1 , \ell)) \subset L_1^-$.  
 We say that the curve  $\lambda(1 , \ell)$ for $\ell\geq 0$ is  based at the point
 $p(2; 1 , \ell)$. For all $x \in \tau(\lambda(1 , \ell))$, we have $n_0(x) = 2$ and $r(x) > 2$.

 Let  $\g'(2 , \ell) \subset \cL_2^-$ for $b \leq \ell < \infty$  denote the curves    corresponding   to the intersection of  $P_{\g}$ with $\cL_2^-$, where 
  $\gamma'(2,\ell)$ for $b\leq \ell <0$ denotes an  interior curve of the intersection with   $\cL_2^-$, assuming that such exists.   
  Let   $\g'(2,\ell)$ for $\ell\geq 0$  denote the boundary curve  with upper endpoint $p'(1; 2 , \ell)$.   
 Set  $\g(2 , \ell) = \sigma_2^{-1}(\g'(2 , \ell)) \subset L_2^-$.  
We say that the curve  $\g(2 , \ell)$ for $\ell\geq 0$ is  based at the point  $p(1; 2
, \ell)$. For all $x \in \tau(\g(2 , \ell))$, we have $n_0(x) = 2$ and $r(x) > 2$.

  Let  $\lambda'(2 , \ell) \subset \cL_2^-$  for $b \leq \ell < \infty$  denote the curves    corresponding   to the intersection of  $P_{\lambda}$ with
$\cL_2^-$, where   $\lambda'(2,\ell)$ for $b\leq \ell <0$ denotes an interior curve of the intersection with   $\cL_2^-$, assuming that such exists. Let 
    $\lambda'(2,\ell)$ for $\ell\geq 0$ denote the boundary curve  with upper endpoint  $p'(2; 2 , \ell)$.    
    Set  $\lambda(2,\ell) = \sigma_2^{-1}(\lambda'(2 , \ell)) \subset L_2^-$.  
We say that the curve  $\lambda(2 , \ell)$ for $\ell\geq 0$ is  based at the point
$p(2; 2 , \ell)$. For all $x \in \tau(\lambda(2 , \ell))$, we have $n_0(x) = 2$ and $r(x) > 2$.

\begin{figure}[!htbp]
\centering
\begin{subfigure}[c]{0.45\textwidth}{\includegraphics[width=70mm]{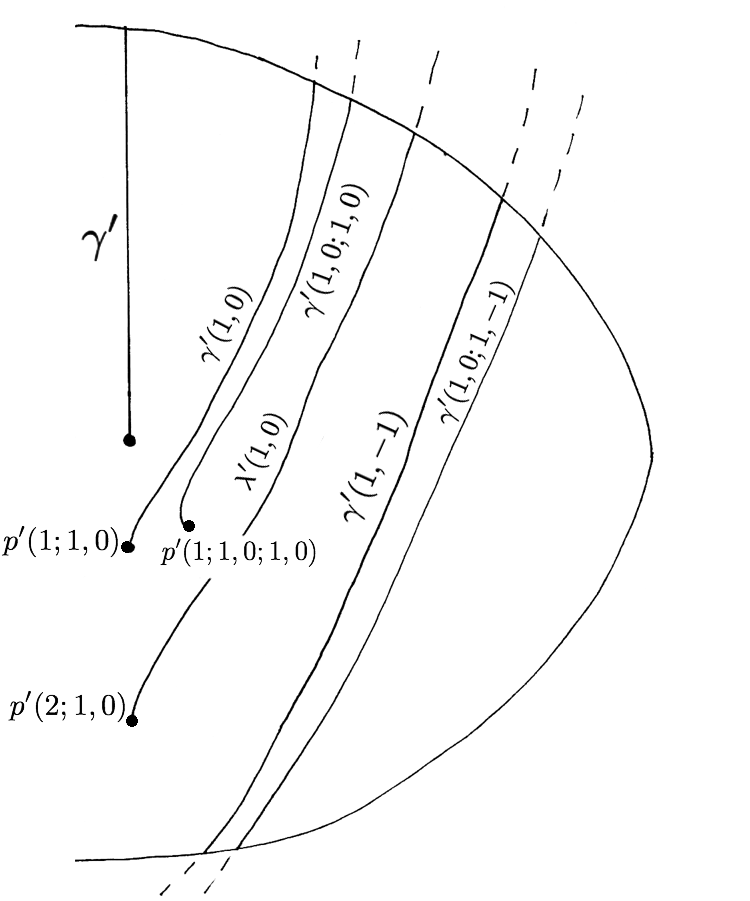}}\end{subfigure}
\begin{subfigure}[c]{0.45\textwidth}{\includegraphics[width=70mm]{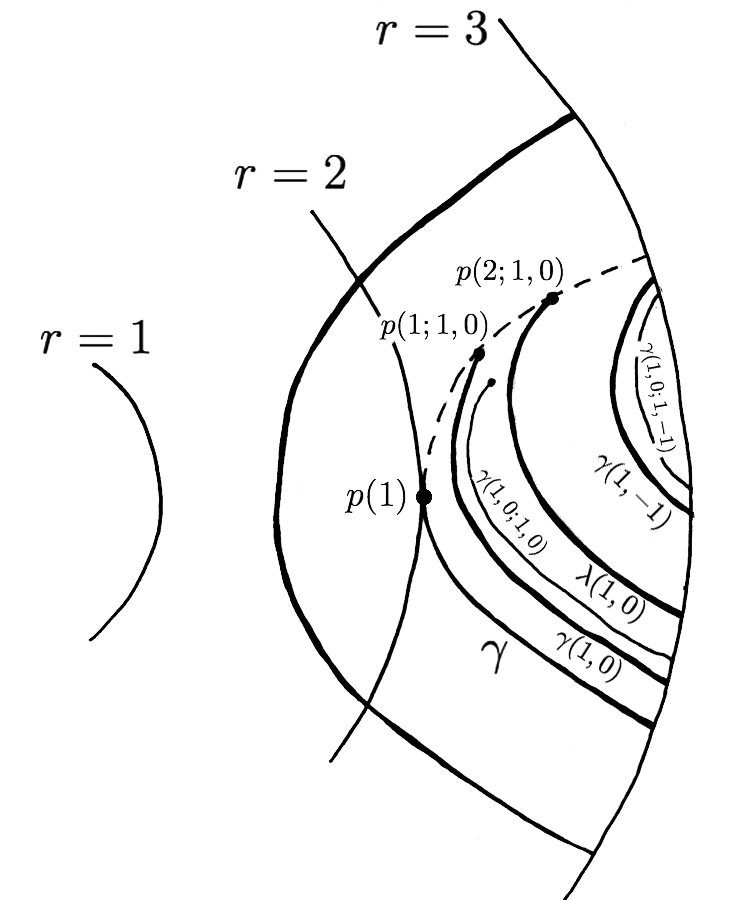}}\end{subfigure}
\caption{\label{fig:intcurvesL1} Curves of levels 1 and 2, in  $\cL_1^-$ and in $L_1^-$}
\vspace{-6pt}
 \end{figure}

\begin{figure}[!htbp]
\centering
\begin{subfigure}[c]{0.45\textwidth}{\includegraphics[width=70mm]{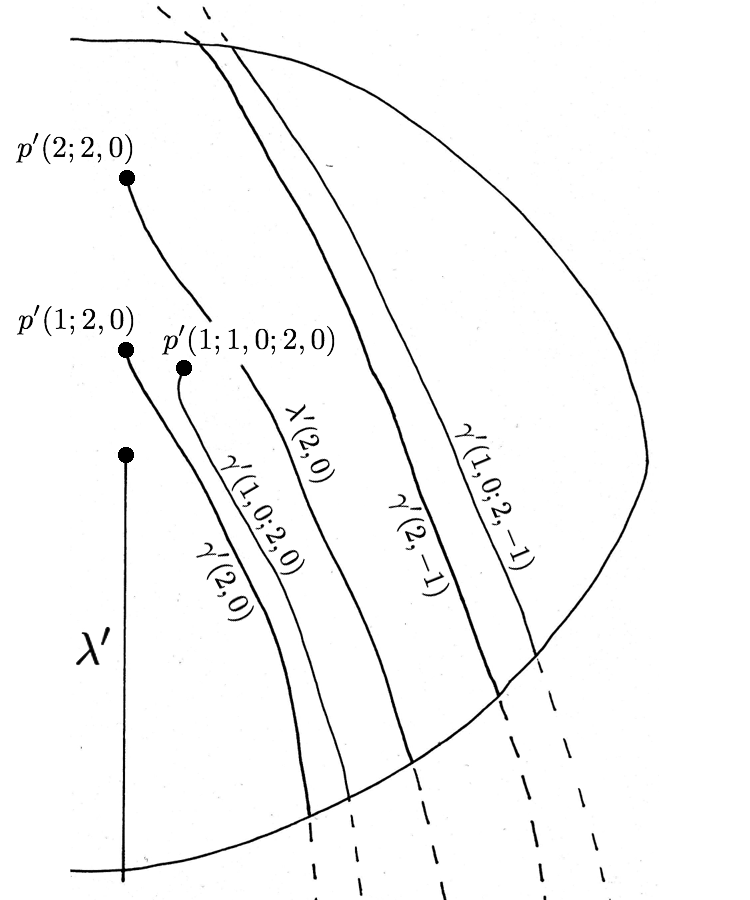}}\end{subfigure}
\begin{subfigure}[c]{0.45\textwidth}{\includegraphics[width=70mm]{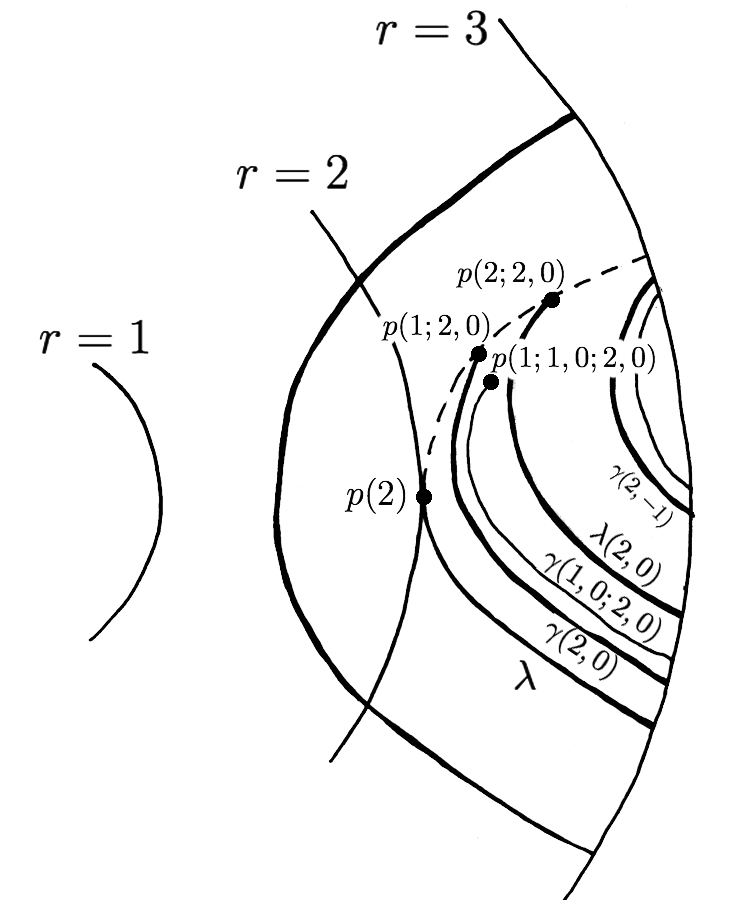}}\end{subfigure}
\caption{\label{fig:intcurvesL2} Curves of levels 1 and 2, in  $\cL_2^-$ and in $L_2^-$}
\vspace{-6pt}
 \end{figure}

Lemmas~\ref{lem-interlaced12} and \ref{lem-interlaced22} show that 
  the families of points  $\{p'(1; i , \ell) \mid \ell \geq 0\}$ and  $\{p'(2; i , \ell') \mid \ell' \geq  0\}$ for $i =1,2$ are
  interlaced, so the same holds for the $\gamma$ and $\lambda$ curves in $\partial_h^- \mW$ and also for  the curves with $\ell<0$. 
  That is, each $\lambda$-curve is between two $\g$-curves, and vice-versa, as  
  illustrated in Figures~\ref{fig:intcurvesL1} and \ref{fig:intcurvesL2}.


We comment on the details of Figures~\ref{fig:intcurvesL1} and
\ref{fig:intcurvesL2}. The curves in the figures come by pairs, in
each pair one curve has level 1 and the other has level 2. We give
details for the level 1 curves, that in the left side figures can be
distinguished by the fact that one of their endpoints is in the
vertical line $\{r=2\}$. In the left side of
Figure~\ref{fig:intcurvesL1}, the closest curve to $\g'$ is
$\g'(1,0)$, the boundary of the first boundary notch of $P_\g'$. The
following level 1 curve is $\lambda'(1,0)$ and the next one is
$\g'(1,-1)$. Observe that $b$ might be $-1$ for $P_\g'$ and zero for
$P_\lambda'$, but the number of internal notches is each propeller differs at
most by 1. On the right side of Figure~\ref{fig:intcurvesL1}, the
image under $\sigma_1^{-1}$ of the curves is illustrated. Observe that
the curves $\g(1,0)$ and $\lambda(1,0)$ have an endpoint in the doted
line, that is the image under $\sigma_1^{-1}$ of the line
$\cC\cap\{z\leq -1\}\cap\cL_1^-$. Analogous considerations apply to
Figure~\ref{fig:intcurvesL2}. We give the description of the level 2
curves in this section.

Finally, to complete the description of the level $1$ set, note that 
  for $i = 1,2$, each curve  $\tau(\g(i , \ell)) \subset E_i$ defines a facing curve  $\tau(\ovg(i , \ell)) \subset S_i$.  These facing curves in $S_i$ are at level $1$ again, so belong to $\fM_0^1$.  
Thus we have
 \begin{equation}\label{eq-level1approx}
\fM_0^1 ~ = ~ \tau(\cR') ~ \cup ~  \tau(P_{\g}') ~ \cup ~
\tau(P_{\lambda}') ~ \cup \left\{ \bigcup_{\ell \geq b} ~ \bigcup_{i =1,2}~ 
\left[ ~ \tau(\ovg(i , \ell))  ~  \cup \tau(\ovl(i , \ell) ) ~ \right] \right\} 
\end{equation}
where the propellers and boundary curves are identified via the
insertion maps $\sigma_i$ as appropriate.   Observe that in the above, we are
  assuming that $b \leq 0$ has a constant value for each family of
  curves. By the interlaced property, the number of internal notches
  in $P_\g'$ and $P_\lambda'$ differs at most by one.

 For $n \geq 1$, the construction of  $\fM_0^{n+1}$ from $\fM_0^n$   follows a procedure similar to the above pattern. 
 
 The basic scheme is that for each image under $\tau$ of a propeller
 at level $n$, the entry curve along each of its notches is a curve at level $n+1$.  The $\Psi_t$-flow of this curve yields a family of propellers   at level $n+1$.  The   resulting propellers at level $n+1$ are contained in the region $ r> 2$, so they are finite.   
 We consider the case of $\fM_0^2$ in some detail, to illustrate the general construction  of $\fM_0^n$ for $n \geq 2$.
 
  Consider the case $b\leq \ell<0$ and let $\g(1 , \ell) =
  \sigma_1^{-1}(\g'(1 , \ell)) \subset L_1^-$. All the points on the
  curve $\g(1,\ell)$ have radius greater than $2$. Thus by
  Proposition~\ref{prop-propexist} the $\cK$-orbit of any point $x\in
  \tau(\g(1,\ell))$ contains the point $\overline{x}\in S_1$ such that
  $x\equiv \overline{x}$. Let $t_x>0$ be such that
  $\Phi_{t_x}(x)=\overline{x}$, then
\begin{equation}\label{eq-Ssurface}
S_{\g(1,\ell)}=\{\Phi_{t}(x) \, |\, x\in \tau(\g(1,\ell)), \, 0\leq
t \leq t_x\},
\end{equation}
is a compact surface embedded in $\mK$. 
Similarly, we obtain  
compact surfaces $S_{\lambda(1,\ell)}$, $S_{\g(2,\ell)}$ and
$S_{\lambda(2,\ell)}$ embedded in $\mK$, for $b\leq \ell <0$. 

The surface $S_{\g(1,\ell)}$ has points of level 2 and might have
points at higher levels, but the level is uniformly bounded above by 
Proposition~\ref{prop-bubbles1}. It contributes to $\fM_0^2$ with a
compact set, a bubble, that we describe in detail in
Section~\ref{sec-bubbles}. 
 To simplify the exposition, we assume for the
rest of this section that $b=0$ and thus that the propellers (at any
level) do not have   internal notches. We discuss the modifications needed to accommodate these internal notches later.

Consider the   curve $\g(1 , \ell) = \sigma_1^{-1}(\g'(1 , \ell)) \subset L_1^-$  where $\ell \geq 0$. It  
goes from the outer boundary of $L_1^- \subset \partial_h^- \mW$ to the point $p(1; 1 , \ell)$ with $r(p(1; 1 , \ell)) >  2$,  and all points on the curve $\g(1 , \ell)$  have radius greater than $2$. 
Thus, the $\Psi_t$-flow of the entire curve $\g(1 , \ell)$ traverses
$\mW$ in finite time, defining a finite propeller $\ds P_{\g(1 ,  \ell)}$ as in   Section~\ref{sec-intropropellers}. 
  By \eqref{eq-Delta}, the propeller $\ds P_{\g(1 , \ell)}$ 
intersects $\cL_1^-$ at most $\max_{x\in \g(1,\ell)}\Delta(r(x))$ times, for each curve in the collection  $\ds \left\{\g(1 , \ell) \right\}_{\ell=0}^\infty$.
In the same manner, we also obtain finite propellers 
  $\ds P_{\lambda(1 , \ell)}$,
  $\ds P_{\g(2 , \ell)}$ and
  $\ds P_{\lambda(2 , \ell)}$ for $\ell \geq 0$.

    Consider the corresponding four collections  of ``notched''  \emph{finite}  propellers
  \begin{equation}\label{eq-propellers2}
\left\{ P_{\g(1,\ell))}' \right\}_{\ell =0}^\infty ~ , ~  
\left\{ P_{\lambda(1,\ell)}' \right\}_{\ell =0}^\infty ~ , ~ 
\left\{ P_{\g(2,\ell)}' \right\}_{\ell =0}^\infty ~ , ~ 
\left\{ P_{\lambda(2,\ell)}' \right\}_{\ell=0}^\infty
\end{equation}
 defined by taking the intersection of the corresponding  propellers with $\mW'$. 
Glue  these notched propellers to $\fM_0^1$ using   $\sigma_1$ and $\sigma_2$, 
where the points added with this gluing are level 2 points, hence are contained in  $\fM_0^2$.

Each propeller in each of the four infinite collections in \eqref{eq-propellers2} yields   a family of curves in
 $\cL_1^-$  and another family of curves in  $\cL_2^-$. 
Since all the propellers we are
considering are finite, the number of curves in the intersection of
each propeller with  $\cL_1^-$ and $\cL_2^-$  is finite. 
Note that the number of notches in a given propeller at  level $n \geq 2$  may   be zero.
Moreover, for $i, j = 1,2$ and $\ell \geq 0$, each  base point  $p(i; j , \ell)$  of   the corresponding generating curve has radius greater than $2$, so the $\cW$-orbit of the labeling points traces out the full boundary of the propeller it defines. 
This is in contrast to the case with the level $1$ propellers, where both forward and backward orbits were required to reach all of the notches. We continue to use the same labeling scheme in levels larger than $1$ for the base points of notches corresponding to the intersections with the faces of $\cL_1^-$ and $\cL_2^-$. 

Following the previous scheme, each  entry region of a notch of
$P_{\g(i_1,\ell_1)}'$ and $P_{\lambda(i_1,\ell_1)}'$ defines a    curve denoted  by 
$\g(i_1 , \ell_1 ; i_2, \ell_2)$ and $\lambda(i_1 , \ell_1 ; i_2,
\ell_2)$ for $i_1, i_2 = 1,2$, where $\ell_1 , \ell_2 \geq 0$ and
$\ell_2$ is bounded above by  
$\max_{x\in\g(i_1,\ell_1)}\Delta(r(x)) + 1 < \infty$ or
$\max_{x\in\lambda(i_1,\ell_1)}\Delta(r(x)) + 1 < \infty$, accordingly.
The index
  $i_2$ indicates that the curve is in $\cL_{i_2}^-$, while the
  indices $(i_1, \ell_1)$ indicate in which notched propeller of level $2$ they are contained. 
  Hence, for example, the curve $\g(i_1,\ell_1 ; 2, \ell_2)$ is in
$\cL_2^-$ and belongs to $P_{\g(i_1,\ell_1)}'$.
  Some level 2 curves are represented in Figures~\ref{fig:intcurvesL1} and
\ref{fig:intcurvesL2}. The details of these pictures are considered further 
  in Section~\ref{sec-doublepropellers}. 
 
Corresponding to each curve defined by the intersection with an entry face, is the facing curve defined by the intersection with an exit face,   
denoted by $\ovg(i_1 , \ell_1 ; i_2, \ell_2)$ and $\ovl(i_1 , \ell_1 ; i_2, \ell_2)$.
Then   $\fM_0^2$ is obtained by attaching these exit curves to $\fM_0^1$, along with the level $2$ finite propellers in \eqref{eq-propellers2}.

The previous steps are now   repeated recursively. Given $\fM_0^{n}$ for $n \geq 2$, we introduce families of curves defined by the entry curves in the notches of the propellers at level $n$, 
\begin{equation}\label{eq-propellers-n}
 \g(i_1 , \ell_1 ; i_2, \ell_2; \cdots ; i_n, \ell_n) \quad , \quad 
  \lambda(i_1 , \ell_1 ; i_2, \ell_2; \cdots ; i_n, \ell_n)
\end{equation}
and their corresponding facing curves defined by the exit curves, for
$i_1, i_2 , \ldots , i_{n} = 1,2$ and   $\ell_i \geq 0$ which are
bounded,  except for $\ell_1$. The base points of these curves are
$$p(1;i_1 , \ell_1 ; i_2, \ell_2; \cdots ; i_n, \ell_n) ~ {\rm and ~} ~ p(2;i_1 , \ell_1 ; i_2, \ell_2; \cdots ; i_n, \ell_n) , ~ {\rm respectively} .$$
 As before, $i_{n}$ indicates that the
curve (or point) is in $\cL_{i_{n}}^-$ and the previous indices
$(i_1,\ell_1 ; i_2, \ell_2; \cdots ; i_{n-1}, \ell_{n-1})$ indicate the propeller that contains the curve. 

The curves in \eqref{eq-propellers-n} generate finite propellers and
the notched finite propellers associated to them. These are hence at
level $n+1$. Then $\fM_0^{n+1}$ is obtained from $\fM_0^{n}$ by attaching the    infinite families of finite propellers at level $n+1$ to each of the previously attached propellers in $\fM_0^n$, along with the exit curves at level $n+1$ in these attached propellers. Thus, we obtain the nested compact sets 
\begin{equation}\label{eq-nestedfamilies}
\fM_0^{0}\subset \fM_0^{1} \subset \cdots \fM_0^n \subset \fM_0^{n+1} \subset \cdots \subset \fM_0 \subset   \overline{\fM_0} =  \fM
\end{equation}

Observe that in the construction of $\fM_0^2$
 we added $4$ countable families of finite propellers, the images under $\tau$ of those in \eqref{eq-propellers2}. In general, in stage   $n$ of the construction, 
we add $2^n$ countable families of finite propellers, where each propeller is indexed by its base  point 
$p(i_0; i_1 , \ell_1 ; i_2, \ell_2; \cdots ; i_{n-1}, \ell_{n-1})$  and
  each family has the indices   $i_0, i_1,i_2,\ldots, i_{n-1}$ in
  common. However, for $(i_0; i_1 , \ell_1 ; i_2, \ell_2; \cdots ;
  i_{n-1})$ fixed, the number of  $\ell_{n-1}$ for which there exists
  propellers with base curve $\g(i_1 , \ell_1 ; i_2, \ell_2; \cdots ;
  i_{n-1}, \ell_{n-1})$ if $i_0=1$ or $\lambda(i_1 , \ell_1 ; i_2,
  \ell_2; \cdots ; i_{n-1}, \ell_{n-1})$ if $i_0=2$, is bounded above
  by $\Delta(r)$ where $r$ is the maximum radius of the points in the
  curve $\g(i_1 , \ell_1 ; i_2,
  \ell_2; \cdots ; i_{n-2}, \ell_{n-2})$ or $\lambda(i_1 , \ell_1 ; i_2,
  \ell_2; \cdots ; i_{n-2}, \ell_{n-2})$, accordingly.

\begin{figure}[!htbp]
\centering
{\includegraphics[width=120mm]{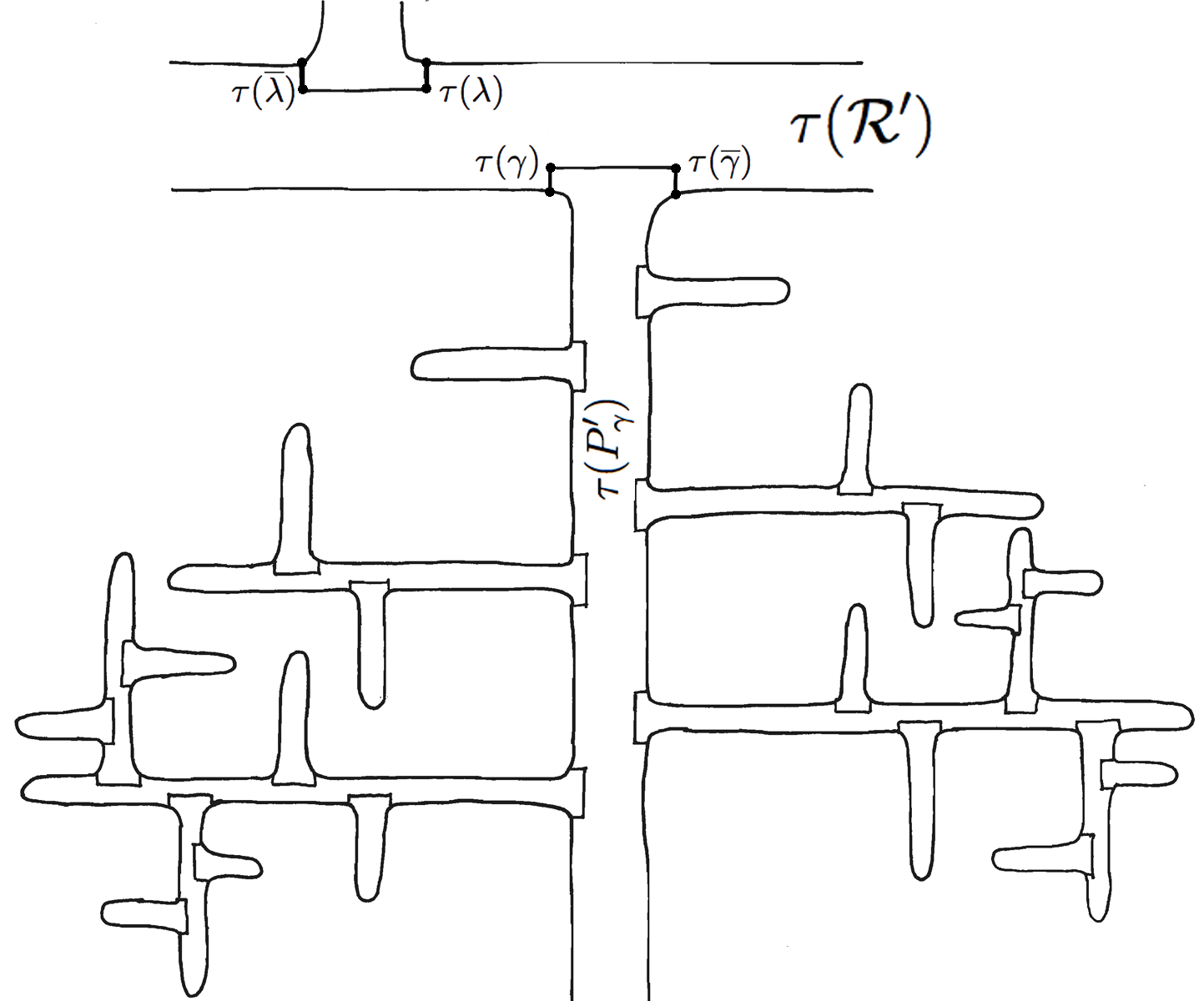}}
\caption{\label{fig:choufleur}  Flattened part of $\fM_0$}
\vspace{-6pt}
\end{figure}

As observed above, 
the set $\tau(\cR^\prime)$ is a cylinder, minus two rectangles, embedded as a ``folded figure eight'' as in 
 Figure~\ref{fig:notched8}. One then adds   two infinite
propellers that wrap around $\tau(\cR^\prime)$ to obtain $\fM_0^1$. 
The following steps in the construction add finite propellers to the infinite
ones in $\fM_0^1$. The boundary of these finite propellers contain all of the finite
\emph{chou-fleurs of Siebenmann}, introduced in \cite{Ghys1995}.
The term ``chou-fleur'' comes from the diagram of flattened
propellers. 
 In fact, if we draw $\cR'$ as a rectangle, then we can add
the flattened propellers as in Figure~\ref{fig:choufleur}. The boundary represents
$\cK$-orbits of the special points $p_1^-$ and $p_2^-$, and any finite
part lying between two facing transition points is a chou-fleur.

 We comment on the details of Figure~\ref{fig:choufleur}.
The higher horizontal band in Figure~\ref{fig:choufleur} represents the notched Reeb cylinder
$\cR'$ to which the two infinite propellers in $\fM_0^1$ are
attached. Part of the infinite propeller $\tau(P'_\g)$ is the main
vertical branch and only the base of the infinite propeller $\tau(P'_\lambda)$ is pictured 
in this diagram. The propeller $\tau(P'_\g)$ generates by its
intersections with the insertions, the finite propellers
in $\fM_0^2$, and so on. 
In $\mK$, the width of the infinite propellers is the same as the
width of $\cR'$. The finite propellers are roughly the same width as
the infinite ones, at least for a certain period of time.
Hypothesis~\ref{hyp-SRI}, implies that the finite propellers attached to
$\tau(P_\g')$ become longer as we move downwards and also that the branching
structure gets more and more complicated.

  \bigskip
 
\section{Double propellers and pseudogroup   dynamics}\label{sec-doublepropellers}

 In this section,   we   begin our analysis of  the dynamics of the Kuperberg flow on $\fM_0$ and the   dynamics of the pseudogroup $\cGK$ acting on $\bRt$. 
A key technique introduced in this section is the  the concept  of \emph{double propellers}, which are
obtained from the $\cK$-orbit of  special  ``parabolic''  curves in $\partial_h^- \mW$.   The double propellers  define  families of nested topological circles  in $\bRt$ which play a fundamental role  in the study of the   dynamics of the pseudogroup $\cGK$ acting on $\bRt$.

We first  introduce the general notion of double propellers in $\mW$,  then
 discuss     their properties. We then define a  labeling system for the
 double propellers generated by the flow $\Phi_t$,  which is based on the labeling system for propellers 
  introduced in the last section.  Finally, we consider the intersections of the double propellers with the rectangle $\bRt$, which generates families of  nested ``ellipses'', and introduce a modified labeling system for these curves which corresponds to the action of the generators of the $\psg$ $\cGK^*$ on $\bRt$.
 
  In the following, we  assume that the constant $b=0$,  for $b$ as introduced in
Section~\ref{sec-proplevels} to index the internal notches in the propeller and postpone the discussion    of the case $b <  0$ to
  Section~\ref{sec-bubbles}.

\begin{figure}[!htbp]
\centering
\begin{subfigure}[c]{0.4\textwidth}{\includegraphics[width=50mm]{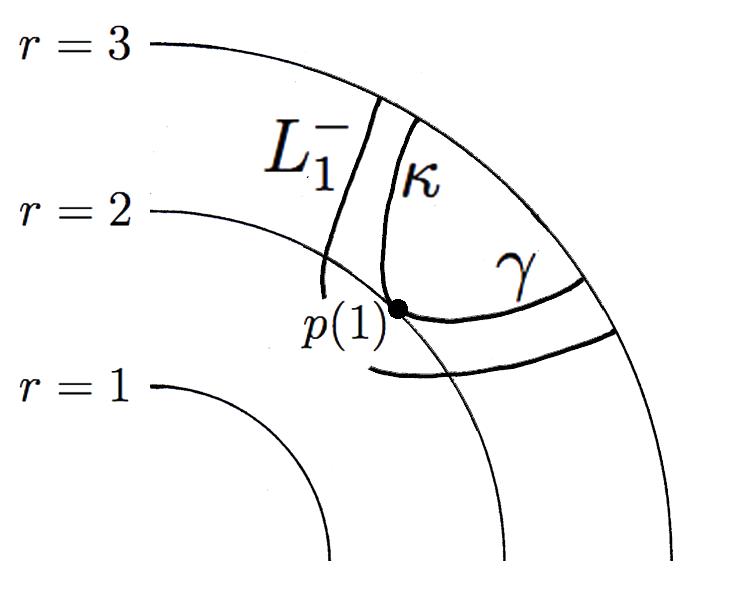}}\end{subfigure}
\begin{subfigure}[c]{0.4\textwidth}{\includegraphics[width=50mm]{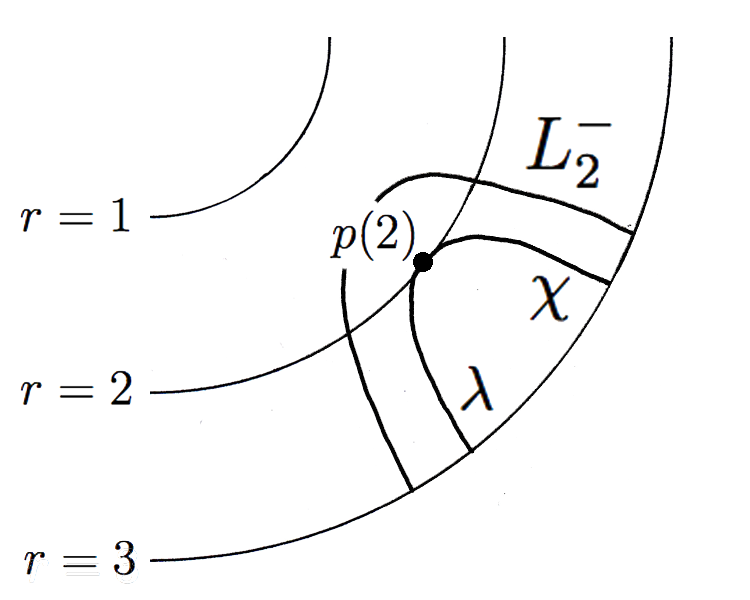}}\end{subfigure}
\caption{\label{fig:Gamma1}  The curves $\G=\g\cup \kappa$ and $\Lambda=\lambda\cup\chi$ in $\partial_h^-\mW$} 
\vspace{-6pt}
 \end{figure}

  Consider a smooth curve $\G \subset L_i^-$ parametrized by $u \colon [0,2] \to \partial_h^- \mW$, where $u_s \equiv u(s)$, such that:  
  \begin{enumerate}
\item $r(u_s) \geq 2$ for all $0 \leq s \leq 2$; 
\item  $r(u_0) = r(u_2) = 3$, so that both endpoints lie in the boundary $\partial_h^- \mW \cap \partial_v \mW$; 
 \item $\G$ is  topologically transverse to the   cylinders $\cC(r)$ for $2 \leq r \leq 3$,  except   at the midpoint $u_1$.
\end{enumerate}
It then follows that   $r(u_s) \geq  r(u_1) = 2+\e$ for all $0 \leq s \leq 2$, and some $\e\geq 0$.
 See Figure~\ref{fig:Gamma1} for an illustration in the case when $\e = 0$.

 Assume that $\e > 0$, so that $r(u_s) > 2$  for all $0 \leq s \leq 2$,   then    the $\cW$-orbit of each $u_s$  traverses $\mW$. 
 Define $T_s$ as  the exit time for the $\cW$-orbit of   $u_s$. 
 The $\cW$-orbits of the points in $\G$   form a    surface embedded in $\mW$,   whose boundary is contained in the boundary of $\mW$, and thus the surface separates $\mW$ into two connected components. This surface is denoted   $P_{\G}$ and called  the \emph{double propeller} defined by the $\Psi_t$-flow of $\G$.

Consider the   curves $\g, \kappa \subset   \partial_h^- \mW$ obtained by dividing 
  the  curve $\G$ into two  segments at the midpoint $s=1$.  Parametrize these curves as follows: 
\begin{eqnarray*}
\g = \G \ | \ [0,1] ~ & , &  ~   u(s) ~ \text{for} ~ 0 \leq s \leq 1 \\
\kappa= \G \  | \ [1,2] ~ & , &  ~  u(2-s) ~ \text{for} ~ 0 \leq s \leq 1 
\end{eqnarray*}
The orbit  $\ds \{\Psi_t(u_1) \mid 0 \leq t \leq T_1\}$ forms the \emph{long boundary} of the propellers $P_{\g}$ and $P_{\kappa}$ generated by the $\Psi_t$-flow of these curves. Then $P_{\G}$ is viewed as the gluing of   $P_{\g}$ and $P_{\kappa}$ along the long boundary, which forms a ``zipper'' joining the two surfaces together, hence  the notation ``double propeller'' for $P_{\G}$. Note that the length of the zipper  tends to infinity as $r(u_1) \to 2$.

 We comment on the details of Figure~\ref{fig:Gamma1}.
 Comparing Figure~\ref{fig:Gamma1} with Figure~\ref{fig:gamma'}, the $\g$ and
$\lambda$ curves are the first half arcs in $\G$ and
$\Lambda$ (first in the $\mS^1$ direction), while $\kappa$ and $\chi$
are the second half arcs in $\G$ and
$\Lambda$, respectively.

If $\e =0$, define   two infinite propellers   $P_{\g}$ and $P_{\kappa}$ 
as in Definition~\ref{def-infpropeller}, and then define 
    $\ds P_{\G} = P_{\g} \cup P_{\kappa}$, 
 where   the  $\Psi_t$-orbit $\cZ_{\g}$  of the midpoint $u_1$,    defined as in \eqref{eq-Zset},  is again common to both $P_{\g}$ and $P_{\kappa}$, and  $P_{\G}$ is viewed as the gluing of the two propellers along an ``infinite zipper''.

\begin{defn}\label{def-doublepropeller}
Let $\G$ be as above with $r(u_1)=2$. Let $\g^{\e}$  and $\kappa^{\e}$
for $0 <\e\leq 1$ be the curves as defined in \eqref{eq-sparem}. The  \emph{infinite double propeller} is the union:
\begin{equation}\label{eq-doublepropeller}
P_{\G} ~ \equiv ~ \cZ_{\g} ~ \cup ~ \bigcup_{\e > 0} ~ \left\{P_{\g^{\e}} \cup P_{\kappa^{\e}}\right\}  
\end{equation}
\end{defn}
Observe that $\overline{P_\G}=P_\G\cup\cR$ as in
Proposition~\ref{prop-propellerclosure}.

 We next apply these ideas to the study of the invariant set $\fM_0$.
  Recall     the curves $\g$ and $\lambda$ defined by
  \eqref{eq-curves} were obtained from the intersection of the   Reeb cylinder $\cR$ in $\mW$  with  the entry faces $\cL_1^-$ and $\cL_2^-$.  
Now   consider   the intersections of the    \emph{full cylinder} $\cC \equiv\{r=2\}$ with the entry faces $\cL_1^-$ and $\cL_2^-$, which  yields   the     curves:
\begin{equation}\label{eq-Gamma}
\G' = \cC \cap \cL_1^- \subset \mW \quad , \quad \G = \sigma_1^{-1}(\G') \subset  L_1^-
\end{equation}
 \vspace{-20pt} 
\begin{equation}\label{eq-Lambda}
\Lambda' = \cC \cap \cL_2^- \subset \mW \quad , \quad \Lambda = \sigma_2^{-1}(\Lambda') \subset  L_2^-
\end{equation}

Hypothesis~\ref{hyp-SRI}  implies that  the curves $\G$  and $\Lambda$
are  topologically transverse to the cylinders $\{r=const.\}$, except at
their middle points, where they are tangent to the cylinder $\cC$, as in
Figure~\ref{fig:Gamma1}.   The endpoints of each of the curves $\G$ and $\Lambda$  have $r$-coordinates equal to 3. 
For the curves $\g'$ and $\lambda'$ as defined in Section~\ref{sec-proplevels}, note that $\g' = \G' \cap \cR$ and $\lambda' = \Lambda' \cap \cR$. Also, set
\begin{equation}\label{eq-kappa}
\kappa = \sigma_1^{-1}(\G' \cap \{z \leq -1\}) \quad , \quad \chi = \sigma_2^{-1}(\Lambda' \cap \{z \geq 1\}) .
\end{equation}
  The graphs of the curves $\G$, $\g$ and $\kappa$ near the special
  point  $\tau^{-1}(p_1^-) \in L_1^-$ are   illustrated  in
  Figure~\ref{fig:Gamma1}, as well as the graphs of $\Lambda$, $\lambda$ and $\chi$ near the special point $\tau^{-1}(p_2^-) \in L_2^-$.

Let  $P_{\G}$ and $P_\Lambda$ be the   infinite double propellers associated to these curves by \eqref{eq-doublepropeller}.
The  curves $\G$ and $\Lambda$ are disjoint in $\partial_h^- \mW$, thus the  infinite double propellers  $P_{\G}$ and $P_{\Lambda}$    are     disjoint in $\mW$.
Consider then the  notched infinite double propellers,
\begin{equation}
P_{\G}' =P_{\G} \cap \mW^\prime \qquad
\mbox{and}\qquad P_{\Lambda}' =P_{\Lambda}\cap
\mW^\prime.
\end{equation}

\begin{lemma}\label{lem-notchedprops}
The notched propellers $P_{\G}'$ and $P_\Lambda' $ are tangent to
$\cC ' = \cC \cap \mW'$ along the $\cW$-orbit of the special points
$\tau^{-1}(p_1^\pm) \in L_1^{\pm}$ and $\tau^{-1}(p_2^\pm)    \in L_2^{\pm}$, respectively.
\end{lemma}
\proof 
The curves $\G$, $\overline{\G}$, $\Lambda$ and
$\overline{\Lambda}$ are tangent to the cylinder
$\cC'$ at their middle points that are the points $\tau^{-1}(p_1^\pm)$
and $\tau^{-1}(p_2^\pm)$, respectively. Hence $P_{\G}'$ and $P_\Lambda ' $ are tangent to
$\cC '$ along the orbits of these points.
\endproof

The bases of the propellers  $P_{\G}'$ and $P_{\Lambda}'$ are then ``glued'' to the notched cylinder $\cC'$ using the maps
$\sigma_1$ and $\sigma_2$, respectively. This   adds to the notched cylinder $\cC'$  two infinitely long ``tubes with notched holes''  that wrap
around $\cC'$,  accumulating on the notched Reeb cylinder
$\cR'$. Moreover, the tubes are tangent to $\cC'$ along the $\cW$-arcs
of the special orbits that are in $\cC'- \cR'$. 

 Consider the $\Phi_t$-flow of the image $\tau(\cC')\subset \mK$:
\begin{equation}\label{eq-doublelamination}
\widehat{\fM}_0 ~ \equiv ~  \{\Phi_t(\tau(\cC')) \mid -\infty < t < \infty \}.
\end{equation}
The level function $n_0$ introduced in Proposition~\ref{prop-levels} extends to a well-defined function on $\widehat{\fM}_0$,  using   the same method of proof, where $\tau(\cC')$ has level $0$ by definition. Thus,     $\widehat{\fM}_0$ can be decomposed into its level
sets:
\begin{equation}
\widehat{\fM}_0^n ~ = ~ \{ x \in \widehat{\fM}_0 \mid n_0(x) \leq n\} ~, ~ n = 0,1,2, \ldots
\end{equation}
Then $\widehat{\fM}_0^0=\tau(\cC'\cup \overline{\Gamma}\cup
\overline{\Lambda})$ as defined above. The set $\widehat{\fM}_0^1$ is
obtained by attaching to the notches of $\tau(\cC')$ the double
infinite propellers $\tau(P_\Gamma')$ and $\tau(P_\Lambda')$, as well
as part of the boundaries of their notches. 
In complete analogy with the construction of $\fM_0$,  repeat this process recursively  for each family of notches, to obtain an  embedded surface $\widehat{\fM}_0 \subset \mK$. Note that $\fM_0 \subset \widehat{\fM}_0$.

Next, we use the   labeling system for $\fM_0$  to describe the
intersection of $\widehat{\fM}_0$  with  $E_1$ and $E_2$. For this
section, we keep the
assumption $b=0$ that will be dropped in Section~\ref{sec-bubbles}.
Recall the notation convention of Remark~\ref{rmk-labels}.
 The  infinite sequence of points $p'(i_0; i_1 , \ell_1)$ for $\ell_1 \geq 0$ introduced in Section~\ref{sec-proplevels} all lie on the $\cW$-orbits of the points $\tau^{-1}(p_{i_0}^\pm)$ which are contained in the open ``half-cylinders'',  either  $\ds \cC \cap \{z < -1\}$ or $\ds \cC \cap   \{z > 1\}$.

The   $\cK$-orbits  of the points $\tau(p(i_0; i_1 , \ell_1))$ for $i_0, i_1 = 1,2$ and $\ell_1 \geq 0$,  yield the families of points
 labeled as $p(i_0; i_1 , \ell_1 ; i_2, \ell_2; \cdots ; i_n,
 \ell_n)$, which are the  basepoints for the     curves introduced in
 Section~\ref{sec-proplevels}. The double propellers also share these
 same basepoints, the basepoint of a double propeller is the ``middle''
 point in the curve generating it. We   adopt the same notation system.

For $i=1,2$, apply $\sigma_i^{-1}$  to  the intersections   of   $P_{\G}$ and $P_{\Lambda}$  with the faces $\cL_i^-$   to obtain, in $L_1^- \cup L_2^- \subset \partial_h^- \mW$, four countable collections of curves, labeled    in a manner corresponding to that used for the curves $\g(i,\ell)$ and $\lambda(i,\ell)$ in \eqref{eq-propellers2}. For $\ell \geq 0$, we set:
\begin{itemize}
\item $\ds  \G(1,\ell)  = \g(1 , \ell)  \cup \kappa(1 , \ell) ~ \subset ~  L_1^-$ ~, ~ based at $p(1; 1 , \ell)$  and corresponding to $P_{\G}$; 
\item $\ds  \G(2,\ell)  = \g(2 , \ell)  \cup \kappa(2 , \ell) ~ \subset ~  L_2^-$ ~, ~ based at $p(1; 2 , \ell)$  and corresponding to $P_{\G}$; 
\item $\ds  \Lambda(1,\ell)  =  \lambda(1 , \ell) \cup \chi(1 , \ell)  ~ \subset ~  L_1^-$~, ~ based at $p(2; 1 , \ell)$  and corresponding to $P_{\Lambda}$; 
\item $\ds  \Lambda(2,\ell)  =  \lambda(2 , \ell) \cup \chi(2 , \ell)  ~ \subset ~  L_2^-$~, ~ based at $p(2; 2 , \ell)$  and corresponding to $P_{\Lambda}$. 
\end{itemize}
The endpoints of each   of the curves $\G(i,\ell)$ and  $\Lambda(j,\ell) $ are contained in the boundary  of  $\partial_h^- \mW$, while  the midpoints are   endpoints for   $\g(i, \ell)$ and $\kappa(i, \ell)$,  or $\lambda(j , \ell)$ and $\chi(j, \ell)$, accordingly. 

Note  that   the curves $\G'(i,\ell)$ and  $\Lambda'(i ,\ell)$ in the faces $\cL_i^-$ are
the result of applying the $\Psi_t$-flow to the curves $\G$ and
$\Lambda$. 
 We thus   obtain
four countable collections of   double propellers
$P_{\G(i,\ell)}$ and $P_{\Lambda(j,\ell)}$, for $i, j=1,2$, and
  notched  double propellers $P_{\G(i,\ell)}'$ and
$P_{\Lambda(j,\ell)}'$.

  As in Section~\ref{sec-proplevels},   this creation and labeling of
  double propellers proceeds recursively, though for levels $n \geq
  2$, the propellers obtained are finite as the base curves are
  contained in the region $r > 2$. 
 There  is a nuance that arises with the construction of the corresponding
  double propellers, however. The   curves  $\kappa(i_1 , \ell_1 ; i_2,
  \ell_2; \cdots ; i_n, \ell_n) \subset L_{i_n}^-$ and $\g(i_1 , \ell_1 ; i_2, \ell_2; \cdots ; i_n,
  \ell_n) \subset L_{i_n}^-$ do not satisfy the radial monotonicity
  hypothesis.  In particular,   their intersection point, the
  ``midpoint'' $u_1$ will not   be the point where $r(u_s)$ is
  minimal.  However, the $\kappa$ and  $\g$  curves  are endpoint isotopic
  to a curve satisfying the transversality condition, hence the
  conclusion that the corresponding curve $\G$ separates the region $L_i^-$ will remain true, so that the double propeller $P_{\G(i_1,\ell_1;\cdots;i_n,\ell_n)}$ is  isotopic to a standard double propeller, hence   will separate $\mW$ into two connected components as well.  Similar remarks also hold for the  double propeller $P_{\Lambda(i_1,\ell_1;\cdots;i_n,\ell_n)}$, and we 
 label the base curves resulting from the $\Phi_t$-flow:
\begin{itemize}
\item $\ds  \g(i_1 , \ell_1 ; i_2, \ell_2; \cdots ; i_n, \ell_n) ~
  \subset ~ \G(i_1 , \ell_1 ; i_2, \ell_2; \cdots ; i_n, \ell_n)
  ~ \subset ~   L_{i_n}^-$ ~,\\  \quad which is based at $p(1 ; i_1 , \ell_1 ; i_2, \ell_2; \cdots ; i_n, \ell_n)$; 
\item $\ds  \kappa(i_1 , \ell_1 ; i_2, \ell_2; \cdots ; i_n, \ell_n) ~
  \subset ~ \G(i_1 , \ell_1 ; i_2, \ell_2; \cdots ; i_n, \ell_n)
  ~ \subset ~   L_{i_n}^-$ ~, \\  \quad which is based at $p(1 ; i_1 , \ell_1 ; i_2, \ell_2; \cdots ; i_n, \ell_n)$; 
\item $\ds  \lambda(i_1 , \ell_1 ; i_2, \ell_2; \cdots ; i_n, \ell_n)
  ~ \subset ~  \Lambda(i_1 , \ell_1 ; i_2, \ell_2; \cdots ; i_n,  \ell_n) ~ \subset ~       L_{i_n}^-$~, \\  \quad which is based at $p(2 ; i_1 , \ell_1 ; i_2, \ell_2; \cdots ; i_n, \ell_n)$; 
\item $\ds  \chi(i_1 , \ell_1 ; i_2, \ell_2; \cdots ; i_n, \ell_n)
  ~ \subset ~  \Lambda(i_1 , \ell_1 ; i_2, \ell_2; \cdots ; i_n,  \ell_n) ~ \subset ~       L_{i_n}^-$~, \\  \quad which is based at $p(2 ; i_1 , \ell_1 ; i_2, \ell_2; \cdots ; i_n, \ell_n)$. 
\end{itemize}
  The shape of $\kappa$, $\g$, $\chi$ and $\lambda$ curves will be
 important in latter sections, and analyzed  in 
 Section~\ref{sec-geometry}.

There is a tangency relation between double propellers with
consecutive levels. Consider the curve $\G\subset L_1^-$ and the curves
$\G(1,\ell)\subset L_1^-$. The base point $p'(1;1,\ell)$ of
$\G'(1,\ell)\subset \cL_1^-$ has radius equal to 2 and thus belongs to
$\cC\cap \cL_1^-$. Then $p(1;1,\ell)$ is a point in $\kappa$ and every
$\G(1,\ell)$ curve is tangent to $\kappa$ at
$p(1;1,\ell)$. Analogously,
\begin{itemize}
\item $\Lambda(1,\ell)$ is tangent to $\kappa$ at $p(2;1,\ell)$;
\item $\G(2,\ell)$ is tangent to $\chi$ at $p(1;2,\ell)$;
\item $\Lambda(2,\ell)$ is tangent to $\chi$ at $p(2;2,\ell)$.
\end{itemize}
Thus the level 2 propellers $P_{\G(i_1,\ell_1)}$ and
$P_{\Lambda(i_1,\ell_1)}$ are tangent along the $\cW$-orbit of
$p(1;i_1,\ell_1)$ and $p(2;i_1,\ell_1)$, respectively, to the level 1
propeller $P_\G$ if $i_1=1$ and $P_\Lambda$ if $i_1=2$.
This is illustrated in   Figure~\ref{fig:level2tangencies} below.

  \begin{figure}[!htbp]
\centering
{\includegraphics[width=100mm]{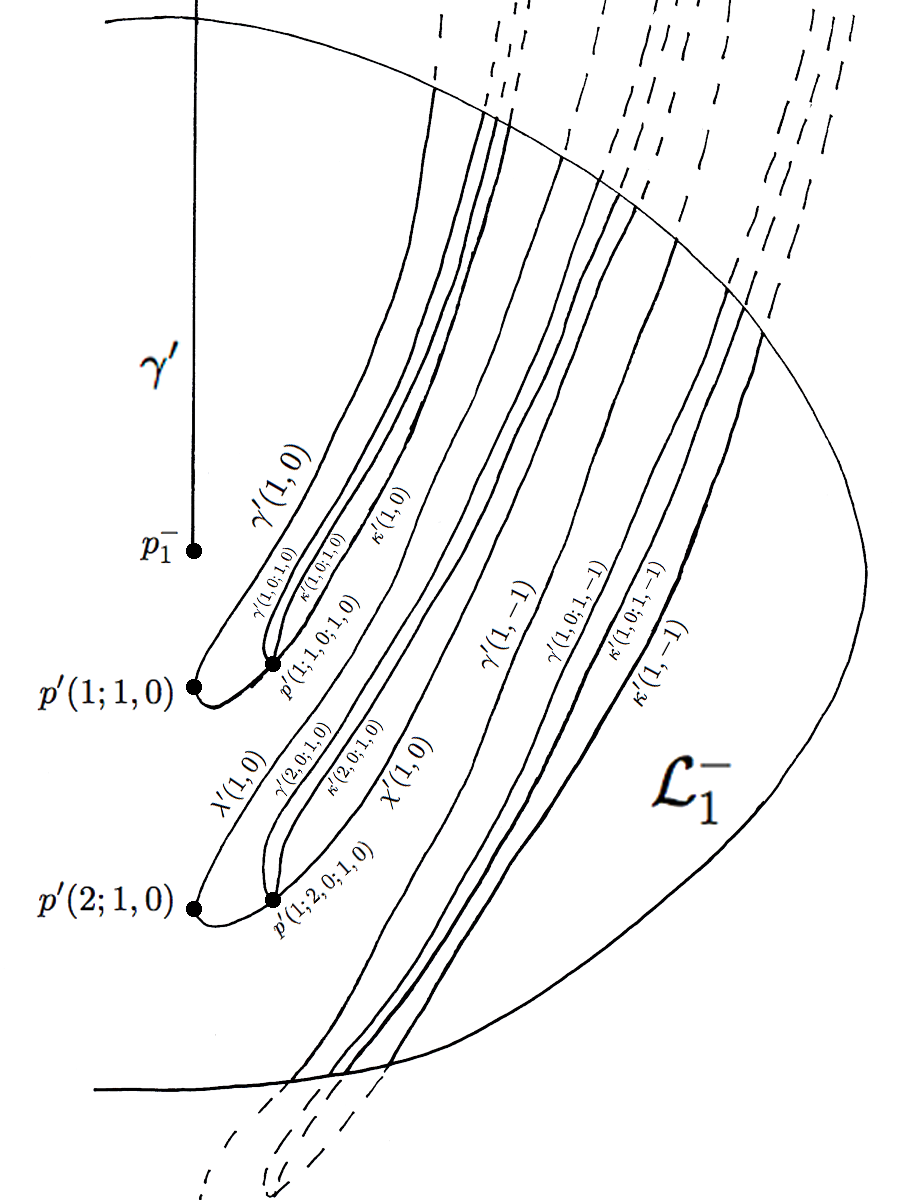}}
\caption{\label{fig:level2tangencies}  Tangencies of level 2 curves with level 1 curves in   $\cL_1^-$.}
\vspace{-6pt}
\end{figure}

In general, a $\G$ or $\Lambda$ curve at level $n$ is tangent to a
$\kappa$ or $\chi$ curve at level $n-1$, more precisely:
\begin{itemize}
\item $\G(1,\ell_1;i_2,\ell_2;\cdots;i_n,\ell_n)\subset L_{i_n}^-$ is
  tangent to $\kappa(i_2,\ell_2;\cdots;i_n,\ell_n)$ at
  $p(1;1,\ell_1;i_2,\ell_2;\cdots;i_n,\ell_n)$;
\item $\G(2,\ell_1;i_2,\ell_2;\cdots;i_n,\ell_n)\subset L_{i_n}^-$ is
  tangent to $\chi(i_2,\ell_2;\cdots;i_n,\ell_n)$ at
  $p(1;2,\ell_1;i_2,\ell_2;\cdots;i_n,\ell_n)$;
\item $\Lambda(1,\ell_1;i_2,\ell_2;\cdots;i_n,\ell_n)\subset L_{i_n}^-$ is
  tangent to $\kappa(i_2,\ell_2;\cdots;i_n,\ell_n)$ at $p(2;1,\ell_1;i_2,\ell_2;\cdots;i_n,\ell_n)$;
\item $\Lambda(2,\ell_1;i_2,\ell_2;\cdots;i_n,\ell_n)\subset L_{i_n}^-$ is
  tangent to $\chi(i_2,\ell_2;\cdots;i_n,\ell_n)$ at $p(2;2,\ell_1;i_2,\ell_2;\cdots;i_n,\ell_n)$.
\end{itemize}
Recall that the points $p(i_0;i_1,\ell_1)$ converge as
$\ell_1\to\infty$ to $p(i_1)$ as described in
Section~\ref{sec-proplevels}. This type of convergence is repeated at
any level. The level $n$ points
$p(i_0;i_1,\ell_1;i_2,\ell_2;\cdots;i_n,\ell_n)$ that belong to the
curve $\kappa(i_2,\ell_2;\cdots;i_n,\ell_n)$  if $i_1=1$ and to the
curve $\chi(i_2,\ell_2;\cdots;i_n,\ell_n)$ if $i_1=2$, converge to
$p(i_1;i_2,\ell_2;\cdots;i_n,\ell_n)$ as $\ell_1\to \infty$. Observe
that the last point is the base point of the $\kappa$ or $\chi$ curve.

 We next consider the intersections of the double propeller surfaces with the rectangle $\bRt$. Note that the $\Psi_t$-flow of a parabolic curve $\Gamma \subset \{r > 2\} \cap \partial_h^- \mW$ with endpoints in $\{r=3\}$ forms a compact surface in $\mW$ with boundary contained in $\partial \mW$, and the vector field $\cW$ is tangent to this surface. Thus, its intersections with the rectangle $\bRt$ are always transverse, so must be a finite union of closed curves. Moreover, the $z$-symmetry of the vector field $\cW$ implies these intersections are symmetric with respect to the horizontal line $\{z=0\} \cap \bRt$.

\begin{figure}[!htbp]
\centering
{\includegraphics[width=100mm]{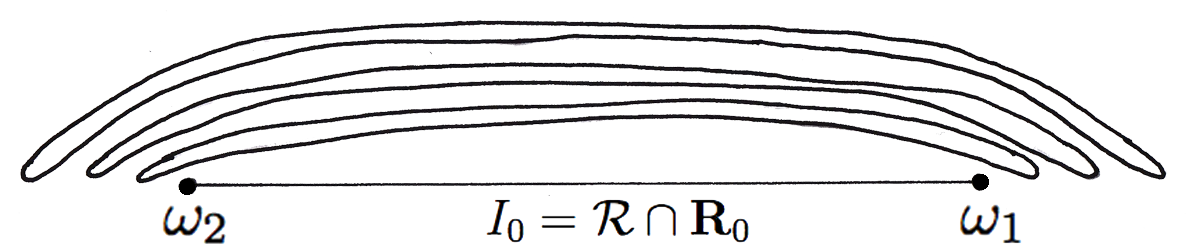}}
\caption[{Trace of an infinite double propeller in $\bRt$}]{\label{fig:ellipsespropeller} Trace of an infinite double propeller in $\bRt$ (viewed sideways)}
\end{figure}

 In the case of the parabolic curves $\G = \g \cup \kappa \subset L_1^-$ and $\Lambda = \lambda \cup \chi \subset L_2^-$,    the intersections of the   propellers generated by the curves $\g$ and $\kappa$ with the set $\bRt$ yields a family of arcs,  as illustrated in  Figure~\ref{fig:arcspropeller}. 
Thus, the intersections with $\bRt$  of the double propeller formed from their $\Psi_t$-flows   yields a family of ``thin circles'', 
as illustrated by Figure~\ref{fig:ellipsespropeller}. We call  these curves ``ellipses'' for  short.

\vspace{-6pt}

We next  develop a labeling system for  the collection of ellipses in  $P_{\Gamma} \cap \bRt$ of the form illustrated in Figure~\ref{fig:ellipsespropeller}. 
First, introduce the following   line segments contained in the intersection $\cC_0\equiv\cC \cap \bRt$:
\begin{equation}\label{eq-intervals}
I_{0}  = \left\{(2,\pi ,z) \mid  -1 \leq z \leq 1 \right\}   ,    J_{0}  = \left\{(2,\pi,z) \mid  -2 \leq z \leq -1 \right\}   ,    K_{0}  = \left\{(2,\pi,z) \mid  1 \leq z \leq 2 \right\} .
\end{equation}
$$N_0 = J_0 \cup I_0 \quad , \quad M_0 =  I_0 \cup K_0$$
The endpoints of the interval $I_0$ are the two special points, $\omega_i = \cO_i \cap \bRt$, for $i=1,2$. The arcs  $N_0$ and $M_0$ are connected curves in $\bRt$ which contain  these special points   as midpoints.

 Proposition~\ref{prop-propexist} implies that the $\Phi_t$-flow of $I_0$
 contains the ``notched'' Reeb cylinder $\tau(\cR')$, thus the
 $\Phi_t$-flow of $I_0$ equals    $\fM_0$. 
The $\Phi_t$-flow of $\cC_0$ contains the double propellers.

The forward $\Phi_t$-flow of $N_0$ intersects  $E_1$ in
$\tau(\G)$ and   generates the embedded double propeller $\tau(P_{\G}')\subset \mK$. 
The intersection $\tau(P_{\G}')\cap
\bRt$ is   a countable family of closed curves labeled as 
$\G_0(\ell)$   that are tangent to
$J_0$ along the forward $\Phi_t$-orbit of $\omega_1$ and to $K_0$ along
the backward $\Phi_t$-orbit of $\omega_1$.

The  forward $\Phi_t$-flow of $M_0$  intersects  $E_2$ in
$\tau(\Lambda)$ and   generates the embedded double propeller $\tau(P_{\Lambda}')\subset \mK$. 
The intersection $\tau(P_{\Lambda}')\cap
\bRt$ is   a countable family of
closed curves $\Lambda_0(\ell)\subset \bRt$  that are tangent to
$J_0$ along the forward $\Phi_t$-orbit of $\omega_2$  and to $K_0$ along
the backward $\Phi_t$-orbit of $\omega_2$. 

 Note that in the above and later, the subscript $``0"$    indicates by convention that a curve is in $\bRt$, and so may be considered as belonging to $\mW$ or $\mK$, according to the context.

The indexing $\ell$ of the closed curves $\G_0(\ell)$  and  $\Lambda_0(\ell)$  and their images under the   $\Phi_t$-flow, 
involves a   subtlety which is analogous to that arising in the proof of  Proposition~\ref{prop-generators}.
 It is possible that some  of the resulting curves from this flow will intersect in forward time $E_1$ and $E_2$, before
returning to $\bRt$. 
Thus,  the labeling for curves in these faces as constructed in
Section~\ref{sec-intropropellers} and above, may include some which do
not intersect $E_1\cup E_2$ before returning to $\bRt$.
This is the motivation for the following labeling convention.

Let    $\G_0(a)$, for    $a \leq 0$, denote   the first such  curve $\G_0(\ell)$  that intersects $\bRt$, 
whose intersection with the vertical line $\{r=2\}\subset \bRt$ is
composed of the two vertex points $p_0(1;1,a)$ and $p_0(1;2,a)$. 

Let   $\G_0(0)$ be the first such curve for which the $\cK$-orbit of its lower vertex point $p_0(1;1,0)$ 
intersects $E_1$ before  intersecting   $\bRt$ again. 
We assume, without loss of generality, that the $\cK$-orbit of $p_0(1;2,0)$ intersects $E_2$ before  intersecting   $\bRt$ again. 
 This symmetry assumption  simplifies the  labeling system.

 \begin{remark}\label{rmk-notation} 
The indexing of the curves starts with  a value of $\ell=a$,
possibly negative, such that for $a \leq \ell < 0$,  the $\cK$-orbit  of the vertices of the curves $\G_0(\ell)$ do not 
intercept $E_1\cup E_2$ in forward time,  before returning to $\bRt$. 
The number of such non-positive  indices $\ell$  is
finite, or might be zero.  In the latter case,  the indexing starts at
$\ell =0$, and it signifies that the forward $\Phi_t$-flow  of the
vertices of $\G_0(\ell)$   intersects the entry regions $E_1$ or $E_2$ in curves $\G(1,\ell)$ or $\G(2,\ell)$, respectively,
before intersecting $\bRt$ again. 
Similar comments apply for the indexing  of  the $\Lambda_0$ curves.
\end{remark}

Observe that $a\leq b\leq 0$ for $b$ as introduced in
Section~\ref{sec-proplevels}. As we   assume that $b=0$, then for $a\leq \ell<0$,  the $\Phi_t$-orbit segment starting at
any point of $\G_0(\ell)$ and ending at the corresponding point of
$\G_0(\ell+1)$ will not contain transition points.

Next consider the embedded propellers
$\tau(P_{\G(1,\ell_1)}')$ which are the image under $\tau$ of compact surfaces traversing $\mW$ from the bottom face $\partial_h^- \mW$ to the top face $\partial_h^+ \mW$.
They  intersect $\cL_1^-$ and $\cL_2^-$ in   finite families of closed curves, denoted by $\G'(1,\ell_1; 1, \ell_2)$ and $\G'(1,\ell_1 ; 2, \ell_2)$, and whose 
 inverse image in $L_i^-$ is denoted by $\G(1,\ell_1 ; i,
\ell_2)$, for $i=1,2$ and $\ell_2 \geq 0$. 
We form  the corresponding finite families of closed curves in $\bRt$ denoted by $\G_0(1, \ell_1 ; \ell_2)$,  
whose positive $\Phi_t$-flows  intersect $E_i$   in the curves $\tau(\Gamma(1,\ell_1;i,\ell_2))$ for $i=1,2$. The index $i$ signifies whether we follow the $\Phi_t$-flow of the curve $\G_0(1, \ell_1 ; \ell_2)$ in the region $\{ z < 0\}$ to the surface $E_1$, or we follow the $\Phi_t$-flow into the region $\{z > 0\}$ to the surface   $E_2$.

\begin{figure}[!htbp]
\centering
{\includegraphics[width=140mm]{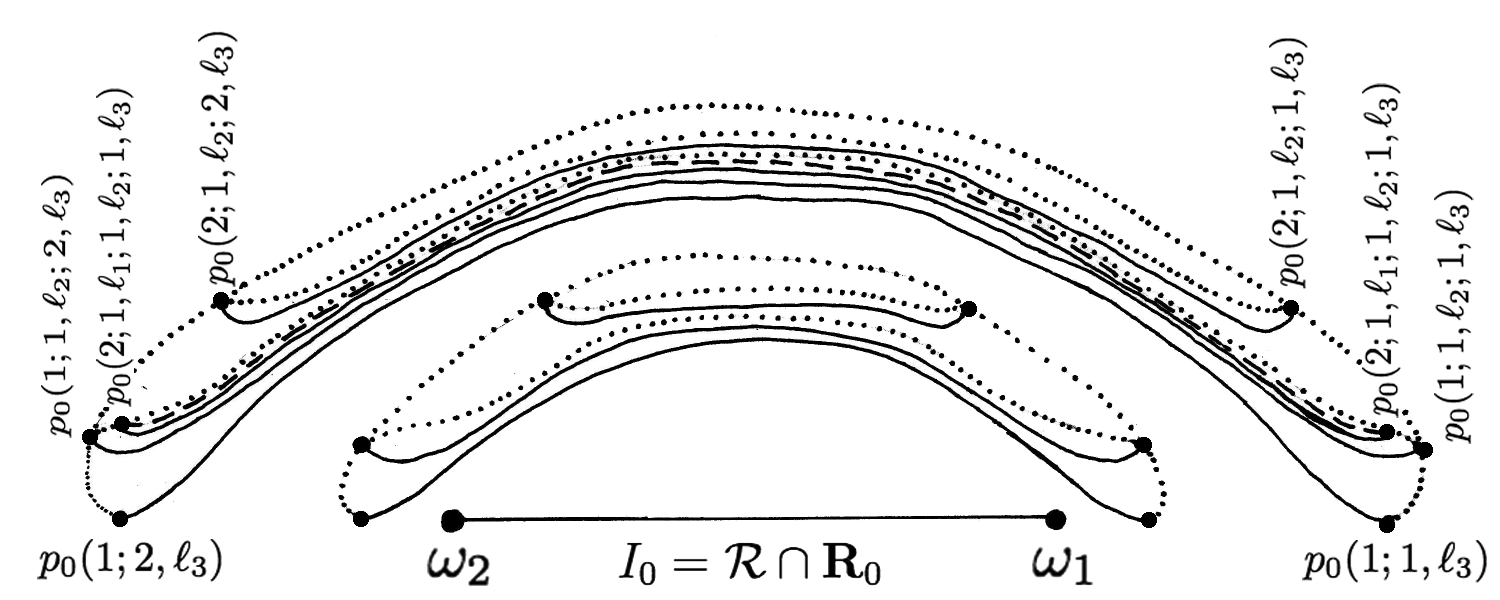}}
\caption[{The $\kappa_0$ and $\chi_0$ curves of levels 1 and 2}]{\label{fig:curvesR0} Curves in $\bRt$ (viewed sideways).
The $\kappa_0$ and $\chi_0$ curves of levels 1 and 2 are represented by dotted lines,
while the dashed line is a level 3 $\kappa_0$-curve.}
\vspace{-6pt}
\end{figure}

 Similarly, the    embedded propellers $\tau(P'_{\Gamma(2,\ell_1)})$  intersect $\bR_0$ in   finite families of closed
curves $\G_0(2, \ell_1 ; \ell_2)$, as illustrated in    Figure~\ref{fig:curvesR0}.
The  positive $\Phi_t$-flows of these curves in $\bRt$  intersects $E_1$ and $E_2$
in   curves denoted by $\tau(\G(2,\ell_1; 1, \ell_2))$ and $\tau(\G(2,\ell_1 ; 2, \ell_2))$
for $\ell_2\geq 0$,
respectively. 
These curves define  in turn, propellers that intersect $\bRt$ along closed curves, which are recursively defined as:
\begin{enumerate}
\item \label{item-sectionlabel1} $\ds  \G_0(i_1, \ell_1 ; i_2, \ell_2; \cdots ; i_{n-1},
  \ell_{n-1} ; \ell_n) ~
  \subset ~ \bRt$ ~   containing the ``vertex points'' \\
  $p_0(1 ; i_1 , \ell_1 ;  i_2, \ell_2; \cdots ; i_n, \ell_n)$  for $i_n =1,2$ and with   level   $n$; \\
\item \label{item-sectionlabel2}  $\ds   \Lambda_0(i_1 , \ell_1 ; i_2, \ell_2; \cdots ; i_{n-1},
  \ell_{n-1} ;  \ell_n) ~  \subset ~ \bRt $~  containing the ``vertex points''  \\
  $p_0(2 ; i_1 , \ell_1 ;  i_2, \ell_2; \cdots ; i_n, \ell_n)$   for $i_n =1,2$ and with  level   $n$. 
\end{enumerate}

The points 
$p_0(i_0 ; i_1 , \ell_1 ;  i_2, \ell_2; \cdots ; i_n, \ell_n)$ correspond to the points 
$p(i_0 ; i_1 , \ell_1 ;  i_2, \ell_2; \cdots ; i_n, \ell_n)$ in the same manner as the curves 
$\G_0(i_1 , \ell_1 ;  i_2, \ell_2; \cdots ; \ell_n)$ correspond to both 
$$\G(i_1 , \ell_1 ;  i_2, \ell_2; \cdots ; 1, \ell_n) ~ {\rm and } ~ \G(i_1 , \ell_1 ; i_2, \ell_2; \cdots ; 2, \ell_n). $$
 We also introduce the curves in $\bRt$:
\begin{enumerate}
\item \label{item-sectionlabelkappa1} 
$\ds \g_0(i_1, \ell_1 ; i_2, \ell_2; \cdots ; i_{n-1},  \ell_{n-1} ; \ell_n)$  and  
$\kappa_0(i_1, \ell_1 ; i_2, \ell_2; \cdots ; i_{n-1},  \ell_{n-1} ; \ell_n)$  contained in\\ 
$\G_0(i_1, \ell_1 ; i_2, \ell_2; \cdots ; i_{n-1}, \ell_{n-1} ; \ell_n)  \subset  \bRt$,  with  
$p_0(1 ; i_1 , \ell_1 ;  i_2, \ell_2; \cdots ; i_n, \ell_n)$  as
common boundary points for $i_n=1,2$; \\
\item \label{item-sectionlabelkappa2} 
$\ds \lambda_0(i_1, \ell_1 ; i_2, \ell_2; \cdots ; i_{n-1}, \ell_{n-1} ; \ell_n)$  and  
$\chi_0(i_1, \ell_1 ; i_2, \ell_2; \cdots ; i_{n-1},  \ell_{n-1} ; \ell_n)$  contained in \\ 
$\Lambda_0(i_1, \ell_1 ; i_2, \ell_2; \cdots ; i_{n-1}, \ell_{n-1} ; \ell_n)  \subset  \bRt$,  
with  $p_0(2 ; i_1 , \ell_1 ;  i_2, \ell_2; \cdots ; i_n, \ell_n)$  as
common boundary  points for $i_n=1,2$.
\end{enumerate}

\begin{remark}
For points of level at least 2, the number of positive indices
$\ell_n$ such that \\ $p_0(i_0;i_1,\ell_1;\cdots;i_n,\ell_n)$ exists is
greater or equal to the number of indices such that the point \\
$p(i_0;i_1,\ell_1;\cdots;i_n,\ell_n)$ exists, since the curves in the
the intersection of propellers at level at least 2 get shorter as we
approach the tip of the propeller and might not intersect the insertions.
\end{remark}

The $\g_0$-curves and $\lambda_0$-curves are in the $\Phi_t$-flow of
$I_0$, while the $\kappa_0$ and $\chi_0$ curves are in
the $\Phi_t$-flow of $J_0$ and $K_0$, respectively. Each
$\G_0(\ell)$ is tangent to $J_0$ at $p_0(1;1,\ell)$  and to $K_0$ at
$p_0(1;2,\ell)$. Then, we obtain that the $\G_0$
 curves at level $n$ are tangent to a $\kappa_0$
or a $\chi_0$ curve of level $n-1$. Likewise, each $\Lambda_0$
curve at level $n$ is tangent to a $\kappa_0$
and a $\chi_0$ curve of level $n-1$, as for curves in $L_i^-$.

Observe that this tangency relation implies that the curves
$\G_0(1,\ell_1;\ell_2)$ and $\Lambda_0(1,\ell_1;\ell_2)$ are tangent
to $\kappa_0(\ell_2)$ at their vertex points. Hence
$\G_0(1,\ell_1;\ell_2)$ and $\Lambda_0(1,\ell_1;\ell_2)$ are inside the
region bounded by $\G_0(\ell_2)$, for any $\ell_1\geq a$. Likewise,
$\G_0(2,\ell_1;\ell_2)$ and $\Lambda_0(2,\ell_1;\ell_2)$ are tangent
to $\chi_0(\ell_2)$ at their vertex points and lie inside the region
bounded by $\Lambda_0(\ell_2)$ for any $\ell_1\geq a$. Iterating this
relation we have that:
\begin{itemize}
\item the curves $\G_0(1,\ell_1;i_2,\ell_2;\cdots;\ell_n)$ and
  $\Lambda_0(1,\ell_1;i_2,\ell_2;\cdots;\ell_n)$ are tangent at their
  vertex points to $\kappa_0(i_2,\ell_2;\cdots;\ell_n)$ and lie inside
  the region bounded by $\G_0(i_2,\ell_2;\cdots;\ell_n)$;
\item the curves $\G_0(2,\ell_1;i_2,\ell_2;\cdots;\ell_n)$ and
  $\Lambda_0(2,\ell_1;i_2,\ell_2;\cdots;\ell_n)$ are tangent at their
  vertex points to $\chi_0(i_2,\ell_2;\cdots;\ell_n)$ and lie inside
  the region bounded by $\Lambda_0(i_2,\ell_2;\cdots;\ell_n)$.
\end{itemize}
Thus the $\G_0$ and $\Lambda_0$ curves form families of nested
ellipses. By construction, inside the region bounded by $\G_0(\ell_n)$
are all the $\G_0$ and $\Lambda_0$ curves whose first $i$-index is
equal to 1 and whose last $\ell$-index is equal to $\ell_n$. This
relation, up to level 3, is illustrated in Figure~\ref{fig:curvesR0}.

 Thus, the  intersection $\widehat{\fM}_0 \cap \bRt$ contains   collections of ``ellipses'' in  $\bRt$, indexed as above by their vertex points.
We next  examine the action of  $\cGK^*$ on  these collections of ellipses in  $\bRt$. 
  The  generators $\ds \{ \phi_1^\pm, \phi_2^\pm, \psi \}$ of $\cGK^*$ have domains which are subsets of $\bRt$ that do not contain all of the ellipses themselves,  so   we  define their actions on the ellipses  in terms of the actions on the vertex points of these    curves.
Observe that all these curves are contained in $\bRt\cap \{r\geq 2\}$. 
A corollary of the proof of Proposition~\ref{prop-sigmadyn} is the following:

\begin{lemma}\label{lem-invariantcurves} 
For each family of curves defined by a member  of    $\ds \{\G_0 , \gamma_0 , \kappa_0 , \Lambda_0 , \lambda_0 ,  \chi_0\}$, 
and for each      $\vp \in \cGK^*$, the action of $\vp$ maps   the family of  curves   to a subset of itself.
\end{lemma}

For $n >1$, recall that the last index $\ell_n$ for a curve $\G_0(i_1,\ell_1 ; \cdots ; i_{n-1},\ell_{n-1} ;\ell_n)$
is bounded above, with the bound depending on the previous indices, and by abuse of notation, 
   let  $\ell_n' = \ell_n'(i_1,\ell_1 ; \cdots ; i_{n-1},\ell_{n-1}) $ denote its maximum value. 
   Recall that there is no restriction on the first index, so    $\ell_1' = \infty$.

Proposition~\ref{prop-sigmadyn}  
implies that each of the sets of points $\{p_0(1;\cdots)\}$ and
$\{p_0(2;\cdots)\}$ is invariant, thus  the action on these points   determines which curve in each of these families is mapped to which curve, and correspondingly, the action on the curves determines the action on these vertices.

\begin{figure}[!htbp]
\centering
{\includegraphics[height=66mm]{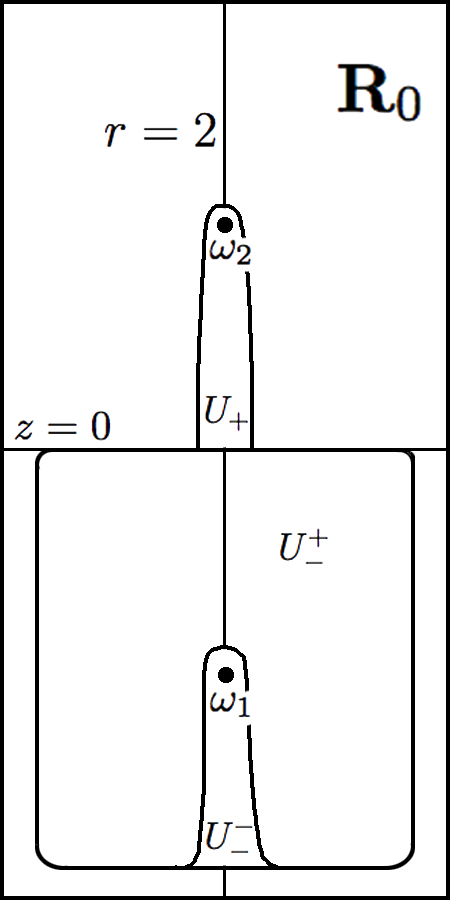}}
\caption{\label{fig:domainspsi333} Domains of continuity for $\psi$}
\vspace{-6pt}
\end{figure}

We first  analyze the action of the map $\psi$. Recall from
Section~\ref{sec-pseudogroup} that the map $\psi$ is obtained from the
return map
of the Wilson flow to $\bRt$, hence $\psi$ preserves the level
function $n_0$. The intervals $I_0$, $J_0$ and $K_0$ are
mapped to themselves, with   fixed-points $\omega_i$, for
$i=1,2$. Consider the open sets in the domain of $\psi$ defined by $U_+\subset \bRt\cap
\{r\geq 2, z> 0\}$ and $U_-\subset \bRt\cap \{r\geq 2, z< 0\}$. The set   $U_-$ further decomposes into open sets defined by the range of $\psi$,  
\begin{equation}\label{eq-domainsPsi}
\psi \colon U_-^- \to  U_- \quad , \quad  \psi \colon U_-^+ \to \bRt\cap \{r\geq 2, z\geq 0\} ,
\end{equation}
as illustrated in Figure~\ref{fig:domainspsi333}.

 The   following is a direct consequence of   our labeling system: 
\begin{lemma}\label{lem-psi_action}
The action $\psi$   on the $\G_0$ curves is as follows:
\begin{enumerate}
\item The map $\psi:U_-^- \to U_-$, sends a subset of the curve
$\G_0(i_1,\ell_1 ; \cdots ; i_{n-1},\ell_{n-1} ;\ell_n)$ to a subset of
$\G_0(i_1,\ell_1 ; \cdots ; i_{n-1},\ell_{n-1} ;\ell_n+1)$, for $0 \leq \ell_n < \ell_n'$.
 For $n=1$,  this operation is always allowed. \\
 
 \item  The map $\psi:U_-^+ \to \bRt\cap \{z\geq 0\}$, sends a subset of the curve
$\G_0(i_1,\ell_1 ; \cdots ; i_{n-1},\ell_{n-1} ;\ell_n)$ to a subset
of itself.  \\

\item  The map $\psi:U_+ \to \bRt\cap \{z\geq 0\}$, sends a subset of the curve
$\G_0(i_1,\ell_1 ; \cdots ; i_{n-1},\ell_{n-1} ;\ell_n)$ to a subset of
$\G_0(i_1,\ell_1 ; \cdots ; i_{n-1},\ell_{n-1} ;\ell_n-1)$. For every $n \geq 1$, the index 
$\ell_n$   is bounded below, hence the map $\psi$ can only be applied to
such a curve a finite number of times.
\end{enumerate}
\end{lemma}
The action on the other families of curves is analogous, but it is slightly different on the base points. The difference is made explicit in Lemma~\ref{lem-actionpoints}.1(b).

Next, we analyze the action of the maps $\phi_i^+\colon U_{\phi_i^+}
\to \bRt$ for $i=1,2$. The domain $U_{\phi_i^+}$ consists of   points $x \in \bRt$
whose forward $\cK$-orbit intercepts the   face $E_i$ and then continues on to
intercept $\bRt$ again, hence the level is increased by one under
these maps. 
  Hypothesis~(K3) on the embeddings $\sigma_i$  implies   that their domains $U_{\phi_i^+} \subset \bRt$ satisfy 
\begin{equation}\label{eq-domainpsi+}
U_{\phi_1^+}   \subset \bRt \cap \{z<0\} \quad , \quad U_{\phi_2^+}  \subset \bRt \cap \{z>0\}
\end{equation}
and that  their ranges satisfy 
\begin{equation}\label{eq-rangepsi+}
\phi_1^+  ( U_{\phi_1^+}   )=V_{\phi_1^+} \subset \bRt \cap \{z<0\} \quad , \quad \phi_2^+  ( U_{\phi_2^+}    )=V_{\phi_2^+} \subset \bRt \cap \{z<0\} ,
\end{equation}
as represented in Figure~\ref{fig:GKdomains}.
 Then \eqref{eq-rangepsi+} and the definition of the  indexing system yields:

\begin{lemma}\label{lem-phi+_action}
For $i=1,2$, the map  $\phi_i^+$ sends a subset of $\G_0(i_1,\ell_1 ; \cdots ;
i_{n-1},\ell_{n-1} ;\ell_n)$ for $\ell_n\geq 0$ 
to a subset of $\G_0(i_1,\ell_1 ; \cdots ; i_{n-1},\ell_{n-1} ;i,\ell_n ; a)$
where $a\leq 0$ is the first index such that the curve $\G_0(a)$ exists. 
\end{lemma}
Note that $\phi_i^+$ sends any point in $U_{\phi_i^+}\cap \{r\geq 2\}$
to a point that is contained in the closed region bounded by the curve
$\G_0(a)$ if $i=1$ and by $\Lambda_0(a)$ if $i=2$.

\begin{figure}[!htbp]
\centering
{\includegraphics[width=110mm]{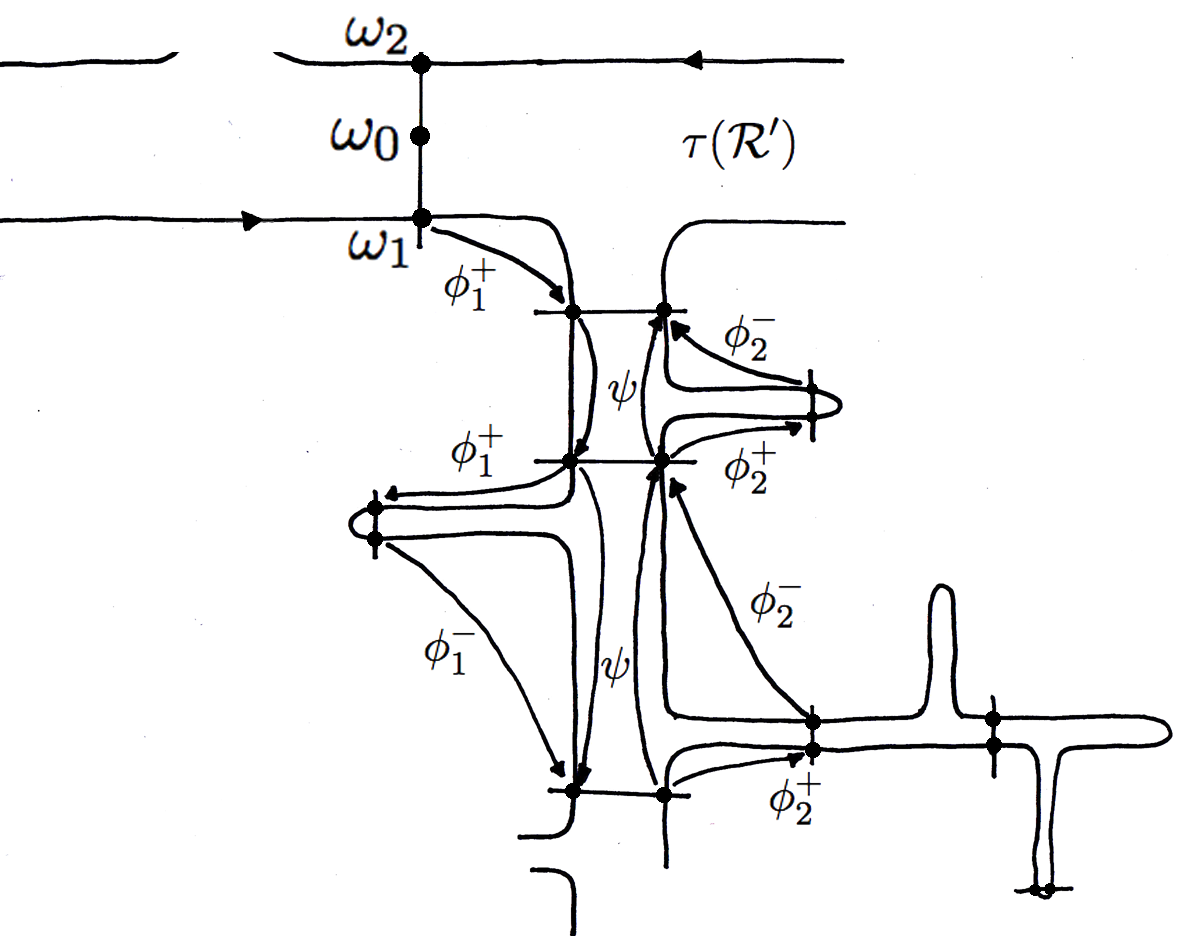}}
\caption[{Intersection of flattened $\fM_0$    with $\bRt$ with the action of $\cG_K^*$}]{\label{fig:generatorsM0R0} Intersection of flattened $\fM_0$    with $\bRt$ with the action of the generators of $\cG_K^*$}
\vspace{-6pt}
\end{figure}

  We comment on Figure~\ref{fig:generatorsM0R0},
  it is useful to compare it to Figure~\ref{fig:choufleur}. The
  horizontal upper band represents the Reeb cylinder $\tau(\cR')$. The
  transverse lines to the cylinder and propellers represent their
  intersections with the rectangle $\bRt$. Since only $\tau(P_\g')$
  and its ramifications are illustrated, these lines correspond to
  $\g_0$-curves. 

The generator $\psi$ preserves level, hence it moves points along
propellers. The points in $\tau(\cR')\cap \bRt$ are
$\psi$-invariant. For points in the left side boundary of
$\tau(P_\g')$, the action of $\psi$ pushes them downwards to the
following $\g_0$-curve, while for points in the right side boundary
$\tau(P_\g')$ it moves points upwards. The maps $\phi_i^+$ add 1 to
the level, so they correspond to jumping to the first curve at higher
level. For example, for the curve $\g_0(0)$, $\phi_1^+$ moves points
to $\g_0(1,0;0)$ and $\phi_2^+$ moves points to $\g_0(2,0;0)$. Observe
that the meaning of ``first'' curve is in the direction of the flow,
that runs in opposite directions near the boundaries of
$\tau(P_\g')$. The actions of $\phi_i^-$ are also illustrated.

\begin{figure}[!htbp]
\centering
\begin{subfigure}[c]{0.45\textwidth}{\includegraphics[width=60mm]{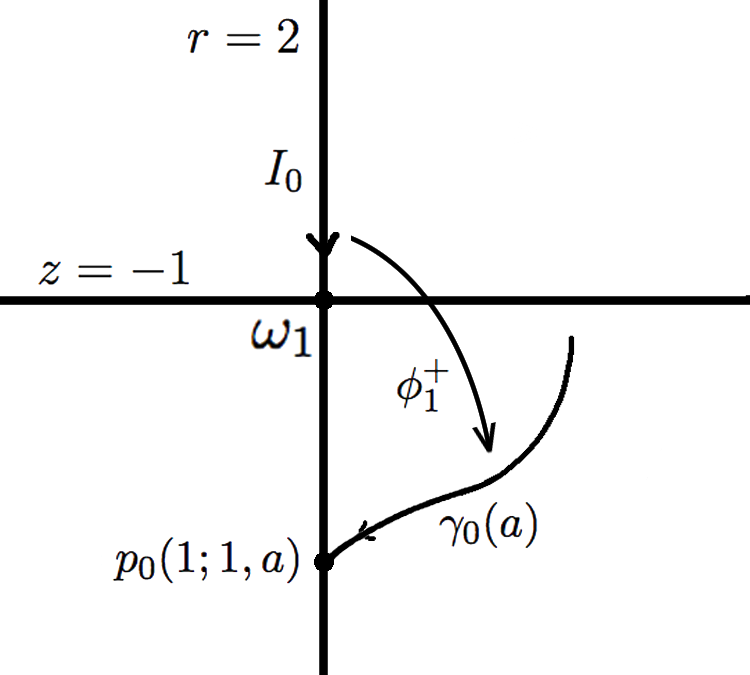}}\end{subfigure}
\begin{subfigure}[c]{0.45\textwidth}{\includegraphics[width=60mm]{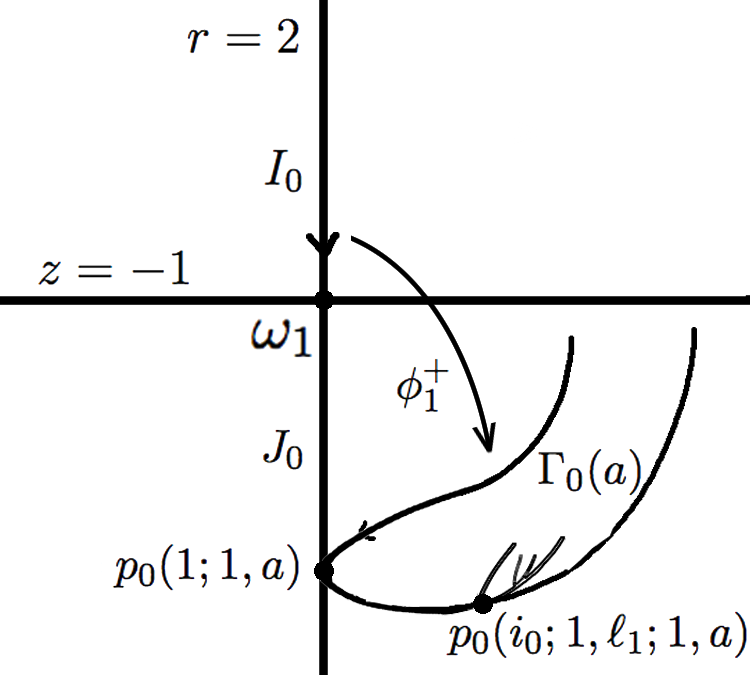}}\end{subfigure}
\caption{\label{fig:phi1} The image under $\phi_1^+$ of $I_0$ and   $N_0 = I_0 \cup J_0$}
\vspace{-10pt}
 \end{figure}

 The action of a general element $\vp \in \cGK^*$ on a curve in one of the families  
  $\ds \{\G_0 , \gamma_0 , \kappa_0 , \Lambda_0 , \lambda_0 ,  \chi_0\}$, 
  can be quite complicated and will be considered further in Section~\ref{sec-entropylamination}.
 For now, consider the simplest case of iterations of the map $\phi_1^+$. 
The image $\phi_1^+(N_0 \cap U_{\phi_1^+}) \subset \bRt$  is    a parabolic curve, twisting upwards as illustrated on the right side of  Figure~\ref{fig:phi1}. 
If the image curve $\phi_1^+(N_0 \cap U_{\phi_1^+})$ is again
(partially) contained in the domain of $\phi_1^+$, then  we can
consider its (restricted) image under $\phi_1^+$ and repeat this
process inductively as long as it is defined. This yields    a family
of parabolic curves, as illustrated on the right side of
Figure~\ref{fig:phi1}, which shows    the image  $\phi_1^+(I_0 \cap
U_{\phi_1^+})$ on the left-hand-side   and three iterations of this
map applied to $N_0$ on the right-hand-side.

The  iterates of $\psi$ map $\G_0(i_1,\ell_1 ; \cdots ; i_{n-1},
  \ell_{n-1} ; a)$ to any $\G_0(i_1,\ell_1 ; \cdots ; i_{n-1},
  \ell_{n-1} ; \ell_n)$. That is,  the action of $\psi$
  generates the family $\G_0(i_1,\ell_1 ; \cdots ; i_{n-1},  \ell_{n-1} ; \ell_n)$ from $\G_0(i_1,\ell_1 ; \cdots ; i_{n-1},  \ell_{n-1} ; a)$.
 
 \medskip

Finally, we analyze the action of the maps   $\phi_i^-\colon U_{\phi_i^-} \to \bRt$ for $i=1,2$, which are defined using the   $\Phi_t$-flow.
The domain $U_{\phi_i^-}$ consists of   points $x \in \bRt$
whose forward $\cK$-orbit intercepts the   face $S_i$ and then continues on to
intercept $\bRt$ again.  Hence,  the level is decreased by $1$ under
these maps,  and $U_{\phi_i^-}\subset \bRt\cap\{z>0\}$ for $i=1,2$.

\begin{lemma}\label{lem-phi-_action}
For $i=1,2$ and $n \geq 2$,  $(\phi_i^-)^{-1}$ maps a subset of  $\G_0(i_1,\ell_1 ; \cdots ; i_{n-1},
  \ell_{n-1} ; \ell_n)$ to a subset of $\G_0(i_1,\ell_1 ; \cdots  ; i_{n-1},
  \ell_{n-1} ; i, \ell_n-1; a)$, where $a$ is defined as in Lemma~\ref{lem-phi+_action}. 
  Hence, $\phi_i^-$ is only defined on the curves $\G_0(i_1,\ell_1 ; \cdots  ; i,
  \ell_{n-1} ; a)$.
\end{lemma}

We now comment on the behavior of the points  $p_0(i_0;
i_1,\ell_1;\cdots ; i_n,\ell_n)$. The points  $p'(i_0;1,\ell)$ are in
$\cL_1^-$ for $\ell \geq 0$, with $r(p'(i_0;1,\ell))=2$ and
$z(p'(i_0;1,\ell))<-1$. Thus,  $p(i_0;1,\ell)\in
\kappa\subset \G\subset L_1^-$ and $p_0(i_0;1,\ell)\in J_0$. Since
$p_0(i_0;1,\ell)\to \omega_1$ as $\ell\to\infty$, for $\ell$ large
enough, the points belong to $U_{\phi_1^+}$.  

If $p_0(i_0;1,\ell_1)\in U_{\phi_1^+}$, Lemma~\ref{lem-phi+_action}
implies that $\phi_1^+(p_0(i_0;1,\ell_1))=p_0(i_0;1,\ell_1;1,a)$
that belongs to the curve $\kappa_0(a)\cap \{z\leq 0\}$. Thus the
points $p_0(i_0;1,\ell_1; 1,a)$ accumulate on
$p_0(1;1,a)=\phi_1^+(\omega_1)$ as $\ell_1\to\infty$. Observe that
$p_0(1;1,a)$ is the vertex point of the curve $\kappa_0(a)$.

Analogously,   the points $p_0(i_0;1,\ell_1;2,a)$ belong to the curve
$\kappa_0(a)$ and  accumulate on $p_0(1;2,a)$, the 
other vertex point of $\kappa_0(a)$, as $\ell_1\to \infty$.
The points $p_0(i_0;2,\ell_1;i_2,a)$ belong to the curve
$\chi_0(a)$ and  accumulate on $p_0(2;i_2,a)$ as $\ell_1\to \infty$.

Lemmas~\ref{lem-psi_action} and \ref{lem-phi-_action} are interpreted
in terms of points in the following way.

\begin{lemma}\label{lem-actionpoints}
The maps $\psi$, $\phi_1^+$ and $\phi_2^+$ act on $p_0$ points as
follow:
\begin{enumerate}
\item For the map $\psi$ we have the three following possibilities:
\begin{enumerate}
\item For $1\leq \ell_n < \ell_n'$,  $\psi:U_-^- \to U_-$  maps  
$p_0(i_0;i_1,\ell_1 ; \cdots ; 1,\ell_n)$ to  
$p_0(i_0,i_1,\ell_1 ; \cdots ; 1,\ell_n+1)$.
 For $n=1$,  this operation is always allowed. \\
 
 \item   $\psi:U_-^+ \to \bRt\cap \{z\geq 0\}$ maps   
$p_0(i_0;i_1,\ell_1 ; \cdots ; 1,\ell_n)$ to   $p_0(i_0;i_1,\ell_1 ; \cdots ; 2,\ell_n)$.  \\

\item    $\psi:U_+ \to \bRt\cap \{z\geq 0\}$ maps  
$p_0(i_0;i_1,\ell_1 ; \cdots ; 2,\ell_n)$ to $p_0(i_0;i_1,\ell_1 ; \cdots ; 2,\ell_n-1)$. For every $n \geq 1$, the index 
$\ell_n$   is bounded below by $a$, hence the map $\psi$ can only be applied to
such a curve a finite number of times.

\end{enumerate}

\item  $\phi_1^+$ maps  a point $p_0(i_0;i_1,\ell_1 ; \cdots ;
  1,\ell_n)$ in its domain to   $p_0(i_0;i_1,\ell_1 ; \cdots ;
  1,\ell_n;1,a)$.

\item  $\phi_2^+$ maps a point $p_0(i_0;i_1,\ell_1 ; \cdots ;
  2,\ell_n)$ in its domain to $p_0(i_0;i_1,\ell_1 ; \cdots ; 2,\ell_n;1,a)$.
\end{enumerate}
\end{lemma}

We have the following consequences of the above results.

\begin{lemma}\label{lem-paccumulationpoints}
Every point $p_0(i_0; i_1,\ell_1;\cdots ; i_n,\ell_n)$ is an
accumulation point of the set of $p_0$ points.
\end{lemma}

\proof The points $p_0(i_0;1,\ell_1;\cdots; i_n,\ell_n)$ lie    in the
  curve $\kappa_0(i_2,\ell_2;\cdots;i_{n-1},\ell_{n-1};\ell_n)$, and
 as $\ell_1\to \infty$ they accumulate   on the point  $p_0(1;i_2,\ell_2;\cdots;i_{n-1},\ell_{n-1};i_n,\ell_n)$.
  
Also,  the points $p_0(i_0;2,\ell_1;\cdots; i_n,\ell_n)$ lie in the
  curve $\chi_0(i_2,\ell_2;\cdots;i_{n-1},\ell_{n-1};\ell_n)$, and as $\ell_1\to \infty$ they 
  accumulate  on the point  $p_0(2;i_2,\ell_2;\cdots;i_{n-1},\ell_{n-1};i_n,\ell_n)$.
\endproof

  \bigskip
 
\section{Normal forms and tree structure of $\fM_0$}\label{sec-normal}

A major theme of this work is to use the pseudogroup $\cGK$    generated by the  map    $\whPhi$ acting  on    $X = \bRt$, as defined in Definition~\ref{def-pseudogroupK}, to study the dynamics of the flow $\Phi_t$. The results of Section~\ref{sec-doublepropellers} interpret the local action of $\cGK$ in terms of the propellers introduced in Section~\ref{sec-proplevels}. The application of these two concepts associated to $\Phi_t$ will be made throughout the remaining sections, and many of the results we establish often rely on highly technical arguments, with various nuances. In this section, we pause to give a broader overview of how these concepts, of pseudogroups and propellers, are   related. To this end, we introduce a ``tree structure'' on the space $\fM_0$ as defined in \eqref{eq-lamination}  and illustrated in Figure~\ref{fig:choufleur}. Then in Proposition~\ref{prop-standardform} we give   an algebraic ``normal form''   for elements of the $\psg$ $\cGK^*$. These two results are   closely related, as we explain below. Finally, in Proposition~\ref{prop-asymptotic}  we   estimate   the number of normal forms as a function of the word length, and show that this   function has subexponential growth. This is an intrinsic property of the flow $\Phi_t$.

\begin{figure}[!htbp]
\centering
{\includegraphics[width=140mm]{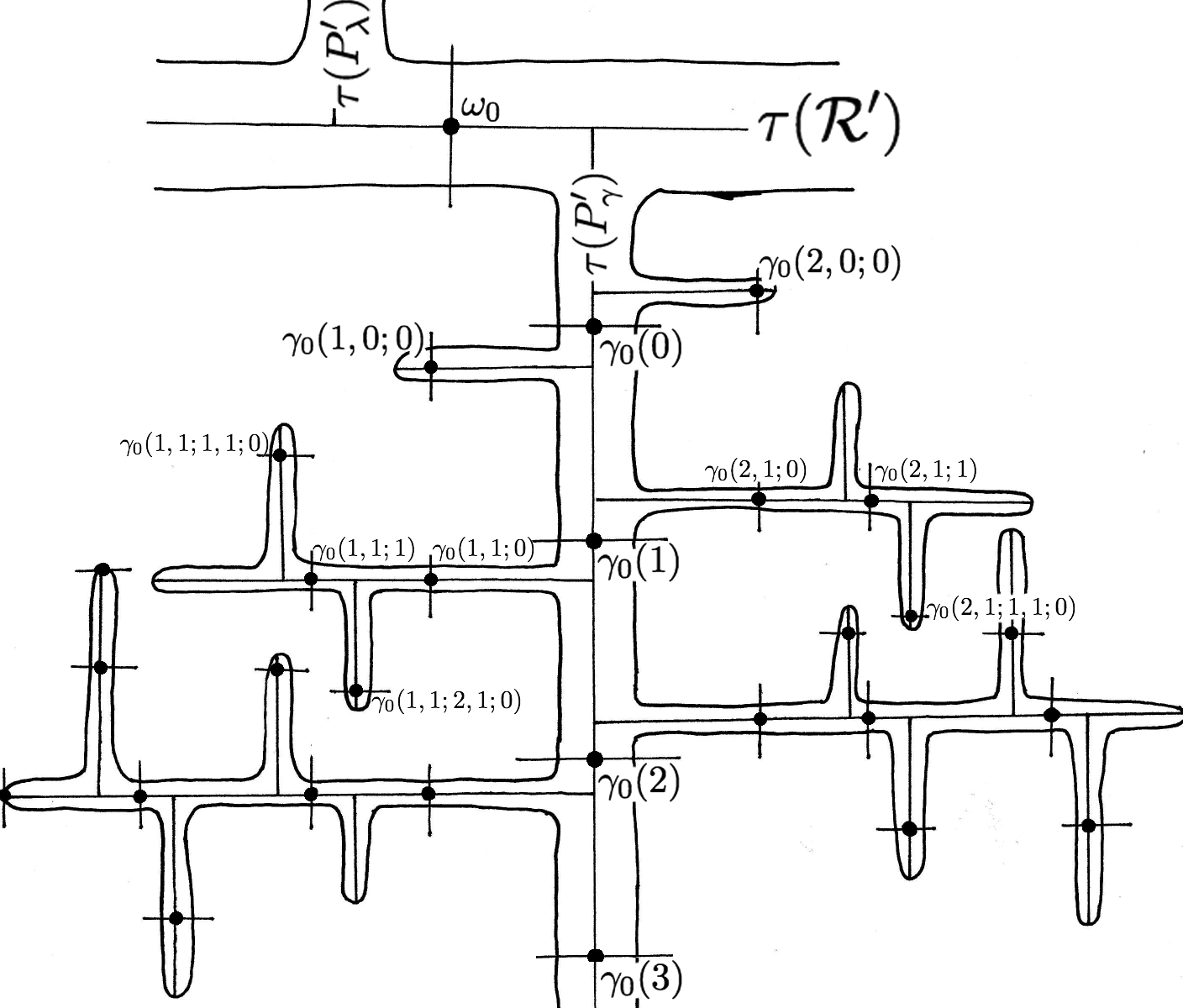}}
\caption{\label{fig:fM0cT}   $\fM_0$ with the  embedded tree $\TP$}
\vspace{-6pt}
\end{figure}

 We first construct the tree $\TP\subset \fM_0$.   Recall that the level decomposition of $\fM_0$   in \eqref{eq-nestedfamilies} expresses the space as a union of propellers, and each propeller is defined by the $\Psi_t$-flow of a curve segment in the face $\partial_h^- \mW'$.   Each propeller then intersects the center annulus 
 $\{z=0\} = \tau(\cA)$. 
   There is a nuance in this picture, as each propeller in the decomposition of $\fM_0$ may give rise to  ``bubbles'' which correspond to  double propellers of uniformly bounded complexity attached on the interior of the given propeller, as in \eqref{eq-Ssurface}. This  will be discussed in Section~\ref{sec-bubbles} and is illustrated in Figure~\ref{fig:choufleurbubbles}. 
   As they do not influence the dynamics of the flow $\Phi_t$,  they do not appear in  Figure~\ref{fig:fM0cT}.

 In the surface $\fM_0$, the intersection of each propeller with the annulus  $\{z=0\}$ forms an embedded   ``center line segment''.  Define $\TP' = \cA \cap \fM_0$ which consists of the union of these embedded line segments. Let $\TP'' \subset \TP'$ consist of the line segments in simple propellers.  
   The connected tree $\TP$ is formed by
adding continuous curve segments in $\fM_0$ joining the line segments in $\TP''$, as 
 is  illustrated in  Figure~\ref{fig:fM0cT}.    
 
 Define the center line   $\cT=\{z=0\} \cap  \bR_0$,   and set     $\fC_0' = \cT \cap \fM_0$ and $\fC' = \fM \cap \cT$. We will show in Proposition~\ref{prop-cantor}  
   that $\fC'$ is    a Cantor set, and observe there that  it has a   decomposition into two disjoint Cantor sets, $\fC' = \fC \cup \fC^1$. The points in the dense subset $\fC_0^1 = \fM_0 \cap \fC_0^1$ of $\fC^1$ correspond to the intersections of the double propellers in bubbles with $\cT$, while  the points in the dense subset $\fC_0 = \fM_0 \cap \fC$   of $\fC$ correspond to the intersections of $\g_0$ or $\lambda_0$ curves with $\cT$.  
The   ``fat dots''  in Figure~\ref{fig:fM0cT} correspond to the  points contained  in the set   $\fC_0$, and the lines through the dots represent the intersection with $\bRt$,

The   basepoint $\omega_0=\cT\cap\tau(\cR')$ is contained in the core
annulus $\cR = \tau(\cR') \subset \fM_0$, illustrated  in the picture as the top horizontal strip, which lies at level  $0$. There
is one edge containing $\omega_0$, the top horizontal line in Figure~\ref{fig:fM0cT}, which
corresponds to $\tau(\cR' \cap \cA)$ and defines  a loop containing $\omega_0$.
Thus, strictly speaking, $\TP$ is a tree provided the  loop containing $\omega_0$ is erased.

  Recall that the $\psg$   $\cGK^*$ denotes the  collection of all maps  formed by compositions of the maps $\{Id, \phi_1^+, \phi_1^-, \phi_2^+, \phi_2^-,  \psi\}$ and their restrictions to open subsets in their domains. The action of these generators on the endpoints of the $\gamma_0$ and $\lambda_0$ curves in $\bRt$ formed by the intersection $\fM_0 \cap \bRt$ were discussed in Section~\ref{sec-doublepropellers}. Note that each such curve intersects the line $\cT$ in a point of $\fC_0$, and the   actions of the generators $\{\psi, \phi_1^+ , \phi_2^+\}$ and their inverses have simple interpretations as actions on the vertices of $\TP$.  We introduce the notation $\{\opsi, \ophi_1, \ophi_2\}$ for the induced action on vertices. 
Later, in  Section~\ref{sec-entropylamination}, we consider the $\psg$ $\cGM$ generated by these actions in   detail. 

   First,  the action of   $\opsi$ fixes the basepoint $\omega_0$, and can be thought of as flowing in a loop around the   horizontal line in the core annulus $\cR'$.
The     action   of   $\opsi$ on points in $\fC_0$,  corresponding  to $\gamma_0$   curves in   $\bRt$,  is  determined by the action of $\psi \colon U_-^- \to U_-$ as described in  Lemma~\ref{lem-psi_action}.   We have: 
\begin{equation}\label{eq-defpsibar2}
\opsi ~ \colon ~  \gamma_0(i_1,\ell_1 ;\cdots ; i_{n-1},\ell_{n-1} ;\ell_n) \cap \cT ~ \mapsto ~  \gamma_0(i_1,\ell_1 ; \cdots ; i_{n-1},\ell_{n-1} ; \ell_n+1) \cap \cT .
\end{equation} 
  Since $\opsi$ does not change the level, we conclude that $\opsi$ moves a point in
$\fC_0$ to the ``next point'' along the   line segment in $\TP''$ containing it.   Here, ``next'' means
moving in the positive $\theta$ direction along the propeller, hence
making one turn around the circle $\cS = \{r=2 ~ \& ~  z=0\}$ in $\mK$.
Thus for $\xi \in \fC_0$ in the domain of the mapping $\opsi$, we have $r(\opsi(\xi)) < r(\xi)$ if $\xi\neq \omega_0$.   
The  action of the map $\opsi$   on the  points $\xi \in \fC_0$  defined by $\lambda_0$ curves is analogous.

We next show that  the paths connecting different level points, which were added to the segments in $\TP''$ to form $\TP$,  correspond to
  the actions of the elements $\ophi_k$  for $k =1,2$, analogous to the action of $\phi_k^+$. 
 Then by Lemma~\ref{lem-phi+_action} we have that:   
\begin{equation}\label{eq-defphik2}
\ophi_k ~ \colon ~  \gamma_0(i_1,\ell_1 ; \cdots ; i_{n-1},\ell_{n-1} ; \ell_n) \cap \cT ~ \mapsto ~  \gamma_0(i_1,\ell_1 ; \cdots ; i_{n-1},\ell_{n-1} ;k,\ell_n ; a) \cap \cT
\end{equation}
where   $\ell_n \geq 0$ as the curve must intersect the surface $E_k$ in a boundary notch, and $a\leq 0$ is the first index such that the corresponding curve in $\bRt$ exists, as explained in 
Remark~\ref{rmk-notation}. 
A similar result holds for the action of $\ophi_k$ on $\lambda_0$ curves. 
Note that  both    maps $\ophi_1$ and $\ophi_2$ increase the level function by $1$.

\begin{figure}[!htbp]
\centering
{\includegraphics[width=60mm]{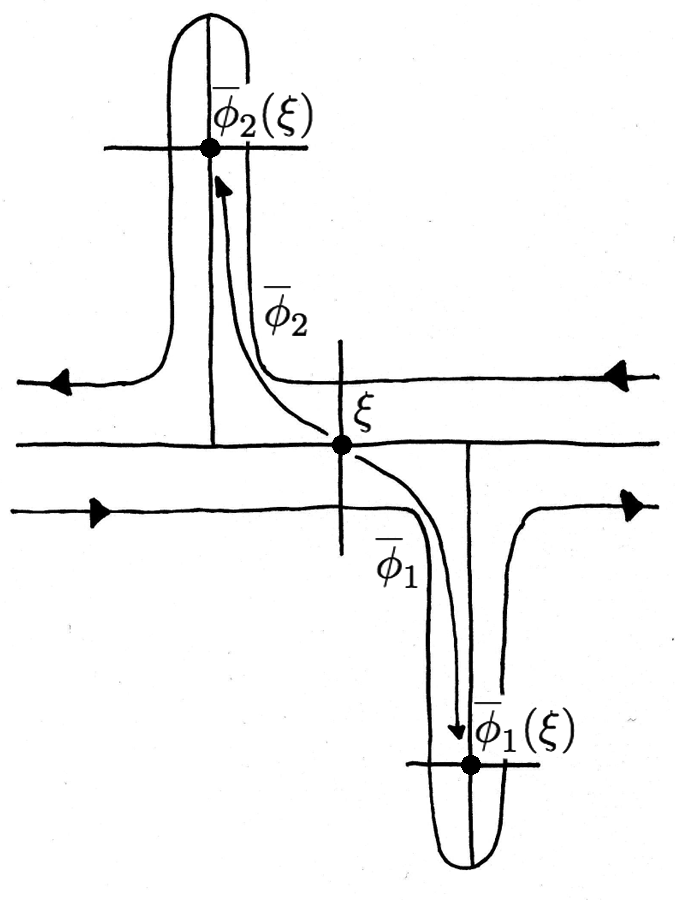}}
\caption{\label{fig:fM0cTa}   The action of $\ophi_k$}
\vspace{-6pt}
\end{figure}

The map $\ophi_1$ acts on $\omega_0$ via the insertion $\phi_1^+$ acting on the point $\omega_1 = (2, \pi, -1)$ to yield    $p_0(1 ; 1 , a) = \phi_1^+(\omega_1)$ which corresponds to  the first vertex $\ophi_1(\omega_0)$ at level $1$. Graphically, the  action of $\ophi_1$  is to ``turn the corner'' to the right. This results in adding a curved segment to $\TP''$ between $\omega_0$ and  the point $\ophi_1(\omega_0)$ in the downward propeller.

The map  $\ophi_2$  acts on $\omega_0$ via the insertion $\phi_2^+$ acting on the point $\omega_2 = (2, \pi, 1)$ to yield   $p_0(2 ; 1 , a) = \phi_2^+(\omega_2)$
which corresponds to the first vertex  $\ophi_2(\omega_0)$ at level $1$.
Graphically, the  action of $\ophi_2$  is to ``turn the corner'' to
the left. This results in adding a curved segment to $\TP''$ between
$\omega_0$ and  the point   $\ophi_2(\omega_0)$  in the upward propeller.

 Note that the remaining vertices in $\TP$ at level $1$ which descend
 the vertical main propeller $\tau(P_\g)$ in Figure~\ref{fig:fM0cT}
 correspond to the points $\opsi^{\ell} \circ \ophi_1(\omega_0)$.   The
 points $\opsi^{\ell} \circ \ophi_2(\omega_0)$ , also at level $1$,
 ascend on the propeller $\tau(P_\lambda)$ above the core annulus, and are not pictured in Figure~\ref{fig:fM0cT}.

The remaining vertices of the tree $\TP$ are at level at least $2$.
The actions of the maps $\{\opsi, \ophi_1, \ophi_2\}$ on these vertices follows the rules above as defined by \eqref{eq-defpsibar2} and \eqref{eq-defphik2}. The action of $\opsi$  is   translation along the line in $\TP''$ containing it, while the actions of the maps $\ophi_1$ and $\ophi_2$ turn the corner into the propeller at the next higher level, to the right or left accordingly. Accordingly, as mentioned in the definition of $\TP$ above, we  add a curved segment to $\TP''$ between the vertex $\xi$ and $\ophi_k(\xi)$,  whenever $\ophi_k(\xi)$  is defined.

We next consider the structure of the $\psg$ $\cGK^*$. The analog results for the pseudogroup generated by $\{\opsi, \ophi_1 , \ophi_2\}$ are studied in Section~\ref{sec-entropylamination}.
 The collection  
 \begin{equation}\label{eq-generators}
 \cGK^{(1)} = \{Id, (\psi)^{\pm 1}, (\phi^+_1)^{\pm 1}, (\phi^-_1)^{\pm 1}, (\phi^+_2)^{\pm 1},  (\phi^-_2)^{\pm 1}\}
\end{equation}
is a symmetric generating set   for the $\psg$ $\cGK^*$.
   Let  $\cGK^{(n)} \subset \cGK^*$ be  the collection of maps defined by   the restrictions of  compositions of at most $n$ elements of $\cGK^{(1)}$. 
The ``word metric''   on $\cGK^*$  is defined by setting    $\|\vp \|  \leq  n$ if  $\vp \in \cG_X^{(n)}$.

Recall that $\fM_0$ is invariant under the flow $\Phi_t$ and  so each $x \in \fM_0$ has infinite orbit. Set $\fMR = \fM \cap \bRt$.    We use the level function defined in Section~\ref{sec-fM} to show that the words for the restricted action $\cGK^* | \fMR$   have a normal form. First, we define:

\begin{defn}\label{def-monotone}
A word in $\cGK^*$ is said to be \emph{monotone (increasing)} if  it has the form
\begin{equation}\label{eq-productformula}
\vp ~ = ~\psi^{{\ell_m}} \circ \phi^+_{i_m}   \circ \cdots \circ \psi^{\ell_2} \circ \phi^+_{i_2} \circ \psi^{\ell_1} \circ \phi^+_{i_1} \circ \psi^{\ell_0} 
\end{equation}
where each $i_{\ell} = 1,2$, each $\ell_{k} \geq 0$. Let  $\cM(n)$ be the set of monotone words of length at most $n$, with $\cM(0) = \{Id\}$. Then let  $\cM(\infty)$ be the union of all collections $\cM(n)$ for $n \geq 0$.
\end{defn}
Observe that for $\vp \in \cM(\infty)$ as in \eqref{eq-productformula}, its length   satisfies  $\| \vp\| = m + \ell_0 + \ell_1 + \cdots + \ell_m$. 

Also, note that the collection of maps $\cM(\infty) \subset \cGK^*$  forms  a monoid, as composition of maps of the form \eqref{eq-productformula} is again of that form, assuming that the composition has non-empty domain.

\begin{remark}
Observe that if the level along the action of a word in $\cGK^*$ is monotone increasing, the word is written as a composition of the maps $\phi_i^+$, $(\phi_i^-)^{-1}$ and $\psi$, for $i=1,2$. By Lemma~\ref{lem-phi-_action}, the action of $(\phi_i^-)^{-1}$ can be replaced by $\phi^k\circ \phi_i^+\circ \psi^{-1}$ for some $k>0$. The word obtained by this substitution is still monotone increasing in level.

Even if the substitution made the expression of an element in $\cGK^*$ longer, the counting method below does not consider this nuance and the estimation on the number of {\it monotone words} is thus an estimation of words along which the level is monotone increasing.
\end{remark}

\begin{prop}\label{prop-standardform}
Let $\vp \in \cGK^*$ with $\|\vp \| \leq n$ and $Dom(\vp) \cap \fMR \ne \emptyset$. Then there exists a factorization $\vp = \vp^+ \circ \vp^-$,  where $\vp^+ \in \cM(n')$   and $(\vp^-)^{-1} \in \cM(n'')$ for integers $n', n''$ with  $n' + n'' \leq n$. Moreover, we have     $Dom(\vp) \subset Dom(\vp^+ \circ \vp^-)$. The factorization $\vp = \vp^+ \circ \vp^-$ is said to be the \emph{normal form} for the word $\vp$.
\end{prop}
\proof
By hypothesis we have
that $\vp=\vp_k\circ\cdots\circ\vp_1$ for some
$k\leq n$ and $\vp_\ell\in\cGK^{(1)}$ for every $1\leq \ell\leq k$. 
 Since $\fM$ is the closure of $\fM_0$, and we assume that $Dom(\vp) \cap \fMR \ne \emptyset$, there exists a
point $\xi_0 \in Dom(\vp) \cap \fM_0 \cap \bRt$. 
Set $\xi_\ell=\vp_\ell\circ\cdots\circ \vp_1(\xi_0)$.

The idea behind the proof is best described in terms of the  tree $\TP$ illustrated in Figure~\ref{fig:fM0cT}.
Given $\xi_0\in \fM_0\cap \bRt$ there exist a unique point $x\in
\cR'\cap \bRt$ such that $\Phi_t(x)=\xi_0$ for some $t > 0$. The $\cK$-orbit segment $[x, \xi_0]_{\cK}$ is a path in $\fM_0$ so can be deformed into a path in the tree $\TP$ from $\omega_0$ to the closest endpoint of the $\gamma_0$ or $\lambda_0$ curve containing $\xi_0$. This tree path define a 
unique monotone word $\vp_{\xi_0}\in \cGK^*$ such that $\vp_{\xi_0}(x)=\xi_0$. 
Analogously, there exist a unique monotone word $\vp_{\xi_k}\in \cGK^*$ and a unique point $y\in
\cR'\cap \bRt$ such that $\vp_{\xi_k}(y)=\xi_k$. The expression we are
looking for is then the simplified word obtained from
$\vp_{\xi_k}\circ(\vp_{\xi_0})^{-1}$, that is   a product of
a monotone word with the inverse of another monotone word.

Let us start the proof. By Proposition~\ref{prop-levels} the level
function $n_0:\fM_0\to\mN$ is well defined. Let $n_*(\vp, \xi_0) =
\min\{ n_0(\xi_{\ell}) \mid 0 \leq \ell \leq k\}$ and  $n^*(\vp, \xi_0)
= \max\{n_0(\xi_{\ell}) \mid 0 \leq \ell \leq k\}$. We consider the
following cases:
\begin{enumerate}
\item $n_0(\xi_0)\leq n_0(\xi_k)$ and $n_0(\xi_0)=n_*(\vp,\xi_0)$.
\item $n_0(\xi_0)\leq n_0(\xi_k)$ and $n_0(\xi_0)>n_*(\vp,\xi_0)$.
\item $n_0(\xi_0)> n_0(\xi_k)$ and $n_0(\xi_k)=n_*(\vp,\xi_0)$.
\item $n_0(\xi_0)> n_0(\xi_k)$ and $n_0(\xi_k)>n_*(\vp,\xi_0)$.
\end{enumerate}

\begin{figure}[!htbp]
\centering
{\includegraphics[width=76mm]{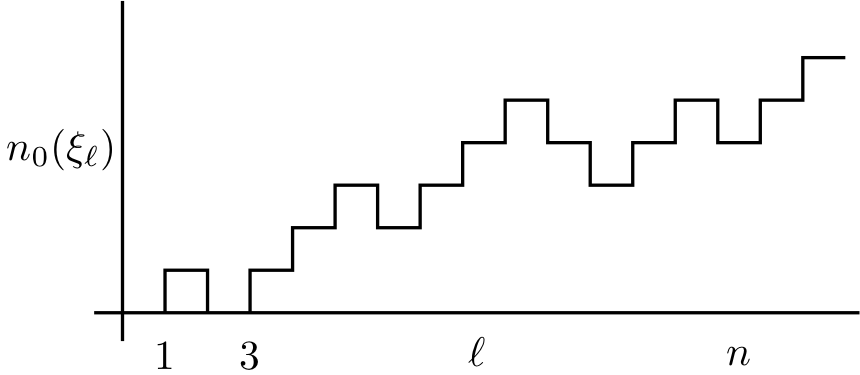}}
\caption{\label{fig:levelfunction} Plot of the level function $n_0(\xi_{\ell})$}
\vspace{-6pt}
\end{figure}

Consider  case (1), for which the plot of the function appears as in
Figure~\ref{fig:levelfunction}. 
Let $\ell_1\geq 0$ be the least index
such that there exists $\ell_2>\ell_1$ with
$n_0(\xi_{\ell_1})=n_0(\xi_{\ell_2})$ and
$n_0(\xi_\ell)>n_0(\xi_{\ell_1})$ for all $\ell_1<\ell<\ell_2$. By
Proposition~\ref{prop-shortcut4} combined with
Lemma~\ref{lem-actionpoints}, we can replace the subword
$\vp_{\ell_2}\circ\vp_{\ell_2-1}\circ \cdots \circ \vp_{\ell_1}$ by the
map $\psi$, thus changing a non-monotone subword of the expression of
$\vp$ by $\psi$. The new word has length less than $n$ and, since the
domain of $\psi$ contains the domain of the other generators in
$\cGK^{(1)}$, its domain contains the domain of $\vp$. 
Repeating this process a finite number of times we obtain a monotone
word $\vp^+$, and putting $\vp^-=Id$ we obtain the conclusion of Proposition~\ref{prop-standardform}.

The subword substitution operations in the proof of case (1) above have alternate   interpretations.  In terms of the sample graph of the level function in Figure~\ref{fig:levelfunction}, we are replacing any   parts of the graph which are not increasing   with a horizontal line, which corresponds to substituting in some power of the generator $\psi^{\ell}$ for the subword. In terms of the tree  $\TP$   in Figure~\ref{fig:fM0cT}, we are following a path that takes a short-cut across the base  of each   propeller, so that the resulting path follows a curve whose distance from the root is monotonically increasing. 
 
For  case (2), let $\ell_*$ be the first index for which
$n_0(\xi_{\ell_*})=n_*(\vp,\xi_0)$. We can apply the method of the first
case from $\xi_0$ to $\xi_{\ell_*}$, to change $\vp$ for a word of the
form $\vp_k\circ\cdots\circ \vp_{\ell_*+1}\circ\vp_*$ with
$(\vp_*)^{-1}$ monotone. The word $\vp_k\circ\cdots\circ
\vp_{\ell_*+1}$ satisfies the hypothesis of the first case and hence
can be expressed as a monotone word. This concludes the second case.

 In terms of the tree  $\TP$, the point $\xi_{\ell_*}$ corresponds to a vertex point in the tree which is the closest to the root vertex in the path in $\fM_0$ between $\xi_0$ and $\xi_k$.

We are left with the cases where $n_0(\xi_0)>n_0(\xi_k)$. In  case (3), let $\ell_*\leq k$ be the first index for which
$n_0(\xi_{\ell_*})=n_*(\vp,\xi_0)$. As in case (2), we can change
the expression of $\vp$ for $\vp_k\circ\cdots\circ \vp_{\ell_*+1}\circ\vp_*$ with
$(\vp_*)^{-1}$ monotone. Since $n_*(\vp,\xi_0)=n_0(\xi_k)$,
Proposition~\ref{prop-shortcut4} implies that the word
$\vp_k\circ\cdots\circ\vp_{\ell_*+1}$ can be replaced by a certain
number of consecutive $\psi$-maps, the number being less or equal to
$k-\ell_*$. Thus $\vp$ can be expressed as the inverse of a monotone
word and the conclusion follows.

The proof for case (4) is analogous to that of case (2).
\endproof

 We next develop an estimate for the  cardinality of the set $\cM(n)$ as a function of $n$. This will be used in later sections when we consider various entropy invariants of the flow $\Phi_t$. 
 
For  $\vp \in \cM(n)$ with $\| \vp \| = n$, assume that $\vp$ has the normal  form \eqref{eq-productformula}. Let   $\xi_0 \in \fM_0  \cap Dom(\vp)$, then 
 let $\xi_{\ell}$ for $0 < \ell \leq n$ denote the images of $\xi_0$ under the partial products of the factors  in \eqref{eq-productformula}. 
   Note  that  $r(\xi_{\ell+1}) \geq r(\xi_{\ell})$, with equality if
   $\xi_{\ell+1} =  \psi(\xi_{\ell})$, and with  strict inequality if
   $\xi_{\ell+1} =  \phi^+_j(\xi_{\ell})$ and
   $\xi_\ell\neq \omega_j$ for $j=1,2$. 
   
Observe that in the case where $\xi_{\ell+1} =  \phi^+_j(\xi_{\ell})$, the point
 $\xi_{\ell}$ must lie in the domain of $\phi^+_j$. Since $r(\xi_{\ell}) \geq 2$, if $j=2$ and   $\xi_{\ell} = \psi^i(\xi_{\ell -i})$, then $i > 0$ must be  sufficiently large so that the orbit transverses bottom half of the  rectangle $\bRt$ in order to enter the domain of $\phi^+_2$. 
By definition \eqref{eq-Delta},   the integer $\Delta(r) = \lfloor \Theta(r)/2\pi \rfloor $ is the number of intersections for the Wilson flow to reach the line $\{z=0\} \cap \bRt$, where $\Delta(r) \to \infty$ as $ r \to 2$. Thus,  $i > \Delta(r)$ where $r = r(\xi_{\ell -i})$. 

\begin{lemma}\label{lem-exponents}
For each $b \geq 1$, there exists integers $N_{b} > 0$ and $L_b >0$ such that for $\vp \in \cM(n)$   of the normal form \eqref{eq-productformula}
with $\ell_0 = 0 $,  $\ell_1 \leq b$ and sub-index   $m$ equaling the number of insertion maps in $\vp$,  then $m \leq N_b$ and each $\ell_k \leq L_b$ for $1 \leq k \leq b$.
Moreover,  $N_b \to \infty$  as $b \to \infty$.  
\end{lemma}
\proof 
The assumption that no factor in \eqref{eq-productformula} is a map $\phi^-_{j}$ for $j=1,2$ implies that the action of $\vp$ on $\bRt$ can be described in terms of the action of   the generators $\{\psi, \phi_1^+, \phi_2^+\}$  on the families of nested ellipses in $\bRt$ 
given by the intersection of $\widehat{\fM}_0 \cap \bRt$. This is   described in Lemmas~\ref{lem-invariantcurves}, \ref{lem-psi_action}, \ref{lem-phi+_action} and \ref{lem-actionpoints}.
The assumption $\ell_0 = 0$ implies that $\vp$ begins with the action of $\phi_{i_1}^+$ and thus   the special point $\omega_{i_1} \in Dom(\vp)$.

Recall that the function $N(r)$ for $r > 2$ was defined   in   Lemma~\ref{lem-maxstep},  and  is the maximum increase in the level function along the orbit of $\cK$ starting at an entry point $\xi$ with $r(\xi) = r$.  As shown in  Section~\ref{sec-global}, the function     $N(r) \to \infty$  as $r \to 2$. 
We are given  $b \geq 1$, so set 
\begin{equation}\label{eq-rbdef}
\rho_b=\min\{r(\phi_{i_2}^+\circ \psi^{\ell_1}\circ \phi_{i_1}^+(\xi)) \mid \xi \in Dom(\vp) ~ \text{and} ~ 0\leq \ell_1\leq b\}>2 ~ .
\end{equation}
 Then $N(\rho_b)$ is the maximal increase in the level function  starting from the point $\ds \phi_{i_2}^+\circ \psi^{\ell_1}\circ \phi_{i_1}^+(\xi)$, and thus $m \leq N(\rho_b) +2$. Set $N_b = N(\rho_b) +2$. Note that $N_b \to \infty$  as $b \to \infty$. 

Let $T_b > 0$ be such that  for all   $x \in \partial_h^- \mW$ with $\rho_b \leq r(x) \leq 3$, then the Wilson flow $\Psi_t(x)$ exits $\partial_h^+ \mW$ with $t \leq T_b$. Then $L_b$ can be taken to be the greatest integer less than $T_b/4\pi$.
\endproof

We use Lemma~\ref{lem-exponents} to obtain an estimate on the growth of the function $ \# \cM(n)$.

 \begin{prop}\label{prop-asymptotic} 
   For each $b \geq 1$, there is a polynomial function $P_b(n)$  of $n$ such that 
the  cardinality of the set $\cM(n)$   satisfies
  \begin{equation}\label{eq-productestimate}
\#\cM(n) \leq P_b(n) \cdot  2^{(n/b)} .
\end{equation}
\end{prop}
 \proof
 For  $\vp \in \cM(n)$ with $\| \vp \| = n$, assume that $\vp$ has the normal  form \eqref{eq-productformula}.

   For the given value of $b$,  let $i(\vp, b) \geq 1$ be the index such that $\ell_i \geq b$ for all $1 \leq i < i(\vp, b)$, and $\ell_i < b$ for $i = i(\vp, b)$.   Factor $\vp = \vp^{(b)} \cdot \vp_{(b)}$ where $\vp^{(b)}$ starts with the map $\ds \phi_{i_{i(\vp, b)}}^+$ and $\vp_{(b)}$ starts with the map $\psi^{\ell_0}$.

 Let $k(\vp, b)$ denote the number of factors in  $\vp^{(b)}$ of the form $\phi_i^+$.  By Lemma~\ref{lem-exponents}, $0 \leq k(\vp, b) \leq N_b$. We have $\| \vp^{(b)}\| \leq n$, so the indices of the factors of the form $\phi_i^+$ appearing in $\vp^{(b)}$ gives a choice of $k(\vp, b)$ indices out of the maximum of $n$ possibilities, and for each choice of index, let $i = 1,2$. Thus, the number of such words has an upper bound ${n \choose k} \cdot 2^k$ for $k = k(\vp,b)$. Then set 
 \begin{equation}\label{eq-productformula1}
B'(n,b) = {n \choose 0} \cdot 2^0 + {n \choose 1} \cdot 2^1 + \cdots + {n \choose N_b} \cdot 2^{N_b}
\end{equation}
Observe that $B'(n,b)$ is a polynomial function of $n$ of degree at most $N_b$. It follows   that $B'(n,b)$   is an upper bound on the number of possible words  $\vp^{(b)}$ which can arise for $\vp \in \cM(n)$ and the given value of $b$.

 Next, consider the number of possible choices for the words $\vp_{(b)}$ which can arise. For $p = i(\vp, b)-1$, we then can write
 \begin{equation}\label{eq-productformula2}
\vp_{(b)} ~ = ~    \psi^{\ell_p} \circ \phi_{i_{p}}^+ \circ \psi^{\ell_{p-1}} \circ \cdots \circ \psi^{\ell_1} \circ \phi_{i_1}^+ \circ \psi^{\ell_0}
\end{equation}
where  each index $\ell_i > b$ for $1 \leq i \leq  p$, and so $p  \leq n/(b+1) < n/b$. 
Observe that  for $1 \leq k \leq p$,  there are at most  $2^p$ possible choices of $i_{k}$, so the number of such choices is bounded above by $2^{(n/b)}$.
 
The placement of the terms $\phi_{i_k}^+$, or equally the choices of the values $\ell_i \geq b$, is given by a more complicated choice function.
  Observe that $b \leq \ell_0 \leq n$, and $b \leq \ell_1 \leq n - \ell_0 -1\leq n-b-1$. Thus, there are at most $n-2b$ possible values for  $\ell_1$. 
 Next, we have  $  b \leq \ell_2 \leq n-\ell_0-\ell_1-2\leq n-2b-2$, so that there are  at most $n - 3b-1$ possible values for   $\ell_2$. We continue in this way up to the choice of $\ell_p$. The number of possible choices of the indices $(\ell_0 , \ell_1 , \ldots , \ell_p)$   is then bounded above by the products of the maximal number   of values  for each $\ell_i$ for $0 \leq i \leq p$,    so is a polynomial in $n$ of degree at most $p+1$, which we denote by $B''(n,b)$. Set 
 \begin{equation}\label{eq-productformula3}
P_b(n) =  B'(n,b) \cdot B''(n,b)  
\end{equation} 
The estimate   \eqref{eq-productestimate} follows.
    \endproof

  \begin{cor}\label{cor-subexponentialgrowth} 
The  function $n \mapsto  \# \cM(n)$  has subexponential growth. That is, we have   
 \begin{equation}\label{eq-productformula4}
\lim_{n\to \infty} ~ \frac{\ln (\# \cM(n))}{n} = 0 .
\end{equation} 
\end{cor} 
\proof
For each $b \geq 1$, the estimate \eqref{eq-productestimate} implies that the limit in \eqref{eq-productformula4} is bounded above by $\ln(2)/b$ hence equals $0$.
\endproof

  \bigskip

 \section{Internal notches and bubbles}\label{sec-bubbles}

In this section, we analyze the properties of the ``bubbles'' that arise when the interior of a propeller intersects an insertion region, resulting in an internal notch, as discussed in Section~\ref{sec-proplevels} and illustrated in Figure~\ref{fig:intnotches}. Recall that an internal notch in a propeller $P_{\g}$ is a rectangular hole whose
boundary is disjoint from the boundary of $P_{\g}$. 
The $\Phi_t$-flow of an internal notch
generates a compact surface, as defined by   (\ref{eq-Ssurface}) and called a ``bubble''.
A key result of this section is that  the bubbles obtained from internal notches  all admit a uniform bound on
their complexity: the difference of level between any two points in a
bubble is uniformly bounded. This is shown in Proposition~\ref{prop-bubbles1}.

Assume that the construction of the plug $\mK$ and the flow $\cK$
 satisfies Hypotheses~\ref{hyp-SRI} and \ref{hyp-genericW}. 

In Sections~\ref{sec-proplevels} and \ref{sec-doublepropellers}, the
orbits of the special points $p_1^{\pm}$ and $p_2^{\pm}$ were   used to label the
boundary notches of the propellers $P_{\g}'$, $P_{\lambda}'$, $P_\G'$ and
$P_\Lambda'$, that are bounded by $\cW$-arcs spaced along the edges of
the propellers. A propeller formed from any of these four families and at any
level $n \geq 1$, might also have internal notches as described in Section~\ref{sec-proplevels}.
 In this section, we give an analogous  labeling for the internal notches and bubbles of the
propellers $P_{\g}'$ and $P_{\lambda}'$, which    then also applies to the double propellers
$P_{\G}'$ and $P_{\Lambda}'$.

First, recall the data that is given. 
 The curve $\g' = \cR \cap \cL_1^-$ denotes the intersection of the face $\cL_1^{-}$ with the Reeb cylinder $\cR$, as in \eqref{eq-curvesprime}   and $\ovg' = \cR \cap \cL_1^+$ is its facing curve. Define  $\g = \sigma_1^{-1}(\g') \subset L_1^- \subset \partial_h^- \mW$ with facing curve $\ovg = \sigma_1^{-1}(\ovg') \subset L_1^+ \subset \partial_h^+ \mW$.  
  The curve $\kappa' \subset \cL_1^-$ is the intersection  of the $\Psi_t$-flow of the interval $J_0$ with $\cL_1^-$ and $\kappa = \sigma_1^{-1}(\kappa') \subset L_1^-$,
  as illustrated in  Figure~\ref{fig:Gamma1}.

The curve $\lambda' = \cR \cap \cL_2^-$,  denotes the intersection of the face $\cL_2^{-}$ with the Reeb cylinder $\cR$, as in \eqref{eq-curvesprime}   and $\ovl' = \cR \cap \cL_2^+$ is its facing curve.  
Define   $\lambda = \sigma_2^{-1}(\lambda') \subset L_2^- \subset \partial_h^- \mW$ with facing curve $\ovl = \sigma_2^{-1}(\ovl') \subset L_2^+ \subset \partial_h^+ \mW$.
 The curve $\chi' \subset \cL_2^-$ is the intersection  of the $\Psi_t$-flow of the interval $K_0$ with $\cL_2^-$ and $\chi = \sigma_1^{-1}(\chi') \subset L_2^-$, 
 again  as illustrated in  Figure~\ref{fig:Gamma1}.

Consider first the propeller $P_{\g}\subset \mW$ generated by the
$\Phi_t$-flow of $\g$,  and let $P_\g'=P_\g\cap \mW'$ be the notched
propeller. The
intersection of $P_\g'$ with $\bRt$ consists of  an infinite
family of arcs with lower endpoints in the vertical line segment $J_0$ that
accumulate on $\omega_1$ and upper endpoints that
accumulate on $\omega_2$, as in Figure~\ref{fig:arcspropeller}. This
implies, as observed in Section~\ref{sec-proplevels}, that $P_\g'$ has
an infinite number of (boundary) notches and a finite number $|b|\geq 0$ of
internal notches. 
We assume for the rest of this section that $b< 0$ and investigate the properties of the bubbles which result from their $\Phi_t$-flows. 

Let $\gamma'(i,\ell)\subset \cL_i^-$ for $b\leq \ell<0$ denote the
intersection of the boundaries of the internal notches of $P_\g'$ with
$\cL_i^-$. Observe that both endpoints of such a curve are contained in the boundary  $\partial
\cL_i^-$, as for the curve $\g'(1,-1)$ in
Figure~\ref{fig:intcurvesL1} and the the curve $\g'(2,-1)$ in
Figure~\ref{fig:intcurvesL2}. The $\Psi_t$-flow of each
$\gamma'(i,\ell)$, from $\cL_i^-$ to $\cL_i^+$, defines a rectangular region in the interior of $P_{\g}$. 

Set $\gamma(i,\ell)=\sigma_i^{-1}(\g'(i,\ell))\subset L_i^-$ and parametrize $\gamma(i,\ell)$ by $\gamma_{(i,\ell)}:[0,2]\to L_i^-$ in
such a way that:
\begin{itemize}
\item $2<r(\g_{(i,\ell)}(1))\leq r(\g_{(i,\ell)}(s))$ for every $s\in
  [0,2]$;
\item $r(\g_{(i,\ell)}(0))=r(\g_{(i,\ell)}(2))=3$,  so both endpoints
  lie in the boundary $\partial_h^-\mW\cap\partial_v\mW$.
\end{itemize}
Analogously define $\lambda(i,\ell)\subset L_i^-$ for $b\leq
\ell<0$. Since the $\g$ and $\lambda$ curves
are interlaced the number of internal notches in $P_\g'$ is equal to
$\pm 1$ the number of internal notches in $P_\lambda'$. To simplify the notation, we assume  that there is the same number of internal
notches in $P_\lambda'$ and $P_\g'$.

\begin{prop}\label{prop-r_b}
There exists $r_b > 2$ such that for any $b\leq \ell <0$ and $i=1,2$, then $r(x') \geq r_b$ for all $x' \in \g(i,\ell)$ and all 
$x' \in \lambda(i,\ell)$.
\end{prop}

\proof For a curve $\zeta \colon [0,1] \to \mW$, set $r(\zeta) = \min \, \{ r(\zeta(s)) \mid 0 \leq s \leq 1\}$. Then set
\begin{equation}\label{eq-rb2}
r_b=\min \{r(\kappa(1,0)), r(\kappa(2,0)), r(\chi(1,0)), r(\chi(2,0))\} >2. 
\end{equation}
 
\begin{figure}[!htbp]
\centering
{\includegraphics[width=60mm]{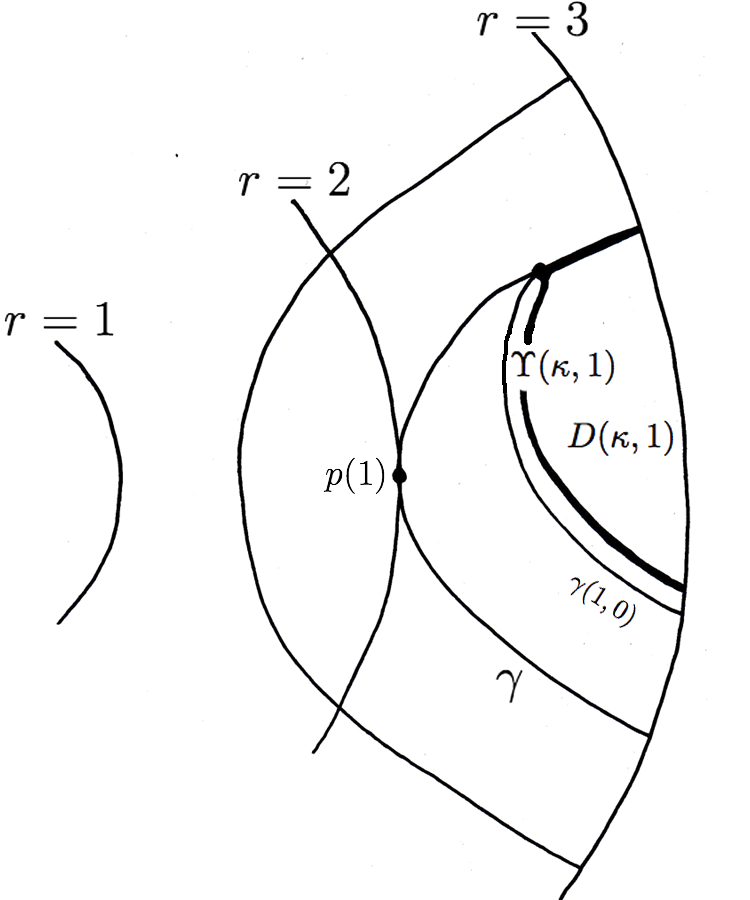}}
\caption{\label{fig:proofr_b}  The curve $\Upsilon(\g,1)$ and the disc  $D(\kappa,1)$ in $L_1^-$}
\vspace{-6pt}
\end{figure}

Consider first the case $\g(1,\ell)$ for $b\leq \ell <0$.  
Recall that  $p(1;1,0) \in \kappa$ is the endpoint of the curve $\gamma(1,0) \subset L_1^-$ and satisfies $r(p(1;1,0)) > 2$. 
Define a curve $\Upsilon(\kappa,1) \subset  L_1^-$ with endpoints in $\{r =3\}$ which first follows the path $\kappa(1,0)$ from its endpoint in $\{r =3\}$ to its endpoint $p(1;1,0)$, then follows the curve $\kappa$ back  to its boundary point in $\{r =3\}$. This is illustrated in  Figure~\ref{fig:proofr_b}. 
Note that  $\Upsilon(\kappa,1)$ divides $L_1^-$ into two topological discs and consider the closure $D(\kappa,1)$ of the disk  contained in the region $r> 2$. 
The radius function restricted to $D(\kappa,1)$ has a minimum value
greater or equal to $r_b$.

 Note that for $b \leq \ell < 0$, we have  $\g(1,\ell)  \subset D(\kappa,1)$ and hence  $r(\g(1,\ell)) \geq r_b$.

The radius estimate for the cases $\g(2,\ell)$, $\lambda(1,\ell)$ and $\lambda(2,\ell)$ follow in the similar manner, using the corresponding curves 
$\Upsilon(\kappa,2)$, $\Upsilon(\chi,1)$ or $\Upsilon(\chi, 2)$. 
\endproof

In what follows, we restrict our attention to the surfaces derived from the $\Psi_t$-flow of the curves $\g(1,\ell)$   for $b\leq \ell<0$, which form the ``bubbles'' in the propeller $P_{\g}$. 
An analogous  discussion applies to the curves $\g(2,\ell)$, $\lambda(1,\ell)$ and $\lambda(2,\ell)$. 

The endpoints of  $\g(1,\ell_1)$ are contained in $\partial_v \mW$ and thus their $\cW$-orbits escape through $\partial_v^+ \mW$ without intersecting $\bRt$.
Thus,  the intersection    $P_{\g(1,\ell_1)} \cap \bRt$ is   a finite family of circles contained
in the region $r\geq r(\g_{(1,\ell_1)}(1))>r_b$. 
The number of circles is bounded by $\Delta( r(\g_{(1,\ell_1)}(1)))+1$, where $\Delta$ is defined by \eqref{eq-Delta}.
The number of circles  thus  admits the uniform upper bound $\Delta(r_b)+1$.

Since $\g(1,\ell_1)$ is contained in the region of $L_1^-$ bounded by
the curve $\Upsilon(\kappa,1)$ introduced in the proof of
Proposition~\ref{prop-r_b}, thus is in the region of $L_1^-$ bounded by
$\Gamma$. The
circles formed by  the intersection   $P_{\g(1,\ell_1)} \cap \bRt$ lie in the discs bounded by a certain $\G_0(\ell_2)$, for $\ell_2\geq
a$. Assuming that $P_{\g(1,\ell_1)}\cap\bRt$ is non-empty, the first
curve in the intersection is in the region of $\bRt$ bounded by
$\G_0(a)$, with $a$ as in Remark~\ref{rmk-notation} and the following curves are obtained by applying the
Wilson map $\psi\in\cGK$. Thus if the intersection is formed by $n$
curves, the first one is in the region bounded by $\G_0(a)$, the
second one in the region bounded by $\G_0(a+1)$ and the $k$th curve is
in the region bounded by $\G_0(a+k)$, for $0\leq k\leq n$.

We extend the labeling to $P_{\g(1,\ell_1)}\cap\bRt$ by naming the
curves in the intersection
$\g_0(1,\ell_1;\ell_2)$, for $a\leq \ell_2\leq a+n$, in such a way that $\g_0(1,\ell_1;\ell_2)$ is
contained in the disc bounded by $\G_0(\ell_2)$. Observe that $\ell_1$
is negative and  that the connected component of
$P_{\g(1,\ell_1)}\cap\bRt$ closer to the boundary $r=3$ of $\bRt$ is
denoted by $\g_0(1,\ell_1;a)$. This convention
agrees with the labeling system introduced in the previous sections.
Analogously we have:
\begin{itemize}
\item The curves in the intersection $P_{\lambda(1,\ell_1)}\cap\bRt$
  are labeled $\lambda_0(1,\ell_1;\ell_2)$ and
  are contained in regions bounded by $\Gamma_0$ curves. Thus
  $\lambda_0(1,\ell_1;\ell_2)$ lies inside the region bounded by
  $\Gamma_0(\ell_2)$.
\item The curves in the intersection $P_{\g(2,\ell_1)}\cap\bRt$
  are labeled $\g_0(2,\ell_1;\ell_2)$ and
  are contained in regions bounded by $\Lambda_0$ curves. Thus
  $\g_0(2,\ell_1;\ell_2)$ lies inside the region bounded by
  $\Lambda_0(\ell_2)$.
\item The curves in the intersection $P_{\lambda(2,\ell_1)}\cap\bRt$
  are labeled $\lambda_0(2,\ell_1;\ell_2)$ and
  are contained in regions bounded by $\Lambda_0$ curves. Thus
  $\lambda_0(2,\ell_1;\ell_2)$ lies inside the region bounded by $\Lambda_0(\ell_2)$.
\end{itemize}

For the notched propellers $P_{\g(i, \ell_1)}'$ and
$P_{\lambda(i,\ell_1)}'$, for $i=1,2$ and $b\leq \ell_1<0$, obtained by
intersecting the corresponding propeller with $\mW'$, their images $\tau(P_{\g(i,\ell_1)}')$ and
$\tau(P_{\lambda(i,\ell_1)}')$ form part of $\fM_0$. 

We now describe the compact surfaces $S_{\g(1,\ell_1)}\subset \mK$ for
$b\leq \ell_1<0$
introduced in \eqref{eq-Ssurface}.  Consider the notched propellers
$P_{\g(1,\ell_1)}'= P_{\g(1,\ell_1)}\cap\mW'$ for $b\leq
\ell_1<0$. 
If
$P_{\g(1,\ell_1)}\cap \cL_1^-$ is non-empty, let
$\g'(1,\ell_1;1,\ell_2)$ be the curves in the intersection. Observe
that $\ell_2$ admits finitely many possible values, at most $\Delta(r_b)+1<\infty$.
In the same way if
$P_{\g(1,\ell_1)}\cap \cL_2^-$ is non-empty, let
$\g'(1,\ell_1;2,\ell_2)$ be the curves in the intersection. Observe
that $\ell_2$ admits finitely many possible values, at most   
$\Delta(r_b)+1<\infty$. Since $\g_0(1, \ell_1;1,b-1)\subset \bRt$ is the region bounded by
$\G_0(b-1)=\g_0(b-1)\cup\kappa_0(b-1)$, the $\cK$-orbits of points in
$\g_0(1, \ell_1;1,b-1)$ come back to $\bRt$ before hitting the
insertions. Hence the first curve in $P_{\g(1,\ell_1)}\cap\cL_1^-$ is either
$\g(1, \ell_1;1,b)$ or $\g(1, \ell_1;1,b+1)$.

\begin{figure}[!htbp]
\centering
{\includegraphics[width=80mm]{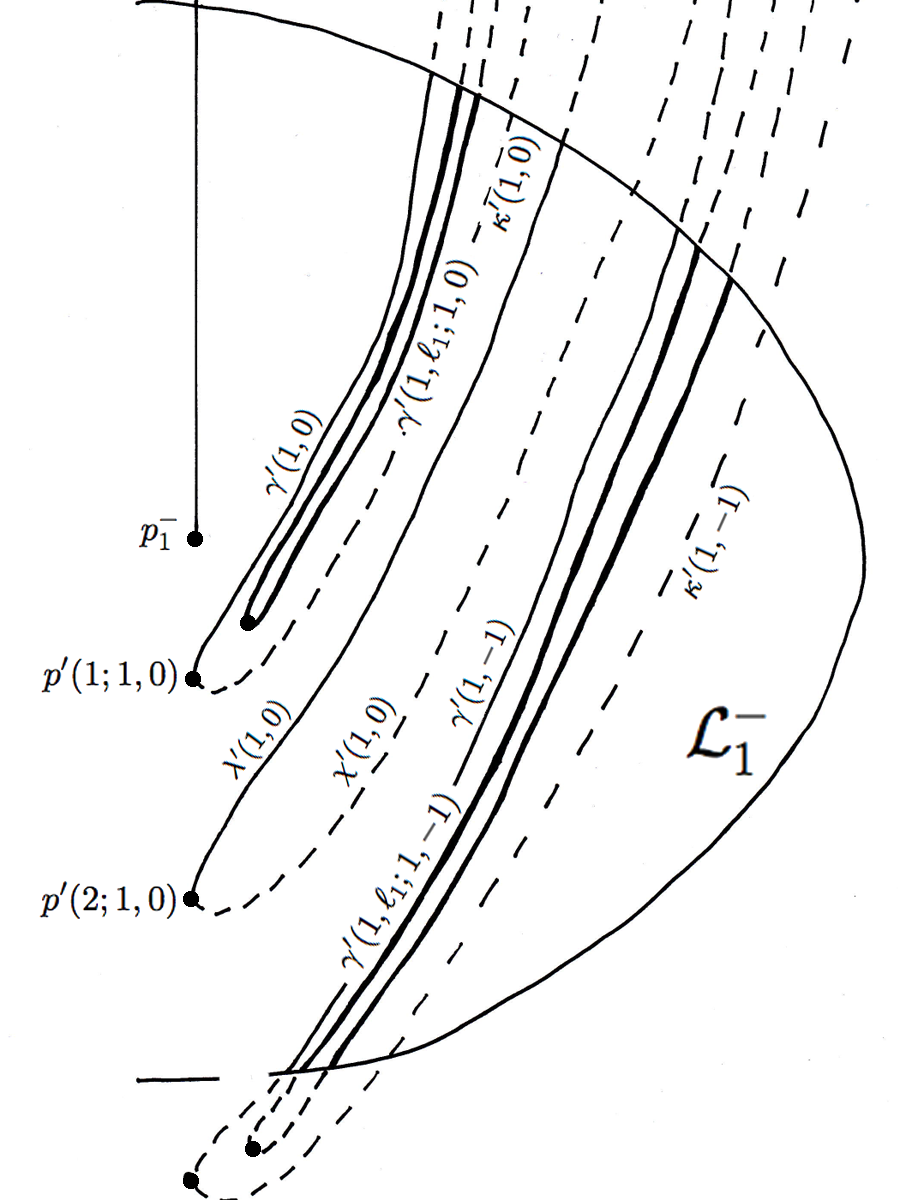}}
\caption{\label{fig:bubblesL-}  Possible intersections of $P_{\g(1,\ell_1)}$ with $\cL_1^-$ for $b\leq \ell_1<0$}
\vspace{-6pt}
\end{figure}

\begin{figure}[!htbp]
\centering
{\includegraphics[width=120mm]{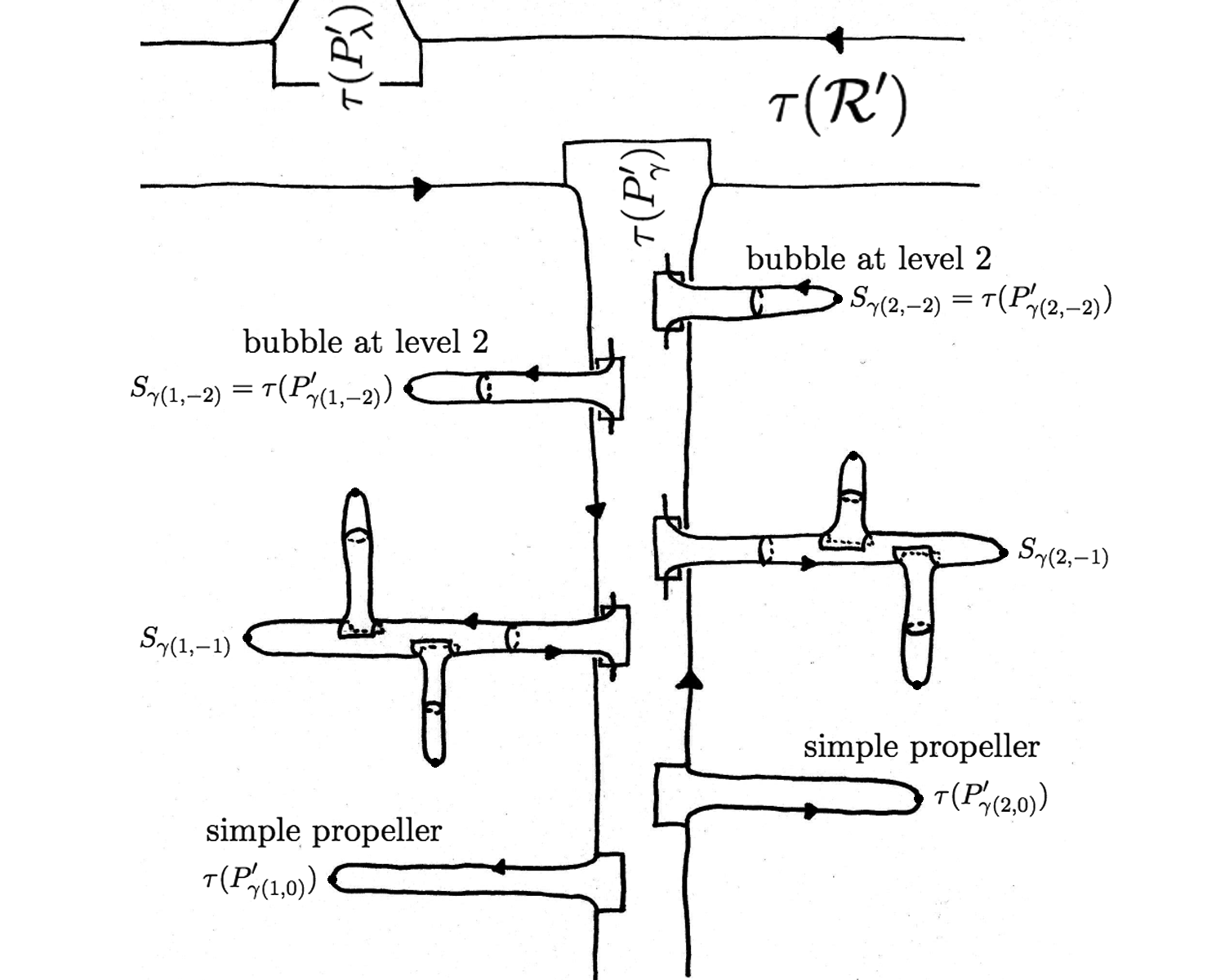}}
\caption{\label{fig:choufleurbubbles}  Flattened part of $\fM_0$ with  bubbles}
\vspace{-6pt}
\end{figure}

Consider now the curves
$\g(1,\ell_1;1,\ell_2)=\sigma_1^{-1}(\g'(1,\ell_1;1,\ell_2))\in
L_1^-$ with $b\leq \ell_1<0$ and $b\leq \ell_2$. Two possible situations arise, as illustrated in
Figure~\ref{fig:bubblesL-}:
\begin{itemize}
\item $\g'(1,\ell_1;1,\ell_2)$ contains a point of the $\cW$-orbit
  of $\g_{(1,\ell_1)}(1)$. In this case the curve is an arc with
  endpoints in $\partial \cL_1^-$. Thus $\g(1,\ell_1;1,\ell_2)$
  generates a finite double propeller
  $P_{\g(1,\ell_1;1,\ell_2)}\subset \mW$. The Radius Inequality
  implies that the minimum radius along $\g(1,\ell_1;1,\ell_2)$ is
  strictly bigger than $r(\g_{(1,\ell_1)}(1))>r_b$, and thus the number
  of circles in $P_{\g(1,\ell_1;1,\ell_2)}\cap\bRt$ is less or equal
  that the number of circles in the intersection
  $P_{\g(1,\ell_1)}\cap\bRt$.
\item $\g'(1,\ell_1;1,\ell_2)$ does not contains points in the
  $\cW$-orbit of $\g_{(1,\ell_1)}(1)$ and thus it is formed by the
  union of two connected components, each being an arc with endpoints
  in $\partial \cL_1^-$. In this case we obtain from
  $\g(1,\ell_1;1,\ell_2)$ two double finite propellers whose union we
  denote by $P_{\g(1,\ell_1;1,\ell_2)}\subset \mW$. Observe that the Radius Inequality
  implies that the minimum radius along any component of $\g(1,\ell_1;1,\ell_2)$ is
  strictly bigger than $r(\g_{(1,\ell_1)}(1))>r_b$, and thus the number
  of circles in $P_{\g(1,\ell_1;1,\ell_2)}\cap\bRt$ is less or equal
  that twice the number of circles in the intersection
  $P_{\g(1,\ell_1)}\cap\bRt$.
\end{itemize}
The construction terminates if $P_{\g(1,\ell_1;1,\ell_2)}\cap \cL_i^- = \emptyset$ for $i=1,2$; otherwise,  it continues in a recursive
manner.

We comment on the details of Figure~\ref{fig:bubblesL-}. There are
two $\Gamma'$ curves, with the corresponding $\kappa'$ in dotted
lines, and one $\Lambda'$ curve, with the corresponding $\chi'$ in dotted
lines. These are level 1 curves. At level two, two curves in the
intersection  $P_{\g(1,\ell_1)}\cap\cL_1^-$ are illustrated:
$\g'(1,\ell_1;1,0)$ with one connected component, and 
$\g'(1,\ell_1;1,-1)$ formed of the union of two connected components.
These curves are contained in the regions of $\cL_1^-$ bounded by
$\G'(1,0)$ and $\G'(1,-1)$, respectively.

Analogous considerations apply to the curves
$\g(1,\ell_1;2,\ell_2)$  thus we obtain two
families of level 3 propellers $P_{\g(1,\ell_1;i_2,\ell_2)}$ for
  $i_2=1,2$ and $b\leq \ell_1<0$ that are part of $S_{\g(1,\ell_1)}$. Let
  $P'_{\g(1,\ell_1;i_2,\ell_2)}=P_{\g(1,\ell_1;i_2,\ell_2)}\cap \mW'$.
Observe that $\tau(P'_{\g(1,\ell_1;i_2,\ell_2)})\subset
S_{\g(1,\ell_1)}$. Thus for $b\leq \ell_1<0$,
\begin{equation}
\bigcup_{i_2=1,2}\bigcup_{\ell_2}\tau(P'_{\g(1,\ell_1;
  i_2,\ell_2)})\subset S_{\g(1,\ell_1)},
\end{equation}
where the union is taken over all the possible values of $\ell_2$,
once $i_2$ is given.

\begin{lemma}\label{lemma-bubblelevel}
There exists $N>0$ such that for any $x\in S_{\g(1,\ell_1)}$ with $x\notin \overline{\g}(1,\ell_1)$, then $2\leq n_0(x)\leq N+1$.
\end{lemma}

\proof The boundary of $S_{\g(1,\ell_1)}$ is composed by
$\tau(\g(1,\ell_1))\subset E_1$, the facing curve
$\tau(\overline{\g}(1,\ell_1))\subset S_1$ and the two $\cK$-orbit segments
going from the endpoints of the entry curve to the endpoints of the
facing exit curve. Since $r(\g_{(1,\ell_1)}(1))>2$ and $\g(1,\ell_1)$
is a level 2 curve, then for any $x\in
S_{\g(1,\ell_1)}$ we know that
$n_0(x)\geq 2$. 

Consider now $\g(1,0)$. Since $r(\g_{(1,0)}(1))>r_b>2$,  the $\cK$-orbit of every point in $\tau(\g(1,0))$
passes through the facing point, and thus  attains a maximum level
$N+1$. Thus $n_0(x)\leq N+1$.
\endproof

This implies that the construction process above terminates after at most
$N$ steps. Thus $S_{\g(1,\ell_1)}$ is the union of
$$\bigcup_{1\leq n\leq
  N+1}\tau(P'_{\g(1,\ell_1;i_2,\ell_2;\cdots
  ;i_n,\ell_n)}),$$
for any combinations of indices, $(i_2,\ell_2;\cdots;i_{n},\ell_n)$, where 
$b\leq \ell_1<0$.
 
In a similar way we obtain the compact surfaces $S_{\g(2,\ell)}$,
$S_{\lambda(1,\ell)}$ and $S_{\lambda(2,\ell)}$. We call these surfaces
the level 2 bubbles in $\fM_0$, since they are associated to a level 2
curve.

Recall that in (\ref{eq-level1approx})
the set $\fM_0^1$ was described for any value of $b$. The
construction of $\fM_0^{n+1}$ from $\fM_0^n$, for $n\geq 1$, follows
the lines of the construction of $\fM_0^1$ from
$\fM_0^0=\tau(\cR')$. To obtain $\fM_0^2$ we have to add the level 2 points in
the $\cK$-orbits of points in $\g(i,\ell)$ and $\lambda(i,\ell)$ for
$i=1,2$ and $b\leq \ell$  and unbounded. The case $\ell\geq 0$ was
described in Section~\ref{sec-proplevels}, here we consider the case
$b\leq \ell<0$ and $i=1,2$, thus the level 2 points in the bubbles $S_{\g(i,\ell)}$,
$S_{\lambda(i,\ell)}$. If $b\neq 0$, we add to the set $\fM_0^1$
described in Section~\ref{sec-proplevels} the notched propellers
$\tau(P'_{\g(i,\ell_1)})$ and $\tau(P'_{\lambda(i,\ell_1)})$ for $i=1,2$
and $b\leq \ell_1<0$ and the finite collection of exit curves
$\tau(\overline{\g}(i_1,\ell_1;i_2,\ell_2))$ and
$\tau(\overline{\lambda}(i_1,\ell_1;i_2,\ell_2))$.

\begin{remark}
The propeller $P_{\g(i_1,\ell_1)}$ with $\ell_1\geq 0$ might also have
internal notches, that is the intersection $P_{\g(i_1,\ell_1)}\cap
\cL_{i_2}^-$, for $i_2=1,2$, might have curves that are arcs having
both endpoints in $\partial \cL_{i_2}^-$. In this case we denote these curves
by $\g(i_1,\ell_1;i_2,\ell_2)$ with $\ell_2<0$ and $\ell_1\geq
0$. Again, such a curve generates a level 3 bubble
$S_{\g(i_1,\ell_1;i_2,\ell_2)}$ that is part of $\fM_0$. 

In general a  propeller $P_{\g(i_1,\ell_1;\cdots;i_n,\ell_n)}$ forms
part of a bubble if at least one of the indices is negative. Assuming
that $\ell_k$ is the first negative index, then
$P_{\g(i_1,\ell_1;\cdots;i_n,\ell_n)}\subset
S_{\g(i_1,\ell_1;\cdots;i_k,\ell_k)}$, that is a level $k+1$ bubble.
\end{remark}

The following result generalizes the above discussion to higher level
bubbles.

\begin{prop}\label{prop-bubbles1}
Let $\g'(i_1,\ell_1;\cdots;i_n,\ell_n)\subset \cL_{i_n}^-$ be any curve
in the construction of $\fM_0$ with at least one negative index. Let
$\ell_k<0$ be the first negative index. Then
\begin{enumerate}
\item
  $\g(i_1,\ell_1;\cdots;i_n,\ell_n)=\sigma_{i_n}^{-1}(\g'(i_1,\ell_1;\cdots;i_n,\ell_n))$
  generates a finite double propeller or a pair of finite double
  propellers. Each propeller intersects $\bRt$ along at
  most $K$ closed curves, for a fixed number $K$;
\item for every $x\in S_{\g(i_1,\ell_1;\cdots;i_n,\ell_n)}$ we have
  that $n+1\leq n_0(x)\leq N+n$, for $N$ as in Lemma~\ref{lemma-bubblelevel}.
\end{enumerate}
\end{prop}

\proof
For the first conclusion, observe that
$\g(i_1,\ell_1;\cdots;1,\ell_n)$ lies in the  region
$D(\kappa,1) \cup D(\chi,1) \subset L_1^-$
described in the proof of Proposition~\ref{prop-r_b}, and $\g(i_1,\ell_1;\cdots;2,\ell_n)$ lies in the
analogous  region $D(\kappa,2) \cup D(\chi,2) \subset L_2^-$. Thus the number of circles in the
intersection of $P_{\g(i_1,\ell_1;\cdots;i_n,\ell_n)}$ with $\bRt$ is
uniformly bounded by $\Delta(r_b)+1$.

The second conclusion follows from the arguments above and the proof
of Lemma~\ref{lemma-bubblelevel}.
\endproof

Observe that if $b\neq 0$ the double propellers $P_\G$ and $P_\Lambda$
will also have internal notches, that is notches that are not
intersecting the orbits of $\omega_1$ or $\omega_2$, accordingly. Thus
they will generate bubbles in the same manner as the simple propellers
$P_\g$ and $P_\lambda$ did. 

We consider the case of $P_\G$ in some detail. The intersection
$P_\G\cap \cL_1^-$ consists of the curves $\G(1,\ell_1)$ for $b\leq
\ell_1$ and unbounded. For $b\leq \ell_1<0$ the curve $\G(1,\ell_1)$
consists of two connected components each having its two endpoints in
$\partial \cL_1^-$, one that corresponds to $\g(1,\ell_1)$ and the
other one $\kappa(1,\ell_1)$ that belongs to
$P_\kappa\cap\cL_1^-$. Thus $\G(1,\ell_1)$ is the disjoint union of
$\g(1,\ell_1)$ and $\kappa(1,\ell_1)$, when $b\leq \ell_1<0$. The
curves $\g(1,\ell_1)$ generate the bubbles $S_{\g(1,\ell_1)}$ that
form part of $\widehat{\fM}_0$ as defined in \eqref{eq-doublelamination}. Analogously, the curves $\kappa(1,\ell_1)$ generate the bubbles $S_{\kappa(1,\ell_1)}$ that
form part of $\widehat{\fM}_0$. The results above apply without any
changes to the bubbles $S_{\kappa(1,\ell_1)}$. Let
$S_{\Gamma(1,\ell_1)}=S_{\g(1,\ell_1)}\cup S_{\kappa(1,\ell_1)}\subset
\widehat{\fM}_0$.

  \bigskip

\section{Wandering points and propellers} \label{sec-wandering}

The  open set of  wandering points $\fW \subset \fM$ for the flow $\Phi_t$  was defined by \eqref{eq-dynamicdecomp} in Section~\ref{sec-minimalset}, as a union of the four classes:
$$\fW    ~ = ~  \fW^{0}    \cup \fW^+ \cup \fW^- \cup \fW^{\infty} $$
where $\fW^0$ consist of all finite orbits, $\fW^+$ is contained in the set of forward trapped orbits,  $\fW^-$ is contained in the set of backward trapped orbits,  and 
$\fW^{\infty}$ is contained in the infinite orbits.

For example, $\fW$  contains the orbits of all points in the region $\{r < 2\}$ by Proposition~\ref{prop-wanderinginterior}. The orbits of points in $\fW^{\infty}$ exhibit the most subtle dynamical properties of the four classes. Here is the main result of the section:

   \begin{thm}\label{thm-wandering}
 Let $x \not\in \fM$, then  $x \in \fW$. 
 \end{thm} 
 
The proof of this result  requires  the introduction of  new classes
of double propellers which are used to describe the orbits of points  $x \in \bRt$ which satisfy $r(x) \geq 2$. 
Recall that a key point in the proofs of
Propositions~\ref{prop-wander} and \ref{prop-cluster-} was to analyze
the   orbits which intersect an entry region $E_i$   in the sets
$\tau(\cE_i^{-,+})$, as illustrated in Figure~\ref{fig:Upsilon}, which
consist  of points that are mapped from the region $\{r > 2\}$ to the
region $\{r < 2\}$ by the insertion maps $\sigma_i$. The corresponding
regions $\sigma_i^{-1}(\cE_i^{-,+})\subset L_i^-$ for $i=1,2$ are
bounded by what we call ``$G-L$'' curves. We construct the propellers
generated by these curves and consider their intersections with the
rectangle $\bRt$, generating families of ``$G_0-L_0$'' curves. These
curves   divide the region $\bRt\cap \{r\geq 2\}$ according to their
dynamical behavior, and the regions thus defined are used to encode
the trapped wandering orbits for $\Phi_t$.  One important application
of this method is Proposition~\ref{prop-r>2R0}, which implies that there are no trapped wandering   orbits which are strictly contained in the region $\{r > 2\}$.

\begin{figure}[!htbp]
\centering
{\includegraphics[width=60mm]{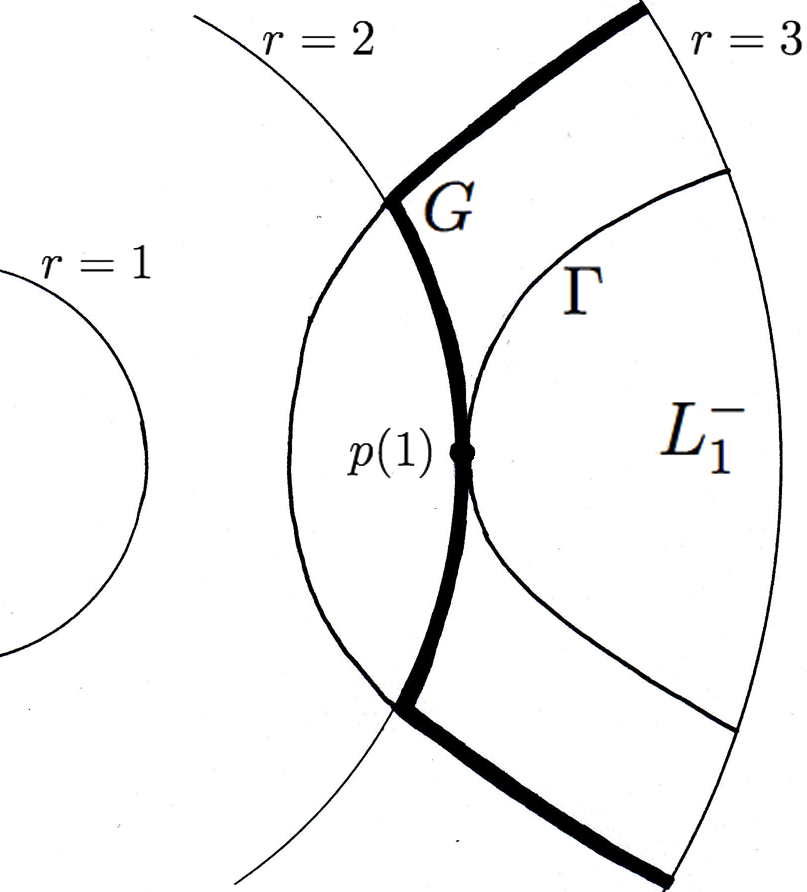}}
\caption{\label{fig:GL1} The   curve $G$ in $L_1^-$}
\vspace{-6pt}
\end{figure}

 We  introduce  the   curves $G \subset L_1^-$ and $L\subset L_2^-$,
 starting with the curve $G$. Recall from Section~\ref{sec-kuperberg} that the boundary of $L_i^-$
is composed of two arcs, $\alpha_i$ and $\alpha_i'$, as illustrated in Figure~\ref{fig:insertiondisks}, with $\alpha_i$  contained in the boundary of $\partial_h^-\mW$. 
The curve $G$ starts  at one of the points in the intersection
$\alpha_1\cap\alpha_1'$, following an arc in $\alpha_1'$   to the
first point with $r$-coordinate equal to 2, then follows the arc $L_1^-\cap
\{r=2\}$, and finally follows the curve  $\alpha_1'$ to the other point in
$\alpha_1\cap\alpha_1'$. Thus, the curve $G$, as illustrated in Figure~\ref{fig:GL1},  is composed of three  smooth arcs, where the first and last are contained in $\alpha_1'$ and  satisfy the radial monotonicity
assumption in Section~\ref{sec-intropropellers}, and the  middle arc  lies on a segment of the circle $\{ r=2\}$, so is tangent in its interior Êto the curve $\Gamma$ at $\sigma_1^{-1}(p_1^-)$ as illustrated in Figure~\ref{fig:GL1}. 
 Note that $G$ bounds   a region in $L_1^-$ that contains the pre-image of the set  $\cE_1^{-,+}$ under the insertion map $\sigma_1$. 
We define the curve $L \subset L_2^-$ using the curve $\alpha_2'$ in the analogous manner, with details omitted.

Next, we  form two infinite propellers  by considering the flow under $\cW$ of the curves $G$ and $L$, as for the curves $\Gamma$ and $\Lambda$ in Section~\ref{sec-intropropellers}.
The $\Psi_t$-flow of the middle arc in $G$ is an infinite strip which spirals around the cylinder $\cC = \{ r= 2\}$ in the region $\{-2 \leq z < -1\}$, and contains  the spiraling curve  $\cZ_{\gamma}^- = \{\Psi_t(\sigma_1^{-1}(p_1^-)) \mid t \geq 0 \}$     in the interior of the strip.
Let $\overline{G}\subset L_1^+$ be the facing curve to $G$, then similar observations apply to the reverse $\Psi_t$-flow of $\overline{G}$ and the spiraling curve 
 $\cZ_{\gamma}^+ = \{\Psi_t(\sigma_1^{-1}(p_1^+)) \mid t \leq 0
 \}$. The $\Psi_t$-flow of the two arcs of $G$ which are contained in  $\alpha_1'$ generate two infinite propellers.
Form the   infinite double   propeller $P_G \subset \mW$ from the union of these two propellers and the two infinite strips. Then $P_G$  contains  
$\ds \cZ_{\g} = \cZ_{\g}^- ~ \cup ~ \cZ_{\g}^+$ defined in  \eqref{eq-Zset}.

Let $P_L$ be the infinite double   propeller  corresponding to $L$   defined in the same way, which then contains $\cZ_{\lambda} = \cZ_{\lambda}^- \cup \cZ_{\lambda}^+$.   As the curves $G$ and $L$ are disjoint, the propellers $P_G$ and $P_L$ are also disjoint. 
 We then have the analogous result to Proposition~\ref{prop-propellerclosure}.

\begin{prop} \label{prop-GLclosure}
The closures $\overline{P}_{G}$ and $\overline{P}_{L}$ of   the infinite propellers $P_G$ and $P_L$   contain the Reeb cylinder $\cR$,    and we have
\begin{equation}
\overline{P}_{G} ~ = ~ P_{G} ~ \cup ~ \cR   ~ , ~ \overline{P}_{L} ~ = ~ P_{L} ~ \cup ~ \cR   .
\end{equation}
\end{prop}
\proof 
Let $\{x_n \in P_G \mid n = 1,2, \ldots\}$ converge to a point
$x_*$. If $r(x_*) > 2$, then we can assume that $r(x_n)>2$ and as the
Wilson flow preserves the radius coordinate $r$, there is a
corresponding sequence of points $\{y_n \in \alpha_1' \cap G \subset
L_1 \mid n = 1,2, \ldots\}$, with $r(y_n) = r(x_n)$ and $y_n$, $x_n$
in the same $\cW$-orbit. Then $r(y_n) \to r(x_*) > 2$, so $x_* \in
P_G$. In the case that $r(x_*) = 2$, then the sequence $\{y_n\}$
converges is either in $G$ or in $\overline{G}$. Assuming it belongs
to $G$, it converges to a point in $G \cap
\{r=2\}$. The forward Wilson flow of these points  converge
to the periodic orbit $\cO_1$. Then the points $x_n$ converge to a
point in the cylinder $\cR$, so $x_* \in \cR$. A similar argument holds for the backward Wilson flow of the facing curve $\overline{\alpha_i'}$, for which the inner endpoints limit on the periodic orbit $\cO_2$.

The analysis of the closure of the propeller $P_L$ proceeds similarly.
\endproof

Observe that $G$ divides $\partial_h^-\mW$ into two open regions. 
Let $\cU$ denote the region in the complement of $G$ which is   contained in $L_1^-$, so that points in $\cU$ have $r$-coordinate bigger than $2$. Let $U$ denote the open region in $\mW$ obtained from the $\Psi_t$-flow of $\cU$ from $L_1^-$ to $L_1^+$. 
It follows  that $P_G$ is contained in the closure of $U$, and so  $\overline{P_G}\subset \overline{U}$. 
  Since $\G$ is contained
in the closure of $\cU$ and $\G\setminus \{\sigma_1^{-1}(p_1^-)\}$
is contained in $\cU$, we conclude that    $P_{\G} - \cZ_{\gamma} \subset U$ as well.
We say that  the double propeller  $P_G$   \emph{envelops}   the double propeller
$P_{\G}$.

Similarly, let $\cV$ denote the region in   the complement of $L$ which is   contained in $L_2^-$, so that points in $\cV$ have $r$-coordinate bigger than $2$.
Let $V$ denote the open region in $\mW$ obtained from the $\Psi_t$-flow of $\cV$ from $L_2^-$ to $L_2^+$. 
Then $P_L$ is contained in the closure of $V$, and $\overline{P_L}\subset \overline{V}$. 
Correspondingly, we say that  $P_L$ envelops $P_\Lambda$.

 Next, consider the notched double propellers
$$P_G'=P_G\cap \mW' \qquad \mbox{and} \qquad P_L'=P_L\cap\mW'.$$ 
Since $\cZ_\g\subset P_G$, the propeller $P_G'$ intersects $\cL_1^-$
infinitely many times, once for each intersection of  
the $\cW$-orbit of the special point $\sigma_1^{-1}(p_1^-)$ with $\cL_1^-$. 
In the same way, $P_G$ intersects $\cL_2^-$ infinitely many times.

\begin{figure}[!htbp]
\centering
\begin{subfigure}[c]{0.4\textwidth}{\includegraphics[height=80mm]{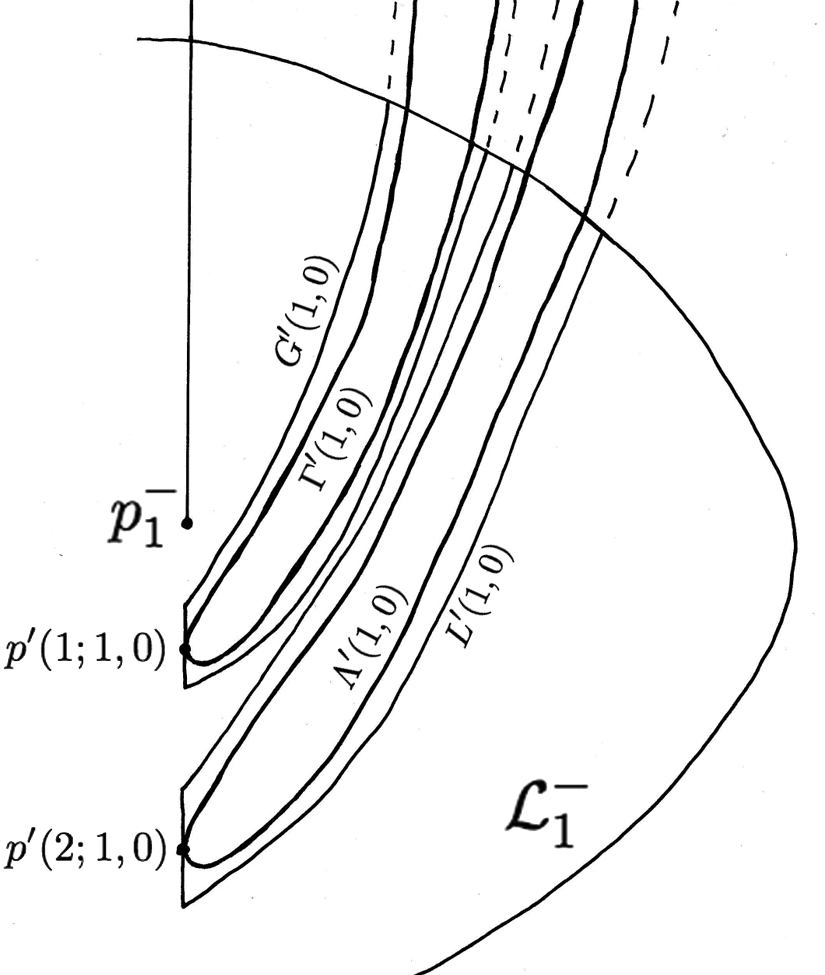}}\caption{$\Gamma'(1,0) \subset G'(1,0)$ and $\Lambda'(1,0) \subset L'(1,0)$}\end{subfigure}
\begin{subfigure}[c]{0.4\textwidth}{\includegraphics[height=80mm]{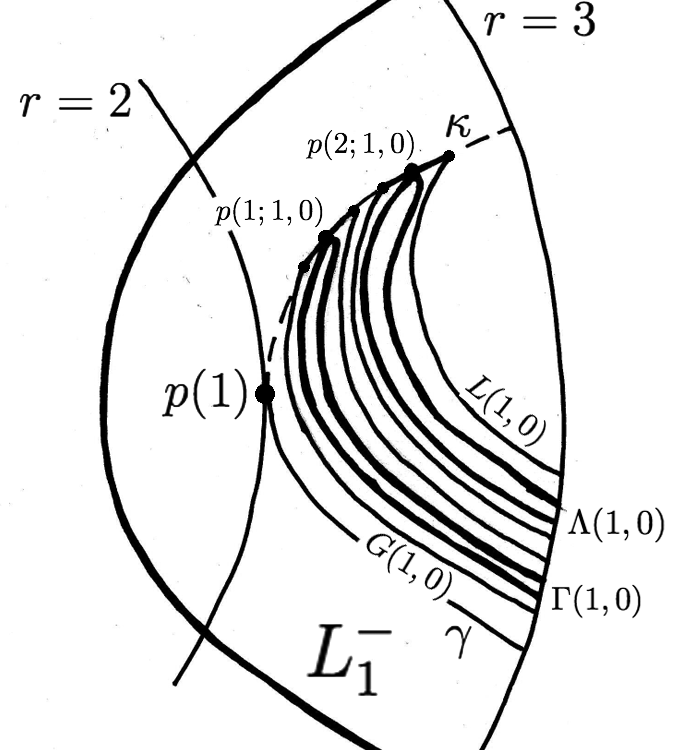}}\caption{$\Gamma(1,0) \subset G(1,0)$ and $\Lambda(1,0) \subset L(1,0)$}\end{subfigure}
\caption{\label{fig:GLcurvesL1} $G$ and $L$ curves at level $1$ in the regions $\cL_1^-$ and $L_1^-$} 
\vspace{-6pt}
\end{figure}

For each $i =1,2$, the intersections  $P_G \cap \cL_i^-$   and $P_L \cap \cL_i^-$  form two
infinite collections of curves $G'(i,\ell)$ and $L'(i,\ell)$, respectively, for
$b \leq \ell $ and unbounded.  
Apply $\sigma _i^{-1}$ to these curves, to obtain in
$L_1^-\cup L_2^-\subset \partial_h^-\mW$, four countable collections
of curves labeled, for $\ell \geq b$, as illustrated in   Figure~\ref{fig:GLcurvesL1}:
\begin{itemize}
\item $G(1,\ell)$ tangent at $p(1;1,\ell)$ to $\G(1,\ell)$ in
  $L_1^-$;
\item $L(1,\ell)$ tangent at $p(2;1,\ell)$ to $\Lambda(1,\ell)$ in
  $L_1^-$;
\item $G(2,\ell)$ tangent at $p(1;2,\ell)$ to $\G(2,\ell)$ in
  $L_2^-$;
\item $L(2,\ell)$ tangent at $p(2;2,\ell)$ to $\Lambda(2,\ell)$ in $L_2^-$.
\end{itemize}
For $b\leq \ell<0$  as introduced in Section~\ref{sec-proplevels} and  illustrated in
Figure~\ref{fig:bubblesL-}, 
 the curves $G'(i,\ell)$ and $L'(i,\ell)$ have two
connected components, which follows from the discussion   in Section~\ref{sec-bubbles}. As in
Proposition~\ref{prop-bubbles1}.1 the bubbles generated by these
curves have uniformly bounded level difference, and thus do not change
the discussion below. Thus, without loss of generality we assume for the rest of the section that $b=0$.

Observe that since the   $\Psi_t$-flow of the   middle arc in $G$ spirals around the cylinder $\cC \cap \{z < -1\}$,  it intersects the interval $J_0$ defined in \eqref{eq-intervals} in a sequence of closed intervals.  Hence  each curve $G(1,\ell)$ contains a closed arc contained in $\kappa\subset \G$, for $\kappa$ as defined in \eqref{eq-kappa}. In
the same way, each curve $L(1,\ell)$ has an arc in $\kappa$. Similarly, 
each of the curves $G(2,\ell)$ and $L(2,\ell)$ intersects $\chi \subset \Lambda$
along a closed  arc.

The endpoints of each one of the $G(i,\ell)$ and $L(i,\ell)$ curves are contained in the
boundary of $\partial_h^-\mW$, while the middle arcs are in
$\kappa$ or $\chi$, accordingly. By flowing these curves
in $\mW$, we obtain four countable collections of finite double
propellers, denoted by  $P_{G(i,\ell)}$ and $P_{L(i,\ell)}$. The curves $G(i,\ell)$ and $L(i,\ell)$ do not satisfy the radial monotonicity
assumption in Section~\ref{sec-intropropellers}. However, the two arcs that are
not contained in $\kappa$ or $\chi$, are endpoint isotopic
  to a curve satisfying the transversality condition, hence the
  conclusion that each of the curves $G(i,\ell)$ and $L(i,\ell)$ separates the
  region $L_i^-$ into two open regions remains true. The finite double propellers
  $P_{G(i,\ell)}$ and $P_{L(i,\ell)}$ are  isotopic to a standard
  finite double propeller and separate $\mW$ into two connected components.

Continue the above process recursively,  to obtain collections of  curves
\begin{itemize}
\item $\ds G(i_1,\ell_1;i_2,\ell_2;\cdots;i_n,\ell_n) \subset L_{i_n}^- \subset \partial_h^-\mW$, 
tangent to $\ds \G(i_1,\ell_1;i_2,\ell_2;\cdots;i_n,\ell_n)$ at\\ $p(1;i_1,\ell_1;i_2,\ell_2;\cdots;i_n,\ell_n)$;\\ 
\item $\ds L(i_1,\ell_1;i_2,\ell_2;\cdots;i_n,\ell_n) \subset L_{i_n}^-  \subset \partial_h^-\mW$ 
tangent to  $\ds \Lambda(i_1,\ell_1;i_2,\ell_2;\cdots;i_n,\ell_n)$ \\
at $p(2;i_1,\ell_1;i_2,\ell_2;\cdots;i_n,\ell_n)$.
\end{itemize}
Let
$\ds    P'_{G(i_1,\ell_1;i_2,\ell_2;\cdots;i_n,\ell_n)}\subset \mW'$  
and  $\ds P'_{L(i_1,\ell_1;i_2,\ell_2;\cdots;i_n,\ell_n)}\subset \mW'$
be the corresponding notched propellers.

We next discuss topological properties of the families of propellers formed from the $G$ and $L$ curves. 
Recall that the points $p(1;i_1,\ell_1;i_2,\ell_2;\cdots;i_n,\ell_n)$ divide the parabolic $\G$-curves  in two arcs, 
$$\G(i_1,\ell_1;i_2,\ell_2;\cdots;i_n,\ell_n) ~ = ~ \g(i_1,\ell_1;i_2,\ell_2;\cdots;i_n,\ell_n) ~ \cup ~\kappa(i_1,\ell_1;i_2,\ell_2;\cdots;i_n,\ell_n) .$$ 
The level $n$ curve $\G(i_1,\ell_1;i_2,\ell_2;\cdots;i_n,\ell_n)$ is tangent at
$p(1;i_1,\ell_1;i_2,\ell_2;\cdots;i_n,\ell_n)$ to a level $n-1$
curve. For $i_1=1$, it is tangent to the curve
$\G(i_2,\ell_2;\cdots;i_n,\ell_n)$ and  
$$p(1;i_1,\ell_1;i_2,\ell_2;\cdots;i_n,\ell_n) ~ {\rm  belongs ~ to ~ } ~ \kappa(i_2,\ell_2;\cdots;i_n,\ell_n). $$
  For $i_1=2$ the curve $\G(i_1,\ell_1;i_2,\ell_2;\cdots;i_n,\ell_n)$  is tangent to the curve
$\Lambda(i_2,\ell_2;\cdots;i_n,\ell_n)$ and  
$$p(1;i_1,\ell_1;i_2,\ell_2;\cdots;i_n,\ell_n) ~ {\rm  belongs ~ to ~ } ~ \chi(i_2,\ell_2;\cdots;i_n,\ell_n) .$$ 
Thus:
\begin{itemize}
\item  $p(1;1,\ell_1;i_2,\ell_2;\cdots;i_n,\ell_n)$ and
$p(2;1,\ell_1;i_2,\ell_2;\cdots;i_n,\ell_n)$ lie in
$\kappa(i_2,\ell_2;\cdots;i_n,\ell_n)$ for any $\ell_1$;

\item  $p(1;2,\ell_1;i_2,\ell_2;\cdots;i_n,\ell_n)$ and
$p(2;2,\ell_1;i_2,\ell_2;\cdots;i_n,\ell_n)$ lie in
$\chi(i_2,\ell_2;\cdots;i_n,\ell_n)$ for any $\ell_1$.
\end{itemize}

\begin{lemma} \label{lem-disjointGL}
Any two distinct $G$ and $L$ curves of the collection above, contained in $L_1^-\cup L_2^-$,
are disjoint.
\end{lemma}

\proof
The curve $\G$ divides $L_1^-$ in two regions, one of which is
contained in $\cU$ and its closure contains all the $G(1,\ell)$ and
$L(1,\ell)$ curves. The closure of this region intersects
$G$ only at the special point $\sigma_1^{-1}(p_1^-)$. Hence, $G$ is disjoint from all $G(1,\ell)$ and
$L(1,\ell)$. In the same way, $L$  is disjoint from all $G(2,\ell)$ and
$L(2,\ell)$. Then proceed inductively, to obtain that any two different curves in the collection
of $G$ and $L$ curves are disjoint.
\endproof

\begin{figure}[!htbp]
\centering
{\includegraphics[width=120mm]{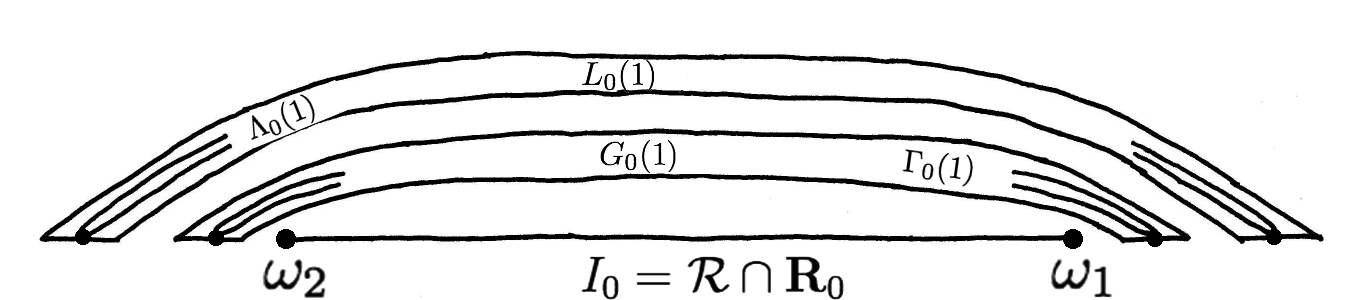}}
\caption[{Endparts of $\G_0$ and $\Lambda_0$ curves inside $G_0$ and $L_0$}]{\label{fig:GL2} Endparts of $\G_0$ and $\Lambda_0$ curves inside $G_0$ and $L_0$ curves in $\bRt$ (viewed sideways)}
\vspace{-6pt}
\end{figure}

We now consider the intersection of the notched propellers
$$P'_{G(i_1,\ell_1;i_2,\ell_2;\cdots;i_n,\ell_n)} \quad , \quad  P'_{L(i_1,\ell_1;i_2,\ell_2;\cdots;i_n,\ell_n)}$$
 with the rectangle
$\bRt$. Following the   notation convention  developed  in Section~\ref{sec-doublepropellers}, 
the intersections
$P_G'\cap\bRt$ and $P_L'\cap\bRt$ yield  countable collections of closed curves, labeled $G_0(\ell)$ and $L_0(\ell)$, 
as illustrated in  Figure~\ref{fig:GL2}. Observe that for each
$\ell\geq a$, where   $a\leq 0$   was defined   in
Remark~\ref{rmk-notation}, the curve $G_0(\ell)$ contains a segment of $J_0$, a segment of
$K_0$, and is
tangent at the two points $p_0(1;i,\ell)$  to the curve 
$\G_0(\ell)$ for $i=1,2$. In the same way, each $L_0(\ell)$ contains a segment of $J_0$, a segment of
$K_0$ and is
tangent at the two points $p_0(2;i,\ell)$ to the curve 
$\Lambda_0(\ell)$  for $i=1,2$.

Each $G_0(\ell)$ curve divides $\bRt$ in two open connected regions,
one of which contains $\G_0(\ell)$ and is contained in $U \cap \bRt$. We call this
region the interior of $G_0(\ell)$ and we
 denote it  by $\cU_0(\ell)$. 
 Analogously, the interior of each 
$L_0(\ell)$ curve is denoted by $\cV_0(\ell)$.

As in the previous section, we define families of $G_0$ and $L_0$ curves   in $\bRt\cap\{r\geq 2\}$  recursively.

Recall that $U$ denotes the interior region of $P_G \subset \mW$, so that $\cU$ as defined previously satisfies 
$\cU = U \cap \partial_h^- \mW \subset L_1^-$.
Also,   $V$ denotes the interior region  of $P_L \subset \mW$, so that $\cV$ as defined previously satisfies 
$\cV = V \cap \partial_h^- \mW \subset L_2^-$.

In a corresponding manner, for the families of $G$ and $L$ curves, let $U(i_1, \ell_1, \cdots, i_n, \ell_n)$ denote the interior region of  $P_{G(i_1, \ell_1, \cdots, i_n, \ell_n)} \subset \mW$, and $V(i_1, \ell_1, \cdots, i_n, \ell_n)$ denote the interior region of  $P_{L(i_1, \ell_1, \cdots, i_n, \ell_n)} \subset \mW$.
Then we define open regions in $\partial_h^- \mW$    as follows:

\begin{eqnarray}
\cU(i_1, \ell_1; \cdots ; i_n, \ell_n) & = &  U(i_1, \ell_1; \cdots ; i_n, \ell_n)   \cap \partial_h^- \mW   \label{eq-Ui1-elln} \\
\cV(i_1, \ell_1; \cdots ; i_n, \ell_n) & = &  V(i_1, \ell_1; \cdots ; i_n, \ell_n)   \cap \partial_h^- \mW .  \label{eq-Vi1-elln}
\end{eqnarray}

The intersection of the regions $U(i_1, \ell_1; \cdots ; i_n, \ell_n)$ and $V(i_1, \ell_1; \cdots ; i_n, \ell_n)$ with $\bRt$ form two finite collections of connected open regions,
which are indexed as follows:
\begin{eqnarray}
\cU_0(i_1, \ell_1; \cdots ;  i_n, \ell_n; \ell_{n+1}) & \subset &  \tau(U(i_1, \ell_1; \cdots ; i_n, \ell_n)   \cap \bRt) ~, ~ \ell_{n+1} \geq a  \label{eq-U0i1-elln}\\
\cV_0(i_1, \ell_1; \cdots ; i_n, \ell_n ; \ell_{n+1}) & \subset &  \tau(V(i_1, \ell_1; \cdots ; i_n, \ell_n)   \cap \bRt)  ~, ~ \ell_{n+1} \geq a,  \label{eq-V0i1-elln}
\end{eqnarray}
 which extends the notation convention introduced above for
 $\cU_0(\ell)$ and $\cV_0(\ell)$. Observe that the index $\ell_{n+1}$
 admits finitely many values since the boundary propellers $P_{G(i_1,
   \ell_1; \cdots ; i_n, \ell_n)}$ and $P_{L(i_1, \ell_1; \cdots ; i_n, \ell_n)}$ are finite.
 
 The  forward $\Psi_t$-flow of the region in \eqref{eq-U0i1-elln} to $\cL_1^-$ yields the region $\sigma_1(\cU(i_1, \ell_1; \cdots ; i_n, \ell_n; 1, \ell_{n+1})) \subset \cL_1^-$, and  its forward $\Psi_t$-flow to $\cL_2^-$ yields the region $\sigma_2(\cU(i_1, \ell_1; \cdots ; i_n, \ell_n; 2, \ell_{n+1})) \subset \cL_2^-$. Similar comments apply to the regions in \eqref{eq-V0i1-elln}.

Recall from  Lemma~\ref{lem-disjointGL} that the collection of all $G$ and $L$ curves in $\partial_h^- \mW$ are disjoint. This implies a containment property between their interior regions:
\begin{lemma}\label{lem-inclusions}
Suppose that either
\begin{itemize}
\item $\cU_0(i_1, \ell_1 ; \cdots ; i_n, \ell_n; \ell_{n+1}) \cap
  \cU_0(i_1', \ell_1'; \cdots ; i_{n'}', \ell_{n'}'; \ell_{n'+1}')
  \ne \emptyset$,
\item $\cU_0(i_1, \ell_1 ; \cdots ; i_n, \ell_n; \ell_{n+1}) \cap
  \cV_0(i_1', \ell_1'; \cdots ; i_{n'}', \ell_{n'}'; \ell_{n'+1}')
  \ne \emptyset$,
\item $\cV_0(i_1, \ell_1 ; \cdots ; i_n, \ell_n; \ell_{n+1}) \cap
  \cV_0(i_1', \ell_1'; \cdots ; i_{n'}', \ell_{n'}'; \ell_{n'+1}')
  \ne \emptyset$,
\end{itemize}
 then $n\ne
 n'$ and $\ell_{n+1} = \ell_{n'+1}'$. 

If $n'>n$, let $n'=n+m$. Then $\cU_0(i_1', \ell_1';
\cdots ; i_{n'}', \ell_{n'}'; \ell_{n+1})$ is contained in 
\begin{itemize}
\item $\cU_0(i_1, \ell_1 ; \cdots ; i_n, \ell_n; \ell_{n+1})$ if $i_{n'-n}=1$ and for
  any $1\leq k\leq n$
\begin{equation}\label{eq-coefficients1}
i_{k} = i_{k+m}'=i_{(n'-n)+k}' ~ {\rm and } ~ \ell_{k} = \ell_{k+m}'=\ell_{(n'-n)+k}' ~.
\end{equation}
\item $\cV_0(i_1, \ell_1 ; \cdots ; i_n, \ell_n; \ell_{n+1})$ if $i_{n'-n}=2$ and for
  any $1\leq k\leq n$
\begin{equation}\label{eq-coefficients2}
i_{k} = i_{k+m}'=i_{(n'-n)+k}' ~ {\rm and } ~ \ell_{k} = \ell_{k+m}'=\ell_{(n'-n)+k}' ~.
\end{equation}
\end{itemize}
The same conclusions hold for $\cV_0(i_1', \ell_1';
\cdots ; i_{n'}', \ell_{n'}'; \ell_{n+1})$.
\end{lemma}
\proof
The boundaries of the $\cU_0$ and $\cV_0$ regions are disjoint by Lemma~\ref{lem-disjointGL}, so the inclusion of one into the other follows, if the regions are not disjoint. 
The identities \eqref{eq-coefficients1} and \eqref{eq-coefficients2} involving the indices follow from the construction of the curves.
\endproof

 We describe two cases to illustrate the conclusion of Lemma~\ref{lem-inclusions}.

 \noindent $\bullet$~ For  $i_1=i_2=1$, the region $\cU_0(1,\ell_1;1,\ell_2;\cdots;\ell_n)$ 
  is contained in 
   $\cU_0(1,\ell_2;\cdots;\ell_n)$ and $\cU_0(i_3,\ell_3;\cdots;\ell_n)$.

 \noindent $\bullet$~ For $i_1=2 , i_2=1$,  the region
 $\cU_0(2,\ell_1;1,\ell_2;\cdots;\ell_n)$ is contained in 
   $\cV_0(1,\ell_2;\cdots;\ell_n)$ and  $\cU_0(i_3,\ell_3;\cdots;\ell_n)$.

Observe that   each propeller $\ds P_{G(i_1,\ell_1;\cdots;i_n,\ell_n)}$ envelops the propeller
$\ds P_{\G(i_1,\ell_1;\cdots;i_n,\ell_n)}$.  This yields a nested relation between   the regions bounded by the $G_0$ and $L_0$ curves, and the $\Gamma_0$ and $\Lambda_0$ curves in $\bRt$ .  This is illustrated in Figure~\ref{fig:GL2}, and we describe the relation in two special cases.

 \noindent $\bullet$~  The region $\cU_0(i_1,\ell_1;i_2,\ell_2;\cdots;\ell_n)$   contains   the curve  $\G_0(i_1,\ell_1;i_2,\ell_2;\cdots;\ell_n)$
in its  interior,  and the intersection of the two associated curves consists of the points
$$G_0(i_1,\ell_1;i_2,\ell_2;\cdots;\ell_n) \cap \G_0(i_1,\ell_1;i_2,\ell_2;\cdots;\ell_n) = p_0(1;i_1,\ell_1;i_2,\ell_2;\cdots;i_n,\ell_n) ~ , ~ i_n=1,2 . $$

For $i_1=1$,  the curve $G_0(1,\ell_1;i_2,\ell_2;\cdots;\ell_n)$ intersects $\kappa_0(i_2,\ell_2;\cdots;\ell_{n})$ in
  two closed arcs, one containing 
  $p_0(1;1,\ell_1;i_2,\ell_2;\cdots;1,\ell_n)$ and the other
  containing    $p_0(1;1,\ell_1;i_2,\ell_2;\cdots;2,\ell_n)$.

  For $i_1=2$, the curve $G_0(2,\ell_1;i_2,\ell_2;\cdots;\ell_n)$ intersects $\chi_0(i_2,\ell_2;\cdots;\ell_{n})$ in
  two closed arcs,  one containing 
  $p_0(1;2,\ell_1;i_2,\ell_2;\cdots;1,\ell_n)$ and the other
  containing    $p_0(1;2,\ell_1;i_2,\ell_2;\cdots;2,\ell_n)$.

\medskip
  
 \noindent $\bullet$~  The region $\cV_0(i_1,\ell_1;i_2,\ell_2;\cdots;\ell_n)$   contains  the curve  $\Lambda_0(i_1,\ell_1;i_2,\ell_2;\cdots;\ell_n)$
 in its interior, and the intersection of the two associated curves consists of the points
 $$ L_0(i_1,\ell_1;i_2,\ell_2;\cdots;\ell_n) \cap \Lambda_0(i_1,\ell_1;i_2,\ell_2;\cdots;\ell_n) = p_0(2;i_1,\ell_1;i_2,\ell_2;\cdots;i_n,\ell_n) ~, ~ i_n=1,2$$

For $i_1=1$, the curve $L_0(1,\ell_1;i_2,\ell_2;\cdots;\ell_n)$ intersects $\kappa_0(i_2,\ell_2;\cdots;\ell_{n})$ in
  two closed arcs,  one containing 
  $p_0(2;1,\ell_1;i_2,\ell_2;\cdots;1,\ell_n)$ and the other
  containing    $p_0(2;1,\ell_1;i_2,\ell_2;\cdots;2,\ell_n)$.

 For $i_1=2$, the curve $L_0(2,\ell_1;i_2,\ell_2;\cdots;\ell_n)$ intersects $\chi_0(i_2,\ell_2;\cdots;\ell_{n})$ in
  two closed arcs,  one containing 
  $p_0(2;2,\ell_1;i_2,\ell_2;\cdots;1,\ell_n)$ and the other
  containing    $p_0(2;2,\ell_1;i_2,\ell_2;\cdots;2,\ell_n)$.

  With these preparations, we can  make the following definition for the points in $\bRt$ with $r(x) > 2$, which extends the notion of ``level'' for those points   not contained in $\fM$. Recall that $\fMR = \fM \cap \bRt$.
  
   \begin{defn}\label{def-level}
   Let $x \in \bRt$ with $r(x) > 2$, and  $x \not\in \fMR$. 
   Say that $x$ has \emph{level at least $n$}, and write $L(x) \geq
   n$, if either there exists an open region $\cU_0$ or $\cV_0$  of level $n$ such that $x$ is contained in the closure of either such regions:
   $$x \in \overline{\cU_0(i_1, \ell_1 ; \cdots ; i_{n-1}, \ell_{n-1}; \ell_{n})} \quad {\rm or} \quad  x \in \overline{\cV_0(i_1, \ell_1 ; \cdots ; i_{n-1}, \ell_{n-1}; \ell_{n})}.$$  
   The \emph{level} of $x$ is the greatest $n \geq 1$ such that $L(x) \geq n$. 
   If $L(x) \geq n$ for all $n$, then $x$ is said to have infinite level, and otherwise has finite level.
   If $x$ is not contained in any such region,  set $L(x) =0$.
      \end{defn}
   
   The following is a key property of the level.
   
   \begin{prop}\label{prop-depthfinite}
    Let $x \in \bRt$ with $r(x) > 2$, and  $x \not\in \fMR$. Then $x$ has finite level.
   \end{prop}
\proof
Let $\fG_0$ denote the union of the $G_0$-curves in $\bRt$, and $\fL_0$ the union of the $L_0$-curves.
Recall that $\fM$ is the closure of $\fM_0$ which is the union of all
$\gamma_0$ and $\lambda_0$ curves in $\bRt$.

Set $\fY = \fMR \cup \fG_0 \cup \fL_0$. The set $\fY$ is bounded in $\bRt$ so its closure is compact, but more is true.
\begin{lemma}\label{lem-closure}
The set $\fY$ is closed, hence is compact.
\end{lemma}
 \proof
 By Proposition~\ref{prop-GLclosure},  the limit of curves $G_0(\ell)$ and $L_0(\ell)$ as $\ell \to \infty$ is contained in $\fY$.
 Thus, for a fixed index $(i_1,*;i_2,\ell_2;\cdots;\ell_n)$, the limit of curves
 $G_0(i_1,\ell_1;i_2,\ell_2;\cdots;\ell_n)$ or $L_0(i_1,\ell_1;i_2,\ell_2;\cdots;\ell_n)$ as $\ell_1 \to \infty$ is contained in $\fY$.

Given a sequence of points in $\fY$ which are contained in $G_0$ or $L_0$ curves with indices $(i_1,\ell_1;i_2,\ell_2;\cdots;\ell_n)$ where the degree $n$ tends to infinity,   we observe that the sequence must be nested by Lemma~\ref{lem-inclusions}, and each interior region 
 $\cU_0(i_1,\ell_1;i_2,\ell_2;\cdots;\ell_n)$ or $\cV_0(i_1,\ell_1;i_2,\ell_2;\cdots;\ell_n)$ contains a corresponding $\gamma_0$ or $\lambda_0$ arc which is contained in $\fM_0$ hence in $\fY$. It follows that any convergent sequence of points in $\fG_0$ or $\fL_0$ of this form must have limit in the set $\fM$.
 \endproof
 
 Let $\e >0$ be chosen small enough, so that $\e < d_{\bRt}(x, \fMR)$ for the metric $d_{\bRt}$ defined in Section~\ref{sec-pseudogroup}. That is, the distance from every point of $\fMR$ to $x$ is greater than $\e$. Let 
 $$\cU_0(\fM, \e) = \{ y \in \bRt \mid d_{\bRt}(y, \fMR) < \e\}$$
Recall that the  $G_0$ and $L_0$ curves in $\bRt$ are disjoint closed curves   by Lemma~\ref{lem-disjointGL}, so for each index there exists an $\e_* > 0$ (where we abuse notation and do not indicate the precise index on $\e_*$) so that each of the sets
 \begin{eqnarray}
\cU_0(i_1,\ell_1;i_2,\ell_2;\cdots;\ell_n, \e_*) & = & \{ y \in \bRt \mid d_{\bRt}(y, G_0(i_1,\ell_1;i_2,\ell_2;\cdots;\ell_n)) < \e_*\}  \label{eq-Gnbhd}\\
\cV_0(i_1,\ell_1;i_2,\ell_2;\cdots;\ell_n, \e_*) & = & \{ y \in \bRt \mid d_{\bRt}(y, L_0(i_1,\ell_1;i_2,\ell_2;\cdots;\ell_n)) < \e_*\} \label{eq-Lnbhd}
\end{eqnarray}
contains exactly one $G_0$ or $L_0$ curve. We comment on these
definitions. 

\begin{figure}[!htbp]
\centering
{\includegraphics[width=66mm]{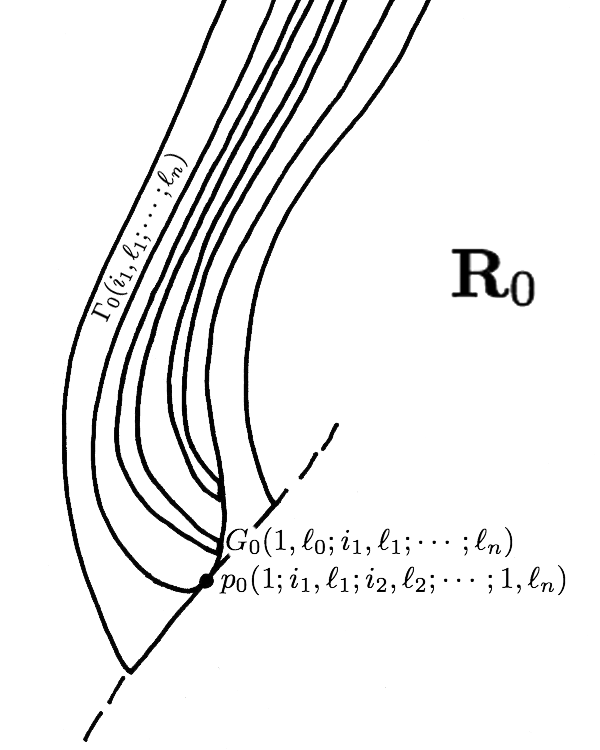}}
\caption[{Lower part of $\G_0$ and $G_0$ curves}]{\label{fig:convergeGL}  Lower part of $\G_0(i_1,\ell_1;i_2,\ell_2;\cdots;\ell_n)$ inside region bounded by $ G_0(i_1, \ell_1; i_2,\ell_2;\cdots;\ell_n)$}
\vspace{-6pt}
\end{figure}

Observe that the curve
$G_0(i_1,\ell_1;i_2,\ell_2;\cdots;\ell_n)$ contains the points
$p_0(1;i_1,\ell_1;i_2,\ell_2;\cdots;i_n,\ell_n)$ for $i_n=1,2$. Since
the curves $G_0(1,\ell;i_1,\ell_1;i_2,\ell_2;\cdots;\ell_n)$ limit as
$\ell\to\infty$ to $\g_0(i_1,\ell_1;i_2,\ell_2;\cdots;\ell_n)$, for
$\ell$ big enough the open set
$\cU_0(i_1,\ell_1;i_2,\ell_2;\cdots;\ell_n, \e_*)$ intersects these
curves, but it does not contains them, as illustrated in Figure~\ref{fig:convergeGL}.

Form an open covering of $\fY$ which consists of the open neighborhood $\cU_0(\fM, \e)$ of $\fMR$, and all open sets of the form 
\eqref{eq-Gnbhd} and \eqref{eq-Lnbhd}. Then by Lemma~\ref{lem-closure}, there is a finite subcovering, which must consist of $\cU_0(\fM, \e)$ and a finite collection of open sets of the form \eqref{eq-Gnbhd} and \eqref{eq-Lnbhd}. It follows that for all but a finite number of exceptions, the $G_0$ and $L_0$ curves are contained in $\cU_0(\fM, \e)$, and thus the closures of their interiors are also contained in $\cU_0(\fM, \e)$. As $x \not\in \cU_0(\fM, \e)$, there are only a finite number of $G_0$ or $L_0$ curves whose interior closures   contain $x$. If $L(x) \ne 0$, then $L(x) $ is the maximum value of the level for all such interior regions which contain $x$, which is  finite.
\endproof

 We next use the   constructions and results above    to describe  the dynamical properties of the $\Phi_t$-orbits of  points $x \in \bRt$  with $r(x) \geq 2$ and $x \not\in \fM_{\bRt}$. First,   recall that the   map  $\tau \colon \mW' \to \mK$   identifies 
the surfaces $L_i^-$ and $\cL_i^-$ via the map $\sigma_i$,  to yield the entry region $E_i$, for $i =1,2$. 
  The map  $\sigma_i$ sends
$\alpha_i'$ to the arc $\beta_i'\subset \partial_h^-\mW\setminus
(L_1^-\cup L_2^-)$. The points in
$\sigma_1(G\cap\alpha_1')$ and $\sigma_2(L\cap\alpha_2')$ are technically secondary  
entry points of $\mK$ as defined in Section~\ref{sec-radius}, but can also be considered as primary entry points as they lie in the closure of the primary entry points. The two arcs $\sigma_1(G\cap\{r=2\})$ and
$\sigma_2(L\cap \{r=2\})$ are exactly the secondary entry points in
$E_1$ and $E_2$, respectively, whose $r$-coordinate is equal to 2. The forward orbit of these points has $r$-coordinate bigger or equal to 2 by 
Corollary~\ref{cor-A_0infinite}.

\begin{prop}\label{prop-L0andG0}
Suppose that  $x\in \bRt$  is not in   the $\cK$-orbit of $\omega_1$ or $\omega_2$, and  is   contained in a $L_0$ or $G_0$ curve, so that for some index $\{i_1,\ell_1;\cdots;\ell_n\}$ we have
$$
x \in G_0(i_1,\ell_1;\cdots;\ell_n) ~ {\rm or}~ x \in   L_0(i_1,\ell_1;\cdots;\ell_n) ~ .
$$
 
 \begin{enumerate}
\item If   $x\notin   \kappa_0(i_2,\ell_2;\cdots;\ell_{n})\cup
\chi_0(i_2,\ell_2;\cdots;\ell_{n})\cup J_0\cup K_0$,
 then there exists
  $t,t'\in \mR$ such that $\Phi_t(x) \in \tau(\beta_j'\times\{2\})$ and $\Phi_{t'}(x) \in \tau(\beta_j'\times\{-2\})$, for $j=1$ or $2$. Hence the orbit of $x$ escapes $\mK$ in positive and negative  time. \\
\item If   $x\in \kappa_0(i_2,\ell_2;\cdots;\ell_{n})  \cup
  \chi_0(i_2,\ell_2;\cdots;\ell_{n})\cup J_0\cup K_0$,
  then there exists $t\in \mR$ such that $y=\Phi_t(x)$ is a transition
  point with $\rho_x(t)=2$. If $y$ is a secondary entry point, its
  forward orbit limits to $\omega_1$, while if $y$ is a
  secondary exit point its backward orbit limits to $\omega_2$.
 \end{enumerate}
\end{prop}

\proof
We consider the case where $x\in G_0(i_1,\ell_1;\cdots;\ell_n)$, with the   
  case $x \in   L_0(i_1,\ell_1;\cdots;\ell_n) $ being analogous. 

If $n=1$, then $x$ is contained in one of the curves $G_0(\ell_1)$, as illustrated in Figure~\ref{fig:GL2}.
Let  $x' = \tau^{-1}(x) \in G_0(\ell_1)\subset \bRt$, then flow
backwards in $\mW$, to obtain a  point $x'_{-1}\in G$, where $G$ is
illustrated in Figure~\ref{fig:GL1}. As $x$ does not lie on the flow
of the special points, the point $x'_{-1}$ lies either on the curve
$G\cap\alpha_1'$ with $r(x'_{-1}) > 2$ or in $G\cap\{r=2\}$. The first
case corresponds to (1) and we have that $x_{-1}'$ is a secondary
entry point on the boundary of the primary entry points. It follows
that  the backward $\Phi_t$-orbit of $x_{-1} = \tau(x'_{-1})$  escapes
$\mK$ in finite time. By Proposition~\ref{prop-propexist} the forward
$\Phi_t$-orbit of $x_{-1}$ passes through the facing point
$\overline{x_{-1}}\in \tau(\beta_1'\times\{2\})$  and thus the forward
orbit exits $\mK$ as claimed. The second case corresponds to (2) and $x_{-1}$ is the secondary entry point $y$, whose forward orbit limits to $\omega_1$ by Proposition~\ref{prop-r=2any}.

Next, assume that $n\geq 2$, then $r(x)>2$. Since $x\in G_0(i_1,\ell_1;\cdots;\ell_n)$ it belongs to the
propeller $\tau(P'_{G(i_1,\ell_1;\cdots;i_{n-1},\ell_{n-1})})$.
Then $x'=\tau^{-1}(x)$ is in the propeller $P_{G(i_1,\ell_1;\cdots;i_{n-1},\ell_{n-1})}$. 
Flow  $x'$ backwards in $\mW$, to obtain a
  point $x'_{-1}\in L_{i_{n-1}}^-$ that belongs to the curve
  $G(i_1,\ell_1;\cdots;i_{n-1},\ell_{n-1})$. By
  Proposition~\ref{prop-propexist}, $x_{-1}=\tau(x'_{-1})$ is a secondary
  entry point in the
  $\cK$-orbit of $x$. 

The point $x_{-1}$ can be identified with a point $x_{-1}^0\in
G_0(i_1,\ell_1;\cdots;\ell_{n-1})$: consider
$x_{-1}$ and flow it backwards from $E_{i_{n-1}}$ to
$\bRt$. Then $r(x_{-1}^0)\geq 2$, and repeat the process
inductively to obtain $x'_{-(n-1)}\in G(i_1,\ell_1)\in L_{i_1}^-$ and
$x_{-(n-1)}^0\in G_0(\ell_1)$. We have two possible situations:

\smallskip

(1)~ If $r(x_{-(n-1)}^0)>2$, we have that $x_{-(n-1)}^0\in   G_0(\ell_1)\setminus (J_0\cup K_0)$ and 
$$
x\in G_0(i_1,\ell_1;\cdots;\ell_n)\setminus \left[
  \kappa_0(i_2,\ell_2; \cdots ;\ell_n) \cup  \chi_0(i_2,\ell_2; \cdots
  ;\ell_n) \cup J_0\cup K_0\right].
$$
Then the $\cW$-orbit of $\tau^{-1}(x_{-(n-1)}^0)$ is finite and
intersects $\partial_h^-\mW$ in a point $x_{-n}'$ contained in
$G\cap\alpha_1'$. Thus $\sigma_1(x'_{-n})$ is in the boundary of the primary entry points of
$\mK$, implying that the orbit of $x$
escapes in negative time. Since $r(\sigma_1(x_{-n}'))>2$, the orbit of $x$ escapes in forward time.

(2)~  If $r(x_{-(n-1)}^0)=2$, then either $z(x_{-(n-1)}^0)<-1$ or   $z(x_{-(n-1)}^0)>1$ and
$$x\in G_0(i_1,\ell_1;\cdots;\ell_n)\cap \left[ \kappa_0(i_2,\ell_2; \cdots ;\ell_n)  \cup \chi_0(i_2,\ell_2; \cdots ;\ell_n) \cup J_0\cup K_0\right] .
$$ 
Assume first that $z(x_{-(n-1)}^0)<-1$, or
  equivalently that $x_{-(n-1)}^0\in J_0$. Consider the point
  $\tau^{-1}(x_{-(n-1)}^0)$ and flow it backwards in $\mW$. We obtain an
  entry point $x_{-n}'$ in $G\cap\{r=2\}$ and $\sigma_1(x_{-n}')$ is a
  secondary entry point with $r$-coordinate equal to 2 in the $\cK$-orbit of $x$. By
  Proposition~\ref{prop-r=2any} its forward $\cK$-orbit accumulates on
  $\omega_1$.

If on the contrary $z(x_{-(n-1)}^0)>1$, or
  equivalently that $x_{-(n-1)}^0\in K_0$, consider the point
  $\tau^{-1}(x_{-(n-1)}^0)$ and flow it forward in $\mW$. We obtain an
  exit point $x_{-n}'$ in $\overline{G}\cap\{r=2\}$ and $\sigma_1(x_{-n}')$ is a
  secondary exit point with $r$-coordinate equal to 2  in the $\cK$-orbit of $x$. By
  Proposition~\ref{prop-r=2any} its backward $\cK$-orbit accumulates on
  $\omega_2$.
\endproof

  \medskip

We next investigate the behavior of the $\cK$-orbit of points in $\bRt$ with $r$-coordinate
bigger than 2 that do not belong to $\fY$.

\begin{prop}\label{prop-R0regions}
Let $x\in\bRt$ with $r(x) > 2$,    assume that $x \not\in \fY$. Let $0 \leq L(x) < \infty$ be the level of $x$, as in Definition~\ref{def-level}.  
  Then we have the following possibilities:
\begin{enumerate}
\item Let $L(x) =0$, so that $x$ is outside every $G_0$ and every $L_0$ curve, then
  $\rho_x(t)>2$ for all $t$, and the orbit of $x$ escapes $\mK$ in
  positive and negative time.\\

\item Let $L(x) = n > 0$, so that $x$ is  contained in an open $\cU_0$ or $\cV_0$ region of level $n$, and assume that $x$ is also in the interior of the corresponding $\G_0$ or   $\Lambda_0$ curve whose vertex lies on the boundary of the region. 
 Then $\rho_x(t)>2$
  for all $t$,  and the orbit of $x$ escapes $\mK$ in
  positive and negative time.\\

\item Let $L(x) = n > 0$, and assume that $x$ lies in  the interior  of 
  $\cU_0(i_1,\ell_1;i_2,\ell_2;\cdots;\ell_n)$ and in the exterior of
  $\G_0(i_1,\ell_1;i_2,\ell_2;\cdots;\ell_n)$, then there exist
  $s<0$ such that $\Phi_{s}(x)$ is a secondary entry point, and 
  for $\e>0$ sufficiently small we have $\rho_x(s-\e)<2$.\\

\item Let $L(x) = n > 0$, and assume that $x$ lies in  the interior  of 
  $\cV_0(i_1,\ell_1;i_2,\ell_2;\cdots;\ell_n)$ and in the exterior of
  $\Lambda_0(i_1,\ell_1;i_2,\ell_2;\cdots;\ell_n)$, then there exist
  $s<0$ such that $\Phi_{s}(x)$ is a secondary entry point,  and 
  for $\e>0$ sufficiently small we have $\rho_x(s-\e)<2$.
\end{enumerate}
\end{prop}

\proof Let $x'\in\bRt\subset \mW$ such that $x=\tau(x')$. By assumption, 
$r(x')>2$ thus its Wilson orbit contains an entry point $x'_{-1}\in \partial_h^-\mW$ with $r(x'_{-1}) = r(x) > 2$. 
Then  Proposition~\ref{prop-propexist} implies
that $\tau(x'_{-1})=x_{-1}$ is in the $\cK$-orbit of $x$. If $x_{-1}$  is a  secondary entry point, we find $x_{-1}^0\in \bRt$ by flowing
$x_{-1}$ backwards.  If $r(x_{-1}^0)>2$ we can repeat this process for as
long as these conditions are satisfied. This reverse flow process  stops when either
$x_{-n}$ is a primary entry point,  or $r(x_{-n}^0)\leq 2$.

\smallskip

We  now analyze the four cases in the theorem:

(1)~ Since $x$ is outside every $G_0$ and every $L_0$ curve,
  $x'_{-1}\in \partial_h^-\mW\setminus (L_1^-\cup L_2^-)$ and hence
  $x_{-1}$ is a primary entry point in the $\cK$-orbit of $x$ with
  $r(x_{-1})=r(x)>2$. Then Proposition~\ref{prop-min} implies that
  $\rho_x(t) >2$ for all $t$, and that the orbit of $x$ escapes $\mK$ in
  positive and negative time.

 \smallskip
   
(2)~ We discuss   the case where $x$ is contained in a $\cU_0$ region. The alternate case where $x$ is contained in a $\cV_0$ region follows similarly.
The assumption that $L(x) = n$, and that  $x$   does not belong to any $G_0$ or $L_0$ curve, implies there is a unique region with 
$x \in \cU_0(i_1,\ell_1;i_2,\ell_2;\cdots;\ell_n)$.  It follows  that $x$ is contained 
 in the interior of the curve  $\G_0(i_1,\ell_1;i_2,\ell_2;\cdots;\ell_n)$, which is a connected component of  
 the  intersection of $\tau(P'_{\G(i_1,\ell_1;i_2,\ell_2;\cdots;i_{n-1},\ell_{n-1})}) \cap \bRt$.
 It then follows that the  point $x'_{-1}$ constructed above lies     inside the region bounded by 
  the curve $\G(i_1,\ell_1;i_2,\ell_2;\cdots;i_{n-1},\ell_{n-1}) \subset L_{i_{n-1}}^-$. 
    Thus,  $x_{-1}^0$ is in the interior of $\G_0(i_1,\ell_1;i_2,\ell_2;\cdots;\ell_{n-1})$,
and we have that  $r(x_{-1}^0)>2$.

We repeat the process recursively $n-1$ times, to obtain $x_{-(n-1)}^0\in
\bRt$ with $r(x_{-(n-1)}^0)>2$ and  contained  in the interior of $\G_0(\ell_1)$. The assumption that $L(x) =n$ implies that $x_{-(n-1)}^0$ lies in the exterior of  each $G_0(1,k;\ell_1)$ and $L_0(1,k;\ell_1)$ curve, for any $a \leq k< \infty$.

The backward $\cW$-orbit of  $x_{-(n-1)}^0$ yields the entry point  $x'_{-n}\in \partial_h^-\mW$. Observe
that  $r(x'_{-n})>2$  since the Wilson flow preserves the radius, and $x_{-n}'$
is  in the region inside $\G$. Hence $\tau(x_{-n}')\in
E_1\cap\{r>2\}$ and the backward $\cK$-orbit of this point to $\bRt$ yields 
the point $x_{-n}^0$. 
Then $r(x_{-n}^0)>2$ and $x_{-n}^0$ is outside every
$G_0$ and $L_0$ curve, so by case  (1) we conclude that $\rho_x(t)>2$ for
all $t$, and the orbit of $x$ escapes $\mK$ in positive and negative time.

(3)~  Since
  $G_0(i_1,\ell_1;i_2,\ell_2;\cdots;\ell_n)$ is in the
  intersection of $\tau(P'_{G(i_1,\ell_1;i_2,\ell_2;\cdots;i_{n-1},\ell_{n-1})})$ and
  $\bRt$, the point $x'_{-1} \in L_{i_{n-1}}^-$ is in the region between the curves
 $$G(i_1,\ell_1;i_2,\ell_2;\cdots;i_{n-1},\ell_{n-1}) \quad {\rm and} \quad \G(i_1,\ell_1;i_2,\ell_2;\cdots;i_{n-1},\ell_{n-1}).$$
  Thus, 
  $x_{-1}^0$ is in the region between the curves
  $$G_0(i_1,\ell_1;i_2,\ell_2;\cdots;\ell_{n-1})\quad {\rm and} \quad   \G_0(i_1,\ell_1;i_2,\ell_2;\cdots;\ell_{n-1}), $$
from which it follows that  $r(x_{-1}^0)>2$.

We repeat the process recursively $n-1$ times to obtain $x_{-(n-1)}^0\in
\bRt$ in the region between the curves $G_0(\ell_1)$ and
$\G_0(\ell_1)$. Then $r(x_{-(n-1)}^0)>2$ and the entry point
$x'_{-n}$ of its $\cW$-orbit has radius bigger than 2. Moreover,
since $G_0(\ell_1)$ is in the intersection of the propeller
$\tau(P_G')$ and $\bRt$, we conclude that $x'_{-n}\in
\cU \subset L_1^-$. Analogously, $x'_{-n}$ is outside the curve $\G \subset L_1^-$. Hence
$x_{-n}$ is a secondary entry point with $r(x_{-n})>2$. Then by the definition of the curve $\G$, for $\e>0$
small $\Phi_{-\e}(x_{-n})$ has $r$-coordinate less than 2. Since
$x_{-n}$ is in the backward $\cK$-orbit of $x$, there exists $s<0$ be
such that $\Phi_{s}(x)=x_{-n}$ and the conclusion follows.

(4)~ The proof of this case proceeds exactly  as for case (3).
\endproof

  We next  give   applications of the properties of $G-L$ curves  as developed in this section, to the study of the trapped  orbits   of points $x \in \mK$ in the complement of the compact subspace $\fM$.

\begin{prop}\label{prop-r>2R0}
Let $x\in \mK$ such that $x \not\in \fM$ and $\rho_x(t)>2$ for all $t$. Then the orbit of
$x$ escapes $\mK$ in positive and negative time.
\end{prop}

\proof
Let $x \in \mK$ with $x \not\in \fM$ and $r(x) > 2$. If the forward
$\cK$-orbit of $x$ exits $\mK$ through an exit point $z$, by
hypothesis $r(z)>2$. Hence
  the $\cK$-orbit of $x$ escapes
$\mK$ in positive and negative time by Proposition~\ref{prop-min}. 
Thus, we need to consider the case where the forward $\cK$-orbit of
$x$  is trapped.

Now assume that the forward $\cK$-orbit of $x$ is trapped, then  Proposition~\ref{prop-cluster+} implies that there exists a subsequence 
$\{ x_{\ell_i} \mid i=1,2, \ldots \} \subset \bRt$   such that $\ds \lim_{i\to \infty} ~ r(x_{\ell_i}') = 2$.  
Hence we can assume that the forward $\cK$-orbit of $x$ intersects
$\bRt$ in a point $x_1$ with $r(x_1)>2$ and arbitrarily close to $2$. 

By hypothesis, $x \not\in \fM$ which implies that $x_1 \not\in
\fM_{\bRt}$, so there exists   $\e > 0$ such that the open ball
$B_{\bRt}(x_1, \e) \subset \bRt$ about $x_1$  is disjoint from
$\fM_{\bRt}$, and hence is disjoint from every   $\g_0$ and
$\lambda_0$ curve. Moreover, since $\rho_x(t)>2$ for all $t$, then we
can choose $\e$ such that $B_{\bRt}(x_1, \e)$ is disjoint from every
$\G_0$ and $\Lambda_0$ curve.
  
  Then $x_1$ satisfies either the hypothesis of Proposition~\ref{prop-L0andG0}.1, the
hypothesis of Proposition~\ref{prop-R0regions}.1,  or the
hypothesis of Proposition~\ref{prop-R0regions}.2. In any case, the
$\cK$-orbit of $x_1$ escapes $\mK$ in positive and negative time, contradicting the assumption that the forward $\cK$-orbit of $x$ is trapped.
\endproof

Recall that the Matsumoto region for a Matsumoto constant $\delta_M > 0$, as defined in Definition~\ref{def-matsumoto},  is the set 
$\cU(\delta_M) = \tau(\{ y' \in \mW \mid 2-\delta_M < r(y') < 2\}) \subset    \mK$. 
We use the above results to   describe the intersection of the $\Phi_t$-flow of  the Matsumoto
region $\cU(\delta_M)$ with the rectangle $\bRt$, in terms of the $G_0$ and $L_0$ curves introduced above.

\begin{prop}\label{prop-wandering}
Suppose  that   $x\in \cU(\delta_M) \subset \mK$  is a secondary entry point, with 
 $y'\in (\cL_1^-\cup \cL_2^-)$ and $x'\in (L_1^-\cup L_2^-)$ such that $\tau(y')=\tau(x')=x$.
Assume that  $2 - \delta_M < r(y') < 2$  and $r(x')>2$.

\begin{enumerate}
\item If $z(y') < -1$, then the forward orbit of $x$ is trapped, and the
  points in the the intersection of its forward orbit and $\bRt \cap \{r\geq 2\}$  are in the regions between
  $G_0(i_1,\ell_1;\cdots;\ell_n)$ and
  $\kappa_0(i_1,\ell_1;\cdots;\ell_n)$, for any collection of indices. 
\item If $z(y')>1$, then the backward orbit of $x$ is trapped and the
  points in the the intersection of its backward orbit and $\bRt \cap \{r\geq 2\}$  are in the regions between
  $L_0(i_1,\ell_1;\cdots;\ell_n)$ and
  $\chi_0(i_1,\ell_1;\cdots;\ell_n)$, for any collection of indices. 
\item If $-1<z(y')<1$, then the orbit of $x$ is an infinite orbit and the
  points in the the intersection of its orbit and $\bRt \cap \{r\geq 2\}$  are in the regions between
  $L_0(i_1,\ell_1;\cdots;\ell_n)$ and
  $\lambda_0(i_1,\ell_1;\cdots;\ell_n)$,  or in the regions between
  $G_0(i_1,\ell_1;\cdots;\ell_n)$ and
  $\g_0(i_1,\ell_1;\cdots;\ell_n)$, for any collection of indices.
\end{enumerate}
\end{prop}

\proof
The first conclusion in each item above, that the orbit of $x$ is trapped, follows from the corresponding case  of 
Proposition~\ref{prop-wander}. For the other claims, we analyze  cases as follows.
 
(1)~ If  $z(y')<-1$, then  
  $x'$ is in the open region of $L_1^-$ bounded by $G$ and $\kappa$. The
  conclusion then follows by Proposition~\ref{prop-R0regions}.3.

(2)~ If  $z(y')>1$, then   
  $x'$ is in the open region of $L_2^-$ bounded by $L$ and $\chi$. The
  conclusion then follows by Proposition~\ref{prop-R0regions}.4.

(3)~ The case $z(y') = 0$ is impossible for a secondary entry point,  by the assumptions on the insertion maps. 
 If $-1<z(y')<0$, then 
  $x'$ is in the open region of $L_1^-$ bounded by $G$ and $\g$. If
  $0<z(y')<1$, then  $x'\in L_2^-$ in the open region bounded by
  $L$ and $\lambda$.    The
  conclusion then follows by Propositions~\ref{prop-R0regions}.3 and \ref{prop-R0regions}.4.
\endproof

We conclude this section with  the proof of Theorem~\ref{thm-wandering}. Recall that we assume $x \not\in \fM$.

First,   if $x$ is a primary entry or exit point, then $x \in \fW$ by  Lemma~\ref{lem-wanderingboundary}. In particular, $x \in \fW$ for every $x$ whose $\cK$-orbit is finite.
 Thus, we need only consider the case for $x$ such that its $\cK$-orbit is infinite.
 We analyze cases,    based on the behavior of the radius function $\rho_x(t) = r(\Phi_t(x))$.

If $x$ satisfies $\rho_x(t) > 2$ for all $t \in \mR$, then  $x$ has finite orbit by Proposition~\ref{prop-r>2R0}, so $x \in \fW$.

 If $\rho_x(t)< 2$ for some $t\in \mR$, then   $x$ is
 wandering by Proposition~\ref{prop-wanderinginterior} and hence $x \in \fW$.

The final case to consider is when $x$ is an infinite orbit with $\rho_x(t) \geq 2$ for all $t$.
Then  there is a point $y \in \bRt$ on the orbit of $x$, with $r(y) = 2$.
As $y \not\in \fM$, it must lie in the interior of one of the
segments $J_0$ or $K_0$ in \eqref{eq-intervals}.
Choose an open neighborhood $U(y,\e)$ for $\e > 0$ sufficiently small, which is disjoint from $\fM$.

If $y$ does not lie in   one of the $G_0$ or $L_0$ curves,  as
illustrated  in Figure~\ref{fig:GL1}, then  we can assume that $U(y,\e)$ is disjoint from these curves.
Consider the cases for the orbit of $z \in U(y,\e) \cap \bRt$.
If $r(z) < 2$ then $z$ is wandering by Proposition~\ref{prop-wanderinginterior}.
For $r(z) \geq 2$, then the point $z$ is outside any $G_0$ or $L_0$
curve. Proposition~\ref{prop-R0regions}.1 implies that the
$\cK$-orbit of $z$ is finite. Thus, $y \in \fW$ and hence $x \in \fW$.

If $y$ is contained in one of  the $G_0$ or $L_0$ curves, it is in a
level 1 curve. For $z \in U(y,\e) \cap \bRt$ we need only consider
the case when $r(z) \geq 2$. If $r(z) =2$ then
Proposition~\ref{prop-L0andG0}.2 implies its $\cW$-orbit is
asymptotic to a special point $\omega_i$. If $r(z) >  2$, then either
$z$ is in $\cU_0(\ell)$ and outside $\G_0(\ell)$ for
some $\ell$, or $z$ is in $\cV_0(\ell)$ and outside
$\Lambda_0(\ell)$ for some $\ell$, thus
Proposition~\ref{prop-R0regions}.3 or ~\ref{prop-R0regions}.4 applies. Then the orbit of $z$ has points with $r$-coordinate less than 2 and thus is a wandering point by Proposition~\ref{prop-wanderinginterior}.

 This covers the possible cases, showing that $\mK - \fM \subset \fW$, which is the assertion of  Theorem~\ref{thm-wandering}.
  
  In the next section, we give a condition that implies no points of $\fM$ are wandering, hence $\mK - \fM = \fW$.

  \bigskip
 
\section{The minimal set   for generic flows}\label{sec-generic}

 The construction of the flow $\Phi_t$ on the Kuperberg plug $\mK$  involves various choices,  in particular   the choice of the  vector field $\cW$ on $\mW$ satisfying the hypotheses in Section~\ref{sec-wilson}, and the choice of the insertions $\sigma_i$ satisfying the hypotheses in Section~\ref{sec-kuperberg}. For any such choice, the $\Phi_t$-flow   preserves  the compact set  $\fM$   defined by the closure of the orbits of the Reeb cylinder $\cR$. 
In this section, we consider  an additional  regularity  hypothesis  on the insertion maps $\sigma_i$,   which is a uniform version of Hypothesis~\ref{hyp-SRI}, and is used to estimate the behavior of the flow $\Phi_t$ near  the special orbits of $\Phi_t$. 

We   say that $\Phi_t$ is \emph{generic} if it satisfies the   hypotheses in Definition~\ref{def-generic}.
The  main result of this section is then:
   
\begin{thm}\label{thm-density}
If the flow  $\Phi_t$  is generic, then $\Sigma = \fM$.  
\end{thm}

The  proof of Theorem~\ref{thm-density}  will occupy the remainder of this section, and requires a detailed metric  analysis of the orbit structure of the  $\Psi_t$-flow near the Reeb cylinder $\cR$, which is used to obtain estimates on the $\Phi_t$-flow near the special orbits.

The equality  $\Sigma = \fM$  has been   previously observed in special cases.  Ghys gave an argument in   \cite[Th\'eor\`eme, p. 302]{Ghys1995} that  this conclusion holds  for   certain generic classes of insertions. The Kuperbergs constructed  in \cite{Kuperbergs1996},    an example  using polynomial vector fields     for which they sketch the proof of that  $\Sigma = \fM$. Our result,  which is motivated by these examples,   yields this conclusion in more generality. The   proof  makes full use of the body of techniques developed in this paper, 
and especially of the properties of the system of double propellers constructed in Section~\ref{sec-doublepropellers}. 
  For the more general case of the construction of flows on $\mK$ which do not satisfy the generic regularity hypotheses,  it seems possible that the inclusion $\Sigma \subset \fM$ may be proper, as was discussed in \cite{Kuperbergs1996}, in which  case, the minimal set $\Sigma$ will be a ``Denjoy-type'' invariant set for a flow on a surface lamination.
 
The additional hypotheses   is formulated using  the notation of  Hypothesis~\ref{hyp-SRI} on the insertion maps $\sigma_i$ for $i=1,2$. 
The projection along the $z$-coordinate in $\mW$ is denoted by   $\pi_z(r, \theta, z) = (r, \theta,-2)$, and 
 we assume that $\sigma_i$ restricted to the   face, $\sigma_i \colon L_i^- \to \mW$,  has image  transverse to the vertical fibers of $\pi_z$. 
 This condition is implicit in the illustrations Figures~\ref{fig:twisted}, \ref{fig:K} and \ref{fig:cWarcs}. Given this assumption, there is a well-defined inverse map   $\ds \vartheta_i = (\pi_z \circ \sigma_i)^{-1} \colon \fD_i \to L_i^-$ with   domain $\fD_i \subset \partial_h^- \mW$, which we recall was given in coordinates  in \eqref{eq-coordinatesVT} by
 \begin{equation}\label{eq-quadraticcurves}
\vartheta_i(r',\theta',-2)  = \left( R_{i,r'}(\theta'), \Theta_{i,r'}(\theta'),-2 \right) =   \left( r(\vartheta_i(r',  \theta',-2)), \theta(\vartheta_i(r', \theta',-2)),-2 \right) .
\end{equation}
 Also, recall that   Hypothesis~\ref{hyp-SRI}.\ref{item-SRI-4} assumes that $R_{i,2}(\theta') = r(\vartheta_i(2,  \theta',-2))$ has non-vanishing derivative, 
  except at $\theta_i'$   defined by  $\vartheta_i(2,\theta_i',-2)= (2,\theta_i,-2)$, and that  $\e_0> 0$ is the constant defined  after the property  \eqref{eq-generic1}.

For $i=1,2$, consider the curves $\g_{i,r}(\theta) =
\vartheta_i(r,\theta,-2)$,   defined for        $1 \leq r \leq 3$ and
$\theta$ such that $(r,\theta,-2) \in \fD_i$. 
The transversality assumption on $\sigma_i$ implies that each  curve $\g_{i,r}$ is non-singular; that is, 
$\frac{d}{d\theta} \g_{i,r}(\theta) \ne 0$.
 We impose a  hypothesis  on the shape of these curves for values of the coordinate $r$ near $2$, which implies that they  have parabolic shape, up to second order.

\begin{hyp}\label{hyp-parabolic}   
For  $i = 1,2$,   $2 \leq r_0 \leq 2 + \e_0$ and $\theta_i -\e_0 \leq \theta \leq \theta_i + \e_0$, assume that   
\begin{equation}\label{eq-quadratic}
\frac{d}{d\theta} \Theta_{i,r_0}(\theta) > 0 \quad, \quad  \frac{d^2}{d\theta^2} R_{i,r_0}(\theta) > 0  \quad, \quad \frac{d}{d\theta} R_{i,r_0}(\theta_i') = 0 
\end{equation}
where $\theta_i'$ satisfies $\vartheta_i(2,\theta_i',-2)= (2,\theta_i,-2)$. Thus for $2 \leq r_0 \leq 2+ \e_0$,  the graph of $R_{i,r_0}(\theta')$ is   parabolic with vertex       $\theta' = \theta_i'$. 
\end{hyp}
Note that \eqref{eq-quadratic}   is a consequence of Hypothesis~\ref{hyp-SRI} when  $r_0 = 2$.  
   
 \begin{defn}\label{def-generic}
 A Kuperberg flow $\Phi_t$ is \emph{generic}   if the construction of $\mK$ and $\cK$ satisfies Hypotheses~\ref{hyp-SRI}, \ref{hyp-genericW} and \ref{hyp-parabolic}. That is,  the singularities for the vanishing of the vertical vector field $\cW$ are of quadratic type, and the insertion yields a quadratic-type  radius function near the special points.
 \end{defn}

 Recall from Proposition~\ref{prop-wilsonproperties}   that for each   point $x \in \cR$ with $-1 < z(x) < 1$, its $\cW$-orbit is forward asymptotic to the periodic  orbit $\cO_2$ and backward asymptotic to the periodic orbit $\cO_1$.
The periodic orbits intersect $\bRt$ in the points $ \omega_i  =  \cO_i \cap \bRt$      for $i=1,2$.  
  Also, recall from \eqref{eq-intervals} that $I_{0}  = \cR \cap \bRt =  \left\{(2,\pi ,z) \mid  -1 \leq z \leq 1 \right\}$ is the line segment in $\bRt$ between $\omega_1$ and $\omega_2$. 
   Recall that the action of the generator $\psi \in \cGK^*$ defined by the Wilson flow preserves $I_0$, and for $\xi \in I_0$ with $-1 < z(\xi) <  1$ we have $z(\xi) < z(\psi(\xi))$. Then set
\begin{equation}\label{eq-funddomain}
I_{\xi}  =   \left\{(2,\pi ,z) \mid  z(\xi)   < z \leq   z(\psi(\xi)) \right\} .
\end{equation}

 \begin{lemma}\label{lem-density1}
Let $\xi \in I_0$ with $-1 < z(\xi) < z(\psi(\xi)) < 0$ and suppose that $I_{\xi} \subset \Sigma$,  then $\Sigma = \fM$.
\end{lemma}
\proof
The hypotheses imply that    $I_{\xi}$ is a fundamental domain for the $\Psi_t$-flow on the interior of $\cR$. That is, 
  the $\cW$-orbit of each interior point $x \in \cR$ intersects the interval $I_{\xi}$.  
If $I_{\xi} \subset \Sigma$ then $\cR' \subset \Sigma$, and hence  $\fM_0 \subset \Sigma$. As $\fM_0$ is dense in $\fM$,   the conclusion follows. 
\endproof

The strategy of the proof of Theorem~\ref{thm-density} is   to show that for some $\xi \in I_0$ sufficiently close to $\omega_1$, the $\cK$-orbit of $\omega_1$ contains $I_{\xi}$ in its closure.  To this end,      we establish some    estimates on the orbit  under $\cGK^*$ of   $\xi \in \bRt$ with     $2 \leq r(\xi) < 2 + \delta$ where   $\delta > 0$ is sufficiently small.
 
We first obtain estimates for the metric behavior of the orbits of the Wilson generator $\psi \in \cGK^*$.    Recall the functions $f$ and $g$ chosen in Section~\ref{sec-wilson}, which are constant in the coordinate $\theta$, with
\begin{equation}\label{eq-wilson2}
\cW =g(r, \theta, z)  \frac{\partial}{\partial  z} + f(r, \theta, z)  \frac{\partial}{\partial  \theta} ~ .
\end{equation}
  Hypothesis~\ref{hyp-genericW}  and condition \eqref{eq-generic1} imply there exists constants $A_g,B_g,C_g$ such that  the quadratic form  $Q_g(u,v) = A_g \cdot u^2 + 2B_g \cdot uv + C_g \cdot v^2$ defined by the Hessian of $g$ at $\omega_1$ is positive definite. As a consequence, 
for  $Q_0(r,z) = (r-2)^2 + (z+1)^2$,    there exists   $D_g > 0$ such that 
  \begin{equation}\label{eq-quadraticest1}
| g(r, \theta, z)  - Q_g(r-2, z+1) | ~ \leq ~   D_g \cdot (|r-2|^3 + |z+1|^3) \quad {\rm for} ~   Q_0(r,z)  \leq \e_0^2
\end{equation}
where $\e_0$ is the constant defined in \eqref{eq-generic1}. The condition \eqref{eq-quadraticest1} implies that for $(r,z)$ sufficiently close to $(2,-1)$,  the error term on the right-hand-side can be made arbitrarily small relative to the distance squared $Q_0(r,z)$ from the special point $(2,-1)$. We also observe that    \eqref{eq-quadraticest1}   implies there exists constants $0 < \lambda_1 \leq \lambda_2$ such that
  \begin{equation}\label{eq-quadraticest2}
\lambda_1 \cdot Q_0(r,z)  \leq  g(r, \theta, z) \leq \lambda_2 \cdot Q_0(r,z)   \quad {\rm for} ~    Q_0(r,z)  \leq \e_0^2 . 
\end{equation}

 Next, consider the action of the maps $\psi^{\ell}$ for $\ell > 0$.
   Let $\xi  \in \bRt$ with $2 \leq r(\xi) \leq 2 + \e_0$ and $-7/4 \leq z(\xi) \leq -1/4$, such that $\psi(\xi)$ is defined and   $z(\psi(\xi))< 0$.
   Let $T(\xi) > 0$ be defined by $\psi(\xi) = \Psi_{T(\xi) }(\xi)$. Then the $z$-coordinate of $\psi(\xi)$ is given by 
   \begin{equation}\label{eq-coordinates}
  z(\psi(\xi)) -  z(\xi) ~ = ~ \int_0^{T(\xi) } ~ g(\Psi_s(\xi)) ~ ds  ~ \geq ~ 0  .
\end{equation}
If $\xi \ne \omega_1$ then $g(\Psi_s(\xi))$ is positive along the orbit segment for $0 \leq s \leq T(\xi)$, hence  $ z(\psi(\xi)) -  z(\xi) > 0$.

Note that if the orbit of $x \in \mW$ avoids the $\e_0$-tube around the periodic orbits $\cO_i$ then $g \equiv 1$ along the orbit, so the $z$-coordinate along the orbit increases at constant rate 1. 
As the cylinder $\cC$ has height $4$, this means that such an orbit traverses an angle of at most $4$, hence it does not   complete a full turn around the cylinder $\cC$. In particular, for $\xi \in \bRt$ with $r(\xi) \geq 2 + \e_0$ and $z(\xi) < 0$ with   $\psi$   defined at $\xi$, this implies the $\cW$-orbit crosses the annulus $\cA$ where $z=0$, then returns to intercept $\bRt$ with $z(\psi(\xi)) > 0$. 
On the other hand, for $\xi \in \bRt$ with $2 < r(\xi) < 2+ \delta$ for $\delta \ll \e_0$ and $z(\xi) < -1$, the $\cW$-orbit traverses a region near $\cO_1$ where the slope is close to $0$, and hence repeatedly  traverses the rectangle, with the number of revolutions increasing   as $\delta \to 0$. In particular, for such $\xi$,  the powers $\psi^{\ell}(\xi)$ form a sequence of points in the vertical line $r= r(\xi)$ with increasing $z$-coordinates, as has been noted previously, especially in the proofs of Propositions~\ref{prop-wander} and \ref{prop-cluster-}.

We next  combine \eqref{eq-quadraticest1} with \eqref{eq-coordinates} to obtain   metric estimates on the   orbit of $\omega_1$ under the action of    $\cGK^*$. 
  Recall that the first transition point for the forward orbit of    $\omega_1$ is the special entry point  $p_1^{-} = \tau(\cL_1^- \cap \cO_1)$ with $r(p_1^-) = 2$, as   illustrated in   Figure~\ref{fig:KR}. 
Section~\ref{sec-proplevels} introduced the alternate notation  $p'(1) = p_1^- \in E_1$ and  $p(1) = \tau^{-1}(p'(1))  \in L_1^-$,  where     $r(p(1)) = 2$ by the Radius Inequality.

  The forward $\cW$-orbit of  $p(1)$  is trapped in the region $\cC \cap \{z < -1\}$, and thus  intercepts  $\cL_1^- \cap \cC$ in an infinite  sequence of points with   increasing $z$-coordinates between $-2$ and $-1$.  For all  $\ell \geq b$, 
 these points   are labeled $p'(1; 1 , \ell) \in  \cL_1^-$,  where  $r(p'(1; 1 , \ell)) =2$ and 
$\ds  z(p'(1; 1 , \ell)) <  z(p'(1; 1 , \ell +1)) < -1$. Moreover,  $z(p'(1; 1 , \ell)) \to -1$ as $\ell \to \infty$.

For each  $\ell \geq 0$,   set $p(1; 1 , \ell) = \sigma_1^{-1}(p'(1; 1 , \ell))  \in L_1^-$ then $r(p(1; 1 , \ell)) > 2$. Note that    $r(p(1; 1 , \ell)) \to 2$ as $\ell \to \infty$, and the sequence $p(1; 1 , \ell)$ accumulates    on  $p(1)$ in $L_1^-$.

\begin{figure}[!htbp]
\centering
{\includegraphics[height=70mm]{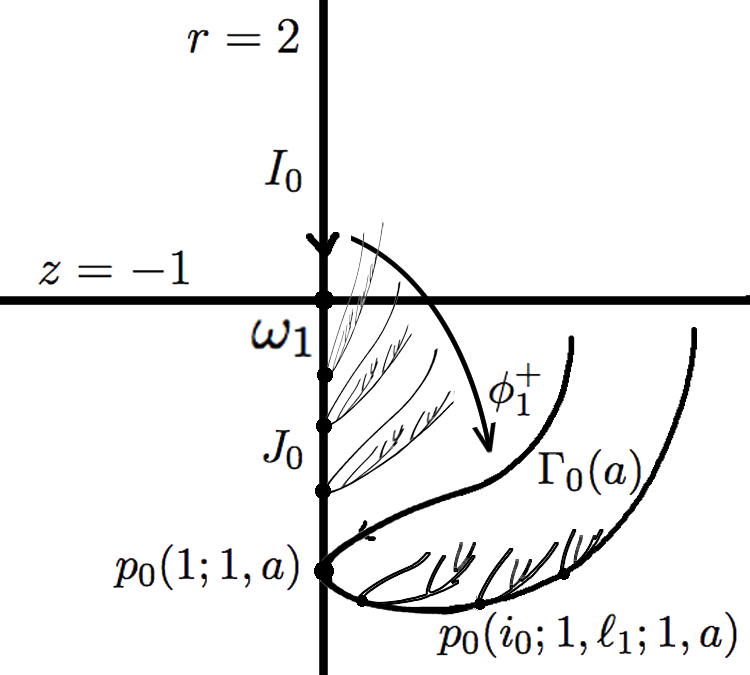}}
\caption{\label{fig:transform2}  Iterations of $\G_0$ and $\Lambda_0$ in $\bRt$ under  $\cGK^*$}
\vspace{-6pt}
\end{figure}

Recall from Section~\ref{sec-doublepropellers} that corresponding to the sequence   $\{ p'(1; 1 , \ell) \mid  \ell \geq 0\} \subset  \cL_1^-$ is a sequence $\{ p_0(1; 1 , \ell) \mid  \ell \geq 0\} \subset  \bRt$ where  $p_0(1; 1 , \ell)$  is the point in $\bRt$ whose  forward  $\cW$-orbit has $p'(1; 1 , \ell)$ as its first transition point. 
Thus,  $r(p_0(1; 1 , \ell) ) = 2$, and  $p_0(1; 1 , \ell +1) = \psi(p_0(1; 1 , \ell) )$ for $\ell \geq 1$, with  $z(p_0(1; 1 , \ell)) \to -1$
(see  Figure~\ref{fig:transform2}).   
Define  $0  < t_1 < t_2 < \cdots $   such that  $\Psi_{t_{\ell}}(p_0(1 ; 1 , 0)) = p_0(1 ; 1 , \ell)$.

\begin{lemma}\label{lem-density2}
There exists constants $C_2, C_3 > 0$ such that  for all $\ell \geq \ell_0$, where  $\ell_0>0$   is such that $z(p_0(1 ; 1 , \ell_0)) > -(1+\e_0)$, then we have:
  \begin{equation}\label{eq-spacing1}
  \frac{-1}{4\pi \lambda_1 \ell + C_2} <   z(p_0(1 ; 1 , \ell +1)) +1 \leq    \frac{-1}{4\pi \lambda_2 \ell + C_3} .
\end{equation}
Moreover, for the constants $\lambda_1 \leq \lambda_2$ introduced in \eqref{eq-quadraticest2}, there exists constants $C_4, C_5, C_6 > 0$   such that for all $\ell \geq \ell_0$, 
\begin{equation}\label{eq-spacing2}
\frac{C_4}{(4\pi \lambda_1 \ell )^2 + C_5}\leq  z(p_0(1 ; 1 , \ell +1))    -   z(p_0(1 ; 1 , \ell))  \leq \frac{C_4}{(4\pi \lambda_2 \ell )^2 + C_6}    .
\end{equation} 
\end{lemma}
 \proof
Note that $z(p_0(1 ; 1 , \ell_0)) > -(1+\e_0)$ implies the $\cW$-orbit
of $p_0(1;1,\ell_0)$ satisfies $$-(1+\e_0) < z(\Psi_t(p_0(1 ; 1 , \ell_0))) < -1 ~ {\rm  for ~ all} ~  t \geq 0, $$
and thus $f(\Psi_t(p_0(1 ; 1 , \ell_0))) = 1$. It then follows from \eqref{eq-wilson2} that $t_{\ell} = t_{\ell_0}  + 4\pi (\ell - \ell_0)$ for all $\ell \geq \ell_0$.

Set $z(t) = z(\Psi_t(p_0(1 ; 1 , \ell_0)))$ for $t \geq 0$, so that  $-(1+\e_0) < z_0 = z(p_0(1 ; 1 , \ell_0)) < -1$. Then we have   
\begin{equation}\label{eq-spacing3}
\lambda_1 \cdot \left(1+ z(t)  \right)^2  ~ \leq ~ \frac{dz}{dt} ~ \leq ~ \lambda_2 \cdot \left(1+ z(t)  \right)^2 ~. 
\end{equation}
As $z(0) = z_0$ set  $C_1 = -1/(z_0 +1) > 0$ then we have
\begin{equation}\label{eq-solution1}
  \frac{-1}{\lambda_1 t + C_1}   ~ \leq ~   z(t) +1 ~ \leq ~   \frac{-1}{\lambda_2   t + C_1 } ~.
\end{equation}
Substitute  $t_{\ell} = t_{\ell_0}  + 4\pi (\ell - \ell_0)$ into \eqref{eq-solution1} and collect constants to obtain \eqref{eq-spacing1}.

 The estimate \eqref{eq-spacing2} follows by subtracting terms in \eqref{eq-spacing1} for $\ell$ and $\ell+1$, and gathering constants.
 Alternately, one can use the formula \eqref{eq-quadraticest2} and the estimate \eqref{eq-quadraticest1} to obtain \eqref{eq-spacing2}.
 \endproof

The estimates in Lemma~\ref{lem-density2} have   counterparts for $\xi \in I_0$ with $-1 < z(\xi) < 0$ for which $z(\psi(\xi)) < -1  + \e_0$. In particular, the length of the interval $I_{\xi}$ defined in \eqref{eq-funddomain} has   quadratic upper and lower bound estimates as a function of $|z(\xi) +1|$.

We next obtain a more precise estimate than \eqref{eq-quadraticest1} for the restriction  of $g(2,\theta_0,z)$   to the line $I_0$.  

Note that for some $\lambda_1 \leq \lambda_g \leq \lambda_2$, the Taylor expansion for   $g(2,\theta_0,z)$   gives the approximation
$\ds g(2,\theta_0,z) = \lambda_g (z+1)^2 + O(|z+1|^3)$.  
Next, choose a scale $0< \e_3 <\e_0/2$ sufficiently small so that this approximation is strong enough to control the dynamics of the $\Psi_t$-flow near the cylinder $\cC$.

\medskip

 \begin{defn}\label{def-e3}
Let $\ds 0 < \e_3 \leq  \min \left\{\e_0/2 ,  1/100,  1/(300   \lambda_g ) \right\}$ be such that 
 \begin{equation}\label{eq-renormalizedest1}
|g(2,\theta_0,z) - \lambda_g (z+1)^2| ~ \leq ~  \frac{\lambda_g (z+1)^2}{100} \quad {\rm for} ~ |z+1| \leq \e_3 ~ .
\end{equation}
 \end{defn}
As a consequence, we have
\begin{equation}\label{eq-renormalizedest1a}
.99 \cdot \lambda_g (z+1)^2  \leq g(2,\theta_0,z) \leq  1.01 \cdot \lambda_g (z+1)^2  \quad {\rm for} ~ |z+1| \leq \e_3 ~ .
\end{equation}

 A key point is that we can now  choose a $\delta_3 > 0$ so that a weaker form of the  estimate \eqref{eq-renormalizedest1a} holds for some neighborhood of a fundamental domain $I_{\xi}$ in $I_0$ as defined by \eqref{eq-funddomain}. Let  $0 < \delta_3 < \e_0/2$ be such  that the following two estimates hold:
 \begin{eqnarray}
|r-2| \leq \delta_3 ~ {\rm and }~ |z+1| \leq \e_3    ~ & ~ \Longrightarrow ~&~ 0 ~ \leq ~ g(r,\theta_0,z) ~ \leq ~ 1.02 \cdot \lambda_g (\e_3)^2 \label{eq-boxest1}\\
|r-2| \leq \delta_3  ~ {\rm and }~ \e_3/4 \leq z+1 \leq \e_3  ~ & ~ \Longrightarrow ~&~ .98 \cdot \lambda_g (z+1)^2  \leq g(r,\theta_0,z) \leq  1.02 \cdot \lambda_g (z+1)^2 \label{eq-boxest2}
\end{eqnarray}
 Note that the assumption that $g$ is non-degenerate in Hypothesis~\ref{hyp-genericW} and the estimate \eqref{eq-quadraticest2} implies that for $\e_3$ sufficiently small,  there exists a constant $C_g > 0$ such that $\delta_3 =  \e_3/C_g$.

 We next  introduce the ``target box'',  which is key to the proof of Theorem~\ref{thm-density},
 \begin{equation}\label{eq-targetbox}
R_0(\delta_3, \e_3) = \left\{ (r,z) \in \bRt \mid |r-2| \leq \delta_3  ~ , ~ \e_3/4 \leq z+1 \leq \e_3 \right\} ~ .
\end{equation}

 \begin{defn}\label{def-R}
We say that  $\xi  \in \bRt$ and $\ell > 0$ satisfy \emph{Condition $R(\delta_3, \e_3, \ell)$} if 
\begin{equation}\label{eq-endzone1}
z(\xi) < -1 ~ , ~ 2 < r(\xi) \leq 2+\delta_3
\end{equation} 
\begin{equation}\label{eq-endzone2}
\e_3/3 < z(\psi^{\ell}(\xi)) +1 < z(\psi^{\ell +1}(\xi)) +1 < \e_3/2 ~ .
\end{equation}
 \end{defn}

We then have an  analog for $r > 2$ of Proposition~\ref{prop-wander}. 

 \begin{lemma}\label{lem-density3} 
 For $\delta_3, \e_3$ as above  and $\xi \in \bRt$ satisfying \eqref{eq-endzone1},   there exists $\ell > 0$ such that   \eqref{eq-endzone2} is satisfied.
  \end{lemma}
 \proof
As $z(\xi) < -1$ is assumed,  there is a unique  $t_1 > 0$ so that   $z(\Psi_{t_1}(\xi)) = -1 + \e_3/4$, and let $t_2 \geq t_1$ be the first subsequent time for which $\Psi_{t_2}(\xi) \in \bRt$. Note that  $t_2 - t_1 \leq 2\pi \, r(\xi) < 6 \pi$, and so by \eqref{eq-coordinates} and \eqref{eq-boxest1} and Definition~\ref{def-e3}, we obtain
\begin{equation}\label{eq-boundedstep1}
z(\Psi_{t_2}(\xi)) - z(\Psi_{t_1}(\xi)) ~ \leq ~  6 \pi \cdot 1.02 \cdot \lambda_g (\e_3)^2 ~ \leq ~ \e_3 \cdot 20/300 = \e_3/15 < \e_3/12
\end{equation}
 so that   $z(\Psi_{t_2}(\xi)) + 1< \e_3/3$.  The same estimates as in \eqref{eq-boundedstep1} show that 
 $$z(\psi^{\ell +1}(\Psi_{t_2}(\xi))) - z(\psi^{\ell}(\Psi_{t_2}(\xi))) \leq \e_3/15 < \e_3/12$$
 for all $\ell > 0$ for which $z(\psi^{\ell}(\Psi_{t_2}(\xi))) +1 < \e_3/2$. The estimate \eqref{eq-endzone2} then follows. 
 \endproof

The last metric estimate required for the action of $\psi$, and  a key point in the proof of Theorem~\ref{thm-density},  is to obtain estimates on   the dependence of $z(\psi^{\ell}(\xi))$ for small changes of $\xi \in \bRt$.  We  first estimate   $z(\psi^{\ell}(\xi+ \Delta z)) - z(\psi^{\ell}(\xi))$ for $\ell > 0$,  where     $\Delta z$ represents a  small vertical increment.

 Recall that  the Wilson flow $\Psi_t$ is rotationally invariant, and for $x_1 = (r_1, \theta_1, z_1)$ we  adopt the     notation $R_{\theta'}(x_1) = x_1 + \theta' = (r_1, \theta_1 + \theta', z_1)$,  so that the rotation  invariance becomes  $\Psi_t(x_1 +   \theta') = \Psi_t(x_1) +   \theta'$.

\begin{lemma}\label{lem-density4}
Let  $x_1 = (r_1, \theta_1, z_1)$ with $-7/4 \leq z_1 \leq -5/4$ and $|r_1-2|  \leq \e_0$, so that $f(x) = g(x) = 1$ for $x$ sufficiently close to $x_1$.  For  $\Delta z$ sufficiently small, set $x_1' = x_1 + \Delta z = (r_1, \theta_1, z_1 + \Delta z)$. For $t \geq 0$, set  $y_1 = \Psi_{t}(x_1)$,  then
\begin{equation}\label{eq-verticalslope}
z(\Psi_t(x_1')) - z(\Psi_t(x_1)) ~ \approx  ~ g(y_1)  \cdot \Delta z  
\end{equation}
\end{lemma}
 \proof
 Let $s_1$ be defined by $z(\Psi_{s_1}(x_1 ')) = z(x_1)$. 
  As the slope of the vector field $\cW$ near $x_1$ is $g(x_1)/f(x_1) = 1$, we have $s_1 =   -\Delta z$, and hence  $\theta(\Psi_{s_1}(x_1')) = \theta_1 -\Delta z$. 
Thus,  $x_1 = R_{s_1}(\Psi_{s_1}(x_1'))$ and so   
\begin{equation} \label{eq-verticalslopecalc}
z\left( \Psi_t(x_1') \right)   =    z\left( \Psi_t(  \Psi_{-s_1}(  R_{-s_1}(x_1))) \right) 
 =    z\left( R_{-s_1}( \Psi_{-s_1}(\Psi_{t}( x_1))) \right)  
  =    z\left(  \Psi_{-s_1}(y_1)\right)
\end{equation}
The slope of $\cW$ at $y_1$ is   $g(y_1)$ as $f(y_1) = 1$, so that  $z\left(  \Psi_{-s_1}(y_1)\right) \approx  z(y_1) + g(y_1) \cdot  \Delta z$, 
   which yields \eqref{eq-verticalslope}.
 \endproof

  In the case where  $x_1 = \xi \in \bRt$ and $\Psi_{t}$ is the map defining  $\psi^{\ell}$,  then 
the  estimate \eqref{eq-verticalslope} can   be made    precise.  
For    $\ell \geq 0$,  set $\xi_{\ell} = \psi^{\ell}(\xi)$. 
Define $T(\xi_0,\ell) > 0$ by  $\psi^{\ell}(\xi) = \Psi_{T(\xi,\ell)}(\xi)$.

 \begin{cor}\label{cor-density4}
Let $\xi_0  \in \bRt$ with $-7/4 \leq z(\xi_0) \leq -1-\e_0$ 
and $\ell > 0$ satisfy \emph{Condition $R(\delta_3, \e_3,\ell)$}.
 Then for $\Delta z$ sufficiently small so that $\xi_0' = \xi_0 + \Delta z$ again satisfies \emph{Condition $R(\delta_3, \e_3, \ell)$}, we have
 \begin{equation}\label{eq-zincreasedz}
{  |} z(\psi^{\ell}(\xi_0')) - z(\psi^{\ell}(\xi_0)) { |} \leq    1.02 \cdot  \lambda_g (\e_3)^2 \cdot \Delta z  
\end{equation}
 \end{cor}
 \proof
The point $\xi_0$ satisfies the  conditions of Lemma~\ref{lem-density4}, and then use \eqref{eq-boxest2} in the proof.
\endproof

\medskip
 
   We  next develop an estimate  for  $z(\psi^{\ell}(\xi_0')) - z(\psi^{\ell}(\xi_0))$,  where $\xi_0' = \xi_0 +    \Delta r$.
   
     Assume that  $\xi_0  \in \bRt$ with $-7/4 \leq z(\xi_0) \leq -1-\e_0$ 
and $\ell > 0$ satisfy condition $R(\delta_3, \e_3, \ell)$.
 Also assume that    $\Delta r$ is sufficiently small so that $\xi_0' = \xi_0 + \Delta r$ again satisfies condition $R(\delta_3, \e_3, \ell)$. 
 
  Then the $z$-coordinate $z(\Psi_{t}(\xi_0))$  increases at constant rate $1$ in the region $\{z \leq -1-\e_0\}$, and subsequently  increases at a possibly slower rate in the region $\{ -1 - \e_0 < z < -1+\e_0\}$. 
  The $\cW$-orbit segment $\{ \Psi_t(\xi_0) \mid 0 \leq t \leq T(\xi_0, \ell)\}$ is geometrically a ``coiled spring'' in $\mW$ which wraps $\ell$-times around the cylinder $\{r=r(\xi)\}$, with the coils most ``tightly spaced'' near the periodic orbit $\cO_1$. If $\Delta r > 0$ then the orbit 
$\{ \Psi_t(\xi_0') \mid 0 \leq t \leq T(\xi_0', \ell)\}$ is slightly less tightly spaced, so the ends of the orbit segment are   spaced further apart. 
More precisely, we have 
  \begin{equation}\label{eq-zvariation1}
  z( \psi^{\ell}(\xi_0))  =  z(\xi_0) ~ + ~  \int_{0}^{T(\xi_0, \ell)} ~   g(\Psi_s(\xi_0)) ~ ds 
\end{equation}
Obtaining a sharp estimate for $z(\psi^{\ell}(\xi_0')) - z(\psi^{\ell}(\xi_0))$ requires estimating this integral as  $\Delta r$ varies in the expression $r(\xi_0') = r(\xi_0) +  \Delta r$, and this requires more detailed estimates of the integrand $g(\Psi_s(\xi_0))$ than is given. However, a sufficient estimate can be obtained by noting that 
$z( \psi^{\ell}(\xi_0))$ is a smooth function with bounded derivatives  on any compact set $K \subset \bRt$ contained in the domain of $\psi^{\ell}$. For $0 < \delta_4 < \delta_3$,  define:
\begin{equation}\label{eq-partialboundregion}
K(\e_3,\delta_3, \delta_4) ~ = ~  \{ (r,z) \mid -2 \leq z \leq -1-\e_3 ~{\rm and} ~ 2 +\delta_4 \leq r \leq 2 + \delta_3\}
\end{equation}
which is a compact set contained in the domain of $\psi^{\ell}$. Then for the function $g$ chosen, define:
\begin{equation}\label{eq-partialbound}
M(g,\ell, \e_3,\delta_3, \delta_4) ~ = ~ \max \left\{ \frac{\partial \{ z( \psi^{\ell}(r,z))\}}{\partial r} \mid   (r,z) \in K(\e_3,\delta_3, \delta_4) \right\}
\end{equation}
Then by the Mean Value Theorem, we have:
 \begin{lemma}\label{lem-density5}
Let  $\xi_0  \in \bRt$ and $\Delta r$ be such that   $\xi_0, \xi_0' \in K(\e_3,\delta_3, \delta_4)$, where $\xi_0' = \xi_0 + \Delta r$. Then
 \begin{equation}\label{eq-zincreasedr}
\left| z(\psi^{\ell}(\xi_0')) - z(\psi^{\ell}(\xi_0)) \right|  ~ \leq ~ M(g,\ell, \e_3,\delta_3, \delta_4) \cdot \Delta r
\end{equation}
\end{lemma}

  The idea of the   proof of Theorem~\ref{thm-density},   is to show that for $\e_3 > 0$ as in Definition~\ref{def-e3}, and for any $\delta_3$ satisfying \eqref{eq-boxest1} and \eqref{eq-boxest2}, the 
 orbits of $\omega_1$ under $\cGK^*$ yields a set of points which have dense $z$-coordinates in a rectangle defined by \eqref{eq-boxest2}.   As $\delta_3$ can be chosen arbitrarily small, this   implies that  the closure of the orbit of $\omega_1$ contains a fundamental domain   for the $\Psi_t$-flow on $\cR$  as in \eqref{eq-funddomain}, so the claim that $\Sigma = \fM$ follows from Lemma~\ref{lem-density1}.   The derivations of the estimates used in the proof of these statements depend fundamentally on Hypothesis~\ref{hyp-parabolic}.

The $\cK$-orbit of $\omega_1 \in \bRt$ defines the sequence of points 
$\{ p'(1; 1 , \ell)   \mid  \ell \geq 0\} \subset  \cL_1^-$ for $\ell \geq 0$, and   corresponding points  
  $\{ p_0(1; 1 , \ell) \mid  \ell \geq 0\} \subset  \bRt$  with $r(p_0(1; 1 , \ell) ) = 2$  and $z(p_0(1 ; 1 , \ell)) \to -1$,  
where   the convergence of $z(p_0(1 ; 1 , \ell))$   is estimated by Lemma~\ref{lem-density2}.

Assume that $a=b=0$ for $a$ and $b$ as defined in Sections~\ref{sec-doublepropellers} and
\ref{sec-proplevels}, respectively.
For each $\ell \geq 0$, we then have the curve $ \kappa_0(\ell)
\subset \bRt$  with lower endpoint $p_0(1 ; 1 , \ell)$. 
For $\ell =0$, the lower part of the curve $\kappa_0(0)$ is the image of the segment $J_0 \subset \bRt$ under the map $\phi_1^+ \in \cGK^*$, and $\kappa_0(\ell)$ is the image of $\kappa_0(0)$ under the map $\psi^{\ell}$ for $\ell \geq 0$. 

The curve $\kappa_0(0)$ is ``parabolic'' by
Hypothesis~\ref{hyp-parabolic}. That is, the $r$-coordinate of the
graph is approximated by a quadratic function of the $z$-coordinate,
for $z$ near $p_0(1 ; 1 , 0)$. In particular, the $r$-value of points
on $\kappa_0(0)$ increases as $-(1+z)$ increases for $z \leq -1$.
As the map $\psi$ preserves the radius coordinate, the same monotonicity property is true for each of the curves $\kappa_0(\ell)$, as illustrated in the graphs in  Figure~\ref{fig:curvesR0}.

Note that the  image under the map  $\phi_1^+$ of the points $\{p_0(1; 1 , \ell) \mid  \ell \geq 0\}$ 
yields a sequence of points  $p_0(1 ; 1 , \ell ; 1, 0) \subset \kappa_0( 0)$ 
for which $r(p_0(1 ; 1 , \ell ; 1, 0)) > 2$ and $p_0(1 ; 1 , \ell ; 1,
0) \to p_0(1;1,0)$ as $\ell \to \infty$. Thus $r(p_0(1;1,\ell;1,0))\to
2$ as $\ell\to \infty$.

Choose $m_1 > 0$ such that  $r(p_0(1;1,m_1;1,0)) < 2+ \delta_3$.   
By Lemma~\ref{lem-density3}, there exists $\ell_1 > 0$ such that 
\begin{equation}\label{eq-targetedrange}
\e_3/3 < z(p_0(1;1,m_1;1,\ell_1)) +1 < z(p_0(1;1,m_1;1,\ell_1+1)) +1 < \e_3/2 ~ ,
\end{equation}
since $p_0(1;1,m_1;1,\ell_1)=\psi^{\ell_1}(p_0(1;1,m_1;1,0))$.

Choose $n_1 > m_1$ such that    $r(p_0(1;1,n_1;1,0)) \leq r(p_0(1;1,m_1;1,\ell_1))$, and  set $\delta_4 = r(p_0(1;1,n_1;1,0)) -2$.

Observe that the $r$-coordinate along   the   curve $\kappa_0(0)$ between $p_0(1;1,n_1;1,0)$ and $p_0(1;1,m_1;1,0)$ is monotone increasing by Hypothesis~\ref{hyp-parabolic}, hence have $r$ values ranging between $2 + \delta_4$ and $2+\delta_3$. We next select a collection of sequences of points in the $\cGK^*$-orbit of $\omega_1$ which are sufficiently closely spaced, and which ``shadow'' this curve segment. Applying the map $\psi^{\ell_1}$ will then yield a collection of points in the region  $R_0(\delta_3, \e_3)$ which have arbitrarily dense $z$-coordinates. 

\begin{figure}[!htbp]
\centering
{\includegraphics[height=70mm]{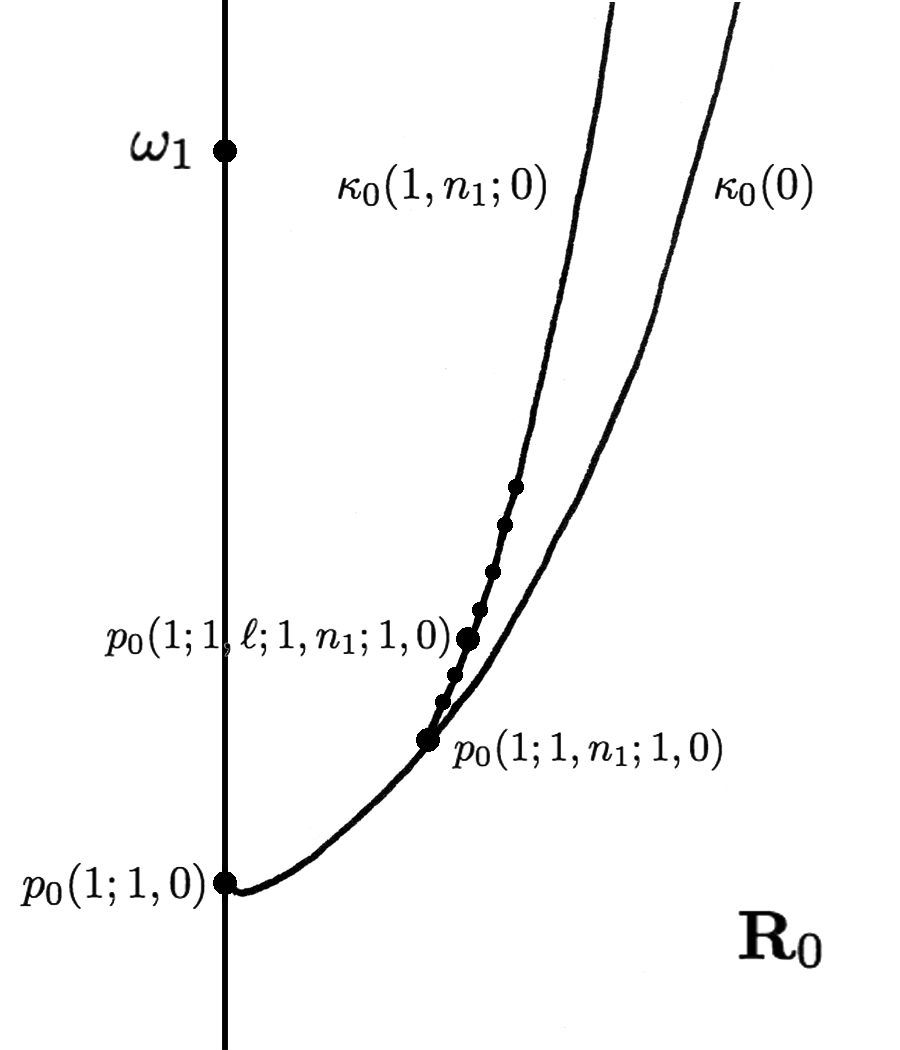}}
\caption[{Points   in the   $\kappa_0$ curve}]{\label{fig:kappa1ell}  The points $p_0(1;1,\ell;1,n_1;1,0)$  in the curve $\kappa_0(1,n_1;0)$ in $\bRt$}
\vspace{-6pt}
\end{figure}

Introduce the constant
\begin{equation}
\mu = 2 \cdot \max ~ \left\{ M(g,\ell_1, \e_3,\delta_3, \delta_4)  ~ , ~  1.02 \cdot  \lambda_g (\e_3)^2 \right\}
\end{equation}
where  $M(g,\ell_1, \e_3,\delta_3, \delta_4)$ is defined by   \ref{eq-partialbound}, and the second term is introduced in \ref{eq-zincreasedz}.
Then $\mu > 2$, and for   $N \geq 2$, set 
\begin{equation}
\e_N = \frac{z(p_0(1;1,m_1;1,\ell_1+1)) - z(p_0(1;1,m_1;1,\ell_1))}{N} < 1/8 \quad , \quad  \delta_N = \e_N /\mu  .
\end{equation}
 The image under $\phi_1^+$ of the lower part of the parabolic curve
 $\kappa_0(n_1)$ with lower endpoint $p_0(1;1,n_1)$ is part of the
 curve $\kappa_0(1,n_1; 0)$ with lower endpoint $p_0(1;1,n_1;
 1,0)$. Thus, the lower part of $\kappa_0(1,n_1; 0)$ is  the image of
 the interval $J_0$ under the composition $\phi_1^+ \circ \psi^{n_1}
 \circ \phi_1^+$. Applying this map to the sequence
 $p_0(1;1,\ell)= \psi^{\ell}(p_0(1;1,0))$ yields a collection of
 points on the curve $\kappa_0(1,n_1; 0)$ which converge to $p_0(1;1,
 n_1;1,0)$, as in Figure~\ref{fig:kappa1ell}. 
 By Lemma~\ref{lem-density2}, there exists $n_2 > 0$ so that for the collection
 \begin{equation}\label{eq-Omega2}
 O(\delta_N ; n_1,n_2) = \left\{ \phi_1^+ \circ \psi^{n_1} \circ \phi_1^+ \circ \psi^{\ell} (p_0(1;1,0)) \mid \ell \geq n_2 \right\} \subset \kappa_0(1,n_1;  0) , 
\end{equation}
 the differences of both the $z$ and $r$-values   of successive points
 for $\ell, \ell +1 \geq n_2$ are bounded above by $\delta_N$. By
 making $n_2$ bigger, we can assume that the point $p_0(1;1,n_2; 1,n_1;1,0)$  has the largest $r$-coordinate   for all  points in $O(\delta_N ; n_1,n_2)$.

 We now proceed inductively. Assume that the integers $\{n_1, n_2, \ldots n_k\}$ have been chosen, with associated sets 
 $O(\delta_N ; n_1,n_2, \ldots, n_i) \subset \cGK^*(\omega_1)$    for
 $i \leq k$, so that $p_0(1;1,n_i; 1,n_{i-1};\ldots ; 1,n_1 ; 1,0)$ 
satisfies that 
$$r(p_0(1;1,n_i; 1,n_{i-1};\ldots ; 1,n_1 ; 1,0))$$ 
is the maximum of the $r$-values for all points in  $O(\delta_N ; n_1,n_2, \ldots, n_i)$.
 
 Suppose that    $r(p_0(1;1,n_k; 1,n_{k-1};\ldots ; 1,n_1 ; 1,0)) <
 r(p_0(1;1,m_1;1,0))$, then the previous points chosen lie inside the
 region defined by the region bounded by the parabola 
 $\G_0(1,n_1; 0)$ and   the region $\{r \leq r(p_0(1;1,m_1;1,0))\}$ which lies inside the $\e_0$-ball in $\bRt$ about $\omega_1$. Thus, we can continue the above process, and   choose $n_{k+1}$ as follows.  
 The curve  $\kappa_0(1, n_k; 1,n_{k-1};\ldots; 1, n_1;0)$   with
 lower endpoint $p_0(1;1,n_k;1,n_{k-1};\cdots;1,0)$   is the image of $J_0$ under the composition of maps
 \begin{equation}\label{eq-compositionk}
\phi_1^+ \circ \psi^{n_1} \circ \phi_1^+ \circ \psi^{n_2} \circ \phi_1 \circ \cdots \circ   \psi^{n_k} \circ \phi_1^+
\end{equation}
 Thus, applying the  map  in \eqref{eq-compositionk}  to the sequence
 $\{ p_0(1;1,\ell) \mid \ell \geq 0\}$ yields a collection of points on the curve   $\kappa_0(1,n_k; \ldots; 1, n_1; 0)$ which converge to $p_0(1;1,n_k;1,n_{k-1};\cdots;1,0)$.  
 Then by Lemma~\ref{lem-density2}, there exists $n_{k+1} > 0$ so that for the collection
 \begin{eqnarray} 
 \lefteqn{O(\delta_N ; n_1,n_2, \ldots, n_{k+1}) }  \label{eq-Omega2k+1} \\
 & = &  \left\{ \phi_1^+ \circ \psi^{n_1} \circ \phi_1^+ \circ \psi^{n_2} \circ \phi_1 \circ \cdots \circ   \psi^{n_k} \circ \phi_1^+ (p_0(1;1,\ell)) \mid \ell \geq n_{k+1} \right\} \nonumber
\end{eqnarray}
which is contained in the curve $\ds  \kappa_0(1,n_k; \ldots; 1, n_1; 0)$,  
the differences of both the $z$ and $r$-values   of successive points
for $\ell, \ell +1 \geq n_{k+1}$ are bounded above by
$\delta_N$. We can assume that $r(p_0(1;1,n_{k+1};1,n_k;\cdots;1,n_1;0))$ is the maximum of the $r$-values for all points in  $O(\delta_N ; n_1,n_2, \ldots, n_{k+1})$.

The collections of points  defined by \eqref{eq-Omega2k+1} are
contained in the region bounded by the parabola  
$\G_0(1,n_1; 0)$.
 Hypothesis~\ref{hyp-parabolic} implies the composite functions in
 \eqref{eq-compositionk}   have   non-zero derivatives which are
 bounded away from zero for points in the region $\delta_4 \leq r(\xi)
 \leq \delta_3$ and $z(\xi) \leq -1-\e_0$. If
 $r(p_0(1;1,n_{k+1};1,n_k;\cdots;1,n_1;0) < r(p_0(1;1,m_1;1,0))$ then
 the $\kappa_0$ curve with lower endpoint at 
 $p_0(1;1,n_{k+1};1,n_k;\cdots;1,n_1;0)$ is defined, and so we can repeat the process to increase $r(p_0(1;1,n_{k+1};1,n_k;\cdots;1,n_1;0))$, unless    
 $$r(p_0(1;1,n_{k+1};1,n_k;\cdots;1,n_1;0)) \geq r(p_0(1;1,m_1;1,0)),$$
  in which case the  inductive selection of   points terminates.

Now apply the map $\psi^{\ell_1}$ to each of the sets $O(\delta_N ; n_1,n_2, \ldots, n_{k})$ to obtain points in the region  $R_0(\delta_3, \e_3)$.
Then Corollary~\ref{cor-density4} and  Lemma~\ref{lem-density5} imply the   difference of the $z$-values of successive points in the sequence are bounded above by $\mu \cdot \delta_N = \e_N$.  Finally, by construction the upper $r$-value of points in $O(\delta_N ; n_1,n_2, \ldots, n_{k})$ is approached by the descending $r$-values of points in  $O(\delta_N ; n_1,n_2, \ldots, n_{k+1})$.  

Thus, the image of $\cGK^*(\omega_1)$ contains points which are in the strip $2 < r \leq 2+ \delta_3$ and whose $z$-values increase from 
$z(p_0(1;1,m_1;1,\ell_1))$ to $z(p_0(1;1,m_1;1,\ell_1+1))$ in increments at most $\e_N$. As $\delta_3 > 0$ and $N \gg 0$ were arbitrary, this implies the closure of $\cGK^*(\omega_1)$ contains $\cR$, as was to be shown to establish Theorem~\ref{thm-density}. \hfill $\Box$

We conclude with a   remark on the proof of Theorem~\ref{thm-density}. The quadratic assumption in Hypothesis~\ref{hyp-genericW}  was not used in any essential way, except to obtain precise estimates on the spacing of the orbits of $\cW$. It seems like that a much weaker assumption than this hypothesis should suffice for the proof.

  \bigskip
 
\section{Geometry of curves in $\bRt$}\label{sec-geometry}
 
The aim of this section is to investigate the geometry of the
collection of curves formed by the intersection $\fM_0 \cap \bRt$,
thus of $\g_0$ and $\lambda_0$ curves. 
Up to now,  we have only used their topology and the fact that the $\Gamma_0$ and $\Lambda_0$ curves form
families of nested ellipses in $\bRt$. For a generic Kuperberg flow, as defined by Definition~\ref{def-generic}, we obtain estimates on the geometry of these curves, including the following  result which has applications  in Sections~\ref{sec-zippered} and \ref{sec-growth}:

\begin{thm}\label{thm-boundedlength}
Let $\Phi_t$ be a generic Kuperberg flow on $\mK$. Then there exists
$L > 0$ such that each connected component of $\fM_0 \cap \bRt$ is a 
curve with length bounded above by $L$.
\end{thm}

For $\g_0$ and $\lambda_0$ curves of level 1, the $r$-coordinate along
the lower half of the curve is monotone increasing. If such condition
was true for any curve in $\fM_0\cap \bRt$, the theorem will be
trivial. As we explain below this monotonicity condition does not
holds for curves at level 2 or higher. To describe how $\g_0$ and $\lambda_0$ curves
fold, we introduce a new set of curves, the $A_0$-curves.

For simplicity, we will restrict our discussion to $\g_0$ curves. Recall that the union of $\g_0$ and $\kappa_0$ curves
gives $\Gamma_0$ curves,  that at high levels are very  thin (as a
consequence of the nesting property). Hence the shape of a $\G_0$
curve, that is a closed curve, can be seen as the shape of a
simple curve. These are the $A_0$-curves described below.

As mentioned above, the lower half of a $\g_0$ curve of level at least
2 has no monotone radius. Then, the curve in $E_1$
generating the corresponding propeller has no monotone radius which
implies that the propeller can have internal notches that do not
correspond to the ones introduced in
Section~\ref{sec-proplevels}. These new internal
notches generate also compact surfaces, to which we refer as bubbles
of second type (the bubbles of first type being the ones introduced in
Section~\ref{sec-bubbles}). This phenomena is discussed at the end of
this section.

The two aspects we are interested in are the
following. First, we show that the lengths of $\g_0$ and $\lambda_0$ curves having
one endpoint in the domain of $\phi_1^+$ and the other endpoint in the
domain of $\phi_2^+$ has an upper bound. Observe that such curves admit
$2-2\e_0$ as a trivial lower bound. Second, these curves are ``nice
flat'' arcs away from their endpoints, in the sense that if we
consider a propeller generated by a $\gamma$-curve truncated near its
tip, the trace that we obtain on the annulus $\cA=\{z=0\}$ will be a spiral
that is never tangent to the lines $\{\theta=const\}$.

To prove Theorem~\ref{thm-boundedlength}, let us start by analyzing level 1 curves. Recall that the curves
$\Gamma_0(\ell)$ for $\ell\geq a$ and unbounded, form the trace of the
propeller $\tau(P_\Gamma')$ on $\bRt$. Hypothesis~\ref{hyp-SRI} implies
that $\Gamma\subset L_1^-$  is a parabolic curve with vertex at $p(1)=\sigma_1^{-1}(p_1^-)$
as illustrated in Figure~\ref{fig:Gamma1}. Thus the curves
$\Gamma_0(\ell)$ are tangent to the vertical line $r=2$ at two points,
their vertices, 
with parabolic shape near these points, as in
Figure~\ref{fig:ellipsespropeller}.

The points in $\Gamma_0(\ell)\subset \mW$ belong to Wilson orbits that
intersect at least $a+\ell$ times $\bRt$, thus orbits that turn at
least $a+\ell-1$ times around the cylinder $\cC$ before hitting $\cA$. It follows that the radius coordinate along these curves, for $\ell>a$, is
bounded above by $2+\e_0$. Thus the curves are contained in the
rectangle $[2,2+\e_0]\times [-2,2]$. Moreover, outside the $\e_0$ tube
around the periodic orbits $\cO_i$, for $i=1,2$, we have that $g=1$
and thus a Wilson orbit intersects at most twice $\bRt$ in this
region, once in each direction. Hence  the curves $\Gamma_0(\ell)$
are almost vertical outside the $\e_0$-neighborhood of the points
$\omega_i$. These considerations apply to higher level $\Gamma_0$-curves.

\begin{figure}[!htbp]
\centering
{\includegraphics[height=70mm]{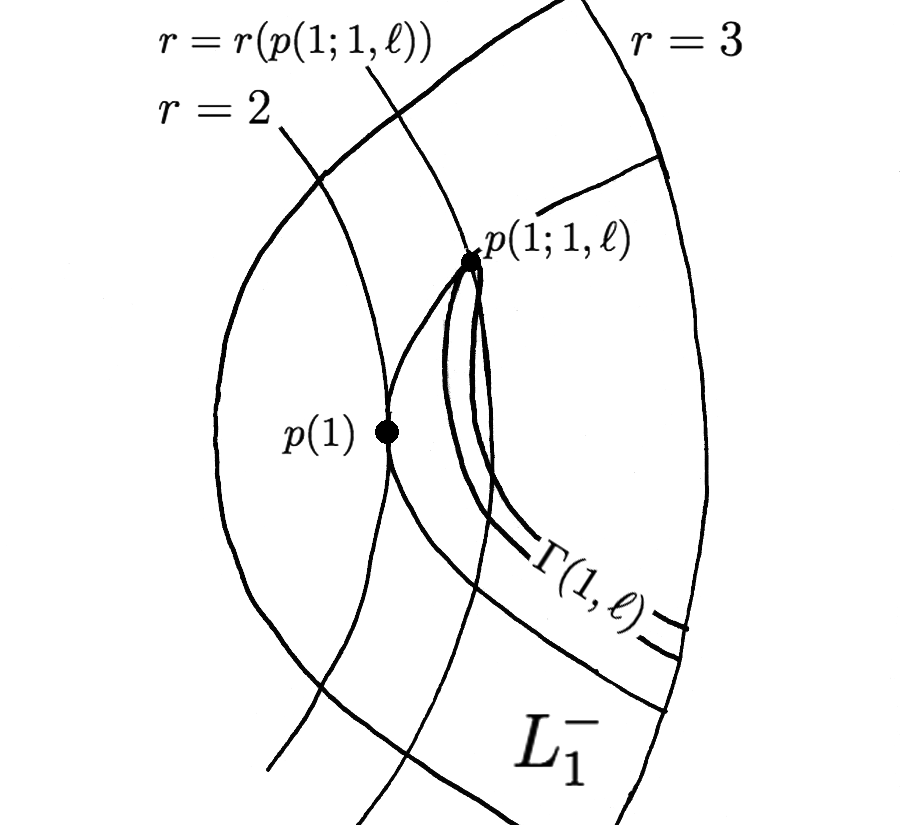}}
\caption[{Intersection of $\Gamma(1,\ell)$ with the circle $r=r(p(1;1,\ell))$}]{\label{fig:nonmonotoneradius} Intersection of $\Gamma(1,\ell)$ with the circle $r=r(p(1;1,\ell))$ in $L_1^-$}
\vspace{-6pt}
\end{figure}

The problem of describing higher level $\Gamma_0$ curves comes from
the fact that the curves $\Gamma(i,\ell)\subset L_i^-$ for $i=1,2$ do not have
monotone radius along $\gamma(i,\ell)$ and the same might apply to
the curves $\kappa(i,\ell)$. To illustrate this situation, let us
consider the curve $\Gamma_0(\ell)$ for $\ell\geq 0$. The hypothesis
$\ell\geq 0$ implies that the positive orbit of its points hit $E_1$
along $\tau(\Gamma(1,\ell))$, that is tangent to the vertical line
$r=2$. Thus $\Gamma(1,\ell)\in L_1^-$ is tangent to $\Gamma$ at
$p(1;1,\ell)$ and is composed of the union of $\g(1,\ell)$ and
$\kappa(1,\ell)$. The circle with radius $r=r(p(1;1,\ell))$ is
transverse to $\gamma$ at two points, one being $p(1;1,\ell)$. Thus
$\gamma(1,\ell)$ can not have monotone radius, as illustrated in
Figure~\ref{fig:nonmonotoneradius}.

Observe that Hypothesis~\ref{hyp-parabolic} implies
that a curve in $\cL_1^-$ crossing the line $z=-1$ is bended by the
inverse of the
insertion map $\sigma_1^{-1}$ at the point with $z$-coordinate equal to $-1$. If the original curve is
vertical, we obtain a parabolic curve in $L_1^-$. For example, $\Gamma'\subset \cL_1^-$ is the vertical straight
line and becomes the parabolic curve $\Gamma\in L_1^-$ with axis
$\sigma_1^{-1}(\{z=-1\})$. 

Let $A'=\{z=-1\}\subset \cL_1^-$, then the
curve $A=\sigma_1^{-1}(A')\subset L_1^-$ is a straight line connecting
the boundary of $\partial_h^-\mW$ to the point $p(1)$ with
radius 2. Observe that $A$ is inside the region of $L_1^-$ that is
bounded by $\Gamma$.

The propeller $P_A$ is an infinite propeller as in
Definition~\ref{def-infpropeller}. Consider the notched propeller
$P_A'=P_A\cap \mW'$ and its image $\tau(P_A')$ in $\mK$. The trace of
this propeller on the rectangle $\bRt$ is an infinite family of arcs
$A_0(\ell)$, each one of them being in the region of $\bRt$ that is
bounded by $\Gamma_0(\ell)$, for the same index $\ell$. Since $A$
intersects $\Gamma$ at the special point $p(1)\in L_1^-$, for each
$\ell$ the curve $A_0(\ell)$ intersects $\Gamma_0(\ell)$ at the two
points $p_0(1;1,\ell)$ and $p_0(1;2,\ell)$. Recall that these are the
points where $\Gamma_0(\ell)$ is tangent to the vertical line $r=2$.

From the curve $A\subset L_1^-$ we can construct the infinite family of $A$-curves and
the corresponding
\mbox{$A_0$-curves} in $\bRt$. Observe that a level $n$ curve in the later family
is contained in the region of $\bRt$ bounded by a $\Gamma_0$-curve at
the same level, in fact, by the $\Gamma_0$-curve having exactly
the same label. We will say that the $\Gamma_0$-curve {\it envelops} the
$A_0$-curve. To be more precise $\Gamma_0(i_1, \ell_1;\cdots ;\ell_n)$
envelops $A_0(i_1, \ell_1;\cdots ;\ell_n)$ and these two curves
intersect at the points $p_0(1;i_1, \ell_1;\cdots ; 1,\ell_n)$ and
$p_0(1;i_1, \ell_1;\cdots ;2,\ell_n)$. Finally, observe that at the
intersection points 
the curves are transverse to each
other and locally $\Gamma_0(i_1, \ell_1;\cdots ;\ell_n)$ has parabolic
shape with axis $A_0(i_1, \ell_1;\cdots ;\ell_n)$.

 As we comment at the beginning of the section we are interested in
 finding upper bounds to the length of $\g_0$-curves and describing
 the end parts of these curves. The length of such curves is
 approximated by the length of $A_0$-curves, thus it is enough to find
 a uniform upper bound to the length of $A_0$-curves.

We consider next $A$-curves in detail. Let $A'(i,\ell)\subset \cL_i^-$
be the entry curves of the notches of the propeller
$P_A'$ and $A(i,\ell)=\sigma_i^{-1}(A'(i,\ell))\subset
L_i^-$, for $i=1,2$. Hypothesis~\ref{hyp-genericW}, on the Wilson flow, implies
that for $r>2+\e_0$ the vertical component of the vector field is
equal to 1 and thus the lines $A_0(\ell)$ are steep in this region. 
Then there exists some positive number $B$ such that for $\ell\geq B$
the curve $A'(1,\ell)$ intersects transversally $A'$. Let $q'(1;1,\ell)$ be the
intersection point. 
The points in the curve $A'(1,\ell)$ belong to $\cW$-orbits that
intersect at least $(a+\ell)$-times the rectangle $\bRt$ before
intersecting the   annulus $\cA$ and for $\ell>0$ the endpoint
$p'(1;1,\ell)\in\cL_1^-$. The genericity assumptions imply thus that
the value of $B$ is small. Assume, without loss of generality, that
$B=b$ for $b$ as defined in Section~\ref{sec-proplevels}.

\begin{figure}[!htbp]
\centering
\begin{subfigure}[c]{0.4\textwidth}{\includegraphics[height=70mm]{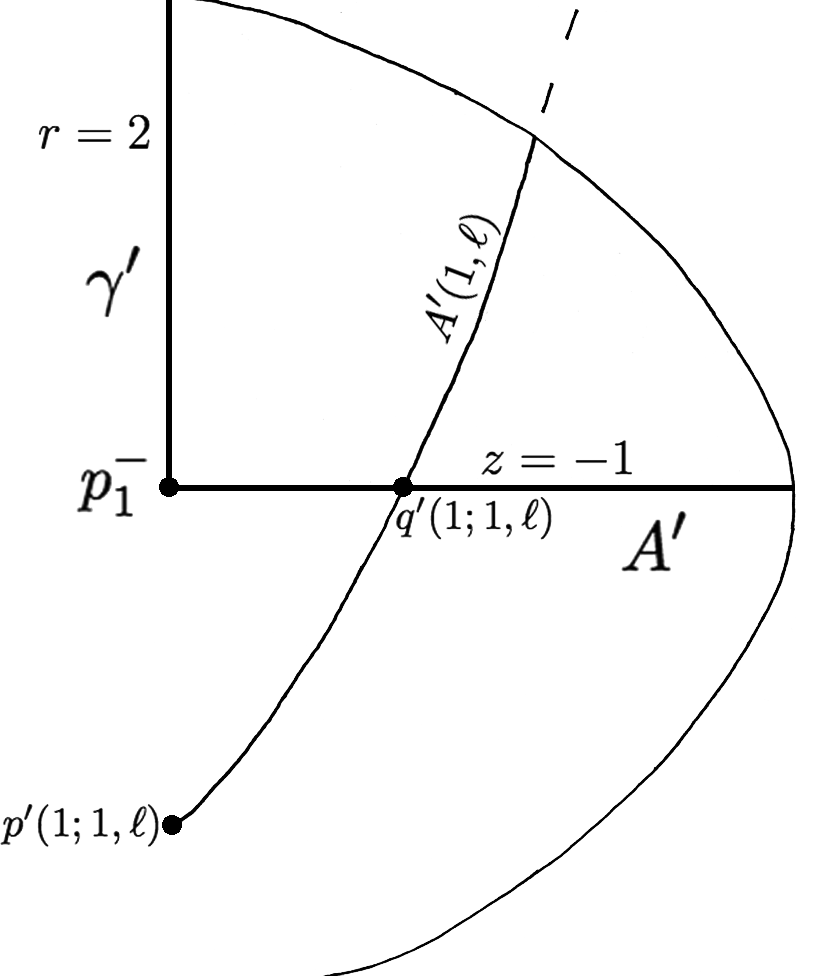}}\end{subfigure}
\begin{subfigure}[c]{0.4\textwidth}{\includegraphics[height=70mm]{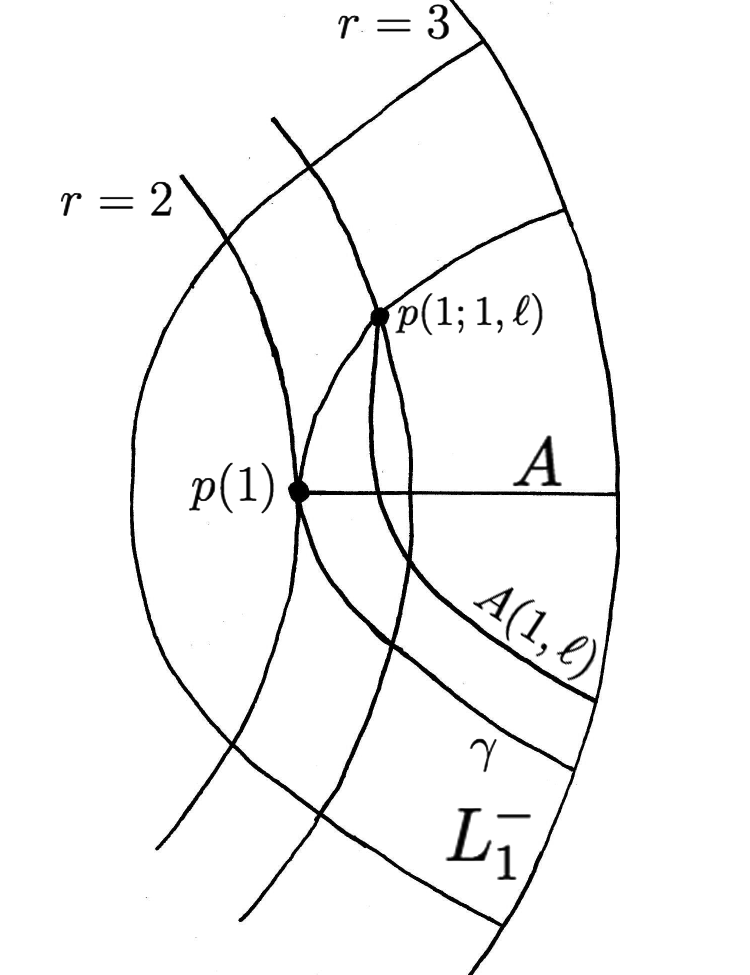}}\end{subfigure}
\caption{\label{fig:firstfolding} Curves $A'$ and $A'(1,\ell)$ in  $\cL_1^-$and the corresponding curves in $L_1^-$}
 \vspace{-6pt}
\end{figure}

Thus, for $\ell \geq 0$,  $A(1,\ell)\subset L_1^-$
is a curve going from the point $p(1;1,\ell)\in \kappa$ to a point in
the boundary of $\partial_h^-\mW$ and crossing $A$ at
$q(1;1,\ell)=\sigma_1^{-1}(q'(1;1,\ell))$, then $A(1,\ell)$ is folded
  at the point $q(1;1,\ell)$, as illustrated in
  Figure~\ref{fig:firstfolding}. Observe that $q(1;1,\ell)$ is not
  necessarily the point of $A(1,\ell)$ with minimum radius, the point
  with minimum radius of $A(1,\ell)$ is in the arc between
  $q(1;1,\ell)$ and $p(1;1,\ell)$.

Since $A$ and $A(1,\ell)$ intersect, the propellers $P_A$ and
$P_{A(1,\ell)}$ intersect along the $\cW$-orbit of the point
$q(1;1,\ell)$. Observe that $P_{A(1,\ell)}$ is a finite propeller. 

\begin{figure}[!htbp]
\centering
{\includegraphics[height=100mm]{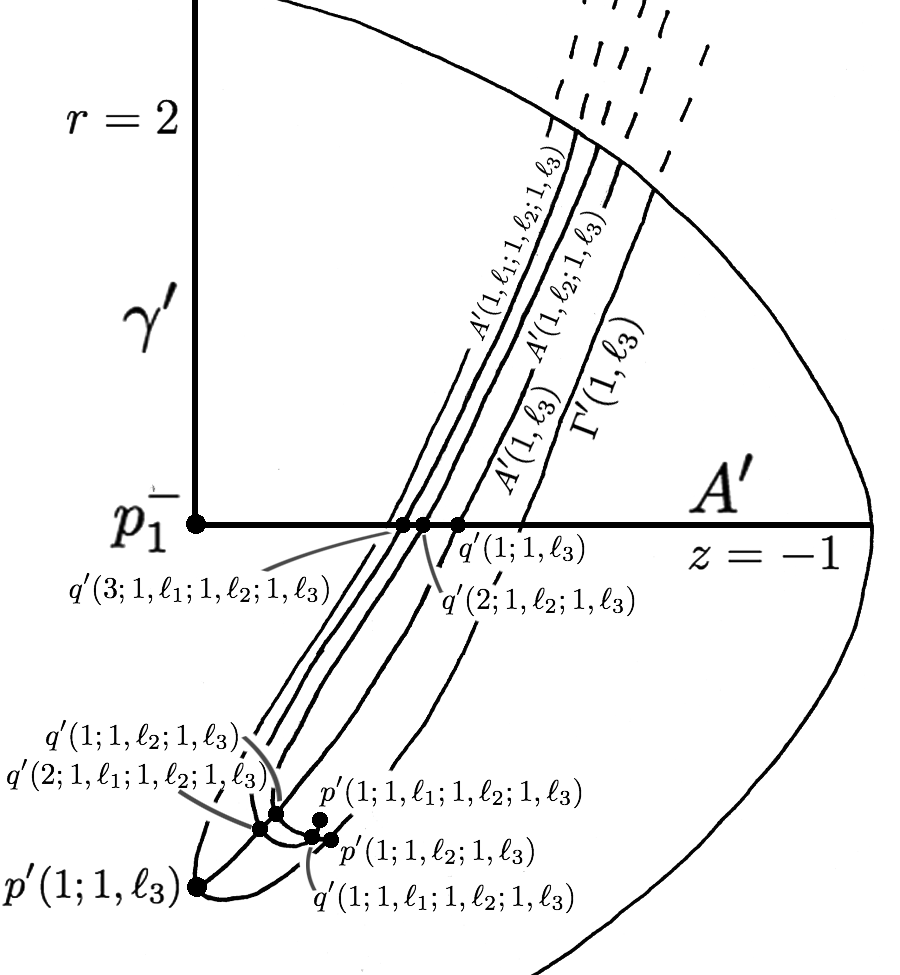}}
\caption{\label{fig:qprimepoints} Three levels of $A'$-curves inside the region bounded by $\Gamma'(1,\ell_3)$ in $\cL_1^-$}
\vspace{-6pt}
\end{figure}

We consider now the intersection of $P_{A(1,\ell_1)}$ with
$\cL_1^-$, for $\ell_1 \geq b$.
Since $P_{A(1,\ell_1)}$ is a
finite propeller the intersection consists of finitely many curves
$A'(1,\ell_1;1,\ell_2)$ for $\ell_1$ fixed and $\ell_2\geq b$ and
bounded. For $\ell_1\geq 0$ and $\ell_2\geq 0$, these curves have one endpoint at
$p'(1;1,\ell_1;1,\ell_2)\subset \kappa'(1,\ell_2)$ and its other
endpoint on the boundary of $\cL_1^-$. Moreover
$A'(1,\ell_1;1,\ell_2)$ is in the region bounded by
$\Gamma'(1,\ell_1;1,\ell_2)$, the later being in the region bounded by
$\Gamma'(1,\ell_2)$, by the results of the previous
sections. For $\ell_1$ big enough and $\ell_2\geq 0$ the curve
$A'(1,\ell_1;1,\ell_2)$ intersects $A'$ and at the intersection the
curves are transverse to each other. Let
$q'(2;1,\ell_1;1,\ell_2)$ be the intersection point. Observe that $A'(1,\ell_1;1,\ell_2)$ also intersects transversally
$A'(1,\ell_2)$ at the point $q'(1;1,\ell_1;1,\ell_2)$ that is in the $\cW$-orbit of
$q(1;1,\ell_1)$.
The points $q'(1;1,\ell_2)=A'(1,\ell_2)\cap A'$ and
$q'(2;1,\ell_1;1,\ell_2)$ are on $A'$ and in the region bounded by
$\G'(\ell_2)$, as illustrated in
Figure~\ref{fig:qprimepoints}. 

We comment on the assumption $\ell_1$ big enough. The curve
$A'(1,\ell_1;1,\ell_2)$ intersects $A'$ if its endpoint
$p'(1;1,\ell_1;1,\ell_2)$ has $z$-coordinate smaller than $-1$,
condition that is satisfied for $\ell_1$ big enough since the sequence
$p'(1;1,\ell_1;1,\ell_2)$ converges to $p'(1;1,\ell_2)$ as $\ell_1\to
\infty$ and $z(p'(1;1,\ell_2))<-1$ for any $\ell_2\geq a$. Observe that since we want to describe curves at high levels,
we just care for curves with $\ell_1$ big. This assumption will be
used several times in what follows.

\begin{remark}\label{rmk-endpointA0}
$A_0$-curves at level 1 have their lower endpoint below the line
$z=-1$, since the lower endpoint belongs to a $\cW$-orbit with radius
coordinate equal to 2. At higher levels this is not true. However, we
will only consider $A_0$-curves whose lower endpoint is below
$z=-1$. The reason for this is that a higher level $A_0$-curve with lower endpoint above $z=-1$
  intersects the domain of the map $\psi^{-1}$. Hence the trace of
  the propeller containing this $A_0$-curve has curves whose lower
  endpoint is below $z=-1$. The latter curve is longer, thus an upper
  bound to the length of the longer curve  suffices.
\end{remark}

\begin{lemma}\label{lem-radiusqpoints}
For a generic $\mK$ plug, $r(q'(2;1,\ell_1;1,\ell_2))<r(q'(1;1,\ell_2))$.
\end{lemma}

\proof
Observe that
$$
r(q'(1;1,\ell_1;1,\ell_2))=r(q(1;1,\ell_1))\leq r(p(1;1,\ell_1))=r(p'(1;1,\ell_1;1,\ell_2),
$$ 
where the two equalities follow since the points are in the same
$\cW$-orbits and the inequality in the middle follows from
Hypothesis~\ref{hyp-parabolic}. Thus Hypothesis~\ref{hyp-genericW}
implies that the $\cW$-orbit of $p(1;1,\ell_1)$ climbs faster than the
$\cW$-orbit of $q(1;1,\ell_1)$ and thus
$$z(q'(1;1,\ell_1;1,\ell_2))\leq
z(p'(1;1,\ell_1;1,\ell_2)).$$ 
For $\ell_1$ big enough and $\ell_2$ bounded, these two
$z$-coordinates are less than $-1$. Thus the intersection point of
$A'(1,\ell_1;1,\ell_2)$ and $A'(1,\ell_2)$ has $z$-coordinate less
than $-1$, implying that $r(q'(2;1,\ell_1;1,\ell_2))<r(q'(1;1,\ell_2))$.
\endproof

Recall that $A'(1,\ell_1;1,\ell_2)$ intersects transversally $\kappa'(1,\ell_2)\subset\Gamma'(1,\ell_2)$
at the point $p'(1;1,\ell_1;1,\ell_2)$. Since the intersection of $A'(1,\ell_1;1,\ell_2)$ with $A'(1,\ell_2)$
is also transverse, $A'(1,\ell_1;1,\ell_2)$ bends at the intersection
point  
$q'(1;1,\ell_1;1,\ell_2)$. This implies that the curve
$\Gamma'(1,\ell_1;1,\ell_2)$, that envelops $A'(1,\ell_1;1,\ell_2)$
and is tangent at $p'(1;1,\ell_1;1,\ell_2)$ to $\kappa'(1,\ell_2)$,
bends at its intersection with $A'(1,\ell_2)$. The intersection between
$\Gamma'(1,\ell_1;1,\ell_2)$ and $A'(1,\ell_2)$ is transverse and
consists of one or two
points: $\g'(1,\ell_1;1,\ell_2)$ always intersects $A'(1,\ell_2)$ and 
$\kappa'(1,\ell_1;1,\ell_2)$ might intersect it to.

Observe that if we fix the index $\ell_2$ and we let $\ell_1\to
\infty$, both sequences $q'(1;1,\ell_1;1,\ell_2)$ and
$p'(1;1,\ell_1;1,\ell_2)$ converge to the point $p'(1;1,\ell_2)$.

We now iterate the process above, considering only the first
insertion for simplicity. The idea is that a level $n$ curve in $\cL_1^-$, with $\ell_1$-index big enough, intersects
a $A'$-curve at each of the previous levels, and these points are the
bending points of the curve. 

\begin{lemma}
For $\ell_1$ sufficiently large, the curve $A'(1,\ell_1;1,\ell_2;\cdots;1,\ell_n)$
intersects $A'$ and the curves $A'(1,\ell_{n-k};\cdots;1,\ell_n)$ for all
$1\leq k<n$. 
\end{lemma}

For
$1\leq k<n$, the
intersection point of $A'(1,\ell_1;1,\ell_2;\cdots;1,\ell_n)$ with $A'(1,\ell_{n-k};\cdots;1,\ell_n)$ is denoted
by $q'(k;1,\ell_1;\cdots;1,\ell_n)$, and the intersection with $A'$ is
denoted by $q'(n;1,\ell_1;\cdots;1,\ell_n)$. Observe that the first
index in the label of $q$-points is the difference of levels between
the $A'$-curves that pass through it.

\proof
For $\ell_1$ big enough we can assume that the endpoint
$p'(1;1,\ell_1;\cdots;1,\ell_n)$ of
$A'(1,\ell_1;1,\ell_2;\cdots;1,\ell_n)$ is in $\cL_1^-$ and has $z$-coordinate smaller
than $-1$, thus the curve intersects $A'$.

We now proceed by induction on $n$. For $n=1$ the conclusion is
straight forward. Assume that the lemma holds up to $n-1$ and consider
the case $n$. For $k=1$ observe that
$A'(1,\ell_1;1,\ell_2;\cdots;1,\ell_n)$ intersects $\kappa'(1, \ell_2;
\cdots;1,\ell_n)$ at $p'(1;1,\ell_1;\cdots;1,\ell_n)$. Following the
proof of
Lemma~\ref{lem-radiusqpoints} we have that 
$$r(q'(n;1,\ell_1;\cdots;1,\ell_n))<r(q'(n-1;1,\ell_2;\cdots;1,\ell_n)).$$
Thus at the line $z=-1$ the curve
$A'(1,\ell_1;1,\ell_2;\cdots;1,\ell_n)$ is on the left of the curve
$A'(1,\ell_2;\cdots;1,\ell_n)$, while near its endpoint it is on the
right. These curves intersect.

By hypothesis, the level $n-1$ curve $A'(1,\ell_2;\cdots;1,\ell_n)$
intersects all the curves of smaller levels. The same holds for the
curve $\Gamma'(1,\ell_2;\cdots;1,\ell_n)$; that is, 
$\Gamma'(1,\ell_2;\cdots;1,\ell_n)$ intersects the curves $A'(1,\ell_{n-k};\cdots;1,\ell_n)$ for all
$2\leq k<n$. Since 
$A'(1,\ell_1;1,\ell_2;\cdots;1,\ell_n)$ is in the region bounded by
$\Gamma(1,\ell_2;\cdots;1,\ell_n)$, we obtain that it must intersect
all the curves $A'(1,\ell_{n-k};\cdots;1,\ell_n)$.
\endproof

Now we consider the level 3 case in
some detail. For $\ell_1$ big enough and $\ell_2$ bounded, the curve
$A(1,\ell_1;1,\ell_2)=\sigma_1^{-1}(A'(1,\ell_1;1,\ell_2))\subset
L_1^-$  intersects transversally:
\begin{itemize}
\item $A(1,\ell_2)$ at the point $q(1;1,\ell_1;1,\ell_2)$.
\item $A$  at the point $q(2;1,\ell_1;1,\ell_2)$.
\end{itemize}
Thus the finite propeller $P_{A(1,\ell_1;1,\ell_2)}$ intersects:
\begin{itemize}
\item $P_{A(1,\ell_2)}$ along the $\cW$-orbit of $q(1;1,\ell_1;1,\ell_2)$.
\item $P_A$ along the $\cW$-orbit of $q(2;1,\ell_1;1,\ell_2)$.
\end{itemize}
Consider the trace of $P_{A(1,\ell_1;1,\ell_2)}$ on $\cL_1^-$, that
forms a finite collection of curves $A'(1,\ell_1;1,\ell_2;1,\ell_3)$
with $\ell_3\geq b$ and bounded. For $\ell_3\geq 0$ we have that:
\begin{itemize}
\item $A'(1,\ell_1;1,\ell_2;1,\ell_3)$ is contained in the region
  bounded by $\Gamma'(1,\ell_1;1,\ell_2;1,\ell_3)$, which is\linebreak contained
  in the region bounded by $\Gamma'(1,\ell_2;1,\ell_3)$, that is in
  the region bounded by $\Gamma'(1,\ell_3)$, by the nesting property.
\item $A'(1,\ell_1;1,\ell_2;1,\ell_3)$ intersects
  $A'(1,\ell_2;1,\ell_3)$ at $q'(1;1,\ell_1;1,\ell_2;1,\ell_3)$;
  $A'(1,\ell_3)$ at\linebreak $q'(2;1,\ell_1;1,\ell_2;1,\ell_3)$ and $A'$ at
  $q'(3;1,\ell_1;1,\ell_2;1,\ell_3)$.
Applying the proof of Lemma~\ref{lem-radiusqpoints} we have that 
$$r(q'(3;1,\ell_1;1,\ell_2;1,\ell_3))<r(q'(2;1,\ell_2;1,\ell_3))<r(q'(1;1,\ell_3)).$$
\end{itemize}
Thus the curve $A'(1,\ell_1; 1,\ell_2;1,\ell_3)$ folds twice: once at
$q'(1;1,\ell_1;1,\ell_2;1,\ell_3)$ and once at\linebreak
$q'(2;1,\ell_1;1,\ell_2;1,\ell_3)$. We deduce that the corresponding
$A_0$-curve in $\bRt$, $A_0(1,\ell_1; 1,\ell_2;\ell_3)$, folds twice near each of
its endpoints due to the symmetry condition on the Wilson flow. Recall that $A_0(1,\ell_1; 1,\ell_2;\ell_3)\subset
\bRt$ is the curve containing points whose forward orbits intersect
$E_1$ along $\tau(A(1,\ell_1;1,\ell_2;1,\ell_3))$. 

Iterating this analysis we get that each $A'$-curve at level $n$ folds $(n-1)$-times, and thus the lower end part of the corresponding $A_0$-curve coils in a ``small''
spiral turning in the counter-clockwise direction.

We can now start the proof
of Theorem~\ref{thm-boundedlength}, that is a direct consequence of
the following proposition.

\begin{prop}\label{prop-lengthgcurves}
The length of an $A_0$-curve is
upper bounded by the length of $A_0(a)+D$ for some uniform constant $D$.
\end{prop}

\proof
As explained in Remark~\ref{rmk-endpointA0} we can restrict to curves
whose lower endpoint is below $z=-1$. 
We start by comparing the length of $A_0(1,\ell_1;\ell_2)$ with the
length of $A_0(\ell_2)$. The idea is to prove first that 
$$\mbox{length}(A_0(1,\ell_1;\ell_2))<\mbox{length}(A_0(\ell_2))+D,$$
that is clearly upper bounded by the length of $A_0(a)+D$.

The curve $A'(1,\ell_1;1,\ell_2)$ in $\cL_1^-$ can be decomposed into
two parts: the one between $p'(1;1,\ell_1;1,\ell_2)$ and
$q'(1;1,\ell_1;1,\ell_2)$ and the rest. Observe that the later part
has approximately the same length as the segment of $A'(1,\ell_2)$
that lies above (in the $z$-direction) the point
$q'(1;1,\ell_1;1,\ell_2)$. The reason for this is that both curves are
in the region bounded by $\G'(1,\ell_2)$ and have no bending points
above $q'(1;1,\ell_1;1,\ell_2)$.

We can
thus bound the length of  $A'(1,\ell_1;1,\ell_2)$ by the sum of the
length of $A'(1,\ell_2)$ plus the distance between the points 
$p'(1;1,\ell_1;1,\ell_2)$ and
$q'(1;1,\ell_1;1,\ell_2)$. We compare the distances in $\bRt$.
Translating the discussion above to $\bRt$, we can bound the length of
$A_0(1,\ell_1;\ell_2)$ by the sum of the length of $A_0(\ell_2)$ plus
twice the distance between $p_0(1;1,\ell_1;1,\ell_2)$ and
$q_0(1;1,\ell_1;1,\ell_2)$. The multiplication by two is due
to the symmetry of these curves with respect to the line
$\{z=0\}$. Let
$\varphi_{\ell_2}=\psi^{a+\ell_2}\circ\phi_1^+\in \cGK^*$, then using
Lemmas~\ref{lem-psi_action}, \ref{lem-phi+_action} and \ref{lem-actionpoints} we obtain
\begin{equation}\label{eqn-difference}
d_{\bRt}(p_0(1;1,\ell_1;1,\ell_2),q_0(1;1,\ell_1;1,\ell_2))  =  d_{\bRt}(\varphi_{\ell_2}(p_0(1;1,\ell_1)),\varphi_{\ell_2}(q_0(1;1,\ell_1))).
\end{equation}

The next lemma is the key ingredient in the proof. Recall that
$U_{\phi_1^+}$ is the domain of the map $\phi_1^+$. Let 
$\cD'(1,\ell)$ be the points in
the region bounded by $\G_0(\ell)$ that belong to  $U_{\phi_1^+}\cap
\{z\leq -1\}$. Define 
$$\cD(1,\ell)=\{x\in \cD'(1,\ell)\,|\, \varphi_\ell(x)\in
\cD'(1,\ell)\},$$
for $\varphi_\ell=\psi^{(a+\ell)}\circ \phi_1^+$.

\begin{lemma}\label{lem-constantdistance}
For a generic Kuperberg plug $\mK$, the norm of the differential of $\varphi_\ell=\psi^{(a+\ell)}\circ \phi_1^+$
is less than 1 for points in $\cD(1,\ell)$.
\end{lemma}

\proof
Observe that $\psi^{-(a+\ell)}(\cD(1,\ell))$ is contained in the
region of $\bRt$ bounded by $\G_0(a)$. The map $\psi$ preserves
the $r$-coordinate. 
The map $\phi_1^+$ maps $\cD(1,\ell)$ to the interior of the region
bounded by $\G_0(a)$, and thus
$$
\mbox{Image}(\phi_1^+(\cD(1,\ell)))\subset
\mbox{Image}(\psi^{-(a+\ell)}(\cD(1,\ell))).
$$
For every $\ell'\geq \ell$ such that the region $\cD(1,\ell')\neq
\emptyset$, we have that 
$$
\mbox{Image}(\phi_1^+(\cD(1,\ell')))\subset
\mbox{Image}(\psi^{-(a+\ell)}(\cD(1,\ell))).
$$
Thus for points in $\cD(1,\ell)$ the map $\phi_1^+$ contracts
distances more than $\psi^{-(a+\ell)}$.  Hence, the norm of the differential of $\phi_1^+$ is smaller than the
norm of the differential of $\psi^{-(a+\ell')}$ for any $\ell'$ such
that $\cD(1,\ell')$ is non-empty. Hence, for $\xi\in
\cD(1,\ell)$ we have that 
\begin{eqnarray*}
\|D\varphi_\ell(\xi)\|& \leq & \|D\psi^{(a+\ell)}(\phi_1^+(\xi))\|\times
\|D\phi_1^+(\xi)\|\\
&< &\|D\psi^{(a+\ell)}(\phi_1^+(\xi))\|\times
\|D\psi^{-(a+\ell)}(\psi^{a+\ell}\circ\phi_1^+(\xi))\|\\
& =&\|D\psi^{(a+\ell)}(\eta)\|\times
\|D\psi^{-(a+\ell)}(\psi^{a+\ell}(\eta))\|=1.
\end{eqnarray*}
where $\eta=\phi_1^+(\xi)$. The last equality uses the fact that
$\psi$ preserves the $r$-coordinate.
\endproof

Observe that for $\xi\in\cD(1,\ell)$ and $\ell'\geq \ell$ we have also
that $\|D\varphi_{\ell'}(\xi)\|<1$, provided
$z(\varphi_{\ell'}(\xi))\leq -1$. 

We come back to bounding the length of the level 2 curve
$A_0(1,\ell_1;\ell_2)$. From \eqref{eqn-difference} and the previous lemma, we
obtain a constant $C<1$ such that 
$$
d_{\bRt}(p_0(1;1,\ell_1;1,\ell_2),q_0(1;1,\ell_1;1,\ell_2))
\leq C\,d_{\bRt}(p_0(1;1,\ell_1),q_0(1;1,\ell_1)).
$$
Let $a_1=\max_{\ell_1\geq 0} \{d_{\bRt}(p_0(1;1,\ell_1),q_0(1;1,\ell_1))\}$
then, 
$$\mbox{length}(A_0(1,\ell_1;\ell_2))\leq \mbox{length}(A_0(\ell_2))+2Ca_1,$$
concluding the proof of Proposition~\ref{prop-lengthgcurves} for level
2 curves for any $D\geq 2Ca_1$. 

We know consider the level 3 case. Again the curve
$A'(1,\ell_1;1,\ell_2;1,\ell_3)\subset \cL_1^-$ can be decomposed into
two parts: the one between $p'(1;1,\ell_1;1,\ell_2;1,\ell_3)$ and
$q'(1;1,\ell_1;1,\ell_2;1,\ell_3)$ and the rest. The later part has
approximately the same length as the part of $A'(1,\ell_2;1,\ell_3)$
that lies above the point $q'(1;1,\ell_1;1,\ell_2;1,\ell_3)$. Thus the
length of $A'(1,\ell_1;1,\ell_2;1,\ell_3)$ is upper bounded by the
length of $A'(1,\ell_2;1,\ell_3)$ plus the distance between the points $p'(1;1,\ell_1;1,\ell_2;1,\ell_3)$ and
$q'(1;1,\ell_1;1,\ell_2;1,\ell_3)$.

In $\bRt$,   the
length of $A_0(1,\ell_1;1,\ell_2;\ell_3)$ is   bounded above by the
length of $A_0(1,\ell_2;\ell_3)$ plus twice the distance between the points $p_0(1;1,\ell_1;1,\ell_2;1,\ell_3)$ and
$q_0(1;1,\ell_1;1,\ell_2;1,\ell_3)$. Observe that for $\varphi_{\ell_3}=\psi^{(a+\ell_3)}\circ
\phi_1^+$ we have that
$p_0(1;1,\ell_1;1,\ell_2;1,\ell_3)=\varphi_{\ell_3}(p_0(1;1,\ell_1;1;\ell_2))$
and $q_0(1;1,\ell_1;1,\ell_2;1,\ell_3) =\varphi_{\ell_3}(q_0(1;1,\ell_1;1;\ell_2))$.
By construction the point $p_0(1;1,\ell_1;1,\ell_2)$ belongs to the
  curve $\kappa_0(\ell_2)\subset \G_0(\ell_2)$ and $q_0(1;1,\ell_1;1,\ell_2)$ belongs to
  $A_0(\ell_2)$.
Thus these two points are in  $\cD(1,\ell_2)$. By
Lemma~\ref{lem-constantdistance} and  $C<1$ as above, we obtain
\begin{eqnarray*}
d_{\bRt}\left(p_0(1;1,\ell_1;1,\ell_2;1,\ell_3),q_0(1;1,\ell_1;1,\ell_2;1,\ell_3)\right)
& \leq & C\,d_{\bRt}(p_0(1;1,\ell_1;1,\ell_2),q_0(1;1,\ell_1;1,\ell_2))\\
& \leq & C^2\,d_{\bRt}(p_0(1;1,\ell_1),q_0(1;1,\ell_1))\\
& \leq & C^2a_1.
\end{eqnarray*}
Thus, 
\begin{eqnarray*}
\mbox{length}(A_0(1,\ell_1;1,\ell_2;\ell_3))& <
&\mbox{length}(A_0(1,\ell_2;\ell_3)) + \\
& &+2d_{\bRt}(p_0(1;1,\ell_1;1,\ell_2;1,\ell_3),q_0(1;1,\ell_1;1,\ell_2;1,\ell_3))\\
& \leq & \mbox{length}(A(\ell_3))+2Ca_1+2C^2a_1\\
& = & \mbox{length}(A(\ell_3))+2Ca_1(1+C).
\end{eqnarray*}
The first inequality follows just by decomposing the lower half of
$A_0(1,\ell_1;1,\ell_2;\ell_3)$ it two parts and using the symmetry of
this curve with respect to $\{z=0\}$. The second inequality follows
from the estimation of the length of level 2 curves and the
computation above.

Iterating this argument for level 4 curves we have that
\begin{eqnarray*}
\mbox{length}(A_0(1,\ell_1;1,\ell_2;1,\ell_3;\ell_4))& <
&\mbox{length}(A_0(1,\ell_2;1,\ell_3;\ell_4))+\\
&
&+2d_{\bRt}(p_0(1;1,\ell_1;1,\ell_2;1,\ell_3;1,\ell_4),q_0(1;1,\ell_1;1,\ell_2;1,\ell_3;1,\ell_4))\\
& <&
\mbox{length}(A_0(\ell_4)) +2Ca_1(1+C)+\\
& & + 2Cd_{\bRt}(p_0(1;1,\ell_1;1,\ell_2;1,\ell_3),q_0(1;1,\ell_1;1,\ell_2;1,\ell_3))\\
& \leq & \mbox{length}(A_0(\ell_4)) +2Ca_1(1+C)+2C^3a_1\\
& = & \mbox{length}(A_0(\ell_4)) +2Ca_1(1+C+C^2).
\end{eqnarray*}

Generalizing to a level $n$ curve, we get 
\begin{eqnarray}
\mbox{length}(A_0(1,\ell_1;\cdots;1,\ell_n)) & < &
\mbox{length}(A_0(\ell_n))+2Ca_1(1+C+\cdots +C^{n-2}) \label{eq-approxlength}\\
& = &
\mbox{length}(A_0(\ell_n))+2a_1\left(\frac{C-C^{n}}{1-C}\right). \nonumber
\end{eqnarray}

The last amount is upper bounded by $\frac{C}{1-C}$, and thus taking
$D=max\{2a_1\left(\frac{C}{1-C}\right), 2a_1\}$ we obtain the desired bound.
\endproof

This finishes the proof of Theorem~\ref{thm-boundedlength}.
We can now describe a more accurate picture of the trace on $\bRt$ of
a level $n\geq 2$ propeller in $\fM_0$ and thus of the trace of $\fM_0$. In
Figure~\ref{fig:arcspropeller} we had made the assumption that the
curve generating a finite propeller has monotone radius and thus that
the longest $\cW$-orbit in the propeller is the orbit of the endpoint of the
curve. As described above, the propellers in $\fM_0$ of level at least
2 are not generated by curves with monotone radius.

Consider a curve
$\g(i_1,\ell_1;\cdots;i_n,\ell_n)\subset L_{i_n}^-$ with $n\geq 2$,
$i_k=1,2$ for $1\leq k\leq n$ and
non monotone
radius. Once more, we restrict the discussion to $\g$-curves for
simplicity, but it applies also to $\lambda$-curves. 

Let $q^\g(n;i_1,\ell_1;\cdots;i_n,\ell_n)$ be the point in
$\g(i_1,\ell_1;\cdots;i_n,\ell_n)$ with smaller radius. 
The trace of the propeller $P_{\g(i_1,\ell_1;\cdots;i_n,\ell_n)}\subset \mW$ on
$\bRt$ consists of at most
$$\Delta(r(q^\g(n;i_1,\ell_1;\cdots;i_n,\ell_n)))+1$$ 
curves, one for each
intersection of the $\cW$-orbit of $q^\g(n;i_1,\ell_1;\cdots;i_n,\ell_n)$
with $\bRt\cap\{z\leq 0\}$. Call this number
$n(q^\g(n;i_1,\ell_1;\cdots;i_n,\ell_n))$.
Recall that the $\cW$-orbit of $p(1;i_1,\ell_1;\cdots;i_n,\ell_n)$
intersects $\bRt\cap\{z\leq 0\}$ in at most
$\Delta(r(p(1;i_1,\ell_1;\cdots;i_n,\ell_n)))+1$ points. Call this number
$n(p(1;i_1,\ell_1;\cdots;i_n,\ell_n))$. Hence the trace of
$P_{\g(i_1,\ell_1;\cdots;i_n,\ell_n)}$ on $\bRt$ are the curves
$\g_0(i_1,\ell_1;\cdots;i_n,\ell_n;\ell_{n+1})$ with $a\leq \ell_{n+1}\leq
n(q^\g(n;i_1,\ell_1;\cdots ;i_n,\ell_n))$.

If
$n(q^\g(n;i_1,\ell_1;\cdots;i_n,\ell_n))-n(p(1;i_1,\ell_1;\cdots;i_n,\ell_n))=0$,
all the curves in $P_{\g(i_1,\ell_1;\cdots;i_n,\ell_n)}\cap \bRt$ are arcs
as in Figure~\ref{fig:arcspropeller}, except that near the endpoints the
arcs coil as explained earlier in this section.

\begin{figure}[!htbp]
\centering
{\includegraphics[width=140mm]{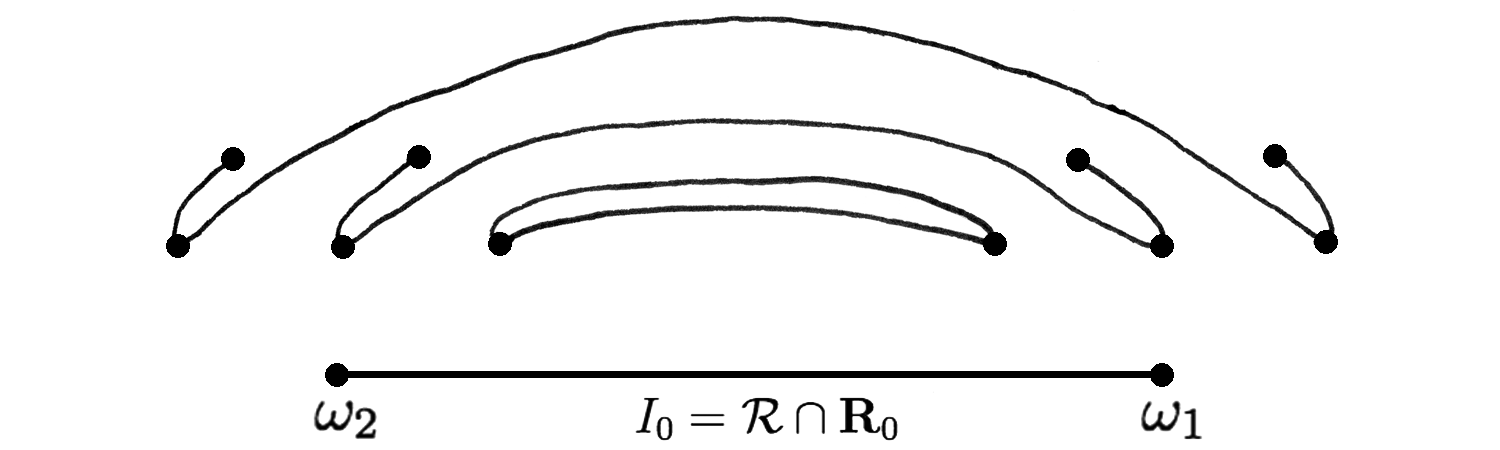}}
\caption[{Trace of the finite propeller  $P_{\g(i_1,\ell_1;\cdots;i_n,\ell_n)}$  in $\bRt$}]{\label{fig:arcsnonmonotone} Trace of the finite propeller  $P_{\g(i_1,\ell_1;\cdots;i_n,\ell_n)}$  in $\bRt$ with  $n(q^\g(n;i_1,\ell_1;\cdots;i_n,\ell_n))-n(p(1;i_1,\ell_1;\cdots;i_n,\ell_n))=1$}
\vspace{-6pt}
\end{figure}

If
$n(q^\g(n;i_1,\ell_1;\cdots;i_n,\ell_n))-n(p(1;i_1,\ell_1;\cdots;i_n,\ell_n))>0$,
then for 
$$n(p(1;i_1,\ell_1;\cdots;i_n,\ell_n))<\ell_{n+1}\leq
n(q^\g(n;i_1,\ell_1;\cdots;i_n,\ell_n))$$
 the curves
$\g_0(i_1,\ell_1;\cdots;i_n,\ell_n;\ell_{n+1})$ are closed, as in
Figure~\ref{fig:arcsnonmonotone}. The closed curves are contained in
the region 
$$\{r(q^\g(n;i_1,\ell_1;\cdots;i_n,\ell_n))\leq
r<r(p(1;i_1,\ell_1;\cdots;i_n,\ell_n))\}\cap \bRt.$$
In conclusion, the trace of $P_{\g(i_1,\ell_1;\cdots; i_n,\ell_n)}$ on
$\bRt$ can have a non-zero but finite number of closed curves.

\begin{figure}[!htbp]
\centering
{\includegraphics[width=120mm]{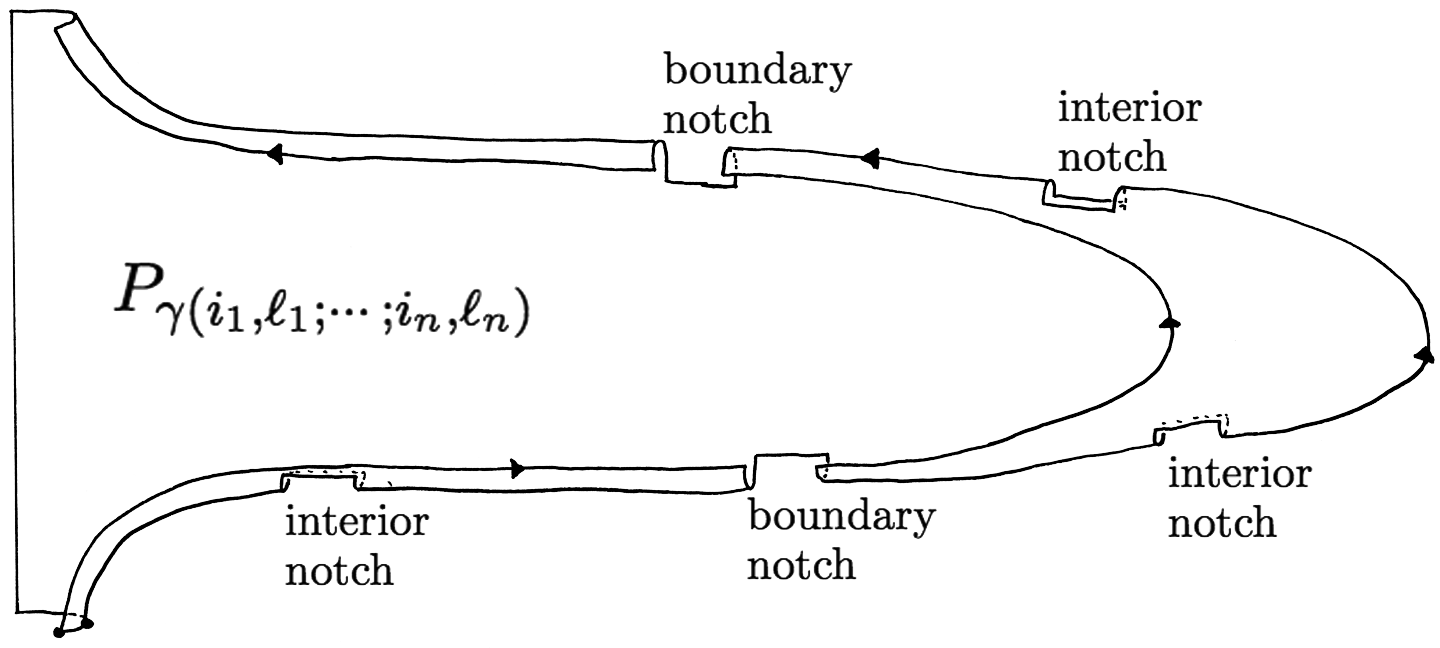}}
\caption[{Notches in the propeller $P_{\g(i_1,\ell_1;\cdots;i_n,\ell_n)}$}]{\label{fig:notches2} Notches  with  $n(q^\g(n;i_1,\ell_1;\cdots;i_n,\ell_n))-n(p(1;i_1,\ell_1;\cdots;i_n,\ell_n))=1$}
\vspace{-6pt}
\end{figure}

We now turn our attention once more to the intersection of
$P_{\g(i_1,\ell_1;\cdots; i_n,\ell_n)}$ with $\cL_i^-$. As explained in
Section~\ref{sec-proplevels}, the intersection consists of a 
number of internal notches, followed by boundary notches. After
the description in this section two new possibilities arise. First,
since the radius is non monotone along the curves in
$P_{\g(i_1,\ell_1;\cdots; i_n,\ell_n)}\cap \bRt$, some of these curves might create
internal notches when flowed to $\cL_i^-$. Second,  the
fact that the trace of  $P_{\g(i_1,\ell_1;\cdots; i_n,\ell_n)}$ on $\bRt$
can have closed curves adds the possibility of having internal notches
near the tip of the propeller, as illustrated in Figure~\ref{fig:notches2}.

The two level 1 propellers $\tau(P_\g')$ and $\tau(P_\lambda')$ in
$\fM_0$ are generated by curves with monotone radius. Thus they have a
finite number $|b|\geq 0$, maybe zero, of internal notches and an infinite number
of (boundary) notches, as described in
Section~\ref{sec-proplevels}. If any, the internal notches generate
bubbles as in Section~\ref{sec-bubbles}.

Consider the case of propellers of level at least 2. Let
$P_{\g(i_1,\ell_1;\cdots;i_n,\ell_n)}$ be a propeller of level
$n+1\geq 2$. Let $\g'(i_1,\ell_1;\cdots;i_n,\ell_n;1,\ell_{n+1}) \subset \cL_1^-$ be the curves
in the intersection of $P_{\g(i_1,\ell_1;\cdots;i_n,\ell_n)}$ with
$\cL_1^-$ and $\g_0(i_1,\ell_1;\cdots;i_n,\ell_n;\ell_{n+1})$ the corresponding
curve in $\bRt$. Internal notches arise when
$\g'(i_1,\ell_1;\cdots;i_n,\ell_n;1,\ell_{n+1})$ has its two
endpoints in $\partial \cL_1^-$.  

Assume that this is the case. If the
corresponding curve
$\g_0(i_1,\ell_1;\cdots;i_n,\ell_n;\ell_{n+1})\subset \bRt$ is not a closed curve,
then the forward orbit of its lower endpoint
$p_0(1;i_1,\ell_1;\cdots;i_n,\ell_n;1,\ell_{n+1})$ misses $\cL_1^-$. Thus the radius coordinate along the curve
$\g(i_1,\ell_1;\cdots;i_n,\ell_n;1,\ell_{n+1})\subset L_1^-$ is big,
meaning that it is uniformly bounded below. The
arguments in Section~\ref{sec-bubbles} apply to this case and thus the
surface generated by the $\cK$-orbits of the points in
$\tau(\g(i_1,\ell_1;\cdots;i_n,\ell_n;1,\ell_{n+1}))$ is compact and
admits a uniform upper bound to the difference of level between points
in it. We leave the details to the reader.

For the second case, assume that $\g'(i_1,\ell_1;\cdots;i_n,\ell_n;1,\ell_{n+1})$ has its two
endpoints in $\partial \cL_1^-$ and that the corresponding curve
$\g_0(i_1,\ell_1;\cdots;i_n,\ell_n;\ell_{n+1})\subset \bRt$ is a closed curve. This means that
$$n(q^\g(n;i_1,\ell_1;\cdots;i_n,\ell_n))-n(p(1;i_1,\ell_1;\cdots;i_n,\ell_n))>0$$
and $n(p(1;i_1,\ell_1;\cdots;i_n,\ell_n))<\ell_{n+1}\leq
n(q^\g(n;i_1,\ell_1;\cdots;i_n,\ell_n))$. 

\begin{lemma}\label{lem-unknown}
There exists $r_B>2$ such that if
$n(q^\g(n;i_1,\ell_1;\cdots;i_n,\ell_n))-n(p(1;i_1,\ell_1;\cdots;i_n,\ell_n))>0$
the radius coordinate on the curves $\g(i_1,\ell_1;\cdots;i_n,\ell_n;1,\ell_{n+1})\in L_1^-$
with
$$n(p(1;i_1,\ell_1;\cdots;i_n,\ell_n))<\ell_{n+1}\leq
n(q^\g(n;i_1,\ell_1;\cdots;i_n,\ell_n))$$
is bounded below by $r_B$.
\end{lemma}

\proof
Observe that the closed curve $\g_0(i_1,\ell_1;\cdots;i_n,\ell_n;\ell_{n+1})\in
\bRt$ is contained in the region bounded by $\G_0(i_2,\ell_2;\cdots;\ell_{n+1})$
if $i_1=1$ and in the region bounded by $\Lambda_0(i_2,\ell_2;\cdots;\ell_{n+1})$
if $i_1=2$. Also, it is far from
$\G_0(i_2,\ell_2;\cdots;\ell_{n+1})$ or $\Lambda_0(i_2,\ell_2;\cdots;\ell_{n+1})$ since the point $p_0(i_1;i_2,\ell_2;
\cdots;1,\ell_{n+1})$ is in the closure of the endpoints of the curves
$\g_0(i_1,m;\cdots;i_n,\ell_n;\ell_{n+1})$ that are arcs, for $m>\ell_1$
and unbounded. That is $\g_0(i_1,\ell_1;\cdots;\ell_{n+1})$ is separated from $\G_0(i_2,\ell_2;\cdots;\ell_{n+1})$ or $\Lambda_0(i_2,\ell_2;\cdots;\ell_{n+1})$ by the curves $\g_0(i_1,m;\cdots;\ell_{n+1})$ and $\lambda_0(i_1,m;\cdots;\ell_{n+1})$ for all $m>\ell_{n+1}$

Assume the conclusion in the lemma is false, then for every $\e>0$
there exists a curve \linebreak $\g(i_1,\ell_1;\cdots;i_n,\ell_n;1,\ell_{n+1})\in
L_1^-$ with
$$n(p(1;i_1,\ell_1;\cdots;i_n,\ell_n))<\ell_{n+1}\leq
n(q^\g(n;i_1,\ell_1;\cdots;i_n,\ell_n))$$
intersecting the ball centered at $p(1)$ of radius $\e$. 
Equivalently, for every $\e>0$ there exists a closed curve  $\g_0(i_1,\ell_1;\cdots;i_n,\ell_n;\ell_{n+1})\in
\bRt$ with
$$n(p(1;i_1,\ell_1;\cdots;i_n,\ell_n))<\ell_{n+1}\leq
n(q^\g(n;i_1,\ell_1;\cdots;i_n,\ell_n))$$
intersecting the ball centered at $\omega_1$ of radius $\e$. Thus
there is a sequence of closed curves accumulating on $\omega_1$ that
belong to the trace of the propellers. This is the contradiction that
we needed.
\endproof

\begin{remark}\label{rmk-secondbubbles}
An immediate consequence of   Lemma~\ref{lem-unknown} is that the 
surfaces given by the $\cK$-orbits of the points in
$\tau(\g(i_1,\ell_1;\cdots;i_n,\ell_n;1,\ell_{n+1}))$, with $\ell_{n+1}$ satisfying the hypotheses of the lemma,  are compact and the
difference of level between points in them is uniformly bounded as in
Proposition~\ref{prop-bubbles1}. We refer to these surfaces as bubbles
of the second type.
\end{remark}

  \bigskip

\section{Zippered laminations}\label{sec-zippered}

  The   notion of an $n$-dimensional  lamination is well-known, and
  can be summarized as an $n$-dimensional   {foliated space} which
  is transversally modeled on a totally disconnected space. 
In this section,  we  introduce the notion of   \emph{zippered laminations}, which are a type of ``stratified'' laminations with possibly pathological behavior for the strata.  After introducing this and some other   preliminary notions, we prove the main result of this section:
  
   \begin{thm} \label{thm-zippered} 
If $\Phi_t$ is  a generic Kuperberg flow on $\mK$, then $\fM$ is a \emph{zippered lamination}. 
 \end{thm}

We begin by recalling the notion of a foliated space \cite{CandelConlon2000,MS2006}.
A \emph{continuum}   is   a non-empty, compact, connected metric space.  
\begin{defn} \label{def-fs}
A  \emph{foliated space of dimension $n$} is a   continuum $\fZ$ with a foliated structure. That is,   there exists: 
\begin{enumerate}
\item a separable metric space $\fX$, and a collection of compact subsets $\fT_i \subset \fX$  for $1 \leq i \leq k$;
\item a finite collection of homeomorphisms  (the charts)  $\{ \vp_i \colon \fT_i \times [-1,1]^n   \to \oU_i  \subset \fZ \mid 1 \leq i \leq k \}$;
\item the ``interiors'' $\{ U_i = \vp_i(\fT_i \times (-1,1)^n) \mid 1 \leq i \leq k \}$ form an   covering for $\fZ$;
\item the charts $\{\vp_i \mid 1 \leq  i \leq k\}$   satisfy the compatibility condition \eqref{eq-overlaps}.
\end{enumerate}
  \end{defn} 
 
Let $\pi_i \colon \oU_i \to \fT_i$ denote the composition of $\vp_i^{-1}$ with projection onto the first factor.  

For $\xi \in \fT_i$ the set $\cP_i(\xi) = \vp_i(\{\xi\} \times [-1,1]^n) \subset \oU_x$ is called a \emph{plaque} for the coordinate chart $\vp_i$. 
For $z \in \oU_i$, we adopt the notation $\cP_i(z) = \cP_i(\pi_i(z))$, so that $z \in \cP_i(z)$. 

For $\xi \in \fT_i$ the plaque $\cP_i(\xi)$ is given the topology so that the restriction $\vp_i \colon   \{\xi\} \times  [-1,1]^n \to \cP_i(\xi)$ is a homeomorphism, hence $int (\cP_i(\xi)) = \vp_i(\{\xi\} \times (-1,1)^n)$.
Then we require, in addition, 
\begin{equation}\label{eq-overlaps}
{\rm for ~ all} ~ z \in U_i \cap U_j ~, ~  \interior(\cP_i(z)) \cap \interior( \cP_j(z)) ~ \mbox{is~ an ~ open ~ subset ~ of both} ~
 \cP_i(z) ~ {\rm and}~ \cP_j(z) ~ .
 \end{equation} 
 
 The collection of sets
$\ds \cV = \{ \vp_i(\{\xi\} \times V) \mid 1 \leq i \leq k ~, ~ \xi \in \fT_i ~, ~ V \subset (-1,1)^n ~ {\rm open}\}$ 
forms the basis for the \emph{fine topology} of $\fZ$. The connected components of the fine topology are called leaves, and define the foliation $\F$ of $\fZ$.
For $x \in \fZ$, let $L_x \subset \fM$ denote the leaf of $\F$ containing $x$. The subsets  $\fT_i$ of $\fX$ are called the  \emph{local transverse models}.

 Many of the techniques of foliation theory extend to foliated spaces, including the existence of holonomy transformations defined by ``parallel transport'' along paths in the leaves,  which generate a pseudogroup model for the dynamics of the foliated space. 
Foliated spaces are considered in greater detail in \cite[Chapter 2]{MS2006} and  \cite[Chapter 11]{CandelConlon2000}.

An \emph{$n$-dimensional lamination} is a    foliated space $\fZ$   of dimension $n$, such that  the transverse model space $\fX$ is totally disconnected.   
The   leaves of $\F$ are then the  path components of $\fZ$. 

The definition of a \emph{lamination with boundary} is obtained by modifying the above definition, to   allow  foliation charts of the form $\vp_i \colon \fT_i \times  ([0,1) \times (-1,1)^{n-1}) \to U_i$. The boundary $\partial \fZ$ then consists of the points which are  images under a chart   of the chart boundary $\fT_i \times (\{0\} \times (-1,1)^{n-1})$.

The notion of a ``zippered lamination'' is motivated by our study of the space $\fM$, which is  the closure of the  non-compact manifold with boundary $\fM_0$. 
The unusual property of $\fM$ is that   the boundary curves of $\fM_0$ are also dense in $\fM$, which is impossible for a usual lamination with boundary. 
The definition of a  zippered lamination below is similar to that of a lamination with boundary, in that each of its leaves is a manifold with boundary,  but   the boundaries of the leaves do not have to align themselves in a way that they are regularly covered by foliation charts. In place of the covering of $\fM$ by standard foliation coordinate charts  as in Definition~\ref{def-fs}, the definition below covers $\fM$ by foliation coordinate charts whose plaques have varying sizes, and cannot be ``standardized'' in a continuous manner, though their overlaps do define continuous holonomy transformations.

 \begin{defn} \label{def-zl}
 An \emph{$n$-dimensional zippered lamination} $\fZ$ is a   continuum whose path components  are $n$-dimensional manifolds with boundary.
Denote the union of the boundaries of the leaves   by $\partial_{\F} \fZ$,  and    the complement by  $\interior(\fZ) =   \fZ  -    \partial_{\F} \fZ$. Then we require that $\interior(\fZ)$   is dense in $\fZ$, and the following data are given, which satisfies the accompanying conditions: there is given a   compact  metric space $\fX$, and for $1 \leq i \leq k$, 
\begin{enumerate}
\item a subset $\fT_i \subset \fX$  whose closure in $\fX$ is totally disconnected;
\item a Borel subset $\cB_i \subset \fT_i \times [-1, 1]^n$, where for each $\xi \in \fT_i$ the slice
$\ds  \cB_i \cap \{\xi\} \times [-1,1]^n$ is a compact subset whose interior contains $\{\xi\} \times \{0\}$, and is homeomorphic to $[-1,1]^n$; 
\item a    homeomorphism   $\ds \vp_i \colon  \cB_i \to \oU_i   \subset \fZ$ onto its image $\oU_i$ with the induced topology from $\fZ$;
\item the interiors $\{U_i = \interior(\oU_i) \mid 1 \leq i \leq k\}$ form a   covering for $\interior(\fZ)$;
\item the collection of charts $\{\vp_i \mid 1 \leq  i \leq k\}$   satisfy the compatibility condition \eqref{eq-overlaps}.
\end{enumerate}
  \end{defn} 
  The interior $\interior(\oU_i) $ for the Borel set $\oU_i$ is defined as the union over the points $\xi \in \fT_i$ of the interiors of the plaques $\cP_i(\xi)=\vp_i(\cB_i \cap \{\xi\} \times [-1,1]^n)$ in $\oU_i$.

Note that the Borel  property in \ref{def-zl}.2 allows the plaques $\cP_i(\xi_{\ell}) \subset \fZ$ to degenerate in size  for a sequence $\xi_{\ell} \in \fT_i$ which converge to a point $\xi_* \in \fX$ with $\xi_* \not\in \fT_i$.  In fact, the limit of such a sequence  must then be   contained in  $\partial_{\F} \fZ$. 
 
The first step in showing that $\fM$ satisfies the conditions of Definition~\ref{def-zl} is to construct the model spaces $\fT_i$.  
Recall that   $\cT=\{z=0\} \cap  \bR_0$, and set:
\begin{equation}\label{eq-cantorset1}
\fC_0' = \cT \cap \fM_0 \quad , \quad \fC' = \cT \cap \fM ~. 
\end{equation}
The    points in $\fC_0'$ are classified according to the properties of the curves in $\fM_0 \cap \bRt$ which define them. 
The intersection of $\cT$ with the Reeb cylinder $\tau(\cR')$ defines the point   $\omega_0 \in \fM_0 \cap \bRt$ with $r(\omega_0) =2$. 
The remaining points in the intersection correspond to the intersections with $\g_0$ or $\lambda_0$ curves, and these can be of  two types:
 
\begin{enumerate}
\item The $\g_0$ or $\lambda_0$ curve is an arc. In the case where the
  forward orbit of 
  points in such a curve hit the entry regions, the curve corresponds
  to boundary notches in the propeller, one for each entry region. 
\item The $\g_0$ or $\lambda_0$ curve is  a closed curve, hence there are   two points in $\fC_0$ corresponding to the same curve.   In the case where the
  forward orbit of 
  points in such a curve hit the entry regions, the curve corresponds
  to interior notches in the propeller, either one or two for each
  entry region. These closed curves correspond either to notches
  generating ``bubbles'' or to intersections of ``bubbles'' with
$\bRt$ as described in Sections~\ref{sec-bubbles} and \ref{sec-geometry}.
\end{enumerate}

Let $\fC_0^1 \subset \fC_0'$ be the intersection points corresponding to closed curves in  $\fM_0 \cap \bRt$.  
Define: 
\begin{equation}\label{eq-cantorset}
\fC_0 = \fC_0' - \fC_0^1 \qquad , \qquad \fC=\overline{\fC_0}.
\end{equation}

\begin{prop}\label{prop-cantor}
The set  $\fC$ is a perfect and totally disconnected subset of $\fC'$.
\end{prop}
\proof
Each point in $\fC_0$ is the intersection of $\cT$ with   a certain
$\gamma_0$ or $\lambda_0$ curve, and thus is defined by the labeling of the endpoints of this curve as in Section~\ref{sec-proplevels}. 
Moreover,  the symmetry of the Wilson flow implies that each such $\gamma_0$ and  
$\lambda_0$ curve in $\fM_0\cap \bR_0$ intersects $\cT$. Thus, 
a sequence of curves whose intersections with $\cT$ defines points in $\fC_0$ which converge, also   defines a sequence of converging endpoints, whose corresponding limit curve defines a point in $\fM \cap \cT$. Thus, $\fC = \overline{\fC_0} \subset \fM \cap \cT = \fC'$. 
Lemma~\ref{lem-paccumulationpoints} implies that each point of $\fC_0$ is a limit point of the set, hence the closure $\fC$ is a perfect set. It is totally disconnected  by the results of Section~\ref{sec-wandering}. 
\endproof

\begin{remark}
It  can also be proved that $\fC'$ is a Cantor set, since any
$q^\g$ or $q^\lambda$ point is the accumulation point of a  
sequence of $q^\g$ and $q^\lambda$ points of lower level.
\end{remark}

The strategy for constructing foliation charts for $\fM$ using the Cantor set $\fC$ is simple in outline, and we next develop. But there is a nuance that arises near points of $\fC$ corresponding to tips of propellers, that requires covering these ``singular points'' by a separate construction.

Consider $\xi \in \fC - \fC_0$, then there exists $\{ x_{k} \mid k \geq 1\}
\subset \fC_0$ with $x_{k} \to \xi$. 
By the results of
Section~\ref{sec-proplevels}, each $x_{k}$ lies on a propeller curve
$\g_0(i_1 , \ell_1 ; i_2, \ell_2; \cdots ; \ell_n)$ or $\lambda_0(i_1
, \ell_1 ; i_2, \ell_2; \cdots ; \ell_n)$ where
the multi-index depends on $k$. 
By passing to a subsequence, we can
assume that they are all of the same type, either all $\gamma_0$ curves, or all $\lambda_0$ curves.
Each such curve is  then
defined by the labeling of its endpoints, denoted simply by
$p_0^{i_n}(x_k)$  to shorten the notation introduced in
Section~\ref{sec-doublepropellers}. 
That is,
$p_0^{i_n}(x_k)=p_0(i_0;i_1,\ell_1; i_2, \ell_2; \cdots ;
i_{n-1},\ell_{n-1}; i_n, \ell_n)$ for $i_n=1,2$, where $i_0 = 1$ for a $\gamma_0$ curve, and $i_0 =2$ for a $\lambda_0$ curve. 

We consider the case where  we have a sequence of   $\g_0$ curves. The case for a sequence of $\lambda_0$ curves is analogous. 
 As $\bRt$ is compact, we can pass to a subsequence, and assume that
 the lower endpoints $p_0^1(x_{k})$ converge to a point $p_0^1(\xi)$ where $z(p_0^1(\xi)) \leq 0$,  and similarly  $p_0^2(x_{k}) \to p_0^2(\xi)$ where $z(p_0^2(\xi)) \geq 0$.

The sequence of endpoints $p_0^1(x_{k})$ all lie on the $\cK$-orbit
of the special point $\omega_1$, so each has a well-defined level by
Proposition~\ref{prop-levelsM}. If there is a uniform bound on this
level, then the limit must be contained in $\fM_0$ contrary to the
choice of $\xi$. Thus, passing to a subsequence yet again, we can
assume that the level $n_0(x_{k})$ is monotone increasing.

It then follows that the number $n$ of indices of the sequence $\g_0(i_1 ,
\ell_1 ; i_2, \ell_2; \cdots ; i_n, \ell_n)$ must tend to infinity,
and by the results of Section~\ref{sec-doublepropellers}, the sequence
$\{x_k\}$ defines a nested sequence of ellipses in $\bRt$, each ellipse  bounded by
the union of the associated $\g_0$ and $\kappa_0$ curves. 
Moreover, the width of these ellipses must tend to zero. It follows that the intersections of the  interiors of these nested ellipses define an arc in $\bRt$,  denoted by 
$[p_0^1(\xi), p_0^2(\xi)]$ where the point $\xi = [p_0^1(\xi), p_0^2(\xi)]  \cap \{z=0\}$.

The arc $[p_0^1(\xi), p_0^2(\xi)]$ is the limit of the boundary
$\g_0$ curves of the nested ellipses, as a point set. Each of these
boundary $\g_0$ curves  is an embedded curve with uniformly bounded length by
Theorem~\ref{thm-boundedlength}. Thus, 
these curves converge  to an immersed curve of finite length. That is,  
in a generic Kuperberg plug, for $\xi \in \fC$  the curve $[p_0^1(\xi), p_0^2(\xi)]$ has uniformly bounded    arclength.

We next repeat these arguments for each  of four rectangular sections to $\fM$, defined as follows, and the associated Cantor transversals.  Introduce the rectangles
\begin{equation}
\bT_i = \{\xi = (r, i \cdot\pi/2, z) \mid ~ 1 \leq r \leq 3 ~,   ~ -2 \leq z \leq 2\} ~, ~ 1 \leq i \leq 4
\end{equation}
so that $\bRt = \bT_2$. We add a fifth rectangle 
$$
\bT_5 = \{\xi = (r,3\pi/8, z) \mid ~ 1 \leq r \leq 3 ~,   ~ -2 \leq z
\leq 2\}.
$$
Each of these rectangles  is disjoint from both the regions $D_i$ and their insertions $\cD_i$ for $i =1,2$ by the choices made in 
Section~\ref{sec-kuperberg}. It is helpful to review the illustrations in 
 Figures~\ref{fig:insertiondisks} and  Figure~\ref{fig:K}. The key property of these sections that we require, is that for any $\cK$ orbit segment $[x,y]_{\cK}$ with $x, y \in \bT_i$, then either the interior segment $(x,y)_{\cK}$ intersects the annulus $\cA$, or it must intersect one of the other sections $\bT_j$ for $j \ne i$.
 
For each of $1 \leq i \leq 5$, introduce the set $\fC_{0,i} \subset
\bT_i \cap \cA \cap \fM_0$   consisting of  the points   that do not belong to a
bubble of any type, and correspond to arcs in $\bT_i \cap \fM_0$. 
Note that      $\fC_{0,2}=\fC_0$. 
Define $\fC_i$ to be the closure of  $\fC_{0,i}$.
We again conclude that each $\xi \in \fC_i$ defines an
immersed arc $[p_i^{1}(\xi), p_i^{2}(\xi)]$, each of whose length admits a uniform upper bound.

For each $1 \leq i \leq 5$, define   functions $S_i^{\pm} \colon \fC_i
\to [0,\infty)$ such that for $\xi \in \fC_i$
\begin{itemize}
\item $S_i^+(\xi)$ is the length of the arc   
 $[\xi, p_i^{2}(\xi)] \subset  [p_i^{1}(\xi), p_i^{2}(\xi)]$ containing $\xi$ as the
 lower endpoint, where $[\xi, p_i^{2}(\xi)] \subset \{z \geq 0\}$; 
\item   $S_i^-(\xi)$ is the length of the arc 
 $[p_i^{1}(\xi), \xi] \subset [p_i^{1}(\xi), p_i^{2}(\xi)]$ containing $\xi$ as the
 upper endpoint, where $[p_i^{1}(\xi), \xi] \subset \{z \leq 0\}$. 
\end{itemize}

 Since $\tau^{-1}(p_i^2(\xi))$ and $\tau^{-1}(p_i^1(\xi))$ belong to the
 same $\cW$-orbit, the symmetry of the $\Psi_t$-flow implies that
 $S_i^+(\xi)=S_i^-(\xi)$.  Let $S_i(\xi) = S_i^+(\xi)$. Then by the above, for each $i$,  the function $\xi \mapsto S_i(\xi)$ admits a uniform
 upper bound for all $\xi \in \fC_i$.  Let $L_i = \sup \{S_i(\xi) \mid \xi \in \fC_i\}$ denote the upper bound of these lengths for $\xi \in \fC_i$.

 The function   $S_i(\xi)$   is not continuous, but we note that the proof of  \eqref{eq-approxlength}
  shows that for $\e > 0$, there exists $N_{\e}$ which is independent of $\xi$, such that there exists $x \in \fM_0 \cap \fC_i$ for which the curve $\sigma_x$ containing $x$ has length $|S_i(x) - S_i(\xi) | < \e$ and moreover, the endpoint $p_i^1(x)$ has level  $n \leq N_\e$.  It follows that the arclength function   $S_i \colon \fC_i \to [0,L_i]$ is the uniform limit of   step-functions on $\fC_i$ so is a Borel function of $\xi$.

Define $\fC_i^+ = \{ \xi \in \fC_i \mid S_i^+(\xi) > 0\}$.   For   $\xi \in \fC_i^+$, let $\gamma_{\xi}(s)$ 
denote the parametrization of $[p_i^{1}(\xi), p_i^{2}(\xi)]$ by arclength,   with $\gamma_{\xi}(0) = \xi$.

 Fix $1 \leq i \leq 5$ and consider a point   $\xi \in \fC_i$. The proof of Proposition~\ref{prop-sigmadyn} carries over to each of the surfaces $\bT_i$, so there exists a constant $\nu_*$ such that 
the set  $\ds \cS_{i,\xi} = \{ s \mid \Phi_s(\xi) \in \fM \cap \bT_i \}$ is syndetic for the constant $\nu_*$. Moreover, the constant $\nu_*$  can be chosen   independent of $i$ and $\xi$.    There is a nuance in the estimation of the return time, $\nu_*$ which bounds  the return time of a point in  $\bT_i$ to the same transversal.   The return point is contained in $\fC_i' = \fM \cap \bT_i$ but not necessarily in the subset $\fC_i$. 
However, Lemma~\ref{lemma-bubblelevel}, Proposition~\ref{prop-bubbles1} and the
discussion in Section~\ref{sec-geometry} imply
that there exists a uniform constant which is an upper bound on the
length of time any orbit spends in a bubble. That is, the time
required to flow across each of the internal notches as illustrated in
Figures~\ref{fig:intnotches} and \ref{fig:notches2} has a uniform bound. Thus, replacing $\nu_*$ with a fixed multiple of itself, we may assume that $\nu_*$ is a syndetic constant for $\ds \cS_{i,\xi} = \{ s \mid \Phi_s(\xi) \in \fC_i \}$ for each $1 \leq i \leq 5$.

 Consider the annulus $\cA\subset \mW$ and its image $\tau(\cA)$
  in $\mK$. For each
$\theta\in \mS^1$ consider the sets
\begin{eqnarray*}
{\bf R}_\theta & = & \{\tau(r,\theta,z)\,|\, (r,\theta, z)\in \mW', \,
2\leq r\leq 3,\, -2\leq z\leq 2\}\\
\cT_\theta & = & \{\tau(r,\theta,0)\,|\,2\leq r\leq 3\}
\end{eqnarray*}
with $\cT_\theta\subset {\bf R}_\theta$. Observe that for those values
of $\theta$ for which all the points $(r,\theta,z)$ are not in $\cD_i$
for $i=1,2$, ${\bf R}_\theta$ is a rectangle in $\mK$. On $\cT_\theta$
consider the set of points $\fC_\theta '=\cT_\theta\cap \fM_0$. This
set contains two types of points: those that belong to a closed curve in
${\bf R}_\theta\cap \fM_0$ that we denote by $\fC_\theta^1$ and those
that belong to an arc in ${\bf R}_\theta\cap \fM_0$ . Consider the set
$\fC_{0,\theta}=\fC_\theta '-\fC_\theta ^1$ and its closure $\fC_\theta$. Observe that $\fC_{0,\pi/2}=\fC_0$. Let $\cA_\fC$ be the union over $\theta$ of
these Cantor sets $\fC_\theta$. Observe that $\cA_\fC$ contains the
intersection of $\tau(\cA)$ with  all the simple propellers in $\fM_0$.

Let $\fM_{\fC_i} \subset \fM\cap\bT_i$ be the path components that contain a point in
$\fC_i$ for $1\leq i\leq 5$, and define $W = \fM_{\fC_1} \cup \cdots \cup \fM_{\fC_5} \cup \cA_\fC$ be the  union of these closed sets in $\mK$. 
For each $i$ and for $x \in \fM_{\fC_i}$,  define      $T_{i}^\pm(x)$,   which  is the  forward and backwards ``return time'' to $W$   for the flow $\Phi_t$:
\begin{equation}\label{eq-parameterlimits}
T_{i}^+(x)   =   \inf \, \{s > 0\mid \Phi_s(x) \in W \} \quad , \quad 
T_{i}^-(x)   =    \sup \, \{s < 0\mid \Phi_s(x) \in W \} ~ .
\end{equation}
For $\xi\in\fC_i$ we have that $ | T_{i}^\pm (\g_\xi(s)) | < \nu_*$ uniformly. The values of these functions are not necessarily continuous in the variables $\xi$ and  $s$. The   discontinuities may arise,  for example, when the endpoint of the flow segment starting at $\g_\xi(s)$  jumps from a section $\bT_j$ to the annulus $\cA$.

For each $1 \leq i \leq 5$,     define the subsets of $\fC_i^+ \times [-L_i, L_i] \times [-\nu_*, \nu_*]$, 
\begin{eqnarray*} 
D_{i,1} ~ & = & ~  \{  (\xi, s,t) \mid \xi \in \fC_i^+ ~ ,  ~ - S_{i}(\xi) \leq s \leq 0~, ~ T_{i}^-(\gamma_{\xi}(s))  \leq t \leq  T_{i}^+(\gamma_{\xi}(s))  \}   \label{eq-domainxi1}\\
D_{i,2} ~ & = & ~  \{  (\xi, s,t) \mid \xi \in \fC_i^+ ~ ,  ~ 0   \leq s \leq  S_{i}(\xi)  ~, ~ T_{i}^-(\gamma_{\xi}(s)) \leq t \leq T_{i}^+(\gamma_{\xi}(s))\}  ~ . \label{eq-domainxi2}
\end{eqnarray*}
Let $\e_* > 0$ be sufficiently small so that the $\e_*$ neighborhood of $W$ in $\mK$ is disjoint from the surfaces $E_i$ and $S_j$ for $1 \leq i,j \leq 2$ as defined in \eqref{eq-sections}. 
Then for each $1 \leq i \leq 5$,  also   define 
 \begin{eqnarray*} 
D_{i,1}^* ~ & = & ~  \{  (\xi, s,t) \mid \xi \in \fC_i^+ ~ ,  ~ -  S_{i}(\xi) \leq s \leq 0~, ~ T_{i}^-(\gamma_{\xi}(s)) -\e_* \leq t \leq  T_{i}^+(\gamma_{\xi}(s)) +\e_* \}   \label{eq-domainxi1e}\\
D_{i,2}^* ~ & = & ~  \{  (\xi, s,t) \mid \xi \in \fC_i^+ ~ ,  ~ 0   \leq s \leq  S_{i}(\xi)  ~, ~ T_{i}^-(\gamma_{\xi}(s)) -\e_* \leq t \leq T_{i}^+(\gamma_{\xi}(s)) +\e_* \} ~ . \label{eq-domainxi2e}
\end{eqnarray*}
Introduce the continuous maps:
\begin{eqnarray*} 
\vp_{i,1} \colon D^*_{i,1} \to \fM ~ & , & ~ \vp_{i,1}(\xi, s,t) = \Phi_t(\gamma_{\xi}(s))   \label{eq-vpdomainxi1}\\
\vp_{i,2} \colon D^*_{i,2} \to \fM ~ & , & ~  \vp_{i,2}(\xi, s,t) = \Phi_t(\gamma_{\xi}(s))    ~ . \label{eq-vpdomainxi2}
\end{eqnarray*}

\begin{prop} \label{prop-injective}
For each $1 \leq i \leq 5$ and $j=1,2$, the map $\vp_{i,j} \colon D^*_{i,j} \to \fM$ is one-to-one. Moreover, we have
$$\fM = \bigcup_{1 \leq i \leq 4} ~ \left\{ \vp_{i,1}(D_{i,1}^*) \cup \vp_{i,2}(D_{i,2}^*)\right\} .$$
\end{prop}
\proof
Suppose that $(\xi, s,t) \ne (\xi', s',t') \in D^*_{i,j}$ satisfy
$\vp_{i,j}(\xi, s,t) = \vp_{i,j}(\xi', s',t')$ for some $i,j$. 
Let $x = \gamma_{\xi}(s)$  and $x' = \gamma_{\xi'}(s')$ then $x' =
\Phi_{-t'}(\vp_{i,j}(\xi,
s,t))=\Phi_{-t'}(\Phi_t(x))=\Phi_{t-t'}(x)$. Assume without loss of
generality that $t-t'>0$. Then there is a
$\cK$-orbit segment $\sigma$ from $x\in \bT_i$ to $x'\in \bT_i$ whose
interior intersects at most once
$W$. Let $x_1, x_2,\ldots$ be the transition
points in the forward $\cK$-orbit of $x$. Observe that $\sigma$ must
contain at least one transition point.

Assume first that $\sigma\cap \cA_\fC$ is empty. By the choices of $\bT_i$, we have that if:
\begin{enumerate}
\item $x\in \bT_1\cap\{z<0\}$ then $[x,x_1]_\cK$ intersects
  $\bT_2\cap\{z<0\}$.
\item $x\in \bT_2\cap\{z<0\}$ then either $[x,x_1]_\cK$ intersects
  $\bT_3\cap\{z<0\}$ or $x_1$ is in $E_1$.
\item $x\in \bT_3\cap\{z<0\}$ then $[x,x_1]_\cK$ intersects $\bT_4\cap\{z<0\}$.
\item $x\in \bT_4\cap\{z<0\}$ then $[x,x_1]_\cK$ intersects
  $\bT_5\cap\{z<0\}$.
\item $x\in \bT_5\cap\{z<0\}$ then $[x,x_1]_\cK$ intersects $\bT_1\cap\{z<0\}$.
\item $x\in \bT_1\cap\{z>0\}$ then $[x,x_1]_\cK$ intersects
  $\bT_5\cap\{z>0\}$.
\item $x\in \bT_2\cap\{z>0\}$ then either $[x,x_1]_\cK$ intersects
  $\bT_1\cap\{z>0\}$ or $x_1$ is in $E_2$.
\item $x\in \bT_3\cap\{z>0\}$ then $[x,x_1]_\cK$ intersects $\bT_2\cap\{z>0\}$.
\item $x\in \bT_4\cap\{z>0\}$ then either $[x,x_1]_\cK$ intersects
  $\bT_3\cap\{z>0\}$ or $x_1$ is in $S_2$.
\item $x\in \bT_5\cap \{z>0\}$ then either $[x,x_1]_\cK$ intersects
  $\bT_4\cap\{z>0\}$ or $x_1$ is in $S_1$. 
\end{enumerate}

If $\sigma$ intersects $W$ only at its endpoints, in cases (1), (3),
(4), (5), (6) and (8) we obtain a contradiction. We next
analyze the other four possible cases.
\begin{itemize}
\item[(2)]
Assume that $x\in \bT_2\cap\{z<0\}$, then $x_1\in E_1$. 
From $x_1$, the orbit segment $\sigma$ either continues to a point
$y_1\in \bT_5$ or it flows to $x_2\in S_1$ before intersecting $W$,
with $x_1\equiv x_2$. In
the first case, if 
  $y_1\in \fM_{\fC_5}$ we obtain a contradiction. 
  If not, $r(x_1)>2$ and $\sigma$ contains the
  point $\overline{x}_1\in S_1$ such that $x_1\equiv \overline{x}_1$
  (this is the case when $\sigma$ contains a segment of $\cK$-orbit
  in a bubble). Thus there exists $k\geq 2$ such that $x_k\in S_1$
  with $x_1\equiv x_k$. From $x_k$ the segment $\sigma$ continues
  and hits $\fM_{\fC_3}\subset \bT_3$, giving us a
  contradiction. Observe that in the second case we obtain the same
  contradiction with $k=2$.

\item[(7)] If $x\in \bT_2\cap\{z>0\}$  then $x_1\in E_2$. From $x_1$, the orbit segment $\sigma$ either continues to a point
$y_1\in \bT_4$ or it flows to $x_2\in S_2$ before intersecting $W$,
with $x_1\equiv x_2$. In
the first case, either $y_1\in \fM_{\fC_4}$ and we obtain a
contradiction, or there exists $k\geq 2$ such that $x_k\in S_2$
  with $x_1\equiv x_k$ and $[x_1,x_k]_{\cK}\cap W=\emptyset$. From $x_k$ the segment $\sigma$ continues
  and hits $\fM_{\fC_1}\subset \bT_1$, giving us a
  contradiction. Observe that in the second case we obtain the same
  contradiction.

\item[(9)] If $x\in \bT_4\cap\{z>0\}$ then $x_1\in S_2$. From $x_1$,
  the orbit segment $\sigma$ continues to a point $y_1\in
  \bT_1$. Since $x_1$ is a secondary exit point, 
  $y_1\in \fM_{\fC_1}$, a contradiction.

\item[(10)] If $x\in \bT_5\cap\{z>0\}$ then $x_1\in S_1$. From $x_1$,
  the orbit segment $\sigma$ continues to a point $y_1\in
  \bT_3$. Since $x_1$ is a secondary exit point, 
  $y_1\in \fM_{\fC_3}$, a contradiction.
\end{itemize}

We are left with the case where the interior of $\sigma$ intersects
$W$, that we separate in two situations, first when the intersection
point is in $W-\cA_\fC$ and second when it is in $\cA_\fC$. In the
first situation, concatenating the cases above, we conclude that there
is no $\cK$-orbit segment with both endpoints in one of  transversals and
intersecting $W-\cA_\fC$.

In the case where $\sigma\cap \cA_\fC\ne \emptyset$, the only
possibility for such a segment is to flow from $x$ to a transition
point $x_1$. If $[x_1,x_2]_\cK$ does not intersects $W$, then
$x_1\equiv x_2$ and we can replace it with a short cut. Then, $x\in
\bT_2$ and, by (2) and (7) above, the arc $[x_2,x_3]_\cK$ intersects
$W-\cA_\fC$. Thus the first intersection point of $\sigma$ with $W$ 
is either in $\bT_1$ or $\bT_3$. A contradiction.

Hence, $[x_1,x_2]_\cK$ must intersect
$\cA_\fC$. Let $y$ be the intersection point. 
If $x_1$ is a secondary
entry point, then $x_1\equiv x_2$, $x\in \bT_2$ and $x'$
is either in $\bT_1$ or $\bT_3$. If $x_1$ is a secondary exit point,
then $x\in \bT_4\cap \{z>0\}$ or $x\in \bT_5\cap\{z>0\}$. In the first
case, $x_1\in S_2$ and thus has $z$-coordinate positive. Case (9)
above implies that its orbit intersects $W-\cA_\fC$ before
intersecting $\cA_\fC$. In the second case, $x_1\in S_1$ and by
assumption its
forward orbit intersects $\cA_\fC$ before intersecting the transversal
$\bT_3$. Then it flows in the negative $\theta$-direction to the
transversal $\bT_2$, and since $x_1$ is a secondary exit point, the
last intersection point belongs to $\fM_{\fC_2}$. The last
contradiction we needed.

Thus, the maps $\vp_{i,j}$ are homeomorphisms. By definition they cover $\fM$.
\endproof

  \bigskip

   \section{Entropy of the flow}\label{sec-entropyflow}

  A celebrated theorem  of Katok    \cite{Katok1980} on the
 entropy for $C^2$-flows      on compact $3$-manifolds implies that  the   topological entropy of the Kuperberg flow $\Phi_t$  is   zero.
  In this section, we give a  geometric proof of this conclusion, based on an
 analysis of the dynamics of the flow $\Phi_t$ restricted to the closed invariant set\ $\fM$.  
The key idea is to relate the flow entropy to   another type of entropy invariant,  which is derived from the pseudogroup dynamics for $\cGK$ acting  on the rectangle $\bRt$.   
 
  We first prove, in Proposition ~\ref{prop-entropywhPhi}, that the entropy of the flow vanishes if the entropy of the return map $\whPhi$ vanishes. Then in Proposition~\ref{prop-entropysection},  we relate the entropy of $\whPhi$   to the entropy of the $\psg$ $\cGK^*$, proving that they are proportional. Finally, in Theorem~\ref{thm-entropyvanish}, we show that the entropy  $h_{GLW}(\cGK^* | \fMR)$ of $\cGK^*$ vanishes.
 The study of the various $\psg$ entropies associated to the flow
 $\Phi_t$   reveals the geometrical principles behind the vanishing of its entropy, and provides further insights into the dynamics of this flow. 
 
We define the entropy of the flow $\Phi_t$  using a variation of the Bowen formulation of   topological entropy  \cite{Bowen1971,Walters1982} for a flow   on a compact metric space $(X, d_X)$. The definition we adopt    is symmetric in the role of the time variable $t$.    Two points $p,q\in X$ are said to be 
\emph{$(\vp_t , T, \e)$-separated} if 
$$ d_X(\vp_t(p),\vp_t(q))>\e \quad \mbox{for some} \quad -T\leq t\leq T ~ .$$
A set $E \subset X$ is \emph{$(\vp_t , T, \e)$-separated} if  all pairs of distinct 
points in $E$ are $(\vp_t , T, \e)$-separated. Let $s(\vp_t , T, \e)$ be the maximal
cardinality of a $(\vp_t , T, \e)$-separated set in $X$. Then the topological entropy
is defined by
\begin{equation}\label{eq-separated}
h_{top}(\vp_t)= \frac{1}{2} \cdot \lim_{\e\to  0} \left\{ \limsup_{T\to\infty}\frac{1}{T}\log(s(\vp_t , T, \e)) \right\} .
\end{equation}
Moreover, for a compact space $X$, the entropy $h_{top}(\vp_t)$ is independent of the choice of metric $d_X$.
 
A relative form of the topological entropy for a flow $\vp_t$ can be defined for any   subset $Y \subset X$, by requiring that the collection of distinct $(\vp_t , T, \e)$-separated  points used in the definition \eqref{eq-separated}  be contained  in $Y$. The restricted topological entropy $h_{top}(\vp_t | Y)$ is bounded above by $h_{top}(\vp_t)$.

   Katok   proved in \cite[Corollary~4.4]{Katok1980} that for a $C^2$-diffeomorphism of a compact surface, or for a $C^2$-flow on  a compact $3$-manifold, its topological entropy is bounded above by   the exponent of the  rate of growth of its periodic orbits.  In particular, Katok's result implies:
\begin{thm}[Katok 1980]\label{thm-katok}
Let $\vp_t \colon M \to M$ be a $C^2$-flow on a compact 3-manifold. Suppose that $\vp_t$ has no periodic orbits, then $h_{top}(\vp_t) = 0$.
\end{thm}

Suppose that $\vp_t$ is an aperiodic flow obtained by inserting the  Kuperberg flow $\Phi_t$ on $\mK$ into a flow box for some flow on a compact $3$-manifold.
 Then $h_{top}(\vp_t) = 0$ by Theorem~\ref{thm-katok}, and thus   we have $h_{top}(\vp_t | \fM) =0$   for the flow $\vp_t$ restricted to the invariant set  $\fM$. By construction, the flows $\vp_t$ and $\Phi_t$ agree on $\fM$, so  that $h_{top}(\Phi_t | \fM) = h_{top}(\vp_t | \fM) = 0$. 
 That is, the topological entropy vanishes for the restriction of the flow $\Phi_t$ to the   space $\fM$ with the induced metric. 
In this section, we give a proof that   $h_{top}(\Phi_t | \fM) = 0$ based on the   dynamics of the restricted flow $\Phi_t \colon \fM \to \fM$ and the geometry of the invariant set $\fM$.
 
The idea of our   approach to the calculation  of $h_{top}(\Phi_t | \fM)$,   is to consider the   return map $\whPhi$ on the rectangle $\bRt$, and relate the entropy of the flow $\Phi_t$ restricted to $\fM$ with the entropy $h_{top}(\whPhi | \fMR)$ for $\whPhi$   restricted to its invariant set $\fMR = \fM \cap \bRt$. It  is standard to relate the entropy of the flow with the entropy of an induced return map to a section to the flow, assuming that the flow is everywhere transverse to the section.   However, the flow $\Phi_t$ is not everywhere transverse to $\bRt$ as it is tangent to $\bRt$ along its center line $\cT = \{z=0\} \cap \bRt$. As discussed below, these tangencies  result in discontinuities for the induced map $\whPhi$, which make the 
 relationship between $h_{top}(\Phi_t | \fM)$ and  the entropy $h_{top}(\whPhi | \fMR)$    more subtle. 
We subsequently relate the entropy $h_{top}(\whPhi | \fMR)$   with the entropy $h_{GLW}(\cGK^* | \fMR)$  for the  $\psg$ $\cGK^*$ defined in Definition~\ref{def-pseudogroup+}, and then show that $h_{GLW}(\cGK^* | \fMR) =0$.

We first  recall the definition   of    entropy for $\psg$s  as  introduced by   Ghys, Langevin and Walzcak in \cite[Section~2]{GLW1988}.
Let $(X, d_X)$ be a compact metric space, and $\cG_X^{(1)} = \{\vp_0 = Id, \vp_1, \vp_1^{-1}, \ldots, \vp_k, \vp_k^{-1}\}$ be a set of local homeomorphisms of $X$, with their inverses. That is, for each $1 \leq i \leq k$ there are open sets $U_i, V_i \subset X$ such that $\vp_i \colon U_i \to V_i$ is a homeomorphism. We require that  each map $\vp_i$ admits an extension to a homeomorphism $\ovp_i$ of  an open neighborhood of the closure $\oU_i$ of the domain in $X$. Let $\cG_X$ denote the      pseudogroup  generated by the collection of maps $\cG_X^{(1)}$, so that the axioms of Definition~\ref{def-pseudogroup} are satisfied.   
 Let $\cG_X^* \subset \cG_X$ denote the $\psg$ generated by the compositions of maps in $\cG_X^{(1)}$, so that $\cG_X$ satisfies conditions (1) to (4) of Definition~\ref{def-pseudogroup}, but not necessarily the condition (5).
  Let  $\cG_X^{(n)} \subset \cG_X^*$ be  the collection of maps defined by   the restrictions of  compositions of at most $n$ elements of $\cG_X^{(1)}$. 

 For $\e > 0$,  say that    $x, y \in X$ are \emph{$(\cG_X^*, n,\e)$-separated} if there exists $\vp \in \cG_X^{(n)}$ such that $x, y$ are in the domain of $\vp$ and  $\ds  d_X (\vp(x) , \vp(y))  > \e$. In particular, if the identity map is the only element of $\cG_X^{(n)}$ which contains both $x$ and $y$ in its domain, then they are $(\cG_X^*, n,\e)$-separated if  and only if $d_X(x,y) > \e$.

A finite set $E \subset X$ is said to be   \emph{$(\cG_X^*, n,\e)$-separated} if   each distinct pair $x , y \in E$  is $(\cG_X^*, n,\e)$-separated.
 Let $s(\cG_X^*, n,\e)$ be the maximal cardinality of a $(\cG_X^*, n,\e)$-separated subset of $X$.

 Define the  \emph{entropy} of $\cG_X^*$ by: 
\begin{equation}\label{eq-pseudogroupentropy}
h_{GLW}(\cG_X^*) ~ = ~ \lim_{\e \to 0} \left\{ \limsup_{n\to\infty}\frac{1}{n}\ln(s(\cG_X^*, n,\e)) \right\} .
\end{equation}
The property $0 < h_{GLW}(\cG_X^*) < \infty$ is independent of the choice of metric on $X$.
Note that the denominator in the expression \eqref{eq-pseudogroupentropy} is the word length $n$. Thus,  $h_{GLW}(\cG_X^*)$ differs from other definitions of entropy for group actions which scale by volume growth, rather than length. 

A key    observation in \cite{GLW1988} is  that  if the generators $\cG_X^{(1)}$ are the restrictions of $C^1$-diffeomorphisms defined on the compact closures of open subsets of  a manifold $M$, then $h_{GLW}(\cG_X^*)  < \infty$.  A discussion of some of the dynamical implications of $h_{GLW}(\cG_X^*)  > 0$ for  $C^1$-pseudogroups is given in \cite{CandelConlon2000,GLW1988,Hurder2014,Walczak2004}.

Given a subset $Y \subset X$, we can restrict consideration to   subsets $E \subset Y$ which are  $(\cG_X^*, n,\e)$-separated.  Let 
$s(\cG_X^*, Y, n,\e)$ be the maximal cardinality of such an $(\cG_X^*, n,\e)$-separated subset of $Y$. Define  the restricted entropy by
 \begin{equation}\label{eq-pseudogroupentropyE}
h_{GLW}(\cG_X^* | Y) ~ = ~ \lim_{\e \to 0} \left\{ \limsup_{n\to\infty}\frac{1}{n}\ln(s(\cG_X^*,Y, n,\e)) \right\} ~ \leq ~ h_{GLW}(\cG_X^*) ~.
\end{equation}

Before considering the entropy invariants associated to the return map $\whPhi$, we  consider      the points of \emph{discontinuity}  for its powers $\whPhi^n$,  and use this to give a description of the domains of continuity of the maps $\whPhi^n$.  We restrict   to   the region $\bRt \cap \{r \geq 2\}$ which contains the set $\fMR$. In Section~\ref{sec-pseudogroup}, the  continuity properties for the induced return map $\whPsi$ of the Wilson flow were  analyzed, which resulted in the description of its domains of continuity in  \eqref{defPsi}. The analysis of the powers of the map $\whPhi$ extends this analysis.

Recall that by Condition (K3) in Section~\ref{sec-kuperberg}, the compact annular region $\cA(2) = \{(r,\theta, 0) \mid r \geq 2\} \subset \cA$ is disjoint from the images of the insertions $\sigma_i$ for $i=1,2$,  thus the vector fields $\cW$ and $\cK$ agree on an open neighborhood of $\cA(2)$ and so 
the flows $\Psi_t$ and $\Phi_t$   agree near $\cA(2)$   for $t$ near to $0$.

Consider a point $\xi_0 \in Dom(\whPhi) \cap \{r \geq 2\}$ with $\xi_1 = \whPhi(\xi_0)$. Then there exists a $\cK$-orbit segment $[\xi_0,\xi_1]_{\cK}$    which intersects $\bRt$ only in its endpoints. If the $\cK$-flow $\Phi_t(\xi_0)$ intersects $\bRt$ transversally at $\xi_1$, then there is an open neighborhood in $\bRt$ of $\xi_0$ consisting of points whose flow is also transverse to $\bRt$ near $\xi_1$, and thus $\whPhi$ is continuous at $\xi_0$. 
On the other hand, if the   $\cK$-flow $\Phi_t(\xi_0)$ intersects $\bRt$ tangentially  at $\xi_1$, then in every open neighborhood of $\xi_0$ there is some point whose flow does not intersects $\bRt$ near $\xi_1$, and  the map $\whPhi$ will be discontinuous at $\xi_0$.
   Thus,   $\xi_0$ is   a point of discontinuity for $\whPhi$ precisely when it  is contained in  the inverse image under $\whPhi$ of the center line  $\cT = \{x \in \bRt \mid z(x)=0\}$. That is,   the set of discontinuities  for   $\whPhi$ is the set   $\whPhi^{-1}(\cT)$. Similarly,     the set of discontinuities  for the inverse map $\whPhi^{-1}$ is the set   $\whPhi(\cT)$.
 
These remarks for $\whPhi$ generalize to its powers as follows. Let  $n > 1$ and  $\xi_0 \in   Dom(\whPhi^n)$, then set $\xi_{\ell} = \whPhi^{\ell}(\xi_0)$ for $1 \leq \ell \leq n$.
 If the $\cK$-orbit segment  $[\xi_0,\xi_n]_{\cK}$  intersects $\bRt$   transversally at each point $\xi_{\ell}$  then there is an open neighborhood  of $\xi_0$ in $\bRt$  consisting of points whose flow is also transverse to $\bRt$, and thus $\whPhi^n$ is continuous at $\xi_0$.
 On the other hand, if the $\cK$-orbit segment  $[\xi_0,\xi_n]_{\cK}$  intersects $\bRt$ tangentially at some point $\xi_{\ell}$,  that is $\xi_{\ell} \in \cT \subset \bRt$, then in every open neighborhood of $\xi_0$ there is some point whose flow does not intersects $\bRt$ near $\xi_{\ell}$, and thus $\xi_0$ is a point of discontinuity for $\whPhi^n$.  For each integer $\ell$, introduce the   subset of $Dom(\whPhi^{\ell})$, 
\begin{equation}
\cL_{\ell} = \whPhi^{-\ell}(\cT) = \whPhi^{-\ell}(\cT \cap Dom(\whPhi^{-\ell}))~ .
\end{equation}

The above discussion implies that for $n \geq 1$, $\xi \in Dom(\whPhi^n)$  is a point of discontinuity for $\whPhi^n$ precisely when there exists $1 \leq \ell \leq n$ such that  $\xi \in \cL_{\ell}$. The analogous statement also holds for $\xi \in Dom(\whPhi^{-n})$ with $n > 0$ and $-n \leq \ell < 0$.
For $n \geq 1$, define  
\begin{equation}\label{eq-discontinuities}
 D^{(n)}(\whPhi)  =    \bigcap_{-n \leq \ell \leq n} ~ \left\{ Dom(\whPhi^{\ell}) - \cL_{\ell} \right\} ~ \cap ~ \{r \geq 2\} ~.
\end{equation}
Let $U \subset  D^{(n)}(\whPhi)$ be a connected component, then $\whPhi^{\ell}|U$ is continuous for all $-n \leq \ell \leq n$.
Thus, we are interested in  the  finite collection of     domains of continuity, 
 \begin{equation}\label{eq-setofdiscontinuities}
\cD_{\fM}^{(n)}(\whPhi)   =    \{ U \subset D^{(n)}(\whPhi) \mid U ~ {\rm is
  ~ a  ~ connected  ~ component ~of}~D^{(n)}(\whPhi) ~{\rm such ~that}~ U \cap \fMR \ne \emptyset \} .
\end{equation}
 Corollary~\ref{cor-min=min} implies that $\fMR \subset Dom(\whPhi^n) \cap  Dom(\whPhi^{-n})$ for all $n$. 
 For $U \in \cD_{\fM}^{(n)}(\whPhi)$, let $\oU \subset \bRt$ denote its closure, then we have:
\begin{lemma}\label{lem-domains}
For each $n \geq 1$, $\ds \fMR \subset \bigcup ~ \{ \oU \mid U \in \cD_{\fM}^{(n)}(\whPhi)\}$.
\end{lemma}

The collection  $\cD_{\fM}^{(n)}(\whPhi)$ thus gives a decomposition of $\whPhi$ into sets on which the function has nice regularity properties. The number of sets in these partitions is an important property of the flow $\Phi_t$. 
   \begin{defn}\label{def-complexity}
 The \emph{complexity  function} of  $\whPhi$ is defined by  $C_{\Phi}(n) = \#\cD_{\fM}^{(n)}(\whPhi)$. That is, for 
      $n \geq 1$,   $C_{\Phi}(n)$ is the number of   connected components of $D^{(n)}(\whPhi)$.
 \end{defn}
 We   develop a geometric interpretation of the sets $\cL_{\ell}$ and $\ds U \in \cD_{\fM}^{(n)}(\whPhi)$, in terms of the flow $\Phi_t$ restricted to  the space $\fM_0$. We use this to  relate the    function $C_{\Phi}(n)$    to the function $\# \cM(n)$,  where the collection $\cM(n) \subset \cGK^{(n)}$ was  defined in Proposition~\ref{prop-standardform}.

 We require some preliminary observations. 
 We first obtain an estimate   relating   the exponent $\ell$ and the time $T$, for the maps $\whPhi^{\ell}$ and the flow $\Phi_T$ restricted to $\fM$.  It is straightforward to obtain  such an estimate      for a complete transversal to the flow $\Phi_t$ with return times bounded away from zero, but the estimate is more subtle for the return map  $\whPhi$.  
 \begin{lemma}\label{lem-returntoR}
There exist  $L_{\Phi} > 0$  so that for $x \in \fMR$, and $T > 0$ such that   $\Phi_T(x) \in \fMR$, and  $n_{x} > 0$ for which $\whPhi^{2n_{x}}(x) = \Phi_T(x)$, then  $0 < L_{\Phi} \, n_{x} \leq  T$. 
\end{lemma}
\proof
We use the properties of the return map $\whPhi$ established in Section~\ref{sec-pseudogroup}.

Let $x \in \fMR$ with $r(x) \geq 2$, $z(x) < 0$ and $z(\whPhi(x)) < 0$. If  the $\cK$-orbit of $x$ does not intersect an insertion region $\cL_i^-$ it must   complete a revolution around the core cylinder $\cC$ before returning to $\bRt$, hence the segment $[x,\whPhi(x)]_{\cK}$ has length bounded below by $L_{\Phi}' = 4\pi$. If  the $\cK$-orbit of $x$ does   intersect an insertion region $\cL_i^-$ then the segment $[x,\whPhi(x)]_{\cK}$ has length at least $L_{\Phi}''$ for some $L_{\Phi}'' > 0$. 
Let $L_{\Phi} = \min \{L_{\Phi}' , L_{\Phi}''\}$.

In the case where $x \in \fMR$ with  $z(x) > 0$ and $z(\whPhi(x)) > 0$, then the segment $[x,\whPhi(x)]_{\cK}$ traverses $\mK$ in the clockwise direction, but is otherwise analogous to the  case above. We assume that the estimates $L_{\Phi}'$ and $L_{\Phi}''$ also apply in this case. 

In the case where $x \in \fMR$ with   $z(x) \leq 0$ and $z(\whPhi(x)) \geq 0$, then the flow $\Phi_t(x)$ crosses the annulus $\{z =0 \}$. In particular, $\whPhi(x) = \psi(x)$ for the map $\psi$ of the pseudogroup $\cGK$. Recall that the action of $\psi$ was described in  Lemma~\ref{lem-psi_action}, and we observe that  the length of the segment $[x,\whPhi(x)]_{\cK}$ need not have a positive lower bound. But this segment must be preceded and followed by a segment of the above two types.
It follows that the composition $\whPhi \circ \whPhi$ applied to $x \in \fMR$ is realized by a $\cK$-orbit segment of length at least $L_{\Phi}$, for any $x \in \fMR$.

It follows that $[x,\whPhi^{2n_{x}}(x)]_{\cK}$ has length at least $L_{\Phi} \cdot n_{x}$.  
\endproof

   Recall that by Propositions~\ref{prop-r_b} and \ref{prop-bubbles1}, the constant $r_b > 2$, introduced in    Definition~\eqref{eq-rb2}, provides a lower bound for the radius of points in bubbles of the first type. Also, Lemma \ref{lem-unknown} provides a constant $r_B>2$ for which bubbles of the second type do not arise.
Finally,   $r_e$ is the exceptional radius defined in \eqref{eq-eradius}. Let $r_u = \min\{r_b, r_B, r_e\}$, so that $2 < r_u \leq r_e$. Observe that the techniques in the proof of Proposition~\ref{prop-syndeticgk*} apply    for any  value of the radius upper bound which is less than $r_e$. In particular, we   conclude   the following:

\begin{prop}\label{prop-syndeticanyr}
Let $\xi \in \bRt$ have   infinite orbit for the flow $\Phi_t$ with
$r(\xi)<r_u$ and  $\rho_{\xi}(t)\geq 2$ for all $t$. Then the set $\ds
\cS^*(\xi, r_u)  = \{ s \mid \Phi_s(\xi) \in \cGK^*(\xi)\subset \bRt ~, ~ r(\Phi_s(\xi)) < r_u\}$ is syndetic in $\mR$ for a constant $\nu_u^*$ which is independent of $\xi$.
\end{prop}

 \begin{cor}\label{cor-syndeticanyr}
 Let $\xi \in \fMR = \fM \cap \bRt$. Then  $\ds
\cS^*(\xi, r_u) = \{ s \mid \Phi_s(\xi) \in \cGK^*(\xi)\subset \bRt ~, ~ r(\Phi_s(\xi)) < r_u\}$ is syndetic in $\mR$ for the constant $\nu_u^*$.
 \end{cor}
 \proof   
Theorem~\ref{thm-minimal} implies that the $\Phi_t$-orbit of $\xi$ contains points arbitrarily close to the special point $\omega_1 \in \bRt$, and thus there exists $t$ such that with $2 \leq r(\Phi_t(\xi))  < r_u$. The claim then follows from Proposition~\ref{prop-syndeticanyr}. 
 \endproof

\begin{cor}\label{cor-returntoe}
 Let $M_u > 0$ be the greatest integer with $M_u \leq 2\nu_u^*$. 
For  each    $x\in\fMR$ there exists $ 0 < \ell_{x}', \ell_x'' \leq M_u$ such that  $r(\whPhi^{\ell_x'}(x))<r_u$ and $r(\whPhi^{-\ell_x''}(x))<r_u$.
\end{cor}
 \proof
This follows from Lemma~\ref{lem-returntoR} and Corollary~\ref{cor-syndeticanyr}. 
  \endproof

  We next characterize the discontinuities of $\whPhi^n$ on $\oU$, for $U \in \cD_{\fM}^{(n)}(\whPhi)$. We first reduce this problem to a ``standard form''.
Recall from \eqref{eq-intervals}  that the  interval $I_0 = \cR \cap \bRt$ is the closed vertical segment  whose $\Phi_t$-flow generates $\fM_0$, and the intersection of $\fM_0$ with $\bRt$ yields the connected components of $\fM_{\bRt}$. 
Also,  the $\Phi_t$-flows  of the closed vertical segments $N_0 =  J_0 \cup  I_0 \subset \bRt$ and $M_0 = I_0 \cup K_0 \subset \bRt$ generate the double propellers discussed in Section~\ref{sec-doublepropellers}. The $\Phi_t$-flows of $N_0$ and $J_0$ generate   families of nested double propellers, whose intersections with $\bRt$   are the topological circles illustrated in Figures~\ref{fig:ellipsespropeller}, \ref{fig:curvesR0} and \ref{fig:bubblesL-}. 
  By the nesting properties of the double propellers,    the discontinuities for $\whPhi^n$ restricted to  $\fMR$ are determined by those for $\whPhi^n$ restricted to $I_0$.

Suppose that $x \in I_0$ is a point of continuity for $\whPhi^n$ with $n \geq 1$, so the   restriction $\whPhi^n$ to $I_0$ is continuous in an open neighborhood of $x$. 
Let $x',x'' \in I_0$ be such that the open connected segment $(x',x'') \subset I_0$ is the maximal open subset of $I_0$ containing $x$ 
such that the restriction $\whPhi^n | (x',x'')$ is continuous. Set $\cJ(x,0)= (x',x'')$.
Then  for each $z \in \cJ(x,0)$ and $w = \whPhi^n(z)$, the $\cK$-orbit segment $[z,w]_{\cK}$   intersects $\bRt$ transversally, while for
  the endpoints $\{x',x''\}$, the $\cK$-orbit segments $[x', \whPhi^n(x')]_{\cK}$ and $[x'', \whPhi^n(x'')]_{\cK}$ must each contain a point of tangency with $\cT$.  The collection of such points of tangency  for the flow $\Phi_t$     partition $I_0$ into maximal subintervals on which the restriction of $\whPhi^n$ is continuous. 
 
  We   develop an estimate on the number of intervals of continuity for $\whPhi^n$ in $I_0$, where we consider the case $n > 0$ first.
Let $x \in I_0$ be a point of continuity. As $r(x) =2$, by Corollary~\ref{cor-returntoe}, there exists an infinite sequence $\{0 < \ell_1 < \ell_2 < \cdots\}$ with $(\ell_{i+1} - \ell_i) < M_u$ for all $i \geq 1$, such that for   $x_i = \whPhi^{\ell_i}(x)$ we have $r(x_i) < r_u$. 
Set   $\cI(x,0) = I_0$, and for $ i > 0$ let $\cI(x, i) \subset \fM_0 \cap \bRt$ denote the maximal connected    component containing  $x_i$.
Then $\cI(x, i)$  is a simple arc by the choice of $r_u$. In particular,  the endpoints of $\cI(x,i)$ are contained in the boundary of $\fM_0$.    
Let $p_0^1(\cI(x,i))$ denote the endpoint of $\cI(x,i)$ in the region $\{z < 0\}$, and     $p_0^2(\cI(x,i))$  denote the endpoint of $\cI(x,i)$ in the region $\{z > 0\}$. In particular,   $p_0^1(\cI(x,0)) = \omega_1$ and $p_0^2(\cI(x,0)) = \omega_2$ are the special points defined in \eqref{eq-omegas}.
 
For  $i > 0$,  denote $\cJ(x,i) = \whPhi^{\ell_i}(\cJ(x,0)) \subset \cI(x,i)$. Then   $\cJ(x,i)$ is an open subinterval,   disjoint from the center line $\cT$. 
 For each $i \geq 0$,  if $\cJ(x,i)$ lies below the line $\cT$, let $\zeta_{i} =  p_0^1(\cI(x,i)) \in \bRt$ denote the lower endpoint of the interval $\cI(x,i)$. 
If $\cJ(x,i)$ lies above the line $\cT$, then  $\zeta_{i} =  p_0^2(\cI(x,i)) \in \bRt$ denotes the upper endpoint of the interval $\cI(x,\ell)$.   
Then Proposition~\ref{prop-syndeticanyr} implies that for each $i \geq 1$ and  $n_i = \ell_{i} - \ell_{i-1}$ the  restriction 
\begin{equation}\label{eq-intervalshift}
\whPhi^{n_i} \colon \cJ(x,i-1) \to \cJ(x,i)
\end{equation}
corresponds to a map $\vp(x,i) \in \cGK^{(n_i)}$ which maps $\zeta_{i-1}$ to $\zeta_{i}$. 
 Let $\Upsilon(x, i) = \vp(x,i) \circ \cdots \circ \vp(x,2)   \circ \vp(x,1) \in \cGK^{(\ell_i)}$.

We next give a geometric interpretation to the maps in $\cGK^*$ defined by the actions \eqref{eq-intervalshift},   in terms of the action of $\cGK^*$ on $\TP$ defined in Section~\ref{sec-normal}. 

The  intersection $\fM_0 \cap \bRt$ consists of simple arcs arising from the intersections with propellers, as discussed in Section~\ref{sec-proplevels}, and double arcs or ellipses which arise from the intersection of $\bRt$ with bubbles in $\fM_0$ as discussed in Section~\ref{sec-bubbles}.
 The simple arcs   in the intersection of $\fM_0 \cap \bRt$ correspond to   transverse line segments to the propellers, as illustrated in    Figure~\ref{fig:fM0cT}.  The double arcs resulting from the intersections with bubbles are   illustrated in Figure~\ref{fig:choufleurbubbles}, but for simplicity are not pictured in Figure~\ref{fig:fM0cT}.   In particular, each arc   $\cI(x,i) \subset \fM_0 \cap \bRt$ is   an arc transverse to the flow $\Phi_t$ in $\fM_0$. 

 For each (simple) propeller $P_{\gamma}$ or $P_{\lambda}$  attached in the level decomposition of  $\fM_0$, as described in Section~\ref{sec-intropropellers},  the intersections  $P_{\gamma} \cap \cA$ and $P_{\lambda} \cap \cA$ form the center line segments of the propeller, as illustrated in     Figure~\ref{fig:fM0cT}.
  The union of all these lines form the subset 
$\TP'' \subset \TP' = \fM_0 \cap \cA$, as defined in Section~\ref{sec-normal}. The   intersections of $\TP''$ with $\cT = \bRt \cap \cA$ are the vertices in the set $\TP''$, as    illustrated  in   Figure~\ref{fig:fM0cT}. 

For each $i \geq 1$,  the map $\vp(x,i)$ induces a map $\overline{\vp(x,i)}$ from the vertex of $\TP$ defined by $\cI(x,i-1) \cap \cT$ to the vertex defined  by $\cI(x,i) \cap \cT$. Thus the map $\Upsilon(x, i) \colon \cJ(x,0) \to \cJ(x,i)$ maps an open set  in the simple arc $I_0$ to an open set in the simple arc $\cI(x,i)$.
We have  $\ell_i = n_1 + \cdots + n_i$ hence $\Upsilon(x, i) = \whPhi^{\ell_i}|\cJ(x,0)$, so that $\cJ(x,0)$ is contained in the domain of $\whPhi^{\ell_i}|I_0$.

We return to the problem of estimating the number of domains of continuity for $\whPhi^n | I_0$.
Let $x \in I_0$ be  a point of continuity, let the indices $\ell_i$ be defined as above, then   let  $i_n > 0$ be the least index $i$ such that  $\ell_{i_n} \geq n$. Note that $i_n \leq n$,  and   Corollary~\ref{cor-syndeticanyr} implies that $\ell_{i_n} \leq n + M_u$.

By  Proposition~\ref{prop-standardform},   $\Upsilon(x, i_n)$ can be written in normal form $\Upsilon(x, i_n) = \vp^+ \circ \vp^- \in \cGK^{(\ell_{i_n})}$, where appropriate subwords in $\Upsilon(x, i_n) $ are replaced by shortcuts, or powers of the generator $\psi$.  The interval $\cI(x,i_n)$ is thus defined by the word $\vp^+ \circ \vp^-$. 
  Proposition~\ref{prop-asymptotic} implies  the number of   normal forms possible for word length $\ell_{i_n}$ is bounded by $\#\cM(\ell_{i_n}) \leq \#\cM(n + M_u)$.  By Corollary~\ref{cor-subexponentialgrowth}, the  function $\# \cM(n + M_u)$
has subexponential growth as a function of $n$. It follows that the total number of simple arcs in $\fM_0 \cap \bRt$ which are in the images of the map $\whPhi^n | I_0$ is bounded above by  $\# \cM(n + M_u)$.

 Let $x \in I_0$ be  a point of continuity, and continuing the notation as above, 
for $y = \whPhi^n(x)$,   the $\cK$-orbit segment $[x, y]_{\cK}$ is   transverse to $\bRt$ at every point of intersection, and the orbit contains a point in each of the sequence of simple arcs $\{\cI(x,i) \mid 0 \leq i \leq i_n\}$.  

 By the choice of the sequence $\{\ell_i\}$, for each $1 \leq i \leq i_n$,  the forward $\cK$-flow in $\fM_0$  of the points in the interval $\cJ(x,i-1)$ to the interval  $\cJ(x,i)$ has length at most $\nu_u^*$.
 Therefore, there is a uniform upper bound $M_u'$  on the number of points   $\xi \in \cJ(x,i-1)$ 
 for which the orbit segment $[\xi, \whPhi^{n_i}(\xi)]_{\cK}$ is  tangent to $\cT$.

 Recall that  $\{x',x''\} \in \cI(x,0)$ are the boundary points of the  interval $\cJ(x,0) \subset \cI(x,0)$, and
  the $\cK$-orbit segments $[x', \whPhi^n(x')]_{\cK}$ and $[x'', \whPhi^n(x'')]_{\cK}$  each contains a point of tangency with $\cT$,  
 which must be one of the points of tangency for the $\Phi_t$-flow between $\cJ(x,i-1)$ and $\cJ(x,i)$ for some $1 \leq i \leq i_n$. 
 Thus, the number of possible   points of tangency with $\bRt$ which can define the lower boundary point $x'$ for $\cJ(x,0)$ has an upper bound estimate  $i_n \cdot M_u' \leq n\cdot M_u'$, and similarly for the upper boundary $x''$.

Thus, if $\cJ \subset I_0$ is a maximal domain of continuity for $\whPhi^n | I_0$, then there exists a connected component    $\cI(U,n) \subset \fM_0 \cap \bRt$ for which  $\whPhi^n(\cJ) \subset \cI(U,n)$. The number of possible such domains is bounded above by  $\# \cM(n + M_u)$. Then the number of maximal connected intervals $\cJ \subset I_0$ for which $\whPhi^n(\cJ) \subset \cI(U,n)$ is bounded above by the number of possible endpoints for the domains $\cJ$, which   is bounded above by  $(n\cdot M_u')^2$. It follows that   the number of domains of continuity for $\whPhi^n$ restricted to $I_0$   has subexponential growth.

    Similar considerations apply to the analysis of $\whPhi^{-n} | I_0$, so that the number of domains of continuity for the collection of maps $\{\whPhi^{\ell}|I_0 \mid -n \leq \ell \leq n \}$ is a function of  subexponential growth in $n$.
    We thus obtain:
    
       \begin{prop}\label{prop-complexity}
The   function $C_{\Phi}(n)$  has subexponential growth for    a generic Kuperberg flow $\Phi_t$.
 \end{prop}

We  now return to the definition and calculation of the entropy invariants associated with $\whPhi$.   Recall from Section~\ref{sec-pseudogroup} that $\cGK$ is the pseudogroup  generated by   $\whPhi$ acting  on    $X = \bRt$.

 Let  $V_+ =  \{ x \in Dom(\whPhi) \mid z(x) > 0\}$, let $V_-$ be the region of $Dom(\whPhi)$ which lies below the curve $\cL_1$ and let $V_0$ be the region between $V_-$ and $V_+$. Then set 
$$\whPhi_- ~ \equiv ~  \whPhi | V_- \quad , \quad \whPhi_0 ~ \equiv ~  \whPhi | V_0 \quad , \quad \whPhi_+ ~ = ~ \whPhi | V_+ ~ .$$
Let $\cG_{\whPhi}^* \subset \cGK$ be the $\psg$ defined as in Definition~\ref{def-pseudogroup+} by the maps $\{\whPhi_- , \whPhi_0 , \whPhi_+\} \subset \cGK$. 
Also, introduce the symmetric generating set   for $\cG_{\whPhi}^*$, 
\begin{equation}\label{eq-gensetforPwhPsi}
 \cG_{\whPhi}^{(1)} = \{Id, (\whPhi_-)^{\pm 1}, (\whPhi_0)^{\pm 1}, (\whPhi_+)^{\pm 1}\} ~.
\end{equation}
Then  a map  $\vp \in \cG_{\whPhi}^{(n)}$ is a composition of at most $n$ maps in the collection $\cG_{\whPhi}^{(1)}$, which by definition of these generating maps is just a composition of $\whPhi$, or its inverse, restricted to a connected domain of continuity.
 In particular, for   each open set $U \in  \cD_{\fM}^{(n)}(\whPhi)$, and $-n \leq \ell \leq n$, the restriction    $\vp = \whPhi^{\ell} | U \in \cG_{\whPhi}^{(n)}$.

Recall that  $d_{\bRt}$ denotes the   metric on $\bRt$ defined in \eqref{eq-metricR}, and     the metric $d_{\mK}$ on $\mK$ was defined  in   Section~\ref{sec-radius} so that $\cK$ has unit length, and thus $\Phi_t$ is a unit speed flow. 

We   restrict  the action of  $\cG_{\whPhi}^*$  to the invariant set $\fMR \subset \bRt$.  Let 
$s(\cG_{\whPhi}^*, \fMR, n,\e)$ be the maximal cardinality of   a $(\cG_{\whPhi}^*, n,\e)$-separated subset of $\fMR$. Define  the restricted entropy by
 \begin{equation}\label{eq-pseudogroupentropyfMR}
h_{GLW}(\cG_{\whPhi}^* | \fMR) ~ = ~ \lim_{\e \to 0} \left\{ \limsup_{n\to\infty}\frac{1}{n}\ln(s(\cG_{\whPhi}^*,\fMR, n,\e)) \right\} ~.
\end{equation}

 Here is the first result.
 
 \begin{prop}\label{prop-entropywhPhi}
$h_{GLW}(\cG_{\whPhi}^* | \fMR)  = 0$ if and only if  $h_{top}(\Phi_t | \fM) = 0$.
\end{prop}
 \proof

   We require   some   preliminary   remarks. 
   As  $\bRt$ is compact submanifold of $\mK$, there exists a constant  $A_1 \geq 1$ so that for $x, y \in \bRt$  we have
\begin{equation}\label{eq-entbounds2}
A_1^{-1} \, d_{\mK}(x,y)  ~ \leq d_{\bRt}(x, y) \leq A_1 \, d_{\mK}(x,y) .
\end{equation}

The submanifold $\bRt$  was chosen to be disjoint from the embedded regions $\cD_i \subset \mW$ for $i=1,2$, so there exists some $\e_2 > 0$ so that the $\e_2$-closed neighborhood of $\bRt$ is   disjoint from these surfaces. That is, for $i=1,2$, we have
$$C_{\mK}(\bRt, \e_2) = \{x \in \mK \mid d_{\mK}(x, \bRt) \leq \e_2 \} \quad , \quad C_{\mK}(\bRt, \e_2)  \cap \cD_i = \emptyset .$$
Then  the flows $\Phi_t$ and $\Psi_t$ agree on the set $C_{\mK}(\bRt, \e_2)$. 

Define a function $\Delta(\e) > 0$ for $\e > 0$ as follows. For $y \in \fMR$ let $s_y' < 0 < s_y''$ be the minimum and maximum values such that 
$\{\Phi_t(y) \mid s_y' \leq t \leq s_y''\} \subset C_{\mK}(\bRt, \e_2)$. Then for $x \in \fMR$ with $x \ne y$ let 
$\Delta(x,y,\e_2) = \min \, \{d_{\mK}(x,\Phi_t(y)) \mid  s_y' \leq t \leq s_y''\}$. Then $0 < \Delta(x,y,\e_2) \leq d_{\mK}(x,y)$, and we set
\begin{equation} \label{eq-Deltaepsilon}
\Delta(\e) = \min \, \left\{ \Delta(x,y,\e_2) \mid x, y \in \fMR ~ {\rm and} ~ d_{\bRt}(x, y) \geq \e \right\}
\end{equation}
Observe that $\Delta(x,y,\e_2)$ is a continuous function of $x\ne y$, and  the set $\{(x,y) \mid x, y \in \fMR ~ {\rm and} ~ d_{\bRt}(x, y) \geq \e\}$ is compact, so $0 < \Delta(\e)  \leq \e$.

Recall that $\nu_{\cK}$ is the constant introduced in Proposition~\ref{prop-sigmadyn}.
As the flow $\Phi_t$ is smooth,     there exists $A_2 > 1$ such that  for   $x,y \in \fM$ and   $-\nu_{\cK}  \leq t\leq \nu_{\cK}$,   
\begin{equation}\label{eq-entbounds1}
A_2^{-1} \, d_{\mK}(x,y)  ~ \leq d_{\mK}(\Phi_t(x),\Phi_t(y)) \leq A_2 \, d_{\mK}(x,y) .
\end{equation}

 Now, assume that $h_{GLW}(\cG_{\whPhi}^* | \fMR)  = \lambda > 0$.   
 Then   for $0 < \e_3 \leq \e_2$ sufficiently small,  and  for   $n$ sufficiently large, there is a $(\cG_{\whPhi}^*, n,\e_3)$-separated subset   $\whE_n \subset \fMR$ with cardinality $\# \whE_n > \exp(n \, \lambda/2)$. Fix $n$ large enough so that this estimate holds.

  We have by definition \eqref{eq-discontinuities} and Lemma~\ref{lem-domains} that for each $0\leq \ell \leq n$,
$$Dom(\whPhi^{\ell}) \cap \fMR \subset \overline{D^{(\ell)}(\whPhi)} \cap \fMR \subset \bigcup ~ \left\{ \oU \mid U \in \cD_{\fM}^{(\ell)}(\whPhi) \right\}  .$$   
 Recall that   $C_{\Phi}(n)$ is the number of connected domains in $\cD_{\fM}^{(n)}(\whPhi)$, so by   the Pigeonhole principle, there exists some connected component $U \in \cD_{\fM}^{(n)}(\whPhi)$ such that $\#  (\whE_n \cap \oU) > \exp(n \, \lambda/2)/C_{\Phi}(n)$. For $n$ sufficiently large, we can assume that $C_{\Phi}(n) \leq \exp(n \, \lambda/4)$ by   Proposition~\ref{prop-complexity}.  Then for such $n$, set  $\whE_n' = \whE_n \cap \oU$, and we have that $\# \whE_n' > \exp(n \, \lambda/4)$.     Let $y\neq x$ in $\whE_n'$, then $x$ and $y$ are $(n,\e_3)$-separated implies there exists 
 $-n\leq \ell_{x,y}\leq n$ such that $d_{\bRt}(\whPhi^{\ell_{x,y}}(x),\whPhi^{\ell_{x,y}}(y))>\e_3$. 

 It may happen that the restricted set  $\whE_n'$ satisfies $\whE_n' \cap (\oU - U) \ne \emptyset$. That is, there exists  $x\in \cL_\ell$ for some $-n\leq \ell\leq n$. For such $x$, there exists $x' \in U$ which is sufficiently close so that the $(n,\e_3)$-separated condition is again satisfied. Thus, by taking a sufficiently small perturbation of each point in $\whE_n'$, we can assume without loss of generality that $\whE_n' \subset U$. 
Recall the definition of $A_1$ in \eqref{eq-entbounds2}, then 
 for  each distinct pair $x,y \in \whE_n'$ there exists $-n  \leq \ell_{x,y} \leq n$ such that 
 \begin{equation}\label{eq-newseparation}
d_{\mK}(\whPhi^{\ell_{x,y}}(x) , \whPhi^{\ell_{x,y}}(y)) > \e_3/A_1  .
\end{equation}

By  Proposition~\ref{prop-sigmadyn}, the forward and backward  return time to $\bRt$  for the $\cK$-orbit of any point  $x \in \fMR$ is bounded   by  $\nu_{\cK}$. Consequently, for any $x \in \fM$ there exists a least $0 \leq t_x < \nu_{\cK}$ such that $\Phi_{t_x}(x) \in \fMR$.

   Set $T_n = n \, \nu_{\cK}$.   
For a distinct pair $x,y \in \whE_n'$ let $-n  \leq \ell_{x,y} \leq n$ be as above, and let $t_x$ be the value such that $x' = \Phi_{t_x}(x) = \whPhi^{\ell_{x,y}}(x)$, and $t_y$ such that $y' = \Phi_{t_y}(y) = \whPhi^{\ell_{x,y}}(y)$. 
  Then $0 < |t_x| , |t_y| \leq \ell_{x,y} \, \nu_{\cK} \leq T_n$.    
  
  Set  $\e_4 = \e_3/A_1$, and $\e_5 = \Delta(\e_4) \leq \e_4$ for the function $\Delta(\e)$ defined by   \eqref{eq-Deltaepsilon}.

  We claim that the set $\whE_n' \subset \fM$ is $(\Phi_t, T_n,\e_5)$-separated. If not, then there exists a distinct pair $x,y \in \whE_n'$ such that
  $d_{\mK}(\Phi_t(x), \Phi_t(y)) \leq \e_5$ for all $-T_n \leq t \leq T_n$, and this holds for $t = t_x$ and $t = t_y$ in particular. 
  It is given that $d_{\mK}(x', y') > \e_4$, so by the definition of $\e_5$  and the function $\Delta(\e)$, we have $\Delta(x,y,\e_2) \geq \e_4$. That is, for all $t$ with $s_y' < t < s_y''$ we have that $\ds  d_{\mK}(x',\Phi_t(y')) \geq \e_4$. If $s_y' < t_1 - t_2 < s_y''$ then this implies that 
   $d_{\mK}(\Phi_{t_1}(x), \Phi_{t_1}(y)) \geq \e_4 \geq \e_5$ contrary to assumption. If $s_y' < t_1 - t_2 < s_y''$ is not satisfied, then 
   $\Phi_{t_1}(y) \not\in C_{\mK}(\bRt, \e_2)$, so that  $d_{\mK}(\Phi_{t_1}(x), \Phi_{t_1}(y)) \geq \e_2 \geq \e_5$, again contrary to assumption.     

It   follows that  $s(\Phi_t | \fM , T_n,\e_5) \geq \exp(n \, \lambda/4)$. 

As $T_n$ is a linear function of $n$, this implies that $h_{top}(\Phi_t | \fM) > 0$, as was to be shown.

 The reverse implication  of Proposition~\ref{prop-entropywhPhi} follows by  showing  
 that    $h_{GLW}(\cG_{\whPhi}^* | \fMR)  =   0$ implies    $h_{top}(\Phi_t | \fM) = 0$.  For this purpose, we
 use an alternate interpretation of the topological entropy $h_{top}(\Phi_t | \fM)$, in terms of the  {minimum cardinality} of a $(\Phi_t, T,\e)$-{spanning set} for $\fM$. The properties of spanning sets used to define  the entropy of a map are  discussed in greater detail in \cite[Chapter~7.2]{Walters1982}.  
 
 For $T \geq 0$, introduce the metric $d_{\mK}^T$ on $\fM$ where for $x, y \in \fM$, 
 $$d_{\mK}^T(x,y) = \max \, \{d_{\mK}(\Phi_t(x), \Phi_t(y)) \mid -T \leq t \leq T\} ~ .$$
 
 A subset $F \subset \fM$ is said to be $(\Phi_t | \fM, T, \e)$-\emph{spanning} if for all $y \in \fM$ there exists $x \in F$ such that $d_{\mK}^T(x,y) \leq \e$.
 Note that for $T' > T$ the metric $d_{T'}$ is finer than $d_T$, so a $(\Phi_t | \fM, T', \e)$-spanning set $F$ is also a $(\Phi_t | \fM, T,\e)$-spanning set.
 Let $r(\Phi_t | \fM, T,\e)$ be the minimal cardinality of a $(\Phi_t | \fM, T,\e)$-spanning set in $\fM$. Then $r(\Phi_t | \fM, T,\e) < \infty$ as $\fM$ is compact, and we have the inequalities 
 $r(\Phi_t | \fM, T', \e) \geq r(\Phi_t | \fM, T,\e)$ for $T' > T$, and $r(\Phi_t | \fM, T,\e') \geq r(\Phi_t | \fM, T,\e)$ for $0 < \e' < \e$.
 Then we have 
\begin{equation}\label{eq-spanning}
h_{top}(\Phi_t | \fM)  = \lim_{\e\to  0} \left\{ \limsup_{T\to\infty}\frac{1}{T}\log(r(\Phi_t | \fM, T,\e)) \right\} .
\end{equation}
As an application, we have:
\begin{lemma}\label{lem-spanningbound}
Suppose that for all $\lambda > 0$, we have that:  for $\e >0$ sufficiently small and all $T > 0$ sufficiently large, there exists a $(\Phi_t | \fM, T, \e)$-spanning set $F(T,\e) \subset \fM$ with $\# \, F(T,\e)  \leq \exp (T \, \lambda)$.  Then $h_{top}(\Phi_t | \fM) =0$.
\end{lemma}

Now, assume  that $h_{GLW}(\cG_{\whPhi}^* | \fMR)  =   0$, and let   $\lambda > 0$. Then  we construct 
   $(\Phi_t | \fM, T, \e)$-spanning sets  which have growth rate less than $\exp(T \, \lambda)$.  As $\lambda$ can be chosen arbitrarily small, then Lemma~\ref{lem-spanningbound}  implies that   $h_{top}(\Phi_t | \fM) =0$, which completes the proof of Proposition~\ref{prop-entropywhPhi}.

Fix $\lambda > 0$ and  $0 < \e < \delta_{\cK}$ then  set  $\e_6 =    \e/8$.
Then   there exists $n_0 > 0$ such that for $n \geq n_{0}$ we have $s(\cG_{\whPhi}^* | \fMR , n, \e_6) \leq \exp(n \, \lambda /2)$. 
Thus,  for each $n \geq n_{0}$  there exists a $(\cG_{\whPhi}^*, n, \e_6)$-separated subset $\whF_n \subset \fMR$ with maximal cardinality which satisfies $\# \whF_n \leq \exp (n \, \lambda/2)$.
 Then for each pair $x \ne y \in \whF_n$ there exists $-n \leq \ell_{x,y} \leq n$ such that $x,y \in Dom(\whPhi^{\ell_{x,y}})$ and 
  $d_{\bRt}(\whPhi^{\ell_{x,y}}(x),\whPhi^{\ell_{x,y}}(y))>\e_6$.

For each $n \geq n_{0}$, the set $\whF_n$ was chosen to have maximal cardinality. Thus, for each $y \in \fMR$,   there exists  $x \in \whF_n$   such that   $d_{\bRt}(\whPhi^{\ell}(x),\whPhi^{\ell}(y))\leq \e_6$ for all $-n \leq \ell \leq n$.
  That is,   $\whF_n$ is $(\whPhi | \fMR, n, \e_6)$-spanning. 
  
 We consider the restricted set $\whF_n \cap U$. The set $\whF_n \subset \fMR$  is $(\cG_{\whPhi}^*, n, \e_6)$-separated for $\fMR$, hence is $(\whPhi | \fMR, n, \e_6)$-spanning for $\fMR$, but this need not imply that  $\whF_n \cap U$ is $(\cG_{\whPhi}^*, n, \e_6)$-spanning for $\fMR \cap \oU$.
  If so, then set    $\whF_n^U = \whF_n \cap U$.  
Otherwise,  if $\whF_n \cap U$ is not $(\cG_{\whPhi}^*, n, \e_6)$-spanning for $\fMR \cap \oU$, then we  define $\whF_n^U \subset U$ 
  by adding to the set $\whF_n \cap U$ sufficiently many points, obtained from $\e_6$-perturbations of points in $\whF_n - (\whF_n \cap U)$,  so that $\whF_n^U$ is $(\cG_{\whPhi}^*, n, \e_6)$-spanning for $\fMR \cap \oU$.  This is done as follows.
A point  $x \in \whF_n - (\whF_n \cap U)$  is said to be   in the $(\cG_{\whPhi}^*, n, \e_6)$-shadow of $U$, if there exists $y_x \in U$ such that 
  $d_{\bRt}(\whPhi^{\ell}(x),\whPhi^{\ell}(y_x))\leq \e_6$ for all $-n \leq \ell \leq n$. 
Then define $\whF_n^U$ as follows:
  $$\whF_n^U = (\whF_n \cap U) \cup \{y_x \mid  x \in \whF_n -(\whF_n \cap U) , ~ x ~ {\rm  is ~ a ~ shadow ~ point ~ for}~  U\} . $$
Set $\e_7 = 2 \e_6 = \e/4$, so that the $\e_7$ ball centered at a point $x \in \bRt$ contains the $\e_6$-ball centered at $\zeta$, for any point $\zeta$ in the $\e_6$-ball centered at $x$. 
Then by construction, $\whF_n^U$  is  $(\whPhi | \fMR, n, \e_7)$-spanning for $\fMR \cap \oU$. 
Note that    $\# \whF_n^U \leq \# \whF_n$.

 \bigskip

The next step in the proof is to use the  flow $\Phi_t$ of the sets $\whF_n^U$, for $U \in \cD_{\fM}^{(n)}(\whPhi)$, to obtain a $(\Phi_t | \fM, T_n,\e)$-spanning set in $\fM$, for appropriate $T_n$. 
There is a technical difficulty that arises in this procedure, due to the fact that the return time for the flow to the space $\bRt$ is not constant. 
  We require the following observation,   that there is a ``modulus of continuity'' for the time-change along an orbit, with respect to  variations of the initial points for   flow segments of $\Phi_t$ with bounded lengths.

 \begin{lemma}\label{lem-boundedovershoot}
 There exists $M_{\cK} > 0$ and  $\delta_{\cK} > 0$ so that for all $0 < \e' < \delta_{\cK}$:

 Let $x \in \fMR$ and $|T| \geq 1$ be such that $\xi = \Phi_T(x) \in \fMR$.
 Let $n_x$ be such that $\xi = \whPhi^{n_x}(x)$. 
For  $y,z \in B_{\bRt}(x,\e') \cap \fMR$, assume that $\ds  d_{\bRt}(\whPhi^{\ell}(x), \whPhi^{\ell}(\eta)) < \e'$ for all $0 \leq \ell \leq n_x$, 
for  both $\eta = y$ and   $\eta = z$.
Let   $T'$ be such that $\Phi_{T'}(y) = \whPhi^{n_x}(y)$, and     $T''$ be such that $\Phi_{T''}(z) = \whPhi^{n_x}(z)$. 
Then   if $d_{\bRt}(y,z) < \e'/(M_{\cK} \,|T|)$, we have 
\begin{equation}\label{eq-uniformtimechange}
d_{\mK}(\Phi_t(y), \Phi_t(z)) ~ \leq ~ \e' \quad {\rm for ~ all} ~ |t| \leq  |T|  ~ .
\end{equation}
 \end{lemma}
 \proof
The curves $\{ \Phi_t(\xi) \mid 0 \leq t \leq T\}$ for $\xi \in \{x,y,z\}$ are    integral curves for the vector field $\cK$. The idea of the proof is that for appropriately chosen $\delta_{\cK}$ the $\cK$-orbits segments for $y$ and $z$   are contained in a suitably small neighborhood of the $\cK$-orbit segment for $x$, chosen so that the length of the vector $\cK$ along these curves are uniformly close. It then follows that the return times for $y$ and $z$ are also uniformly close. The details of the proof use  standard methods of differential equations, and are left to the reader.  
 \endproof

We use Lemma~\ref{lem-boundedovershoot} to construct a set $\whF_n'$ whose points are sufficiently close together, so that the $\Phi_t$ flow of $\whF_n'$ 
 yields a $(\Phi_t | \fM, T_n', \e)$-spanning set, for a sequence of values $T_n'$ which tend to infinity.

 Let $T_n = L_{\Phi} \cdot n$, and set   $\e_8 = \e_7/(2 M_{\cK} \,T_n) = \e_7/(2 n \, L_{\Phi} \, M_{\cK})$ for   $L_{\Phi}$  as defined in Lemma~\ref{lem-returntoR} and   $M_{\cK}$ as defined in Lemma~\ref{lem-boundedovershoot}. 
 Note that there  exists a constant $A_7 > 0$, which  depends only on the geometry of $\bRt$,  such that the number of   points required for an $\e_8$-spanning set of the disk $D_{\bRt}(x, \e_7) \subset \bRt$, \emph{ for the metric $d_{\bRt}$}, is bounded above by 
 $A_7 \, (\e_7/\e_8)^2 = A_7 \, (2 n \, L_{\Phi} \, M_{\cK})^2$.

 For each $U \in \cD_{\fM}^{(n)}(\whPhi)$ and  each $x \in \whF_n^U$,   choose   
    an $\e_8$-spanning set  for $\fMR \cap D_{\bRt}(x, \e_7) \cap \oU$, 
$$S(U,x,\e_8) = \{ x(U,x,i) \mid i \in \cI(U,x,\e_8) \} \subset \fMR \cap D_{\bRt}(x, \e_7) \cap U ~ ,$$
where the index set $\cI(U,x,\e_8)$ has cardinality $\# \, \cI(U,x,\e_8) \leq  A_7 \, (2 n \, L_{\Phi} \, M_{\cK})^2$.
 Then set
 \begin{equation}\label{eq-spanningset}
\whF_n' = \bigcup_{U \in \cD_{\fM}^{(n)}(\whPhi)} ~ \bigcup_{x \in \whF_n^U} ~ S(U,x,\e_8) ~ .
\end{equation}
Observe that $\# \, \whF_n' \leq C_{\Phi}(n) \cdot \exp (n \, \lambda/2) \cdot A_7 \, (2 n \, L_{\Phi} \, M_{\cK})^2$.

 We next construct a set $F_n$ which is contained in the $\Phi_t$-flow of the set $\whF_n'$, for a sequence of sufficiently small increments of the time parameter. 
Let $M_n$  be the greatest  integer with $M_n \leq T_n/\e_7$, and  for each $0 \leq \ell \leq M_n$   set   $t_{\ell} = \ell  \cdot \e_7$. Then define 
\begin{equation}\label{eq-spanningsetPhi}
F_n = \bigcup_{0 \leq \ell \leq M_n} ~ \Phi_{\ell}(\whF_n') ~ \subset ~ \fM ~ .
\end{equation}
   Observe that 
   $$\# F_n = M_n \cdot \# \whF_n' \leq  L_{\Phi} \cdot n/\e_7 \cdot C_{\Phi}(n) \cdot \exp (n \, \lambda/2) \cdot A_7 \, (2 n \, L_{\Phi} \, M_{\cK})^2 .$$
 Thus, for $n_0' \geq n_0$ sufficiently large, we have that $n \geq n_0'$ implies that $\# \, \whF_n \leq \exp (n \, \lambda)$.
Set $T_n' = T_n - \nu_{\cK}$. 

\begin{lemma}\label{lem-subexpspan}
  $F_n$   is $(\Phi_t | \fM, T_n', \e)$-spanning. 
\end{lemma}
\proof
 Let $\eta \in \fM$, then we must show that there exists $\xi \in F_n$ such that 
 $$d_{\mK}(\Phi_t(\eta) , \Phi_t(\xi)) \leq \e \quad {\rm for ~ all}~ -T_n \leq t \leq T_n ~ . $$
 
 First,   there exists $0 \leq t_{\eta} < \nu_{\cK}$ such that $y = \Phi_{-t_{\eta}}(\eta) \in \fMR$. Then there exists $U \in \cD_{\fM}^{(n)}(\whPhi)$ with $y \in \oU$, and $x \in \whF_n^U$ for which $\ds d_{\bRt}(\whPhi^{\ell}(x), \whPhi^{\ell}(y)) \leq \e_7$   for all $-n  \leq \ell \leq n$. 
By the choice of the set $N(U,x)$, there exists $z\in N(U,x)$ for which   $d_{\bRt}(y,z) \leq \e_8$. Note that $z \in \whF_n'$.

 We may assume that $T_n \geq  \nu_{\cK} + 1$, as we consider $T_n \to \infty$. Then observe that   the    hypotheses of  Lemma~\ref{lem-boundedovershoot} are satisfied for $\e' = \e_7$,  $T = T_n$,   $x \in \whF_n$ and the pair of points $\{y,z\}$ as above, and so  the estimate \eqref{eq-uniformtimechange} holds for these points.
 Thus we have:
 \begin{equation}\label{eq-approximation}
d_{\mK}(\Phi_t(y), \Phi_t(z)) ~ \leq ~ \e_7 =    \e/4\quad {\rm for ~ all} ~  |t| \leq  T_n ~ .
\end{equation}
Choose $0 \leq \ell \leq  M_n$   so that $|t_{\eta} - t_{\ell}| \leq  \e_7$. Then $d_{\mK}(\Phi_{t_{\eta}}(z) , \Phi_{t_{\ell}}(z)) \leq  \e_7$ 
 as the flow $\Phi_t$ has unit speed.

Then for $-T_n' \leq t \leq T_n'$, note that  $|t+t_{\eta}| \leq  |T_n - \nu_{\cK}| + |\nu_{\cK}| \leq T_n$, so  we have the estimates:
\begin{eqnarray*}
d_{\mK}(\Phi_{t}(\eta) , \Phi_t(\xi))  & = &  d_{\mK}(\Phi_{t} \circ \Phi_{t_{\eta}}(y) , \Phi_{t} \circ \Phi_{t_{\ell}}(z)) \\
& \leq &  d_{\mK}(\Phi_{t + t_{\eta}}(y) , \Phi_{t+t_{\eta}}(z)) +  d_{\mK}(\Phi_{t + t_{\eta}}(z) , \Phi_{t+t_{\ell}}(z)) \\
& \leq  &  \e_7 + \e_7 = \e/2 ~ .
\end{eqnarray*}
as was to be shown. 
\endproof

 Finally, to complete the proof of Proposition~\ref{prop-entropywhPhi}, note that  $T_n = L_{\Phi} \cdot n$ is a linear function of $n$, so the growth rate of the sets $F_n$ has an exponential bound  $\exp(n \, \lambda) = \exp(T_n \, \lambda/L_{\Phi})$. As the choice of    $\lambda > 0$ was arbitrary,   by Lemma~\ref{lem-spanningbound} we conclude that  $h_{top}(\Phi_t | \fM) =0$.
 \endproof

 \bigskip

There is   another entropy invariant associated to $\Phi_t$, which is defined in terms of the action  of   the   $\psg$  $\cGK^*$ on $\bRt$  introduced in  Definition~\ref{def-pseudogroup+}, and equipped with the symmetric generating set  
\begin{equation}\label{eq-symmetricgenset}
\cGK^{(1)} = \{Id, (\psi)^{\pm 1}, (\phi^+_1)^{\pm 1}, (\phi^-_1)^{\pm 1}, (\phi^+_2)^{\pm 1},  (\phi^-_2)^{\pm 1}\} .
\end{equation}
Let $h_{GLW}(\cGK^* | \fMR)$ denote the pseudogroup entropy  associated to the generating set $\cGK^{(1)}$ acting on $\fMR$.

Recall that   the map $\psi \in \cGK^{(1)}$    represents the return map for the Wilson flow, and corresponds to the ``short-cuts'' introduced in Lemma~\ref{lem-shortcut}.
 This additional generator raises the possibility that the growth rate of $(\cGK^*, n,\e)$-separated  points may be greater the growth rate of $(\cG_{\whPhi}^*, n,\e)$-separated  points. On the other hand, the $\psg$ $\cG_{\whPhi}^*$ is the ``full pseudogroup'' generated by the return map $\whPhi$, while the generators for $\cGK^*$ form just  a subset of these. We use Proposition~\ref{prop-syndeticgk*} to compare these growth rates, and prove the following:

\begin{prop}\label{prop-entropysection}
 $h_{GLW}(\cG_{\whPhi}^* | \fMR) > 0$ implies that $h_{GLW}(\cGK^* | \fMR) > 0$.
\end{prop}
\proof

Assume that $h_{GLW}(\cG_{\whPhi}^* | \fMR) =  \lambda > 0$. For $\e > 0$ sufficiently small, and    $n$ sufficiently large, 
let $E_0(n, \e) \subset \fMR$  be  a $(\cG_{\whPhi}^*, n,\e)$-separated set   with  cardinality 
$\# E_0(n, \e) \geq \exp(n \, \lambda/2)$.  
Then for each distinct pair $x,y \in E_{0}(n, \e)$, there exists $-n \leq \ell_{x,y} \leq n$ such that 
$\ds d_{\bRt}(\whPhi^{\ell_{x,y}}(x) , \whPhi^{\ell_{x,y}}(y) ) > \e$.

 By Corollary~\ref{cor-syndeticanyr}, for each $x \in \fMR$, the orbit  $\cGK^*(x) \subset \fMR$ is syndetic for the constant $\nu_{u}^*$.  Recall that $M_u > 0$ is the greatest integer with $M_u \leq 2\nu_u^*$ then 
by Corollary~\ref{cor-returntoe}, there exists $0 \leq  m_{x} \leq M_u$ such that  $\whPhi^{-m_{x}}(x) \in \cGK^*(x)$. 
Also by the proof of Proposition~\ref{prop-syndeticgk*}, there exists $0 \leq \ell_x \leq M_u/2$ such that    $r(\whPhi^{\ell_x + \ell_{x,y}}(x)) < r_e$.  Thus, module changing $\ell_x$ for $\ell_x+1$, we can further assume that $|z(\whPhi^{\ell_x + \ell_{x,y}}(x))|> \delta$   for some $\delta > 0$ independent of the choice of $x$ and  $0 \leq \ell_x \leq M_u$.

For each $x \in  E_0(n, \e)$ we have  associated two indices, $0 \leq \ell_x \leq M_u$ and $0 \leq m_{x} \leq M_u$. 
The total number of pairs $\{\ell_x, m_x\}$ is at most $M_u^2$, so there exists values $0 \leq \ell' \leq M_u$ and $0  \leq m' \leq M_u$
such that the set
$$E_0(n, \e)(\ell', m') = \{x \in E_0(n, \e) \mid \ell_x = \ell' ~ {\rm and}~ m_x = m'\}$$
has cardinality at least $\# E_0(n, \e)/M_u^2$. Fix such $\ell', m'$ and set $E_0'(n, \e) = E_0(n, \e)(\ell', m')$.
 Set  $\cE_0'(n, \e) = \whPhi^{-m'}(E_0'(n, \e) ) \subset \fMR$, which is then   $(\cG_{\whPhi}^*, m' + n,\e)$-separated.

  Proposition~\ref{prop-sigmadyn} shows that for each $\xi  \in \fMR$ and $\eta = \whPhi(\xi)$, the $\cK$-orbit segment  $[\xi , \eta]_{\cK}$ has length bounded above by  $\nu_{\cK}$. Then there is a uniform norm on the derivative $D_{\xi}\Phi_t \colon T_{\xi}\mK \to  T_{\eta}\mK$, so there  is a uniform upper bound $M_{\Phi} \geq 1$ on the norm of the derivative matrix of $D_{\xi}\whPhi$ and on its inverse $(D_{\xi}\whPhi)^{-1}$. Thus, for   $\xi, \xi' \in \fMR \cap Dom(\whPhi^{\ell'})$, we have 
\begin{equation}\label{eq-normsbounded}
(M_{\Phi})^{-\ell'} \cdot d_{\bRt} (\xi , \xi') ~ \leq ~ d_{\bRt}(\whPhi^{\ell'}(\xi) , \whPhi^{\ell'}(\xi') ) ~ \leq ~ (M_{\Phi})^{\ell'} \cdot d_{\bRt} (\xi , \xi') .
\end{equation}

For $x \ne y \in E_0'(n, \e) $, set $\xi = \whPhi^{-m'}(x) \in \cE_0'(n, \e)$ and  $\xi' =  \whPhi^{-m'}(y)\in \cE_0'(n, \e)$. 

Let  $\ell_{x,y}' = \ell_{x,y} + \ell' + m' \leq n + 2M_u$, then  by the choice of $\ell'$ and $m'$,  Corollary~\ref{cor-syndeticanyr}  implies  the map $\ds \whPhi^{\ell_{x,y}'} \in \cGK^{(n + 2M_u)}$. Also, we  have the estimate
\begin{equation}\label{eq-normsbounded2}
d_{\bRt}(\whPhi^{\ell_{x,y}' }(\xi) , \whPhi^{\ell_{x,y}' }(\xi')) = d_{\bRt} (\whPhi^{\ell' + \ell_{x,y}}(x) , \whPhi^{\ell' + \ell_{x,y}}(y)) \geq M_{\Phi}^{-m'} \cdot d_{\bRt} (\whPhi^{\ell_{x,y}}(x) , \whPhi^{\ell_{x,y}}(y)) > M_{\Phi}^{-M_u} \cdot \e ~.
\end{equation}
Set $\e' =  M_{\Phi}^{-M_u} \cdot \e$.  Then    $\cE_0'(n, \e)$ is  $(2 (n + 2M_u) \cdot  \nu_{\cK}, \e')$-separated   for the action of $\cGK^*$ on $\fMR$.  

These estimates hold for a sequence of integers $n \to \infty$, so it follows that $h_{GLW}(\cGK^* | \fMR) \geq \lambda/2 > 0$, as was to be shown. 
\endproof

 \bigskip
 
The rest of this section is devoted to the study of the entropy invariant  $h_{GLW}(\cGK^* | \fMR)$. One approach to the study of the restricted entropy $h_{top}(\Phi_t | \fM)$ is to calculate the invariant measures for the flow, then estimate the Lyapunov spectrum for the flow  on $\fM$ and apply the Margulis-Ruelle Inequality  \cite{Mane1983,Ruelle1978} to estimate it from above. To show  the entropy $h_{top}(\Phi_t | \fM) = 0$ it then suffices to  show that the Lyapunov spectrum must be trivial for any invariant measure, which is done by analyzing the asymptotic  properties of the derivative along the flow. In essence,  the idea behind the proof of the following result is analogous, in that we  obtain uniform estimates for the expansive properties of the dynamics    of the  $\psg$  $\cGK^*$.

\begin{thm}\label{thm-entropyvanish}
  $h_{GLW}(\cGK^* | \fMR) = 0$.
\end{thm}
\proof

Assume that  $h_{GLW}(\cGK^* | \fMR) = \lambda > 0$.
Let $\e > 0$ and $n_{\e}> 0$ be such that 
\begin{equation}\label{eq-entestimate1}
s(\cGK,\fMR, n,\e) > \exp(n \, \lambda/2) \quad {\rm for ~ all}~ n \geq n_{\e} .
\end{equation}
Let  $E_{n} \subset \fMR$ be a sequence of  $(\cG_{\whPhi}^*, n,\e)$-separated   sets such that $\# E_n \geq \exp(n \,\lambda/2)$.  
Actually, it may be necessary to pass to a subsequence $\{n_i \mid i=1, 2, \ldots \}$ to obtain the estimate \eqref{eq-entestimate1}, but for simplicity of notation, without loss of generality we assume there exists sets $E_n$ for all $n \geq n_{\e}$ with this property.
Then 
for each pair $\xi \ne \eta \in E_n$, either $d_{\bRt}(\xi, \eta) \geq  \e$, or else there exists   a word $\vp \in \cGK^{(n)}$ with   $\xi, \eta \in Dom(\vp)$,   and $d_{\bRt}(\vp(\xi) , \vp(\eta)) \geq \e$.  
We show this leads to a contradiction.

Recall that the metric $d_{\bRt}$ on  $\bRt$ is the usual euclidean flat metric, where its sides have lengths   $2 \times 4$. 
Choose $m$ so that 
 $$\frac{1}{\sqrt{18}} \cdot \exp(n \, \lambda/4) < m < \frac{1}{4} \cdot \exp(n \, \lambda/4) ,$$
and divide $\bRt$  into $2m \times 4m$ uniform squares with sides of length   $\delta = 1/m$ and diameter $\sqrt{2}/m^2$.
The number of squares is thus bounded above by $\frac{1}{2} \cdot \exp(n \, \lambda/2)$, 
so by   the  {Pigeonhole principle},  for some square region must contain two distinct points of $E_n$. 
 Hence,   for each $n$ there exists     $x_n, y_n \in E_n$ with  $d_{\bRt}(x_n, y_n) <  \frac{1}{3} \cdot \exp(-n \, \lambda/4)$, and 
   $\vp_n \in \cM(n)$ with $d_{\bRt}(\vp_n(x_n), \vp_n(y_n)) > \e$. 
   
 For $\xi \in Dom(\vp_n)$,  let $D_{\xi}(\vp_n)$ denote the $2 \times 2$ Jacobian matrix of first derivatives of $\vp_n$ at the point $\xi$, and let $\| D_{\xi}(\vp_n)\|$ denote its matrix norm. 
   Then  by the Mean Value Theorem, 
   there exist $w_n \in   Dom(\vp_n)$ in this same square region containing $x_n, y_n$ such that $\| D_{w_n}(\vp_n) \| \geq \e/3 \cdot \exp(n \, \lambda/4)$.   
It follows that   the norms $\| D_{w_n}(\vp_{n}) \|$ grow exponentially as a function of $n$.

 We next show that the uniform norms  $\| D_{\xi}(\vp)\|$ of  the derivative  matrices of $\vp \in  \cGK^{(n)}$ admit arbitrarily small exponential bounds as functions of $n$, which yields a contradiction. 
 We remark that this calculation implies   the   Lyapunov spectrum of the flow $\Phi_t$ for any invariant measure supported on $\fM$ is trivial.
  
For an invertible matrix $A$,   introduce the ``symmetric   norm''   $|\| A \|| = \max \{\|A\|, \| A^{-1}\|\}$, where $\| A \|$ is the usual sup-norm on the linear transformation defined by $A$.

Let $D_{\xi}(\phi_i^+)$ denote the $2 \times 2$ Jacobian matrix of first derivatives of $\phi^+_i$ at the point $\xi \in Dom(\phi^+_i)$. 
Define 
 \begin{equation}\label{eq-entestimate2}
C(\phi) = \sup \left\{  |\| D_{\xi}(\phi^+_i)\| |  \mid i = 1,2 ~ \& ~ \xi \in Dom(\phi^+_i) \right\} .
\end{equation}

For   $b \geq 1$,   introduce the  upper bound
 \begin{equation}\label{eq-entestimate3}
C(\psi, b) ~ = ~   \sup ~\left\{  |\| D_{\xi}(\psi^{\ell})\||   \mid 0 \leq \ell \leq b \quad {\rm and} \quad  \xi \in Dom(\psi^{\ell})  \right\} .
\end{equation}

 As the matrix $D_{\omega_i}(\psi)$ is the identity at the fixed points $\omega_i$ for $i=1,2$,  for $\xi \in \bRt$ with $r(\xi) = 2$, we have that  $\ds \lim_{\ell \to \pm \infty} ~ \frac{\ln \| D_{\xi}(\psi^{\ell})\|}{|\ell |} =0$.   Thus, for every $\mu > 1$, there exists $n(\psi, \mu) > 0$ such that 
\begin{equation}\label{eq-entestimate4}
1 \leq |\| D_{\xi}(\psi^{\ell})\|| \leq \mu^{\ell} \quad {\rm for} ~ \ell \geq n(\psi, \mu), ~ \xi \in Dom(\psi^{\ell}).
\end{equation}
 
 Recall from Section~\ref{sec-normal} the notion of  monotone words $\cM(n)$ of length at most $n$ in Definition~\ref{def-monotone}, and that by Proposition~\ref{prop-standardform}  every $\vp \in \cGK^{(n)}$ admits   a factorization $\vp = \vp^+ \circ \vp^-$,  where $\vp^+ \in \cM(n')$   and $(\vp^-)^{-1} \in \cM(n'')$ for integers $n', n''$ with  $n' + n'' \leq n$. 
For $\vp \in \cM(n)$  we   use this factorization and  the estimates   \eqref{eq-entestimate2}, \eqref{eq-entestimate3}, and \eqref{eq-entestimate4} to obtain a uniform estimate for the norm  $|\| D(\vp)\| |$  as a function of $n$.  
    
For $\mu = \exp(\lambda/32)$, choose $\ds b \geq   \max \{ n(\psi, \mu) , (32/\lambda) \cdot \ln ( C(\phi)) \}$, so  that  $C(\phi)^{1/b}  <  \exp(\lambda/32)$.
Also,  for   $\ell \geq b \geq n(\psi, \mu)$,  we have $\ds 1 \leq |\| D_{\xi}(\psi^{\ell})\|| \leq \exp(\ell \, \lambda/32)$  for all $\xi \in Dom(\psi^{\ell})$.

Given  $\vp \in \cM(n)$ recall  the product representation in \eqref{eq-productformula}. For the $b$ chosen above, there exists $i(\vp, b) \geq 1$ such that $\ell_i \geq b$ for all $1 \leq i < i(\vp, b)$, and $\ell_i \leq b$ for $i = i(\vp, b)$.  Write $\vp = \vp^{(b)} \cdot \vp_{(b)}$ where $\vp^{(b)}$ starts with the map $\ds \phi^+_{i(\vp, b)}$ and $\vp_{(b)}$ starts with  $\psi^{\ell_0}$.
Assume that  $n \geq \max \{n_{\e} , n(\psi, \mu)\}$,  where  $n_{\e}$ was introduced in the paragraph after \eqref{eq-entestimate1}.

The factor $\vp^{(b)}$ contains at most $N_b$ generators of the form $\phi^+_i$, and between each such map is a term $\psi^{\ell_i}$, so there are at most $N_b+1$ such factors with exponents either $\ell_i \leq b$ or $\ell_i > b$. Then
\begin{equation}\label{eq-entestimate5}
|\|D(\vp^{(b)})\|| ~ \leq ~  C(\phi)^{M_b} \cdot C(\psi, b)^{(N_b+1)} \cdot   \exp( n \, \lambda/32)   ~ .
\end{equation}
The factor $\vp_{(b)}$ contains at most $n/b$ generators of the form $\phi^+_i$, and between each such generator is a term $\psi^{\ell_i}$ where $\ell_i > b$ by the definition of $\vp_{(b)}$. Then
\begin{equation}\label{eq-entestimate6}
|\| D(\vp_{(b)})\|| ~ \leq ~  C(\phi)^{n/b} \cdot \exp( (n \, \lambda/32) ~ \leq ~  \exp(n \, \lambda/32) \cdot \exp( (n \, \lambda/32) ~ =  ~\exp(n \, \lambda/16) ~ .
\end{equation}
 Combining these estimates,  we have:
 \begin{lemma}\label{lem-upperbound} 
For $\lambda > 0$,  let $\ds b \geq   \max \{ n(\psi, \mu) , (32/\lambda) \cdot \ln ( C(\phi)) \}$. For all $n \geq \max \{n_{\e} , n(\psi, \mu)\}$, we have
\begin{eqnarray}
\max\{ \| D(\vp)\| , \| D(\vp)^{-1}\|\} ~ & \leq & ~  |\|D(\vp^{(b)})\|| \cdot |\|D(\vp_{(b)})\|| \label{eq-entestimate7} \\
 & \leq &  ~ \exp( n \, \lambda/32) \cdot \exp(n \, \lambda/16)  < \exp(n \, \lambda/8) ~ . \nonumber
\end{eqnarray}
\end{lemma}

Now recall that with the assumption that $h_{GLW}(\cGK | \fMR) = \lambda > 0$, then  for each $n \geq n_{\e}$, there exists  $\vp_{n} \in \cM(n)$ and $w_n \in Dom(\vp_{n})$   such that $\| D_{w_n}(\vp_{n}) \| \geq \e \, \exp(n \, \lambda/4)$. Proposition~\ref{prop-standardform} implies that each such word $\vp_{n}$ admits a factorization into a monotone decreasing factor and an increasing factor,  each of length at most $n$. Apply the estimate \eqref{eq-entestimate7} to both factors to obtain that for $n$ sufficiently large,  
$$\| D_{w_n}(\vp_{n}) \|  < \exp(n \, \lambda/8) \cdot \exp(n \, \lambda/8) = \exp(n \, \lambda/4) ~.$$
For $n$ sufficiently large this   yields a contradiction, which  
 completes the proof of Theorem~\ref{thm-entropyvanish}.
  \endproof
 
  \bigskip
 
\section{Lamination entropy}\label{sec-entropylamination}

The minimal set for a flow is an invariant of topological conjugacy, so represents an ``invariant'' for the flow. Properties of the minimal set can be used to study and classify the flow as a dynamical system. In the case of a  generic Kuperberg flow   $\Phi_t$ the space $\fM$ is the unique minimal set, and is   fundamental  for the study of the dynamics of the flow as has been seen in previous sections.  In this section, we  use the zippered lamination structure of $\fM$ to define      invariants for the lamination itself, derived from the holonomy pseudogroup for $\fM$.  The invariants studied in this section are ``entropy-like''. In the next section, we study growth-type invariants for the leaves, and both types of invariants reveal the beautiful subtleties of the class of generic Kuperberg flows. 

We first give precise definitions of the maps $\{\opsi, \ophi_1 , \ophi_2\}$ introduced in Section~\ref{sec-normal}, and of the $\psg$ $\cGM^*$ they generate. We then study how this $\psg$ is associated to the holonomy of the zippered lamination, and compare it with the $\psg$ $\cGK^*$ introduced in Section~\ref{sec-pseudogroup}. These two $\psg$s are naturally closely related, as the action of $\cGK^*$ on $\bRt$ studied in previous sections induces an action on the families of curves in the intersection $\fM \cap \bRt$  and thus induce maps in $\cGM^*$ acting on the transverse Cantor set $\fC$ which parametrizes these curves. 

 A key difference is that $\cGM^*$   contains maps defined by the holonomy along all paths in the leaves of $\fM$, and not just those paths following a $\cK$-orbit. In particular,  the short-cut maps, which play a role in the study of the entropy for the $\cGK^*$ action, arise naturally from the geometry of the leaves of $\fM$.  
The dynamics of the $\cGM^*$ action on the transverse Cantor set $\fC$ defined in \eqref{eq-cantorset} provides an alternate source of invariants for the dynamics of the flow $\Phi_t$ on $\fM$.  

The first invariant studied is the   entropy $h_{GLW}(\cGM^*)$ of the action of $\cGK^*$ on the transverse Cantor set $\fC \subset \cT$.  We relate  this entropy to   the  entropy for the action  of the $\psg$ $\cGK^*$ on $\fMR$ as studied in the previous section, and   
    prove that $h_{GLW}(\cGM^*) =0$ as is naturally expected. 

The main result of this section is that there is a non-zero entropy invariant for the $\psg$ $\cGM^*$ which is obtained by considering   ``lamination slow entropy'' $h_{GLW}^{\alpha}(\cGM^*)$ for $0 < \alpha \leq 1$. We show in Theorem~\ref{thm-slowentropy} 
  that if the insertion maps $\sigma_i$ used in the construction of $\mK$ have ``slow growth'', as defined in Definition~\ref{def-slowgrowth}, then   $h_{GLW}^{1/2}(\cGM^*) > 0$ for growth constant $\alpha = 1/2$.  Thus, even though the flow entropy $h_{top}(\Phi_t | \fM) = 0$, there is enough ``chaos'' in the orbits of $\Phi_t$ and hence in the holonomy action of $\cGM^*$, to imply that it has  positive slow entropy. We also note, at the end of the section, how one can obtain $h_{GLW}(\cGM^*) > 0$ for the usual entropy by modifying properties of the insertions maps $\sigma_i$, which gives a new insight into the results of the Kuperbergs in \cite[Section~8]{Kuperbergs1996}.

   Recall  that $\cT=\{z=0\} \cap  \bR_0$ is the line segment
 in  $\bRt$ transverse to the interiors of most of the leaves of $\fM$, 
and $\fC\subset \cT$ as defined in \eqref{eq-cantorset} is a Cantor set by Proposition~\ref{prop-cantor}. 
Let $\fMC \subset \fMR = \fM \cap \bRt$ be the union of the path components of $\fMR$ that contain a point in
$\fC$. 

 The holonomy pseudogroup $\cGM$  for  the zippered lamination $\fM$ can be defined using a covering of $\fM$ by foliation charts. The proof of Theorem~\ref{thm-zippered} introduced the sections $\bT_i$ and their associated Cantor sets $\fC_i$  which form the model spaces for the foliation charts constructed there, and so   the holonomy defined using this  covering produces a pseudogroup. 
 On the other hand,   Corollary~\ref{cor-min=min} implies that each leaf in $\fM$ intersects the section $\bRt$ so we can alternately define $\cGM$ using the induced holonomy maps on the   space $\fMC$ associated to the transverse Cantor set $\fC$.
 This reduction yields a pseudogroup $\cGM$ which is most closely related to our previous constructions, so we assume that $\cGM$ is defined using   the space $\fC$. 
 Also, recall that $\fM_0$ is dense in $\fM$ and the intersection $\fM_0 \cap \bRt$ is dense in $\fMR$ so the holonomy maps of $\cGM$ can be defined by paths in the leaf $\fM_0$ with endpoints in the transversal $\fC$.

Recall the    symmetric generating set   for the $\psg$ $\cGK^*$  
$$\cGK^{(1)} = \{Id, (\psi)^{\pm 1}, (\phi^+_1)^{\pm 1}, (\phi^-_1)^{\pm 1}, (\phi^+_2)^{\pm 1},  (\phi^-_2)^{\pm 1}\}$$
which act on the endpoints of the arcs in $\fM_\fC$ and so induce actions of the tree $\TP$ as discussed in Section~\ref{sec-normal}. The vertices of the tree $\TP$ are points in $\fC_0$ so we get an induced action on $\fC_0$.  We   formalize the definitions of the local homeomorphisms defined  on subsets of $\fMC$ induced by these maps. 
We first note  that $\fMC$ ``fibers'' over $\fC$ in the following sense:
 \begin{lemma}\label{lem-structureM}
 There exists a continuous map $\pi_{\fM} \colon \fMC \to \fC$,  such that for
 $\xi \in \fC$, 
 \begin{enumerate}
\item $\pi_{\fM}(\xi) = \xi$; 
\item  $\pi_{\fM}(p_0^k(\xi)) = \xi$ for $k=1,2$,
\end{enumerate}
where $p_0^k$ for $k=1,2$ are the maps introduced in Section~\ref{sec-zippered}. 
The fibers  of the map  $\pi_{\fM} \colon \fMC \to
\fC$ are compact intervals  $\pi_{\fM}^{-1}(\xi) \subset \bRt$ whose
lengths are bounded above.
 \end{lemma}
 \proof
 For each $\xi \in \fC$ there exists a maximal connected closed arc
 in $\fMC$ intersecting $\cT$ at $\xi$. 
 Let $p_0^1(\xi)$ be the lower endpoint of this arc contained in $\{z
 \leq 0\}$, and $p_0^2(\xi)$   the upper endpoint of this arc contained in $\{z \geq 0\}$. 
 For $\eta \in \fMC$ there is a unique   arc with $\eta \in [p_0^1(\xi), p_0^2(\xi)]$, and we set $\pi_{\fM}(\eta) = \xi$. Then properties (1) and (2) follow by definition.
 
 The continuity of the map $\pi_{\fM}$ follows from the proof of Theorem~\ref{thm-zippered}. The bound on the lengths of the fibers is a  consequence of
Theorem~\ref{thm-boundedlength}.
  \endproof

For points in the dense subset  $\fC_0 \subset \fC$,  the   definition of $\pi_{\fM}$  is based on the intersections with $\cT$ of the  families of   $\g_0$ and $\lambda_0$ curves in $\bRt$.  That is,  each $\xi \in \fC_0$  has the form
\begin{equation}\label{eq-xi-labels}
\xi = \gamma_0(i_1,\ell_1 ;\cdots ; i_{n-1},\ell_{n-1} ;\ell_n) \cap \cT  \quad {\rm or}  \quad 
\xi = \lambda_0(i_1,\ell_1 ;\cdots ; i_{n-1},\ell_{n-1} ;\ell_n) \cap \cT .
\end{equation}
Each such curve   is   
defined by the labeling of its endpoints denoted   by $p_0^1(\xi)$ for the lower endpoint and
$p_0^2(\xi)$ for the upper endpoint. Then for $\eta \in\fM_0\cap\bRt$
in a $\g_0$ or $\lambda_0$ curve which is an arc, the application $\pi_{\fM}$   takes the
point $\eta$  to the intersection with $\cT$ of the $\g_0$ or $\lambda_0$ curve which contains it.

 The map $\phi_k^+$ for $k=1,2$ induces maps $\ophi_k$   with domains defined by
\begin{equation}\label{eq-domainophi}
Dom(\ophi_k) ~ = ~ \left\{  \xi \in \fC  \mid  p_0^k(\xi) \in Dom(\phi_k^+) \right\} ~.
\end{equation}
The maps are then formally defined by:
 \begin{eqnarray}
\ophi_1(\xi) ~ & = &~ \pi_{\fM}(\phi_1^+(p_0^1(\xi))) ~ {\rm for} ~  \xi \in Dom(\ophi_1) \nonumber \\
\ophi_2(\xi) ~ & = &~   \pi_{\fM}(\phi_2^+(p_0^2(\xi))) ~ {\rm for} ~  \xi \in Dom(\ophi_2)  ~ . \nonumber
\end{eqnarray}
The map $\ophi_k$ has   inverse $\ophi_k^{-1}$ defined on the image of $\ophi_k$. 
 The action of the maps $\phi_k^+$  on the $\gamma_0$ curves in $\fM_0$  are  illustrated in Figure~\ref{fig:generatorsM0R0}, and the action of the induced maps $\ophi_k$  
on the vertices of the tree $\TP$ are illustrated in Figure~\ref{fig:fM0cTa}, as discussed in Section~\ref{sec-normal}. 
 Note that each map    $\ophi_k$ increases the  level of the curve in $\fMC$ defining a point $\xi \in \fC$ by 1,  for $k=1,2$.

      Next,    we define the map $\opsi \in \cGM$ induced by the map $\psi$.
Recall from the discussion in Section~\ref{sec-pseudogroup} that   the Wilson flow  reverses direction at the annulus $\cA = \{z=0\} \subset \mW$, and is anti-symmetric with respect to $\cA$. 
Define the domain of $\opsi$ by 
\begin{equation}\label{eq-domainopsi}
Dom(\opsi) ~ = ~ \left\{  \xi \in \fC  \mid  z(p_0^1(\xi)) < 0 ~ \& ~ z(\psi_-(p_0^1(\xi)))   leq    0  \right\} ,
\end{equation}
for $\psi_-$ as defined in \eqref{eq-defpsiextended}.
The map $\opsi$ is then formally defined by:
 \begin{eqnarray}
\opsi(\xi)  ~ & = &~ \pi_{\fM}(\psi_-(p_0^1(\xi))) ~ {\rm for} ~  \xi \in Dom(\opsi) ~ .\nonumber 
\end{eqnarray}
 The map $\opsi$ has   inverse $\opsi^{-1}$ defined on the image of $\opsi$. 
 The action of  $\psi$ and its inverse on the $\gamma_0$ curves in $\fM_0$,  is  illustrated in Figure~\ref{fig:generatorsM0R0} as the vertical maps. The action of the induced map $\opsi$  on the vertices of the tree $\TP$ is given by translation along the center lines in $\TP''$, as discussed in Section~\ref{sec-normal}. 
 Note that the map    $\opsi$ preserves the  level.

   \begin{defn}\label{def-cGM}
   Let  $\cGM$ be the   pseudogroup generated by the collection of maps
\begin{equation}\label{eq-generatorsGM}
\cGM^{(1)} \equiv    \{\overline{Id},  \opsi , \opsi^{-1} , \ophi_1 , \ophi_1^{-1}, \ophi_2 , \ophi_2^{-1}\} ~ .
\end{equation}
Let  $\cGM^*$ be the $\psg$ formed  by  the compositions of maps in  $\cGM^{(1)}$
 and the  restrictions of these compositions to open subsets of their domains in $\fC$. 
  \end{defn}

  \begin{remark}
  We do not need to consider the maps $\phi_i^-$, since for $\xi\in \fC$ such that $p_0^i(\xi)\in Dom((\phi_i^-)^{-1})$ we have that $p_0^i(\xi)\in Dom(\phi_i^+\circ\psi^{-1})$ and
  $$\pi_\fM((\phi_i^-)^{-1}(\xi))=\pi_\fM(\phi_i^+\circ\psi^{-1}(p_0^i(\xi)).$$
  \end{remark}

The action of $\cGM$ on $\fC$ has a geometric model. Identify the points of $\fC_0$ with the vertices of $\TP$. Then the discussion in 
Section~\ref{sec-normal} gives the   actions of the generators in \eqref{eq-generatorsGM} on the   tree $\TP$. This remark is the basis for the proof of the following:

 \begin{prop}\label{prop-pseudogroupGM}
  The holonomy pseudogroup $\cGM$  of $\fM$ is the pseudogroup acting on $\fC$ defined by the holonomy maps   associated to  leafwise paths in $\fM$.
 \end{prop}
 \proof
The leaf $\fM_0$ is dense in $\fM$, so it suffices to consider the holonomy   associated to a path $\sigma \colon [0,1] \to \fM_0$ with $\sigma(0), \sigma(1) \in \fC_0$.
 As each $\xi \in \fC_0$ corresponds to a vertex of  $\TP$, we can assume that the path $\sigma \colon [0,1] \to \TP$. 
 Recall from Section~\ref{sec-normal} that $\TP'' \subset \TP' =   \cA \cap \fM_0$ consists of the line segments in $\fM \cap \cA$   which are formed by the intersection of $\cA$   with the simple propellers in the construction of $\fM_0$. 
 Then  $\TP$ is formed by
adding continuous curve segments in $\fM_0$ joining the line segments in $\TP''$, as 
 is  illustrated in  Figure~\ref{fig:fM0cT}.   
Thus, the edges of the tree $\TP$ consists of two types, those which belong to $\TP''$ and the added connecting paths. 

If the path $\sigma$ traverses  a connecting segment, then the holonomy map this induces  corresponds to the action of one of the maps $\ds \{\ophi_1 , \ophi_1^{-1}, \ophi_2 , \ophi_2^{-1}\}$.  If the path $\sigma$ traverses a segment in $\TP''$, then the holonomy map this induces  corresponds to the action of one of the maps $\ds \{\opsi , \opsi^{-1} \}$. Thus, the holonomy maps defined by such $\sigma$ are contained in $\cGM$.

Conversely, given a word in $\cGM$ we can associate a path $\sigma \colon [0,1] \to \TP$ with $\sigma(0), \sigma(1) \in \fC_0$ using the same correspondence of the generators in $\cGM^{(1)}$ with segments in $\TP$, so that every word in $\cGM$ is associated with the holonomy along a leafwise path. 
 \endproof

  As   noted in Section~\ref{sec-normal}, $\TP$ is a tree except for the loop based at the basepoint vertex $\omega_0$. Thus, a 
   path $\sigma \colon [0,1] \to \TP$ is   homotopic with endpoints fixed, to a   path which is monotone in level. 
    We derive from this remark  the existence of normal forms for words in $\cGM^*$ exactly as in Proposition~\ref{prop-standardform}.

 The word length on $\cGM^*$ is defined as before, where $\| \ovp \| \leq n$ if $\ovp$ can be expressed as a composition of at most $n$ maps in $\cGM^{(1)}$.
 Also, define monotone words in $\cGM^*$ in analogy with Definition~\ref{def-monotone}:
\begin{defn}\label{defn-monotoneM}
A word   $\ovp \in \cGM^*$ is said to be \emph{monotone increasing} if  it is written in the form
\begin{equation}\label{eq-productformulaGM}
\ovp=\opsi^{\ell_m} \circ \ophi_{j_m} \circ \opsi^{\ell_{m-1}} \circ \cdots
\circ   \ophi_{j_2} \circ      \opsi^{\ell_1} \circ \ophi_{j_1}
\end{equation}
where each $j_k=1,2$ and $\ell_k\geq 0$. 
Set $\cM_\fM(0) = \{\overline{Id}\}$, and define 
\begin{equation}\label{eq-monotonedefGM}
\cM_\fM(n) ~ = ~  \{\ovp \in \cGM^* \mid \ovp ~{\rm monotone} ~  \&~ \|\ovp \| \leq n\} \quad; \quad \cM_\fM(\infty) ~ = ~ \bigcup_{n \geq 0} ~ \cM_\fM(n) ~ .
\end{equation}
\end{defn}

The proof of the upper bound estimate on the function $\# \cM_\fM(n)$ in Proposition~\ref{prop-asymptotic} applies verbatim to the function $\ds \# \cM_\fM(n)$, and  we have:
\begin{prop}\label{prop-asymptoticGM}
  For each $b \geq 1$, there is a polynomial function $P_b(n)$  of $n$ such that 
the  cardinality of the set $\cM_\fM(n)$   satisfies
  \begin{equation}\label{eq-productestimateGM}
\#\cM_\fM(n) \leq P_b(n) \cdot  2^{(n/b)} .
\end{equation}
\end{prop}

The   following result is analogous to    Corollary~\ref{cor-subexponentialgrowth}, and the proof follows in the same way.          
\begin{cor}\label{cor-asymptoticM} 
The  cardinality of the set $\cM_\fM(n)$ of monotone words of length
at most $n$ in $\cGM$ satisfies
 \begin{equation}\label{eq-productformula4a}
\lim_{n\to \infty} ~ \frac{\ln (\# \cM_\fM(n))}{n} = 0 .
\end{equation} 
\end{cor}

We then have the following result, which is analogous to Proposition~\ref{prop-standardform}, and follows from  Proposition~\ref{prop-pseudogroupGM} and the comments afterwards.
\begin{prop}\label{prop-standardformGM}
Let $\ovp \in \cGM^*$ with $\|\ovp \| \leq n$. Then there exists a factorization $\ovp = \ovp^+ \circ \ovp^-$,  where $\ovp^+ \in \cM_\fM(n')$   and $(\ovp^-)^{-1} \in \cM_\fM(n'')$ for integers $n', n''$ with  $n' + n'' \leq n$. Moreover, we have     $Dom(\ovp) \subset Dom(\ovp^+ \circ \ovp^-)$. 
\end{prop}

The factorization $\ovp = \ovp^+ \circ \ovp^-$ is said to be the \emph{normal form} for the word $\ovp$.
 
 We conclude this discussion of the structure of $\cGM^*$ with some remarks concerning the relation between the actions of $\cGM^*$ and $\cGK^*$.
  First, recall that the domain $D(\whPsi)_-^-$ of the the induced action of the Wilson map $\psi_-$ contains the segment $J_0  \subset \bRt$ as defined in \eqref{eq-intervals}, and $\psi_-$ defines a strict contraction of $J_0$ to the fixed point $\omega_1$ which is the upper endpoint of $J_0$. For the $\gamma_0$ and $\lambda_0$ curves in $\bRt$ with lower endpoint in $J_0$, the action of $\psi_-$ on these curves induces the map $\opsi$ which is a strict contraction to $\omega_0 \in \fC$. Thus, the holonomy pseudogroup $\cGM$ contains a contracting fixed point, 
  which is just the holonomy along the loop in $\TP$ containing the point $\omega_0$.

The second remark concerns the action of $\psi$ on the endpoints of curves in $\fMC$ and the action of $\opsi$ on $\TP$. Let $\xi \in \fC$ be defined by a $\gamma_0$ or $\lambda_0$ curve in $\fMC$ with endpoints $x = p_0^1(\xi)$ and $\ox = p_0^2(\xi)$ where $r(x) = r(\ox) = r_0 > 2$. Let $k_x \geq 0$ be the exponent such that $\ox = \psi^{k_x}(x)$.  Then $x_{\ell} = \psi^{\ell}(x)$ for $0 \leq \ell \leq k_x$ is a sequence of points with $r(x_{\ell}) = r_0$ and $z(x_{\ell +1}) > z(x_{\ell})$ for $0 \leq \ell < k_x$. Let $m_x$ be the greatest integer with $m_x \leq k_x/2$. By the anti-symmetry of the Wilson flow, we have $z(x_{m_x}) < 0$ and $z(x_{m_x + 1}) \geq 0$. That is, the action of $\psi$ on the curves with lower endpoints    
$\{x_{\ell} \mid 0 \leq \ell \leq m_x\}$ produces a sequence of curves whose intersections with $\cT$ have decreasing radius. For $k_x$ even,  the closest curve to $I_0$ in the sequence, the innermost curve, corresponds to the endpoint  $x_{m_x}$. When $k_x$ is odd, the closest curve degenerates to the single point $x_{m_x+1}$ with $z(x_{m_x + 1}) =0$. The action of $\psi$ on the endpoints 
$\{x_{\ell} \mid m_x < \ell < k_x\}$ then reverses this process, where the endpoints $x_0$ and $x_{k_x}$ correspond to the same curve in $\fMC$ and thus the same point in $\fC$. 

For the induced map $\opsi \in \cGM^*$ acting on $\TP$, its action on the vertex $\xi \in \fC_0$ corresponding to the $\gamma_0$ or $\lambda_0$ curve chosen, defines a sequence of vertices $\opsi^{\ell}(\xi) = \pi_{\fM}(x_{\ell})$ in $\TP$ tending to the tip of the propeller containing $\xi$, as in Figure~\ref{fig:fM0cT}. However, if $k_x$ is even, the point $\xi_{m_x}$ is then the furthest point along this sequence, which is closest to the tip of the propeller. If $k_x$ is odd, $\xi_{m_x+1}$ is the furthest point along the sequence.
In terms of the normal forms for words in $\cGK^*$ defined in Proposition~\ref{prop-standardform} and the normal forms for words in $\cGM^*$ defined in Proposition~\ref{prop-standardformGM}, this implies that a power $\psi^{\ell}$ appearing in \eqref{eq-productformula} can collapse to a  power $\opsi^{m}$ appearing in \eqref{eq-productformulaGM} with $m < \ell$.  We will say that the action of $\cGM^*$ has ``leaf short-cuts'', which correspond to the holonomy of a path which goes from a point on one side of a propeller, to a point on the opposite side, avoiding the trip that the action of $\cGK^*$ must follow to the extremal end of the propeller, to get to the  point on the opposite side.

Finally,  note that   Hypothesis~\ref{hyp-genericW} and the structure of the propellers imply that for  $\xi \in Dom(\opsi)$ we have 
  $r(\opsi(\xi))< r(\xi)$,   as  $r(\psi_-(p_0^1(\xi))) = r(p_0^1(\xi))$.

We next discuss the entropy associated to the action of the $\psg$ $\cGM^*$ on $\fC$. 
Let   $\cT = \cA \cap \bRt$ have the metric defined by the radial coordinate $r$, and endow $\fC \subset \cT$ with the restricted metric, denoted by  $d_{\fC}$. 
 For $\e > 0$,  say that    $\xi_1, \xi_2 \in \fC$ are \emph{$(n,\e)$-separated} if there exists $\ovp \in \cGM^{(n)}$ such that $\xi_1, \xi_2$ are in the domain of $\ovp$, and  $\ds  d_{\fC} (\ovp(\xi_1) , \ovp(\xi_2)) ~ > ~ \e$. 
 A finite set $\cS \subset \fC$ is said to be   \emph{$(n,\e)$-separated} if   each distinct pair $\xi_1 , \xi_2 \in \cS$  is $(n,\e)$-separated.
 Let $s(\cGM^*,n,\e)$ be the maximal cardinality of an $(n,\e)$-separated subset of $\fC$. 

 Then as in \eqref{eq-pseudogroupentropy}, define the entropy of $\cGM^*$ by: 
\begin{equation}\label{eq-lamentropy}
h_{GLW}(\cGM^*) ~ = ~ \lim_{\e \to 0} \left\{ \limsup_{n\to\infty}\frac{1}{n}\ln(s(\cGM^*,n,\e)) \right\} .
\end{equation}

The entropy $h_{GLW}(\cGM^*)$  is closely related to the entropy $h_{GLW}(\cGK^* | \fMC)$, as the  action of  $\cGK^*$ on $\bRt$ and $\cGM^*$ on $\fC$ are ``almost intertwined''  by the projection map  $\pi_{\fM} \colon \fMC \to \fC$ defined in Lemma~\ref{lem-structureM}.
The relation between the entropy for a group action under a factor map  suggests that $h_{GLW}(\cGM^*)$ should be  bounded above by
$h_{GLW}(\cGK^* | \fMC)$,  which vanishes by Theorem~\ref{thm-entropyvanish}.

 The standard argument for factor maps used to prove this result does not actually suffice in the situation we consider, as the   ``leaf short-cuts''  in   $\cGM^*$   (discussed above) imply there are actions in $\cGK^*$ which are collapsed  by the map $\pi_{\fM}$.  We will use instead a more straightforward approach to the estimation of $h_{GLW}(\cGM^*)$  based on the fact that the space $\fC$ is contained in an interval, and the function $\ds \# \cGM^{(n)}$ has subexponential growth.

 \begin{thm}\label{thm-entropyGLW0}
Let $\mK$ be a generic Kuperberg plug, then $h_{GLW}(\cGM^*) = 0$. 
\end{thm}
\proof
  Suppose that $h_{GLW}(\cGM^*) = \lambda > 0$, then there exists $\e > 0$ and a subsequence of sets $\whE_{n_i} \subset \fC$ which are $\e$-separated by elements of $\cGM^{(n_i)}$ and $\# \whE_{n_i} > \exp(n_i \, \lambda/2)$. We show this yields a contradiction.
   Fix a choice of $n = n_i$ such that the set $\whE_n$ satisfies $\# \whE_{n} > \exp(n \, \lambda/2)$ and is $\e$-separated by elements of $\cGM^{(n)}$. 

Then for each pair $\xi \ne \xi' \in \whE_n$ there exists $\ovp \in \cGM^{(n)}$  such that $\ds d_{\fC}(\ovp(\xi), \ovp(\xi')) > \e$. We may assume that the word $\ovp$ has the normal form
$\ds \ovp = \ovp^+ \circ \ovp^-$ where each of $\ovp^+$ and $\ds (\ovp^-)^{-1}$ is a monotone word as in  \eqref{eq-productformulaGM}. 

 Proposition~\ref{prop-asymptoticGM} and    Corollary~\ref{cor-asymptoticM}  imply that the function $\ds \# \cGM^{(n)}$ has subexponential growth. It follows that for $n$ sufficiently large, there exists $\ovp \in \cGM^{(n)}$ and a subset $\ds \whE_{\ovp} \subset \whE_n \cap Dom(\ovp)$ with $\# \whE_{\ovp} > \exp(n \, \lambda/4)$ and for each $\xi \ne \xi' \in \whE_{\ovp}$ we have $\ds d_{\fC}(\ovp(\xi), \ovp(\xi')) > \e$.

Thus we conclude that the set $\ds \ovp(\whE_{\ovp}) \subset \fC \subset \cT$ is $\e$-separated. This is clearly impossible for $n$ sufficiently large, as the length of $\cT$ is finite.
\endproof

 The above proof that  $h_{GLW}(\cGM^*) = 0$  depends fundamentally on the estimate of the growth rate of the reduced word function $\ds \# \cGM^{(n)}$. This estimation follows from  Proposition~\ref{prop-asymptoticGM}, which implies that this function has exponential growth of arbitrarily small exponent, hence must have subexponential growth. 
This  suggests considering a more sensitive entropy-type invariant, which detects growth rates that are subexponential, yet faster than any polynomial function, in order to obtain non-vanishing dynamical invariants of the action of   $\cGM^*$ on $\fC$. It turns out that such invariants exists   in the literature.

 The \emph{slow entropy} of a map was introduced in the works of Katok and Thouvenot \cite{KatokThouvenot1997} and Cheng and Li \cite{ChengLi2010}, and we adapt this idea for the action of $\cGM^*$.
 For $0 < \alpha < 1$,  define  the $\alpha$-entropy of $\cGM^*$, or just the slow entropy, by
\begin{equation} \label{eq-slowentropy}
h_{GLW}^{\alpha}(\cGM^*) ~   =   ~ \lim_{\e \to 0} \left\{ \limsup_{n\to\infty}\frac{1}{n^{\alpha}}\ln(s(\cGM^*,n,\e)) \right\}   
\end{equation}
where $0 \leq h_{GLW}^{\alpha}(\cGM^*)  \leq \infty$. If  $h_{GLW}^{\alpha}(\cGM^*)  > 0$ and $0 < \beta < \alpha$, then $h_{GLW}^{\beta}(\cGM^*)  = \infty$.

Note that for $0 < \alpha < 1$, the function $\exp (n^{\alpha})$
grows faster than any polynomial function, but is slower than any exponential function, so the invariant $h_{GLW}^{\alpha}(\cGM^*)$ has the right character for measuring the complexity of the dynamics of $\cGM^*$

The \emph{entropy dimension} for a continuous  transformation $T \colon X \to X$ of a compact metric space $X$ was introduced by de Carvalho \cite{deCarvalho1997}, and studied further by Cheng and Li \cite{ChengLi2010}. We define an analogous invariant for the $\psg$ action of $\cGM^*$ on $\fC$, 
given by the number $0 \leq D (\cGM^*) \leq 1$ defined by 
\begin{equation}
D(\cGM^*)  = \inf \{ \alpha \mid 0 < \alpha \leq 1 ~{\rm and } ~ h_{GLW}^{\alpha}(\cGM^*)  = 0\} ~ . 
\end{equation}

The proof of Theorem~\ref{thm-entropyGLW0} suggests that to show $h_{GLW}^{\alpha}(\cGM^*) \ne 0$ for some $\alpha$, it suffices to 
 estimate the growth rate of the function $\ds \# \cGM^{(n)}$ more precisely. This, in turn, requires a more precise accounting for what monotone words of the form \eqref{eq-productformulaGM} actually exist in $\cGM^{(n)}$. This is an extremely difficult question to answer in general, but we next describe further hypotheses on the construction of $\mK$, the  notion of ``slow growth''   in  Definition~\ref{def-slowgrowth} and ``fast growth''  in  Definition~\ref{def-fastgrowth}  below,  which makes a lower bound estimate possible.  We obtain three results,   Theorems~\ref{thm-slowentropy}, \ref{thm-fastgrowthentropy} and \ref{thm-hyperbolicentropy} below. Due to the highly technical nature of the proofs and the estimates required, we present in detail only the proof of the following:

\begin{thm}\label{thm-slowentropy}
Let $\Phi_t$ be a generic Kuperberg flow. If the insertion maps $\sigma_j$  have ``slow growth''   in the sense of  Definition~\ref{def-slowgrowth}, then $h_{GLW}^{1/2}(\cGM) > 0$, and thus $1/2 \leq D(\cGM^*)  \leq 1$.
\end{thm}

  The strategy of the proof of Theorem~\ref{thm-slowentropy} is to develop an ``admissibility  criterion'' for strings $I =  (\ell_1, \ldots , \ell_m)$ and $J = (j_1, \ldots , j_m)$ such that for each admissible pair $(I,J)$, we obtain a point $\xi_{(I,J)} \in \fC_0$ by the expression
  \begin{equation}\label{eq-xialtlabels}
\xi_{(I,J)} = \ovp_{(I,J)}(\omega_0) = \opsi^{\ell_m} \circ \ophi_{j_m} \circ \opsi^{\ell_{m-1}} \circ \cdots \circ   \ophi_{j_2} \circ      \opsi^{\ell_1} \circ \ophi_{j_1}(\omega_0) .
\end{equation}

Then observe that the  images of the maps $\ophi_j$ on their domains in $\fMC$    are contained in
  the disjoint compact regions bounded by the parabolic curves $\G_0(a)$ or $\Lambda_0(a)$, according to whether $j=1$ or $2$. 
  The parabolic curves $\G_0(a)$ and $\Lambda_0(a)$ define disjoint compact subsets   $I(\G_0), I(\Lambda_0) \subset \fC$,  where $I(\G_0)$ consists of the points $\xi \in \fC$ whose corresponding path component in $\fMR$ is contained in the closure of the region bounded by $\G_0$, and similarly for $I(\Lambda_0)$. Then choose  $\ve_0 > 0$ such that 
 \begin{equation}\label{eq-defve0}
\ve_0 < \inf \{ d_{\fC}(\pi_\fM(\xi), \pi_\fM(\xi')) \mid \xi \in I(\G_0) ~ , ~ \xi' \in  I(\Lambda_0) \} ~.
\end{equation}
 Two  points $\xi = \ovp_{(I,J)}(\omega_0), \xi' = \ovp_{(I',J')}(\omega_0) \in \fC$ where $J$ and $J'$ have $\ell_m = \ell_{m}' =0$, and terminate with  distinct indices $j_m \ne j'_{m'}$ will then be $\ve_0$-separated in $\fC$.  The strategy is then to construct collections of points  via formula \eqref{eq-xialtlabels} which are $(\cGM^*, \ve_0 , n)$-separated.   The key technical problem is   to   estimate the number of admissible strings which give rise to $(\cGM^*, \ve_0 , n)$-separated  points in this way, using words in $\cGK^*$.

 Recall that the difference between $\cGM^*$ and $\cGK^*$ is that in the former we    allow the leaf short-cuts which replace terms $\psi^{\ell}$ in \eqref{eq-productformula} with   terms $\opsi^{m}$ in \eqref{eq-productformulaGM} where $m \leq \ell$. The leaf short-cuts only arise when applying the map $\psi$ repeatedly to a point  $\xi \in \bRt$ with $z(\xi) < 0$, such that we eventually have $z(\psi^{\ell}(\xi)) > 0$. In particular, they do not arise for the orbits
$p_0(j ; 1 , \ell) = \psi^{\ell}(\phi^+_{j}(\omega_{j}))$ for $j =1,2$.

Introduce the involution $\iota \colon \bRt \to \bRt$ defined by $\iota(r, \pi, z) = (r, \pi, -z)$. For a point $\xi \in \fC$, the action of $\iota$ switches the endpoints 
$p_0^1(\xi)$ and $p_0^2(\xi)$ of the $\gamma_0$ or $\lambda_0$ arc in $\bRt$ through $\xi$, as in the proof of Lemma~\ref{lem-structureM}. In particular, $\iota$ induces the identity map on $\fC$. 

 Extend the symmetric generating set  in \eqref{eq-symmetricgenset} for the $\psg$ $\cGK^*$   by adding the element $\iota$ restricted to $\fMR$, to obtain
 \begin{equation}\label{eq-symmetricgensetaugmented}
\wcGK^{(1)} = \{Id, \iota,  (\psi)^{\pm 1}, (\phi^+_1)^{\pm 1}, (\phi^-_1)^{\pm 1}, (\phi^+_2)^{\pm 1},  (\phi^-_2)^{\pm 1}\} ~.
\end{equation}
Let $\wcGK^*$ denote the augmented $\psg$ generated by this set.  Then the problem is to obtain criteria for the existence for sufficient numbers of  words in $\wcGK^*$
that generate $(\cGM^*, \ve_0 , n)$-separated sets in $\fC$ with sufficient growth rates.

 The proof that the composition $\ovp_{(I,J)}$ in \eqref{eq-xialtlabels} is defined at $\omega_0$ requires technical estimates analogous to those used in   the proof of Theorem~\ref{thm-density}. The existence of a point $\xi_{(I,J)}$
  is interpreted as a statement about the composition of generators of the pseudogroup $\wcGK^*$, which requires a   careful analysis of the dynamics of $\wcGK^*$ near the special points $\omega_i \in \bRt$. This leads to estimates which give   sufficient conditions for $(I,J)$ to be admissible, so that the point $\xi_{(I,J)}$  is well-defined. 
  
For simplicity we assume that $a=0$, for $a$ as in Remark~\ref{rmk-notation}.

Observe that for a word $\ovp_{(I,J)}$ as in \eqref{eq-xialtlabels}, the initial composition    
$ \opsi^{\ell_1} \circ \ophi_{j_1}(\omega_0)$  corresponds to a $\g_0$
or $\lambda_0$ curve whose lower endpoint is the point $p_0(j_1; 1, \ell_1) \in J_0$. Then  $\pi_{\fM}(\phi_{j_2}^+(p_0(j_1; j_2, \ell_1)))$ must lie in the domain of the remaining factor     $\ovp_{(I',J')}$ defined by 
$\ovp_{(I,J)} =  \ovp_{(I',J')} \circ \ophi_{j_2} \circ \opsi^{\ell_1} \circ \ophi_{j_1}$. However, if $j_2 = 2$, it is necessary to apply the involution $\iota$ to the point $p_0(j_1;1,\ell_1)$ to obtain a point in the domain of $\phi_2^+$, as will be seen in the subsequent  construction of separated sets. 

    Recall the   integer valued function $N(r)$ for $2 < r < 3$ introduced in Section~\ref{sec-global}, which is defined using the function $\delta(r)$ in   \eqref{eq-minradineq}.  The function $N(r)$ is an upper bound on the number of insertion maps that can be applied to a point $\xi \in \bRt$ with $r(\xi) = r$, and $N(r)$ is unbounded as $r$ decreases to $r=2$. 
      In particular, $N(r)$ for $r=r(\phi_{j_2}^+(p_0(j_1;j_2,\ell_1)))$ provides an upper bound on the number of subsequent maps $\ophi_{j_k}$ which can appear in the remaining term $\ovp_{(I',J')}$ for an admissible word $\ovp_{(I,J)}$.    
    Thus, a more precise estimate on the growth rate of the function  $N(r)$ for $r> 2$ will yield     estimates on the growth rates of sets of $(\e,n)$-separated points for $\cGM^*$.

   The growth rate of the function  $N(r)$ for $r> 2$ is closely related to the geometry of the embedding maps $\phi^+_i$ for $i=1,2$ near the special points $\omega_i \in \bRt$ as defined in    \eqref{eq-omegas}, which we next consider. We use the assumptions in Hypotheses~\ref{hyp-genericW} and \ref{hyp-parabolic} to analyze the properties of the maps $\phi_i^+$ on sufficiently small neighborhoods of the points $\omega_i$.
   
 Recall that  the first transition point for the forward $\cK$-orbit of    $\omega_i$ is the special entry point  $p_i^{-} \in E_i$  and its backward $\cK$-orbit   is the special exit point  $p_i^{+}   \in S_i$ as     illustrated in  Figure~\ref{fig:KR}. For  $\e > 0$, define the closed squares in $\bRt$ centered on the special points $\omega_i$ for $i=1,2$,
$$ S_{\bRt}(\omega_i ,\e) = \{(r,\pi, z) \mid |r-2| \leq \e  ~ {\rm and }~ |z - (-1)^i| \leq \e  \} ~ ,$$
and let  $\ds  D_{\bRt}(\omega_i, \e) \subset \bRt$ be the closed ball centered at $\omega_i$ with radius $\e$. Then 
$\ds S_{\bRt}(\omega_i ,\e/\sqrt{2}) \subset D_{\bRt}(\omega_i ,\e)$.

We next choose $\e_0' > 0$ so that the various generic hypotheses apply to the points in $S_{\bRt}(\omega_i ,\e_0')$. Recall that  the constant  $\e_0 > 0$ was chosen   in Hypothesis~\ref{hyp-genericW}  so that  the estimate \eqref{eq-quadraticest2} 
holds    for   the Wilson vector field $\cK$ on the disk $\ds  D_{\bRt}(\omega_1, \e_0)$. Recall that for the same $\e_0$,  we   assumed that   Hypothesis~\ref{hyp-parabolic}   on   the insertion maps   $\sigma_i$ for $i = 1,2$, holds for  $2 \leq r_0 \leq 2 + \e_0$ and $\theta_i -\e_0 \leq \theta \leq \theta_i + \e_0$.  

Choose    $0 < \e_0' \leq \e_0/\sqrt{2}$   sufficiently small so  that    $S_{\bRt}(\omega_i, \e_0') \subset Dom(\phi^+_i)$ for $i=1,2$.

Moreover, require that the forward $\Phi_t$-flow of $S_{\bRt}(\omega_i, \e_0')$  to the surface $E_i$ is contained  in the rectangular domain for Hypothesis~\ref{hyp-parabolic}.  This   automatically holds for the $r$ coordinate, as the flow $\Phi_t$ preserves the radius coordinate along these trajectories, but the condition $\theta_i -\e_0 \leq \theta \leq \theta_i + \e_0$ on the image in $E_i$   imposes a restraint on the choice $\e_0'$.

Thus, Hypothesis~\ref{hyp-genericW} applies for the Wilson return map $\psi$ near $\omega_i$ and Hypothesis~\ref{hyp-parabolic} applies for the radial coordinates of the return map $\phi_i^+$ near $\omega_i$.

Consider first  the map $\phi^+_1$.  By       Hypothesis~\ref{hyp-SRI},  the ``parabolic function'' $z \mapsto r(\phi^+_1(2,  \pi, z))$ has a minimum value $2$ at  $z=-1$.  Then
 \begin{equation}\label{eq-increasingr}
2 = r(\phi^+_1(\omega_1))   ~ < ~ r(\phi^+_1(r,\pi, z)))  ~ {\rm for ~ all} ~  (r,\pi, z) \ne \omega_1 ~ {\rm with}~ r \geq 2 .
\end{equation}
  Moreover,    Hypothesis~\ref{hyp-parabolic} implies   that  the images of the vertical lines $r=c$ for $2 \leq c \leq 2 +\e_0'$  in the neighborhood of  $\omega_1$ are   mapped by $\phi^+_1$ to parabolic curves, so that  $z \mapsto r(\phi_1^+(c,  \pi, z))$   has a   minimum value at a unique value  $z= \zeta_1(c)$.  
Thus the forward orbit of $(c,\phi,\zeta_1(c))$ intersects $E_1$ in a point with $z$-coordinate equal to $-1$. 
  Set $\rho_1(c) = r(\phi_1^+(c,  \pi, \zeta_1(c)))> 2$.
  
For  each $2 \leq c \leq 2 +\e_0'$ the function $z \mapsto r(\phi_1^+(c,  \pi, z))$ has vertex point $(\zeta_1(c) , \rho_1(c))$.  The  graph of this function  near the point $(\zeta_1(c) , \rho_1(c))$  has upward parabolic shape,  where   $\zeta_1(c)$  gives the ``offset'' of the parabolic vertex along the $z$-axis.

Hypothesis~\ref{hyp-parabolic} implies that  the function   $c \mapsto \zeta_1(c)$  is a smooth   function of $c$ with $\zeta_1(2) = -1$ and $-1-\e_0' \leq \zeta_1(c) \leq -1 + \e_0'$. We also have   that the function  $c \mapsto \rho_1(c)$     is   smooth,   with 
   $\rho_1(2) =2$ and  $c \mapsto \rho_1(c)$  is strictly increasing for $c > 2$. 
  Moreover, there    exists   $0 < \alpha_1 \leq \beta_1$   such that for  $2 \leq c \leq 2 +\e_0'$ we have the quadratic bounds
 \begin{equation}\label{eq-localest1}
   \rho_1(c) + \alpha_1  \cdot  (z- \zeta_1(c))^2  ~ \leq ~ r( \phi^+_1(c, \pi, z) )  \leq  \rho_1(c) + \beta_1  \cdot (z- \zeta_1(c))^2   ~ .
\end{equation}

The   map $\phi^+_2$ near $\omega_2$ admits a similar analysis, yielding functions $c \mapsto \zeta_2(c)$ and $c \mapsto \rho_2(c)$ with the properties 
$1-\e_0' \leq \zeta_2(c) \leq 1 + \e_0'$ for $2 \leq c \leq 2 +\e_0'$ and   constants  $0 < \alpha_2 \leq \beta_2$ for which there is  the quadratic estimate
\begin{equation}\label{eq-localest2}
   \rho_2(r) + \alpha_2  \cdot   \, (z- \zeta_2(c))^2  ~ \leq ~ r( \phi^+_2(r, \pi, z) )  \leq  \rho_2(r) + \beta_2  \cdot   \, (z- \zeta_2(r))^2   ~ .
\end{equation}

Set $\alpha_{\Phi} = \min\{ \alpha_1, \alpha_2\}$ and $\beta_{\Phi} = \max \{\beta_1, \beta_2\}$.

We also have that the function  $c \mapsto \zeta_i(c)$ is smooth with $\zeta_i(2) = (-1)^i$, for $i =1,2$,  so there exists a constant $\whzeta \geq 0$ so that 
\begin{equation}\label{eq-zetaest}
|\zeta_i(c) - (-1)^i| \leq \whzeta \cdot |c-2| \quad {\rm for}~ 2 \leq c \leq 2 + \e_0' .
\end{equation}
That is, the vertical offset of each parabolic curve  $z \mapsto r(\phi_1^+(c,  \pi, z))$ has a linear bound as a function of $(c-2)$,  for $2 \leq c \leq 2 + \e_0'$.

Recall that we require an estimate on the function $N(r)$   introduced in Section~\ref{sec-global}, for $r$ near $2$, and this is defined in terms of the function $\delta(r)$   defined in   \eqref{eq-minradineq}. Using the normal forms  \eqref{eq-localest1} and \eqref{eq-localest2},  observe    that $\ds \delta(r) = \min \{\rho_1(r), \rho_2(r) \}$ for $2 \leq r \leq 2 +\e_0'$. We introduce two classes for the functions $\rho_i$ which are determinant in the behavior of the function $N(r)$.
Note that each function $\rho_i$ is smooth and strictly increasing, so its derivative $\rho_i'(r) \geq 1$ for all $2 \leq r \leq 2 +\e_0'$.   
 \begin{defn}\label{def-slowgrowth}
The insertions $\sigma_j$ are said to have \emph{slow growth} if   the derivatives  $\rho_i'(2) =1$ for $i=1,2$.
In this case,  there exists $C_{\rho} > 0$ and   uniform estimates on $\rho_{i}(r)$, 
\begin{equation}\label{eq-slowgrowthest}
 r < \rho_i(r) ~ \leq ~ r + C_{\rho} \, (r-2)^2 \quad {\rm for} ~ 2 < r \leq 2 + \e_0' .
\end{equation}
\end{defn}

 \begin{defn}\label{def-fastgrowth}
The insertions $\sigma_i$ are said to have   \emph{fast growth} if  there exists $\lambda > 1$ such that derivatives  $\rho_i'(r) \geq \lambda$ for $2 \leq r \leq 2 + \e_0'$ with $i=1,2$. 
In this case, we have the uniform estimate, 
\begin{equation}\label{eq-fastgrowthest}
  \rho_i(r) ~  \geq  2 + \lambda \, (r-2)  \quad {\rm for} ~ 2 \leq r \leq 2 + \e_0' .
\end{equation}
\end{defn}
 
  These two properties of the insertion maps are sufficient to obtain the required estimates  on the domains for the maps   $\ovp_{(I,J)}$ as discussed above. 
We also   require  some preliminary estimates   which describe the quantitative behavior  of the map $\psi$  near the special point $\omega_1$. The derivation of these  estimates is similar to those in Section~\ref{sec-generic}, so the arguments are only briefly sketched.

Let $\xi = (r,\pi,z) \in \bRt$.  Define times  $0  = T_0(\xi) <  T_1(\xi) <  T_2(\xi) <  \cdots $   where  $\psi^{\ell}(\xi) = \Psi_{T_{\ell}(\xi)}(\xi)$ for $ \ell \geq 0$ such that $\psi^{\ell}(\xi)$ is defined.
Then $r(\psi^{\ell}(\xi)) = r(\xi)$, and by \eqref{eq-coordinates} we have
\begin{equation}\label{eq-generalsteps}
z(\psi^{\ell+1}(\xi)) -  z(\psi^{\ell}(\xi)) ~ = ~ \int_{T_{\ell}(\xi)}^{T_{\ell+1}(\xi)} ~ g(\Psi_s(\xi)) ~ ds    .
\end{equation}

The return time for the flow $\Psi_t$ at $\xi$    is $2 \pi \cdot r(\xi)$,  so for  $\xi$ with  $2 \leq r(\xi) \leq 2 +\e_0$,    the domain of the integral in \eqref{eq-generalsteps}  satisfies   $4 \pi \leq T_{\ell+1}(\xi) - T_{\ell}(\xi) \leq (4+2\e_0') \pi$.   
Moreover,   Hypothesis~\ref{hyp-genericW} implies that  $g(r,\theta,z)$ is a non-decreasing function of $r \geq 2$, hence      \eqref{eq-generalsteps} yields 
\begin{equation}\label{eq-increasing}
z(\psi^{\ell}(r',\pi, z)) \geq z(\psi^{\ell}(r, \pi, z)) \quad {\rm  for} ~ r' > r \geq 2 .
\end{equation}
  Recall that $g(\xi) =1$ if $r(\xi) \geq 2 +\e_0$, so for $\xi \in \bRt$ with $r(\xi) \geq 2+\e_0$, the orbit $\Psi_t(\xi)$ escapes the plug $\mW$ at a time 
  $t = 4 - z(\xi) \leq 6$, hence does not complete a full revolution. 
 Moreover,  if  $\xi \in \bRt$ satisfies $r(\xi) =2$ and  $-2 \leq z \leq -1-\e_0$, then $\psi(\xi)$ is defined and satisfies $-(1+\e_0) < z(\psi(\xi)) < -1$.

       In particular, let $\xi_0 = \phi_j^+(\omega_j)$ for   $j=1,2$,  so that  $r(\xi_0) = 2$, then  we have    $$z(p_0(j ; 1 , 0)) = z(\psi(\xi_0)) > -(1+\e_0) .$$     
  We can thus take  $\ell_0 = 1$   in  Lemma~\ref{lem-density2}  so that for $\ell \geq 1$, 
  the estimates \eqref{eq-spacing1} and \eqref{eq-spacing2}  hold for $ \psi^{\ell}(\xi_0) = p_0(j ; 1 , \ell)$.     
   
\begin{lemma}\label{lem-r=2verticalapproach}
There exists   constants $0 < \mu_1 = 4 \pi  \lambda_1 \leq 4\pi  \lambda_2 < \mu_2 $ and $b_0 > 0$ such that 
for $\ell \geq b_0$, for $j=1,2$, the point $\ds p_0(j ; 1 , \ell) $  satisfies 
  \begin{equation}\label{eq-spacing1aa}
   -1 -   1/(\mu_1 \, \ell) ~ < ~   z(p_0(j ; 1 , \ell) )  ~  <   ~ -1 -  1/(\mu_2 \, \ell) ~ .
\end{equation}
\end{lemma}
 \proof
 Following the notation as in Lemma~\ref{lem-density2}, set $\mu_1  = 4\pi \lambda_1$ and choose $4 \pi \lambda_2 < \mu_2 < 8 \pi\lambda_2$. Then for 
 $b_0 \geq C_3/(\mu_2   - 4\pi \lambda_2)$, \eqref{eq-spacing1aa} follows from \eqref{eq-spacing1} of Lemma~\ref{lem-density2}. 
 \endproof

    Consider next the   case where $\xi \in \bRt$ with $-2 \leq z(\xi) < -1$ and $2 < r(\xi) < 2 + \e_0'$.   Then the    $\cW$-orbit   $\Psi_t(\xi)$ escapes from $\mW$ in finite time, and thus   $\psi^{\ell}(\xi)$ is only defined for a finite range of $\ell$, where $\ell \to \infty$ as $r(\xi)$ approaches $2$.  The next two results give an estimate for the range of $\ell$ for which $z(\psi^{\ell}(\xi)) \leq -1+\e_0'$, using methods analogous to those  used in the proofs of  Lemmas~\ref{lem-density3}, \ref{lem-density4} and \ref{lem-density5}.

  \begin{lemma}\label{lem-spacing}
There exists   constants $U_g > 0$ and    $0 < \e_1 < \e_0'$ such that for all $0 < \e \leq \e_1$ and $\xi \in \bRt$ with $z(\xi) \leq -(1+\e)$ and $2 \leq r(\xi) \leq 2+\e$, then there exists $\ell_{\xi} > 0$ so that 
\begin{equation}\label{eq-spacingest}
-1 - \e ~ \leq ~  z(\psi^{\ell}(\xi)) ~ \leq ~ -1 + \e \quad {\rm for ~ all} \quad \ell_{\xi} \leq \ell \leq \ell_{\xi} + U_g/\e ~ .
\end{equation}
  \end{lemma}
 \proof
 Recall from Section~\ref{sec-proplevels} that  
 Hypothesis~\ref{hyp-genericW}  and condition \eqref{eq-generic1} imply there exists constants $A_g,B_g,C_g$ such that  the quadratic form  $Q_g(u,v) = A_g \, u^2 + 2B_g \, uv + C_g \,  v^2$ defined by the Hessian of $g$ at $\omega_1$ is positive definite. 
  The remainder term for the quadratic Taylor approximation to $g(r,\theta,z)$ is dominated by a scalar multiple of   $\ds  d_{\mW}(\xi, \cO_1)^3$, where   $\ds d_{\mW}(\xi, \cO_1) = \sqrt{(r-2)^2 + (z+1)^2}$ denotes the distance from $\xi = (r,\theta,z)$ to the periodic orbit $\cO_1$. 
 As   $Q_g(u,v)$ is positive definite,  there exists   $0 < \e_2 \leq \e_0'$ so that for $\xi = (r,\theta,z)$ 
 \begin{equation}\label{eq-quadraticest}
| g(\xi) - Q_g(r-2,z+1) | ~ \leq ~ Q_g(r-2,z+1)/4 \quad {\rm for} \quad d_{\mW}(\xi, \cO_1) \leq \e_2 .
\end{equation}

Let $\lambda_2$ denote the maximum eigenvalue of the quadratic form $Q_g$. Then for $0 < \e \leq \e_2/\sqrt{2}$  we have  by  \eqref{eq-quadraticest}     
  \begin{equation}\label{eq-smallestimate}
\max   \left\{ g(\xi) \mid \xi = (r, \theta, z) ~ , ~ | z+ 1| \leq  \e  ~ {\rm and} ~ | r -2| \leq \e  \right\} ~ \leq ~ 4 \lambda_2 \, \e^2 ~ .
\end{equation}
 
 Let $\e_1 = \min \{ \e_2 /\sqrt{2},   1/(24   \pi \, \lambda_2) \}$.    
 
 Given   $0 < \e \leq \e_1$  and $\xi \in \bRt$ with $z(\xi) \leq -(1+\e)$ and $2 \leq r(\xi) \leq 2+\e$, we     show    \eqref{eq-spacingest} holds.
 
 First, note   there exist  $t_1 \geq 0$ so that   $z(\Psi_{t_1}(\xi)) = -(1 + \e)$.  
 Let $t_2 \geq t_1$ be the first subsequent time for which $\Psi_{t_2}(\xi) \in \bRt$.
 Then there exists $\ell_{\xi} \geq 0$ so that $\ds \Psi_{t_2}(\xi) = \psi^{\ell_{\xi}}(\xi)$. Set $\xi_* = \Psi_{t_2}(\xi)$.

Note that $z(\xi_*) \geq -(1+\e)$ and $t_2 - t_1 \leq 2\pi \, r(\xi) < 6 \pi$,    so by \eqref{eq-generalsteps} and \eqref{eq-smallestimate} we have
  \begin{equation}\label{eq-firststepmax}
0 ~  \leq  ~ z(\Psi_{t_2}(\xi))  - z(\Psi_{t_1}(\xi)) ~ < ~  4\lambda_2 \, \e^2 \cdot  6 \pi ~ = ~  24 \pi \lambda_2 \, \e^2 ~ .
\end{equation}
As $z(\Psi_{t_1}(\xi)) = -(1 + \e)$ this yields 
$$- (1 + \e) ~ \leq ~  z(\xi_*) ~  < ~ - (1 + \e) + 24 \pi \lambda_2 \, \e^2   ~ \leq ~  - (1 + \e) +   \e ~ = ~ -1 $$
 The same reasoning   yields, for $\ell > 0$ such that $\ds z(\psi^{\ell}(\xi_*)) \leq -1 + \e$, then
 \begin{equation}\label{eq-secondstep}
  z(\psi^{\ell}(\xi_*)) - z(\psi^{\ell -1}(\xi_*)) ~ < ~      24   \pi \lambda_2 \, \e^2  ~ \leq  ~   \e  ~ .
 \end{equation}
    Set $U_g = 1/(24   \pi \lambda_2)$. 
  It then follows by applying \eqref{eq-secondstep} recursively, that for $\ell \leq U_g/\e$ we have
  \begin{equation}\label{eq-finalstep}
0 <  z(\psi^{\ell}(\xi_*)) - z(\xi_*) \leq \ell \cdot    24   \pi \lambda_2 \, \e^2 \leq   (U_g/\e) \cdot  (24   \pi \lambda_2 \, \e^2)   ~ = ~ \e  ~ .
 \end{equation}
It follows that for $0 \leq \ell \leq U_g/\e$ we have
\begin{equation*}
 - (1 + \e) ~ \leq ~  z(\xi_*) ~  < ~ z(\psi^{\ell}(\xi_*))   ~   \leq  ~     z(\xi_*) +   \e ~  \leq ~ -1 + \e 
\end{equation*}
which  was to be shown.  
  \endproof

  Next, we give a form of   density estimate for the flow $\Phi_t$.
  Recall that for $j=1,2$, the function  $\zeta_j(c)$ is the ``offset'' of the vertex of the parabolic graph  $z \mapsto \phi_j^+(c,\pi, z)$. The function satisfies      $\zeta_j(2) = (-1)^j$ and there is  $\whzeta > 0$ so that the bounds 
\eqref{eq-zetaest} hold for $j=1,2$.

 \begin{lemma}\label{lem-spacing3} Let $0 < \e_1 < \e_0'$ be the constant of Lemma~\ref{lem-spacing}. Then there exists $L_g > 0$ so that 
  for all $0 < \e \leq \e_1$ and $\xi \in \bRt$ with $-(1+\e) \leq z(\xi) \leq -(1+\e/2)$ and $2 \leq r(\xi) \leq 2+\min\{1,  1/\whzeta\}\cdot \e/2$, then there exists $0 \leq \ell_* \leq L_g/\e$ 
  so that 
\begin{equation}\label{eq-spacingest2}
0 <  | z(\psi^{\ell_*}(\xi)) - (-1)^j \, \zeta_j(r(\xi))|  ~ < ~  \e ~ .
\end{equation}
\end{lemma}
\proof
 If $\zeta_1(r(\xi)  \leq  -1$ or  $\zeta_2(r(\xi)  \geq  1$, then
 there is nothing to show. In fact, we need only consider the  case
 where $r_0 = r(\xi)$ is such that    $ 2 < r_0 \leq \e/2\cdot \min\{1,  1/\whzeta\}$ and $-1 + \e \leq \zeta_1(r_0)  <   -1 + \whzeta \cdot (r_0-2) \leq -1 + \e/2$.
It then suffices to estimate the least value of $\ell_* > 0$ such that $\psi^{\ell_*}(\xi) > -1 -\e/2$.

Let $\lambda_1$ denote the minimum eigenvalue of the quadratic form $Q_g$. 
By  \eqref{eq-quadraticest}    there is a lower bound
 \begin{equation}\label{eq-smallestimate2}
\min ~ \left\{ g(\xi) \mid \xi = (r, \theta, z) ~ , ~ \e/2 \leq | z+ 1| \leq  \e  ~ {\rm and} ~ | r -2| \leq \e  \right\} ~ \geq ~  \lambda_1 \, \e^2/(4 \cdot 1.01) \geq \lambda_1 \, \e^2/5 ~ .
\end{equation}
By \eqref{eq-generalsteps} and \eqref{eq-smallestimate2}, and using that  the flow $\Phi_t$ has return time at least $4\pi$ for $r_0 > 2$, 
if $-1-\e \leq z(\xi) < z(\psi(\xi)) \leq -1-\e/2$, then  we have
  \begin{equation}\label{eq-firststepmin}
0 ~  \leq  ~ z(\psi(\xi))  - z(\xi) ~ \geq  ~ 4 \pi \, \lambda_1 \, \e^2/5   ~ .
\end{equation}
Thus, the least   $\ell_*$ satisfies the estimate $\ds \ell_* \cdot 4 \pi \, \lambda_1 \, \e^2/5 \leq  \e/2$,  so for $L_g = 5/(8 \pi \, \lambda_1)$ there exists $\ell_* \leq L_g/\e$ such that \eqref{eq-spacingest2} holds. 
\endproof

After  these preliminary considerations, we return to the proof of Theorem~\ref{thm-slowentropy}.
We consider monotone words in $\cGM^*$ of the form \eqref{eq-productformulaGM}, and their lifts to monotone words in $\wcGK^*$.

 We assume that the  insertion maps $\sigma_j$ satisfy Definition~\ref{def-slowgrowth}, with $C_{\rho}$ the constant so that \eqref{eq-slowgrowthest} is satisfied. 
    Also, recall $\beta_{\Phi} = \max \{\beta_1, \beta_2\}$ for the constants $\beta_j$ as defined by \eqref{eq-localest1}, \eqref{eq-localest2},  and $\whzeta \geq 0$ was defined so that 
the estimate \eqref{eq-zetaest} holds. Then define $\ds C_{\Phi} = \max \{1, C_{\rho} , \beta_{\Phi}\}$.

 Let   $0 < \e < \min \{\e_1 , 1/C_{\Phi}\}$ for   $\e_1$   as  defined   in Lemma~\ref{lem-spacing}.

By Lemma~\ref{lem-r=2verticalapproach}, there exists an integer 
   $b_1 \geq 1$   such that for $\ell \geq b_1$, for $k=1,2$, the point $\ds p_0(k ; 1 , \ell) $  satisfies the estimates in \eqref{eq-spacing1aa}. 
       For $m \geq 1$,  let  $b_m  = m \, b_1$. Then   $1/(\mu_1 \, b_m)  \leq \e/m \leq 1/(m C_{\Phi})$.

Now consider strings $I =  (\ell_1, \ldots, \ell_m)$ and $J
= (j_1, \ldots , j_m)$. We develop a criteria for when the point
$\xi_{(I,J)} = \ovp_{(I,J)}(\omega_0)$ as in \eqref{eq-xialtlabels} is defined.   Take $p_0(j_1;1,0) = \phi^+_{j_1}(\omega_{j_1})$ for $j_1 =1,2$.
Assume that  $\ell_1 \geq b_m$ and set       $(r_1,\pi, z_1)  =
\psi^{\ell_1}(p_0(j_1;1,0))= p_0(j_1;1,\ell_1)$. Then $r_1 =2$ and
$z_1<-1$. Set $\upsilon_1 = |z_1 +1|$ so that $0 < \upsilon_1 \leq \e/m$ by the choice of $\ell_1$.

  For  $j_2 = 1$, we have $\phi^+_{1}(p_0(j_1;1,\ell_1))=p_0(j_1;1,\ell_1;1,0)$, 
  while for $j_2 = 2$, we use the involution $\iota$ to obtain  
  $$\phi^+_{2}(\iota(p_0(j_1;1,\ell_1)))=\phi_2^+(p_0(j_1;2,\ell_1))=p_0(j_1;  2,\ell_1;1,0).$$ 
These points are well-defined by the choice of $\ell_1$. 

Set  $r_2 = r(p_0(j_1;j_2,\ell_1; 1,0))$,  then by  the quadratic estimates \eqref{eq-localest1} or  \eqref{eq-localest2}  for $r_1 =2$ and so $\rho_{j_1}(r_1) =2$,   we have
\begin{equation}\label{eq-bounds2}
2 +  \alpha_{\Phi}  \cdot    \upsilon_1^2  ~ \leq ~ r_2   \leq  2 +  \beta_{\Phi}  \cdot     \upsilon_1^2 \leq 2 + \beta_{\Phi} \cdot (\e/m)^2 < 2 +  \e/m^2
\end{equation}
since $\beta_{\Phi} \, \e \leq C_{\Phi} \, \e < 1$. Thus, $r_2  < 2 + \e/m$. 

Now set $b_m' = b_m   + m\, L_g/\e$ for $L_g = 5/(8 \pi \, \lambda_1)$ as defined in Lemma~\ref{lem-spacing3}, so $b_m' \leq m \, \{ 1/\mu_1 + 5/(8 \pi \, \lambda_1)\}/\e$. 
 Then by Lemmas~\ref{lem-spacing} and \ref{lem-spacing3} applied for $\e/m$ we can  choose $\ell_2 \leq b_m'$ so that 
 $$(r_2,\pi, z_2) = p_0(j_1;j_2,\ell_1; 1,\ell_2) = \psi^{\ell_2}(p_0(j_1;j_2,\ell_1; 1,0))$$ 
 is  defined with $|z_2 + 1| \leq \e/m$, and  the 
 vertical ``offset'' along the line $r=r_2$ is given by 
\begin{equation}\label{eq-bounds2a}
\upsilon_2 = |z_2 - \zeta_{j_2}(r_2)| \leq \e/m  ~ .
\end{equation}
  For  $j_3 = 1,2$, we get
  \begin{eqnarray*}
p_0(j_1;j_2,\ell_1;1, \ell_2;  1,0)   &  = &  \phi^+_{1}(p_0(j_1;j_2,\ell_1;1, \ell_2))  \\
p_0(j_1; j_2, \ell_1; 2, \ell_2; 1,0)  &  = &    \phi^+_{2}(\iota(p_0(j_1;j_2,\ell_1;1, \ell_2)))  
\end{eqnarray*}

 Let $r_3 = r(p_0(j_1; j_2,\ell_1;j_3, \ell_2; 1,0))$, then by  the quadratic estimates \eqref{eq-localest1} or  \eqref{eq-localest2},   the ``slow estimate''  \eqref{eq-slowgrowthest}, and the inductive estimates \eqref{eq-bounds2} and \eqref{eq-bounds2a}, we have
\begin{eqnarray*}\label{eq-bounds3}
  r_{3}   ~ & \leq &  ~  \rho_{j_2}(r_2) + \beta_{\Phi}  \cdot     \upsilon_2^2   \\
   ~ & \leq & ~  [r_2 + C_{\rho}\cdot (r_2 -2)^2] + [\beta_{\Phi}  \cdot   (\e/m)^2] \\
   ~ & \leq & ~ [2 + \beta_{\Phi} \cdot (\e/m)^2 ]   +  [ C_{\rho} \cdot  (\beta_{\Phi} \cdot (\e/m)^2)^2 ]  + [ \beta_{\Phi}  \cdot   (\e/m)^2 ] \\ 
   ~ & = & ~ 2   + 2 \beta_{\Phi}  \cdot   \e^2/m^2 +    C_{\rho} \cdot \beta_{\Phi}^2 \cdot \e^4/m^4    \\
   ~ & \leq & ~ 2   + 2  \, \e/m^2 +  C_\rho \cdot  \e^2/m^4 < 2   + 3 \, \e/m^2          
\end{eqnarray*}
where the last inequality follows from $C_\rho\cdot \e^2/m^4\leq C_\Phi\cdot \e^2/m^4<\e/m^4\leq \e/m^2$.

Then by Lemma~\ref{lem-spacing} applied for $\e/m$, we can  choose $\ell_3 \leq b_m'$ so that 
 $$(r_3, \pi, z_3) = p_0(j_1; j_2,\ell_1; j_3,\ell_2; 1,\ell_3) =
 \psi^{\ell_3}(p_0(j_1;j_2,\ell_1; j_3,\ell_2; 1,0))$$ 
 is  defined with $|z_3 + 1| \leq \e/m$, and  the   vertical ``offset'' along the line $r=r_3$ is given by 
 $$\upsilon_3 = |z_3 - \zeta_{j_3}(r_3)| \leq \e/m  ~  . $$

We continue recursively, assuming the definitions as above  for $1 \leq i \leq k$, and that the following inductive assumptions hold for $1 < i \leq k$:  
\begin{equation}\label{eq-recursive}
\ell_i \leq b_m' \quad , \quad |z_{i} + 1| < \e/m \quad , \quad 2 \leq r_{i} \leq ~ 2 + (1 + 2i) \, \e/m^2 \quad , \quad  \upsilon_{i} \leq \e/m
\end{equation}
 Also assume that $r_k \leq  \min\{1,  1/\whzeta\}\cdot \e/2$. 
Then for  $j_{k+1} = 1,2$, we get 
 \begin{eqnarray*}
p_0(j_1;j_2,\ell_1;j_3,\ell_2; \cdots ; 1,\ell_k;  1,0)  &  = &
\phi^+_{1}(p_0(j_1; j_2,\ell_1;j_3,\ell_2; \cdots ; 1,\ell_k))   \\
p_0(j_1;j_2,\ell_1;j_3,\ell_2;\cdots; 2,\ell_k;1,0)   &  = &    \phi^+_{2}(\iota(p_0(j_1;j_2,\ell_1;j_3,\ell_2; \cdots ; 1,\ell_k)))    
\end{eqnarray*}
and set 
$r_{k+1}   =  r(p_0(j_1;j_2,\ell_1;j_3,\ell_2; \cdots ; j_{k+1},\ell_k; 1,0))$.

We require one additional assumption and a small calculation to complete the inductive step of the construction. Suppose that $(1+2k) \leq m$, then $C_{\rho}\, \e < 1$ by choice of $\e$, so $(1+2k) < m/\sqrt{C_{\rho}\, \e}$.
Thus $(1+2k)^2 < m^2/\e C_{\rho}$ and so $(1+2k)^2 (\e/m^2) < 1/C_{\rho}$  
which yields $C_{\rho} ((1+2k)\e/m^2)^2 <  \e/m^2$.

Finally,  use  the quadratic estimates \eqref{eq-localest1} or  \eqref{eq-localest2},   the ``slow estimate''  \eqref{eq-slowgrowthest},   the inductive estimates \eqref{eq-bounds2} and \eqref{eq-bounds2a}, and the small calculation above to obtain
\begin{eqnarray*}\label{eq-boundsk}
r_{k+1}  ~ & \leq &  ~  \rho_{j_k}(r_k)  + \beta_{\Phi}  \cdot     \upsilon_k^2   \\
   ~ & \leq & ~  [  r_k + C_{\rho}\cdot (r_k -2)^2 ] + [ \beta_{\Phi}  \cdot   (\e/m)^2  ] \\
   ~ & \leq & ~  2 + (1 + 2k) \, \e/m^2   + \e/m^2  +   C_{\rho} ((1+2k)\e/m^2)^2  \\
   ~ & < & ~   2  (1 + 2k) \, \e/m^2  + \e/m^2  +   \e/m^2 \\
   ~ & = & ~    2  (1 + 2(k+1)) \, \e/m^2  
\end{eqnarray*}
Then by Lemma~\ref{lem-spacing} applied for $\e/m$, we can  choose $\ell_{k+1} \leq b_m'$ so that 
 $$(r_{k+1}, \pi,  z_{k+1})   =  p_0(j_1;j_2,\ell_1;j_3,\ell_2; \cdots ; j_{k+1},\ell_k;  1,\ell_{k+1}) = \psi^{\ell_{k+1}}(p_0(j_1;j_2,\ell_1;j_3,\ell_2; \cdots ; j_{k+1},\ell_k;  1,0))$$ 
 is  defined with $|z_{k+1} + 1| \leq \e/m$, and  the   vertical ``offset'' along the line $r=r_{k+1}$ is given by 
 $$\upsilon_{k+1} = |z_{k+1} - \zeta_{j_{k+1}}(r_{k+1})| \leq \e/m  ~  . $$ 
which completes the recursive step.

  We use the constructions above to obtain lower bound estimates on
  the number $s(\cGM, n, \varepsilon)$ of $(n,\varepsilon)$-separated
  words for the action of $\cGM^*$ on $\fC$.

For $0 < \e < \min \{\e_1 , 1/C_{\Phi}\}$ as before  and     $\delta =  \min\{1,  1/\whzeta\}\cdot \e/2$, 
we construct orbits in the rectangular regions in $\bRt$ centered on the special orbits $\omega_j$ for $j=1,2$, 
\begin{equation}
\left\{ (r,\pi, z) \in \bRt \mid |z - (-1)^j| < \e ~ , ~ 2 \leq r \leq \delta \right\}
\end{equation}
As the value of $\e>0$ tends to $0$, the density of such points increases as well, so that one observes the slow entropy of $\cGM^*$ is concentrated in these   regions around the special orbits.

Let   $b_1 \geq 1/(\mu_1 \e)$   be such that for all $\ell \geq b_1$, for $j=1,2$, the point $\ds p_0(j ; 1 , \ell) $  satisfies the estimates in \eqref{eq-spacing1aa}.

    Choose $m \geq 1$,  and    set $b_m = m \, b_1$. As  in  Lemma~\ref{lem-spacing3},  let $b_m'$ be the greatest integer satisfying
$$b_m' ~ \leq ~  b_m   + m\, L_g/\e ~ \leq ~   \{ 1/\mu_1 + 5/(8 \pi \, \lambda_1)\} \cdot m/\e  ~ .$$ 

 Let $k_m$ be the greatest integer for which  $k +1\leq m/4$, then we have the bound $ 2 (1 + 2(k_m+1)) \, \e/m^2 < \e/m$.     So by the recursive procedure above, for $k \leq k_m$       we can realize the point  
    $ \xi_{(I,J)} = \ovp_{(I,J)}(\omega_0) $ defined by \eqref{eq-xialtlabels}, where there are $2^k$ choices of the string $J = (j_1, \ldots , j_k)$, 
          for a fixed strong 
    $I =  (\ell_1, \ldots , \ell_k)$ satisfying $\ell_1 \leq b_m$ and $\ell_i \leq  b_m'$ for $1 < i \leq k$. 
        Note that such a word   has the length estimate
   $$ \| \ovp_{(I,J)}\| \leq k \cdot   b_m'  \leq b_1\,  k_m^2/4  ~ .$$
Recall that $\varepsilon_0$ was defined in \eqref{eq-defve0}.
\begin{lemma}\label{lem-sepppoints}
For $n \leq b_1 \, k^2/4$ we have $s(\cGM^*, n, \varepsilon_0) \geq 2^k$. 
\end{lemma}    
\proof
Let $J = (j_1, \ldots , j_k)$ and $J' = (j_1', \ldots , j_k')$, and let $I = I' = (\ell_1, \ldots , \ell_k)$.
Suppose that $J \ne J'$ then let $1 \leq \nu \leq k$ be the greatest integer such that $j_{\nu} \ne j_{\nu}'$.  Set 
\begin{equation}\label{eq-partialstring}
 \ovp_{(I,J,\nu)} = \opsi^{\ell_k} \circ \ophi_{j_k} \circ \opsi^{\ell_{k-1}} \circ   \ophi_{j_{k-1}} \circ \cdots \circ   \opsi^{\ell_{\nu+1}} \circ \ophi_{j_{\nu+1}}\circ   \opsi^{\ell_{\nu}} 
\end{equation}
where $\|  \ovp_{(I,J,\nu)}\| \leq \| \vp_{(I,J)}\|$.   Assume that $\xi_{(I,J)} = \vp_{(I,J)}(\omega_0)$ and $\xi_{(I,J')} = \vp_{(I,J')}(\omega_0)$ are defined, then 
 \begin{eqnarray*}
(\ovp_{(I,J,\nu)})^{-1}(\xi_{(I,J)}) ~ & = & ~    \ophi_{j_{\nu}} \circ \opsi^{\ell_{\nu-1}} \circ \cdots \circ        \opsi^{\ell_1} \circ \ophi_{j_1}(\omega_0) \\
(\ovp_{(I,J,\nu)})^{-1}(\xi_{(I,J')}) ~ & = & ~   \ophi_{j_{\nu}'} \circ \opsi^{\ell_{\nu-1}} \circ \cdots \circ        \opsi^{\ell_1} \circ \ophi_{j_1'}(\omega_0) 
\end{eqnarray*}
As $j_{\nu} \ne j_{\nu}'$ we have 
$\ds d_{\fM}\left( (\ovp_{(I,J,\nu)})^{-1}(\xi_{(I,J)}) , (\ovp_{(I,J,\nu)})^{-1}(\xi_{(I,J')}) \right) \geq \varepsilon_0$. Thus the collection of points  $\xi_{(I,J)}$ constructed above corresponding to the initial choices of $\e$ and $b_1$ yields a collection of at least $2^k$ points which are $(n, \varepsilon_0)$-separated for  $n \leq b_1 \, k^2/4$.
 \endproof

Then by Lemma~\ref{lem-sepppoints} we have for all $k > 0$ and $n \leq b_1 \, k^2/4$ the estimate
$$\frac{\ln(s(\cGM^*,n,\varepsilon_0))}{\sqrt{n}} ~ \geq ~ \frac{\ln(2^k)}{\sqrt{b_1 \, k^2/4}} ~ = ~ \frac{2 \, \ln(2)}{\sqrt{b_1}}  $$
from which we conclude that $\ds  h_{GLW}^{1/2}(\cGM^*) > 0$,   completing the proof of Theorem~\ref{thm-slowentropy}. \hfill $\Box$

We conclude this section with two further results concerning the lamination entropy for Kuperberg flows. We only sketch the proofs, which follow the same approach as the proof of Theorem~\ref{thm-slowentropy} above.

  \begin{thm}\label{thm-fastgrowthentropy}
Let $\Phi_t$ be a generic Kuperberg flow. If the insertion maps $\sigma_j$  have ``fast growth''   in the sense of     Definition~\ref{def-fastgrowth}, 
then the number $s(\cGM^*, n, \ve_0)$ of $(\cGM^*, n, \ve_0)$-separated points for the $\cGM^*$ action on $\fC$ is asymptotically proportional to   $n$. 
 \end{thm}
 \proof
Suppose we are given a  finite set $\cS_n \subset \fC$ of $(\cGM^*, n, \ve_0)$--separated points, and the corresponding  subset $\cE_n = p_0^1(\cS_n) \subset \fMC$ which  is $(\wcGK^*, 2n,\varepsilon_0')$-separated, as in the proof of Theorem~\ref{thm-slowentropy}.

Let $\xi_1 \ne \xi_2 \in \cS_n$, and suppose that $\ovp \in \cGM^{(n)}$ satisfies $d_{\fM}(\ovp(\xi_1), \ovp(\xi_2)) \geq \e$. 
Then set  $\eta_i = p_0^1(\xi_i)$ and let $\vp \in \wcGK^*$ be the
lifted word which satisfies $d_{\bRt}(\vp(\eta_1), \vp(\eta_2)) \geq
\e'$. For the purposes of the estimates below, we can assume that $\vp
\in \cM(n)$, for $\cM(n)$ as in Definition~\ref{def-monotone}, and   $\eta_1 , \eta_2 \in Dom(\vp)$.

For an analysis of $\vp$ as in the proof of Theorem~\ref{thm-slowentropy}, we set the offset  distances    $\upsilon_k = 0$ for all $k \geq 2$, so the  lower bound estimates on the values $r_k$ of the $k-th$ image of  the special point $\omega_0$ 
are given by a recursive estimate using the estimates \eqref{eq-fastgrowthest} in Definition~\ref{def-fastgrowth}. 
Start with Lemma~\ref{lem-r=2verticalapproach} for some  point $p_0(j ; 1 , \ell) $ where $\ell \geq b_0$. Then apply the  fast growth estimate \eqref{eq-fastgrowthest} to obtain  a recursive estimate for $r_{k+1}$ in terms of $r_k$, which for  $k \geq 3$ yields 
\begin{equation}\label{eq-bounds5}
r_k - 2 ~ \geq ~   \lambda^{k-2} \, (r_2 -2) ~  \geq ~   \lambda^{k-2} \, \alpha_{\Phi}/(\mu_2 \, \ell_1)^2 .
  \end{equation}

Apply the estimate \eqref{eq-bounds5} to  $\vp$ which is a monotone word ending in $\phi^+_{j_m}$. 
Assume that    $r(\vp(\eta_i))   \geq 2 +\delta'$, then   $2 < r(\eta_i) \leq 2 + (\delta'/\lambda^m)$. Use that $0 < \delta' \leq 1$, then  the   power $\ell_1$ of the initial term $\psi$ of $\vp$ satisfies
 \begin{equation}\label{eq-bounds6}
\ell_1 ~ \geq ~     \sqrt{\lambda^{m-2} \, \alpha_{\Phi}/ \mu_2^2 \, \delta'}  ~ \geq ~     \sqrt{\lambda^{m-2} \, \alpha_{\Phi}/ \mu_2^2}
\end{equation}
so the power $\ell_1$ grows exponentially with $m$, at the rate approximately $\lambda^{m/2}$. Thus, to obtain $2^m$ words in $\cM(n)$ the length of the first segment $\psi^{\ell_1}$ must be approximately $\lambda^{m/2}$, and   we obtain the  asymptotic  estimate on the word length
$$n \sim  m + 1 + \sqrt{1 + \lambda + \lambda^2 + \cdots + \lambda^m} \sim \lambda^{m/2} .$$
That is, the word length required to obtain an $(n, \varepsilon_0)$--separated collection of points  grows exponentially if the set of points is assumed to grow exponentially. 
 \endproof

The work \cite[Section~8]{Kuperbergs1996} of Greg and Krystyna
Kuperberg introduces  piecewise linear (PL) versions of  her   flows, and studies
properties of their minimal sets. Our last result contributes another
insight to the dynamical properties of the piecewise-smooth flows. It
is based on a simple observation that if the Wilson flow is allowed to
have a discontinuity in its defining vector field $\cW$ along the periodic orbits, then   we can obtain  the   opposite conclusion to that of 
Theorem~\ref{thm-entropyGLW0}. The assumption that the Wilson flow is smooth forces the holonomy of $\Psi_t$ along the periodic orbits to be unipotent, 
and the generic Hypothesis~\ref{hyp-genericW} yields the estimates in Lemma~\ref{lem-r=2verticalapproach} which play a fundamental role above. However, for a piecewise-smooth flow $\cW$ the derivative of the transverse holonomy of $\Psi_t$ along the periodic orbits need not equal $1$, and in fact can be constructed so that the points $\omega_i$ are   hyperbolic attracting  for the map $\psi$.  It is then a long exercise in the methods of the this section to show:
   \begin{thm}\label{thm-hyperbolicentropy}
Let $\Phi_t$ be a   Kuperberg flow constructed from a piecewise-smooth Wilson flow $\Psi_t$ whose holonomy along the  periodic orbits is hyperbolic, then  
$h_{GLW}(\cGM^*) > 0$.  
 \end{thm}
 \proof
   The proof of Theorem~\ref{thm-entropyGLW0} constructs collections of $(\varepsilon_0', 2n)$-separated points for the action of the augmented  pseudogroup $\wcGK^*$ on $\fMC$.  If the map $\psi \in \cGK^*$ has hyperbolic attracting points $\omega_i$ for $i=1,2$, then the estimate \eqref{eq-spacing1aa} becomes exponential, which implies that the number of insertions $\phi_j^+$ that can be realized grows exponentially fast with the length of the initial word $\psi^{\ell_1}$.  
\endproof  

The discussion and results of this section suggest two problems to consider.
\begin{quest} \label{quest-entropyflow1}
 Suppose that $\Phi_t$ is a Kuperberg flow with a hyperbolic singularity at the special orbits, as discussed above, so that $h_{GLW}(\cGM^*) > 0$. Is it also true that   $h_{top}(\Phi_t) > 0$ for these flows?
 \end{quest}
 It seems likely that a careful consideration of the methods of this section will provide an affirmative answer to the above question. The second problem is in regards to the calculations used in the proof of Theorem~\ref{thm-slowentropy}.
 
\begin{quest} \label{quest-entropyflow2}
 Let $\Phi_t$ be a generic Kuperberg flow. If the insertion maps $\sigma_j$  have ``slow growth''   in the sense of  Definition~\ref{def-slowgrowth}, does the flow $\Phi_t$ have slow entropy $h_{top}^{1/2}(\Phi_t) > 0$? \end{quest}
 
  \bigskip

\section{Growth of leaves}\label{sec-growth}

In this section, we study the growth type of the surface $\fM_0$ considered as a leaf in $\fM$. The growth type is a natural invariant of the flow $\Phi_t$ and we will show is closely related to the slow entropy of the $\psg$ $\cGM^*$ introduced in the previous section. This provides a result analogous to   Manning's Theorem   in \cite{Manning1979},  that the volume growth rate of the universal cover for a compact manifold $M$ with negative sectional curvature is related to the entropy of the geodesic flow for $M$.

The idea is to use the construction of the tree $\TP$ in Section~\ref{sec-normal},    the action of $\cGM^*$ on
$\fC_0$ and the geometry of $\fM_0$  as discussed previously   in Sections~\ref{sec-normal},
\ref{sec-geometry} and \ref{sec-entropylamination} to calculate the volume growth function of $\fM_0$ which is defined as follows.

 The smoothly embedded zippered lamination $\fM \subset \mK$   inherits a Riemannian metric from $\mK$, and we let $d_{\fM}$ denote the induced distance function on  the \emph{leaves} of $\fM$. The submanifold $\fM_0 \subset \mK$ with boundary is given this distance function, and we let 
$\ds B_{\omega_0}(s)  = \{ x \in \fM_0 \mid d_{\fM}(\omega_0 , x) \leq s\}$  
be the closed ball of radius $s$ about the basepoint $\omega_0 = (2, \pi, 0) = \cR' \cap \cT$. Let $\A(X)$ denote the Riemannian area of a Borel subset $X \subset \fM_0$. Then $\mathrm{Gr}(\fM_0, s) = \A(B_{\omega_0}(s))$ is called the \emph{growth function} of $\fM_0$.

Given functions $f_1, f_2 \colon [0,\infty) \to [0, \infty)$ say that $f_1 \lesssim f_2$ if there exists constants $A, B, C > 0$ such that for all $s \geq 0$, we have that $  f_2(s) ~ \leq ~ A \cdot f_1(B \cdot s) + C$.
 Say that    $f_1 \sim f_2$ if both $f_1 \lesssim f_2$ and $f_2 \lesssim f_1$ hold.    This defines   equivalence relation on functions, which defines their   \emph{growth type}.

The growth function $\mathrm{Gr}(\fM_0, s)$ for $\fM_0$ depends upon the choice of Riemannian metric on $\mK$ and basepoint $\omega_0 \in \fM_0$, however  the growth type $[\mathrm{Gr}(\fM_0, s)]$  is   independent of   these choices, as observed by Milnor  \cite{Milnor1968} for coverings of compact manifolds and Plante  \cite{Plante1975} for leaves of foliations. 

We say that $\fM_0$ has \emph{exponential growth type} if $\mathrm{Gr}(\fM_0, s) \sim \exp(s)$. Note that $\exp(\lambda \, s) \sim \exp(s)$ for any $\lambda > 0$, so there is only one growth class of ``exponential type''.
We say that $\fM_0$ has \emph{nonexponential growth type} if $\mathrm{Gr}(\fM_0, s)  \lesssim \exp(s)$ but $\exp(s) \not\lesssim \mathrm{Gr}(\fM_0, s)$. 
We also have the subclass of   nonexponential growth type,  where $\fM_0$ has  \emph{quasi-polynomial growth type} if there exists $d \geq 0$ such that $\mathrm{Gr}(\fM_0, s)  \lesssim s^d$.  
The growth type of a leaf of a foliation or lamination is an
entropy-type invariant of its dynamics, as discussed in
\cite{Hurder2014}. 

 Here is the main result of this section:
\begin{thm}\label{thm-volumegrowth}
Let $\Phi_t$ be a generic Kuperberg flow. If the insertion maps $\sigma_j$  have ``slow growth''   in the sense of  Definition~\ref{def-slowgrowth},  then the growth type  of $\fM_0$ is nonexponential, and satisfies
\begin{equation}\label{eq-growthestimates}
\exp(\sqrt{s})   \lesssim \mathrm{Gr}(\fM_0, s)   \lesssim \exp(s)  
\end{equation}
In particular, $\fM_0$ does  not have quasi-polynomial growth type.
\end{thm}
\proof
 The proof of this result occupies the rest of this section. The first step is to 
 elaborate on   the relation between  the geometry of the tree $\TP$ and action of $\cGM^*$. 

 Consider the monoid $\cM_\fM(\infty)$ of monotone words in $\cGM^*$
 as defined in Definition~\ref{defn-monotoneM}.  The \emph{Cayley graph} of $\cM_\fM(\infty)$, denoted by $|\cM_\fM|$,  is the graph with:
 \begin{enumerate}
\item vertices  given by the set $\{ \ovp(\omega_0) \mid \ovp \in \cM_\fM(\infty) \}$, and 
\item edges  given by the actions of the maps $\{     \opsi, \ophi_1, \ophi_2 \}$ on the vertices.
\end{enumerate}
To be more precise, for $i=1,2$, 
     there is an  edge $\langle \ophi_i , \omega_0 \rangle$ joining
     $\omega_0$ to   the vertex $\ophi_i(\omega_0)$. For $\ovp \in \cM_\fM(n)$
     suppose that $\opsi \circ \ovp \in \cM_\fM(n+1)$ then we have an edge
     $\langle  \opsi , \ovp(\omega_0) \rangle$ joining $\ovp(\omega_0)$ to $\opsi\circ
     \ovp(\omega_0)$. For $\ovp \in \cM_\fM(n)$ and $i=1,2$ suppose that
     $\ophi_i \circ \ovp \in \cM_\fM(n+1)$ then we have an edge $\langle
     \ophi_i , \ovp(\omega_0) \rangle$ joining $\ovp(\omega_0)$ to $\ophi_i \circ
     \ovp(\omega_0)$. 
     
     All edges of $|\cM_\fM|$ are assigned length $1$
     with the standard metric on each, so $|\cM_\fM|$ becomes a
     complete metric space. Moreover,  each vertex $\ovp(\omega_0)$ has valence
     equal to one plus the number of words in $\{\opsi \circ \ovp, \ophi_1 \circ \ovp, \ophi_2 \circ \ovp\}$ which are well-defined at $\omega_0$.
     
We   compare the geometry of the graph $|\cM_\fM|$ with that of  the tree $\TP$ introduced  in Section~\ref{sec-normal} and that of the manifold $\fM_0$, using the following standard notion:
\begin{defn}\label{def-qi}
A map $f \colon X \to Y$ between metric spaces $(X, d_X)$ and $(Y, d_Y)$ is a \emph{quasi-isometry} is there exists constants $C_f \geq 0$ and $\lambda_f \geq 1$ so that for all $x,x' \in X$ we have
$$\lambda_f^{-1}\cdot d_Y(f(x), f(x')) - C_f ~ \leq ~ d_X(x,x') ~ \leq ~ \lambda_f \cdot d_Y(f(x), f(x')) + C_f  ~.$$
Moreover, for all $y \in Y$ there exists $x \in X$ such that $d_Y(y, f(x)) \leq C_f$.
\end{defn}

The first comparison is a consequence of  our previous observations.
 \begin{prop} \label{prop-cayley}
 There is an  embedding $\wtPhi \colon |\cM_\fM| \to \TP$ which is a quasi-isometry.
\end{prop}

\proof
 Define $\wtPhi$ as follows.  The special vertex point  $\omega_0$ is
 sent to     basepoint $\omega_0  \in \TP$ which is pictured in
 Figure~\ref{fig:fM0cT} as the point in the  upper horizontal
 strip. The other vertices are mapped to the points of $\TP$ defined
 by the intersection with $\cT$ of the curve in $\bRt$ defined by the
 action of $\ovp$ on   $\omega_0$. For $\ovp'$ a generator of
 $\cGM^*$, the edge $\langle \ovp',
 \ovp(\omega_0)\rangle$ of $|\cM_\fM|$ is
 mapped by a constant speed curve to the corresponding branch of $\TP
 \subset \fM_0$ connecting the points $\ovp(\omega_0)$ and
 $\ovp'\circ\ovp(\omega_0)$ that belong to
 $\fC_0$.
We give a    uniform estimate on the lengths of these branches of $\TP$.

 \begin{lemma} \label{lem-lengths}
There exists $0 < L_1 < L_2$ such that for each $\xi, \xi' \in \fC_0$ which are related by the action of a generator $\{\opsi, \ophi_1 , \ophi_2\}$ of $\cGM$, 
the segment $[\xi, \xi']$ of the tree $\TP\subset \fM_0$ joining them has length satisfying $L_1 \leq L([\xi, \xi']) \leq L_2$.
\end{lemma}
\proof
For the case when $\xi' = \opsi(\xi)$,  it was observed that $4\pi \leq L[\xi, \xi'] \leq 6\pi$. Thus, we need consider the case when $\xi' = \ophi_k(\xi)$ for $k =1,2$.
Each point $\xi \in \fC_0$ is joined to the endpoints $\{p_0^1(\xi) ,
p_0^2(\xi)\}$ of the $\g_0$ or $\lambda_0$ corresponding to it. The
segment of $\cK$-orbit between $p_0^i(\xi)$ and $\phi_i^+(p_0^i(\xi))$, for
$i=1,2$, is contained in the concatenation of two $\cW$-arcs and thus
it is upper and lower bounded as a consequence of
Corollary~\ref{cor-segmentlengths}. Let $L>0$ be the upper bound and
$L' > 0$  be the lower bound. 

Taking $L_1=\min\{L, 4\pi\}$ and $L_2=\max\{L',6\pi\}$, the claim follows.
  \endproof
 
To complete the proof of Proposition~\ref{prop-cayley}, note that Lemma~\ref{lem-lengths} shows   the edges of $\TP$ have   lengths which are uniformly bounded above and below, which   implies that $\wtPhi$ is a quasi-isometry. 
 \endproof

    The proof of Theorem~\ref{thm-zippered} shows that every point of $\fM_0$ is a uniformly bounded distance from a point of $\TP$. Thus, $\fC_0$ is a \emph{net} in $\fM_0$, which implies:

 \begin{cor}\label{cor-subexponential}
For a generic Kuperberg Plug, the map $\whPhi \colon |\cM_\fM| \to \fM_0$ obtained from the composition of $\wtPhi$ with the inclusion $\TP\subset \fM_0$, is a quasi-isometry.  
That is,  the graph $\TP$ is a ``tree model'' for the space $\fM_0$.    
 \end{cor}

  Since, by Corollary~\ref{cor-subexponential}, $|\cM_\fM|$   is quasi-isometric to  $\fM_0$ we have reduced the study of the growth properties of the leaf
  $\fM_0$ to those of the monoid $\cM_\fM(\infty)$, so   by  Corollary~\ref{cor-asymptoticM} we have:        
\begin{prop}\label{prop-asymptoticM} 
Both $\TP$ and  $\fM_0$ have  subexponential growth rates.
\end{prop}

We now return to the proof of Theorem~\ref{thm-volumegrowth}. By Proposition~\ref{prop-asymptoticM} we  have $\ds \mathrm{Gr}(\fM_0, s)   \lesssim \exp(s)$. 

To establish the lower bound $\ds \exp(\sqrt{s})   \lesssim \mathrm{Gr}(\fM_0, s)$, 
we assume that   $\Phi_t$ is a generic Kuperberg flow whose insertion maps $\sigma_j$  have ``slow growth''. 
Then the proof of Theorem~\ref{thm-slowentropy} constructs a subset of words in $\cM_\fM(\infty)$ whose action on the basepoint $\omega_0 \in \TP$ yields a set of images which grow at the rate $\exp(\sqrt{s})$. In particular, the number of words in $\cM_\fM(\infty)$ must grow at least this rate, which establishes  \eqref{eq-growthestimates}.  
 \endproof

   \bigskip
 
\section{Shape of the minimal set}\label{sec-shape}

In this   section, we consider the topological properties of the minimal set $\Sigma$ for a generic Kuperberg flow  $\Phi_t$.
The space $\Sigma $ is compact and connected, so is a continuum, but  its definition in terms of the closure of orbits reveals little about its topological nature. The natural framework for the  study of topological properties of spaces such as $\Sigma$ is using   shape theory. For example,   Krystyna Kuperberg posed the question whether $\Sigma $ has stable shape? Stable shape is discussed below, and  is about the nicest property one can expect  for a minimal set that is not a compact submanifold.     Theorem~\ref{thm-stableshape}  below    shows that $\Sigma $ does not have stable shape. This result follows  from the equality   $\Sigma = \fM$ for a generic flow, and the structure theory for $\fM$ developed in the previous sections of this work.

  Shape theory studies the topological properties of a topological space $\fZ$  using a form of \v{C}ech homotopy theory.  The definition of  {shape} for a space $\fZ$  embedded in the Hilbert cube was introduced by Borsuk \cite{Borsuk1968,Borsuk1975}. Later   developments and   results of shape theory are  discussed  in the texts \cite{DydakSegal1978,MardesicSegal1982} and the historical essay  \cite{Mardesic1999}. 
 We give a brief definition of the shape   of a compactum $\fZ$ embedded in a metric space $X$, following the works of Fox~\cite{Fox1972}, Morita~\cite{Morita1975} Marde{\v{s}}i{\'c}~\cite{MardesicSegal1982}, and the suggestions of the referee. 
  \begin{defn}\label{def-shapeapprox}
A sequence  $\fU = \{U_{\ell} \mid \ell =1,2,\ldots\}$ is called a \emph{shape approximation}   of $\fZ \subset X$ if: 
\begin{enumerate}
\item each $U_{\ell}$ is an open neighborhood of $\fZ$ in $X$ which is homotopy equivalent to a compact polyhedron;
\item $U_{\ell +1} \subset U_{\ell}$ for $\ell \geq 1$, and their closures satisfy $\ds \bigcap_{\ell \geq 1} ~ \oU_{\ell} = \fZ$.
\end{enumerate}
 \end{defn} 
Suppose that $X, X'$ are  connected manifolds, that  $\fU$ is a shape approximation for the compact subset $\fZ \subset X$, and $\fU'$ is a  shape approximation for the compact subset $\fZ' \subset X'$.
The compacta $\fZ, \fZ'$ are said to have the same shape if the following conditions are satisfied:
\begin{enumerate}
\item There are an order-preserving map $\phi \colon \mZ \to \mZ$, and for each $n \geq 1$, a continuous map $f_n \colon U_{\phi(n)} \to U_n'$ such that for any pair $n \leq m$, the restriction $f_n|U_{\phi(m)}$ is homotopic to $f_m$ as maps from $U_{\phi(m)}$ to $U_n'$.\medskip
\item There are an order-preserving map $\psi \colon \mZ \to \mZ$, and for each $n \geq 1$, a continuous map $g_n \colon U_{\psi(n)}' \to U_n$ such that for any pair $n \leq m$, the restriction $g_n|U_{\psi(m)}'$ is homotopic to $g_m$ as maps from $U_{\psi(m)}'$ to $U_n$.\medskip
\item  For each $n \geq 1$, there exists $m \geq \max\{n, \phi \circ \psi(n)\}$ such that the restriction of $g_n \circ f_{\psi(n)}$ to $U_m$ is homotopic to the inclusion as maps from $U_m$ to $U_n$. \medskip
\item  For each $n \geq 1$, there exists $m \geq \max\{n, \psi \circ \phi(n)\}$ such that the restriction of $f_n \circ g_{\phi(n)}$ to $U_m'$ is homotopic to the inclusion as maps from $U_m'$ to $U_n'$. \medskip 
\end{enumerate}

\begin{defn}
Let $\fZ \subset X$ be a compact subset of a connected manifold $X$. 
 Then the \emph{shape of $\fZ$ }is defined to be the equivalence class of a shape approximation of $\fZ$ as above. 
\end{defn}   
It is a basic fact of shape theory   that this definition does not depend upon the choice of shape approximations, and that two homotopic compacta have the same shape.
  The references  \cite{DydakSegal1978,MardesicSegal1982}  give complete details and alternate approaches to defining the shape of a space.  
A   concise overview of some of the key aspects of shape theory for continua embedded in Riemannian manifolds        is   given  in \cite[Section 2]{ClarkHunton2012}.

For the purposes of this work, we consider the case where $X$ is a connected compact Riemannian manifold, and $\fZ \subset X$  is an embedded continuum with the induced metric from that on $X$.
The shape of $\fZ$ can then be defined using    a  shape approximation $\fU$     defined by  a descending chain of open $\e$-neighborhoods of $\fZ$ in $X$, given by 
$\ds U_{\ell} = \{ x \in X \mid  d_X(x, \fZ) < \e_{\ell}\}$  where we have  $0 <  \e_{\ell +1} < \e_{\ell}$ for all $\ell \geq 1$, and  $\ds \lim_{\ell \to \infty}  \, \e_{\ell}   =   0$.

  \begin{defn}\label{def-stableshape}
 A compactum $\fZ$ has \emph{stable shape} if it is \emph{shape equivalent} to   a finite polyhedron. That is, there exists a shape approximation $\fU$    such that each inclusion $\iota \colon U_{\ell +1}\hookrightarrow  U_{\ell}$ induces a homotopy equivalence,  and $U_1$ has the homotopy type of a finite polyhedron. 
 \end{defn} 
 Some examples of spaces with stable shape are  compact   manifolds,  and more generally  finite $CW$-complexes. A less obvious
 example is the minimal set for a Denjoy flow on $\mT^2$ whose shape is equivalent to the wedge of two circles.  In particular, 
 the minimal set  of an aperiodic $C^1$-flow on plugs as constructed by Schweitzer in \cite{Schweitzer1974} has stable shape. 
A result of Krasinkiewicz shows that a continuum embedded in a closed orientable surface is either shape equivalent to a finite wedge of circles,  or    has the shape of a ``Hawaiian earring'' \cite{Krasinkiewicz1981,McMillan1975}. 
In higher dimensions,  Clark and Hunton show in \cite{ClarkHunton2012} that for 
 an $n$-dimensional lamination $\fZ$ embedded in an $(n+1)$-dimensional manifold $M$ such that  $\fZ$ is an   attractor  for a smooth diffeomorphism $f \colon M \to M$, and for which the restriction $f \colon \fZ \to \fZ$ is an    expanding map on the leaves of $\fZ$,   then $\fZ$  has stable shape.  
 In contrast, the results of the previous sections are used to show that the shape properties of the minimal set for a generic Kuperberg flow are not so simple. The first result   is the following. 

    \begin{thm} \label{thm-stableshape}
 The  minimal set $\Sigma$ of a generic Kuperberg flow does  not have stable shape.
 \end{thm}

We begin the proof of this result after some further discussions of the shape properties of $\fM$.

  \begin{defn} \label{def-movable}
 A compactum $\fZ \subset  X$ is said to be \emph{movable  in $X$} if for every neighborhood $U$
 of $\fZ$,  there exists a neighborhood $U_0 \subset U$ of $\fZ$ such
 that, for every neighborhood $W\subset U_0$ of $\fZ$, there is a continuous map $\varphi \colon U_0 \times [0,1] \to U$  satisfying the condition $\varphi(x,0)=x$ and $\varphi(x,1) \in W$ for every point $ x \in U_0$.
\end{defn}
The notion of a movable compactum was introduced by Borsuk \cite{Borsuk1969}, as a generalization of spaces having the shape of an  \emph{absolute neighborhood retract} (ANR's), and is an invariant of the shape  of $\fZ$. See    \cite{ClarkHunton2012,DydakSegal1978,MardesicSegal1982} for further discussions concerning movability.  
It is a standard result  that a compactum $\fZ$ with stable shape is movable. 
The movable property distinguishes between the shape of a Hawaiian earring and a Vietoris solenoid; the former is movable and the latter is not.
It  is a more subtle problem to construct compacta  which are invariant sets for dynamical systems, which are movable but do not have stable shape, such as given in   \cite{Sindelarova2007}.  

Showing the movable property for a space  requires the construction of a homotopy retract  $\varphi$ with the properties stated in the definition, whose existence  can be difficult to achieve in practice. There is an alternate condition on homology groups, weaker than  the movable condition,  which is much easier to check.

 \begin{prop}\label{prop-MLmove}
Let $\fZ$ be a movable compacta with shape approximation $\fU$. Then the    homology groups satisfy   the \emph{Mittag-Leffler Condition}:
 For all $\ell \geq 1$, there exists $p \geq \ell$ such that for any $q \geq p$, the maps on   homology groups for $m \geq 1$ induced by the inclusion maps  satisfy
 \begin{equation}\label{eq-shapehomology}
\text{Image}\left\{ H_m(U_p; \mZ) \to H_m(U_{\ell} ; \mZ) \right\}  = \text{Image}\left\{ H_m(U_q; \mZ) \to H_m(U_{\ell} ; \mZ) \right\}  ~ .
\end{equation}
\end{prop}

 This result is a   special case of a more general   Mittag-Leffler condition, as discussed in detail in \cite{ClarkHunton2012}.  For example,   the above   form of the Mittag-Leffler condition can be used to show that the Vietoris solenoid formed from the inverse limit of coverings of the circle is not movable.
 
We can now state   an additional shape property for the minimal set of a generic Kuperberg flow. 
\begin{thm}\label{thm-MLhomology}
Let $\Sigma$ be the minimal set for a generic Kuperberg flow. Then the Mittag-Leffler condition   for homology groups is satisfied. That is, given a shape approximation $\fU = \{U_{\ell}\}$ for $\Sigma$,
then  for any $\ell\geq 1$ there exists $p>\ell$ such that for any $q\geq p$ 
\begin{equation}\label{eq-MLhomology}
Image\{H_1(U_p;\mZ)\to H_1(U_\ell;\mZ)\} = Image \{H_1(U_q;\mZ)\to H_1(U_\ell;\mZ)\}.
\end{equation}
\end{thm}

The proof of Theorem~\ref{thm-MLhomology} follows by exhibiting a generating set of homology classes in $H_1(U_p;\mZ)$ which are represented by   closed loops in $U_p$, and it is shown that the images of these loops  in the   space $U_{\ell}$  become     homologous  to  classes in a fixed $3$-dimensional subspace $G_{\ell} \subset H_1(U_\ell;\mZ)$,  for $p \gg k$ sufficiently large. Moreover, the space $G_{\ell}$ is the image of the  map on homology induced by the inclusion map.
The proof of Theorem~\ref{thm-stableshape}  then follows from    the following   result.

\begin{prop}\label{prop-notstable}
 Let $\fU = \{U_{\ell} \mid \ell =1,2,\ldots\}$ be a  {shape approximation}   of $\fZ \subset X$, such that for $k > 0$
\begin{itemize}
\item the rank of $H_1(U_k;\mZ)\geq 2^k$,
\item for all $\ell > k$ the rank of the image $H_1(U_{\ell} ; \mZ) \to H_1(U_k ; \mZ)$ is 3.
\end{itemize} 
Assume that for any shape approximation of $\fV = \{V_{\ell} \mid \ell =1,2,\ldots\}$  the rank of the homology groups $H_1(V_{\ell} ; \mZ)$ is strictly greater than 3, then $\fZ$ does not have stable shape. 
 \end{prop}
 
\proof
Suppose that $\fZ$ has stable shape, so that there exists a shape approximation $\fV = \{V_{\ell} \mid \ell =1,2,\ldots\}$    such that each inclusion $\iota \colon V_{\ell +1}\hookrightarrow  V_{\ell}$ induces a homotopy equivalence, and thus the inclusion $V_{\ell} \hookrightarrow V_k$ is a homotopy equivalence for all $\ell\geq k$. Let   $3<n_0 = {\rm rank}(H_1(V_1 ; \mZ))$, then $n_0 = {\rm rank}(H_1(V_{\ell} ; \mZ))$ for all $\ell > 0$. Also, $n_0$ is the rank of the image of the map $ H_1(V_{\ell} ; \mZ) \to H_1(V_k ; \mZ)$ for all $\ell\geq k$.

 Then, as both $\fU$ and $\fV$ are shape approximations of $\fZ$, there exists $k_0, k, \ell_0$ and $\ell$ such that $V_{k_0} \subset U_k\subset V_{\ell_0}\subset U_\ell$. Thus we obtain the sequence of maps on homology induced by the inclusions:
\begin{equation}\label{eq-chain}
H_1(V_{k_0} ; \mZ) \to H_1(U_{k} ; \mZ) \to H_1(V_{\ell_0} ; \mZ) \to H_1(U_{\ell} ; \mZ).
\end{equation}
Then   the rank of the image of $H_1(V_{k_0} ; \mZ) \to H_1(U_k ; \mZ)$ must be equal to $n_0$, and the same holds for the map  $H_1(V_{\ell_0} ; \mZ) \to H_1(U_{\ell} ; \mZ)$. Hence the rank of  the image of $H_1(U_{k} ; \mZ) \to H_1(U_{\ell} ; \mZ)$ is equal to $n_0>3$  a contradiction. 
\endproof

 The strategy for the proof of Theorem~\ref{thm-stableshape} is   to construct a   shape approximation  $\fU$ for $\Sigma$ so that the conditions of Proposition~\ref{prop-notstable} are satisfied. This uses the properties of the level function   as developed in  Sections~\ref{sec-fM}, \ref{sec-proplevels} and \ref{sec-doublepropellers}. We then show  in Proposition~\ref{prop-notbound3}  that there is no shape approximation such that the rank of the homology groups is bounded above by 3.

 Theorem~\ref{thm-density}  implies that $\Sigma = \fM$ for a generic flow, and thus we analyze the shape properties of   $\fM$.  
 As $\fM$  is the closure of  $\fM_0$, we have 
 $U_{\mK}(\fM, \e) \subset U_{\mK}(\fM_0, \e')$ for all $0 < \e < \e'$. Thus, it   suffices to consider   shape approximations of $\fM_0$.

The space  $\fM_0$ retracts to the embedded tree    $\TP \subset \fM_0$  defined in Section~\ref{sec-normal} and illustrated in Figure~\ref{fig:fM0cT}. In the following, the    shape approximation $\{U_k \mid k > 0\}$ for   $\fM_0$ is given by  sets $U_k$ for $k \geq 0$, where each $U_k$ is an open neighborhood   in $\mK$ of a compact set $\fN_k$ which is homotopy equivalent  to a $1$-dimensional complex formed by taking quotients of $\TP$.   Each  open  neighborhood $U_k$ of $\fM_0$  thus ``fuses together'' sets of points of the tree $\TP$, so that the shape approximations to $\fM_0$ are homotopy equivalent to a bouquet of circles formed from paths in $\TP$ that  travel out a branch of $\TP$ and then are ``closed  up'' via paths in $\fN_k$.  The systematic description of the   classes in first homology that arise  in this way invokes the labeling system for the double propellers developed in Section~\ref{sec-doublepropellers} and the corresponding labels for the vertices of $\TP$.  
  
  First, to fix   notation,  for  $X \subset \mK$ and $\e > 0$, set:
\begin{itemize}
\item $U_{\mK}(X, \e) = \{ x \in \mK \mid d_{\mK}(x, X) < \e\}$ is the \emph{open} $\e$-neighborhood of $X$ in $\mK$; 
\item $C_{\mK}(X, \e) = \{ x \in \mK \mid d_{\mK}(x, X) \leq \e\}$ is the \emph{closed} $\e$-neighborhood of $X$ in $\mK$.
\end{itemize}

We now begin the recursive construction of a sequence   of   compact  subsets  $\fN_{k} \subset \mK$ for $k \geq 0$, which will be constructed to satisfy the conditions:
  \begin{enumerate}
\item  $\fM \subset \fN_{k}$ and $\fN_{k +1} \subset \fN_{k}$ for all $k \geq 0$;
\item $\fN_{k}$ is homotopy equivalent to a finite wedge of circles;
\item For all $\e > 0$, there exists $k > 0$ such that $\fN_{k} \subset U_{\mK}(\fM_0 , \e)$. 
\end{enumerate}
It then follows that sufficiently small open neighborhoods of the sets in the collection $\{ \fN_{k}  \mid k \geq 0\}$ yield a shape approximation for $\fM$.
 Moreover, each $\fN_k$ is constructed so that the branches of  $\TP$ above level $k$ are collapsed in the set $\fN_k$ to branches at level $k$.  This critical property  is achieved by   introducing  the notion of  ``filled double propellers'' at level $k$, which are compact regions in $\mK$ which contain all the branches of $\TP$ with level at least $k$. This containment property   follows from the      nesting properties of double propellers, as described in Section~\ref{sec-doublepropellers}. Colloquially, these  filled propellers can be thought of as a collection of  ``gloves'' which envelop collections of  branches  of $\TP$ at   level greater than $k$, grouping them together  at the end of a branch at level $k$. The recursive construction of the sets $\fN_k$ is made precise in the following.
We first define $\fN_0$, $\fN_1$ and $\fN_2$,  and establish their properties in detail. For the sets $\fN_k$ with $k > 2$, we provide   fewer details,  as the proofs follow the same outline as for the cases of $\fN_k$ with $k \leq 2$, except as noted. Begin by setting: 
\begin{equation}\label{eq-N0}
\fN_0 ~ = ~    \left\{ x \in \mK \mid r(x) \geq 2 \right\} ~  \subset ~ \mK .
\end{equation}

 \begin{lemma}\label{lem-shape0-1} 
 $\fN_0$ is compact.
 \end{lemma}
\proof
The discontinuities of the radius function $r \colon \mK \to [1,3]$ are contained in the closed subset $\cT_{\cK} \subset \mK$ defined by \eqref{eq-transversal}. By the formulation  \eqref{eq-radius} of the Radius Inequality, the radius value increases at a point of discontinuity, hence $\fN_0$ is closed in $\mK$, and thus is compact.
\endproof

\begin{figure}[!htbp]
\centering
{\includegraphics[width=50mm]{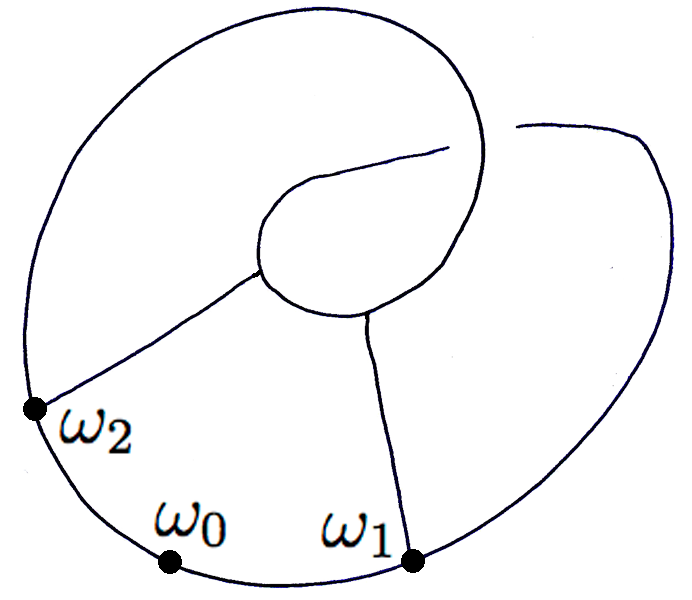}}
\caption{\label{fig:homotopyK} Homotopy type of $\mK$ and $\fN_0$}
\end{figure}

 \begin{lemma}\label{lem-shape0-2} 
 The inclusion $\fN_0 \subset \mK$ is a homotopy equivalence, and hence $\fN_0$      is homotopy equivalent  to a bouquet of three circles. 
 \end{lemma}
\proof
The inclusion $\{x \in \mW \mid r(x) \geq 2\} \subset \mW$ is a homotopy equivalence, with both spaces retracting to the   circle $\{ x=(2, \theta, 0) \mid 0 \leq \theta \leq 2\pi\}$. The identification map $\tau \colon \mW \to \mK$ creates a cross-arc between the circle and itself for each of the insertion maps $\sigma_i$,   as illustrated in Figure~\ref{fig:K}.
Thus, both spaces are homotopy equivalent  to the curve in  Figure~\ref{fig:homotopyK},  and so are homotopy equivalent  to a bouquet of three
circles. 
\endproof

 We fix notation for the generators of   $H_1(\fN_0; \mZ)$ as identified above. 
Let $[R]$ be the class defined by the loop starting at $\omega_0$ and following the Reeb cylinder around in a counterclockwise direction back to $\omega_0$. For $i=1,2$, let $[b_i]$ be the class defined by the straight path from $\omega_0$ to $\omega_i$, then    following the cross-arc created by the insertion $\sigma_i$ and then returning to $\omega_0$.   Then $H_1(\fN_0; \mZ)$ is the group  generated by  $\{ [R], [b_1] , [b_2]\}$. Note that the choice of the generators $[b_1]$ and $[b_2]$ have an ambiguity in terms of how the end of the cross-arc is joined to the basepoint $\omega_0$ to close the path.  Differing choices of   closing paths   result in the addition of an integer multiple of   $[R]$. We make a choice of these connecting paths, and consider the choice fixed in the following.

  \bigskip

Next, we consider the   construction of $\fN_1$ which  introduces the notion of  \emph{filled double propellers}, and also reveals the   mechanisms for the change  in the topology of the spaces $\fN_k$ as $k$ increases.

Recall  the  definitions \eqref{eq-Gamma} and \eqref{eq-Lambda} from Section~\ref{sec-doublepropellers},   
  $$
\G' = \cC \cap \cL_1^- \subset \mW ~ , ~ \G = \sigma_1^{-1}(\G') \subset  L_1^-
\quad , \quad
\Lambda' = \cC \cap \cL_2^- \subset \mW ~ , ~ \Lambda = \sigma_2^{-1}(\Lambda') \subset  L_2^- .
$$
The    \emph{filled double propellers at level $1$} are defined as follows.  Let    $\ds L^-(1)
\subset L_1^-$  denote the closure of the interior   region bounded by $\G$,
and  $\ds L^-(2)  \subset L_2^-$  the closure of interior   region bounded
by  $\Lambda$. Analogously, define $L^+(i)$ in the exit region $L_i^+$
for $i=1,2$.
 The Wilson flow of the curves $\G$ and $\Lambda$ generate  the infinite double propellers    $P_{\G} , P_{\Lambda} \subset \mW$.
 Let $\ds D(1)  \subset \mW$ denote the forward $\Psi_t$-flow of $L^-(1)$.
  We call $D(1)$ the filled double propeller associated to $\G$.   
  Analogously, let    $\ds D(2)  \subset \mW$ denote the forward $\Psi_t$-flow of $L^-(2)$, then $D(2)$ is the filled double propeller associated to $\Lambda$. Note that $D(1)$ and $D(2)$ are disjoint subsets of $\mW$.
  
   Consider the notched Reeb cylinder $\cR'$, as illustrated in
Figure~\ref{fig:notches}.     The filled double propellers in $\mW$ do not intersect
$\cR'$, but they do intersect every neighborhood of it, as the closures of both $D(1)$ and $D(2)$   contain the cylinder $\cR$, as discussed in the proof of Proposition~\ref{prop-propellerclosure}.   Also, recall that  $\tau(\cR') \subset \mK$ is the embedded Reeb cylinder,   illustrated in   Figure~\ref{fig:notched8}.

   For $\delta > 0$, define the   subset of $\mK$, 
  \begin{equation}
 C_{\mK}^+(\tau(\cR'), \delta)  ~ = ~  C_{\mK}(\tau(\cR'), \delta)    \cap    \fN_0 ~ = ~ \{ x \in \mK \mid d_{\mK}(x, \tau(\cR')) \leq \delta\} \cap \fN_0 ~ .
\end{equation}
Then  $\ds  C_{\mK}^+(\tau(\cR'), \delta) $ is a closed $\delta$-neighborhood   of $\tau(\cR')$, contained in the compact subset $\fN_0  \subset \mK$, hence is compact. 
Choose  $\delta_1$ small enough so that the set  $C_{\mK}^+(\tau(\cR'), \delta_1)$ retracts to $\tau(\cR')$.
In terms of the illustration  Figure~\ref{fig:notched8}, we choose $\delta_1$ less than the distance between the insertions and the edges of the gaps, so that 
$C_{\mK}^+(\tau(\cR'), \delta_1)$ has no ``self-intersections''.

For $i=1,2$,  set   $\ds \cD(i)  = \tau(D(i) \cap \wmW) \subset \mK$, where $\wmW$ is  the closure of $\mW'$ as defined in Section~\ref{sec-radius}.  
   Then  define:
\begin{equation}\label{eq-N1}
\fN_1 ~ = ~     C_{\mK}^+(\tau(\cR'), \delta_1)  ~ \cup ~   \cD(1) ~ \cup ~ \cD(2)  ~ .
\end{equation}
Thus,  $\fN_1$   is obtained by  attaching two infinite filled double
propellers to  $C_{\mK}^+(\tau(\cR'), \delta_1)$.

\begin{lemma}\label{lem-shape1-1} 
 $\fN_1$ is compact.
 \end{lemma}
\proof
Consider $\tau^{-1}( C_{\mK}^+(\tau(\cR'), \delta_1))$ which is a thickened  cylinder. For each $i=1,2$, the   region $D(i)$ is a connected solid spiral turning around $\cC$ an
infinite number of times, where after a finite number of turns, with  the number of times   determined by $\delta_1$,   the end of the solid region is contained in  $\tau^{-1}(C_{\mK}^+(\tau(\cR'), \delta_1))$. Thus, $D(i)-(D(i)\cap \tau^{-1}(C_{\mK}^+(\tau(\cR'), \delta_1)))$ has compact closure in $\mW$.  It   follows that   $\fN_1$   is   the union   of the compact set
 $C_{\mK}^+(\tau(\cR'), \delta_1)$  with  the images under $\tau$ of compact subsets of $D(1)$  and $D(2)$,   hence  is compact. 
 \endproof

\begin{lemma}\label{lem-shape1-2}  
$\fM \subset \fN_1 \subset \fN_0$.
\end{lemma}
\proof
It was shown in Section~\ref{sec-proplevels} that $\fM_0$ is the ascending union of the sets $\fM_0^n$ for $n \geq 0$, where $\fM_0^n$ consists of propellers at level at most $n$.  By construction, $\fN_1$ contains the set  $\fM_0^1$. It was   observed in Section~\ref{sec-doublepropellers} that 
   all of the $\G$ and $\Lambda$ curves at level at least $2$ in the face $\partial_h^-\mW$  are contained in the interiors $L^-(i)$ of the parabolic arcs $\G, \Lambda \subset \partial_h^- \mW$. Thus, the double propellers with level at least $2$ are contained in the filled propellers  $D(1)$  and $D(2)$,  hence $\fM_0 \subset \fN_1$. As $\fN_1$ is closed, it also contains $\fM$.
\endproof

\begin{lemma}\label{lem-shape1-3}  
For $\epsilon > 0$ sufficiently small,   $U(\fN_1 , \e)$ retracts to $\fN_1$.
\end{lemma}
\proof
The proof of Lemma~\ref{lem-shape1-1} shows that $\fN_1$ is the union of three compact submanifolds with boundary with corners. Hence,  $\fN_1$ is a submanifold with boundary with corners, from which the conclusion follows.
\endproof

\begin{lemma}\label{lem-shape1-4}  
 $\fN_1$   has the homotopy type of  a finite wedge of circles.
\end{lemma}
\proof
For each $i=1,2$,   attaching $D(i)$ to $\tau^{-1}(C_{\mK}^+(\tau(\cR'), \delta_1))$ is homotopy equivalent to attaching one endpoint of a  line segment to the core circle of the cylinder $\cR$.
 For the image   $\cD(i) \subset \mK$   the other endpoint of the segment is attached to $C_{\mK}^+(\tau(\cR'), \delta_1)$, resulting in a space with the same homotopy type as $\fN_0$. However, there are additional ``handles'' formed in the image which are not as immediate to visualize. These result from the intersections of   the images  $\ds \cD(i)$ for $i=1,2$ with themselves. That is,  
while the filled double propellers $D(1), D(2)\subset\mW$ are  disjoint from
each other,   their images $\cD(1)$ and $\cD(2)$ in $\mK$ under the quotient map $\tau \colon \mW \to \mK$ are not disjoint, and this gives rise to   additional $1$-cycles in $\fN_1$. 
The  loops created by these self intersections   are  called \emph{exceptional cycles}. 
While the number of intersections in $\mK$ of the filled propellers with themselves   will generate an infinite number of homotopy classes of loops,   there are at most finitely many homotopy classes of exceptional cycles   outside of  $\ds  C_{\mK}^+(\tau(\cR'), \delta_1)$. It follows that  $\fN_1$  has the homotopy type of  a finite wedge of circles.
\endproof

Note that the inclusion $\fN_1 \subset \fN_0$ induces a map between homology groups
$$\iota_1\colon H_1(\fN_1; \mZ) \to H_1(\fN_0; \mZ) \cong \langle [R], [b_1], [b_2] \rangle .$$
We claim that this map is surjective. Since $\fN_1$ contains a neighborhood of $\tau(\cR')$, it contains a loop representing the class $[R]\in H_1(\fN_0; \mZ)$. Consider in $\fN_1$ a loop that starts at $\omega_0$ then goes through the face $E_i$ for $i=1,2$,  to the first propeller, and after a certain number of turns depending on the value of $\delta_1$, the propeller intersects $C_\mK^+(\tau(\cR'),\delta_1)$. Thus the loop can be closed by passing back to the Reeb cylinder. The image of this loop is ${m_i}[R]+[b_i]$, for some $m_i>0$. Hence $[b_1]$ and $[b_2]$ are in the image and $\iota_1$ is surjective.

We next describe in more detail the classes in $H_1(\fN_1;\mZ)$  generated by intersections and self-intersections of the regions $\cD(1)$ and $\cD(2)$. We use the labeling of the double propellers  from Sections~\ref{sec-proplevels} and \ref{sec-doublepropellers}  to index   the intersections of these filled propellers.  
Recall that $\cD(1)$ is the interior region of the infinite
double propeller $\tau(P'_\G)$, and the propeller $P_\G\subset \mW$
intersects $\cL_i^-$ along the curves $\G'(i,\ell)$ for $i=1,2$ and
$\ell\geq b$, with $|b|$ the number of internal notches in the propellers, as defined in
Section~\ref{sec-doublepropellers}. 
Analogous observations  apply for $\cD(2)$ and the infinite
double propeller $\tau(P'_\Lambda)$. 
   The curves $\G'$ and $\Lambda'$  at levels $1$ and $2$ 
are illustrated in   Figure~\ref{fig:level2tangencies}. 

Let $L^-(1;i,\ell)\subset L_i^-$ be the
compact region bounded by
$\G(i,\ell)=\sigma_i^{-1}(\G'(i,\ell))$. 
For the exit regions, we have the analogous compact regions $L^+(1;i,\ell)\subset  L_i^+$.  
For $\Lambda$-curves, we have the analogous  entry and exit   compact regions $L^\pm(2;i,\ell)\subset L_i^\pm$ .

Consider their images  under the map $\tau \colon \mW \to \mK$ which identifies $L_i^\pm$ and $\cL_i^\pm$ with the secondary entry and secondary exit regions. We now analyze  the intersections and self-intersections  that are created. 
For $i_0,i_1=1,2$, the regions $\tau(L^-(i_0;i_1,\ell))\subset \cD(i_0)$ are the components of  the set $\cD(i_0)\cap E_{i_1}$, so by     the nesting property of $\G$ and
$\Lambda$ curves, we have that:
\begin{itemize}
\item $\tau(L^\pm(1;1,\ell))\subset \tau(L^\pm(1))$ and thus is
  contained in $\cD(1)$;
\item $\tau(L^\pm(1;2,\ell))\subset \tau(L^\pm(2))$ and thus is
  contained in $\cD(2)$.
\end{itemize}
Analogously, we have:
\begin{itemize}
\item $\tau(L^\pm(2;1,\ell))\subset \tau(L^\pm(1))$ and thus is  contained in $\cD(1)$;
\item $\tau(L^\pm(2;2,\ell))\subset \tau(L^\pm(2))$ and thus is  contained in $\cD(2)$.
\end{itemize} 
Thus,  $\cD(1)\cap\cD(2)$ consists of the regions:
\begin{itemize}
\item $\tau(L^\pm(1;2,\ell))\subset E_2\cup S_2$, for   $\ell\geq b$ and unbounded; 
\item $\tau(L^\pm(2;1,\ell))\subset E_1\cup S_1$, for   $\ell\geq b$ and unbounded.
\end{itemize}

\begin{remark}\label{rmk-elldelta}
There exists an index $\ell(\delta_1)$, which  tends to infinity as   $\delta_1$ tends to zero, such that  $\tau(L^\pm(i_0;i_1,\ell))$ deforms into  $C_{\mK}^+(\tau(\cR'),\delta_1)$ in $\fN_1$ for $\ell \geq \ell(\delta_1)$. It follows that the number  of such intersections which are not deformable in $\fN_1$ to the core $C_{\mK}^+(\tau(\cR'),\delta_1)$  is finite,  and bounded above by  $|b|+\ell(\delta_1)$.
\end{remark}

We complete our description of the topology of $\fN_1$ by giving  
 a set of generators for $H_1(\fN_1;\mZ)$. There are three types of generators: the generator at level 0; the generators that cover the branching of the tree $\TP$ from level 0 to level 1; and the generators that make one turn around the Reeb cylinder along a level 1 propeller. The branches of $\TP$ at higher levels are contained in the filled propellers at level 1 and thus do not contribute to the homology of $\fN_1$. 
 
  \begin{defn}\label{def-generators1}
 Consider the following  generators of $H_1(\fN_1;\mZ)$:

{\rm 

\begin{enumerate}
\item The loop $[R]$ corresponding to the fundamental class of the Reeb cylinder.\\
\item The exceptional loops $E^1(i_0,i_1)$ at level 1, which   are formed as   follows. Consider the point $\tau(p(i_0))\in E_{i_0}$, which we recall  is the intersection of the image of the periodic orbit $\cO_{i_0}$ of the Wilson plug with the entry face of the corresponding insertion. Recall that $\tau(\gamma(i_1,b))$ is the first curve in the intersection of $\tau(P_\g)$ with $E_{i_1}$ and $\tau(\lambda(i_1,b))$ is the first curve in the intersection of $\tau(P_\lambda)$ with $E_{i_1}$. These curves are illustrated in Figures~\ref{fig:intcurvesL1} and \ref{fig:intcurvesL2}.

Connect $\tau(p(i_0))$ to a point $\tau(p(i_0;i_1,b))\in E_{i_1}$ by a path tangent to the propeller $\tau(P_\g)$ if $i_0=1$, and tangent to $\tau(P_\lambda)$ if $i_0=2$. Note that $\tau(p(i_0;i_1,b))$ was not defined before, but for homology purposes we just need a point in $\tau(\gamma(i_1,b))$ if $i_0=1$, or in $\tau(\lambda(i_1,b))$ if $i_0=2$. Then $\tau(p(i_0;i_1,b))\in \tau(L^-(i_1))\subset E_{i_1}$, thus can be connected inside $\fN_1$ to $\tau(p(i_1))$. If $i_1=i_0$ we obtain a loop, and otherwise  close the loop using  a path from $\tau(p(i_1))$ to $\tau(p(i_0))$ contained in $\tau(\cR')$ (the two possible choices differ by $[R]$).
These loops are  illustrated in Figure~\ref{fig:Figure54a}. Observe that $E^1(i_0,i_1)$ does not depend on the choice of $\delta_1$ and that there are $2^2=4$ such loops.\\

\begin{figure}[!htbp]
\centering
{\includegraphics[width=105mm]{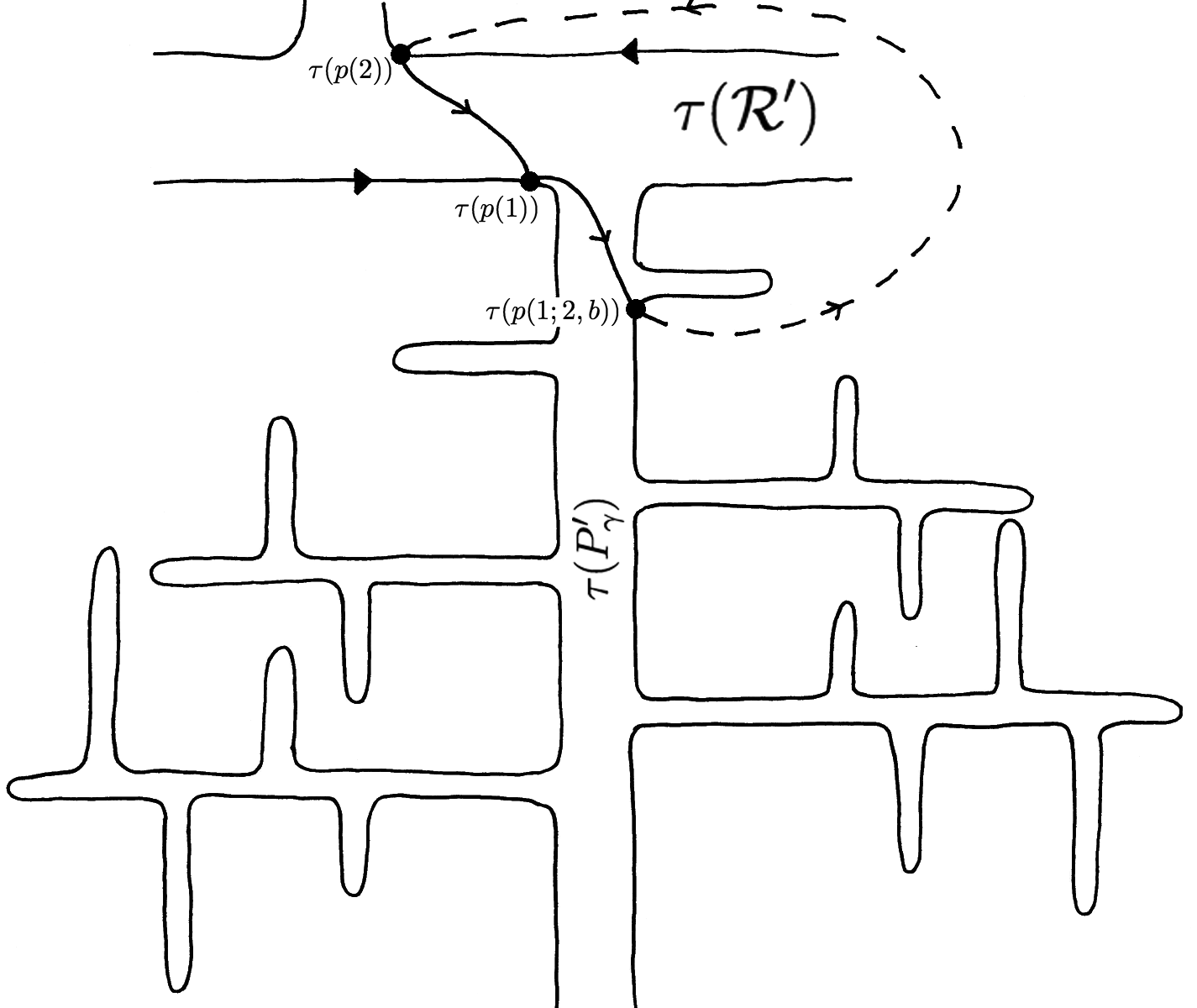}}
\caption{\label{fig:Figure54a}  The exceptional loop $E^1(i_0,i_1)$  at level $k=1$. }
\begin{minipage}[t]{14cm} 
Solid lines represent paths in the set $\fM_0^1$ and   dotted lines represent paths in the intersection of $\fN_1$ with the entry regions $E_1$ and $E_2$. Here,    $i_0 = 1$ and $i_1 = 2$. 
\end{minipage}
\end{figure}

\item The  loops $T^1(i_0,i_1,i_2;\ell)$ of type (2) at level 1, for $b\leq \ell < \ell(\delta_1)$. Consider a loop which starts by connecting $\tau(p(i_1))$ to $\tau(p(i_0;i_1,\ell))$ by a path contained in $E_{i_1}$, then goes from $\tau(p(i_0;i_1,\ell))$ to $\tau(p(i_0;i_2,\ell+1))\in E_{i_2}$ by a path tangent to the corresponding level one propeller (that is, to $\tau(P_\g)$ if $i_0=1$ and to $\tau(P_\lambda)$ if $i_0=2$). This segment makes a full turn around the Reeb cylinder and is contained in a level 1 propeller. Observe that $i_2$ is not necessarily equal to $i_1$. Then connect by a path in $E_{i_2}$ the points $\tau(p(i_0;i_2,\ell+1))$ and $\tau(p(i_2))$. If $i_2=i_1$ we obtain a loop; otherwise  close  the loop by a path from $\tau(p(i_2))$ to $\tau(p(i_1))$ contained in $\tau(\cR')$. These loops are  illustrated in Figure~\ref{fig:Figure55a}. 

Observe that for $\ell=\ell(\delta_1)-1$, the point $\tau(p(i_0;i_2,\ell+1))=\tau(p(i_0;i_2,\ell(\delta_1))$ belongs to $C_\mK^+(\tau(\cR'),\delta_1)$ and hence the loop can be closed by a path segment in $C_\mK^+(\tau(\cR'),\delta_1)$ but outside $\fM_0^1-\tau(\cR')$,  from the corresponding level one propeller to the Reeb cylinder.
Observe that there are $8(|b|+\ell(\delta_1))$ such loops $T^1(i_0,i_1,i_2;\ell)$ that are not homologically trivial or equivalent to $[R]$.
\end{enumerate}
 
}
\end{defn}

\begin{figure}[!htbp]
\centering
{\includegraphics[width=105mm]{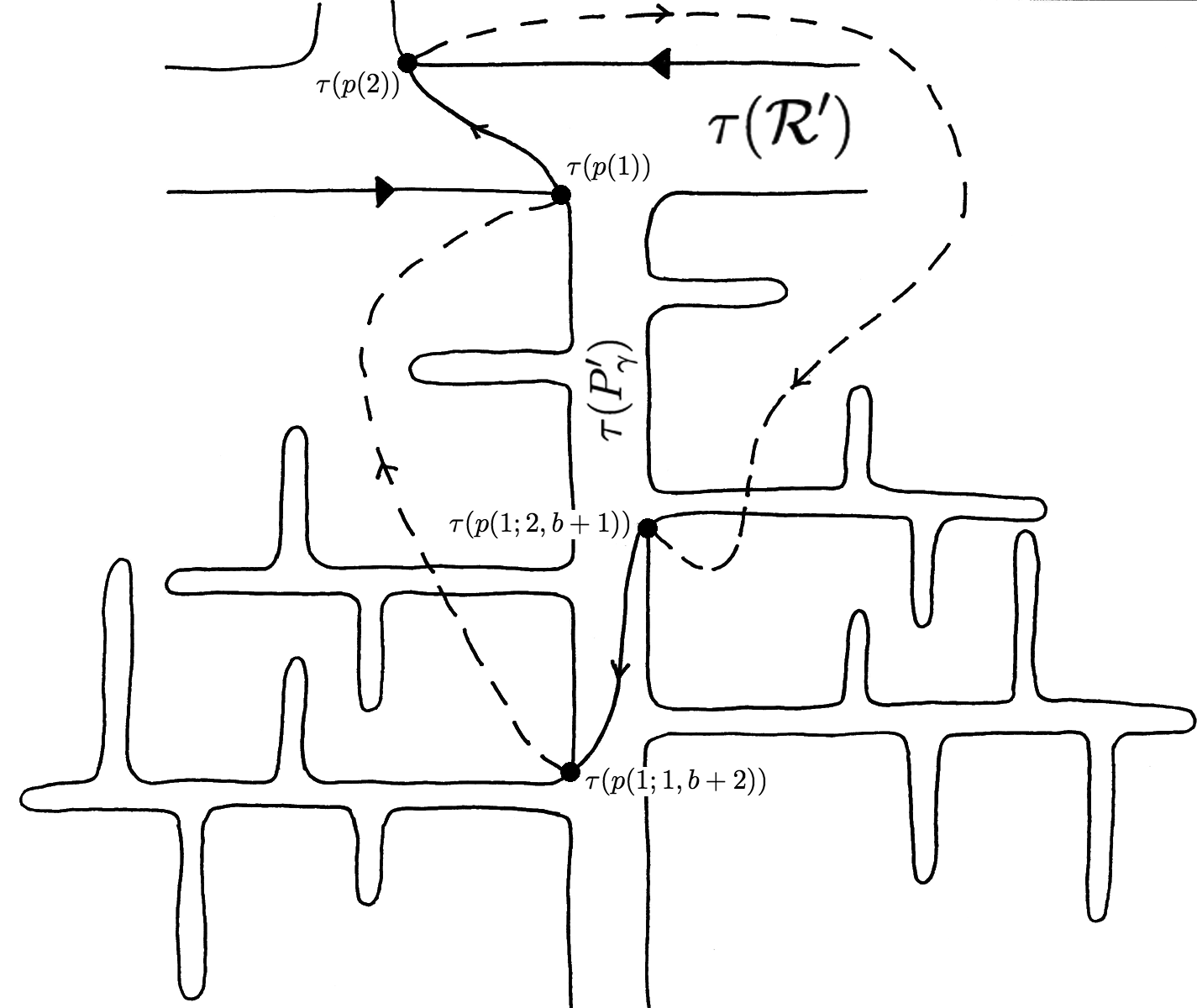}}
\caption{\label{fig:Figure55a} The type (2) loop $T^1(i_0,i_1,i_2;\ell)$   at level $k=1$. }
\begin{minipage}[t]{14cm} 
  Solid lines represent paths in the set $\fM_0^1$ and   dotted lines represent paths in the intersection of $\fN_1$ with the entry regions $E_1$ and $E_2$. Here,  $i_0 =1$, $i_1 = 2$, $i_2 = 1$ and $\ell_1 = b+1$.
\end{minipage}
\end{figure}

 The loops $E^1(i_0,i_1)$ and $T^1(i_0,i_1,i_2;\ell)$ constructed in items (2) and (3) above are based at $\tau(p(i_0))$, and they can be connected to the basepoint $\omega_0$ by a short path in $\tau(\cR')$ to obtain pointed loops.  We consider the image of these loops in $H_1(\fN_0;\mZ)$, and   express the homology classes they determine    in terms of the generators $\{[R], [b_1], [b_2]\}$.  As the resulting homology classes do not depend on the basepoint chosen, we do not mention the addition of the basepoint paths again. 
  
 We now specify the classes in $H_1(\fN_0 ; \mZ)$ determined by the above classes.  For $E^1(i_0,i_1)$,  the image is $|a-b| \cdot [R]+ [b_{i_0}]$,
  for $a$ as introduced in Section~\ref{sec-doublepropellers}. For $T^1(i_0,i_1,i_2;\ell)$ with $b\leq \ell \leq\ell(\delta_1)-1$,  the image is simply $[R]$ or $2[R]$ depending on the choice of closing path from $\tau(p(i_2))$ to $\tau(p(i_1))$. All other classes in $H_1(\fN_1, \mZ)$ are combinations of these classes, as can be seen from   geometric considerations. For example, consider the loop generated as in case (2), given by a path starting at $\tau(p(i_0))$ and going to a point $\tau(p(i_0;i_1,\ell_1))\in E_{i_1}$ for $\ell_1 > b$, by a path tangent to the propeller $\tau(P_\g)$ if $i_0=1$, , and tangent to $\tau(P_\lambda)$ if $i_0=2$.
This will be homologous to a class $E^1(i_0,i_1)$ plus a sum of classes $T^1(i_0,i_1,i_2;\ell)$ for $b < \ell \leq \ell_1$.
 
 The construction of the space $\fN_{k}$ for $k > 1$    is
 analogous to that of $\fN_1$, and follows an inductive procedure. Assume that the space $\fM_0^{k -1}$  has been constructed. Then $\fN_{k}$  is defined as  a union of a closed
 $1$-sided neighborhood of the   space $\fM_0^{k -1}$ with the  filled double propellers at level $k$.  Note that the filled double propellers at level $k-1$  in $\fN_{k-1}$ will no longer be contained in $\fN_{k}$, as they are replaced by filled double propellers at level $k$. 
  Though it is simple to state the recursive construction of $\fN_{k}$ in outline,   the actual construction becomes increasingly complicated to describe,  as it is based on the labeling for components of the filled double propellers given in Section~\ref{sec-doublepropellers}. We discuss this process in further detail for level $k=2$  below.

Recall that the double propellers at level $2$ in $\mW$ are labeled by their
generating curves $\G(i_1,\ell)$ and $\Lambda(i_1,\ell)$ in $L_{i_1}^-$ for $i_1
= 1,2$ and $\ell \geq b$. For each such curve, we define the filled
region $L^-(i_0; i_1 , \ell)$ and the corresponding filled double propeller $D(i_0; i_1 , \ell)$ at level $2$ is given by the  $\Psi_t$-flow of $L^-(i_0; i_1 , \ell)$.
Regarding this notation, if $i_0=1$, then the boundary of $D(i_0; i_1 , \ell)$ is $\G(i_1,\ell)$. If $i_0=2$, then the boundary of $D(i_0; i_1 , \ell)$ is
$\Lambda(i_1,\ell)$. As all $\G$ and $\Lambda$ curves which generate these double propellers are contained in the
region $\{r>2\}\subset \partial_h^-\mW$, each filled double propeller     $D(i_0; i_1 , \ell)$   is a compact region in $\mW$.

 Choose  $\delta_2\leq  \delta_1/2$ sufficiently small, so that the one-sided closed $\delta_2$-neighborhood of $\ds \fM_0^1$ given by 
\begin{equation}\label{eq-retract2}
C_{\mK}^+(\fM_0^1, \delta_2) =C_{\mK}(\fM_0^1, \delta_2) \cap \fN_1
\end{equation}
does not contains all the level 2 propellers. Set   $\cD(i_0; i_1 , \ell) = \tau(D(i_0; i_1 , \ell)  \cap \wmW)$ and define 
\begin{equation}\label{eq-N2}
\fN_2 ~ = ~     C_{\mK}^+(\fM_0^1, \delta_2)  ~ \cup ~  \left\{\bigcup \left\{ \cD(i_0; i_1 , \ell) \mid i_0,i_1 = 1,2 ~ {\rm and}  ~ \ell \geq b \right\}\right\} ~   .
\end{equation}
Thus,    $\fN_2$ is obtained by attaching to the  set $C_{\mK}^+(\tau(\cR'), \delta_2)$  a closed $1$-sided $\delta_2$-neighborhood of the level $1$ propellers to obtain $C_{\mK}^+(\fM_0^1,\delta_2)$, and
 then attaching the filled double propellers at level $2$.

 We  next state  versions  of Lemmas~\ref{lem-shape1-2} and \ref{lem-shape1-3}   for $\fN_2$ whose proofs follow in the same way.

 \begin{lemma}\label{lem-shape2-2} 
$\fM \subset \fN_2 \subset \fN_1$.
\end{lemma}

\begin{lemma}\label{lem-shape2-3}  
For $\epsilon > 0$ sufficiently small,   $U(\fN_2 , \e)$ retracts to $\fN_2$.
\end{lemma}

The description of the homology of $\fN_2$ is analogous   to the description of $\fN_1$ above, and again uses the labeling of the double propellers  from Sections~\ref{sec-proplevels} and \ref{sec-doublepropellers}.  However, in addition to the exceptional cycles in $\fN_2$ that arise from the intersections and self-intersections of the filled double propellers at level $2$, there is an added subtlety in the existence   cycles of type (2) at different levels, which occur for   the spaces $\fN_k$ when $k \geq 2$, since in these cases the added propellers are no longer infinite. 

We begin with some observations, before enumerating a set a generators of $H_1(\fN_2;\mZ)$.  For $\ell$ large enough,  $\cD(i_0;i_1,\ell)$ intersects the entry region $E_{i_2}$ for $i_2=1,2$. 
The intersection $\cD(i_0;i_1,\ell) \cap E_{i_2}$ is   along the regions bounded by  the curves $\tau(\G(i_1,\ell;i_2,\ell'))$ if $i_0=1$,  and $\tau(\Lambda(i_1,\ell;i_2,\ell'))$ if $i_0=2$, where    $\ell'\geq b$ and is bounded above. 
 As before, the regions  $D(i_0;i_1,\ell)$ are disjoint and contractible in
$\mW$, but their images $\cD(i_0;i_1,\ell)$ may have multiple intersections between them, and  also with the neighborhood   $C_\mK^+(\fM_0^1,\delta_2)$,
 so they are not necessarily contractible. 

First define $\ell(\delta_2)$ as in Remark~\ref{rmk-elldelta}, so that the curves $\tau(\g(i_1,\ell))$ and $\tau(\lambda(i_1,\ell))$ are contained in $C_\mK^+(\tau(\cR'),\delta_2)$ if $\ell\geq \ell(\delta_2)$.

 Let $t(i_0;i_1)\geq b$ be the greatest number such that $\cD(i_0;i_1,\ell)$ for $\ell=t(i_0;i_1)$ does not  intersect $E_1$, nor $E_2$. Observe that this phenomena is independent of the choice of $\delta_2$. 
 
We now consider the possible intersections of the filled double propellers at level 2 with $C_\mK^+(\fM_0^1,\delta_2)$.
The filled double propellers $\cD(i_0;i_1,\ell)$  get longer as
$\ell$ increases, thus there exist   constants 
\begin{equation}\label{eq-level2fork=2}
\ell(\delta_2;i_0;i_1) >s(\delta_2;i_0;i_1)\geq b
\end{equation}
such that, for $i_0,i_1=1,2$:
\begin{enumerate}
\item For $b\leq \ell<s(\delta_2;i_0;i_1)$, the filled double propellers
  $D(i_0;i_1,\ell)\subset \mW$ and the pre-image $\tau^{-1}(C_{\mK}^+(\fM_0^1,
  \delta_2))$ are disjoint. If $\ell<t(i_0;i_1)$ then $\cD(i_0;i_1,\ell)$ is homotopically
  trivial and the attachment of $\cD(i_0;i_1,\ell)$ to  $C_{\mK}^+(\fM_0^1,
  \delta_2)$ does not change the homotopy type. \\
\item For $s(\delta_2;i_0;i_1)\leq \ell <\ell(\delta_2;i_0;i_1)$,  the filled double
  propellers $D(i_0;i_1,\ell)\subset \mW$ intersect
  the set $\tau^{-1}(C_{\mK}^+(\fM_0^1, \delta_2))$ and $\cD(i_0;i_1,\ell)$ does not retract in $\fN_2$ to $C_{\mK}^+(\fM_0^1,
  \delta_2)$. In this case, 
$$\cD(i_0;i_1,\ell)\,-\, \left(\cD(i_0;i_1,\ell)\cap C_{\mK}^+(\fM_0^1,\delta_2)\right)$$
is a non-empty   submanifold with compact closure, and the finite filled double propeller
$\cD(i_0;i_1,\ell)$ is such that its base and tip are in $C_{\mK}^+(\fM_0^1,\delta_2)$.  The attachment of $\cD(i_0;i_1,\ell)$ to  $C_{\mK}^+(\fM_0^1, \delta_2)$ adds a handle.\\
\item For $\ell\geq \ell(\delta_2;i_0;i_1)$, the filled double propellers
  $\cD(i_0;i_1,\ell)$ retract in $\fN_2$ to $C_{\mK}^+(\fM_0^1,
  \delta_2)$, and thus adding these propellers does not change the
  homotopy type of $C_{\mK}^+(\fM_0^1,\delta_2)$.
\end{enumerate}

The proof of the following result is analogous to that of   Lemmas~\ref{lem-shape1-1} and \ref{lem-shape1-4}.
\begin{lemma}\label{lem-shape2-14}
$\fN_2$ is compact and has the homotopy type of a finite wedge of circles.
\end{lemma}

The trace of $C_{\mK}^+(\fM_0^1,\delta_2)$ on $\bRt$ contains the 1-sided
$\delta_2$-neighborhood of $\tau(\cR')\cap\bRt$. From \eqref{eq-retract2} we get that this trace contains also   the $\delta_2$-neighborhood  of  the curves
$\g_0(\ell_1)$, $\lambda_0(\ell_1)$ for $\ell_1\geq a$ and
unbounded. For $\ell_1 \geq b$, the trace of $\cD(i_0;i_1,\ell_1)$ is the region bounded by
the curves $\G_0(i_1,\ell_1;\ell_2)$ if $i_0=1$, and by the curves
$\Lambda_0(i_1,\ell_1;\ell_2)$ if $i_0=2$, for $\ell_2\geq a$ and
bounded. 

Recall from Section~\ref{sec-doublepropellers}, that for $i_1$ and
$\ell_2$ fixed, the curves $\g_0(i_1,\ell_1;\ell_2)$ accumulate on a level $1$ curve as $\ell_1\to\infty$,  on 
$\g_0(\ell_2)$ if $i_1=1$, and on $\lambda_0(\ell_2)$ if $i_1=2$. Analogously, the curves $\lambda_0(i_1,\ell_1;\ell_2)$ accumulate  as
$\ell_1\to\infty$,  on
$\g_0(\ell_2)$ if $i_1=1$ and on $\lambda_0(\ell_2)$ if $i_1=2$. Thus,  there exists $\ell'(\delta_2;i_0;i_1)$ such that
for every $\ell_1\geq \ell'(\delta_2;i_0;i_1)$, for $i_0=1$  the curve
$\G_0(i_1,\ell_1;a)$,  and for $i_0=2$  the curve $\Lambda_0(i_1,\ell_1;a)$, 
intersects the 1-sided $\delta_2$-neighborhood of
$\g_0(a)\cup\lambda_0(a)$.

The index  $\ell(\delta_2;i_0;i_1)$ may be required to be chosen larger 
than $\ell'(\delta_2;i_0;i_1)$ in order  that the conditions of case (3) are satisfied: that is,  the
filled double propeller $D(i_0;i_1,\ell_1)\subset \mW$ intersects
$\tau^{-1}(C_{\mK}^+(\fM_0^1,\delta_2))$ along all its length,  and thus $\cD(i_0;i_1,\ell_1)$ retracts in $\fN_2$ to  $C_{\mK}^+(\fM_0^1,\delta_2)$.

We  estimate the value of  $s(\delta_2;i_0;i_1)$ in \eqref{eq-level2fork=2} above. To simplify the discussion, consider the case   $i_0=1$. 
For $\ell_1<\ell'(\delta_2;1;i_1)\leq \ell(\delta_2;1;i_1)$
and $\delta_2$ small,
the curve $\G_0(i_1,\ell_1;a)$ is disjoint from the 1-sided
$\delta_2$-neighborhood of $\g_0(a)\cup\lambda_0(a)$. Take
$\ell_1=\ell'(\delta_2;i_0;i_1)-1$, then there exists $\ell_2>a$ such that
$\G_0(i_1,\ell_1;\ell_2)$ intersects the 1-sided
$\delta_2$-neighborhood of $\g_0(\ell_2)\cup\lambda_0(\ell_2)$. Thus
the filled double propeller $\cD(1;i_1,\ell_1)$ satisfies case (2).

As $\ell_1$ decreases, the filled double propellers $\cD(1;i_1,\ell_1)$
get shorter, and thus there exists $s'(\delta_2;1;i_1)$ such that if
$\ell_1<s'(\delta_2;1;i_1)$, the curves $\G_0(i_1,\ell_1;\ell_2)$ are
disjoint from $C_{\mK}^+(\fM_0^1,\delta_2)$ for any $\ell_2\geq a$. In the
same way we obtain $s'(\delta_2;2;i_1)$, and hence
$s(\delta_2;i_0;i_1)\leq s'(\delta_2;i_0;i_1)$ such that case (1) is satisfied.

The addition of the
propellers $\cD(i_0;i_1,\ell)$ in case (2) creates new handles: each propeller retracts to its intersection with the annulus
$\cA$, thus is an arc whose endpoints are in $C_{\mK}^+(\fM_0^1,\delta_2)$
that is not contained in the previous set. As for $\fN_1$, we obtain exceptional cycles in $\fN_2$ that arise from the intersections of the filled propellers at level 2. 
Thus,  the inclusion $\fN_2 \subset \fN_1$  is not a homotopy equivalence for $\delta_2$ sufficiently small. 
 
 Before giving the descriptions of the classes of generators for $H_1(\fN_2, \mZ)$, we make an observation.   Consider  point   $\tau(p(i_0;i_1,\ell)) \in C_\mK^+(\tau(\cR'),\delta_1)\cap E_{i_1}$ at  level 1,  for $i_1=1,2$,    where  $\ell$ is sufficiently large, but such that     $\tau(p(i_0;i_1,\ell))$ is not contained in $C_\mK^+(\tau(\cR'),\delta_2)$. The closest level 1 point to $\tau(p(i_0;i_1,\ell))$ is $\tau(p(i_2;i_1,\ell+1))$ for $i_2\neq i_0$, since $\g$ and $\lambda$ curves are interlaced as in \eqref{eq-interlaced1}, \eqref{eq-interlaced2}, \eqref{eq-interlaced12} and \eqref{eq-interlaced22}. If the $\delta_2$-neighborhood of $\tau(p(i_2;i_1,\ell+1))$  contains the first point $\tau(p(i_0;i_1,\ell))$, then we can compose these paths to obtain a loop as, as described  in item (3) of the list below.

We can now define a set of generators of $H_1(\fN_2;\mZ)$, which are divided into the following classes: the class $[R]$ at level 0; the ones that cover the branching from level 0 to level 1; the ones that cover the branching from level 1 to level 2; and those that allow to make one turn around the Reeb cylinder along a level 2 propeller. 
 
 \begin{defn}\label{def-generators2}
 Consider the following  generators of $H_1(\fN_2;\mZ)$:

{\rm 

 \begin{enumerate}
\item The loop $[R]$ corresponding to the fundamental class of the Reeb cylinder.\\
\item The   loops $P^2_1(i_0)$ at level 1, which   are formed as   follows. 
 Observe that the level 1 propellers limit to the Reeb cylinder, so that after  $|a|+\ell(\delta_2)$ turns, where $\ell(\delta_2)$ depends on $\delta_2$, the propellers $\tau(P_\g)$ and $\tau(P_\lambda)$ intersect $C_\mK^+(\tau(\cR'),\delta_2)$. Consider a loop starting by a segment that connects $\tau(p(i_0))$ to $\tau(p(i_0;i_1,\ell(\delta_2))) \in C_\mK^+(\tau(\cR'),\delta_2)$ and is tangent to $\fM_0^1$. Then connect the last point to $\tau(p(i_1))$ by a path in $E_{i_1}$ and then back to $\tau(p(i_0))$ along the Reeb cylinder. The value of $\ell(\delta_2)$ is the smallest number for which such a loop exists.
 Observe that the two loops defined by paths
 \begin{eqnarray*}
\tau(p(i_0)) \qquad \xrightarrow{\text{tangent}} \qquad & \tau(p(i_0;1,\ell(\delta_2))) & \qquad \xrightarrow{E_{1}} \qquad \tau(p(1)),\\
\tau(p(i_0)) \qquad  \xrightarrow{\text{tangent}} \qquad &  \tau(p(i_0;2,\ell(\delta_2))) & \qquad \xrightarrow{E_{2}} \qquad \tau(p(2))
\end{eqnarray*}
 are homologous, since the paths from the level 1 points $\tau(p(i_0;1,\ell(\delta_2)))$ and $\tau(p(i_0;2,\ell(\delta_2)))$ to $\tau(p(1))$, as well as the corresponding level 1 propeller, are contained in the $\delta_2$-neighborhood of the Reeb cylinder.\\

 \item  The one turn loops $S_1^2(i_0,i_1;\ell)$ at level 1, which   are formed as   follows. 
 Note that if $\tau(p(i_0;i_1,\ell))$ and $\tau(p(i_2;i_1,\ell+1))$ for $i_2\neq i_0$ do not belong to 
$C_\mK^+(\tau(\cR'),\delta_2)$,
there is loop formed by concatenating the path tangent to $\fM_0$ whose endpoints are $\tau(p(i_0;i_1,\ell))$ and $\tau(p(i_2;i_1,\ell+1))$ with the path joining these two points that is contained in $E_{i_1}\cap \fN_2$.  

Observe that for $\ell$ sufficiently large, the points  $\tau(p(i_0;i_1,\ell))$ and  $\tau(p(i_2;i_1,\ell+1))$ belong to the $\delta_2$-neighborhood of the Reeb cylinder, and thus the loop is homotopic to $P_1^2(i_2)-P_1^2(i_0)+[R]$, as explained in the proof of Proposition~\ref{prop-iota2} below.\\

 \item The exceptional loops $E^2(i_0,i_1,i_2;\ell_1)$ at level $2$ , which   are formed as   follows. 
 As in the case $k=1$, we define $\tau(p(i_0;i_1,\ell_1;i_2,b))$ to be a point in the curve $\tau(\gamma(i_1,\ell_1;i_2,b))$ if $i_0=1$, and in the curve $\tau(\lambda(i_1,\ell_1;i_2,b))$ if $i_0=2$. For $\ell_1\geq t(i_0;i_1)$, consider a path starting at $\tau(p(i_0;i_1,\ell_1))$ and tangent to the corresponding level 2 propeller up to the point $\tau(p(i_0;i_1,\ell_1;i_2,b))\in E_{i_2}$ (that is, tangent to $\tau(P_{\g(i_1,\ell_1)})$ if $i_0=1$ and to $\tau(P_{\lambda(i_1,\ell_1)})$ if $i_0=2$). Since $\tau(p(i_0;i_1,\ell_1;i_2,b))$ is contained in the region bounded by $\tau(\G(i_2,b))$ if $i_1=1$ and by $\tau(\Lambda(i_2,b))$ if $i_1=2$, we can connect $\tau(p(i_0;i_1,\ell_1;i_2,b))$ to $\tau(p(i_1;i_2,b))$  by a short segment in $E_{i_2}\cap \fN_2$ (as follows from Lemma~\ref{lem-paccumulationpoints}). Finally, connect $\tau(p(i_1;i_2,b))$ to $\tau(p(i_0;i_1,\ell_1))$ by a path tangent to $\fM_0^1$. 
 Observe that for the last step, there is only one possible choice up to multiples of $[R]$, since $\fM_0^1$ retracts to a tree (up to the loop corresponding to the Reeb cylinder). 
 These loops are illustrated in Figure~\ref{fig:Figure56a}.

\begin{figure}[!htbp]
\centering
{\includegraphics[width=105mm]{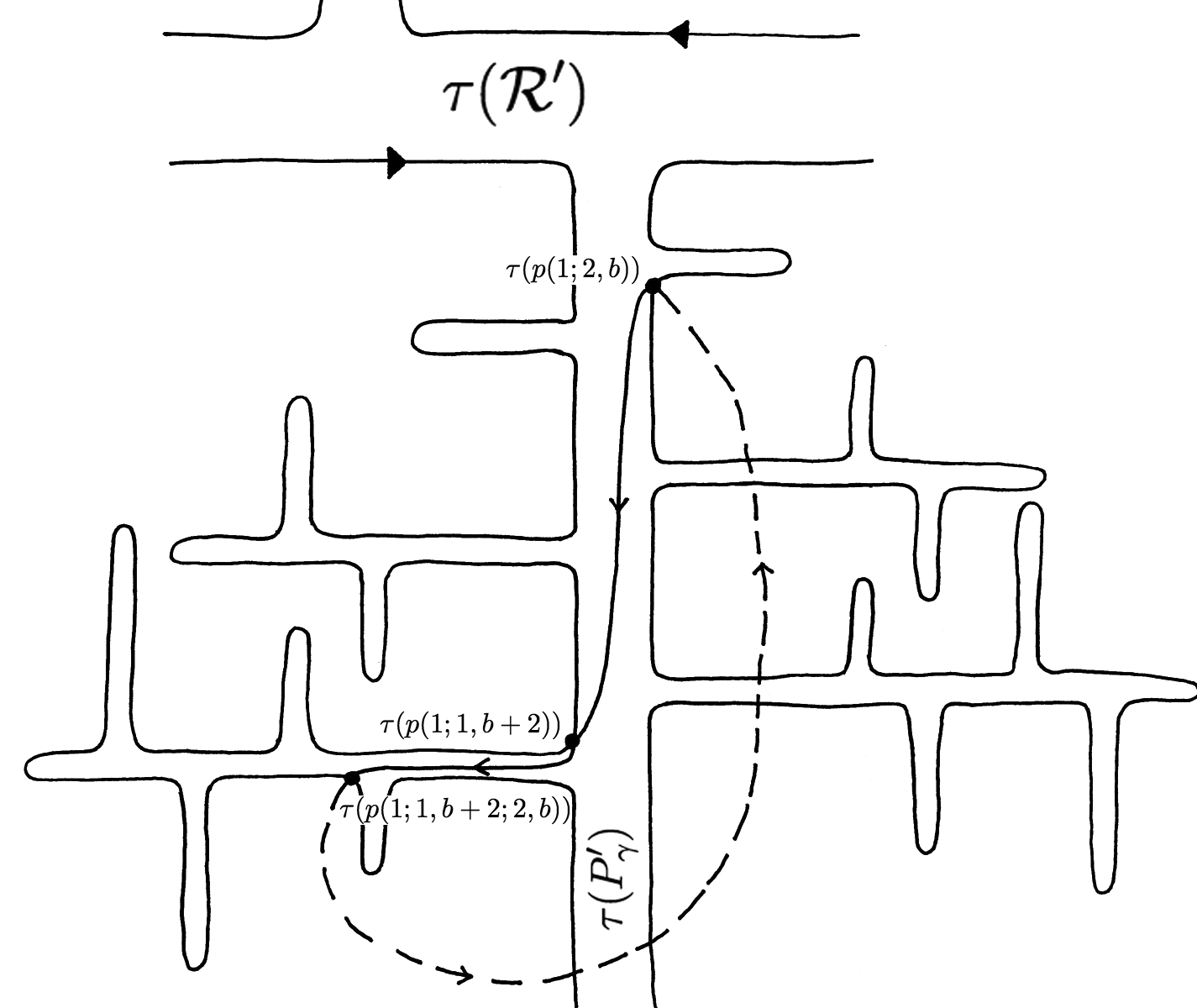}}
\caption{\label{fig:Figure56a} The exceptional loop $E(1,1,2;b+2)$ at level $k=2$}
\begin{minipage}[t]{14cm} 
Solid lines represent paths in the set $\fM_0^2$ and   dotted lines represent paths in the intersection of $\fN_2$ with the entry regions $E_1$ and $E_2$.
Here,  $i_0 =1$, $i_1 = 2$, $i_2 = 2$ and $\ell_1 = b+2$.
\end{minipage}
\end{figure}

 Observe that for $\ell_1 \geq \ell(\delta_2;i_0,i_1)$, these loops retract to $C_\mK^+(\fM_0^1,\delta_2)$. Thus there are at least $2^3(\ell(\delta_2;i_0,i_1)-t(i_0,i_1))$ such loops that are not homologous to loops in $C_\mK^+(\fM_0^1,\delta_2)$.\\

 \item The type (2) loops $T^2(i_0,i_1,i_2,i_3;\ell_1,\ell_2)$ at level 2. Recall that for $\ell_2\geq b$, the point $\tau(p(i_1;i_2,\ell_2))$ lies in  the curve $\tau(\G(i_2,\ell_2))$  if $i_1=1$, and in the curve $\tau(\Lambda(i_2,\ell_2))$ if $i_1=2$. Thus, we can connect this point to $\tau(p(i_0; i_1,\ell_1;i_2,\ell_2))$ for any $\ell_1$ such that the point exists and for $i_0=1,2$. The connecting path can taken to lie in $E_{i_2}\cap \fN_2$.  From $\tau(p(i_0; i_1,\ell_1;i_2,\ell_2))$ take a path tangent to the corresponding propeller to the point $\tau(p(i_0; i_1,\ell_1;i_3,\ell_2+1))$, for $i_3=1,2$. Since $\tau(p(i_0; i_1,\ell_1;i_3,\ell_2+1))\in E_{i_3}$ is in the region bounded by $\tau(\G(i_1,\ell_1;i_3,\ell_2+1))$ or by $\tau(\Lambda(i_1,\ell_1;i_3,\ell_2+1))$, we can connect this point to $\tau(p(i_1;i_3,\ell_2+1))$ by a path in $E_{i_3}\cap \fN_2$ and then back to $\tau(p(i_1;i_2,\ell_2))$ by a path tangent to $\fM_0^1$. These loops are illustrated in Figure~\ref{fig:Figure57a}.

\begin{figure}[!htbp]
\centering
{\includegraphics[width=105mm]{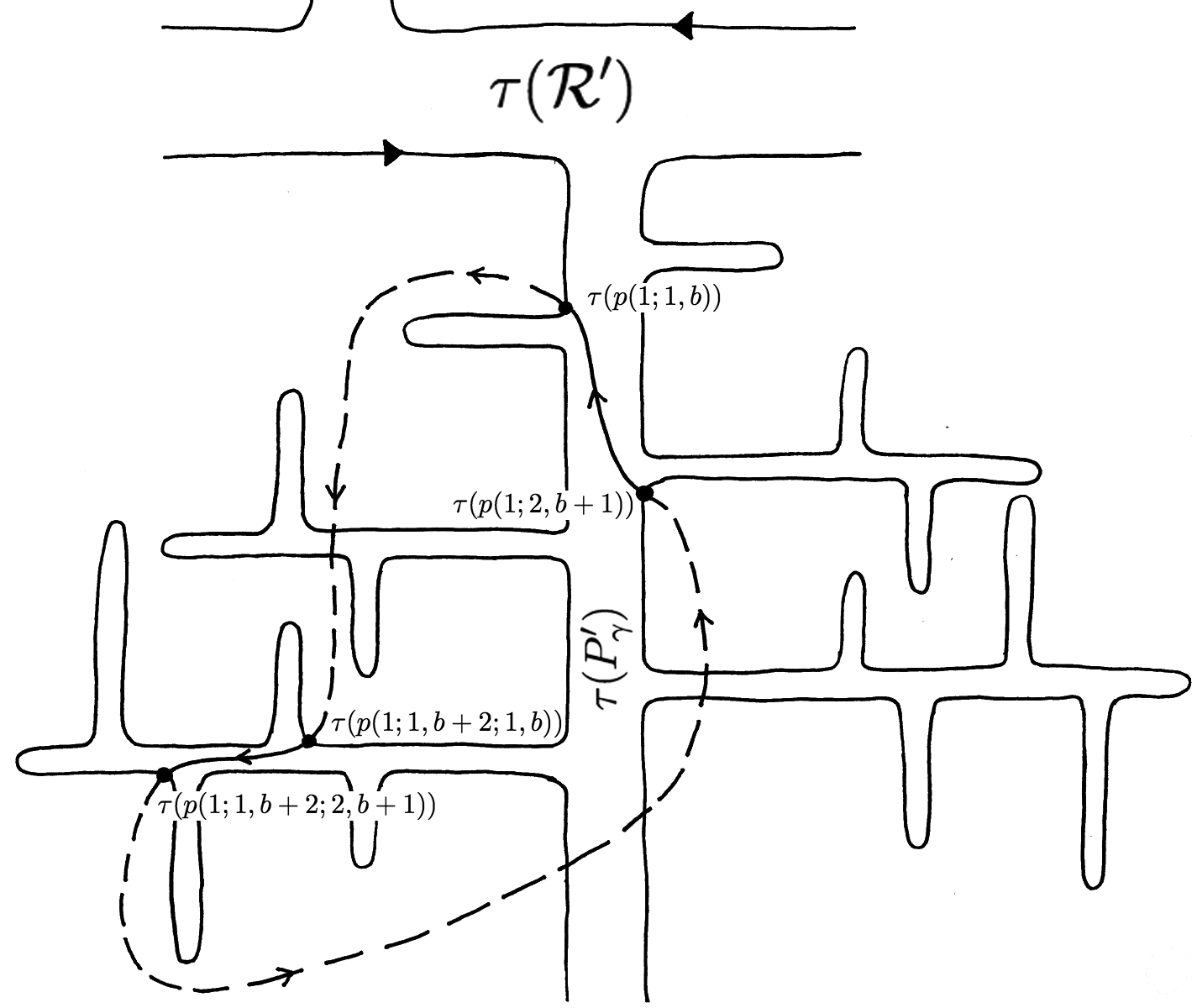}}
\caption{\label{fig:Figure57a}   The type (2)  loop $T^2(1,1,1,2;b+2,b+1)$  at level $k=2$}
\begin{minipage}[t]{14cm} 
Solid lines represent paths in the set $\fM_0^2$ and   dotted lines represent paths in $\fN_2$ with the entry regions $E_1$ and $E_2$.
Here,  $i_1 =1$, $i_2 = 1$, $i_3 = 2$, $\ell_1 = b+2$ and $\ell_2=b$.
\end{minipage}
\end{figure}

 Observe that for $\ell_2\geq \ell(\delta_2; i_0,i_1)$, the paths $T^2(i_0,i_1,i_2,i_3;\ell_1,\ell_2)$  are homologically trivial, as explained  in the proof of  Proposition~\ref{prop-iota2} below.  
 \end{enumerate}
 }
 \end{defn}
 
 In summary, the set of generators for $H_1(\fN_2; \mZ)$ are   
 constructed as follows.  First, we include the classes generated by the loops in $C_\mK^+(\fM_0^1,\delta_2)$, that are either at level 0 or level 1.  These are the loops described in cases (1), (2) and (3) above. Then we must include 
  the homology classes generated by the exceptional loops which arise  from  the branching of the level 2 propellers, as  described in cases  (4) and (5).

 \begin{prop}\label{prop-iota2}
 The inclusion $\fN_2\subset \fN_1$ induces a   map 
 $$\iota_2:H_1(\fN_2;\mZ) \to H_1(\fN_1;\mZ)$$
 such that:
 \begin{enumerate}
 \item $\iota_2(P_1^2(i_0))=T^1(i_0,i_1,i_1;\ell(\delta_2)-1)+ T^1(i_0,i_1,i_1;\ell(\delta_2)-2)+ \cdots + T^1(i_0,i_1,i_1;b)+ E^1(i_0,i_1)$ for any $i_1$;
 
 \item $\iota_2(S_1^2(i_0,i_1;\ell_1))$ is homologous to 
 $$T^1(i_2,i_1,i_1;\ell_1)+\cdots +T^1(i_2,i_1,i_1;b)+E^1(i_2,i_1)-E^1(i_0,i_1)-T^1(i_0,i_1,i_1;b)-\cdots-T^1(i_0,i_1,i_1;\ell_1-1);$$
 
 \item $\iota_2(E^2(i_0,i_1,i_2;\ell_1))=T^1(i_0,i_1,i_1;\ell_1-1)+ T^1(i_0,i_1,i_1;\ell_1-2)+ \cdots + T^1(i_0,i_1,i_1;b)+ E^1(i_0,i_1)$;
 
 \item $\iota_2(T^2(i_0,i_1,i_2,i_3;\ell_1,\ell_2))$ is the trivial element.
 \end{enumerate}
 \end{prop}
 
 In the proof below, there is a basic technique that is used to investigate the image of the different loops. This operation consist in moving by a homotopy   transformation, every   path segment tangent to a propeller at level 2, to a    path segment   tangent to a propeller at level 1, and joining the endpoints of these two segments by paths in the entry regions $E_1$ and $E_2$ accordingly. The reason why this is possible, is that the level 1 propellers in $\fN_1$ are filled,  and thus contain all path segments  tangent to higher level propellers. Once this homotopy operation is performed on a given loop, we then put the resulting deformed loop into    ``simplest form'' and evaluate the resulting homology class.
 
 \proof 
 Consider first an exceptional loop $E^2(i_0,i_1,i_2;\ell_1)$ of the class (4) which is illustrated in Figure~\ref{fig:Figure56a}. We consider its image  in $\fN_1$. Recall that the double propellers at level 1 are filled in $\fN_1$.  Consider the segment of $E^2(i_0,i_1,i_2;\ell_1)$ from $\tau(p(i_0;i_1,\ell_1))$ to $\tau(p(i_0;i_1,\ell_1;i_2,b))$ that is tangent to a level 2 propeller, as written in \eqref{level2-1} below. Since as $\ell_1\to \infty$, we have that the points $\tau(p(i_0;i_1,\ell_1))$ limit to $\tau(p(i_1))$, and the points $\tau(p(i_0;i_1,\ell_1;i_2,b))$ limit to $\tau(p(i_1;i_2,b))$. Thus,  the tangent path segment at level 2 is homotopic to the union of the  segment tangent to $\fM_0$ from $\tau(p(i_1))$ to $\tau(p(i_1;i_2,b))$ and short paths in $E_{i_1}$ and $E_{i_2}$. The first one of these short paths connects $\tau(p(i_0;i_1,\ell_1))$ to $\tau(p(i_1))$, and the second one connects $\tau(p(i_0;i_1,\ell_1;i_2,b))$ to $\tau(p(i_1;i_2,b))$, and is written in \eqref{level2-2} below. 
We can thus write the loop $E^2(i_0,i_1,i_2;\ell_1)$ in the following way:
 \begin{eqnarray}
 \tau(p(i_0;i_1,\ell_1))  & \xrightarrow{\text{tangent}} & \tau(p(i_0;i_1,\ell_1;i_2,b))\label{level2-1}\\
 & \xrightarrow{E_{i_2}} & \tau(p(i_1;i_2,b))\label{level2-2}\\
 & \xrightarrow{\text{tangent}} & \tau(p(i_0;i_1,\ell_1)).\label{level2-3}
\end{eqnarray}  
Then by the previous discussion, this loop is homologous  in $\fN_1$ to the loop 
\begin{equation}\label{level2-4}
\tau(p(i_0;i_1,\ell_1))\xrightarrow{E_{i_1}}  \tau(p(i_1)) \xrightarrow{\text{tangent}} \tau(p(i_1;i_2,b))  \xrightarrow{\text{tangent}} \tau(p(i_0;i_1,\ell_1)).
\end{equation}
Concatenate the   two tangential paths to  obtain the loop:
\begin{equation}\label{level2-5}
\tau(p(i_0;i_1,\ell_1))\xrightarrow{E_{i_1}}  \tau(p(i_1)) \xrightarrow{\text{tangent}} \tau(p(i_0;i_1,\ell_1)).
\end{equation}
 Observe that this loop is independent of $i_2$. Also, if $i_0\neq i_1$ the segment between $\tau(p(i_1))$ and $\tau(p(i_0;i_1,\ell_1))$ is homotopic to a path segment containing    $\tau(p(i_0))$. 
 
 The loop $P_1^2(i_0)$ corresponds to the following concatenation
 \begin{eqnarray}
 \tau(p(i_0)) & \xrightarrow{\text{tangent}} & \tau(p(i_0;i_1,\ell(\delta_2))) \label{level2-6}\\ 
 & \xrightarrow{E_{i_1}} &  \tau(p(i_1)) \label{level2-7}\\
 & \xrightarrow{\text{tangent}} & \tau(p(i_0)), \label{level2-8}
 \end{eqnarray}
 for any $i_1$.
Thus, in $\fN_1$ the loops $E^2(i_0,i_1,i_2;\ell_1)$ and $P_1^2(i_0)$ are homotopic.

In terms of the generators of $H_1(\fN_1;\mZ)$, for $\ell_1\geq b$ we   then have
\begin{equation}\label{level2-9}
\iota_2(E^2(i_0,i_1,i_2;\ell_1))=T^1(i_0,i_1,i_1;\ell_1-1)+ T^1(i_0,i_1,i_1;\ell_1-2)+ \cdots + T^1(i_0,i_1,i_1;b)+ E^1(i_0,i_1).
\end{equation}
  It is helpful to consider Figure~\ref{fig:Figure56a}. 
So while $E^1(i_0,i_1)$ is not necessarily in the image of $\iota_2$,   its composition with  level 1 loops of type (2) is in the image of $\iota_2$ for each pair $(i_0,i_1)$. 

For the  loops $P^2_1(i_0)$ at level 1, described in class (2) above, we obtain by similar reasoning that
\begin{equation}\label{level2-10}
\iota_2(P_1^2(i_0))=T^1(i_0,i_1,i_1;\ell(\delta_2)-1)+ T^1(i_0,i_1,i_1;\ell(\delta_2)-2)+ \cdots + T^1(i_0,i_1,i_1;b)+ E^1(i_0,i_1).
\end{equation}

 Next, consider the one turn loop $S_1^2(i_0,i_1;\ell_1)$, described in class (3) above,  that is formed when $\tau(p(i_0;i_1,\ell_1))$ and $\tau(p(i_2;i_1,\ell_1+1))$ for $i_2\neq i_0$ do not belong to 
$C_\mK^+(\tau(\cR'),\delta_2)$.
Schematically, this loop can written as
\begin{equation}\label{level2-11}
\tau(p(i_0;i_1,\ell_1))\xrightarrow{\text{tangent}} \tau(p(i_2;i_1,\ell_1+1)) \xrightarrow{E_{i_1}} \tau(p(i_0;i_1,\ell_1)).
\end{equation}
Since $i_2\neq i_0$ the tangent part of the loop can be written as the concatenation 
\begin{equation}\label{eq-S12}
\tau(p(i_0;i_1,\ell_1))\xrightarrow{\text{tangent}} \tau(p(i_0)) \xrightarrow{\text{tangent}} \tau(p(i_2))\xrightarrow{\text{tangent}} \tau(p(i_2;i_1,\ell_1+1)).
\end{equation}

If the two endpoints belong to $C_\mK^+(\tau(\cR'),\delta_1)$, we have that $\ell_1\geq\ell(\delta_1)$. Hence, after inclusion in $\fN_1$, we can further decompose the loop as the concatenation
\begin{eqnarray}
\tau(p(i_0;i_1,\ell_1)) & \xrightarrow{\text{tangent}} & \tau(p(i_0;i_1,\ell(\delta_1))) \xrightarrow{E_{i_1}} \tau(p(i_1)) \label{level2-12}\\
 & \xrightarrow{E_{i_1}} & \tau(p(i_0;i_1,\ell(\delta_1))) \xrightarrow{\text{tangent}} \tau(p(i_0)) \xrightarrow{\text{tangent}} \tau(p(i_1)) \label{level2-13}\\
& \xrightarrow{\text{tangent}} & \tau(p(i_2)) \xrightarrow{\text{tangent}} \tau(p(i_2;i_1,\ell(\delta_1)))\xrightarrow{E_{i_1}} \tau(p(i_1)) \label{level2-14}\\
& \xrightarrow{E_{i_1}} & \tau(p(i_2;i_1,\ell(\delta_1))) \xrightarrow{\text{tangent}}  \tau(p(i_2;i_1,\ell_1+1)). \label{level2-15}
\end{eqnarray}
The paths in the first and last lines \eqref{level2-12} and \eqref{level2-15} are segments contained in $C_\mK^+(\tau(\cR'),\delta_1)$. The path in the second line \eqref{level2-13} corresponds to the loop $-P_1^2(i_0)$, and the path in the third line \eqref{level2-14} corresponds to the loop $P_1^2(i_2)$. Let $\approx$ denote homological equivalence, then we conclude 
\begin{equation}\label{level2-16}
\iota_2(S_1^2(i_0,i_1;\ell_1)\approx \iota_2(P_1^2(i_2))-\iota_2(P_1^2(i_0))+m[R],
\end{equation}
for $m=-(\ell_1-\ell(\delta_1))+(\ell_1+1-\ell(\delta_1))=1$.

If at least one of the endpoints in \eqref{eq-S12} is not contained in $C_\mK^+(\tau(\cR'),\delta_1)$,  then in $\fN_1$ we can express the loop $\iota_2(S_1^2(i_0;i_1,\ell_1)$ as a concatenation of paths:
\begin{eqnarray}
\tau(p(i_0;i_1,\ell_1)) & \xrightarrow{\text{tangent}} & \tau(p(i_0;i_1,\ell_1-1)) \xrightarrow{E_{i_1}} \tau(p(i_1)) \xrightarrow{E_{i_1}} \tau(p(i_0;i_1,\ell_1-1))\label{S12-1}\\
&\xrightarrow{\text{tangent}} & \tau(p(i_0;i_1,\ell_1-2)) \xrightarrow{E_{i_1}} \tau(p(i_1)) \xrightarrow{E_{i_1}} \tau(p(i_0;i_1,\ell_1-2))\label{S12-2}\\
& \vdots & \nonumber \\
&\xrightarrow{\text{tangent}} & \tau(p(i_0;i_1,b)) \xrightarrow{\text{tangent}} \tau(p(i_0)) \xrightarrow{\text{tangent}} \tau(p(i_1))\xrightarrow{E_{i_1}}  \tau(p(i_0;i_1,b))\label{S12-3}\\
&  \xrightarrow{E_{i_1}} & \tau(p(i_1)) \xrightarrow{\text{tangent}} \tau(p(i_2)) \xrightarrow{\text{tangent}} \tau(p(i_2;i_1,\ell_1+1)).
\end{eqnarray}
Then the path in  \eqref{S12-1} is homotopic to $-T^1(i_0,i_1,i_1;\ell_1-1)$ followed by a segment contained in $E_{i_1}$,  from $\tau(p(i_0;i_1,\ell_1))$ to $\tau(p(i_0;i_1,\ell_1-1))$. Analogously, the path in line \eqref{S12-2} is homotopic to $-T^1(i_0,i_1,i_1;\ell_1-2)$ followed by the path egment contained in $E_{i_1}$, from $\tau(p(i_0;i_1,\ell_1-1))$ to $\tau(p(i_0;i_1,\ell_1-2))$. Continuing in this way, the path in line \eqref{S12-3} is homotopic to $-E^1(i_0,i_1)$ followed by the segment contained in $E_{i_1}$, from $\tau(p(i_0;i_1,b))$ to $\tau(p(i_1))$. Thus, the path corresponding to the first arrow in \eqref{eq-S12} is homotopic 
 in $\fN_1$ to 
\begin{equation}\label{level2-17}
-E^1(i_0,i_1)-T^1(i_0,i_1,i_1;b)-\cdots-T^1(i_0,i_1,i_1;\ell_1-1).
\end{equation}
Applying the same argument to the the path corresponding to the last arrow in \eqref{eq-S12} we conclude that $\iota_2(S_1^2(i_0,i_1;\ell_1))$ is homotopic to 
\begin{equation}\label{level2-18}
T^1(i_2,i_1,i_1;\ell_1)+\cdots +T^1(i_2,i_1,i_1;b)+E^1(i_2,i_1)-E^1(i_0,i_1)-T^1(i_0,i_1,i_1;b)-\cdots-T^1(i_0,i_1,i_1;\ell_1-1).
\end{equation}
Finally consider a type (2) loop at level 2, $T^2(i_0,i_1,i_2,i_3;\ell_1,\ell_2)$, as described in class (5) above:
\begin{eqnarray}
\tau(p(i_1;i_2,\ell_2)) & \xrightarrow{E_{i_2}} & \tau(p(i_0;i_1,\ell_1;i_2,\ell_2)) \label{level2-19}\\
& \xrightarrow{\text{tangent}} & \tau(p(i_0;i_1,\ell_1;i_3,\ell_2+1)) \label{level2-20}\\
& \xrightarrow{E_{i_3}} & \tau(p(i_1;i_3,\ell_2+1)) \label{level2-21}\\
& \xrightarrow{\text{tangent}} & \tau(p(i_1;i_2,\ell_2)). \label{level2-22}
\end{eqnarray}
In $\fN_1$ the first tangent segment becomes $\tau(p(i_1;i_2,\ell_2)) \xrightarrow{\text{tangent}}  \tau(p(i_1;i_3,\ell_2+1))$, and the points $\tau(p(i_1;i_2,\ell_2))$ and $\tau(p(i_1;i_3,\ell_2+1))$ can be joined to $\tau(p(i_2))$ and $\tau(p(i_3))$, respectively,   by  segments in $E_{i_2}$ and $E_{i_3}$, respectively, to obtain the path
$$\tau(p(i_1;i_3,\ell_2+1))\xrightarrow{\text{tangent}}  \tau(p(i_1;i_2,\ell_2))$$ 
 at levels 0 and 1. Hence, $i_2(T^2(i_0,i_1,i_2,i_3;\ell_1,\ell_2))$ is the loop 
\begin{eqnarray}
\tau(p(i_1;i_2,\ell_2)) & \xrightarrow{E_{i_2}} & \tau(p(i_1;i_2,\ell_2)) \label{level2-23}\\
& \xrightarrow{\text{tangent}} & \tau(p(i_1;i_3,\ell_2+1)) \label{level2-24}\\
& \xrightarrow{E_{i_3}} & \tau(p(i_1;i_3,\ell_2+1)) \label{level2-25}\\
& \xrightarrow{\text{tangent}} & \tau(p(i_1;i_2,\ell_2)), \label{level2-26}
\end{eqnarray}
which  yields the trivial class in $H_1(\fN_1,\mZ)$.  This conclusion can be seen in the illustration  Figure~\ref{fig:Figure57a}, as the two tangential path components become homotopic, but in reverse directions, when these tangential paths are homotoped to tangential paths in $\fN_1$.
 \endproof

The construction of $\fN_2$  has an interpretation in terms of the embedded tree $\TP \subset \fM_0$.
The closed set  $\fN_2$ defines a closed neighborhood of $\fM_0$ and hence of the embedded
tree $\TP \subset \fM_0$. Topologically, this 
corresponds to two operations: first
attaching, at a suitably large distance from the root point $\omega_0$,
the ends of the two level $1$ branches to $\omega_0$. 
Second, selecting a
 finite number, equal to $\sum_{i_0,i_1}\ell(\delta_2;i_0;i_1)-s(\delta_2;i_0;i_1)$,  of
 level 2 branches whose endparts get identified to a certain point in
 the level 1 branches. The branches of the
tree at higher levels are all collapsed into the filled regions, hence
to one of the $\cD(i_0;i_1,\ell)$.

  We  next describe the construction of $\fN_3$, which proceeds  in complete analogy with that of $\fN_2$, and so various repetitive  details are left to the reader.    
  
  For $0 < \delta_3<\delta_2/2$, introduce    the one-sided closed $\delta_3$-neighborhood of $\ds \fM_0^2$ given by  
 \begin{equation}\label{eq-retract3}
 C_{\mK}^+(\fM_0^2,\delta_3)=C_{\mK}(\fM_0^2,\delta_3)\cap \fN_2.
 \end{equation}  
 Choose $\delta_3$ sufficiently small, so that $C_{\mK}^+(\fM_0^2,\delta_3)$ does not contain all the level 3 propellers. 
 
The double propellers at level 3 are given by the collections
\begin{equation}\label{eq-filled3}
\tau(P'_{\G(i_1,\ell_1;i_2,\ell_2)}) ~ {\rm and} ~ \tau(P'_{\Lambda(i_1,\ell_1;i_2,\ell_2)}) \quad {\rm for} ~ \ell_1,\ell_2\geq b ~ {\rm and}  ~ i_1,i_2=1,2 ~ .
\end{equation}
 Denote the corresponding filled double propellers in $\mK$  by $\cD(1;i_1,\ell_1;i_2,\ell_2)$ and $\cD(2;i_1,\ell_1;i_2,\ell_2)$, respectively.  
Observe that the lengths of the double propellers $\ds P'_{\G(i_1,\ell_1;i_2,\ell_2)}$ and $\ds P'_{\Lambda(i_1,\ell_1;i_2,\ell_2)}$ are  not bounded above as $\ell_1, \ell_2 \to \infty$. 
Then set:
\begin{equation}\label{eq-N3}
\fN_3 ~ = ~     C_{\mK}^+(\fM_0^2, \delta_3)  ~ \cup ~  \left\{ ~ \bigcup \left\{ \cD(i_0;i_1,\ell_1;i_2,\ell_2) \mid i_0,i_1,i_2 = 1,2 ~ \& ~ \ell_1,\ell_2 \geq b \right\} ~ \right\} ~  \subset ~ \mK .
\end{equation}

As in the previous analysis of the homotopy types  of $\fN_1$ and $\fN_2$, the filled double propellers in \eqref{eq-N3} have both intersections and self-intersections, so their addition 
to the space $\ds C_{\mK}^+(\fM_0^2, \delta_3)$ adds multiple types of exceptional cycles to the homology of $\fN_3$.

As in the   analysis of the homotopy type  of   $\fN_2$, we define the constant $\ell(\delta_3)$ to be the least integer such that the point  $\tau(p(i_0;i_1,\ell(\delta_3)))$ is contained in the $\delta_3$-neighborhood of the Reeb cylinder, and the constants $\ell(\delta_3;i_0,i_1)>s(\delta_3;i_0,i_1)\geq b$ are defined analogously to the level $k=2$ case   \eqref{eq-level2fork=2}. Moreover, there are constants
$t(i_0;i_1,\ell_1;i_2)$ such that the  filled double propellers
$\cD(i_0;i_1,\ell_1;i_2,\ell)$ 
intersect the regions $E_1$ and $E_2$ for with $\ell\geq t(i_0;i_1,\ell_1;i_2)$, hence they are not homologically trivial. We then have
\begin{equation}\label{eq-level2fork=3}
 \ell(\delta_3;i_0;i_1,\ell_1;i_2) > s(\delta_3;i_0;i_1,\ell_1;i_2)  \geq b
\end{equation}
 and the following conditions are satisfied:
\begin{enumerate}
\item For $b\leq \ell<s(\delta_3;i_0;i_1,\ell_1;i_2)$,  the filled double propellers   $\cD(i_0;i_1,\ell_1;i_2,\ell)$ intersect  $C_{\mK}^+(\fM_0^2, \delta_3)$ only near their generating curves $\ds P'_{\G(i_1,\ell_1;i_2,\ell_2)}$ and $\ds P'_{\Lambda(i_1,\ell_1;i_2,\ell_2)}$. \\
\item For $s(\delta_3;i_0;i_1,\ell_1;i_2)\leq \ell  <\ell(\delta_3;i_0;i_1,\ell_1;i_2)$,  the filled double   propeller $\cD(i_0;i_1,\ell_1;i_2,\ell)$ is such that its base and tip are contained in $C_{\mK}^+(\fM_0^2,\delta_3)$ and do not retract in $\fN_3$ to $C_{\mK}^+(\fM_0^2,\delta_3)$. 
Thus, each of these propellers adds a handle  to $C_{\mK}^+(\fM_0^2,\delta_3)$. \\
\item For $\ell\geq \ell(\delta_3;i_0;i_1,\ell_1;i_2)$, the attachment of  the filled double propellers
  $\cD(i_0;i_1,\ell_1;i_2,\ell)$ does not change the   homotopy type of $C_{\mK}^+(\fM_0^2,\delta_3)$.
\end{enumerate}

It follows that the analogs of  Lemmas~\ref{lem-shape2-2}, \ref{lem-shape2-3} and   \ref{lem-shape2-14}    hold for $\fN_3$.

We next  describe the image of the   map on homology   induced by the inclusion.
A generating set for $H_1(\fN_3;\mZ)$ is constructed in an analogous way as before. 

\begin{defn}\label{def-generators3}
Consider the following  generators of $H_1(\fN_3;\mZ)$:
 
{\rm 

\begin{enumerate}
\item The loop $[R]$ corresponding to the fundamental class of the Reeb cylinder.\\

\item Level 1 and 2  loops. 
\begin{enumerate}
\item Level 1 loops $P^3_1(i_0)$. These are the same as  the loops of  class (2)   generating $H_1(\fN_2;\mZ)$, represented schematically for   $i_1 =1,2$, as
$$\tau(p(i_0)) \xrightarrow{\text{tangent}} \tau(p(i_0;i_1,\ell(\delta_3))) \xrightarrow{E_{i_1}} \tau(p(i_1)) \xrightarrow{\text{tangent}} \tau(p(i_0)).
$$

\item   Level 2 loops $P^3_2(i_0,i_1;\ell_1)$. The level 2 propellers limit to the level 1 propellers, and get longer as the indices $(i_0;i_1,\ell_1)$ increase, 
  for $\ell_1\geq s(\delta_3;i_0;i_1)$, the  propeller $\tau(P_{\g(i_1,\ell_1)})$ if $i_0=1$, or the propellers $\tau(P_{\lambda(i_1,\ell_1)})$ if $i_0=2$, intersect $C_\mK^+(\fM_0^1,\delta_3)$ near their base and tip, and thus represent handles  of $\fN_3$. Then  for $s(\delta_3;i_0;i_1)\leq \ell_1<\ell(\delta_3;i_0,i_1)$,   consider the loops starting at $\tau(p(i_0;i_1,\ell_1))$ tangent to $\fM_0^2$ up to the point $\tau(p(i_0;i_1,\ell_1;i_2,\ell_2))$. For $\ell_2$ sufficiently large, the point $\tau(p(i_0;i_1,\ell_1;i_2,\ell_2))$ is contained in $C_\mK^+(\fM_0^1,\delta_3)$, thus we can join this point to $\tau(p(i_1;i_2,\ell_2))$ and then close the loop by a path tangent to $\fM_0^1$. 

The  homology class of the   loop $P^3_2(i_0,i_1;\ell_1)$ is independent of $i_2 =1,2$, as was the case for  the analogous level 1 loops.
 Given $(i_0,i_1,\ell_1)$, consider the smallest $\ell_2$ such that $\tau(p(i_0;i_1,\ell_1;i_2,\ell_2))$ is contained in $C_\mK^+(\fM_0^1,\delta_3)$ for   $i_2 =1,2$. Since the part of the level 2 propeller that lies between the points $\tau(p(i_0;i_1,\ell_1;i_2,\ell_2))$ and the tip retracts to $C_\mK^+(\fM_0^1,\delta_3)$, for any $\ell_2'>\ell_2$ the loops 
$$
\tau(p(i_0;i_1,\ell_1)) \xrightarrow{\text{tangent}} \tau(p(i_0;i_1,\ell_1;i_2,\ell_2)) \xrightarrow{E_{i_2}} \tau(p(i_1;i_2,\ell_2)) \xrightarrow{\text{tangent}} \tau(p(i_0;i_1,\ell_1))$$ 
$$
\tau(p(i_0;i_1,\ell_1)) \xrightarrow{\text{tangent}} \tau(p(i_0;i_1,\ell_1;i_2,\ell_2')) \xrightarrow{E_{i_2}} \tau(p(i_1;i_2,\ell_2')) \xrightarrow{\text{tangent}} \tau(p(i_0;i_1,\ell_1))$$
are homotopic.\\
\end{enumerate}

\item The one turn loops.
\begin{enumerate}
\item  Level 1 loops $S_1^3(i_0,i_1;\ell)$.  As in the case of the loops of class (3) for $\fN_2$,  consider a pair of points $\tau(p(i_0;i_1,\ell))$ and $\tau(p(i_2;i_1,\ell+1))$ for $i_2\neq i_0$. By the definition in \eqref{eq-retract3}, the two points can be connected by a path in $E_{i_1}\cap \fN_3$ if they belong to 
 $$C_\mK^+(\tau(\cR'),\delta_2)-C_\mK^+(\tau(\cR'),\delta_3),$$
 and the $\delta_3$-neighborhood of $\tau(p(i_2;i_1,\ell+1))$ contains $\tau(p(i_0;i_1,\ell))$. In this case,
there is a loop formed by concatenating the path tangent to $\fM_0$ between the points   $\tau(p(i_0;i_1,\ell))$ and $\tau(p(i_2;i_1,\ell+1))$, with the path joining these two points that is contained in $E_{i_1}\cap \fN_3$.

  Observe that for $\ell$ sufficiently large, the points  $\tau(p(i_0;i_1,\ell))$ and  $\tau(p(i_2;i_1,\ell+1))$ belong to the $\delta_3$-neighborhood of the Reeb cylinder, and thus the homology class of the loop $S_1^3(i_0,i_1;\ell)$ is equal  to $P_1^3(i_2)-P_1^3(i_0)+[R]$, as explained in the proof of Proposition~\ref{prop-iota2}.\\
 
\item  Level 2 one turn loops  $S_2^3(i_0,i_1,i_2;\ell_1,\ell_2)$. For a point $\tau(p(i_0;i_1,\ell_1;i_2,\ell_2))$, the nesting property of the ellipses, as discussed in Section~\ref{sec-doublepropellers}, implies that  $\tau(p(i_3;i_1,\ell_1+1;i_2,\ell_2))$ for $i_3\neq i_0$ is the closest point at the same level. If the $\delta_3$-neighborhood of $\tau(p(i_3;i_1,\ell_1+1;i_2,\ell_2))$ contains $\tau(p(i_0;i_1,\ell_1;i_2,\ell_2))$, then  there is loop formed by concatenating the path tangent to $\fM_0$ whose endpoints are $\tau(p(i_0;i_1,\ell_1;i_2,\ell_2))$ and $\tau(p(i_3;i_1,\ell_1+1;i_2,\ell_2))$, with the path contained in $E_{i_2}\cap \fN_3$ joining these two points.

 Observe that for $\ell_1$ sufficiently large, the points   $\tau(p(i_0;i_1,\ell_1;i_2,\ell_2))$ and  $\tau(p(i_3;i_1,\ell_1+1;i_2,\ell_2))$ belong to the $\delta_3$-neighborhood of $\fM_0^1$ and thus the loop is homotopic to a loop in $C_\mK^+(\fM_0^1,\delta_3)$.\\
\end{enumerate}

\item The exceptional loops $E^3(i_0,i_1,i_2,i_3;\ell_1,\ell_2)$.   As in the case of the loops of class (4) for $\fN_2$,  let $\tau(p(i_0;i_1,\ell_1;i_2,\ell_2;i_3,b))$ be any point in the curve   $\tau(\g(i_1,\ell_1;i_2,\ell_2;i_3,b))\subset E_{i_3}$ if $i_0=1$, and in the curve $\tau(\lambda(i_1,\ell_1;i_2,\ell_2;i_3,b))\subset E_{i_3}$ if $i_0=2$. For $\ell_2\geq t(i_0;i_1,\ell_1;i_2)$, consider a loop starting at $\tau(p(i_0;i_1,\ell_1;i_2,\ell_2))$, then tangent to the corresponding level 3 propeller up to the point $\tau(p(i_0;i_1,\ell_1;i_2,\ell_2;i_3,b))\in E_{i_3}$, for $i_3=1,2$. Then join this last point to $\tau(p(i_1;i_2,\ell_2;i_3,b))$ by a short segment contained in $E_{i_3}\cap \fN_3$, and return tangent to $\fM_0^2$ to the point $\tau(p(i_0;i_1,\ell_1;i_2,\ell_2))$.

Observe that for $\ell_2\geq \ell(\delta_3;i_0;i_1,\ell_1;i_2)$ the loops $E^3(i_0,i_1,i_2,i_3;\ell_1,\ell_2)$ are contained in $\fM_0^2$ and thus are homotopic  to level 1 and level 2 loops. \\

\item The type (2) loops $T^3(i_0,i_1,i_2,i_3,i_4;\ell_1,\ell_2,\ell_3)$ at level 3. 
For  $\ell_3\geq b$, the point $\tau(p(i_1;i_2,\ell_2;i_3,\ell_3))$ is contained in the curve $\tau(\G(i_2,\ell_2;i_3,\ell_3))$  if $i_1=1$, and in the curve $\tau(\Lambda(i_2,\ell_2;i_3,\ell_3))$ if $i_1=2$. Thus,   we can connect this point to  $\tau(p(i_0; i_1,\ell_1;i_2,\ell_2;i_3,\ell_3))$ for any $i_0$, and for $\ell_1$ big enough so that the last point exists. The connecting path can taken to lie in $E_{i_3}\cap \fN_3$.  From $\tau(p(i_0; i_1,\ell_1;i_2,\ell_2;i_3,\ell_3))$, take a path tangent to the corresponding level 3 propeller to the point $\tau(p(i_0; i_1,\ell_1;i_2,\ell_2;i_4,\ell_3+1))$, for $i_4=1,2$. Since 
 $\ds \tau(p(i_0; i_1,\ell_1;i_2,\ell_2;i_4,\ell_3+1))\in E_{i_4}$ 
 is in the region bounded by $\tau(\G(i_2,\ell_2;i_4,\ell_3+1))$ or by $\tau(\Lambda(i_2,\ell_2;i_4,\ell_3+1))$, we can connect this point to $\tau(p(i_1;i_2,\ell_2;i_4,\ell_3+1))$ by a path in $E_{i_4}\cap \fN_3$ and then back to $\tau(p(i_1;i_2,\ell_2;i_3,\ell_3))$ by a path tangent to $\fM_0^2$.
 \end{enumerate}
}
\end{defn}

\begin{remark}\label{rmk-generators}
The loops in  (1), (2a) and (3a)  above   generate the   homology group   $H_1(C_\mK^+(\fM_0^1,\delta_3); \mZ)$. Adding the loops in (2b) and (3b) in the above listing complete the list of generators for the    homology group  $H_1(C_\mK^+(\fM_0^2,\delta_3); \mZ)$. As in the case of $\fN_2$, the exceptional and type (2) loops are given by the intersections of the filled propellers, in this case at level 3.
\end{remark}
 
  \begin{prop}\label{prop-iota3}
 The inclusion $\fN_3\subset \fN_2$ induces a map 
 $$\iota_3:H_1(\fN_3;\mZ) \to H_1(\fN_2;\mZ)$$
 such that:
 \begin{enumerate}
  \item $\iota_3(P^3_1(i_0))=(\ell(\delta_3)-\ell(\delta_2))[R]+ P^2_1(i_0)$;
  
  \item $\iota_3(P^3_2(i_0,i_1;\ell_1))$ is homologous to
 $$T^2(i_0,i_1,i_2,i_2;\ell_1,\ell_2-1)+ T^2(i_0,i_1,i_2,i_2;\ell_1,\ell_2-2)+ \cdots + T^2(i_0,i_1,i_2,i_2;\ell_1,b)+ E^2(i_0,i_1,i_2;\ell_1),$$
 for any $i_2$ and $\ell_2$ the minimum possible value so that the loop $P^3_2(i_0,i_1;\ell_1)$ exists.
 
 \item $\iota_3(S_1^3(i_0;i_1,\ell_1))$ is homologous to $\iota_3(P_1^3(i_2))-\iota_3(P_1^3(i_0))+[R]$;
 
 \item $\iota_3(S_2^3(i_0,i_1,i_2;\ell_1,\ell_2))$ is homologous to 
 $$\iota_3(-P_2^3(i_0,i_1;\ell_1) - P_1^3(i_0)+P_1^3(i_3)+P_2^3(i_3;i_1,\ell_1+1)+[R]),$$
 for $i_3\neq i_0$ and up to the sum of type (2) loops at level 2;
 
 \item $\iota_3(E^3(i_0,i_1,i_2,i_3;\ell_1,\ell_2))$ is homologous to 
 $$T^2(i_0,i_1,i_2,i_2;\ell_1,\ell_2-1)+ T^2(i_0,i_1,i_2,i_2;\ell_1,\ell_2-2)+ \cdots + T^2(i_0,i_1,i_2,i_2;\ell_1,b)+ E^2(i_0,i_1,i_2;\ell_1);$$

 \item $\iota_3(T^3(i_0,i_1,i_2,i_3,i_4;\ell_1,\ell_2,\ell_3))$ is the trivial element.
 \end{enumerate}
 \end{prop}
 
 The idea of the proof is the same as that of Proposition~\ref{prop-iota2}: given a loop in $\fN_3$ each path segment in the loop that is tangent to a propeller at level 3, can be homotoped  in the filled propeller to a path segment    tangent to a propeller at level 2, since the propellers at level 2 are filled in $\fN_2$. We then use segments in the faces of the filled propellers in $\fN_2$ to close up the paths into a loop contained in $\fN_2$. The resulting path is then identified up to homology with a combination of the generators for $H_1(\fN_2;\mZ)$. We sketch the details of the individual cases below.
 
 \proof 
 There are two immediate cases which proceed exactly as in the proof of Proposition~\ref{prop-iota2}.
The image of a level 1 loop of $\fN_3$ is equal in homology to the addition  of the class defined by  a level 1 loop in $\fN_2$  with a  multiple of  $[R]$. 
We give the details for the level 2 one turn loops, and for the level 2 exceptional loops.

Consider a level 2 loop $P^3_2(i_0,i_1,i_2;\ell_1,\ell_2)$. Assume without loss of generality that $i_0=1$. Observe that the loop $P^3_2(i_0,i_1,i_2;\ell_1,\ell_2)$ exists if $\ell_2\geq s(\delta_3;1;i_1,\ell_1;i_2)$ and that $s(\delta_3;1;i_1,\ell_1;i_2)>s(\delta_2;1;i_1)$. Also, the level 2 propeller $\tau(P_{\g(i_1,\ell_1)})$ retracts to $C_\mK^+(\fM_0^1,\delta_2)$ only if $\ell_2>\ell(\delta_2;1;i_1)$. 
We write schematically the loop $P^3_2(1,i_1,i_2;\ell_1,\ell_2)$ as
\begin{eqnarray*}
\tau(p(i_0;i_1,\ell_1)) &\xrightarrow{\text{tangent}}  & \tau(p(i_0;i_1,\ell_1;i_2,\ell_2))\\
& \xrightarrow{E_{i_2}}  & \tau(p(i_1;i_2,\ell_2))\\
&\xrightarrow{\text{tangent}}  & \tau(p(i_0;i_1,\ell_1)).
\end{eqnarray*}
Since $\tau(P_{\g(i_1,\ell_1)})\subset \cD(1;i_1,\ell_1)$, the loop $P^3_2(1,i_1,i_2;\ell_1,\ell_2)$    is homotopic  to the composition of exceptional and type (2) loops in $\fN_2$, so that we obtain 
 \begin{align*}
 \iota_3(P^3_2(i_0,i_1,i_2;\ell_1,\ell_2)) &\approx T^2(i_0,i_1,i_2,i_2;\ell_1,\ell_2-1)+ T^2(i_0,i_1,i_2,i_2;\ell_1,\ell_2-2)+ \cdots\\
 & + T^2(i_0,i_1,i_2,i_2;\ell_1,b)+ E^2(i_0,i_1,i_2;\ell_1).
 \end{align*}
 
 Consider now an exceptional loop $E^3(i_0,i_1,i_2,i_3;\ell_1,\ell_2)$. Since the double propellers at level 2 are filled in $\fN_2$, the tangent segment from $\tau(p(i_0;i_1,\ell_1;i_2,\ell_2))$ to $\tau(p(i_0;i_1,\ell_1;i_2,\ell_2;i_3,b))$ can be identified with the segment from $\tau(p(i_1;i_2,\ell_2))$ to $\tau(p(i_1;i_2,\ell_2;i_3,b))$. Hence the image of an exceptional loop becomes
\begin{eqnarray}
\tau(p(i_0;i_1,\ell_1;i_2,\ell_2)) & \xrightarrow{E_{i_2}}  & \tau(p(i_1;i_2,\ell_2)) \label{eq-block3-line1}\\
& \xrightarrow{\text{tangent}} &\tau(p(i_1;i_2,\ell_2;i_3,b)) \label{eq-block3-line2}\\ 
& \xrightarrow{\text{tangent}}  & \tau(p(i_1;i_2,\ell_2))\\ & \xrightarrow{\text{tangent}} &\tau(p(i_0;i_1,\ell_1;i_2,\ell_2)). \label{eq-block3-line3}
\end{eqnarray}
 This simplifies to 
 $\tau(p(i_0;i_1,\ell_1;i_2,\ell_2))\xrightarrow{E_{i_2}}  \tau(p(i_1;i_2,\ell_2))  \xrightarrow{\text{tangent}} \tau(p(i_0;i_1,\ell_1;i_2,\ell_2))$, 
 from which we obtain 
 \begin{align*}
 \iota_3(E^3(i_0,i_1,i_2,i_3;\ell_1,\ell_2))&= T^2(i_0,i_1,i_2,i_2;\ell_1,\ell_2-1)+ T^2(i_0,i_1,i_2,i_2;\ell_1,\ell_2-2)+ \cdots \\
 &+ T^2(i_0,i_1,i_2,i_2;\ell_1,b)+ E^2(i_0,i_1,i_2;\ell_1).
 \end{align*}

Consider the one turn loop $S_1^3(i_0,i_1;\ell_1)$ that is formed when $\tau(p(i_0;i_1,\ell_1))$ and $\tau(p(i_2;i_1,\ell_1+1))$, for $i_2\neq i_0$, belong to $C_\mK^+(\tau(\cR'),\delta_2) - C_\mK^+(\tau(\cR'),\delta_3)$, thus $\ell_1\geq \ell(\delta_2)$. Hence, as in the proof of Proposition~\ref{prop-iota2} we obtain that $\iota_3(S_1^3(i_0,i_1;\ell_1))$ is homologous to 
 $$\iota_3(- P_1^3(i_0)+P_1^3(i_3)+[R]).$$
 
 For the level 2 one turn loops, we can express $S_2^3(i_0,i_1,i_2;\ell_1,\ell_2)$ as the concatenation of paths, 
 \begin{eqnarray}
 \tau(p(i_0;i_1,\ell_1;i_2,\ell_2)) &\xrightarrow{\text{tangent}}& \tau(p(i_0;i_1,\ell_1)) \label{eq-block3-line4}\\
 &\xrightarrow{\text{tangent}}& \tau(p(i_0)) \label{eq-block3-line5}\\
 &\xrightarrow{\text{tangent}}& \tau(p(i_3;i_1,\ell_1+1)) \label{eq-block3-line6}\\
 &\xrightarrow{\text{tangent}}& \tau(p(i_3;i_1,\ell_1;i_2,\ell_2)) \label{eq-block3-line7}\\
 & \xrightarrow{E_{i_2}} & \tau(p(i_0;i_1,\ell_1;i_2,\ell_2)). \label{eq-block3-line8}
 \end{eqnarray}
In $\fN_2$, the tangent path in \eqref{eq-block3-line4} is homotopic to 
$$-T^2(i_0,i_1,i_2,i_2;\ell_1,\ell_2-1)-\cdots - T^2(i_0,i_1,i_2,i_2;\ell_1,b)- E^2(i_0,i_1,i_2;\ell_1),$$
the the tangent path in \eqref{eq-block3-line5} is homotopic to
\begin{equation}\label{eq-block3-line9}
-(\ell_1-\ell(\delta_2))[R]-P_1^2(i_0),
\end{equation}
the the tangent path in \eqref{eq-block3-line6} is homotopic to 
\begin{equation}\label{eq-block3-line10}
(\ell_1+1-\ell(\delta_2))[R]+P_1^2(i_3),
\end{equation}
and the the tangent path in \eqref{eq-block3-line7} is homotopic to 
\begin{equation}\label{eq-block3-line11}
T^2(i_3,i_1,i_2,i_2;\ell_1+1,\ell_2-1)+\cdots + T^2(i_3,i_1,i_2,i_2;\ell_1+1,b)+ E^2(i_3,i_1,i_2;\ell_1+1).
\end{equation}
Thus $\iota_3(S_2^3(i_0,i_1,i_2;\ell_1,\ell_2))$ is homologous to 
\begin{eqnarray}
-T^2(i_0,i_1,i_2,i_2;\ell_1,\ell_2-1) & -&\cdots - T^2(i_0,i_1,i_2,i_2;\ell_1,b)- E^2(i_0,i_1,i_2;\ell_1) \label{eq-block3-line12}\\
& + & [R] -P_1^2(i_0) +P_1^2(i_3) \label{eq-block3-line13}\\
&+& T^2(i_3,i_1,i_2,i_2;\ell_1+1,\ell_2-1)+\cdots + T^2(i_3,i_1,i_2,i_2;\ell_1+1,b) \label{eq-block3-line14}\\
& + & E^2(i_3,i_1,i_2;\ell_1+1). \label{eq-block3-line15}
\end{eqnarray}

 Analogously to the case of $\iota_2$, the homology class of  $T^3(i_0,i_1,i_2,i_3,i_4;\ell_1,\ell_2,\ell_3)$ is trivial.
 \endproof

 In general,   for $k > 3$, we consider $\delta_k<\delta_{k-1}/2$ such
 that  $ C_{\mK}^+(\fM_0^{k -1}, \delta_{k}) = C_{\mK}(\fM_0^{k -1}, \delta_{k}) \cap \fN_{k -1}$ does not
contain all of the propellers at level $k$.
 
 The space  $\fN_{k}$ is defined by  attaching to $ C_{\mK}^+(\fM_0^{k -1}, \delta_{k})$ the collection of filled
 propellers at level $k$. 
As the number of handles in $\fN_{k-1}$
 depends on $\delta_{k-1}$, the sets $\fN_{k-1}$ and $C_{\mK}^+(\fM_0^{k   -1}, \delta_{k})$ are not  homotopy  equivalent  for $\delta_k$ small enough. 
Analogous versions of Lemmas~\ref{lem-shape2-2}, \ref{lem-shape2-3} and \ref{lem-shape2-14}  hold for $\fN_{k}$. 

\begin{lemma}\label{lem-shapek}
The set $\fN_k$ is compact and satisfies:
\begin{enumerate}
\item $\fM\subset \fN_k\subset \fN_{k-1}$;
\item For $\e>0$ sufficiently small, $U(\fN_k,\e)$ retracts to $\fN_k$;
\item $\fN_k$ has the homotopy type of a finite wedge of circles.
\end{enumerate}
\end{lemma}

We now generalize the description of the generators of the homology and give a lower bound on the rank of the homology groups $H_1(\fN_k;\mZ)$.
\begin{defn}\label{def-generatorsk}
 Consider the following  generators of $H_1(\fN_k;\mZ)$:

{\rm 

\begin{enumerate}
\item The loop $[R]$ corresponding to the fundamental class of the Reeb cylinder.\\

\item Level 1 to $k-1$ loops.
\begin{enumerate}
\item The level 1 loops  $P_1^k(i_0)$ are defined by choosing $i_1$ then $P_1^k(i_0)$ is the loop:
\begin{eqnarray*}
\tau(p(i_0)) & \xrightarrow{\text{tangent}}  &\tau(p(i_0;i_1,\ell(\delta_k)))\\
& \xrightarrow{E_{i_1}}  & \tau(p(i_1))\\
& \xrightarrow{\text{tangent}}  & \tau(p(i_0)). 
\end{eqnarray*}

\item For $2\leq n\leq k-1$, the level $n$ loops $P_n^k(i_0,i_1,\cdots,i_{n-1};\ell_1,\cdots,\ell_{n-1})$ are defined as:
\begin{eqnarray*}
\tau(p(i_0;i_1,\ell_1;\cdots;i_{n-1},\ell_{n-1})) & \xrightarrow{\text{tangent}}  &\tau(p(i_0;i_1,\ell_1;\cdots;i_n,\ell_n))\\
& \xrightarrow{E_{i_n}}  & \tau(p(i_1;i_2,\ell_2;\cdots;i_n,\ell_n))\\
& \xrightarrow{\text{tangent}}  & \tau(p(i_0;i_1,\ell_1;\cdots;i_{n-1},\ell_{n-1})),
\end{eqnarray*}
for any $i_n$ and where $\ell_n$ is  the smallest number for which such a loop exist. \\
\end{enumerate}

\item  The one turn loops at level $n$ for $1\leq n\leq k-1$, denoted by $S_n^k(i_0,i_1,\cdots, i_n;\ell_1,\cdots,\ell_n)$, is represented schematically as
\begin{eqnarray*}
\tau(p(i_0;i_1,\ell_1;\cdots;i_n,\ell_n)) & \xrightarrow{\text{tangent}}  & \tau(p(i_0;i_1,\ell_1+1;\cdots;i_n,\ell_n))\\
& \xrightarrow{E_{i_1}}  & \tau(p(i_0;i_1,\ell_1;\cdots;i_n,\ell_n)).
\end{eqnarray*}

\item  The exceptional loops at level $k$, denoted by $E^k(i_0,i_1,\cdots, i_k;\ell_1,\cdots,\ell_{k-1})$:
\begin{eqnarray*}
\tau(p(i_0;i_1,\ell_1;\cdots;i_{k-1},\ell_{k-1})) & \xrightarrow{\text{tangent}}  &\tau(p(i_0;i_1,\ell_1;\cdots;i_{k-1},\ell_{k-1};i_k,b))\\
& \xrightarrow{E_{i_k}}  & \tau(p(i_1;i_2,\ell_2;\cdots;i_{k-1},\ell_{k-1};i_k,b))\\
& \xrightarrow{\text{tangent}}  & \tau(p(i_0;i_1,\ell_1;\cdots;i_{k-1},\ell_{k-1})).
\end{eqnarray*}

\item The type (2) loop at level $k$, denoted by $T^k(i_0,i_1,\cdots,i_{k+1};\ell_1,\cdots,\ell_k)$:
\begin{eqnarray*}
\tau(p(i_1;i_2,\ell_2;\cdots;i_k,\ell_k)) & \xrightarrow{E_{i_k}}  &\tau(p(i_0;i_1,\ell_1;\cdots;i_k,\ell_k))\\
& \xrightarrow{\text{tangent}}  & \tau(p(i_0;i_1,\ell_1;\cdots;i_{k+1},\ell_k+1))\\
& \xrightarrow{E_{i_{k+1}}}  & \tau(p(i_1;i_2,\ell_2;\cdots;i_{k+1},\ell_k+1))\\
& \xrightarrow{\text{tangent}}  & \tau(p(i_1;i_2,\ell_2;\cdots;i_k,\ell_k)).
\end{eqnarray*}
\end{enumerate}
}
\end{defn}

We then have the analog of  Proposition~\ref{prop-iota3}:
\begin{prop}\label{prop-iotaellk}
The inclusion $\fN_k\subset \fN_{k-1}$ induces a map
$$\iota_k:H_1(\fN_k;\mZ) \to H_1(\fN_{k-1};\mZ)$$
such that:
\begin{enumerate}
 \item $\iota_k(P^k_1(i_0))=(\ell(\delta_k)-\ell(\delta_{k-1}))[R]+ P^{k-1}_1(i_0)$;

 \item For $2\leq n\leq k-2$, $\iota_k(P^k_n(i_0,\cdots,i_{n-1};\ell_1,\cdots,\ell_{n-1}))$ is homologous to
 $$P^{k-1}_n(i_0,\cdots,i_{n-1};\ell_1,\cdots,\ell_{n-1}).$$
 
\item For level $k-1$, $\iota_k(P^k_{k-1}(i_0,\cdots,i_{k-2};\ell_1,\cdots,\ell_{k-2}))$ is homologous to
\begin{align*}
T^{k-1}(i_0,\cdots,i_{k-1},i_{k-1};\ell_1,\cdots,\ell_{k-1}-1) &+ T^{k-1}(i_0,\cdots,i_{k-1},i_{k-1};\ell_1,\cdots,\ell_{k-1}-2)+\cdots\\
 &+ T^{k-1}(i_0,\cdots,i_{k-1},i_{k-1};\ell_1,\cdots,b)\\
 &+ E^{k-1}(i_0,\cdots,i_{k-1};\ell_1,\cdots,\ell_{k-2}),
 \end{align*}
 for any $i_{k-1}$ and the minimum possible value of $\ell_{k-1}$.

\item $\iota_k(S_1^k(i_0,i_1;\ell_1))=-P_1^{k-1}(i_0)+P_1^{k-1}(i_2)+[R]$ for $i_2\neq i_0$.

\item For $2\leq n\leq k-1$, the image  of $S_n^k(i_0,\cdots,i_n;\ell_1,\cdots,\ell_n)$ is homologous to the image under $\iota_k$ of 
$$-\sum_{m=1}^n P_m^k(i_0,\cdots,i_{m-1};\ell_1, \cdots,\ell_{m-1}) +\sum_{m=1}^n P_m^k(i_{n+1},i_1,\cdots,i_{m-1};\ell_1, \cdots,\ell_{m-1}) +[R],$$
for $i_{n+1}\neq i_0$ and modulo adding some type (2) loops at level $k-1$.

\item For the exceptional loops we have
\begin{align*}
\iota_k(E^k(i_0,\cdots,i_k;\ell_1,\cdots,\ell_{k-1})) &= T^{k-1}(i_0,\cdots,i_{k-1},i_{k-1};\ell_1,\cdots,\ell_{k-1}-1)\\
&+ T^{k-1}(i_0,\cdots,i_{k-1},i_{k-1};\ell_1,\cdots,\ell_{k-1}-2)+\cdots\\
 & + T^{k-1}(i_0,\cdots,i_{k-1},i_{k-1};\ell_1,\cdots,b)\\
 &+ E^{k-1}(i_0,\cdots,i_{k-1};\ell_1,\cdots,\ell_{k-2}).
\end{align*}

 \item $\iota_k(T^k(i_0,i_1,\cdots,i_{k+1};\ell_1,\cdots,\ell_k))$ is the trivial element.
 \end{enumerate}
\end{prop}

Computing the exact rank of the groups $H_1(\fN_k;\mZ)$ seems an impossible task, since the values of the $\ell$-indices for which the different types of loops exist, and are not trivial, does not appear to follow any simple pattern. In contrast, we can easily give a lower bound. 

\begin{cor}\label{cor-rankh1}
The rank of $H_1(\fN_k;\mZ)$ is at least $2^{k+2}-1$.
\end{cor}

\proof
We have three distinctive elements in $H_1(\fN_k;\mZ)$: the classes $[R]$, $P_1^k(1)$ and $P_1^k(2)$. Then for each $2\leq n\leq k-1$ and each combination $(i_0,i_1,\cdots, i_{n-1})$ of 1's and 2's,  there is at least one level $n$ loop. Thus we obtain for each $n$, $2^n$ elements in $H_1(\fN_k;\mZ)$. Observe that the actual number of generators  is greater than this.  For example, for $n=2$,  given $(i_0,i_1)$ then for any $\ell_1$ such that $s(\delta_k;i_0,i_1)\leq \ell_1< \ell(\delta_k;i_0,i_1)$ there is a level 2 loop and we are just counting 1 of these. 

The number of exceptional loops is given by the sum of 
$$2\left(\ell(\delta_k;i_0;i_1,\ell_1;\cdots;i_{k-1},\ell_{k-2};i_{k-1})- s(\delta_k;i_0;i_1,\ell_1;\cdots;i_{k-1},\ell_{k-2};i_{k-1})\right)$$ 
over all possible combinations of these indices. Again, counting  two generators ($i_k=1,2$) for each combination $(i_0,i_1,\cdots,i_{k-1})$ of 1's and 2's, we conclude that there are at least $2^{k+1}$ loops. Analogously, there are at least $2^k$ type (2) loops at level $k$.

We obtain that the rank of $H_1(\fN_k;\mZ)$ is lower bounded by $1+2+\cdots+2^{k+1}=2^{k+2}-1$. Observe that we are not counting the one turn loops at any level.
\endproof

The objective is to build a shape approximation $\fU=\{U_\ell \,|\, \ell=1,2,\ldots\}$ of $\fM$ satisfying the hypothesis of Proposition~\ref{prop-notstable}.  That is, we require  that for $k > 0$: 
\begin{itemize}
\item the rank of $H_1(U_k;\mZ)\geq 2^k$,
\item for all $\ell > k$ the rank of the image $H_1(U_{\ell} ; \mZ) \to H_1(U_k ; \mZ)$ is 3.
\end{itemize} 
Observe that taking a sequence $\e_k$ of sufficiently small positive numbers and $U_k=U(\fN_k,\e_k)$, then Lemma~\ref{lem-shapek} and Corollary~\ref{cor-rankh1} imply the first condition. In order to satisfy the second condition, we will extract a subsequence of $\fN_{n_k}$ for which the rank of the image $H_1(\fN_{n_{k+1}} ; \mZ) \to H_1(\fN_{n_k} ; \mZ)$ is 3. This follows from the following result.

  The rank of the image of the map $\iota_\ell^0:H_1(\fN_\ell;\mZ)\to H_1(\fN_0;\mZ)$ is 3, since the generators $[R]$, $[b_1]$ and $[b_2]$ are in the image and they generate $H_1(\fN_0;\mZ)$.
  
As in Remark~\ref{rmk-generators}, the first homology group of $C_\mK^+(\fM_0^1,\delta_k)$ is generated by the loops $[R]$, $P_1^k(1)$, $P_1^k(2)$ and the level 1 one turn loops $S_1^k(i_0,i_1;\ell_1)$. 
 Since the rank of $H_1(\fN_0;\mZ)$ is 3,  what the proposition asserts is that given $k$ there is a number $\ell$ big enough, such that all the homology classes in $\fN_\ell$ are homologous inside $\fN_k$ to loops contained in the $\delta_k$-neighborhood of $\fM_0^1$, and even more in the group generated by $[R]$, $P_1^k(1)$ and $P_1^k(2)$. The idea behind the proof is that if there is a loop in $\fN_\ell$, its image in $\fN_k$ is a loop that travels along propellers of level at must $k$. These propellers have to be sufficiently near the core cylinder $\tau(\cR')$ to guarantee that they branch up to level $\ell$. Thus the image of the loop in $\fN_k$ has to be close to levels 0 and 1, and then are homologous to a combination of the loops  $[R]$, $P_1^k(1)$, $P_1^k(2)$. In order to make this assertion precise we need to consider the images of the exceptional loops and the loops at levels $n$ for any $2\leq n\leq \ell-1$.

\begin{prop}\label{prop-rankkn}
Given $k\geq 2$  there exists $\ell>k$ such that the image of the map
\begin{equation}
\iota_\ell^k:H_1(\fN_\ell;\mZ)\to H_1(\fN_k;\mZ)
\end{equation}
has rank equal to 3.
\end{prop}

 \proof

Fix $k\geq 2$ and consider the subgroup $G_k$ of $H_1(\fN_k;\mZ)$ generated by $[R]$, $P_1^k(1)$ and $P_1^k(2)$. 
We will find $\ell \gg k$ such that the image of $\iota_\ell^k$ is $G_k$. 

In the following arguments,   all the homotopies will be   contained in    $\fN_k$.

First, observe that   case (1) of Proposition~\ref{prop-iotaellk} shows  that the map $\ds \iota_{i} \colon H_1(\fN_{i};\mZ) \to H_1(\fN_{i-1};\mZ)$
satisfies for each $i > k$, 
$$\iota_i^{i-1}(P^i_1(i_0))=(\ell(\delta_i)-\ell(\delta_{i-1}))[R]+ P^{i-1}_1(i_0).$$
It then follows by induction, that for any $\ell > k$,         the images of the homology classes defined by the level 1 loops, $P_1^\ell(i_0)$ for $i_0=1,2$,  are in $G_k$.
Thus, the subgroup 
$G_k$ is contained in the image of $\iota_\ell^k$.

  \bigskip

Next consider   the   one turn loops  $S_n^k(i_0,\cdots,i_n;\ell_1,\cdots,\ell_n)$ defined in case (3) of Definitions~\ref{def-generators2}, \ref{def-generators3} and \ref{def-generatorsk}.  Case (4) of Proposition~\ref{prop-iotaellk} gives a reductive procedure for reducing the homology classes of the 1 turn loops $S_n^k(i_0,\cdots,i_n;\ell_1,\cdots,\ell_n)$ to some combination of $[R]$ and the classes  defined by level n loops, for some  $1 \leq n < \ell$. 
Thus,  the one turn loops   
do not contribute to the rank of the image of $\iota_\ell^k$ for any $\ell>k$. 

 \bigskip

 The Type (2) loops at level $\ell$ were introduced in    case (5) of Definitions~\ref{def-generators2}, \ref{def-generators3} and \ref{def-generatorsk}. For each $\ell > 1$, it was shown that the images of the homology classes of these loops under the map  $\ds \iota_{\ell} \colon H_1(\fN_{\ell};\mZ) \to H_1(\fN_{\ell -1};\mZ)$ are trivial. 
 Thus, for all $\ell > k$, their   images  by  $\ds \iota_{\ell}^k \colon H_1(\fN_{\ell};\mZ) \to H_1(\fN_{k};\mZ)$ are trivial.

 \bigskip

Next, we consider      the level $n$ loops  $P^\ell_n(i_0,\cdots,i_{n-1};\ell_1,\cdots,\ell_{n-1})$ for $2\leq n\leq \ell-1$, as defined in cases (2b) of 
Definitions~\ref{def-generators3} and \ref{def-generatorsk}.
These are loops defined by paths that go out, possibly a ``long'' distance, along a propeller at level $n$, then connect to a path along propellers at a   level less than $n$, and return along   the shortest path tangent to $\fM_0$.

Fix $n$ between 2 and $\ell-1$, and assume without loss of generality that $i_0=1$.   Let $P$ be a propeller satisfying the two conditions: the tip of the level $n$ propeller $P = \tau(P_{\g(i_1,\ell_1;\cdots;i_{n-1},\ell_{n-1})})$ is contained in $C_\mK^+(\fM_0^{n-1},\delta_\ell)$, and $P$ does not retract to this set.  Note that these conditions  depend on the value of $\delta_\ell$. Then the  loop $P^\ell_n(i_0,\cdots,i_{n-1};\ell_1,\cdots,\ell_{n-1})$ exists.

Recall from Lemma \ref{lem-paccumulationpoints} that level $n$ propellers accumulate (as the first $\ell$-index goes to infinity) on level $n-1$ propellers, that respectively accumulate on level $n-2$ propellers, and so forth. Assume for a moment that $\ell>k$ is fixed, with the level $n$ propeller as above.  Then there exists a positive number $\beta_\ell>\delta_\ell$ such that  the propeller $P$ is contained in $C_\mK^+(\fM_0^1,\beta_\ell)$. Moreover, as $\ell$ tends to infinity, $\beta_\ell$ decreases    to zero. 
  Given this observation,  we claim that  for $k$ fixed, there exists $\ell>k$ such that  any propeller $P$ as above is contained in $C_\mK^+(\fM_0^1,\delta_k)$, as it suffices to take any $\ell > k$ such that $\beta_\ell \leq \delta_k$. 
  
  Recall that  $P^\ell_n(i_0,\cdots,i_{n-1};\ell_1,\cdots,\ell_{n-1})$, for $2\leq n\leq \ell-1$ is formed by the concatenation of the paths
  \begin{eqnarray}
  \tau(p(i_0;i_1,\ell_1;\cdots;i_{n-1},\ell_{n-1})) &\xrightarrow{\text{tangent}} & \tau(p(i_0;i_1,\ell_1;\cdots;i_{n-1},\ell_{n-1};i_n,\ell_n)) \label{eq-block0-line1}\\
  & \xrightarrow{E_{i_n}} & \tau(p(i_1;i_2,\ell_2;\cdots ; i_n,\ell_n)) \label{eq-block0-line2}\\
  &\xrightarrow{\text{tangent}} & p(i_0;i_1,\ell_1;\cdots;i_{n-1},\ell_{n-1})). \label{eq-block0-line3}
  \end{eqnarray}
 The tangent paths \eqref{eq-block0-line1} and \eqref{eq-block0-line3} are   contained in $C_\mK^+(\fM_0^1,\delta_k)$, and thus are homotopic  to tangent   paths contained in propellers at levels 0 and 1.  Thus after homotopy, we obtain the loop
  \begin{equation} \label{eq-block0-line4}
\tau(p(i_{n-1})) \xrightarrow{\text{tangent}} \tau(p(i_{n-1};i_n,\ell_n)) \xrightarrow{E_{i_n}} \tau(p(i_n)) \xrightarrow{\text{tangent}} \tau(p(i_{n-1})),
\end{equation}
  that is homologous to a class in $G_k$.

 \bigskip
  
Now  consider   the image of the homology classes defined by the exceptional loops  $E^\ell(i_0,\cdots,i_\ell;\ell_1,\cdots;\ell_{\ell-1})$  which were   introduced in    case (4) of Definitions~\ref{def-generators2}, \ref{def-generators3} and \ref{def-generatorsk}. 
   Recall that  such a loop   can be described as the composition of paths: 
 \begin{eqnarray}
\tau(p(i_0;i_1,\ell_1;\cdots;i_{\ell-1},\ell_{\ell-1})) & \xrightarrow{\text{tangent($\ell$)}}  &\tau(p(i_0;i_1,\ell_1;\cdots;i_{\ell-1},\ell_{\ell-1};i_\ell,b)) \label{eq-block1-line1}\\
& \xrightarrow{E_{i_\ell}}  & \tau(p(i_1;i_2,\ell_2;\cdots;i_{\ell-1},\ell_{\ell-1};i_\ell,b)) \label{eq-block1-line2}\\
& \xrightarrow{\text{tangent($<\ell$)}}  & \tau(p(i_0;i_1,\ell_1;\cdots;i_{\ell-1},\ell_{\ell-1})). \label{eq-block1-line3}
\end{eqnarray}
  Figure~\ref{fig:Figure56a} illustrates such a path for the case $k=2$.  The path in \eqref{eq-block1-line1} starts at a point  on a propeller at level $\ell -1$, entering a propeller at level $\ell$ and continuing to a point   which is the first entry point into a filled propeller. Next,  \eqref{eq-block1-line2} is a short path in the face $E_{i_\ell}$ to a point   at lower level. It then follows the path \eqref{eq-block1-line3} along propellers at level less than $\ell$ back to the starting point. A key point, as seen in the analysis  below, is that while the   path \eqref{eq-block1-line3} at lower level may be ``short'' as in Figure~\ref{fig:Figure56a},   it may just as well traverse the tree $\TP$ from one extreme to another and not be ``short''.

The level $k$ propellers are filled in   $\fN_k$, hence the tangent path segments in \eqref{eq-block1-line1} and \eqref{eq-block1-line3} at level greater or equal to $k$ are homotopic to   tangent segments at level $k$. First, the tangent segment  \eqref{eq-block1-line1}  is homotopic to 
\begin{equation}\label{eq-block2-line1}
\tau(p(i_q;i_{q+1},\ell_{q+1};\cdots;i_{\ell-1},\ell_{\ell-1}))  \xrightarrow{\text{tangent($k$)}}  \tau(p(i_q;i_{q+1},\ell_{q+1};\cdots;i_{\ell-1},\ell_{\ell-1};i_\ell,b)),
\end{equation}
for $q=\ell-k$. 
Note that the result \eqref{eq-block2-line1} illustrates an important aspect of the homotopy operation moving a path segment from a level $\ell$ propeller to one at a lower level.  
This results in the  elimination of the initial stages of the labeling; that is,   the initial indices $(i_0;i_1,\ell_1;\cdots;i_{q-1},\ell_{q-1})$ are deleted   in \eqref{eq-block2-line1}.

Next, the   second tangent segment \eqref{eq-block1-line3}, which is a path
\begin{equation}\label{eq-block4-line1}
\tau(p(i_1;i_2,\ell_2;\cdots;i_{\ell-1},\ell_{\ell-1};i_\ell,b)) \xrightarrow{\text{tangent($<\ell$)}}   \tau(p(i_0;i_1,\ell_1;\cdots;i_{\ell-1},\ell_{\ell-1})),
\end{equation}
can be decomposed as a concatenation of segments tangent to propellers at level less than $\ell$.  This path travels through the tree $\TP$    to connect the two level $\ell -1$ endpoints in $\fM_0^{\ell-1}$.  We denote the paths at each level as follows, where we note that now the final indices are first being removed   in lines \eqref{eq-block5-line1} to \eqref{eq-block5-line4}, as the path  goes to   points at successively lower levels. The path then travels along the set $\fM_0^{k-1}$ as indicted in \eqref{eq-block5-line4}. It    then continues back up the levels to the path in \eqref{eq-block5-line6}. This process yields the following path which is equivalent to   the path   \eqref{eq-block4-line1}.

\begin{eqnarray}
\tau(p(i_1;i_2,\ell_2;\cdots; i_\ell,b)) & \xrightarrow{\text{tangent($\ell-1$)}} & \tau(p(i_1;i_2,\ell_2;\cdots;i_{\ell-1},\ell_{\ell-1})) \label{eq-block5-line1}\\
& \xrightarrow{\text{tangent($\ell-2$)}} & \tau(p(i_1;i_2,\ell_2;\cdots;i_{\ell-2},\ell_{\ell-2})) \label{eq-block5-line2}\\
& \cdots & \nonumber\\
& \xrightarrow{\text{tangent($k=\ell-q$)}} & \tau(p(i_1;i_2,\ell_2;\cdots;i_{k},\ell_{k})) \label{eq-block5-line3}\\
& \xrightarrow{\text{tangent($<k$)}} & \tau(p(i_0;i_1,\ell_1;\cdots;i_{k-1},\ell_{k-1})) \label{eq-block5-line4}\\
& \xrightarrow{\text{tangent($k$)}} & \tau(p(i_0;i_1,\ell_1;\cdots;i_k,\ell_k)) \label{eq-block5-line5}\\
& \cdots& \nonumber\\
& \xrightarrow{\text{tangent($\ell-1$)}} & \tau(p(i_0;i_1,\ell_1;\cdots;i_{\ell-1},\ell_{\ell-1})). \label{eq-block5-line6}
\end{eqnarray}
Note that the reduction in level process stops with the paths in level $k$, as   we can no longer move the curves by a homotopy through the unfilled propellers at level less than $k$. Instead, the curve simply follows the propellers in $\fM_0^k$, as indicated by the path in  \eqref{eq-block5-line4}.

Replacing each of the tangential segments above, \eqref{eq-block5-line1} to \eqref{eq-block5-line6},  with a segment of level at most $k=\ell-q$, we obtain:
\begin{eqnarray}
\tau(p(i_1;i_2,\ell_2;\cdots; i_\ell,b)) & \xrightarrow{E_{i_\ell}} & \tau(p(i_q;i_{q+1},\ell_{q+1};\cdots; i_\ell,b)) \label{eq-block6-line1}\\
& \xrightarrow{\text{tangent($k$)}} & \tau(p(i_q;i_{q+1},\ell_{q+1};\cdots;i_{\ell-1},\ell_{\ell-1})) \label{eq-block6-line2}\\
& \xrightarrow{E_{i_{\ell-1}}} & \tau(p(i_{q-1};i_q,\ell_q;\cdots;i_{\ell-1},\ell_{\ell-1})) \label{eq-block6-line3}\\
& \xrightarrow{\text{tangent($k$)}} & \tau(p(i_{q-1};i_q,\ell_q;\cdots;i_{\ell-2},\ell_{\ell-2})) \label{eq-block6-line4}\\
& \cdots & \nonumber\\
& \xrightarrow{\text{tangent($k$)}} & \tau(p(i_1;i_2,\ell_2;\cdots;i_{k},\ell_{k})) \label{eq-block6-line5}\\
& \xrightarrow{\text{tangent($<k$)}} & \tau(p(i_0;i_1,\ell_1;\cdots;i_{k-1},\ell_{k-1})) \label{eq-block6-line6}\\
& \xrightarrow{\text{tangent($k$)}} & \tau(p(i_0;i_1,\ell_1;\cdots;i_k,\ell_k)) \label{eq-block6-line7}\\
& \xrightarrow{E_{i_k}} & \tau(p(i_1;i_2,\ell_2;\cdots;i_k,\ell_k)) \label{eq-block6-line8}\\
& \cdots & \nonumber\\
& \xrightarrow{\text{tangent($k$)}} & \tau(p(i_{q-1};i_{q},\ell_{q};\cdots;i_{\ell-1},\ell_{\ell-1})). \label{eq-block6-line9}
\end{eqnarray}
As before,   the homotopy operation eliminates the initial indices labeling the endpoints of the path  in $\fM_0^{\ell}$ as the path is moved to a tangential path in a lower level of  $\fM_0^k$.

Now observe that if we follow the paths and homotopies indicated  in the lines
 \eqref{eq-block2-line1} which is homotopic to \eqref{eq-block1-line1}, followed by \eqref{eq-block1-line2}, followed by  \eqref{eq-block6-line1}, followed by  \eqref{eq-block6-line2}, we obtain a loop:
 \begin{eqnarray*}
\tau(p(i_q;i_{q+1},\ell_{q+1};\cdots;i_{\ell-1},\ell_{\ell-1}))  & \xrightarrow{\text{tangent($k$)}} &  \tau(p(i_q;i_{q+1},\ell_{q+1};\cdots;i_{\ell-1},\ell_{\ell-1};i_\ell,b))\\
  &  \xrightarrow{E_{i_\ell}}   & \tau(p(i_1;i_2,\ell_2;\cdots;i_{\ell-1},\ell_{\ell-1};i_\ell,b)) \\
 &  \xrightarrow{E_{i_\ell}} & \tau(p(i_q;i_{q+1},\ell_{q+1};\cdots; i_\ell,b))  \\
   & \xrightarrow{\text{tangent($k$)}} &  \tau(p(i_q;i_{q+1},\ell_{q+1};\cdots;i_{\ell-1},\ell_{\ell-1})) 
\end{eqnarray*}
 which is homotopic to the trivial loop, as the first and last tangential paths are inverses of each other.

Thus, the image in $\fN^k$ of the exceptional loop   $\ds E^\ell(i_0,\cdots,i_\ell;\ell_1,\cdots;\ell_{\ell-1})$ defined by the paths in \eqref{eq-block1-line1}, \eqref{eq-block1-line2} and \eqref{eq-block1-line3} is homotopic to   the loop that starts at the endpoint of   \eqref{eq-block6-line3}, then follows \eqref{eq-block6-line4} up to \eqref{eq-block6-line9}, as indicated in the following:

\begin{eqnarray}
\tau(p(i_{q-1};i_q,\ell_q;\cdots;i_{\ell-1},\ell_{\ell-1})) & \xrightarrow{\text{tangent($k$)}} & \tau(p(i_{q-1};i_q,\ell_q;\cdots;i_{\ell-2},\ell_{\ell-2})) \label{eq-block8-line1}\\
& \cdots & \nonumber\\
& \xrightarrow{\text{tangent($k$)}} & \tau(p(i_1;i_2,\ell_2;\cdots;i_{k},\ell_{k})) \label{eq-block8-line2}\\
& \xrightarrow{\text{tangent($<k$)}} & \tau(p(i_0;i_1,\ell_1;\cdots;i_{k-1},\ell_{k-1})) \label{eq-block8-line3}\\
& \xrightarrow{\text{tangent($k$)}} & \tau(p(i_0;i_1,\ell_1;\cdots;i_k,\ell_k)) \label{eq-block8-line4}\\
& \xrightarrow{E_{i_k}} & \tau(p(i_1;i_2,\ell_2;\cdots;i_k,\ell_k)) \label{eq-block8-line5}\\
& \cdots &\nonumber\\
& \xrightarrow{\text{tangent($k$)}} & \tau(p(i_{q-1};i_{q},\ell_{q};\cdots;i_{\ell-1},\ell_{\ell-1})). \label{eq-block8-line6}
\end{eqnarray}

Now we express the above loop defined by the paths in \eqref{eq-block8-line1} to \eqref{eq-block8-line6} in more detail.
 We   rewrite this loop, starting with the segment tangent in \eqref{eq-block8-line3} to $\fM_0^{k-1}$:
 
 \begin{eqnarray}
 \tau(p(i_1;i_2,\ell_2;\cdots;i_k,\ell_k)) & \xrightarrow{\text{tangent($<k$)}} &  \tau(p(i_0;i_1,
 \ell_1;\cdots;i_{k-1},\ell_{k-1})) \label{eq-cancel1}\\
 & \xrightarrow{\text{tangent($k$)}} &  \tau(p(i_0;i_1,\ell_1;\cdots;i_k,\ell_k))  \label{eq-cancel2}\\
 & \xrightarrow{E_{i_k}} &  \tau(p(i_1;i_2, \ell_2;\cdots;i_k,\ell_k))  \label{eq-cancel3}\\
  & \xrightarrow{\text{tangent($k$)}} &  \tau(p(i_1;i_2,\ell_2;\cdots;i_{k+1},\ell_{k+1}))  \label{eq-cancel4}\\
 & \xrightarrow{E_{i_{k+1}}} &  \tau(p(i_2;i_3, \ell_3;\cdots;i_{k+1},\ell_{k+1}))  \label{eq-cancel5}\\
 & \vdots & \nonumber\\
 & \xrightarrow{E_{i_{\ell-2}}} & \tau(p(i_{q-1};i_q,\ell_q;\cdots;i_{\ell-2},\ell_{\ell-2}))\label{eq-cancel6}\\
 & \xrightarrow{\text{tangent($k$)}} & \tau(p(i_{q-1};i_q,\ell_q;\cdots;i_{\ell-1},\ell_{\ell-1})) \label{eq-cancel7}\\
 & \xrightarrow{\text{tangent($k$)}} & \tau(p(i_{q-1};i_q,\ell_q;\cdots;i_{\ell-2},\ell_{\ell-2})) \label{eq-cancel8}\\
 & \xrightarrow{E_{i_{\ell-2}}} & \tau(p(i_{q-2};i_{q-1},\ell_{q-1};\cdots;i_{\ell-2},\ell_{\ell - 2})) \label{eq-cancel9}\\
 & \vdots & \nonumber\\
 & \xrightarrow{E_{i_{k+1}}} & \tau(p(i_1;i_2,\ell_1;\cdots;i_{k+1},\ell_{k+1})) \label{eq-cancel10}\\
 & \xrightarrow{\text{tangent($k$)}} & \tau(p(i_1;i_2,\ell_2;\cdots;i_k,\ell_k)). \label{eq-cancel11}
 \end{eqnarray}

 Observe that the tangent parts in lines \eqref{eq-cancel2} to \eqref{eq-cancel7} go from a point at level $k-1$ to a point at level $k$, while the tangent parts in lines \eqref{eq-cancel8} to \eqref{eq-cancel11} go from a point at level $k$ to a point at level $k-1$

 The arrows in lines \eqref{eq-cancel7} and \eqref{eq-cancel8} represent paths that are tangent to the same propeller, have the same endpoints and opposite directions, hence their concatenation is homotopic to the trivial loop.

 Also notice that the  arrows in lines \eqref{eq-cancel6} and \eqref{eq-cancel9} represent paths in $E_{i_{\ell-2}}$ with the same endpoints and opposite directions, hence their concatenation is again homotopic to the trivial loop. 
 
 The process continues cancelling pairs of arrows in the representation above, until \eqref{eq-cancel10} cancels with \eqref{eq-cancel5} and \eqref{eq-cancel11} cancels with \eqref{eq-cancel4}.  Hence we are left with the loop:
\begin{eqnarray}
 \tau(p(i_1;i_2,\ell_2;\cdots;i_k,\ell_k)) & \xrightarrow{\text{tangent($<k$)}} &  \tau(p(i_0;i_1,\ell_1;\cdots;i_{k-1},\ell_{k-1})) \label{eq-block9-line21}\\
 & \xrightarrow{\text{tangent($k$)}} &  \tau(p(i_0;i_1,\ell_1;\cdots;i_k,\ell_k)) \label{eq-block9-line22} \\
 & \xrightarrow{E_{i_k}} &  \tau(p(i_1;i_2, \ell_2;\cdots;i_k,\ell_k)). \label{eq-block9-line23}
\end{eqnarray}
The loop defined by the paths in \eqref{eq-block9-line21}, \eqref{eq-block9-line22} and \eqref{eq-block9-line23},  can then be re-arranged (by moving the first arrow to be last)  as the concatenation 
\begin{eqnarray}
\tau(p(i_0;i_1,\ell_1;\cdots;i_{k-1},\ell_{k-1})) & \xrightarrow{\text{tangent($k$)}}  & \tau(p(i_0;i_1,\ell_1;\cdots;i_k,\ell_k)) \label{eq-block9-line1}\\
& \xrightarrow{E_{i_k}} & \tau(p(i_1;i_2,\ell_2;\cdots;i_k,\ell_k)) \label{eq-block9-line2}\\
& \xrightarrow{\text{tangent($<k$)}} & \tau(p(i_0;i_1,\ell_1;\cdots;i_{k-1},\ell_{k-1})). \label{eq-block9-line3}
\end{eqnarray}
Note that the first maps follows a level $k$ propeller around the cylinder, until the point $\tau(p(i_0;i_1,\ell_1;\cdots;i_k,\ell_k))$, then it closes by going down one level in the face $E_{i_k}$, then it closes up by following the tree back at levels less than $k$.

The homology class of the loop defined by \eqref{eq-block9-line1} to \eqref{eq-block9-line2} to \eqref{eq-block9-line3} then equals the class   $$\cE(i_0,\cdots, i_k;\ell_1,\cdots, \ell_k)  \in H_1(\fN_k;\mZ), $$    defined by
\begin{align*}
\cE(i_0,\cdots, i_k;\ell_1,\cdots, \ell_k)= & ~ T^k(i_0,\cdots, i_k,i_k;\ell_1,\cdots,\ell_{k-1},\ell_k-1)+T^k(i_0,\cdots, i_k,i_k;\ell_1,\cdots,\ell_{k-1},\ell_k-2)\\
& +\cdots +T^k(i_0,\cdots, i_k,i_k;\ell_1,\cdots,\ell_{k-1},b) +E^k(i_0,\cdots,i_k;\ell_1,\cdots,\ell_{k-1}).
\end{align*}
It follows that the homology class   in $H_1(\fN_k;\mZ)$  of the exceptional loop   $\ds E^\ell(i_0,\cdots,i_\ell;\ell_1,\cdots;\ell_{\ell-1})$  is equal to $\cE(i_0,\cdots, i_k;\ell_1,\cdots, \ell_k)  \in H_1(\fN_k;\mZ)$. 
Observe that the class $\cE(i_0,\cdots, i_k;\ell_1,\cdots, \ell_k)$ depends on $i_0$  and the first $k$ of the double indices defining the loop
$E^\ell(i_0,\cdots,i_\ell;\ell_1,\cdots,\ell_{\ell-1})$, but is independent of  the indices above $k$.

For simplicity, now assume   that $i_0=1$, then the loop defined by \eqref{eq-block9-line1} to \eqref{eq-block9-line2} to \eqref{eq-block9-line3} becomes the level $k$ loop 
\begin{eqnarray}
\tau(p(1;i_1,\ell_1;\cdots;i_{k-1},\ell_{k-1})) & \xrightarrow{\text{tangent($k$)}}  & \tau(p(1;i_1,\ell_1;\cdots;i_k,\ell_k)) \label{eq-block10-line1}\\
& \xrightarrow{E_{i_k}} & \tau(p(i_1;i_2,\ell_2;\cdots;i_k,\ell_k)) \label{eq-block10-line2}\\
& \xrightarrow{\text{tangent($<k$)}} & \tau(p(1;i_1,\ell_1;\cdots;i_{k-1},\ell_{k-1})),  \label{eq-block10-line3}
\end{eqnarray}
where the tangent part at level $k$ is along the propeller $\tau(P_{\g(i_1,\ell_1;\cdots;i_k,\ell_k)})$. 

 Recall that we are assuming that the index $\ell$ is much larger than $k$.
 
  If we assume that $\ell_k>\ell(\delta_k;i_0;\cdots;i_{k-2},\ell_{k-2};i_{k-1})$, then 
  then the loop defined by concatenating the path in  \eqref{eq-block10-line1}, followed by \eqref{eq-block10-line2}, followed by \eqref{eq-block10-line3},      is then homotopic to a loop in $C_\mK^+(\fM_0^{k-1},\delta_k)$.
  
  Now proceed  recursively. If we assume that  $\ell_{k-1}>\ell(\delta_k;i_0;\cdots;i_{k-3},\ell_{k-3};i_{k-2})$, then the result of the homotopy of the above loop into $C_\mK^+(\fM_0^{k-1},\delta_k)$ is then  homotopic to a loop in $C_\mK^+(\fM_0^{k-2},\delta_k)$.
  
Continuing in this way,  assuming that all  the $\ell$-indices are large enough, we obtain a homotopy of the loop defined by \eqref{eq-block10-line1},   \eqref{eq-block10-line2} and  \eqref{eq-block10-line3}, to   a loop in $C_\mK^+(\fM_0^1,\delta_k)$. 
Let   $\widetilde{E}(i_0,\cdots, i_k;\ell_1,\cdots, \ell_k)$ denote the loop  contained in $C_\mK^+(\fM_0^1,\delta_k)$ which results from the above homotopies.
  Then for $\ell$ sufficiently large,  the homology class  $\cE(i_0,\cdots, i_k;\ell_1,\cdots, \ell_k)$ is defined by the loop $\widetilde{E}(i_0,\cdots, i_k;\ell_1,\cdots, \ell_k)$.

We now claim, that  for $\ell$ sufficiently large, the homology class  $\cE(i_0,\cdots, i_k;\ell_1,\cdots, \ell_k) \in G_k$.

Observe that the endpoint of line  \eqref{eq-block10-line1}, that is the point $\tau(p(1;i_1,\ell_1;i_2.\ell_2;\cdots;i_k,\ell_k))$, can be connected to the level 1 point $\tau(p(i_{k-1};i_k,\ell_k))$ by a path in $E_{i_k}$ that is inside the set $C_\mK^+(\fM_0^1,\delta_k)$. 

Recall that as $\ell_k$ increases to infinity,  the level 1 points $\tau(p(i_{k-1};i_k,\ell_k))$ approach the level 0 point $\tau(p(i_k))$. 

As the index $\ell$ is increased,  the index $\ell_k$ also must increase. We may  thus assume that $\ell_k$ is such that $\tau(p(1;i_1,\ell_1;i_2.\ell_2;\cdots;i_k,\ell_k))$ can be connected to the level 0 point $\tau(p(i_k))$ by a path in $E_{i_k}\cap C_\mK^+(\fM_0^1,\delta_k)$. We choose a value for  $\ell$ so that this condition is satified.

Since the   loop $\widetilde{E}(i_0,\cdots, i_k;\ell_1,\cdots, \ell_k)$  is contained in $C_\mK^+(\fM_0^1,\delta_k)$, the path in line \eqref{eq-block10-line1} can be homotoped to a path tangent at level 1. 
The choice of $\ell$ sufficiently large then implies that $\widetilde{E}(i_0,\cdots, i_k;\ell_1,\cdots, \ell_k)$ is homotopic to the loop
\begin{eqnarray}
\tau(p(i_{k-1})) & \xrightarrow{\text{tangent($1$)}}  & \tau(p(i_{k-1};i_k,\ell_k)) \label{eq-block11-line1}\\
& \xrightarrow{E_{i_k}} & \tau(p(i_k)) \label{eq-block11-line2}\\
& \xrightarrow{\text{tangent($0$)}} & \tau(p(i_{k-1})).  \label{eq-block11-line3}
\end{eqnarray}
Thus $\widetilde{E}(i_0,\cdots, i_k;\ell_1,\cdots, \ell_k)$ is homotopic to a loop that travels along a level 1 propeller as in line \eqref{eq-block11-line1}, then jumps back to the Reeb cylinder  as in line \eqref{eq-block11-line2} and then closes by a path contained in the Reeb cylinder  as in line \eqref{eq-block11-line3}. 

Observe that since line \eqref{eq-block11-line3} is homotopic to line \eqref{eq-block10-line3}, it might be more complicated than the shortest path from $\tau(p(i_{k-1}))$ to $\tau(p(i_k))$. 
In any case, the path represented by lines \eqref{eq-block11-line1} to \eqref{eq-block11-line3} can be written as the sum $P_1^k(i_{k-1})$ and a certain multiple of $[R]$, thus we get a loop in $G_k$.
 \endproof

 To define the shape approximation $\fU = \{U_k \mid k=1,2,\ldots\}$  satisfying the hypothesis of Proposition~\ref{prop-notstable}, let $U_1=U(\fN_2,\e_1)$ for $\e_1>0$ small enough such that $U_1$ retracts to $\fN_2$ (as in Lemma~\ref{lem-shape2-3}). By Proposition~\ref{prop-rankkn} there exists $\ell_2>2$ such that the image of $\iota_{\ell_2}^2$ has rank 3. Set $U_2=U(\fN_{\ell_2},\e_2)$ for $\e_2>0$ small enough such that $U_2$ retracts to $\fN_{\ell_2}$. Recursively define $\ell_k$ such that the image of $\iota_{\ell_k}^{\ell_{k-1}}$ has rank 3 and set $U_k=U(\fN_{\ell_k},\e_k)$ for $\e_k>0$ small enough such that $U_k$ retracts to $\fN_{\ell_k}$.

\medskip

The following result finishes the proof of Theorem~\ref{thm-stableshape}.

\begin{prop}\label{prop-notbound3}
There is no shape approximation $\fV = \{V_k \mid k=1,2,\ldots\}$ of $\fM$ such that the rank of the homology groups $H_1(V_k,\mZ)$ is 3 for any $k$.
\end{prop}

\proof
Let $\fV = \{V_k \mid k=1,2,\ldots\}$   be a given shape approximation of $\fM$. 
Take $k>0$, then there exists $\e>0$ depending on $k$, such that $C_\mK^+(\fM_0^1,\e)\subset V_k$. The first homology group of $C_\mK^+(\fM_0^1,\e)$ has at least 3 generators, one corresponding   to the Reeb cylinder, and analogous to the class $[R]$ introduced above. The other two generators  are associated to paths which travel from the basepoint on the Reeb cylinder,  through an entry region, then back down one of the two level 1 propellers  back to a neighborhood of the Reeb cylinder. These are analogous to the classes $[b_1]$ and $[b_2]$ introduced before.   For $k$ sufficiently large, we can assume that these three generators are also generators of $H_1(V_k,\mZ)$. Thus, to prove the lemma, it suffices to show that for $k$ sufficiently large,  there is at least one more generator of $H_1(V_k,\mZ)$. As before, let $G_k \subset H_1(V_k ; \mZ)$ be the group generated by the three elements above.

Recall from Section~\ref{sec-zippered} that  the connected components of  the intersection    $\fMR = \fM \cap \bRt$  form a (singular) lamination, and the intersection $\fC' =  \fM \cap \cT$ with the transversal  $\cT \subset \bRt$ (introduced in \eqref{eq-cantorset1}) is a Cantor set. Each $V_k$ is an open set containing $\fM$,  thus the set $V_k\cap \bRt$ is a neighborhood of $\fMR$. It follows    that  for $k$ sufficiently large, the number of connected components of  $V_k\cap \bRt$ must tend to infinity.  Moreover, given  distinct points $x,y \in \fC'$, it follows that for $k$ sufficiently large, the points must be contained in distinct connected components of $V_k\cap \bRt$. Also, let $\g_x \subset \fMR$ be the arc-component containing $x$, and likewise let $\g_y \subset \fMR$ be the arc-component containing $y$. Then for a possibly larger value of $k$, the compact arcs $\g_x$ and $\g_y$ are contained in distinct connected components of $V_k \cap \bRt$.

 For each $k$,    let $V_k^0$ denote the connected component of $V_k\cap \bRt$ containing the trace  $\tau(\cR')\cap \bRt$ of the Reeb cylinder. Then   for $k$ sufficiently large, $V_k^0$ does not contain all the level one points $p_0(i_0;i_1,\ell_1)$ in $\bRt$. Even more,  let $n_k$ be the largest integer   such that      for all $a \leq \ell_1\leq n_k$, the point $p_0(i_0;i_1,\ell_1)$ and the arc-component of $\fMR$ it defines are not contained in $V_k^0$, for any pair $(i_0,i_1)$.

Recall that $p_0(i_0;1,\ell)$ denotes the lower endpoint of the curve $\g_0(\ell) \subset \fMR$ if $i_0=1$, and of $\lambda_0(\ell) \subset \fMR$ if $i_0=2$; and $p_0(i_0;2,\ell)$ is the upper endpoint of the curve $\g_0(\ell)$ if $i_0=1$, and of $\lambda_0(\ell)$ if $i_0=2$. 

Let $k$ be sufficiently large so that $n_k > 1$. We set $i_1 = 1$ in the following. For each   level 2 point $p_0(1;1,\ell_1;1,n_k)$,   let 
 $\ds  \g_0(1,\ell_1;n_k)   \subset  \fMR$ be the curve with this as lower endpoint. 
Let $V_k^1$ be the connected component of $V_k\cap \bRt$ containing the point $p_0(1;1,n_k)$, and hence the curve $\g_0(n_k)$. Note that by our choices, $V_k^0 \cap V_k^1 = \emptyset$. 
Then  for $\ell_1$ sufficiently large,   we have  $\ds  \g_0(1,\ell_2;n_k)  \subset V_k^1$.
  
Using the notation for paths as above, 
define the loop $\sigma(\ell_1)$ as 
\begin{eqnarray*}
p_0(1;1,n_k) & \xrightarrow{\text{tangent($<3$)}}  & p_0(1;1,\ell_1;1,n_k)\\
& \xrightarrow{V_k\cap \bRt} & p_0(1;1,n_k).
\end{eqnarray*}
Moreover, this path can be chosen so that it is disjoint for the set $V_k^0$. In terms of the tree $\TP$, we are choosing a path in a level 2 branch that contains a vertex in   $V_k^1 \cap \TP$ but no vertices in $V_k^0 \cap \TP$.

Let $[\sigma(\ell_1)] \in H_1(V_k,\mZ)$ be the homology class of this loop. We claim that it is not contained in  $G_k$, and thus  the rank of $H_1(V_k,\mZ)$ is greater than three as needed. Suppose that $\sigma(\ell_1)$ is homologous to a class in $G_k$, then there exists a connected singular surface $B \subset V_k$ with boundary, such that $\sigma(\ell_1)$ is one boundary component, and the other boundary component represents a class in $G_k$. As $B$ is connected, this implies that $V_k^0 \cap V_k^1 \ne \emptyset$, which is a contradiction.
\endproof

Note that the above proof of Proposition~\ref{prop-notbound3} is a particular case of using the Mayer-Vietoris Theorem for the decomposition of $V_k$ as the union of the sets $V_k - (V_k \cap \bRt)$ and the set $V_k'$ which is a small open neighborhood of $V_k \cap \bRt$ in $V_k$.

We conclude this section with a proof   that the Mittag-Leffler Condition in Theorem~\ref{thm-MLhomology} holds for the first homology groups.

\proof[Proof of Theorem~\ref{thm-MLhomology}]
Let $\fU = \{U_{\ell}\}$ for $\Sigma$ be a shape approximation  of $\fM$. We must show that     for any $\ell\geq 1$ there exists $p>\ell$ such that for any $q\geq p$ 
\begin{equation}\label{eq-MLhomology2}
Image\{H_1(U_p;\mZ)\to H_1(U_\ell;\mZ)\} = Image \{H_1(U_q;\mZ)\to H_1(U_\ell;\mZ)\}.
\end{equation}

Fix $\ell > 0$, and we may assume without loss of generality that $U_{\ell} \subset \fN_0$. Choose $k$ sufficiently large so that $\fN_k \subset U_{\ell}$.  Then by  Proposition~\ref{prop-rankkn}, there exists $m > k$ so that the image of the map $\iota_m^k:H_1(\fN_m;\mZ)\to H_1(\fN_k;\mZ)$ has rank 3. Let $\e>0$ so that $U(\fN_m,\e)$ retracts to $\fN_m$. Then choose $p$ sufficiently large so that $U_p \subset  U(\fN_m,\e)$, then for all $q \geq p$ we   have $U_q \subset  U_p \subset U(\fN_m,\e)$. 

For any such $q \geq p$, there exists some $v$ such that $\fN_v \subset U_q$. Then consider the sequence of homology groups, with  maps induced by the inclusions, 
 \begin{equation}\label{eq-MLproof}
H_1(\fN_v ; \mZ) \longrightarrow  H_1(U_q ; \mZ) \longrightarrow H_1(U(\fN_m,\e) ; \mZ) \cong  H_1(\fN_m ; \mZ)  \xrightarrow{\iota_m^k}  H_1(\fN_k; \mZ) \longrightarrow  H_1(U_{\ell} ; \mZ) 
\end{equation}
The proof of Proposition~\ref{prop-rankkn} shows that the images of $H_1(\fN_v ; \mZ)$ and $H_1(\fN_m ; \mZ)$ 
under these maps are both equal to the subgroup $G_k \subset H_1(\fN_k; \mZ)$, so it follows that the image of $H_1(U_q ; \mZ) \to H_1(U_{\ell} ; \mZ)$ equals the image of $G_k$ in $H_1(U_{\ell} ; \mZ)$. Thus, the image is independent of the choice of $q\geq p$, as was to be shown. 
\endproof
     
\vfill
\eject


\end{document}